\documentclass[a4paper, 10pt, openany]{book}
\usepackage{graphicx,graphics}
\usepackage{calc}
\usepackage[cmtip,arrow,all,2cell]{xy}
\usepackage{pb-diagram,pb-xy}
\usepackage{amsmath,amsthm,amssymb,tabularx}
\usepackage[mathscr]{eucal}
\usepackage{hyperref}

\numberwithin{equation}{section}

\renewcommand{\thesection}{\S\,\arabic{section}}

\theoremstyle{plain}

\newtheorem{tm}{Theorem}[section]

\newtheorem{prop}[tm]{Proposition}
\newtheorem{lm}[tm]{Lemma}
\newtheorem{cor}[tm]{Corollary}

\theoremstyle{definition}

\newtheorem{ex}{Example}[section]

\newtheorem{rem}[ex]{Remark}

\newcommand{\btm}{\begin{tm}}
\newcommand{\etm}{\end{tm}}
\newcommand{\blm}{\begin{lm}}
\newcommand{\elm}{\end{lm}}
\newcommand{\bprop}{\begin{prop}}
\newcommand{\eprop}{\end{prop}}
\newcommand{\bcor}{\begin{cor}}
\newcommand{\ecor}{\end{cor}}
\newcommand{\bex}{\begin{ex}}
\newcommand{\eex}{\end{ex}}
\newcommand{\brem}{\begin{rem}}
\newcommand{\erem}{\end{rem}}

\newcommand{\bpr}{\begin{proof}}
\newcommand{\epr}{\end{proof}}
\newcommand{\beq}{\begin{equation}}
\newcommand{\eeq}{\end{equation}}
\newcommand{\bit}{\begin{itemize}}
\newcommand{\eit}{\end{itemize}}

\newcommand{\norm}[1]{\left\Vert#1\right\Vert}
\newcommand{\abs}[1]{\left\vert#1\right\vert}
\newcommand{\set}[1]{\left\{#1\right\}}

\allowdisplaybreaks

  \newenvironment{list3}[1]
 {
 \begin{list}{}{
 \topsep=5pt %пробел сверху перед первой записью в списке
 \parsep=0pt  %пробел между абзацами в одной записи
 \itemsep=3pt %пробел между записями (дополнительно к \parsep)
 \itemindent=0pt %пробел слева от метки
 \labelwidth=20pt %ширина метки в записи
 \leftmargin=20pt %пробел слева перед первой записью в списке
 \rightmargin=20pt %пробел справа перед первой записью в списке
 \labelsep=10pt %пробел между меткой и записью
 }{#1}
 \end{list}
 }

\def \le {\leqslant}
\def \ge {\geqslant}

\def\e{\varepsilon}
\def\R{{\mathbb{R}}}
\def\C{{\Bbb C}}

\def \N {\mathbb{N}}

\def \Z {\mathbb{Z}}

\def\br{\operatorname{\sf br}}
\def\Ste{\operatorname{\tt Ste}}
\def\SteAlg{\operatorname{\tt SteAlg}}
\def\InvSteAlg{\operatorname{\tt InvSteAlg}}
\def \d {\operatorname{\sf{d}}}

\def\Int{\operatorname{\sf{Int}}}

\def\Spec{\operatorname{\sf{Spec}}}
\def\Sec{\operatorname{\sf{Sec}}}

\def\Hom{\operatorname{\sf Hom}}
\def\Ob{\operatorname{\sf Ob}}
\def\op{\operatorname{\sf op}}
\def\id{\operatorname{\sf id}}
\def\Mor{\operatorname{\sf Mor}}
\def\card{\operatorname{\sf card}}
\def\sp{\operatorname{\sf span}}
\def\Sp{\operatorname{\sf Span}}

\def\supp{\operatorname{\sf supp}}
\def\csp{\operatorname{\overline{\sf span}}}

\def\d{\operatorname{d}}

\def\tr{\operatorname{\sf tr}}
\def\Mono{\operatorname{\sf Mono}}
\def\Epi{\operatorname{\sf Epi}}

\def\Trig{\operatorname{\sf Trig}}
\def\ord{\operatorname{\sf ord}}

\def\SMono{\operatorname{\sf SMono}}
\def\SEpi{\operatorname{\sf SEpi}}
\def\DEpi{\operatorname{\sf DEpi}}
\def\DBim{\operatorname{\sf DBim}}

\def\DiffMor{\operatorname{\sf DiffMor}}

\def\Iso{\operatorname{\sf Iso}}
\def\Bim{\operatorname{\sf Bim}}

\def\Bind{\operatorname{\sf Bind}}

\def\env{\operatorname{\sf env}}
\def\rf{\operatorname{\sf ref}}
\def\Env{\operatorname{\sf Env}}
\def\Rf{\operatorname{\sf Ref}}

\def\Real{\operatorname{\sf Re}}

\def\Coim{\operatorname{\sf Coim}}
\def\coim{\operatorname{\sf coim}}
\def\Im{\operatorname{\sf Im}}
\def\im{\operatorname{\sf im}}
\def\red{\operatorname{\sf red}}
\def\ker{\operatorname{\sf ker}}
\def\Ker{\operatorname{\sf Ker}}
\def\coker{\operatorname{\sf coker}}
\def\Coker{\operatorname{\sf Coker}}
\def\Dom{\operatorname{\sf Dom}}
\def\Ran{\operatorname{\sf Ran}}

\def\Quot{\operatorname{\sf Quot}}

\def\leftlim{\mathop{\varprojlim}\limits}
\def\rightlim{\mathop{\varinjlim}\limits}

\def \N {\mathbb{N}}
\def \Z {\mathbb{Z}}
\def \R {\mathbb{R}}
\def \T {\mathbb{T}}
\def \C {\mathbb{C}}
\def\Diff{\operatorname{\sf Diff}}
\def\Der{\operatorname{\sf Der}}

\def\id{\operatorname{\sf id}}
\def\Int{\operatorname{\sf Int}}
\def\d{\operatorname{\sf d}}

\def\supp{\operatorname{\sf supp}}

\def\exp{\operatorname{\sf exp}}
\def\card{\operatorname{\sf card}}

\def\d{\operatorname{\sf{d}}}

\def \e {\varepsilon}
\def \ph {\varphi}

\def\Jet{\operatorname{\sf{Jet}}}
\def\jet{\operatorname{\sf{jet}}}

\def\boxedast{\raise.8pt\hbox{\rlap{$\mkern2.5mu *$}}\square}

\def\subarr{\subset\kern-7pt\raise-2.4pt\hbox{\resizebox{5pt}{3.5pt}{$\to$}}\kern2pt}
\def\suparr{\reflectbox{$\subarr$}}

\def\osubarr{\subset\kern-6pt\raise-2.4pt\hbox{\resizebox{5pt}{3.5pt}{$\to$}}\kern-3pt\raise3.4pt\hbox{\resizebox{3pt}{3pt}{$\circ$}}\kern2pt}

\def\quarr{\gets\kern-3.7pt\raise0pt\hbox{\resizebox{2pt}{2.5pt}{$\backslash$}}\kern2pt}

\def\qparr{\kern-1pt\reflectbox{$\quarr$}\kern1pt}

\def\oquarr{\gets\kern-3.7pt\raise0pt\hbox{\resizebox{2pt}{2.5pt}{$\backslash$}}\kern-1.6pt\raise-3.5pt\hbox{\resizebox{3pt}{3pt}{$\circ$}}\kern5pt}

\def\ph{\varphi}
\def\l{\left(}
\def\r{\right)}

\makeindex

\addtocounter{section}{-1}

\begin{document}

\title{CONTINUOUS AND SMOOTH ENVELOPES OF TOPOLOGICAL ALGEBRAS}

\author{S.S.Akbarov}

\maketitle

\section{Geometries as categorical constructions}

\paragraph{Observation tools and visible image.}
The part of mathematics that studies the constructions on the objects called {\it manifolds} (or {\it varieties}) can be divided into four domains:
\bit{

\item[--] {\it algebraic geometry} (which can be perceived as a science studying the structure generated on algebraic varieties $M$ by the algebra ${\mathcal P}(M)$ of polynomials),

\item[--] {\it complex geometry} (where the algebras ${\mathcal O}(M)$ of holomorphic functions on complex manifolds $M$ play the same role),

\item[--] {\it differential geometry} (with the algebras ${\mathcal E}(M)$ of smooth functions on smooth manifolds $M$),

\item[--] {\it topology} (where the algebras ${\mathcal C}(M)$ of continuous functions on topological spaces $M$ can be considered as structure algebras).

}\eit
The obvious parallels between these disciplines inspire an idea that there must exist a universal scheme inside mathematics that explains these similarities and allows to discuss the differences in formal terms.

Such a scheme indeed exists, and the idea leading to it is borrowed from physics and can be expressed in the formula:

    \centerline{\it the visible picture depends on the observation tools.}
\vglue9pt

\noindent
As an example, in astronomy the visible image of an object under study that appears in the observer's mind when he uses optical telescope differs from what he sees with his own eyes, or when he uses radio telescope, or X-ray telescope, etc.

It turns out that with certain understanding of the terms ``observation tools'' and ``visible image'' in mathematics one can form a general view at least on the last three disciplines in this list, -- complex geometry, differential geometry and topology, -- and they will be reflections of one common reality, the pictures that appear as results of the choice of a concrete set of tools. This leads to an intriguing picture, where it becomes possible to compare these ``geometries as disciplines'', to find common features, differences, generalizations, new examples, and so on.

It is convenient to assume that the common reality which these geometries reflect is some, enough wide, category of topological associative algebras, for instance, the category $\SteAlg$ of stereotype algebras\footnote{See definition of stereotype algebra below on page \pageref{DEF:ster-alg}.} (possibly, with some supplementary structures, like involution). Then for the formalization of the scheme of {\it observation}, which we discuss here, the following two agreement are sufficient:
\bit{
\item[1)] by the observation tools one means morphisms of a given class $\varPhi$ in this category,

\item[2)] the object in study (an algebra) $A$ and its visible image (another algebra) $E$ are connected through a natural morphism $A\to E$ (like an original and a photo), and the class $\varOmega$ of such morphisms, called ``class of representations'' is given from the very beginning.
 }\eit
Under these assumptions the ``visible image'' $E$ of an object $A$ can be interpreted as its {\it envelope}\footnote{See definition of envelope on page \pageref{DEF:diagr-obolochka}.} $\Env_\varPhi^\varOmega A$, generated by the class of tools (morphisms) $\varPhi$ in the given class of representations (morphisms) $\varOmega$. Every concrete choice of classes $\varPhi$ and $\varOmega$ give birth some ``projection'' of functional analysis (which is understood here as a theory of topological algebras) into geometry (understood as a theory of functional algebras in ``generalized sense'').

Historically one of the first examples of an envelope of a topological algebra was the Arens-Michael envelope, introduced by J.~L.~Taylor in \cite{Taylor-1}. This construction was studied in detail by A.~Yu.~Pirkovaskii in his researches on ``noncommutative complex geometry'' \cite{Pirkovskii}. In the work \cite{Akbarov-stein-groups} it was applied by the author to a generalization of Pontryagin duality to a class of (not necessarily commutative) complex Lie groups. Up to the further terminological corrections (see \cite{Akbarov-env}), these investigations can be considered as an application of this scheme of observation with a result as a {\it categorical construction of complex geometry}.

Later Yu.~N.~Kuznetsova obtained analogous results in \cite{Kuznetsova}, where she constructed a variant of duality theory in topology with homomorphisms into $C^*$-algebras as observation tools. Again, up to later correction in \cite{Akbarov-env} and in this paper, the results of \cite{Kuznetsova} can be interpreted as a way for {\it categorical construction of topology}.

In the present paper we suggest an analogous way for {\it categorical construction of differential geometry}. We describe here a construction of {\it smooth envelope} of stereotype algebra, we study its properties and build a generalization of Pontryagin duality for some class of real Lie groups. Pictorially the results of this paper and the papers we mentioned here can be presented in the following table:

\bigskip

{\small

\begin{tabular}{|c|c|c|c|}
  \hline
  &  &  & \\
   Discipline\label{TABLE:comparison} & {\bf Complex} & {\bf Differential} & {\bf Topology} \\
  & {\bf geometry} & {\bf geometry} & \\
  &  &  & \\
  \hline
  &  &  & \\
 Key  & Algebra ${\mathcal O}(M)$ & Algebra ${\mathcal E}(M)$ & Algebra ${\mathcal C}(M)$ \\
  example & of holomorphic functions & of smooth functions & of continuous functions \\
  of visible & on a complex & on a smooth  & on a topological \\
  image & manifold $M$ & manifold $M$ & space $M$ \\
  &  &  & \\
  \hline
  &  &  & \\
  Observation & Homomorphisms & Differential & Involutive \\
  tools & into Banach & involutive  & homomorphisms \\
 $\varPhi$ & algebras & homomorphisms & into $C^*$-algebras \\
  &  & into $C^*$-algebras & \\
  &  & with joined & \\
  &  & self-adjoint & \\
    &  & nilpotent elements & \\
  &  &  & \\
  \hline
  &  &  & \\
  Class of & Dense & Dense & Dense \\
  representations & epimorphisms & epimorphisms & epimorphisms \\
  $\varOmega$ &  &  & \\
  &  &  & \\
  \hline
  &  &  & \\
 Visible & Holomorphic envelope & Smooth envelope & Continuous envelope \\
 image of an object $A$ & $\Env_{\mathcal O}A$  & $\Env_{\mathcal E}A$ & $\Env_{\mathcal C}A$ \\
  &  &  & \\
  \hline
  &  &  & \\
 Key  & $\Env_{\mathcal O}A={\mathcal O}(M)$  & $\Env_{\mathcal E}A={\mathcal E}(M)$  & $\Env_{\mathcal C}A={\mathcal C}(M)$ \\
 representation & for a subalgebra & for a subalgebra & for a subalgebra \\
  & $A\subseteq{\mathcal O}(M)$ & $A\subseteq{\mathcal E}(M)$ & $A\subseteq{\mathcal C}(M)$ \\
  &  &  & \\
  \hline
  &  &  & \\
 \xymatrix @R=0pt @C=3pt {\\ \text{Reflexivity} \\ \text{diagram}} &
 \xymatrix @R=3pt @C=3pt
 {
 {\mathcal O}^\star(G) & \ar@{|->}[rr]^{\Env_{\mathcal O}} & & &
 {\mathcal O}_{\exp}^\star(G)
 \\
 &&&& \ar@{|->}[dd]^{\star} \\
 &&&& \\
 \ar@{|->}[uu]^{\star} &&&&  \\
  {\mathcal O}(G)
 & & &
 \ar@{|->}[ll]^{\Env_{\mathcal O}}
 & {\mathcal O}_{\exp}(G)
    }
    &
 \xymatrix @R=3pt @C=3pt
 {
 {\mathcal E}^\star(G) & \ar@{|->}[rr]^{\Env_{\mathcal E}} & & &
 \Env_{\mathcal E}{\mathcal E}^\star(G)
 \\
 &&&& \ar@{|->}[dd]^{\star} \\
 &&&& \\
 \ar@{|->}[uu]^{\star} &&&&  \\
  {\mathcal E}(G)
 & & &
 \ar@{|->}[ll]^{\Env_{\mathcal E}}
 & {\mathcal K}_\infty(G)
    }
    &
 \xymatrix @R=3pt @C=3pt
 {
 {\mathcal C}^\star(G) & \ar@{|->}[rr]^{\Env_{\mathcal C}} & & &
 \Env_{\mathcal C}{\mathcal C}^\star(G)
 \\
 &&&& \ar@{|->}[dd]^{\star} \\
 &&&& \\
 \ar@{|->}[uu]^{\star} &&&&  \\
  {\mathcal C}(G)
 & & &
 \ar@{|->}[ll]^{\Env_{\mathcal C}}
 & {\mathcal K}(G)
    }
  \\
  &  &  & \\
  \hline
   &  &  & \\
 \xymatrix @R=0pt @C=3pt { \text{Reflexivity} \\ \text{diagram}\\ \text{for}\\ \text{commutative}\\ \text{groups}} &
 \xymatrix @R=3pt @C=3pt
 {
 {\mathcal O}^\star(G) & \ar@{|->}[rr]^{\mathcal F} & & &
 {\mathcal O}(\widehat{G})
 \\
 &&&& \ar@{|->}[dd]^{\star} \\
 &&&& \\
 \ar@{|->}[uu]^{\star} &&&&  \\
  {\mathcal O}(G)
 & & &
 \ar@{|->}[ll]^{\mathcal F}
 & {\mathcal O}^\star(\widehat{G})
    }
    &
 \xymatrix @R=3pt @C=3pt
 {
 {\mathcal E}^\star(G) & \ar@{|->}[rr]^{\mathcal F} & & &
 {\mathcal E}(\widehat{G})
 \\
 &&&& \ar@{|->}[dd]^{\star} \\
 &&&& \\
 \ar@{|->}[uu]^{\star} &&&&  \\
  {\mathcal E}(G)
 & & &
 \ar@{|->}[ll]^{\mathcal F}
 & {\mathcal E}^\star(\widehat{G})
    }
    &
 \xymatrix @R=3pt @C=3pt
 {
 {\mathcal C}^\star(G) & \ar@{|->}[rr]^{\mathcal F} & & &
 {\mathcal C}(\widehat{G})
 \\
 &&&& \ar@{|->}[dd]^{\star} \\
 &&&& \\
 \ar@{|->}[uu]^{\star} &&&&  \\
  {\mathcal C}(G)
 & & &
 \ar@{|->}[ll]^{\mathcal F}
 & {\mathcal C}^\star(\widehat{G})
    }
  \\
  &  &  & \\
  \hline
\end{tabular}
}

\paragraph{Complex geometry.}

The second column in this table was chronologically first, so it is logical to start the explanation with it.

{\it Complex geometry} studies complex manifolds with supplementary structures like Hermitian metrics, or connextions, or curvature, etc. \cite{Griffiths-Harris}. Usually a complex manifold is defined by its sheaf of holomorphic functions, and for a mathematician with functional-analytic mentality this construction poses a psychological problem. However among complex manifolds there is a subclass, which does not require the notion of sheaf for its description, its objects are called {\it Stein manifold} \cite{Grauert-Remmert},\cite{Taylor}. For our purposes they are good, since, first, they visually illustrate our idea, and, second, the passage from them to the general case (of a manifold defined by sheaf) in the category theory seems to be easy, since a sheaf itself is a simple categorical construction (however, the formal generalization was not constructed yet).

A Stein manifold $M$ is defined by its algebra ${\mathcal O}(M)$ of holomorphic functions (and respectively, all the supplementary structures, like metrics on $M$, also can be defined as constructions on ${\mathcal O}(M)$). As a corollary, it is logical to consider the algebra ${\mathcal O}(M)$ as a key example of visible image in complex geometry. In this science it has been noticed long ago (see \cite{Pirkovskii}, \cite{Akbarov-stein-groups}, \cite{Akbarov-env}) that ${\mathcal O}(M)$ can be restored from a (sufficiently wide) subalgebra $A\subseteq{\mathcal O}(M)$ with the help of the system of Banach quotient algebras $A/p$ for various continuous submultiplicative seminorms $p$ on $A$ -- then ${\mathcal O}(M)$ is just a projective limit of $A/p$:
\beq\label{O(M)=leftlim_p-A/p}
{\mathcal O}(M)=\leftlim_{p}A/p.
\eeq
Of course, this holds not for any $A\subseteq{\mathcal O}(M)$: if we want \eqref{O(M)=leftlim_p-A/p} to be true, $A$ must be sufficiently close to ${\mathcal O}(M)$. An important example was found by A.~Yu.~Pirkovskii \cite{Pirkovskii}: if  $M$ is an affine algebraic manifold, then ${\mathcal O}(M)$ can be restored in this way by the subalgebra $A={\mathcal P}(M)$ of polynomials on $M$.

Up to the details insufficient in the first approximation, the equality \eqref{O(M)=leftlim_p-A/p} can be understood as the proposition that ${\mathcal O}(M)$ is an envelope of $A$ in the class $\varOmega=\DEpi$ of dense epimorphisms\footnote{Dense epimorphisms are defined below on page \pageref{DEF:DEpi}.} of the category ${\SteAlg}$ of stereotpye algebras with respect to the clas $\varPhi=\Mor({\SteAlg},{\tt BanAlg})$ of homomorphisms into Banach algebras:
\beq\label{Env_(Mor-SteAlg,BanAlg)^DEpi-A=O(M)}
{\mathcal O}(M)=\Env_{\Mor({\SteAlg},{\tt BanAlg})}^{\DEpi}A.
\eeq
The projective limit on the right side of \eqref{O(M)=leftlim_p-A/p} is called the {\it Arens-Michael envelope}, and for the arbitrary stereotype algebras $A$ it doesn't coincide with the obect on the right side of \eqref{Env_(Mor-SteAlg,BanAlg)^DEpi-A=O(M)}, that is why in \cite{Akbarov-env} the author introduced a new term for $\Env_{\Mor({\SteAlg},{\tt BanAlg})}^{\DEpi}A$, the {\it holomorphic envelope}. In the table at the page \pageref{TABLE:comparison} this envelope is denoted by $\Env_{\mathcal O}A$. Up to this notation, the equality \eqref{Env_(Mor-SteAlg,BanAlg)^DEpi-A=O(M)} is the proposition in the sixth cell of the second column of the table. We call it {\it key representation} having in mind that it describes the mechanism of discerning objects as a key example in this science, the algebra ${\mathcal O}(M)$.

The main result in \cite{Akbarov-stein-groups} is the diagram in the next to last cell of the column:
\beq\label{DIAGR:holom-dvoistv}
 \xymatrix @R=.5pc @C=1pc
 {
 {\mathcal O}^\star(G) & \ar@{|->}[rr]^{\Env_{\mathcal O}} & & &
 {\mathcal O}_{\exp}^\star(G)
 \\
 &&&& \ar@{|->}[dd]^{\star} \\
 &&&& \\
 \ar@{|->}[uu]^{\star} &&&&  \\
  {\mathcal O}(G)
 & & &
 \ar@{|->}[ll]^{\Env_{\mathcal O}}
 & {\mathcal O}_{\exp}(G)
    }
\eeq
Here ${\mathcal O}(G)$ is the algebra of holomorphic functions on a complex Lie group $G$, ${\mathcal O}_{\exp}(G)$ the algebra of holomorphic functions of exponential type on $G$, ${\mathcal O}^\star(G)$ and ${\mathcal O}_{\exp}^\star(G)$ the dual convolution algebras of analytic functionals, and $\star$ the operation of passing to the dual stereotype space\footnote{See definition at p.\pageref{DEF:star}.}.

Diagram \eqref{DIAGR:holom-dvoistv} shows that the objects in its corners satisfy some reflexivity conditions. For example, if we denote the composition $\star\circ\Env_{\mathcal O}$ by \ $\widehat{}$\, then ${\mathcal O}^\star(G)$ and ${\mathcal O}_{\exp}(G)$ become reflexive with respect to \ $\widehat{}$\ :
$$
\widehat{\widehat{A}}\cong A.
$$
That is why \eqref{DIAGR:holom-dvoistv} is called {\it reflexivity diagram}. The stated proposition (that on the fourth step the chain leads to the initial object) is proved so far only for complex Lie groups $G$ with algebraic component of identity, but, obviously, it must be true for wider class of groups (which namely is a subject of further investigations).

In the special case when $G$ is commutative the reflexivity diagram has the form
\beq\label{DIAGR:holom-dvoistv-F}
 \xymatrix @R=.5pc @C=1pc
 {
 {\mathcal O}^\star(G) & \ar@{|->}[rr]^{\mathcal F} & & &
 {\mathcal O}(\widehat{G})
 \\
 &&&& \ar@{|->}[dd]^{\star} \\
 &&&& \\
 \ar@{|->}[uu]^{\star} &&&&  \\
  {\mathcal O}(G)
 & & &
 \ar@{|->}[ll]^{\mathcal F}
 & {\mathcal O}^\star(\widehat{G})
    }
\eeq
and here $\widehat{G}$ is the group of complex characters, i.e. homomorphisms $\chi:G\to\C^\times$ into the multiplicative group $\C^\times$ of non-zero complex numbers, and ${\mathcal F}$ is the Fourier transform. That is the statement of the last cell of the column.

Diagrams \eqref{DIAGR:holom-dvoistv} and \eqref{DIAGR:holom-dvoistv-F} together define a generalization of Pontryagin duality from commutative complex Lie groups to (non necessarily commutative) Lie groups with algebraic component of identity, and the author sees here a justification of the theory built in \cite{Akbarov-stein-groups}.

\paragraph{Topology.}
After the explanation about the second column the idea of the third and the fourth ones is more or less clear. Following chronology, we then have to explain the fourth column, which presents applications of the described above scheme to topology. The similarity with the second column is that we describe here some already published facts, namely the results of the Yu.~N.~Kuznetsova paper \cite{Kuznetsova}. But the difference is that we retell here the material of \cite{Kuznetsova} in detail, with formulations and proofs. The aim is to repair some inaccuracies and mistakes of \cite{Kuznetsova}. We put this material into $\S 3$ (and partly into $\S 1$) (however, the content of $\S 3$ and $\S 1$ is not exhausted by the results of \cite{Kuznetsova}).

Briefly this theme can be presented as follows.

The class of objects that {\it general topology} studies  -- topological spaces -- is so wide, that for constructing geometries on them one have to pick out subclasses consisting of spaces with given properties of separability, local connectedness, local Euclidean property, etc. Like in complex geometry, the key example of visible image here is a  functional algebra, in this case the algebra ${\mathcal C}(M)$ of continuous functions on a topological space $M$.

If we claim that ${\mathcal C}(M)$ totally defines $M$, then from the very beginning we have to choose spaces $M$ on which the functions from ${\mathcal C}(M)$ separate points. As a corollary, if the scalar field is $\C$ (or $\R$), then $M$ must belong to a standard class satisfying the conditions of functional separability, for example the class of Tikhonov spaces. We consider a narrower class, the class of paracompact locally compact spaces, since the functional algebra ${\mathcal C}(M)$ is relatively simple for them \cite[8.1]{Akbarov} (it is clear, however, that this class can be extended in different directions).

The substantial difference with the case of complex manifold is that, first, the functional algebra ${\mathcal C}(M)$ has natural involution (that allows to restore $M$), and, second, if we want ${\mathcal C}(M)$ to be a visible image, we have to use other observation tools, in this case the involutive homomorphisms into $C^*$-algebras. We denote the corresponding envelope by $\Env_{\mathcal C}$ and call it {\it continuous}. The key representation is described by Theorem \ref{C-obolochka-podalgebry-v-C(M)} below:
$$
\Env_{\mathcal C}A={\mathcal C}(M)
$$
for each unital involutive subalgebra $A$ in ${\mathcal C}(M)$ such that the dual map of spectra $\Spec(A)\gets M$ is a covering (see definition at p.\pageref{DEF:nalozhenie}). The reflexivity diagram is presented below with the label \eqref{chetyrehugolnik-C-C*}, and its commutative analog with the label \eqref{chetyrehugolnik-C-C*-F}. Curiously, so far these diagrams are proved only for a quite narrow class of Moore groups.

\paragraph{Differential geometry.}

This science studies smooth manifolds with supplementary structures like those mentioned in the case of complex manifolds: Riemennian metrics, connectedness, curvature \cite{Alexeevskii-Vinogradov-Lychagin}. They all can be understood as constructions on the algebra ${\mathcal E}(M)$ of smooth functions on a manifold $M$, so the key example of visible image here must be, of course, ${\mathcal E}(M)$.

Up to the last time it was unclear, what are the observation tools in this science. It is clear, for example, that these can't be  arbitrary involutive homomorphisms into $C^*$-algebras (like in topology), since this assumption immediately leads to the continuous envelope $\Env_{\mathcal C}$, which turns ${\mathcal E}(M)$ into ${\mathcal C}(M)$:
$$
\Env_{\mathcal C}{\mathcal E}(M)={\mathcal C}(M)
$$
(this is not good, since the observation tools must not distort the key example of visible image).

The contribution of this paper is that we suggest here the observation tools that allow to recognize the algebras of smooth functions ${\mathcal E}(M)$, and thus suitable for differential geometry. As this is declared in the table at p.\pageref{TABLE:comparison}, these instruments are involutive differential homomorphisms into $C^*$-algebras with joined self-adjoint nilpotent elements (see definitions of $\S 4$). We denote the corresponding envelope by  $\Env_{\mathcal E}$ and call it {\it smooth} (definition at p.\pageref{DEF:C^infty-obolochka}). The key representation is described in Theorem \ref{C^infty-obolochka-podalgebry-v-C^infty}:
$$
\Env_{\mathcal E} A={\mathcal E}(M)
$$
for each involutive unital subalgebra $A\subseteq{\mathcal E}(M)$ such that the dual map of spectra $\Spec(A)\gets M$ is a covering, and in each point $t\in M$ the map of tangent spaces $T_s[A]\gets T_s(M)$ is an isomorphism.
The reflexivity diagrams for the general case and for the case of commutative group have labels  \eqref{chetyrehugolnik-E-E*} and \eqref{chetyrehugolnik-E-E*-F}. The first one is proved only for the class of the groups $C\times K$, where $C$ is an Abelian compactly generated Lie group, and $K$ a compact Lie group.

\paragraph{What is ``geometry as a discipline''?}

The famous F.~Klein Erlangen program, the B.~Riemann field approach, the E.~Cartan construction (see \cite{Alexeevskii-Vinogradov-Lychagin}, \cite{Sharpe}), define the term {\it geometry} in mathematics as a theory obtained by some given schemes inside each of those four disciplines that we mentioned at the beginning. They don't explain the meaning of this word in the very names of those big disciplies  -- ``complex geometry'', ``differential geometry''...

This paper together with the earlier mentioned papers \cite{Akbarov-stein-groups}, \cite{Kuznetsova}, \cite{Akbarov-env}, suggests a way to explain this. As we told before, the scheme we describe gives a purely categorical way to construct ``geometry as a discipline'': {\it each choice of the class of observation tools $\varPhi$ and the class of representations $\varOmega$ in the category of algebras\footnote{We consider stereotype algebras, but nothing prevents to consider other categories of algebras.} generates a category of algebras, which can be naturally interpreted as a geometric discipline, defined by the classes $\varPhi$ and $\varOmega$.}

This must be interesting as an attempt to look at mathematics ``from above'', but this view promises also a practical use. Here are the justifying reasons:

1) As we see in the table at p.\pageref{TABLE:comparison}, each geometry has its own duality theory that generalizes the Pontryagin duality for Abelian groups. The classes of groups covered by these new duality theories are not yet exhaustively described, and moreover in the last two examples -- in differential geometry and in topology -- it is already clear that in this choice of observation tools these classes can't be that wide as we would like them to be (see discussion at p.\pageref{DEF:gruppy-razlich-C^*-algebrami} and p.\pageref{DEF:gruppy-razlich-C^*-algebrami-s-nilp}). Nevertheless, there are some reasons to think that in comparison with the existing general theories these ``geometric duality theories'', though the loss in generality, win in utility, since they describe properties of objects inside these geometric disciplines, without the necessity to go out of their limits. As a comparison we can mention the duality theory constructed by L.~I.~Vainerman, G.~I.~Kac,  M.~Enock and J.-M.~Schwartz (see \cite{Enock-Schwartz}). In this theory the algebras representing a group $G$ are the algebra $L^\infty(G)$ of sufficiently bounded functions on $G$, or the von Neumann algebra ${\mathcal L}(G)$ generated by the left regular representation of $G$ (with the supplementary structures -- analogues of the comultiplication and the antipode, and also a weight). This is a general theory for all locally compact groups, but if we compare it with the duality in complex geometry presented by the reflexivity diagram \eqref{DIAGR:holom-dvoistv} (where the class of groups is much narrower), we can reasonably hope that the algebras used in \eqref{DIAGR:holom-dvoistv} -- ${\mathcal O}(G)$, ${\mathcal O}^\star(G)$, ${\mathcal O}_{\exp}(G)$, ${\mathcal O}_{\exp}^\star(G)$ -- will be more likely relevant in complex geometry than the algebras $L^\infty(G)$ and ${\mathcal L}(G)$, which need supplementary linkage to the theory.

2) Another argument is the following: {\it the described approach binds more tightly functional analysis to algebra and geometry}, since the arising objects are ``much more categorical'', than the classical ones. An illustration: the group algebras ${\mathcal O}^\star(G)$, ${\mathcal E}^\star(G)$, ${\mathcal C}^\star(G)$, mentioned in the next to last line of the table at p.\pageref{TABLE:comparison} are initial objects in the category of holomorphic, smooth and continuous representations of the group $G$ (i.e. they have the properties described by diagrams of the type  \eqref{C*(G)=grup-alg-A} below, see \cite[Theorem 10.12]{Akbarov}; actually that is why ${\mathcal O}^\star(G)$, ${\mathcal E}^\star(G)$, ${\mathcal C}^\star(G)$ have a right to be called group algebras). As a comparison, the algebras $L^\infty(G)$ and ${\mathcal L}(G)$ have similar properties only for discrete groups $G$. Another illustration:  ${\mathcal O}^\star(G)$, ${\mathcal E}^\star(G)$, ${\mathcal C}^\star(G)$ are Hopf algebras. We wrote about this in \cite{Akbarov-stein-groups}: among all the non-stereotype generalizations of Pontryagin duality only in the duality for finite groups the representing algebras are Hopf algebras, in all other theories this is not so.

3) An obvious use of this categorical view on geometry is, no more, no less, than the {\it clarification of the subject of noncommutative geometry} \cite{Connes}. Since in each case the visible images are not nesessarily commutative algebras, the constructions we suggest can be considered as formal definitions of ``noncommutative geometries'' (noncommu\-tative complex geometry, noncommutative differential geometry, noncommutative topology, etc. -- see a little survey in A.~Yu.~Pirkovskii's paper \cite{Pirkovskii}).

\paragraph{Prospects.}\label{PERSPEKTIVY}

It is evident to the author that the facts found so far in this field are just some accidentally caught pictures from a superficial layer, and they hide a deep and substantial science which deserves further investigations. As an advertisement we think it will be useful to sketch some subjects that immediately occur to mind.

1. First of all it must be clear that the classes of groups for which the described variants of duality theories are true, can be extended. For instance, it is likely that in the complex geometry the class of groups with algebraic component of identity can be extended to the class of groups with linear component. In differential geometry this class is very likely the class of Lie-Moore groups.

2. The fact that the obtained classes of groups are nevertheless quite narrow, means evidently that the observation tools found by now are too rough. There must be finer tools which allow to construct duality theory for wider classes of groups. In particular, in topology, -- for the class of all locally compact groups. One of the possible ways for this is the substitution of $C^*$-algebras by von Neumann algebras.

3. Another interesting direction is the substitution of the scalar  fields $\C$ and $\R$ by some other fields. One can, for example consider the field of $p$-adic numbers (of course the class of ``visible groups'' will change very much in this way). Another intriguing alternative is the residue field $\Z_2=\Z/2\Z$. It can be endowed by a topology in such a way that for an arbitrary topological space $M$ the topological separability in $M$  becomes equivalent to the separability by continuous functions $f:M\to\Z_2$. This topology on $\Z_2$ is the topology of the connected colon, where the closed sets in $\Z_2=\{0;1\}$ are exactly $\varnothing$, $\{0\}$ and $\{0;1\}$ (and their complements are open sets). Certainly, the substitution of $\C$ by $\Z_2$ will dramatically change the class of topological spaces on which the geometry can be constructed, but for this one apparently have to build as a fundament a theory analogous to the stereotype one with the field $\Z_2$.

\paragraph{Acknowledgements.}

The author thanks M.~B\"achtold, Y.~Choi, A.~I.~Degtyarev, M.~J.~Dupr\'e, G.~Elencwajg, P.~I.~Katsylo, Yu.~N.~Kuznetsova, B.~McKay, P.~Michor, A.~Yu.~Pirkovskii, V.~L.~Popov, C.~Remling, A.~I.~Shtern, E.~B.~Vinberg and Q.~Yuan for useful advices.

\paragraph{Terminology.}
Everywhere in category theory we use the terminology of textbooks \cite{Bucur-Deleanu}, \cite{Tsalenko-Shulgeifer} and of handbook  \cite{General-algebra}, and as a set-theoretic fundament for the notion of category we choose the Morse-Kelley theory \cite{38}.

The notations $\Mono({\tt K})$, $\Epi({\tt K})$, $\SMono({\tt K})$ and $\SEpi({\tt K})$ mean the classes of monomorphisms, epimorphisms, strong monomorphisms and strong epimorphisms (the last two are defined at p.\pageref{razbienie-kvadrata}) respectively in the category ${\tt K}$.
We say that a category ${\tt K}$ is
 \bit{
\item[---] {\it injectively (projectively) complete}, if each functor $K:{\tt M}\to{\tt K}$ from a small category $\tt M$  (i.e. a category where the class of morphisms is a set) has an injective (projective) limit,

\item[---] {\it complete}, if it is injectively and projectively complete,

 \item[---] {\it linearly complete}\label{DEF:lineino-polnaya-kategoriya}, if any functor from a linearly ordered set to ${\tt K}$ has injective and projective limits.
 }\eit

For any morphism $\ph:X\to Y$ in an arbitrary category the symbols
$\Dom\ph$ and $\Ran\ph$ mean respectively the domain and the range of $\ph$,
i.e. $\Dom\ph=X$ and $\Ran\ph=Y$. If $\tt L$ and $\tt M$ are two classes of objects in $\tt K$, then $\Mor({\tt L},{\tt M})$ means the class of morphisms with domains in $\tt L$ and ranges in $\tt M$:
$$
\Mor({\tt L},{\tt M})=\{\ph\in\Mor({\tt K}):\quad \Dom\ph\in{\tt L}\ \&\ \Ran\ph\in{\tt M}\}.
$$

Let $\varPhi$ be a class of morphisms and ${\tt L}$ a class of objects in a category ${\tt K}$. We say that
\bit{
\item[---]\label{DEF:goes-from} {\it $\varPhi$ goes from ${\tt L}$}, if for any object $X\in{\tt L}$ there is a morphism $\ph\in\varPhi$, going from $X$:
    $$
    \forall X\in{\tt L}\qquad \exists \ph\in\varPhi\qquad \Dom\ph=X;
    $$
in the special case, if $\tt L$ consists of only one object $X$, we say that $\varPhi$ {\it goes from $X$},

\item[---]\label{DEF:goes-to} {\it $\varPhi$ goes to ${\tt L}$}, if for any object $X\in{\tt L}$ there is a morphism $\ph\in\varPhi$, going to $X$:
    $$
    \forall X\in{\tt L}\qquad \exists \ph\in\varPhi\qquad \Ran\ph=X.
    $$
in the special case when $\tt L$ consists of only one object $X$, we say that $\varPhi$ {\it goes to $X$}.
    }\eit

If a topological space $Y$ is imbedded into a topological space $X$ (injectively, but not necessarily in such a way that the topology of $Y$ is inherited from $X$),
and $A$ is a subset in $Y$, then to distinguish the closure of $A$ in $Y$ from its closure in $X$,
we denote the first one by $\overline{A}^Y$, and the second by $\overline{A}^X$\label{DEF:overline^X}.

In the theory of topological vector spaces we follow the textbook by H.~Schaefer \cite{Schaefer}. In particular, we assume that {\it all locally convex spaces (LCS, shortly) are Hausdorff}. By {\it morphism of locally convex spaces} (LCS) we mean an arbitrary linear continuous map $\ph:X\to Y$. The category of all locally convex spaces with these morphisms is denoted by ${\tt LCS}$.

We say that a subset $M$ in a locally convex space $X$ is {\it total (in $X$)}, if its linear span $\Sp M$ is dense in $X$:
$$
\overline{\Sp M}^X=X.
$$

\section{Envelopes and refinements in categories}

Envelopes and refinements in abstract categories were introduced by the author in \cite{Akbarov-env}. We give here main definitions and formulate some facts from \cite{Akbarov-env} that we will need in this text.

\subsection{Nodal decomposition}

\paragraph{Standard classes of monomorphisms and epimorphisms.}
 Recall that a morphism $\ph:X\to Y$ is called
 \bit{
\item[---] a {\it monomorphism}, if any equality
$\ph\circ\alpha=\ph\circ\beta$ implies $\alpha=\beta$;

\item[---] an {\it epimorphism}, if any equality
$\alpha\circ\ph=\beta\circ\ph$ implies $\alpha=\beta$;

\item[---] a {\it bimorphism}, if it is a monomorphism and an epimorphism.
 }\eit

These notions have several important variations, of which two will be useful for us. The first one of these variations is based on the notion of factorisation.

\bit{

\item[$\bullet$]
A {\it factorisation} of a morphism $\ph:X\to Y$ is any its representation as acomposition of an epimorphism and a monomorphism, i.e. any diagram
\beq\label{DEF:faktorizatsiya}
\begin{diagram}
\node{X}\arrow{se,b}{\e}\arrow[2]{e,t}{\ph}\node[2]{Y}\\
\node[2]{M}\arrow{ne,b}{\mu}
\end{diagram}
\eeq
where $\e$ is an epimorphism, and $\mu$ a monomorphism.

\item[$\bullet$]
A monomorphism $\mu:X\to Y$ is said to be {\it immediate}\label{DEF:immediate-mono}, if in any its factorization $\mu=\mu'\circ\e$ the epimorphism $\e$ is automatically an isomorphism. Note that for a monomorphism $\mu$ in any its factorization $\mu=\mu'\circ\e$ the epimorphism $\e$ is automatically a bimorphism. As a corollary, the condition of being immediate monomorphism for $\mu$ is equivalent to the requirement that in any decomposition $\mu=\mu'\circ\e$, where $\e$ is a bimorphism, and $\mu'$ a monomorphism, the morphism $\e$ must be an isomorphism. It is natural to call a monomorphism $\mu'$ in the factorization $\mu=\mu'\circ\e$ a {\it mediator}\label{DEF:posrednik-monomorfizma} of the monomorphism $\mu$, then the epithet ``immediate'' for $\mu$ will mean that there are no non-trivial  mediators for $\mu$ (i.e. mediators, which are not isomorphic to $\mu$ in $\Mono_Y$ -- see below definition \eqref{DEF:varGamma(X)}, here $\varGamma=\Mono$).

\item[$\bullet$]
An epimorphism $\e:X\to Y$ is said to be {\it immediate}\label{DEF:immediate-epi}, if
if $\e$ is an immediate monomorphism in the dual category. In other words, in any factorization $\e=\mu\circ\e'$ the monomorphism $\mu$ must be automatically an isomorphism. Note that for an epimorphism $\e$ in any its factorization $\e=\mu\circ\e'$ the monomorphism $\mu$ is automatically a bimorphism. As a corollary, the condition of being immediate epimorphism for $\e$ is equivalent to the requirement that in any decomposition $\e=\mu\circ\e'$, where  $\mu$ is a bimorphism, and $\e'$ an epimorphism, the morphism $\mu$ must be an isomorphism. It is natural to call an epimorphism $\e'$ in the factorization $\e=\mu\circ\e'$ a {\it mediator}\label{DEF:posrednik-epimorfizma} of the epimorphism $\e$, then the epithet ``immediate'' for  $\e$ will mean that there are no non-trivial mediators for $\e$ (i.e. mediators, which are not isomorphic to $\e$ in $\Epi^X$ -- see below definition \eqref{morphism-v-Epi(X)}, here $\varOmega=\Epi$).

}\eit

The second variation is based on the notion of diagonalization.

\bit{

\item
A pair of morphisms $(\mu,\e)$ is said to be {\it diagonizable}
\cite{General-algebra,Tsalenko-Shulgeifer}, if for all morphisms
$\alpha:\Dom\e\to\Dom\mu$ and $\beta:\Ran\e\to\Ran\mu$ such that $\mu\circ\alpha=\beta\circ\e$
there exists a morphism $\delta:B\to C$, such that the following diagram is commutative:
\beq\label{DIAGR:diagonaliziruemost}
\begin{diagram}
\node{\Dom\e}\arrow{s,l}{\alpha}\arrow{e,t}{\e}\node{\Ran\e}\arrow{s,r}{\beta}\arrow{sw,t,--}{\delta}
\\
\node{\Dom\mu}\arrow{e,b}{\mu}\node{\Ran\mu}
\end{diagram}
\eeq
This is denoted as $\mu\downarrow\e$.
}\eit

\bit{
\item[$\bullet$] A monomorphism $\mu$ is called a {\it strong monomorphism}\label{DEF:SMono}, if it is diagonalisable from the right by any epimorphism, i.e. for any epimorphism $\e$ and for any morphisms $\alpha$ and $\beta$ such that $\beta\circ\e=\mu\circ\alpha$ there is a (necessarily, unique) morphism $\delta$ such that the diagram \eqref{DIAGR:diagonaliziruemost} is commutative. The class of all strong monomorphisms $\mu$ which come to (i.e.  $\Ran\mu=X$) is denoted by $\SMono_X$.

\item[$\bullet$] Dually, an epimorphism $\e$ is called a {\it strong epimorphism}\label{DEF:SEpi}, if it is diagonalizable from the left by any monomorphism, i.e. for any monomorphism $\mu$ and for any morphisms $\alpha$ and $\beta$ such that $\beta\circ\e=\mu\circ\alpha$, there is a (necessarily, unique) morphism $\delta$ such that the diagram \eqref{DIAGR:diagonaliziruemost} is commutative. The class of all strong epimorphisms which go from $X$ (i.e. $\Dom\e=X$) is denoted by $\SEpi^X$.
}\eit

\brem
Formally the relation $\mu\downarrow\e$ does not mean automatically that $\mu\in\Mono$ and $\e\in\Epi$: in the category of vector spaces over $\C$ a pair of morphisms $\mu=0:\C\to 0$ and $\e=0:0\to\C$ is diagonalizable:
$$
\begin{diagram}
\node{0}\arrow{s,l}{\alpha}\arrow{e,t}{\e=0}\node{\C}\arrow{s,r}{\beta}\arrow{sw,t,--}{\delta}
\\
\node{\C}\arrow{e,b}{\mu=0}\node{0}
\end{diagram}
$$
\erem

\brem\label{REM:o-razbienii-kvadrata} If $\mu\in\Mono$ (or $\e\in\Epi$), then $\delta$ is unique: if $\delta'$ is another morphism with the same property, then
$$
\mu\circ\delta=\beta=\mu\circ\delta'\quad\Longrightarrow\quad \delta=\delta'.
$$
Besides this, the commutativity of the upper triangle in \eqref{DIAGR:diagonaliziruemost} implies the commutativity of the lower triangle, and vise versa. For example,
 \beq\label{VERH-treug<=>NIZHN-treug}
\alpha=\delta\circ\e\quad\Longrightarrow\quad \beta\circ\underset{\scriptsize\begin{matrix}\text{\rotatebox{90}{$\owns$}}\\ \Epi\end{matrix}}{\e}=\mu\circ\alpha=\mu\circ\delta\circ\underset{\scriptsize\begin{matrix}\text{\rotatebox{90}{$\owns$}}\\ \Epi\end{matrix}}{\e}\quad\Longrightarrow\quad \beta=\mu\circ\delta
 \eeq
\erem

Obviously these classes of morphisms are connected by the following implications
$$
\text{$\mu$ is a strong monomorphism}\quad\Longrightarrow\quad
\text{$\mu$ is an immediate monomorphism}\quad\Longrightarrow\quad
\text{$\mu$ is a monomorphism},
$$
$$
\text{$\e$ is a strong epimorphism}\quad\Longrightarrow\quad
\text{$\e$ is an immediate epimorphism}\quad\Longrightarrow\quad
\text{$\e$ is an epimorphism}.
$$

Let $\varGamma$ be a class of monomorphisms in a category ${\tt K}$, and all local identities belong to it:
$$
\{1_X;\ X\in\Ob({\tt K})\}\subseteq\varGamma\subseteq\Mono({\tt K})
$$
(the key examples are the classes $\varGamma=\Mono$ and $\varGamma=\SMono$). For each object $X$ in ${\tt K}$ let us denote by $\varGamma_X$ the class of all morphisms in $\varGamma$ with $X$ as range:
 \beq\label{DEF:varGamma(X)}
 \varGamma_X=\{\sigma\in\varGamma:\quad \Ran\sigma=X\}.
 \eeq
It is a category, where a morphism $\rho\overset{\varkappa}{\longrightarrow}\sigma$ from an object $\rho\in\varGamma_X$ into an object $\sigma\in\varGamma_X$, i.e. a monomorphism $\rho:A\to X$ into a monomorphism $\sigma:B\to X$, is an arbitrary morphism $\varkappa:A\to B$ in ${\tt K}$ such that the following diagram is commutative:
\beq\label{morphism-v-Mono(X)}
\xymatrix @R=1pc @C=2pc
{
A\ar[rd]^{\rho}\ar[dd]_{\varkappa} &   \\
  & X \\
B\ar[ru]_{\sigma} &
}
\eeq
Actually, this diagram in the initial category ${\tt K}$ can be considered as a morphism $\rho\overset{\varkappa}{\longrightarrow}\sigma$ in the category  $\varGamma_X$. A composition of such morphisms $\rho\overset{\varkappa}{\longrightarrow}\sigma$ and $\sigma\overset{\lambda}{\longrightarrow}\tau$, i.e. of diagrams
$$
\xymatrix @R=1pc @C=2pc
{
A\ar[rd]^{\rho}\ar[dd]_{\varkappa} &   \\
  & X \\
B\ar[ru]_{\sigma} &
}
\qquad
\xymatrix @R=1pc @C=2pc
{
B\ar[rd]^{\sigma}\ar[dd]_{\lambda} &   \\
  & X \\
C\ar[ru]_{\tau} &
}
$$
is a morphism $\rho\overset{\lambda\circ\varkappa}{\longrightarrow}\tau$, i.e. a diagram
$$
\xymatrix @R=1pc @C=2pc
{
A\ar[rd]^{\rho}\ar[dd]_{\lambda\circ\varkappa} &   \\
  & X \\
C\ar[ru]_{\tau} &
}
$$
One can conceive it as a result of splicing of  the initial diagrams along the common edge $\sigma$, adding the arrow of composition $\varkappa\circ\lambda$, and then throwing away the vertex $B$ together with all its incidental edges:
$$
\xymatrix % @R=1pc @C=2pc
{
A\ar@/^2ex/[rrd]^{\rho}\ar[dd]_{\lambda\circ\varkappa}\ar@{-->}[rd]_{\varkappa} & &  \\
  & B\ar@{-->}[r]_{\sigma}\ar@{-->}[ld]_{\lambda} &  X \\
  C\ar@/_2ex/[rru]_{\tau} & &
}
$$
Of course, local identities in $\varGamma_X$ are diagrams of the form
$$
\xymatrix @R=1pc @C=2pc
{
A\ar[rd]^{\rho}\ar[dd]_{1_X} &   \\
  & X \\
A\ar[ru]_{\rho} &
}
$$

\bit{

\item[$\bullet$] A {\it system of subobjects of the class $\varGamma$} in an object $X$ of a category ${\tt K}$ is an arbitrary skeleton $S$ of the category $\varGamma_X$, such that the morphism $1_X$ belongs to $S$. In other words, a subclass $S$ in $\varGamma_X$ is a system of subobjects in $X$, if
\bit{

\item[(a)]\label{DEF:toch-sist-podobj-a} the local identity of $X$ belongs to $S$:
$$
1_X\in S,
$$

\item[(b)]\label{DEF:toch-sist-podobj-b} every monomorphism $\mu\in\varGamma_X$ has an isomorphic monomorphism in the class $S$:
$$
\forall\mu\in\varGamma_X\qquad\exists\sigma\in S\qquad \mu\cong\sigma.
$$

\item[(c)]\label{DEF:toch-sist-podobj-c} in $S$ an isomorphism (in the sense of category $\varGamma_X$) is equivalent to the identity:
$$
\forall\sigma,\tau\in S\qquad \Big( \sigma\cong\tau\quad\Longleftrightarrow\quad \sigma=\tau \Big)
$$
}\eit
In the Morse-Kelley axiomatics of Set Theory such a class $S$ always exists. When it is chosen, its elements are called {\it subobjects of the class $\varGamma$ of the object $X$}. The class $S$ is endowed with a structure of the full subcategory in $\varGamma_X$.
}\eit

 \bit{
\item[$\bullet$] We say that a category ${\tt K}$ is {\it well-powered in the class $\varGamma$}, if each object $X$ has a system of subobjects $S$ of the class $\varGamma$, which is a set (i.e. not a proper class).
 }\eit

\bex The standard categories frequently used as examples, like the category of sets, groups, vector spaces, algebras (over a given field), topological spaces, topological vector spaces, topological algebras, etc., are, obviously, well-powered in the class $\Mono$.
\eex

 \btm \cite[Theorem 2.24]{Akbarov-env}\label{TH:o-lok-malosti-v-podobjektah}
In the Morse-Kelley axiomatics of Set Theory if a category ${\tt K}$ is well-powered in subobjects of a class $\varGamma$, then there is a map $X\mapsto S_X$ which assigns to each object $X$ in ${\tt K}$ its system of subobjects $S_X$ of the class $\varGamma$ (and $S_X$ is a set).
 \etm

Let $\varOmega$ be a class of epimorphisms in a category ${\tt K}$, and all local identities belong to it:
$$
\{1_X;\ X\in\Ob({\tt K})\}\subseteq\varOmega\subseteq\Epi({\tt K})
$$
(the key examples are the classes $\varOmega=\Epi$ and $\varOmega=\SEpi$). For each object $X$ in ${\tt K}$ we denote by $\varOmega^X$ the class of all morphisms in $\varOmega$ with the domain $X$:
 \beq\label{DEF:varOmega(X)}
 \varOmega^X=\{\sigma\in\varOmega:\quad \Dom\sigma=X\}.
 \eeq
This class forms a category where a morphism $\rho\overset{\varkappa}{\longrightarrow}\sigma$ from an object $\rho\in\varOmega^X$ into an object $\sigma\in\varOmega^X$, i.e. from an epimorphism $\rho:X\to A$ into an epimorphism $\sigma:X\to B$, is an arbitrary morphism $\varkappa:A\to B$ in ${\tt K}$ such that the following diagram is commutative
\beq\label{morphism-v-Epi(X)}
\xymatrix @R=1pc @C=2pc
{
 & A\ar[dd]^{\varkappa}    \\
X\ar[ru]^{\rho}\ar[rd]_{\sigma}  &  \\
 & B
}
\eeq
Actually, this diagram in the initial category ${\tt K}$ can be considered as a morphism $\rho\overset{\varkappa}{\longrightarrow}\sigma$ in $\varOmega^X$. A composition of two such morphisms $\rho\overset{\varkappa}{\longrightarrow}\sigma$ and $\sigma\overset{\lambda}{\longrightarrow}\tau$, i.e. diagrams $$
\xymatrix @R=1pc @C=2pc
{
 & A\ar[dd]^{\varkappa}    \\
X\ar[ru]^{\rho}\ar[rd]_{\sigma}  &  \\
 & B
}
\qquad
\xymatrix @R=1pc @C=2pc
{
 & B\ar[dd]^{\lambda}    \\
X\ar[ru]^{\sigma}\ar[rd]_{\tau}  &  \\
 & C
}
$$
is a morphism $\rho\overset{\lambda\circ\varkappa}{\longrightarrow}\tau$, i.e. a diagram
$$
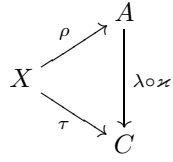
\xymatrix @R=1pc @C=2pc
{
 & A\ar[dd]^{\lambda\circ\varkappa}    \\
X\ar[ru]^{\rho}\ar[rd]_{\tau}  &  \\
 & C
}
$$
One can conceive it as a result of splicing of  the initial diagrams along the common edge $\sigma$, adding the arrow of composition $\lambda\circ\varkappa$, and then throwing away the vertex $B$ together with all its incidental edges:
$$
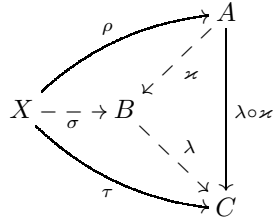
\xymatrix % @R=1pc @C=2pc
{
 & & A\ar[dd]^{\lambda\circ\varkappa}\ar@{-->}[ld]^{\varkappa}    \\
X\ar@/^2ex/[rru]^{\rho}\ar@{-->}[r]_{\sigma}\ar@/_2ex/[rrd]_{\tau}  & B\ar@{-->}[rd]^{\lambda} &  \\
 & & C
}
$$
Of course, local identities in $\varOmega^X$ are diagrams of the form
$$
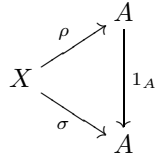
\xymatrix @R=1pc @C=2pc
{
 & A\ar[dd]^{1_A}    \\
X\ar[ru]^{\rho}\ar[rd]_{\sigma}  &  \\
 & A
}
$$

It is useful to define a preorder $\to$ in $\varOmega^X$:
 \beq\label{DEF:le-in-F_X}
\rho\to\sigma\Longleftrightarrow\quad \exists \iota\in\Mor({\tt K})\quad
\sigma=\iota\circ\rho.
 \eeq
The morphism $\iota$, if it exists, is unique, and is an epimorphism (since $\rho$ and $\sigma$ are epimorphisms). As a corollary, there is an operation, which to each pair of morphisms $\rho,\sigma\in\varOmega^X$ with the property $\rho\to\sigma$ assigns the morphism $\iota=\iota^\sigma_\rho$ in
\eqref{DEF:le-in-F_X}:
 \beq\label{DEF:le-in-F_X-*}
\sigma=\iota^\sigma_\rho\circ\rho.
 \eeq
If $\pi\to\rho\to\sigma$, then the chain
$$
\iota^\sigma_\pi\circ\pi=\sigma=\iota^\sigma_\rho\circ\rho=\iota^\sigma_\rho\circ\iota^\rho_\pi\circ\pi,
$$
implies by epimorphy of $\pi$ the equality
\beq\label{iota_rho^tau=iota_rho^sigma-circ-iota_sigma^tau}
\iota^\sigma_\pi=\iota^\sigma_\rho\circ\iota^\rho_\pi.
 \eeq

\bit{

\item[$\bullet$] A {\it system of quotient objects of the class $\varOmega$} on an object $X$ in a category ${\tt K}$ is an arbitrary skeleton $Q$ of the category $\varOmega^X$, such that $1_X$ belongs to $Q$. In other words, a subclass $Q$ in $\varOmega^X$ is called a system of quotient objects on $X$, if
\bit{

\item[(a)] the local identity of $X$ belongs to $Q$:
$$
1_X\in Q,
$$

\item[(b)] every epimorphism $\e\in\varOmega^X$ has an isomorphic epimorphism in $Q$:
$$
\forall\e\in\varOmega^X\qquad\exists\pi\in Q\qquad \e\cong\pi,
$$

\item[(c)] in $Q$ an isomorphism (in the sense of category $\varOmega^X$) is equivalent to the identity:
$$
\forall\pi,\rho\in Q\qquad \Big( \pi\cong\rho\quad\Longleftrightarrow\quad \pi=\rho \Big)
$$
}\eit
In the Morse-Kelley axiomatics of Set Theory such a class $Q$ always exists. When it is chosen, its elements are called {\it quotient objects of the class $\varOmega$ of the object $X$}. The class $Q$ is endowed with the structure of full subcategory in $\varOmega^X$.

}\eit

 \bit{
\item[$\bullet$] We say that a category ${\tt K}$ is {\it co-well-powered in the class $\varOmega$}, if each object  $X$ has a system of quotient objects $Q$ of the class $\varOmega$, which is a set (i.e. not a proper class).
 }\eit

\bex Among the standard categories -- the category of sets, groups, vector spaces, algebras over a given field, topological spaces, topological vector spaces, topological algebras -- some are co-well-powered in the class $\Epi$, but sometimes this is not easy to prove (see \cite{Adamek-Rosicky}). In contrast to this the co-well-poweredness in the class $\SEpi$ is verified much easier. \eex

 \btm \cite[Theorem 2.31]{Akbarov-env} \label{TH:o-lok-malosti-v-faktor-objektah}
In the Morse-Kelley axiomatics of Set Theory if a category ${\tt K}$ is co-well-powered in the class $\varOmega$, then there exists a map $X\mapsto Q_X$ which assigns to any object $X$ in ${\tt K}$ a system of its quotient-objects $Q_X$ of the class $\varOmega$ (and $Q_X$ is a set).
 \etm

\paragraph{Nodal decomposition.} Suppose a morphism $\ph$ in a category $\tt K$ is decomposed into a composition
\beq\label{DEF:uzlovoe-razlozhenie}
\ph=\iota\circ\rho\circ\gamma,
\eeq
where
\bit{

\item[(i)] $\gamma$ is a strong epimorphism,

\item[(ii)] $\rho$ is a bimorphism,

\item[(iii)] $\iota$ is a strong monomorphism.

}\eit
Then the triple $(\iota,\rho,\gamma)$ is called a {\it nodal decomposition} of the morphism $\ph$.

\btm A nodal decomposition, if it exists, is unique: if $(\iota,\rho,\gamma)$ and $(\iota',\rho',\gamma')$ are two nodal decompositions of the morphism $\ph$,
$$
\xymatrix %@R=2.5pc @C=4.0pc
{
X\ar[rrr]^{\ph}\ar[rd]^{\gamma'}\ar[dd]_{\gamma} & & & Y \\
& P'\ar[r]^{\rho'} & Q'\ar[ru]^{\iota'} & \\
P\ar[rrr]_{\rho} & & & Q \ar[uu]_{\iota}
}
$$
then there exist (necessarily, unique) isomorphisms $\sigma:P\to P'$ and $\tau:Q'\to Q$ such that the following diagram is commutative:
\beq\label{otnoshenie-mezhdu-stand-razlozh}
\xymatrix %@R=2.5pc @C=4.0pc
{
X\ar[rrr]^{\ph}\ar[rd]^{\gamma'}\ar[dd]_{\gamma} & & & Y \\
& P'\ar[r]^{\rho'} & Q'\ar[ru]^{\iota'}\ar@{-->}[rd]_{\tau} & \\
P\ar[rrr]_{\rho}\ar@{-->}[ru]_{\sigma} & & & Q \ar[uu]_{\iota}
}
\eeq
\etm

\bit{
\item[$\bullet$]
From the uniqueness (up to isomorphism) of the nodal decomposition $\ph=\iota'\circ\rho'\circ\gamma'$ it follows that one can assign notations to its components. We will further depict a nodal decomposition of a morphism $\ph:X\to Y$ as a diagram
\beq\label{DEF:oboznacheniya-dlya-uzlov-razlozh}
\begin{diagram}
\node{X}\arrow{s,l}{\coim_\infty\ph}\arrow{e,t}{\ph}\node{Y} \\
\node{\Coim_\infty\ph}\arrow{e,t}{\red_\infty\ph}\node{\Im_\infty\ph}\arrow{n,r}{\im_\infty\ph}
\end{diagram}
\eeq
(where elements are defined up to isomorphisms). The proof of Theorem \ref{TH:sush-uzlov-razlozh-v-polnoi-lok-maloi-kateg} below and Remark \ref{REM:struktura-uzlov-razlozh-v-kateg-s-nulem} justify these notations, since they show that $\coim_\infty$, $\red_\infty$ and $\im_\infty$ can be conceived as a sort of ``transfinite induction'' of the usual operation $\coim$, $\red$ and $\im$ in preabelian categories:
\begin{align*}
& \coim_\infty=\lim_{n\to\infty}\underbrace{\coim\circ\coim\circ...\circ\coim}_{\text{$n$ multipliers}} \\
& \red_\infty=\lim_{n\to\infty}\underbrace{\red\circ\red\circ...\circ\red}_{\text{$n$ multipliers}} \\
& \im_\infty=\lim_{n\to\infty}\underbrace{\im\circ\im\circ...\circ\im}_{\text{$n$ multipliers}}
\end{align*}
We will call
\bit{
\item[---] $\im_\infty\ph$ a {\it nodal image},

\item[---] $\red_\infty\ph$ a {\it nodal reduced part},

\item[---] $\coim_\infty\ph$  a {\it nodal coimage}
 }\eit
 of the morphism $\ph$.
}\eit

\bit{
\item[$\bullet$]\label{DEF:strog-epi-razlich-mono} Let us say that in a category ${\tt K}$

 \bit{
\item[---]  {\it strong epimorphisms discern monomorphisms}, if the reverse is true: from the fact that a morphism $\mu$ is not a monomorphism it follows that $\mu$ can be represented as a composition $\mu=\mu'\circ\e$, where $\e$ is a strong epimorphism, which is not an isomorphism,

\item[---] {\it strong monomorphisms discern epimorphisms}, if the reverse is true: from the fact that a morphism $\e$ is not an epimorphism it follows that $\e$ can be represented as a composition $\e=\mu\circ\e'$, where $\mu$ is a strong monomorphism, which is not an isomorphism.
 }\eit

\item[$\bullet$]\label{DEF:kat-s-uzlov-razlozheniem}
    We say also that ${\tt K}$ is a {\it category with a nodal decomposition}, if every morphism $\ph$ in ${\tt K}$ has a nodal decomposition.
 }\eit

\btm\cite[Theorem 2.36]{Akbarov-env} \label{TH:sush-uzlov-razlozh-v-polnoi-lok-maloi-kateg} Let ${\tt K}$ ba a linearly complete, well-powered in strong monomorphisms and co-well-powered in strong epimorphisms category, where strong epimorphisms discern monomorphisms, and, dually, strong monomorphisms discern epimorphisms. Then ${\tt K}$ is a category with nodal decomposition.
 \etm

\paragraph{Connection with  the base decomposition in pre-Abelian categories.}

Recall that in {\it pre-Abelian category} \cite{Bucur-Deleanu,General-algebra} each morphism $\ph:X\to Y$ has a kernel and a cokernel. This implies that for each $\ph$ there exists a unique morphism $\red\ph$ such that the following diagram is commutative:
\beq\label{EX:bazis-razlozh}
\begin{diagram}
\node{X}\arrow{s,l}{\coim\ph}\arrow{e,t}{\ph}\node{Y} \\
\node{\Coim\ph}\arrow{e,t,--}{\red\ph}\node{\Im\ph}\arrow{n,r}{\im\ph}
\end{diagram}
\eeq
where the morphism $\coim\ph=\coker(\ker\ph)$ is called a {\it coimage} of $\ph$, the morphism  $\im\ph=\ker(\coker\ph)$ an {\it image} of $\ph$. The morphism $\red\ph$ is called the {\it reduced part} of $\ph$.

\bit{
\item[$\bullet$] The decomposition \eqref{EX:bazis-razlozh} will be called the {\it basic decomposition} of the morphism $\ph$.
}\eit

{\it If a category ${\tt K}$ is Abelian, then each basic decomposition in it is nodal}. But if ${\tt K}$ is not Abelian, then these composition do not necessarily coincide \cite[Example 4.98]{Akbarov-env}.

\btm\cite[Theorem 2.42]{Akbarov-env}\footnote{In the English text of the paper \cite{Akbarov-env} this theorem has a misprint: the condition of linear completeness is omited there.}\label{TH:faktorizatsija-v-lok-maloi-kategorii} If a category $\tt K$ is pre-Abelian, linearly complete, well-powered in the class $\SMono$ and co-well-powered in the class $\SEpi$, then $\tt K$ is a category with nodal decomposition.
\etm

\btm\cite[Remark 2.44]{Akbarov-env}\label{TH:svayz-bazisnogo-i-uzlovogo-razlozheniya} In a pre-Abelian category if a morphism $\ph$ has a nodal decomposition $\ph=\im_\infty\ph\circ\red_\infty\ph\circ\coim_\infty\ph$, then there exist unique morphisms $\sigma$ and $\tau$ such that the following diagram is commutative:
\beq\label{svayz-bazisnogo-i-uzlovogo-razlozheniya}
\xymatrix  @R=2.5pc @C=4.0pc
{
X\ar[rrr]^{\ph}\ar[rd]^{\coim_\infty\ph}\ar[dd]_{\coim\ph} & & & Y \\
& \Coim_\infty\ph\ar[r]_{\red_\infty\ph} & \Im_\infty\ph\ar[ru]^{\im_\infty\ph}\ar@{-->}[rd]_{\tau} & \\
\Coim\ph\ar[rrr]_{\red\ph}\ar@{-->}[ru]_{\sigma} & & & \Im\ph \ar[uu]_{\im\ph}
}
\eeq
where $\ph=\im\ph\circ\red\ph\circ\coim\ph$ is a basic decomposition of $\ph$.
\etm

\brem
If $\tt K$ is not Abelian, then $\sigma$ and $\tau$ are not necessarily isomorphisms \cite[Example 4.98]{Akbarov-env}.
\erem

\paragraph{Factorization of a category.} Recall that the notion of diagonalizability was defined on page \pageref{DIAGR:diagonaliziruemost}.

\bit{
\item For each class of morphisms $\varLambda$ in $\tt K$

 \bit{
\item[---] its {\it epimorphic conjugate class} is the class
$$
\varLambda^\downarrow=\{\e\in\Epi({\tt K}):\forall\lambda\in\varLambda\quad \lambda\downarrow\e\}.
$$
\item[---] its {\it monomorphic conjugate class} is the class
$$
{^\downarrow\kern-1pt\varLambda}=\{\mu\in\Mono({\tt K}):\forall\lambda\in\varLambda\quad \mu\downarrow\lambda\}.
$$
}\eit
}\eit
Clearly, for each class of morphisms $\varLambda$
\begin{align}
& \Iso\subseteq\varLambda^\downarrow\subseteq\Epi, && \Iso\circ\ \varLambda^\downarrow\subseteq\varLambda^\downarrow
\label{Iso-subseteq-varTheta^downarrow-subseteq-Epi}
\\
& \Iso\subseteq{^\downarrow\kern-1pt\varLambda}\subseteq\Mono, &&
{^\downarrow\kern-1pt\varLambda}\circ\Iso\subseteq{^\downarrow\kern-1pt\varLambda}
\label{Iso-subseteq-^downarrow-varTheta-subseteq-Mono}
\end{align}

\bit{

\item
Let us say that classes of morphisms $\varGamma$ and $\varOmega$ define a {\it factorization of the category\footnote{This construction is also called a {\it bicategory} \cite{General-algebra,Tsalenko-Shulgeifer}.}
$\tt K$}\label{DEF:faktorizatsija-v-kategorii}, if
 \bit{
 \item[F.1] $\varOmega$ is the epimorphic conjugate class for $\varGamma$:
 $$
 \varGamma^\downarrow=\varOmega
 $$

 \item[F.2] $\varGamma$ is the monomorphic conjugate class for $\varOmega$:
 $$
 \varGamma={^\downarrow\varOmega},
 $$

 \item[F.3] the composition of the class $\varGamma$ and $\varOmega$ covers the class of all morphisms:
 $$
 \varGamma\circ\varOmega=\Mor({\tt K})
 $$
 (this means that each morphism $\ph\in\Mor({\tt K})$ can be represented as a composition $\mu\circ\e$, where  $\mu\in\varGamma$, $\e\in\varOmega$).
 }\eit
If these conditions are fulfilled, we write
 \beq\label{K=varGamma-circledcirc-varOmega}
{\tt K}=\varGamma\circledcirc\varOmega.
 \eeq
}\eit

\bex\label{TH:faktorizatsija-v-kategorii-s-uzlov-razlozh-1}
In a category ${\tt K}$ with the nodal decomposition the following classes of morphisms define factorizations:
$$
{\tt K}=\Mono\circledcirc\SEpi=\SMono\circledcirc\Epi.
$$
\eex

\btm\cite[Theorem 8.2]{Tsalenko-Shulgeifer}\label{TH:o-faktorizatsii}
Classes $\varGamma$ and $\varOmega$ define a factorization of $\tt K$
$$
{\tt K}=\varGamma\circledcirc\varOmega
$$
if and only if the following conditions hold:
\bit{
 \item[(i)] $\varGamma\subseteq\Mono({\tt K})$ and $\varOmega\subseteq\Epi({\tt K})$,

 \item[(ii)] $\Iso({\tt K})\subseteq\varOmega\cap\varGamma$,

 \item[(iii)] for each morphism $\ph\in\Mor({\tt K})$ there is a decomposition
 \beq\label{faktorizatsiya-v-kat-s-faktoriz}
 \ph=\mu_{\ph}\circ\e_{\ph},\qquad \mu_{\ph}\in\varGamma,\quad \e_{\ph}\in\varOmega
 \eeq
 \item[(iv)] for any other decomposition with the same properties
$$
 \ph=\mu'\circ\e',\qquad \mu'\in\varGamma,\quad \e'\in\varOmega
$$
there is a morphism $\theta\in\Iso({\tt K})$ such that
$$
\mu'=\mu_{\ph}\circ\theta,\qquad \e'=\theta^{-1}\circ\e_{\ph}.
$$
 }\eit
\etm

 \bit{
\item Let us say that a class of morphisms $\varOmega$ in $\tt K$ is
{\it monomorphically complementable}\label{DEF:klass-monomorfno-dopolnyaem}, if
\beq\label{klass-monomorfno-dopolnyaem}
{\tt K}={^\downarrow\varOmega}\circledcirc\varOmega.
\eeq
In other words, $\varOmega$ must be epimorphic conjugate to its monomorphic conjugate class
$$
\varOmega=(^\downarrow\varOmega)^\downarrow,
$$
and the composition of ${^\downarrow\varOmega}$ and $\varOmega$ must cover the class of all morphisms:
$$
{^\downarrow\varOmega}\circ\varOmega=\Mor(\tt K).
$$
In this case the class ${^\downarrow\varOmega}$ will be called the {\it monomorphuc complement} to $\varOmega$.
}\eit

\brem
From \eqref{Iso-subseteq-varTheta^downarrow-subseteq-Epi} it follows that if a class of morphisms $\varOmega$ is monomorphically complementable, then
\beq\label{Iso-circ-varOmega-subseteq-varOmega}
\Iso\subseteq\varOmega\subseteq\Epi,\qquad \Iso\circ\varOmega\subseteq\varOmega
\eeq
\erem

 \bit{
\item Similarly, we say that the class of morphisms $\varGamma$ in $\tt K$ is
{\it epimorphically complementable}\label{DEF:klass-epimorfno-dopolnyaem}, if
\beq\label{klass-epimorfno-dopolnyaem}
{\tt K}=\varGamma\circledcirc\varGamma^\downarrow.
\eeq
In other words, $\varGamma$ must be the monomorphic conjugate to its epimorphic conjugate class
$$
\varGamma={^\downarrow(\varGamma^\downarrow)},
$$
and the composition of the classes $\varGamma$ and $\varGamma^\downarrow$ must cover the class of all morphisms:
$$
\varGamma\circ\varGamma^\downarrow=\Mor(\tt K).
$$
In this case the class $\varGamma^\downarrow$ will be called the {\it epimorphic complement} to $\varGamma$.
}\eit

\brem
From \eqref{Iso-subseteq-^downarrow-varTheta-subseteq-Mono} it follows that if a class $\varGamma$ is epimorphically complementable, then
\beq\label{varGamma-circ-Iso-subseteq-varGamma}
\Iso\subseteq\varGamma\subseteq\Mono,\qquad \varGamma\circ\Iso\subseteq\varGamma.
\eeq
\erem

\subsection{Envelopes and refinements}

\paragraph{Envelopes.}

 \bit{
\item[$\bullet$] A morphism $\sigma:X\to X'$ in a category $\tt K$ is called an
{\it extension of the object $X\in\Ob({\tt K})$ in the class of morphisms $\varOmega$ with respect to the class of morphisms $\varPhi$}, if $\sigma\in\varOmega$, and for any morphism
$\ph:X\to B$ from the class $\varPhi$ there exists a unique morphism $\ph':X'\to B$ in
${\tt K}$ such that the following diagram is commutative:
 \beq\label{DEF:diagr-rasshirenie}
\begin{diagram}
\node[2]{X} \arrow{sw,t}{\varOmega\owns\sigma} \arrow{se,t}{\forall\ph\in\varPhi}\\
\node{X'}  \arrow[2]{e,b,--}{\exists!\ph'} \node[2]{B}
\end{diagram}
\eeq

\item[$\bullet$] An extension $\rho:X\to E$ of an object $X\in\Ob({\tt K})$ in the class of morphisms $\varOmega$
with respect to the class of morphisms $\varPhi$ is called an {\it envelope of $X$ in $\varOmega$ with respect to $\varPhi$}, if for any other extension $\sigma:X\to X'$ (of $X$ in $\varOmega$ with respect to $\varPhi$) there is a unique morphism $\upsilon:X'\to E$ in $\tt K$ such that the following diagram is commutative:
 \beq\label{DEF:diagr-obolochka}
\begin{diagram}
\node[2]{X} \arrow{sw,t}{\forall\sigma} \arrow{se,t}{\rho}\\
\node{X'}  \arrow[2]{e,b,--}{\exists!\upsilon} \node[2]{E}
\end{diagram}
 \eeq
For the morphism of envelope $\rho:X\to E$ we use the notation
    \beq\label{DEF:env_F^L(X)}
    \rho=\env_{\varPhi}^\varOmega X.
    \eeq
The very object $E$ is also called an {\it envelope} of $X$ (in $\varOmega$ with respect to $\varPhi$), and we use the following notation for it:
    \beq\label{DEF:Env_F^L(X)}
E=\Env_{\varPhi}^\varOmega X.
 \eeq
}\eit

 \bit{
\item[$\bullet$] Let us say that in a category ${\tt K}$ {\it a class of morphisms $\varPhi$ is generated on the inside by a class of morphisms  $\varPsi$}, if
    \beq\label{DEF:morfizmy-porozhdayutsya-iznutri}
    \varPsi\subseteq\varPhi\subseteq\Mor({\tt K})\circ\ \varPsi.
    \eeq
  }\eit

\btm\cite[Theorem 3.5]{Akbarov-env}\label{TH:morfizmy-porozhdayutsya-iznutri}  Suppose that in a category ${\tt K}$ a class of morphisms $\varPhi$ is generated on the inside by a class of morphisms  $\varPsi$. Then for any class of epimorphisms $\varOmega$ (it  is not necessary that $\varOmega$ contains all epimorphisms of ${\tt K}$) and for any object $X$ the existence of envelope $\env^\varOmega_\varPsi X$ is equivalent to the existence of envelope $\env^\varOmega_\varPhi X$, and these envelopes coincide:
 \beq\label{env_Psi=env_Phi}
\env^\varOmega _\varPsi X=\env^\varOmega _\varPhi X.
 \eeq
 \etm

\bit{
\item[$\bullet$]\label{DEF:varPhi-razlich-morfizmy-snaruzhi} Let us say that a class of morphisms $\varPhi$ in a category ${\tt K}$  {\it differs morphisms on the outside}, if for any two different parallel morphisms $\alpha\ne\beta:X\to Y$ there is a morphism $\ph:Y\to M$ from the class $\varPhi$ such that $\ph\circ\alpha\ne\ph\circ\beta$.
}\eit

\btm\cite[Theorem 3.6]{Akbarov-env}\label{TH:Phi-razdel-moprfizmy}
If a class of morphisms $\varPhi$ differs morphisms on the outside, then for any class of morphisms $\varOmega$
 \bit{
\item[(i)] each extension in $\varOmega$ with respect to $\varPhi$ is a monomorphism,

\item[(ii)] an envelope with respect to $\varPhi$ in $\varOmega$ exists if and only if there exists an envelope with respect to $\varPhi$ in the class $\varOmega\cap\Mono$; in this case these envelopes coincide:
$$
\env_{\varPhi}^{\varOmega}=\env_{\varPhi}^{\varOmega\cap\Mono},
$$

\item[(iii)] if the class $\varOmega$ contains all monomorphisms,
 $$
 \varOmega\supseteq\Mono,
 $$
then the existence of the envelope with respect to $\varPhi$ in $\Mono$ automatically implies the existence of envelope with respect to $\varPhi$ in $\varOmega$, and the coincidence of these envelopes:
$$
\env_{\varPhi}^{\varOmega}=\env_{\varPhi}^{\Mono}.
$$

 }\eit\etm

\bit{
\item[$\bullet$] Let us remind that a class of morphisms $\varPhi$ in a category ${\tt K}$ is called a {\it right ideal}, if
$$
\varPhi\circ\Mor({\tt K})\subseteq\varPhi
$$
(i.e. for any $\ph\in\varPhi$ and for any morphism $\mu$ in ${\tt K}$ the composition $\ph\circ\mu$ belongs to $\varPhi$).
}\eit

\btm\cite[Theorem 3.7]{Akbarov-env}\label{TH:Phi-razdel-moprfizmy-*} If a class of morphisms $\varPhi$ differs morphisms on the outside and is a right ideal in the category ${\tt K}$, then for any class of morphisms $\varOmega$
 \bit{
\item[(i)] each extension in $\varOmega$ with respect to $\varPhi$ is a bimorphism,

\item[(ii)] an envelope with respect to $\varPhi$ in $\varOmega$ exists if an only if there exists an envelope with respect to $\varPhi$ in the class $\varOmega\cap\Bim$ of bimorphisms belonging to $\varOmega$; in this case these envelopes coincide:
$$
\env_{\varPhi}^{\varOmega}=\env_{\varPhi}^{\varOmega\cap\Bim}.
$$

\item[(iii)] if the class $\varOmega$ contains all bimorphisms,
 $$
 \varOmega\supseteq\Bim,
 $$
 then an envelope with respect to $\varPhi$ in $\varOmega$ exists if an only if there exists an envelope with respect to $\varPhi$ in $\Bim$, and these envelopes coincide:
$$
\env_{\varPhi}^{\varOmega}=\env_{\varPhi}^{\Bim}.
$$

 }\eit
\etm

A special case of envelope is the construction, where $\varOmega$ and/or $\varPhi$ are classes of all morphisms into the objects from some given subclasses in $\Ob({\tt K})$. The accurate formulation for the case, when both classes $\varOmega$ and $\varPhi$ are defined in such a way is the following. Suppose we have a category ${\tt K}$ and two subclasses ${\tt L}$ and ${\tt M}$ in the class $\Ob({\tt K})$ of objects in ${\tt K}$.

 \bit{
\item[$\bullet$] A morphism $\sigma:X\to X'$ is called an {\it extension of the object
$X\in{\tt K}$ in the class ${\tt L}$ with respect to the class ${\tt M}$}, if $X'\in{\tt L}$ and for any object
$B\in{\tt M}$ and any morphism $\ph:X\to B$ there exists a unique morphism $\ph':X'\to B$ such that the following diagram is commutative:
$$
\put(6,-35){$\begin{matrix}\text{\rotatebox{90}{$\owns$}}\\ {\tt L}
\end{matrix}$} \put(73,-35){$\begin{matrix}\text{\rotatebox{90}{$\owns$}}\\
{\tt M} \end{matrix}$}
\begin{diagram}
\node[2]{X} \arrow{sw,t}{\sigma} \arrow{se,t}{\forall\ph}\\
\node{X'}  \arrow[2]{e,b,--}{\exists!\ph'} \node[2]{B}
\end{diagram}
$$

\item[$\bullet$]\label{DEF:obolochka-otn-klassa} An extension $\rho:X\to E$ of an object $X\in{\tt K}$ in the class ${\tt L}$ with respect to the class ${\tt M}$ is called an {\it envelope of the object $X\in{\tt K}$ in the class ${\tt L}$ with respect to the class ${\tt M}$}, and we denote this by formula
    \beq\label{DEF:env_M^L(A)}
    \rho=\env_{\tt M}^{\tt L} X,
    \eeq
    if for any  other extension $\sigma:X\to X'$ (of the object $X$ in the class ${\tt L}$ with respect to the class ${\tt M}$) there exists a unique morphism $\upsilon:X'\to E$ such that the following diagram is commutative:
 \beq\label{diagr:obolochka}
\put(6,-35){$\begin{matrix}\text{\rotatebox{90}{$\owns$}}\\ {\tt L}
\end{matrix}$} \put(73,-35){$\begin{matrix}\text{\rotatebox{90}{$\owns$}}\\
{\tt L} \end{matrix}$}
\begin{diagram}
\node[2]{X} \arrow{sw,t}{\forall\sigma} \arrow{se,t}{\rho}\\
\node{X'}  \arrow[2]{e,b,--}{\exists!\upsilon} \node[2]{E}
\end{diagram}
 \eeq
The object $E$ is also called an {\it envelope} of the object $X$ (in the class of objects ${\tt L}$ with respect to the class of objects ${\tt M}$), and we will use the following notation for it:
 \beq\label{DEF:Env_M^L(A)}
E=\Env_{\tt M}^{\tt L} X.
 \eeq

}\eit

\bit{
\item[$\bullet$] Let us say that a class of objects ${\tt M}$ in a category ${\tt K}$
 {\it differs morphisms on the outside}, if the class of morphisms with ranges in ${\tt M}$ possesses this property  (in the sense of definition on page \pageref{DEF:varPhi-razlich-morfizmy-snaruzhi}), i.e. for any two different parallel morphisms $\alpha\ne\beta:X\to Y$ there is a morphism $\ph:Y\to M\in{\tt M}$ such that $\ph\circ\alpha\ne\ph\circ\beta$.
}\eit

From Theorem \ref{TH:Phi-razdel-moprfizmy-*} we have

\btm\label{TH:M-razdel-moprfizmy} If a class of objects ${\tt M}$ differs morphisms on the outside, then for any class of objects ${\tt L}$
 \bit{
\item[(i)] each envelope in ${\tt L}$ with respect to ${\tt M}$ is a bimorphism,

\item[(ii)] an envelope in ${\tt L}$ with respect to ${\tt M}$ exists if and only if there exists an anvelope in the class of bimorphisms with the values in ${\tt L}$ with respect to ${\tt M}$; in this case these envelopes coincide:
$$
\env_{\tt M}^{\tt L}=\env_{\tt M}^{\Bim({\tt K},{\tt L})}.
$$
 }\eit
\etm

\bex In the category ${\tt Tikh}$ of Tikhonov spaces the {\it Stone-\v{C}ech compactification} $\beta:X\to\beta X$
is an envelope of the space $X$ in the class ${\tt Com}$ of compact spaces with respect to the same class ${\tt Com}$:
$$
\beta X=\Env^{\tt Com}X.
$$
\eex

\bex\label{EX:popolnenie} {\bf Completion} $X^\blacktriangledown$ of a locally convex space $X$ is an envelope of $X$ in the category ${\tt LCS}$ of all locally convex spaces with respect to the class ${\tt Ban}$ of Banach spaces:
$$
X^\blacktriangledown=\Env_{\tt Ban}^{\tt LCS}X.
$$
\eex

\paragraph{Refinements.}

 \bit{
\item[$\bullet$] A morphism $\sigma:X'\to X$ in the category  $\tt K$ is called an {\it enrichment of the object $X\in{\tt K}$ in the class of morphisms $\varGamma$ by means of the class of morphisms $\varPhi$}, if $\sigma\in\varGamma$, and for any morphism $\ph:B\to X$, $\ph\in\varPhi$, there exists a unique morphism $\ph':B\to X'$ in the category  $\tt K$, such that the following diagram is commutative:
 \beq\label{DEF:suzhenie}
\begin{diagram}
\node[2]{X}  \\
\node{B}\arrow{ne,t}{\forall\ph\in\varPhi}\arrow[2]{e,b,--}{\exists!\ph'}
\node[2]{X'}\arrow{nw,t}{\sigma\in\varGamma}
\end{diagram}
 \eeq

\item[$\bullet$] An enrichment $\rho:E\to X$ of the object $X\in{\sf Ob}(\tt K)$ in the class of morphisms $\varGamma$ by means of the class of morphisms $\varPhi$ is called a {\it refinement of $X$ in the class $\varGamma$ by means of $\varPhi$}, if for any other enrichment $\sigma:X'\to X$ (of $X$ in $\varGamma$ by means of  $\varPhi$) there exists a unique morphism $\upsilon:E\to X'$ in $\tt K$, such that the following diagram is commutative:
 \beq\label{DIAGR:otpechatok}
\begin{diagram}
\node[2]{X}  \\
\node{E}\arrow{ne,t}{\rho}\arrow[2]{e,b,--}{\exists!\upsilon}
\node[2]{X'}\arrow{nw,t}{\forall\sigma}
\end{diagram}
 \eeq
For the morphism of refinement $\rho:E\to X$ we use the notation
    \beq\label{DEF:tr_F^L(A)}
    \rho=\rf_{\varPhi}^\varGamma  X.
    \eeq
The very object $E$ is also called a {\it refinement} of $X$ in $\varGamma$ by means of $\varPhi$, and is denoted by
    \beq\label{DEF:Tr_F^L(A)}
E=\Rf_{\varPhi}^\varGamma  X.
 \eeq
}\eit

 \bit{
\item[$\bullet$]\label{DEF:morfizmy-porozhdayutsya-snaruzhi} Let us say that in a category ${\tt K}$ a {\it class of morphisms $\varPhi$ is generated on the outside by a class of morphisms $\varPsi$}, if
    $$
    \varPsi\subseteq\varPhi\subseteq\varPsi\circ\Mor({\tt K}).
    $$
    }\eit
The following fact is dual to Theorem \ref{TH:morfizmy-porozhdayutsya-iznutri}:

\btm\label{TH:morfizmy-porozhdayutsya-snaruzhi} Suppose in a category ${\tt K}$ a class of morphisms $\varPhi$ is generated on the outside by a class of morphisms $\varPsi$. Then for any class of monomorphisms $\varGamma$ (it is not necessary that $\varGamma$ contains all monomorphisms of the category ${\tt K}$) and for any object $X$ the existence of refinement $\rf^\varGamma_\varPsi X$ is equivalent to the existence of the refinement  $\rf^\varGamma_\varPhi X$, and these refinements coincide:
 \beq\label{imp_Psi=imp_Phi}
\rf^\varGamma_\varPsi X=\rf^\varGamma_\varPhi X.
 \eeq
\etm

\bit{
\item[$\bullet$]\label{DEF:varPhi-razlich-morfizmy-iznutri}  Let us say that a class of morphisms $\varPhi$ in a category ${\tt K}$  {\it differs morphisms on the inside}, if for any two different parallel morphisms $\alpha\ne\beta:X\to Y$ there is a morphism $\ph:M\to X$ from the class $\varPhi$ such that  $\alpha\circ\ph\ne\beta\circ\ph$.
}\eit

The following result is dual to Theorem \ref{TH:Phi-razdel-moprfizmy}:

\btm\label{TH:Phi-razdel-moprfizmy-iznutri} If the class of morphisms $\varPhi$ differs morphisms on the inside, then for any class of morphisms $\varGamma$
 \bit{
\item[(i)] every enrichment in $\varGamma$ by means of $\varPhi$ is an epimorphism,

\item[(ii)] the refinement in $\varGamma$ by means of $\varPhi$ exists if and only if there exists a refinement in $\varGamma\cap\Mono$ by means of $\varPhi$; in that case these refinements coincide:
$$
\rf_{\varPhi}^{\varGamma}=\rf_{\varPhi}^{\varGamma\cap\Epi},
$$

\item[(iii)] if the class $\varGamma$ contains all epimorphisms,
 $$
 \varGamma\supseteq\Epi,
 $$
 then the existence of a refinement in $\Epi$ by means of $\varPhi$ automatically implies the existence of a refinement in $\varGamma$ by means of  $\varPhi$, and the coincidence of these refinements:
$$
\rf_{\varPhi}^{\varGamma}=\rf_{\varPhi}^{\Epi}.
$$
 }\eit
\etm

\bit{
\item[$\bullet$] Let us remind that a class of morphisms $\varPhi$ in a category ${\tt K}$ is called a {\it left ideal}, if
$$
\Mor({\tt K})\circ\varPhi\subseteq\varPhi
$$
(i.e. for any $\ph\in\varPhi$ and for any morphism $\mu$ in ${\tt K}$ the composition $\mu\circ\ph$ belongs to $\varPhi$).
}\eit

The following is dual to Theorem \ref{TH:Phi-razdel-moprfizmy-*}

\btm\label{TH:Phi-razdel-moprfizmy-iznutri-*} If a class of morphisms $\varPhi$ differs morphisms on the inside and is a left ideal in the category  ${\tt K}$, then for any class of morphisms $\varGamma$
 \bit{
\item[(i)] every enrichment in $\varGamma$ by means of $\varPhi$ is a bimorphism,

\item[(ii)] a refinement in $\varGamma$ by means of $\varPhi$ exists if and only if there exists a refinement in  $\varGamma\cap\Bim$ by means of $\varPhi$; in that case these refinements coincide:
$$
\rf_{\varPhi}^{\varGamma}=\rf_{\varPhi}^{\varGamma\cap\Bim}.
$$

\item[(iii)] if $\varGamma$ contains all bimorphisms,
 $$
 \varGamma\supseteq\Bim,
 $$
 then a refinement in $\varGamma$ by means of $\varPhi$ exists if and only if there exists a refinement in $\Bim$ by means of $\varPhi$, and these refinements coincide:
$$
\rf_{\varPhi}^{\varGamma}=\rf_{\varPhi}^{\Bim}.
$$
 }\eit
\etm

A special case of refinement is the situation when $\varGamma$ and/or $\varPhi$ are classes of all morphisms from some given subclass of objects in $\Ob({\tt K})$. An exact formulation for the case when both classes $\varGamma$ and $\varPhi$ are defined in this way is the following: suppose we have a category ${\tt K}$ and two subclasses ${\tt L}$ and ${\tt M}$ in the class $\Ob({\tt K})$ of objects of ${\tt K}$.

 \bit{
\item[$\bullet$] A morphism $\sigma:X'\to X$ is called an {\it enrichment of the object $X\in{\sf Ob}(\tt K)$ in the class of objects ${\tt L}$ by means of the class of objects ${\tt M}$}, if for any object $B\in{\tt M}$ and for any morphism $\ph:B\to X$ there is a unique morphism $\ph':B\to X'$ such that the following diagram is commutative:
$$
\put(5,-35){$\begin{matrix}\text{\rotatebox{90}{$\owns$}}\\ {\tt M}
\end{matrix}$} \put(71,-35){$\begin{matrix}\text{\rotatebox{90}{$\owns$}}\\
{\tt L} \end{matrix}$}
\begin{diagram}
\node[2]{X}  \\
\node{B}\arrow{ne,t}{\forall\ph}\arrow[2]{e,b,--}{\exists!\ph'}
\node[2]{X'}\arrow{nw,t}{\sigma}
\end{diagram}
$$

\item[$\bullet$]\label{DEF:nachinka-otn-klassa} An enrichment $\rho:E\to X$ of the object $X\in{\sf Ob}(\tt K)$ in the class of objects ${\tt L}$ by means of the class of objects ${\tt M}$ is called a {\it refinement of the object $X\in{\sf Ob}(\tt K)$ in the class of objects ${\tt L}$ by means of the class of objects ${\tt M}$}, and we write in this case
    \beq\label{DEF:tr_M^L(A)}
    \rho=\rf_{\tt M}^{\tt L} X,
    \eeq
    if for any other enrichment $\sigma:X'\to X$ (of the object $X\in{\sf Ob}(\tt K)$ in the class of objects ${\tt L}$ by means of the class of objects ${\tt M}$) there is a unique morphism $\upsilon:E\to X'$ such that the following diagram is commutative:
 \beq\label{diagr:sled}
\put(5,-35){$\begin{matrix}\text{\rotatebox{90}{$\owns$}}\\ {\tt L}
\end{matrix}$} \put(71,-35){$\begin{matrix}\text{\rotatebox{90}{$\owns$}}\\
{\tt L} \end{matrix}$}
\begin{diagram}
\node[2]{X}  \\
\node{E}\arrow{ne,t}{\rho}\arrow[2]{e,b,--}{\exists!\upsilon}
\node[2]{X'}\arrow{nw,t}{\forall\sigma}
\end{diagram}
 \eeq
The very object $E$ is also called a {\it refinement} of the object $X\in{\sf Ob}(\tt K)$ in the class of objects ${\tt L}$ by means of the class of objects ${\tt M}$, and we use the following notation for it:
    \beq\label{DEF:Tr_M^L(A)}
E=\Rf_{\tt M}^{\tt L} X.
 \eeq
}\eit

\bit{
\item[$\bullet$] Let us say that a class of objects ${\tt M}$ in the category ${\tt K}$  {\it differs morphisms on the inside}, if the class of all morphisms going from objects of ${\tt M}$ has this property (in the sense of definition on page \pageref{DEF:varPhi-razlich-morfizmy-iznutri}), i.e. for any two different parallel morphisms $\alpha\ne\beta:X\to Y$ there is a morphism $\ph:M\to X$ such that $\alpha\circ\ph\ne\beta\circ\ph$.
}\eit

Theorem \ref{TH:Phi-razdel-moprfizmy-iznutri-*} implies

\btm\cite[Theorem 3.19]{Akbarov-env}\label{TH:M-razdel-moprfizmy-iznutri} If a class of objects ${\tt M}$ differs morphisms on the inside, then for any class of objects ${\tt L}$
 \bit{
\item[(i)]  each domain of convergence in the class ${\tt L}$ by means of the class ${\tt M}$ is a bimorphism,

\item[(ii)] a refinement in the class ${\tt L}$ by means of the class ${\tt M}$ exists if and only if there exists a refinement in the class of bimorphisms going from ${\tt L}$ by means of the class ${\tt M}$; in that case these refinements coincide:
$$
\rf_{\tt M}^{\tt L}=\rf_{\tt M}^{\Bim({\tt L},{\tt K})}.
$$
 }\eit
\etm

\bex {\bf Simply connected covering} used in the theory of Lie groups is from the categorial point of view a refinement in the class of pointed simply connected coverings by means of empty class of morphisms in the category of connected locally connected and semilocally simply connected pointed topological spaces
(see definitions in \cite{Postnikov}). \eex

\bex {\bf Bornologification} (see definition in \cite{Kriegl-Michor}) $X_{\text{\rm born}}$ of a locally convex space $X$ is a refinement of $X$ in the category ${\tt LCS}$ of locally convex spaces by means of the subcategory ${\tt Norm}$ of normed spaces:
$$
X_{\text{\rm born}}=\Rf_{\tt Norm}^{\tt LCS}X
$$
\eex

\bex {\bf Saturation} $X^\blacktriangle$ of a pseudocomplete locally convex space $X$ is a refinement in the category ${\tt LCS}$ of locally convex spaces in its object $X$ by means of the subcategory ${\tt Smi}$ of the Smith spaces (see definitions in \cite{Akbarov}):
$$
X^\blacktriangle=\Rf_{\tt Smi}^{\tt LCS}X
$$
 \eex

\paragraph{Connection with nodal decomposition.}

Envelopes and refinements are connected to nodal decomposition through a series results. We mention the two shortest of them.

\btm\cite[pp.64, 65]{Akbarov-env} Let ${\tt K}$ be a category with nodal decomposition. Then
\bit{
\item[(i)] if ${\tt K}$ is a category with products, and is co-well-powered in the class $\Epi$, then in ${\tt K}$ each object $X$ has an envelope in the class $\Epi$ with respect to an arbitrary class of morphisms $\varPhi$, going from $X$,

\item[(ii)] if ${\tt K}$ is a category with coproducts, and is well-powered in the class $\Mono$, then in ${\tt K}$ each object $X$ has a refinement in the class $\Mono$ by means of the class of morphisms $\varPhi$, going to $X$.

}\eit
\etm

Let us say that in a category ${\tt K}$
 \bit{
 \item[---]\label{DEF:epi-raspoznayut-mono}  {\it epimorphisms discern monomorphisms}, if from the fact that a morphism $\mu$ is not a monomorphism it follows that $\mu$ can be represented as a composition $\mu=\mu'\circ\e$, where $\e$ is an epimorphism, which is not an isomorphism,

 \item[---] {\it monomorphisms discern epimorphisms}, if from the fact that a morphism $\e$ is not an epimorphism it follows that $\e$ can be represented as a composition $\e=\mu\circ\e'$, where $\mu$ is a monomorphism, which is not an isomorphism.
 }\eit

\btm\cite[Theorem 3.31]{Akbarov-env}\label{TH:env+imp=>uzl=razl}  Suppose that in a category ${\tt K}$
 \bit{
 \item[\rm (a)] epimorphisms discern monomorphisms, and, dually, monomorphisms discern epimorphisms,

 \item[\rm (b)] every immediate monomorphism is a strong monomorphism, and, dually, every immediate epimorphism is a strong epimorphism,

 \item[\rm (c)] every object $X$ has an envelope in the class $\Epi$ of all epimorphisms with respect to any morphism, starting from $X$, and, dually, in every object $X$ there is a refinement in the class $\Mono$ of all monomorphisms with respect to any morphism coming to $X$.
   }\eit
Then ${\tt K}$ is a category with nodal decomposition.
 \etm

\subsection{Functoriality}

\bit{
    \item
Let us say that {\it the envelope $\Env^\varOmega_\varPhi$ can be defined as a functor}, if there exist
        \bit{
\item[E.1] a map $X\mapsto (E(X),e_X)$, that to each object $X$ in $\tt K$ assigns a morphism $e_X:X\to E(X)$ in ${\tt K}$, which is an envelope in $\varOmega$ with respect to $\varPhi$:
$$
E(X)=\Env_{\varPhi}^\varOmega X,\qquad e_X=\env_{\varPhi}^\varOmega X
$$
    \item[E.2] a map $\alpha\mapsto E(\alpha)$, that each morphism $\alpha:X\to Y$ in $\tt K$ turns into a morphism  $E(\alpha):E(X)\to E(Y)$ in $\tt K$ in such a way that the following diagram is commutative
    \beq\label{DIAGR:funktorialnost-env-e-E}
\xymatrix @R=2.pc @C=5.0pc % @M=14pt
{
X\ar[d]^{\alpha}\ar[r]^{e_X} & E(X)\ar@{-->}[d]^{E(\alpha)} \\
Y\ar[r]^{e_Y} & E(Y) \\
}
\eeq
}\eit
and the following identities hold
\beq\label{tozhdestva:funktorialnost-env-e-E}
E(1_X)=1_{E(X)},\qquad E(\beta\circ\alpha)=E(\beta)\circ E(\alpha)
\eeq
Clearly, in this case the map $(X,\alpha)\mapsto(E(X),E(\alpha))$ is a covariant functor from ${\tt K}$ into ${\tt K}$, and the map $X\mapsto e_X$ is a natural transformation of the identity functor $(X,\alpha)\mapsto(X,\alpha)$ into the functor $(X,\alpha)\mapsto(E(X),E(\alpha))$.

    \item Let us say that {\it the envelope $\Env^\varOmega_\varPhi$ can be defined as an idempotent functor}, if in addition to E.1 and E.2 one can ensure the condition

 \bit{
 \item[E.3] for each object $X\in\Ob({\tt K})$ the morphism $e_{E(X)}:E(X)\to E(E(X))$ is the local identity:
\beq\label{e_(E(X))=1_(E(X))}
E(E(X))=E(X),\qquad e_{E(X)}=1_{E(X)}\qquad X\in \Ob({\tt K}).
\eeq
 }\eit
 }\eit

\brem\label{REM:E(e_X)=1_(E(X))}
If $\varOmega\subseteq\Epi$, then \eqref{e_(E(X))=1_(E(X))} implies
\beq\label{E(e_X)=1_(E(X))}
E(e_X)=1_{E(X)}\qquad X\in \Ob({\tt K}).
\eeq
Indeed, if we put $\alpha=e_X$ into \eqref{DIAGR:funktorialnost-env-e-E}, we obtain
$$
\xymatrix @R=2.pc @C=8.0pc % @M=14pt
{
X\ar[d]^{e_X}\ar[r]^{e_X} & E(X)\ar[d]^{E(e_X)} \\
E(X)\ar[r]^{e_{E(X)}=1_{E(X)}} & E(E(X))=E(X) \\
}
$$
i.e. $E(e_X)\circ e_X=1_{E(X)}\circ e_X$, and, since $e_X\in\varOmega\subseteq\Epi$, we can cancel it: $E(e_X)=1_{E(X)}$.
\erem

\paragraph{Nets of epimorphisms.}

\bit{ \item[$\bullet$] Suppose that to each object $X\in\Ob({\tt K})$ in a category ${\tt K}$ it is assigned a subset ${\mathcal N}^X$ in the class
$\Epi^X$ of all epimorphisms of the category ${\tt K}$, going from $X$, and the following three requirements are fulfilled:
    \bit{
\item[(a)]\label{AX:set-Epi-a} for each object $X$ the set ${\mathcal N}^X$ is non-empty and is directed to the left with respect to the pre-order \eqref{DEF:le-in-F_X} inherited from $\Epi^X$:
$$
\forall \sigma,\sigma'\in {\mathcal N}^X\quad \exists\rho\in{\mathcal N}^X\quad
\rho\to\sigma\ \& \ \rho\to\sigma',
$$

\item[(b)]\label{AX:set-Epi-b} for each object $X$ the covariant system of morphisms generated by ${\mathcal N}^X$  \beq\label{DEF:kategoriya-svyazyv-morpfizmov}
    \Bind({\mathcal N}^X):=\{\iota_\rho^\sigma;\ \rho,\sigma\in{\mathcal N}^X,\ \rho\to\sigma\}
    \eeq
    (the morphisms $\iota_\rho^\sigma$ were defined in \eqref{DEF:le-in-F_X-*}; by \eqref{iota_rho^tau=iota_rho^sigma-circ-iota_sigma^tau} this system is a covariant functor from the set ${\mathcal N}^X$ considered as a full subcategory in $\Epi^X$ into $\tt K$) has a projective limit in $\tt K$;

\item[(c)]\label{AX:set-Epi-c}  for each morphism $\alpha:X\to Y$ and for each element $\tau\in{\mathcal N}^Y$ there are an element $\sigma\in{\mathcal N}^X$ and a morphism $\alpha_\sigma^\tau:\Ran\sigma\to\Ran\tau$ such that the following diagram is commutative
    \beq\label{DIAGR:set} \xymatrix @R=2.5pc @C=4.0pc {
 X\ar[r]^{\alpha}\ar@{-->}[d]_{\sigma} & Y\ar[d]^{\tau} \\
 \Ran\sigma\ar@{-->}[r]_{\alpha_\sigma^\tau} & \Ran\tau
    } \eeq
    (a remark: for given $\alpha$, $\sigma$ and $\tau$ the morphism $\alpha_\sigma^\tau$, if exists, must be unique, since $\sigma$ is an epimorphism).

 }\eit
Then
 \bit{
\item[---] we call the family of set ${\mathcal N}=\{{\mathcal N}^X;\ X\in\Ob({\tt K})\}$ a {\it net of epimorphisms}\label{DEF:set-epimorf} in the category ${\tt K}$, and the elements of the sets ${\mathcal N}^X$ {\it elements of the net} ${\mathcal N}$,

\item[---] for each object $X$ the system of morphisms $\Bind({\mathcal N}^X)$ defined by equalities \eqref{DEF:kategoriya-svyazyv-morpfizmov} will be called the {\it system of binding morphisms of the net ${\mathcal N}$ over the vertex $X$}, its projective limit (which exists by condition (b)) is a projective cone, whose vertex will be denoted by $X_{\mathcal N}$, and the morphisms going from it by $\sigma_{\mathcal N}=\leftlim_{\rho\in{\mathcal N}^X} \iota^\sigma_\rho: X_{\mathcal N}\to \Ran\sigma$:
 \beq\label{X_F-proektiv-sistema} \xymatrix @R=2.5pc @C=2.0pc {
 & X_{\mathcal N}\ar[dr]^{\sigma_{\mathcal N}}\ar[dl]_{\rho_{\mathcal N}} &  \\
 \Ran\rho\ar[rr]^{\iota^\sigma_\rho} & &  \Ran\sigma
}\qquad (\rho\to\sigma);
 \eeq
in addition, by \eqref{DEF:le-in-F_X-*}, the system of epimorphisms ${\mathcal N}^X$ is also a projective cone of the system $\Bind({\mathcal N}^X)$:
 \beq\label{F_X-proektiv-sistema} \xymatrix @R=2.5pc @C=2.0pc {
 & X\ar[dr]^{\sigma}\ar[dl]_{\rho} &  \\
 \Ran\rho\ar[rr]^{\iota^\sigma_\rho} & &  \Ran\sigma
}\qquad (\rho\to\sigma),
 \eeq
 so there must exist a natural morphism from $X$ into the vertex $X_{\mathcal N}$ of the projective limit of the system $\Bind({\mathcal N}^X)$. We  denote this morphism by $\leftlim{\mathcal N}^X$ and call it the {\it local limit of the net ${\mathcal N}$ of epimorphisms at the object $X$}: \beq\label{DIAGR:sigma-sigma_F}
\xymatrix @R=2.5pc @C=2.0pc {
   X\ar[dr]_{\sigma}\ar@{-->}[rr]^{\leftlim{\mathcal N}^X} & & X_{\mathcal N}\ar[ld]^{\sigma_{\mathcal N}}  \\
   & \Ran\sigma  &
}\qquad (\sigma\in{\mathcal N}^X). \eeq

\item[---] the element $\sigma$ of the net in diagram \eqref{DIAGR:set} will be called a {\it counterfort} of the element $\tau$ of the net.
 }\eit
 }\eit

\btm\cite[Theorem 3.38]{Akbarov-env}\label{TH:funktorialnost-pri-seti-Epi-i-dolonyaemosti}
Let ${\mathcal N}$ be a net of epimorphisms in a category ${\tt K}$, that generates a class of morphisms $\varPhi$ on the inside:
    $$
    {\mathcal N}\subseteq\varPhi\subseteq\Mor({\tt K})\circ {\mathcal N}.
    $$
Then for each monomorphically complementable\footnote{See definition on p.\pageref{DEF:klass-monomorfno-dopolnyaem}.} class of epimorphisms $\varOmega$,
$$
{^\downarrow\varOmega}\circledcirc\varOmega={\tt K},
$$
the following holds:
 \bit{
\item[(a)] for each object $X$ in ${\tt K}$ the morphism $\e_{\leftlim{\mathcal N}^X}$ in the factrization  \eqref{faktorizatsiya-v-kat-s-faktoriz} defined by the classes ${^\downarrow\varOmega}$ and $\varOmega$, is an envelope  $\env_\varPhi^{\varOmega} X$ in $\varOmega$ with respect to $\varPhi$:
\beq\label{im_infty-lim-F_X=env_M^Epi-X-1}
\e_{\leftlim{\mathcal N}^X}=\env_\varPhi^{\varOmega} X,
\eeq

\item[(b)] for each morphism $\alpha:X\to Y$ in ${\tt K}$ and for any choice of envelopes $\env_\varPhi^{\varOmega}X$ and $\env_\varPhi^{\varOmega}Y$ there exists a unique morphism $\Env_\varPhi^{\varOmega} \alpha:\Env_\varPhi^{\varOmega}  X\to \Env_\varPhi^{\varOmega}Y$ in ${\tt K}$ such that the following diagram is commutative:
\beq\label{DIAGR:funktorialnost-env_varPhi^Epi-v-kat-s-uzl-razl-1}
\xymatrix @R=2.pc @C=5.0pc % @M=14pt
{
X\ar[d]^{\alpha}\ar[r]^{\env_\varPhi^{\varOmega}  X} & \Env_\varPhi^{\varOmega}  X\ar@{-->}[d]^{\Env_\varPhi^{\varOmega} \alpha} \\
Y\ar[r]^{\env_\varPhi^{\varOmega}  Y} & \Env_\varPhi^{\varOmega}  Y \\
}
\eeq

\item[(c)] if in addition $\tt K$ is co-well-powered in the class $\varOmega$, then the envelope  $\Env_\varPhi^{\varOmega}$ can be defined as a functor.
}\eit
 \etm

\bex
Let $X$ be a locally convex space. To each closed convex balanced neighbourhood of zero $U$ in $X$ let us assign the closed subspace $\Ker U=\bigcap_{\varepsilon>0}\varepsilon\cdot U$ in $X$ and the quotient space $X/\Ker U$, which we endow with (not the quotient topology, but) the topology of normed space with the unit ball $U+\Ker U$. Then the completion $(X/\Ker U)^\blacktriangledown$ is a Banach space. We denote it by $X/U$ and call it the {\it Banach quotient space of the space $X$ by the neighbourhood of zero $U$}\label{DEF:X/U}. The natural map from $X$ into  $X/U$
$$
 \xymatrix  @R=2.5pc @C=4pc
 {
X \ar@/^5ex/[rr]^{\rho_U}\ar[r]_{\tau_U}  & X/\Ker U \ar[r]_{\blacktriangledown_{X/\Ker U}}
 & (X/\Ker U)^\blacktriangledown=X/U
 }
$$
(where $\tau_U$ is the quotient map, and $\blacktriangledown_{X/\Ker U}$ the completion map) will be called the {\it Banach quotient map of the space $X$ by the neighbourhood of zero $U$}. We denote by $\mathcal B$ the clas of all Banach quotient maps $\{\rho_U:X\to X/U\}$, where $X$ runs over the class of locally convex spaces, and $U$ the class of all closed convex balanced neighbourhoods of zero in $X$.

{\it The class $\mathcal B$ of Banach quotient maps is a net of epimorphisms in the category ${\tt LCS}$ of locally convex spaces, and the pre-order\footnote{The pre-order $\to$ in the class $\Epi^X$ was defined on page \pageref{DEF:le-in-F_X}.}
$\to$  in $\mathcal B$ is equivalent to the embedding of the corresponding neighbourhoods of zero, up to a positive scalar multiplier:}
 \beq\label{rho_V-to-rho_U-<=>-V-subseteq-U}
\rho_V\to\rho_U\quad\Longleftrightarrow\quad\exists\e>0\quad \e\cdot V\subseteq
U.
 \eeq
Certainly, the class $\mathcal B$ generates on the inside the class $\Mor({\tt LCS},{\tt Ban})$ of all morphisms into Banach spaces. This implies that the completion map $X\mapsto X^\blacktriangle$, mentioned in the Example \ref{EX:popolnenie} as an envelope in ${\tt LCS}$ with respect to the class ${\tt Ban}$ of Banach spaces, coincides with the envelope with respect to the class $\mathcal B$ of Banach quotient maps:
$$
X^\blacktriangledown=\Env_{\tt Ban}^{\tt LCS}X=\Env_{\mathcal B}^{\tt LCS}X.
$$
\eex

\paragraph{Regular envelopes.}

\bit{
\item Let us say that {\it a class of morphisms $\varOmega$ pushes a class of morphisms $\varPhi$}, if
\beq\label{Omega-podderzhivaet-Phi}
\forall\psi\in\Mor({\tt K})\qquad \forall \sigma\in\varOmega\qquad \big(\psi\circ\sigma\in\varPhi\quad\Longrightarrow\quad \psi\in\varPhi\big).
\eeq
}\eit

The following result complements Theorem \ref{TH:funktorialnost-pri-seti-Epi-i-dolonyaemosti}.

\btm\cite[Theorems 3.42, 3.48]{Akbarov-env}\label{TH:reg-obolochka} Suppose a category $\tt K$ and classes of morphisms $\varOmega$ and $\varPhi$ in it satisfy the following conditions:
\bit{

\item[RE.1:] $\tt K$ is projectively complete,

\item[RE.2:] $\varOmega$ is monomorphically complementable: ${^\downarrow\varOmega}\circledcirc\varOmega={\tt K}$,

\item[RE.3:] ${\tt K}$ is co-well-powered in the class $\varOmega$,

\item[RE.4:] $\varPhi$ goes from\footnote{In the sense of definition on p.\pageref{DEF:goes-from}.} $\Ob({\tt K})$ and is a right ideal in ${\tt K}$:
$$
\Ob(\tt K)=\{\Dom\ph;\ \ph\in\varPhi\},\qquad \varPhi\circ\Mor({\tt K})\subseteq\varPhi.
$$
\item[RE.5:] the class $\varOmega$ pushes the class $\varPhi$.

}\eit
Then
 \bit{
\item[(a)] there is a net of epimorphisms ${\mathcal N}$ in ${\tt K}$ such that for each object $X$ in ${\tt K}$ the morphism $\e_{\leftlim{\mathcal N}^X}$ in the factorization \eqref{faktorizatsiya-v-kat-s-faktoriz} is an envelope $\env_\varPhi^{\varOmega} X$ in $\varOmega$ with respect to $\varPhi$:
\beq\label{im_infty-lim F_X=env_M^Epi-X-1-*}
\e_{\leftlim{\mathcal N}^X}=\env_\varPhi^{\varOmega} X,
\eeq

\item[(b)] for each morphism $\alpha:X\to Y$ in ${\tt K}$ and for any choice of envelopes   $\env_\varPhi^{\varOmega}X$ and $\env_\varPhi^{\varOmega}Y$ there exists a unique morphism  $\Env_\varPhi^{\varOmega} \alpha:\Env_\varPhi^{\varOmega}  X\to \Env_\varPhi^{\varOmega}Y$ in ${\tt K}$, such that the following diagram is commutative:
\beq\label{DIAGR:funktorialnost-env_varPhi^Epi-v-kat-s-uzl-razl-1-*}
\xymatrix @R=2.pc @C=5.0pc % @M=14pt
{
X\ar[d]^{\alpha}\ar[r]^{\env_\varPhi^{\varOmega}  X} & \Env_\varPhi^{\varOmega}  X\ar@{-->}[d]^{\Env_\varPhi^{\varOmega} \alpha} \\
Y\ar[r]^{\env_\varPhi^{\varOmega}  Y} & \Env_\varPhi^{\varOmega}  Y \\
}
\eeq

\item[(c)] the envelope $\Env_\varPhi^{\varOmega}$ can be defined as an idempotent functor.
}\eit\etm

\bit{

\item If the conditions RE.1-RE.5 hold, then we say that {\it the classes of morphisms $\varOmega$ and $\varPhi$ define a regular envelope in the category $\tt K$}\label{DEF:reg-obolochka}, or that {\it the envelope $\Env_\varPhi^\varOmega$ is regular}.
}\eit

\paragraph{Envelopes coherent with tensor product.}

Let $\tt K$ be a minoidal category \cite{MacLane} with the tensor product $\otimes$ and the unit object $I$.
 \bit{\label{DEF:obolochka-soglasovana-s-tenz-proizv}
\item We say that {\it the envelope $\Env^\varOmega_\varPhi$ is coherent with the tensor product $\otimes$} in $\tt K$, if the following conditions are fulfilled:

 \bit{

\item[T.1]\label{DEF:T.1} The tensor product $\rho\otimes\sigma:X\otimes Y\to X'\otimes Y'$ of any two extensions  $\rho:X\to X'$ and $\sigma:Y\to Y'$ (in $\varOmega$ with respect to $\varPhi$) is an extension (in $\varOmega$ with respect to $\varPhi$).

\item[T.2]\label{DEF:T.2} The local identity $1_I:I\to I$ of the unit object $I$ is the envelope (in $\varOmega$ with respect to $\varPhi$):
\beq\label{E(I)=I}
\env_\varPhi^\varOmega I=1_I
\eeq
 }\eit
 }\eit

 \bit{
\item An object $A$ in a category $\tt K$ is said to be {\it complete} in the class of morphisms $\varOmega\subseteq\Epi$ with respect to the class of morphisms $\varPhi$, if there exists an envelope of $A$ in $\varOmega$ with respect to $\varPhi$, which is an isomorphism: $\env^\varOmega_\varPhi A\in\Iso$.
 }\eit

Let $\Env_\varPhi^\varOmega$ be a regular envelope, coherent with the tensor product in $\tt K$,
$E=\Env_\varPhi^\varOmega$ the idempotent functor built in Theorem \ref{TH:reg-obolochka}, and ${\tt L}$ the (full) subcategory of complete objects in $\tt K$. For any objects $A,B\in{\tt L}$ and for any morphisms
$\ph,\psi\in{\tt L}$ we put
\beq\label{A-overset(E)(otimes)B:=E(A-otimes-B)}
A\overset{E}{\otimes}B:=E(A\otimes B),\qquad \ph\overset{E}{\otimes}\psi:=E(\ph\otimes\psi).
\eeq
The following identity holds:
\beq\label{E(E(X)-otimes-E(Y))=E(X)-overset(E)(otimes)E(Y)}
E(X)\overset{E}{\otimes}E(Y)=E(E(X)\otimes E(Y)),\qquad X,Y\in \Ob(\tt K).
\eeq

\btm\cite[Theorem 3.63]{Akbarov-env}\label{TH:sushestvovanie-tenz-proizv-v-L}
Suppose $\Env_\varPhi^\varOmega$ is a regular envelope, coherent with the tensor product in $\tt K$. Then the formulas \eqref{A-overset(E)(otimes)B:=E(A-otimes-B)} define a structure of monoidal category on $\tt L$ (with  $\overset{E}{\otimes}$ as tensor product and $I$ as unit object), and the functor of envelope $E:{\tt K}\to {\tt L}$ defined in Theorem \ref{TH:reg-obolochka}, is monoidal.
\etm

\bcor\cite[Corollary 3.65]{Akbarov-env}\label{COR:Env-sohranyaet-algebry}
Suppose $\Env_\varPhi^\varOmega$ is a regular envelope coherent with the tensor product in $\tt K$.
The operation $\Env_\varPhi^\varOmega$ turns each algebra (respectively, coalgebra, bialgebra, Hopf algebra) $A$ in $\tt K$ into an algebra (respectively, coalgebra, bialgebra, Hopf algebra) $\Env_\varPhi^\varOmega A$ in $\tt L$.
\ecor

\section{Stereotype spaces}\label{SEC:kateg-Ste}

The theory of stereotype spaces was developed in the series of author's works, of which the main papers are \cite{Akbarov}, \cite{Akbarov-stein-groups} and \cite{Akbarov-env}. We give here some definitions and facts of the theory. We mostly refer the reader to the proofs in \cite{Akbarov} and \cite{Akbarov-env}, but some propositions are new, in particular, the results on stereotype spaces with involution, and we prove them in this text. We also prove  the properties of tensor products in Theorem \ref{TH:svoistva-tenz-proizvedenij} which were mentioned in \cite{Akbarov} and \cite{Akbarov-env} without proofs.

\subsection{Pseudocompletion and pseudosaturation}

\paragraph{Pseudocompleteness and pseudocompletion.}

A set $S$ in a locally convex space $X$ is said to be {\it totally bounded} (or {\it precompact})
\cite{Schaefer}, if for each neighbourhood of zero $U$ in $X$ there is a finite set $A$ such that the shifts of $U$ by the elements of $A$ cover $S$:
$$
S\subseteq U+A.
$$
This is equivalent to the total boundedness of $S$ in the sense of the uniform structure induced from $X$  \cite{Engelking} (i.e. {\it $A$ can be chosen as a subset in $S$}).

\bit{
 \item
A locally convex space $X$ is said to be {\it pseudocomplete}, if every totally bounded Cauchy net in $X$ converges. This is equivalent to the claim that every closed totally bounded set in $X$ is compact.
}\eit
This is connected with the usual completeness and quasicompleteness\footnote{A locally convex space $X$ is said to be {\it quasicomplete}, if every bounded Cauchy net in $X$ converges.} by the implications
\beq\label{implikatsii-dlya-psevdopolnoty}
\text{$X$ is complete $\Longrightarrow$ $X$ is quasicomplete $\Longrightarrow$ $X$ is pseudocomplete}
\eeq
In the metrizable case these properties are equivalent.

Like in the case of completeness, each locally convex space $X$ has pseudocompletion, i.e. the ``outside-nearest'' pseudocomplete space. Formally this construction is described in the following

\btm \label{th-1.5} There exists a map $X\mapsto\triangledown_X$ that assigns to each locally convex space $X$ a linear continuous map ${\triangledown}_X :X \to X^{\triangledown}$ into a pseudocomplete locally convex space
$X^{\triangledown}$ in such a way that the following conditions are fulfilled:
\begin{itemize}
\item[(i)] $X$ is pseudocomplete if an only if ${\triangledown}_X :X\to X^{\triangledown}$ is an isomorphism;

\item[(ii)] for any linear continuous map $\ph :X\to Y$ of locally convex spaces there is a unique linear continuous map $\ph^{\triangledown} :X^{\triangledown} \to Y^{\triangledown}$ such that the following diagram is commutative:
\beq
\begin{diagram}
\node{X} \arrow{e,t}{\triangledown_X} \arrow{s,l}{\varphi}
\node{X^\triangledown} \arrow{s,r,--}{\varphi^\triangledown}
\\
\node{Y} \arrow{e,t}{\triangledown_Y} \node{Y^\triangledown}
\end{diagram}.
\label{eqI.12}
\eeq \end{itemize}
\etm

From $(i),(ii)$ it follows that for any linear continuous map $\ph :X\to Y$ into a pseudocomplete space $Y$ there exists a unique linear continuous map $X^{\triangledown} \to Y$ such that the following diagram is commutative:
\beq\label{eqI.13}
\begin{diagram}
\node{X} \arrow[2]{e,t}{\triangledown_X} \arrow{se,b}{\varphi}
\node[2]{X^\triangledown} \arrow{sw,--}{}
\\
\node[2]{Y}
\end{diagram}.
\eeq
This means by the way that, the morphism $\triangledown_X:X\to X^\triangledown$ is an extension of $X$ in $\Ob({\tt LCS})$ with respect to the object $\C$. Since $\C$ differs morphisms on the outside in ${\tt LCS}$, by Theorem \ref{TH:M-razdel-moprfizmy}, $\triangledown_X:X\to X^\triangledown$ is a bimorphism. This implies in its turn that the morphism  $\triangledown_X:X\to X^\triangledown$ is unique up to an isomorphism in $\Epi^X$.

 \bit{
 \item
The space $X^{\triangledown}$ is called the {\it pseudocompletion}, and the map ${\triangledown}_X :X \to
X^{\triangledown}$ the {\it pseudocompletion map} of the locally convex space $X$. From $(ii)$, it follows also, that the map $\ph \mapsto \ph ^{\triangledown}$ is a covariant functor of the category ${\tt LCS}$ into itself: $(\psi \circ \ph)^{\triangledown}=\psi^{\triangledown} \circ \ph^{\triangledown}$. We call it the {\it pseudocompletion functor}.
 }\eit

\begin{itemize}
\item
A linear continuous map $\ph:X\to Y$ of locally convex spaces will be called
\begin{itemize}
\item[--] an {\it embedding} (respectively, a {\it weak embedding}, a {\it relative embedding}), if it is injective and open (respectively, weakly open, relatively open),

\item[--] a {\it dense embeddding} (respectively, a {\it dense weak embedding}, a {\it dense relative embedding}), if in addition the set of values $\ph(X)$ is dense in $Y$.
\end{itemize}
\end{itemize}

\btm \label{th-1.6} For any locally convex space $X$ the pseudocompletion map
${\triangledown}_X :X \to X^{\triangledown}$ is a dense embedding.
\etm

Like usual completion, the operation of pseudocompletion $X\mapsto X^\triangledown$ adds new elements to $X$, but does not change the topology of $X$.

\paragraph{Pseudosaturateness and pseudosaturation.}

A set $D$ in a locally convex space $X$ is called {\it capacious}, if for any totally bounded set $S\subseteq X$ there is a finite set $A\subseteq X$ such that the shifts of $D$ by the elements of $A$ cover $S$:
$$
S\subseteq D+A.
$$
Note that {\it if $D$ is convex, then $A$ can be chosen as a subset in $S$} (and this gives an equivalent condition on $D$).

\bit{
\item A locally convex space $X$ is said to be {\it pseudosaturated}, if for every closed convex balanced capacious set $D$ in $X$ is a neighbourhood of zero in $X$.
}\eit
In the theory of topological vector spaces this property is connected with the metrizability and barreledness:
\beq\label{implikatsii-dlya-psevdonasyshennosti}
\text{$X$ is metrizable $\Longrightarrow$ $X$ is barreled $\Longrightarrow$ $X$ is pseudosaturated.}
\eeq

\btm[criterion of being pseudosaturated]\label{th-2.5}
For a locally convex space $X$ the following conditions are equivalent:
\begin{itemize}
\item[(i)] $X$ is pseudosaturated,

\item[(ii)] if a set of linear continuous functionals $F\subseteq X'$ is equicontinuous on each totally bounded set $S\subseteq X$, then $F$ is equicontinuous on $X$.

\item[(iii)] if $Y$ is a locally convex space and $\varPhi$ is a set of linear continuous maps $\ph:X\to Y$, equicontinuous on each totally bounded set $S\subseteq X$, then $\varPhi$ is equicontinuous on $X$.
\end{itemize}
\etm

Remarkably, there exists a construction, dual to the construction of pseudocompletion, that assigns to each locally convex space $X$ an ``inside-nearest'' pseudosaturated locally convex space $X^\vartriangle$:

\btm \label{th-1.16} There exists a map $X\mapsto\vartriangle_X$, that assigns to each locally convex space $X$ a linear continuous map $\vartriangle _X :X^\vartriangle \to X$ from a pseudosaturated locally convex space $X^\vartriangle$ in such a way that the following conditions are fulfilled:
\begin{itemize}
\item[(i)] $X$ is pseudosaturated if and only if $\vartriangle _X:X^\vartriangle \to X$ is an isomorphism;

\item[(ii)] for any linear continuous map $\ph:Y\to X$ of locally convex spaces there is a unique linear continuous map $\ph^\vartriangle :Y^\vartriangle\to X^\vartriangle$ such that the following diagram is commutative:
\beq
\begin{diagram}
\node{X} \node{X^\vartriangle} \arrow{w,t}{\vartriangle_X}
\\
\node{Y} \arrow{n,l}{\varphi} \node{Y^\vartriangle} \arrow{w,t}{\vartriangle_Y}
\arrow{n,r,--}{\varphi^\vartriangle}
\end{diagram}.
\label{eqI.25}
\eeq \end{itemize}
\etm

From $(i),(ii)$ it follows that for any linear continuous map $\ph :Y\to X$ from a pseudosaturated locally convex space $Y$ there is a unique linear continuous map $Y\to X^{\vartriangle}$ such that the following diagram is commutative:
\beq\label{eqI.26}
\begin{diagram}
\node{X} \node[2]{X^\vartriangle} \arrow[2]{w,t}{\vartriangle_X}
\\
\node[2]{Y} \arrow{nw,b}{\varphi} \arrow{ne,b,--}{}
\end{diagram}.
\eeq
This means by the way that the morphism $\vartriangle_X:X^\vartriangle\to X$ is an enrichment of $X$ in the class $\Ob({\tt LCS})$ by means of the object $\C$. Since $\C$ differs morphisms on the inside in ${\tt LCS}$, by Theorem \ref{TH:M-razdel-moprfizmy-iznutri}, $\vartriangle_X:X^\vartriangle\to X$ is a bimorphism. This implies in its turn that the morphism $\vartriangle_X:X^\vartriangle\to X$ is unique up to an isomorphism in $\Mono_X$.

\bit{\label{DEF:pseudosaturation}
\item  The space $X^\vartriangle$ is called the {\it pseudosaturation}, and the map $\vartriangle _X :X^\vartriangle \to X$ the {\it pseudosaturation map} of the space $X$. From $(ii)$ it follows that the map $\ph \mapsto \ph ^{\vartriangle}$ is a covariant functor of the category ${\tt LCS}$ into itself: $(\psi\circ\ph)^{\vartriangle}=\psi^{\vartriangle} \circ \ph^{\vartriangle}$. We call it {\it pseudosaturation functor}.
}\eit

\begin{itemize}
\item
A linear continuous map $\ph:X\to Y$ of locally convex spaces will be called
\begin{itemize}
\item[--] a {\it covering} (respectively, a {\it weak covering}, a {\it relative covering}), if it is surjective and closed (respectively, weakly closed, relatively closed),

\item[--] an {\it exact covering} (respectively, an {\it exact weak covering}, an {\it exact relative covering}), if in addition it is injective.
\end{itemize}
\end{itemize}

\brem
If a space $X$ is pseudocomplete and $\ph:X\to Y$ is an exact covering, then for each totally bounded set $S\subseteq X$ the restriction $\ph|_S:S\to\ph(S)$ is a homeomorphism of topological spaces.
\erem

\btm\label{th-1.17} For any locally convex space $X$ the pseudosaturation map $\vartriangle _X :X^\vartriangle \to X$ is an exact covering.
\etm

The pseudosaturation $X^\vartriangle$ can be imagined as a new, stronger topologization of the space $X$, which preserves the system of totally bounded sets and the topology on each totally bounded set in $X$.

\paragraph{Independence and consistency.}
The following examples show that pseudocompleteness and pseudosaturateness are independent conditions.

\begin{ex}\label{ex-4.9} Let $X$ be an infinite-dimensional Banach space, and $Y=X'_\sigma$ its dual space with the $X$-weak topology. The space $Y$ is pseudocomplete, but not pseudosaturated.
\end{ex}

\begin{ex}\label{ex-4.10} An arbitrary non-complete metrizable locally convex space is pseudosaturated, but not pseudocomplete space.
\end{ex}

However, the operations $X\mapsto X^\triangledown$ and $X\mapsto X^\vartriangle$ are consistent in the following sense:

\btm\label{TH:X-psevdopolno=>X^vartriangle-psevdopolno} For a locally convex space $X$
\bit{
\item[---] if $X$ is pseudocomplete, then its pseudosaturation $X^\vartriangle$ is also pseudocomplete,

\item[---] if $X$ is pseudosaturated, then its pseudocompletion $X^\triangledown$ is also pseudosaturated.
}\eit
\etm

\brem
It remains open, if these operations commute:
$$
(X^\triangledown)^\vartriangle\overset{?}{=}(X^\vartriangle)^\triangledown.
$$
\erem

\paragraph{Duality between pseudocompleteness and pseudosaturateness.}

Let $X$ be a locally convex space over the field of complex numbers $\C$. Denote by $X^\star$\label{DEF:X^star} the set of lineat continuous functionals $f:X\to \C$ endowed by the {\it topology of uniform convergence on totally bounded sets} in $X$. We call $X^\star$ the {\it dual space} for the space $X$.

If $B\subseteq X$ and $F\subseteq X^\star$ are arbitrary sets, then by $B^\circ$ and ${^\circ F}$ we denote their  (direct and reverse) polars (in $X^\star$ and in $X$):
$$
B^\circ= \{f\in X^\star: |f|_B:= \sup_{x\in B}|f(x)|\le 1\}, \qquad {^\circ F}=
\{x\in X: |x|_F:= \sup_{f\in F}|f(x)|\le 1\}
$$
Similarly, {\it annihilators} of $B$ and $F$ are the sets
$$
B^\perp= \{f\in X^\star:\quad \forall x\in B\quad f(x)=0\}, \qquad
{^\perp F}=\{x\in X: \quad \forall f\in F\quad f(x)=0\}.
$$

\blm\label{th-2.1} For each locally convex space $X$
\begin{itemize}
\item[(a)] if $B\subseteq X$ is totally bounded, then $B^\circ \subseteq X^\star$ is capacious;

\item[(b)] if $B\subseteq X$ is capacious, then $B^\circ \subseteq X^\star$ is totally bounded;

\item[(c)] if $F\subseteq X^\star$ is totally bounded, then ${^\circ F}\subseteq X$ is capacious;

\item[(d)] if $F\subseteq X^\star$ is capacious, then ${^\circ F}\subseteq X$ is totally bounded.
\end{itemize}
\elm

\btm\label{TH:X-psevdopolno=>X^star-psevdonasysheno} For an arbitrary LCS $X$
\bit{
\item[---] if $X$ is pseudocomplete, then $X^\star$ is pseudosaturated,

\item[---] if $X$ is pseudosaturated, then $X^\star$ is pseudocomplete.
}\eit
\etm

\bit{

\item For a linear continuous map $\ph:X\to Y$ of locally convex spaces the symbol $\ph^\star:Y^\star\to X^\star$ means its dual map
    \beq\label{ph^star(f)}
    \ph^\star(f)=f\circ\ph,\qquad f\in Y^\star.
    \eeq
    Obviously, it is also a linear continuous map, and we call it the {\it dual map} for the map $\ph$.

}\eit

\btm \cite[Theorems 3.2, 3.1]{Akbarov} For a morphism of locally convex spaces $\ph:X\to Y$
 \bit{
\item[--] if $X$ is pseudosaturated, and $\ph:X\to Y$ is a dense embedding, then $\ph^\star:Y^\star\to X^\star$ is an exact covering,

\item[--] if $X$ is pseudocomplete, and $\ph:X\to Y$ is an exact covering, then $\ph^\star:Y^\star\to X^\star$ is a dense embedding.
 }\eit
\etm

\btm\cite[Theorem 3.14]{Akbarov} \label{th-3.14} For any pseudocomplete locally convex spaces $X$ there exists a unique morphism of locally convex spaces
 \beq\label{(X^vartriangle)^star=(X^star)^triangledown}
 \xymatrix @C=5.0pc{  (X^\vartriangle)^\star\ar@{~>}[r]^{\alpha_X} & (X^\star)^\triangledown}
\eeq
such that the following diagram is commutative:
\beq \label{eqIII.18}
 \xymatrix    @R=2.5pc @C=1.5pc
 {
 (X^\vartriangle)^\star \ar@{~>}[rr]^{\alpha_X} & & (X^\star)^\triangledown
\\
 & X^\star \ar[ul]^{(\vartriangle_X)^\star} \ar[ur]_{\triangledown_{X^\star}} &
 }.
\eeq
This morphism $\alpha_X$ is an isomorphism of locally convex spaces, and the map $X\mapsto\alpha_X$ is an isomorphism of functors: for each morphism of locally convex spaces $\ph :X\to Y$ the following diagram is commutative:
   \beq\label{eqIII.19}
 \xymatrix    @R=2.5pc @C=2.0pc
 {
 (X^\vartriangle)^\star \ar@{~>}[rr]^{\alpha_X} & & (X^\star)^\triangledown
 \\
 (Y^\vartriangle)^\star \ar@{~>}[rr]^{\alpha_Y} \ar[u]^{(\varphi^\vartriangle)^\star}
 && (Y^\star)^\triangledown \ar[u]_{(\varphi^\star)^\triangledown}
 }
\eeq
\etm

\btm\cite[Theorem 3.15]{Akbarov} \label{th-3.15} For any pseudosaturated locally convex space $X$ there exists a unique morphism of locally convex spaces
 \beq\label{(X^triangledown)^star=(X^star)^vartriangle}
 \xymatrix @C=5.0pc
 {  (X^\triangledown)^\star\ar@{~>}[r]^{\beta_X} & (X^\star)^\vartriangle}
 \eeq
such that the following diagram is commutative:
\beq \label{eqIII.20}
 \xymatrix  %  @R=2.5pc @C=1.0pc
 {
 (X^\triangledown)^\star \ar@{~>}[rr]^{\beta_X}
 \ar[dr]_{(\triangledown_X)^\star} && (X^\star)^\vartriangle
\ar[dl]^{\vartriangle_{X^\star}}
\\
 & X^\star &
 };
\eeq
This morphism $\beta_X$ is an isomorphism of locally convex spaces, and the map $X\mapsto\beta_X$ is an isomorphism of cunctors: for any morphism of locally convex spaces $\ph :X\to Y$ the following diagram is commutative:
\beq \label{eqIII.21}
 \xymatrix  %  @R=2.5pc @C=1.0pc
 {
 (X^\triangledown)^\star \ar@{~>}[rr]^{\beta_X} & & (X^\star)^\vartriangle
\\
 (Y^\triangledown)^\star \ar@{~>}[rr]^{\beta_Y}
\ar[u]^{(\varphi^\triangledown)^\star} && (Y^\star)^\vartriangle
\ar[u]_{(\varphi^\star)^\vartriangle}
 }.
\eeq
\etm

\subsection{Stereotype spaces}

\paragraph{The map $i_X:X\to X^{\star\star}$.}
\label{SUBSEC:i_X}

The {\it second dual space} $X^{\star\star}$ of a locally convex space $X$ is the dual space to the first dual space:
$$
  X^{\star\star}=(X^\star)^\star
$$
(each star $\star$ means the topology of uniform convergence on totally bounded sets). The formula
\beq\label{DEF:i_X(x)(f)=f(x)}
     i_X(x)(f)=f(x),\qquad x\in X,
\eeq
defines a natural map $i_X:X\to X^{\star\star}$.

\bit{
 \item[$\bullet$]\label{DEF:otkr-otobr}
Let us say that a linear map of locally convex spaces $\ph:X\to Y$ is {\it open}\footnote{We use the notion of {\it open map} in the sense different from the one used in General topology \cite{Engelking}.}, if the image $\ph(U)$ of any neighborhood of zero  $U\subseteq X$ is a neighborhood of zero in the subspace $\ph(X)$ of $Y$ (with the topology inherited from $Y$):
$$
\forall
U\in \mathcal{U}(X) \quad \exists V\in \mathcal{U}(Y) \quad \ph(U)\supseteq
\ph(X)\cap V.
$$
Certainly, it is sufficient here to assume that $U$ is open and absolutely convex. By the obvious formula
\beq\label{ph(X)-cap-V=ph(ph^(-1)(V))}
\ph(X)\cap V=\ph\Big(\ph^{-1}(V)\Big),\qquad V\subseteq Y,
\eeq
(valid for any map of sets $\ph:X\to Y$ and for any subset $V\subseteq Y$), this condition can be rewritten as follows:
$$
\forall
U\in \mathcal{U}(X) \quad \exists V\in \mathcal{U}(Y) \quad \ph(U)\supseteq
\ph\Big(\ph^{-1}(V)\Big).
$$
 }\eit

\btm \cite[Theorem 2.8]{Akbarov}\label{th-2.8} For each LCS $X$ the map $i_X:X\to X^{\star\star}$ is injective, open and has dense set of values in $X^{\star\star}$.
\etm

\btm \cite[Theorem 2.12]{Akbarov} \label{th-2.12} For an arbitrary LCS $X$ the following conditions are equivalent:
\begin{itemize}
\item[(i)] the space $X$ is pseudocomplete;

\item[(ii)] the map $i_X:X\to X^{\star\star}$ is surjective (and hence, bijective).
\end{itemize}
\etm

\bcor\label{cor-2.13} If a locally convex space $X$ is pseudocomplete, then a
(continuous and bijective) map $i_X^{-1}:X^{\star\star} \to X$ is defined, and it is an exact covering.
\ecor

\btm \cite[Theorem 2.14]{Akbarov} \label{th-2.14} For an arbitrary LCS $X$ the following conditions are equivalent:
\begin{itemize}
\item[(i)] the space $X$ is pseudosaturated;

\item[(ii)] the map $i_X:X\to X^{\star\star}$ is continuous.
\end{itemize}
\etm

\bcor\label{cor-2.15} If a locally convex space $X$ is pseudosaturated, then the map $i_X:X\to X^{\star\star}$ is a dense embedding.
\ecor

\paragraph{Definition of stereotype space and examples.}

 \bit{
\item A locally convex space $X$ over $\C$ (one can consider also $\R$) is said to be {\it stereotype}, if its natural map into the second dual space $i_X:X\to X^{\star\star}$ is an isomorphism of locally convex spaces.

\item The class of all stereotype spaces is denoted by ${\tt Ste}$, it forms a category with linear continuous maps as morphisms.

}\eit
Certainly, if $X$ is a stereotype space then its dual space $X^\star$ is also stereotype.
From Theorems \ref{th-2.12} and \ref{th-2.14} we have the following criterion:

\btm\label{TH:kriterij-ster}
A locally convex space $X$ is stereotype if and only if it is pseudocomplete and pseudosaturated.
\etm

This mean in particular, that there are non-stereotype locally convex spaces (since there are non-pesudocomplet and non-pseudosaturated spaces, see Examples \ref{ex-4.9} and \ref{ex-4.10}). Nevertheless, the class of stereotype spaces ${\tt Ste}$ turns out to be amazingly wide. This is seen from the following series of examples, generalizing each other.

\bex All {\sl Banach spaces} are stereotype.
\eex

\bex All {\sl Fr\'{e}chet spaces} are stereotype.
\eex

\bex\label{ex-4.3} All {\sl quasicomplete barreled spaces} are stereotype.
\eex

As a corollary, the place of stereotype spaces among other frequently used classes of spaces can be illustrated by the following diagram:

 \beq\label{DIAGR:ster-prostranstva}
 {\small
 \begin{picture}(380,220)
\put(195,105){\oval(330,200)} \put(110,180){\text{\sf\large STEREOTYPE SPACES}} \put(195,90){\oval(270,140)} \put(110,130){\text{\sf\large quasicomplete barreled spaces}}
 \put(160,80){\oval(140,60)\put(-40,15){\text{\sf Fr\'{e}chet spaces}}}
 \put(160,70){\oval(100,25)\put(-40,-2){\text{\sf Banach spaces}}}
 \put(243,65){\oval(110,60)} \put(245,60){\text{\sf reflexive}}
 \put(249,50){\text{\sf spaces}}
\end{picture}
 }
 \eeq
This picture is supplemented by the examples of spaces, dual to the already mentioned, and having quite unwonted\footnote{Because of the non-standard notion of dual space.} properties:

\bex A locally convex space $X$ is called a {\it Smith space}\footnote{After M.F.Smith \cite{Smith}.}, if it is a complete $k$-space\footnote{\label{DEF:k-space}A topological space $X$ is called {\it $k$-space} or {\it Kelley space}, if every set $M\subseteq X$ having closed trace $M\cap K$ on each compact set $K\subseteq X$ is closed in $X$.} and has a {\it universal compact set}, i.e. a compact set $K\subset X$ that absorbs any other compact set $T\subset X$:  $T\subseteq \lambda K$ for some $\lambda\in \C$. It is known that {\it a locally convex space $X$ is a Smith space if and only if it is stereotype and its dual space $X^\star$ is a Banach space}.
\eex

\bex A locally convex space $X$ is called a {\it Brauner space}\footnote{After K.Brauner \cite{Brauner}.}, if it is a complete $k$-space\footnote{See footnote \ref{DEF:k-space}.} and has a {\it countable fundamental system of compact sets}, i.e. a sequence of compact sets $K_n\subseteq X$ such that every compact set $T\subseteq X$ is contained in some $K_n$. {\it A locally convex space $X$ is a Brauner space if and only if it is stereotype and its dual space $X^\star$ is a Fr\'{e}chet space}.
\eex

The connections between the spaces of Fr\'{e}chet, Brauner, Banach, and Smith are illustrated in the
following diagram (where the 180 degree rotation corresponds to the passage to the dual class):
$$
 \begin{picture}(400,140)
 \put(130,90){\oval(210,80)} \put(80,110){\text{\sf\Large Fr\'{e}chet spaces}}
 \put(170,70){\oval(290,40)[l]} \put(60,68){\text{\sf Banach spaces}}
 \put(260,50){\oval(210,80)} \put(220,23){\text{\sf\Large Brauner spaces}}
 \put(220,70){\oval(290,40)[r]} \put(270,68){\text{\sf Smith spaces}}
 \put(165,73){\text{\sf\footnotesize finite-dimensional}}
 \put(185,65){\text{\sf\footnotesize spaces}}
\end{picture}
$$

\paragraph{Completeness.}

\btm\cite[Theorem 4.21]{Akbarov}\label{TH:STE-polnaya-kategoriya} The category ${\tt Ste}$ is complete: each functor from a smal category into $\tt Ste$ has injective and projective limit.
\etm

In th case of direct sums and direct products these constructions coincide with those in the category ${\tt LCS}$ of locally convex spaces, while in the general case the difference is that the injective limits in ${\tt LCS}$ must be pseudocompleted, and the projective limits must be pseudosaturated:
\begin{align}
& \label{summa-v-Ste}
\overset{\tt Ste}{\bigoplus_{i\in I}} X_i=\overset{\tt LCS}{\bigoplus_{i\in I}} X_i  && \overset{\tt Ste}{\prod_{i\in I}} X_i=\overset{\tt LCS}{\prod_{i\in I}} X_i \\
& \label{rightlim-v-Ste}
\overset{\tt Ste}{\rightlim_{i\to \infty}}X_i=\l \overset{\tt LCS}{\rightlim_{i\to \infty}}X_i\r^\triangledown, &&
\overset{\tt Ste}{\leftlim_{i\to \infty}}X_i=\l \overset{\tt LCS}{\leftlim_{i\to \infty}}X_i\r^\vartriangle.
 \end{align}

\subsection{Nodal decomposition, envelopes and refinements in $\Ste$}

\paragraph{Subspaces and the envelope of a set of vectors.}

\bit{
\item[$\bullet$]
Let $Y$ be a subset in a stereotype space $X$ endowed with the structure of stereotype space in such a way that the set-theoretic enclosure  $Y\subseteq X$ becomes a morphism of stereotype spaces (i.e. a linear continuous map). Then the stereotype space $Y$ is called a {\it subspace} of the stereotype space $X$, and the set-theoretic enclosure $\sigma:Y\subseteq X$ its {\it representing monomorphism}. The record
$$
Y\subarr X
$$
or
$$
X\suparr Y
$$
will mean that $Y$ is a subspace of the stereotype space $X$. If in addition we write
$$
Y=X
$$
then this means that the stereotype spaces $Y$ and $X$ coincide not only as sets but also with their algebraic and topological structure.

}\eit

\bprop\cite[Proposition 4.65]{Akbarov-env}\label{PROP:mu-cong-podpr-v-Ste} For a morphism $\mu:Z\to X$ in the category ${\tt Ste}$ of stereotype spaces the following conditions are equivalent:
 \bit{
\item[(i)] $\mu$ is a monomorphism,

\item[(ii)] there exists a subspace $Y$ in $X$ with the representing monomorphism $\sigma:Y\subarr X$ and an isomorphism $\theta:Z\to Y$ of stereotype spaces such that the following diagram is commutative:
$$
 \xymatrix @R=1.0pc @C=2.5pc
 {
 Z\ar[rd]^{\mu}\ar@{~>}[dd]_{\theta} & \\
  & X\\
 Y\ar@{-->}[ru]_{\sigma} &
 }
$$
 }\eit
\eprop

\bcor\label{COR:Ste-lok-mala-v-Mono}
The category $\tt Ste$ is well-powered in the class $\Mono$.
\ecor

 \bit{
 \item[$\bullet$]
Suppose we have a sequence of two subspaces
$$
Z\subarr Y\subarr X,
$$
and the enclosure $Z\subarr Y$ is a bimorphism of stereotype spaces, i.e. apart from the other requirements, $Z$ is dense in $Y$ (with respect to the topology of $Y$):
$$
\overline{Z}^Y=Y.
$$
Then we will say that the subspace $Y$ is a {\it mediator} for the subspace $Z$ in the space $X$.

 \item[$\bullet$]
We call a subspace $Z$ of a stereotype space $X$ an {\it immediate subspace} in $X$, if it has no non-isomorphic mediators, i.e. for any mediator $Y$ in $X$ the corresponding enclosure $Z\subarr Y$ is an isomorphism. In this case we use the record $Z\osubarr X$:
 $$
 Z\osubarr X\qquad\Longleftrightarrow\qquad \forall Y\quad \bigg( \Big(Z\subarr Y\subarr X \quad\&\quad \overline{Z}^Y=Y\Big)\quad\Longrightarrow\quad Z=Y\bigg).
 $$

 }\eit

\brem In the category of locally convex spaces ${\tt LCS}$ the same construction gives a widely used object: immediate subspaces in a locally convex space $X$ are exactly closed subspaces in $X$ with the topology inherited from $X$. Below in Examples \ref{EX:zamk-neposr-podpr-so-strogo-mazhor-topol} and \ref{EX:nezamk-neposr-podpr} we will see that in the category ${\tt Ste}$ of stereotype spaces the situation becomes sufficiently more complicated.
\erem

\bprop\cite[Proposition 4.68]{Akbarov-env}\label{PROP:mu-cong-neposr-podpr-v-Ste} For a morphism $\mu:Z\to X$ in the category ${\tt Ste}$ the following conditions are equivalent:
 \bit{
\item[(i)] $\mu$ is an immediate monomorphism\footnote{Immediate monomorphisms were defined on page \pageref{DEF:immediate-mono}.},

\item[(ii)] there exists an immediate subspace $Y$ of $X$ with a representing monomorphism $\sigma:Y\subarr X$ and an isomorphism $\theta:Z\to Y$ such that the following diagram is commutative
\beq\label{predstavlenie-monomorfizma}
 \xymatrix @R=1.0pc @C=2.5pc
 {
 Z\ar[rd]^{\mu}\ar@{~>}[dd]_{\theta} & \\
  & X\\
 Y\ar@{-->}[ru]_{\sigma} &
 }
\eeq
 }\eit\noindent
 The subspaces $Y$ and the morphism $\theta$ here are uniquely defined by $Z$ and $\mu$.
\eprop

\medskip
\centerline{\bf Properties of immediate subspaces:\footnote{See proofs in \cite[Section 4.5]{Akbarov-env}}}

\bit{\it

\item[$1^\circ$.] If $Z\subseteq Y\osubarr X$ and $Z\subarr X$, then $Z\subarr Y$.

\item[$2^\circ$.] If $Z\subseteq Y\osubarr X$ and $Z\osubarr X$, then $Z\osubarr Y$.

\item[$3^\circ$.] If $Y\osubarr X$ and $\overline{Y}^X=Y$, then the topology of $Y$ is the pseudosaturation of the topology induced from $X$.

\item[$4^\circ$.] If $Y\osubarr X$ and $\overline{Y}^X=X$, then $Y=X$.

\item[$5^\circ$.] If $Y\osubarr X$ and $T\subseteq Y$ is a compact subset in $X$, then $T$ is compact in $Y$.

}\eit

\bex\cite[Example 4.70]{Akbarov-env}\label{EX:zamk-neposr-podpr-so-strogo-mazhor-topol} There exists a stereotype space $P$ with a closed immediate subspace $Q$, which topology is not inherited from $P$, and, moreover, some continuous functionals $g\in Q^\star$ cannot be continuously extended on $P$.
\eex

\bex\cite[Example 4.71]{Akbarov-env}\label{EX:nezamk-neposr-podpr} There exists a stereotype space $X$ with an immediate subspace $Z$, such that $Z$ is not closed as a set in $X$.
\eex

\bit{

\item[$\bullet$] The {\it envelope} of a set $M\subseteq X$ in a stereotype space $X$ is the minimal immediate subspace in $X$, that contains $M$. It is denoted by $\Env^X M$, or by $\Env M$:
    \beq\label{DEF:Env^XM}
    M\subseteq\Env^XM\osubarr X\quad \&\quad \forall Y\osubarr X\quad (M\subseteq Y\ \Longrightarrow\ \Env^XM\subseteq Y)
    \eeq

    }\eit

\medskip
\centerline{\bf Properties of $\Env^XM$:\footnote{See proofs in \cite[Section 4.5]{Akbarov-env}}}

\bit{\it

\item[$1^\circ$.] The envelope $\Env^XM$ always exists.

\item[$2^\circ$.] If $M\subseteq Y\osubarr X$, then $\Env^XM=\Env^YM\osubarr Y$.

\item[$3^\circ$.] The envelope $\Env^X M$ of each set $M\subseteq X$ is an immediate subspace in $X$, that contains $M$ as a total subset:
\beq\label{M-subseteq-Sp_infty-M^X-osubarr-X}
M\subseteq \Env^X M\osubarr X,
\qquad
\overline{\Sp M}^{\Env^X M}=\Env^X M.
\eeq

\item[$4^\circ$.] Each subspace $Y$ in a stereotype space $X$ is a subspace in its envelope $\Env^X Y$
\beq\label{Y-subarr-X=>Y-subarr-Sp_infty^X(Y)}
Y\subarr X\quad\Longrightarrow\quad Y\subarr \Env^X Y,
\eeq
and $Y$ is an immediate subspace in $X$ if and only if it coincides (with the topology) with its envelope:
\beq\label{Y-osubarr-X=>Y-=-Sp_infty^X(Y)}
Y\osubarr X\quad\Longleftrightarrow\quad Y=\Env^X Y.
\eeq

\item[$5^\circ$.] If $\ph:Y\to X$ is a morphism of stereotype spaces, that maps a set $N\subseteq Y$ into a set $M\subseteq X$,
$$
\ph(N)\subseteq M,
$$
then $\ph$ continuously maps $\Env^Y N$ into $\Env^X M$:
 $$
 \xymatrix{
 Y\ar[r]^{\ph} & X\\
 \Env^Y N\ar[u]\ar@{-->}[r] & \Env^X M \ar[u]
 }
 $$
In the special cases:
 \begin{align}
& \begin{Bmatrix}Y  \subarr  X \\ \text{\rotatebox{90}{$\subseteq$}} \quad \ \text{\rotatebox{90}{$\subseteq$}}\\
N  \subseteq  M \end{Bmatrix}\quad\Longrightarrow\quad \Env^Y N\subarr \Env^X M,
\label{N-subseteq-M,Y=X}
\\
& \begin{Bmatrix}Y  \osubarr  X \\ \text{\rotatebox{90}{$\subseteq$}} \quad \ \text{\rotatebox{90}{$\subseteq$}}\\
N  \subseteq  M \end{Bmatrix}\quad\Longrightarrow\quad \Env^Y N\osubarr \Env^X M,
\label{N-subseteq-M,Y-osubarr-X}
 \end{align}

\item[$6^\circ$.] If $\overline{\Sp M}^X=M$, then the envelope $\Env^XM$ is a pseudosaturation of the space $M$ with respect to the topology induced from $X$.

\item[$7^\circ$.] If $\overline{\Sp M}^X=X$, then $\Env^XM=X$.

}\eit

\paragraph{Quotient spaces and refinements of sets of functionals.}

 \bit{
\item[$\bullet$] Let $X$ be a stereotype space, and
 \bit{

 \item[1)] in $X$ as in a locally convex space we take a closed subspace $E$,

 \item[2)] on the quotient space $X/E$ we consider an arbitrary locally convex topology $\tau$, which is majorated by the natural quotient topology of $X/E$,

 \item[3)] in the completion $(X/E)^\blacktriangledown$ of the locally convex space $X/E$ with the topology $\tau$ we take a subspace $Y$, which contains $X/E$ and is a stereotype space with respect to the topology inherited from $(X/E)^\blacktriangledown$.
 }\eit
Then we call the stereotype space $Y$ a {\it quotient space of the stereotype space} $X$, and the composition $\upsilon=\sigma\circ\pi$ of the quotient map $\pi:X\to X/E$ and the natural enclosure $\sigma:X/E\to Y$ is called the {\it representing epimorphism} of the quotient space $Y$. The record
$$
Y\quarr X
$$
or the record
$$
X\qparr Y
$$
will mean that $Y$ is a quotient space of the stereotype space $X$.

 }\eit

\bprop\label{PROP:e-cong-neposr-faktor-pr-v-Ste} For a morphism $\e:Z\gets X$ in the category ${\tt Ste}$ the following conditions are equivalent:
 \bit{
\item[(i)] $\e$ is an epimorphism,

\item[(ii)] there is a quotient space $Y$ of $X$ with the representing epimorphism $\upsilon:Y\quarr X$, and an isomorphism $\theta:Z\gets Y$ such that the following diagram is commutative:
\beq\label{predstavlenie-epimorfizma}
 \xymatrix @R=1.0pc @C=2.5pc
 {
 Z & \\
  & X\ar[lu]_{\e} \ar@{-->}[ld]^{\upsilon} \\
 Y\ar@{~>}[uu]^{\theta} &
 }
\eeq
 }\eit
\eprop

\bcor\label{COR:Ste-lok-mala-v-Epi}
The category $\tt Ste$ is co-well-powered in the class $\Epi$.
\ecor

The formalization of the idea of quotient object we have presented here has a qualitative shortcoming in comparison with the notion of subspace which we considered above: the problem is that the relation $\quarr$ does not establish a partial order in the system $\Quot(P)$ of quotient spaces of a stereotype space $P$. By the set-theoretic reasons no one of axioms of partial order (reflexivity, antisymmetry and transitivity) holds for $\quarr$. In particular, the first two axioms do not hold since the situation when $Y\quarr X$ and at the same time $Y=X$ is impossible. To explain this, let us agree for simplicity that we do not take into account the necessity to pass to a subspace in the completion which was stated in the step 3 of our definition -- then $Y\quarr X$ (and $Y\ne\varnothing$) implies by the axiom of regularity \cite[Appendix, Axiom VII]{38} that there exists an element  $y\in Y$ such that $y\cap Y=\varnothing$. But if in addition $Y=X$, then the element $y$, being a coset of $X$, i.e. a non-empty subset in $X$, must have non-empty intersection $y\cap Y=y\cap X=y\ne\varnothing$ with $X=Y$. As to the transitivity, in the situation when $Z\quarr Y$ and $Y\quarr X$ the elements of $Z$ are non-empty sets of elements of $Y$, and each such element is a non-empty set of elements of $X$. From the point of view of set theory this is not the same as if elements of $Z$ were sets of elements of $X$, so in this situation the relation $Z\quarr X$ is also impossible. This forces us to introduce a new binary relation.
 \bit{
 \item[$\bullet$] Suppose $Y\quarr X$ and $Z\quarr X$. We will say that the quotient space $Y$ {\it subordinates} the quotient space $Z$, and we write in this situation $Z\le Y$, if there exists a morphism $\varkappa:Y\to Z$ such that the following diagram is commutative:
 \beq\label{DEF:le-dlya-faktor-prostranstv}
\xymatrix @R=1pc @C=2pc
{
Y\ar[dd]_{\varkappa} &   \\
  & X\ar[lu]_{\upsilon_Y}\ar[ld]^{\upsilon_Z} \\
Z &
}
\eeq
(here $\upsilon_Y$ and $\upsilon_Z$ are representing epimorphisms for $Y$ and $Z$). The morphism $\varkappa$, if exists, must be, first, unique, and, second, an epimorphism.
 }\eit

 \bit{
 \item[$\bullet$]
Let $Y$ and $Z$ be two quotient spaces of $X$ such that
$$
Z\le Y,
$$
and the epimorphism $\varkappa:Z\gets Y$ in diagram \eqref{DEF:le-dlya-faktor-prostranstv} is a monomorphism (and hence, a bimorphism) of stereotype spaces. Then we will say that the quotient space $Y$ is a {\it mediator} for the quotient space $Z$ of the space $X$.

 \item[$\bullet$]
We call a quotient space $Z$ of a stereotype space $X$ an {\it immediate quotient space} in $X$, if it has no non-isomorphic mediators, i.e. for any its mediator $Y$ in $X$ the corresponding epimorphism $Z\quarr Y$ is an isomorphism. We write in this case $Z\oquarr X$:
 $$
 Z\oquarr X\qquad\Longleftrightarrow\qquad \forall Y\quad \bigg( \Big(Y\quarr X \quad\&\quad Z\le Y\quad \&\quad  \varkappa\in\Mono\Big)\quad\Longrightarrow\quad Z=Y\bigg).
 $$

 }\eit

\brem In the category of locally convex spaces ${\tt LCS}$ the immediate quotient spaces of a locally convex space $X$ are exactly quotient space of  $X$ by closed subspaces with the usual quotient topologies. Like in the case of subspaces, in the category ${\tt Ste}$ of stereotype spaces the situation becomes more complicated (see below Examples \ref{EX:otkr-neposr-faktor-pr-ne-X/E} and \ref{EX:nezamk-neposr-faktor-pr}).
\erem

 The following statement is dual to Proposition  \ref{PROP:mu-cong-neposr-podpr-v-Ste}, and can be proved by the dual reasoning:

\bprop\label{PROP:mu-cong-neposr-faktor-pr-v-Ste} For a morphism $\e:Z\gets X$ in the category ${\tt Ste}$ the following conditions are equivalent:
 \bit{
\item[(i)]  $\e$ is an immediate epimorphism\footnote{Immediate epimorphisms were defined on page \pageref{DEF:immediate-epi}.},

\item[(ii)] there exists an immediate quotient space $Y$ of the stereotype space $X$ with the representing morphism $\upsilon:Y\quarr X$ and an isomorphism $\theta:Z\gets Y$ such that the following diagram is commutative:
\beq\label{predstavlenie-neposr-epimorfizma}
 \xymatrix @R=1.0pc @C=2.5pc
 {
 Z & \\
  & X\ar[lu]_{\e} \ar@{-->}[ld]^{\upsilon} \\
 Y\ar@{~>}[uu]^{\theta} &
 }
\eeq
 }\eit\noindent
The quotient space $Y$ and the morphism $\theta$ are uniquely defined by $Z$ and $\e$.
\eprop

\medskip
\centerline{\bf Properties of immediate quotient spaces:\footnote{See proofs in \cite[Section 4.6]{Akbarov-env}.}}

\bit{\it

\item[$1^\circ$.] If $Y\oquarr X$, $Z\quarr X$ and $Z\le Y$, then there exists $Z'$ such that $Z\cong Z'\quarr Y$.

\item[$2^\circ$.] If $Y\oquarr X$, $Z\oquarr X$ and $Z\le Y$, then there exists $Z'$ such that $Z\cong Z'\oquarr Y$.

\item[$3^\circ$.] If the representing morphism of an immediate quotient space $\upsilon:Y\oquarr X$ is an open map, then the topology of $Y$ is a pseudocompletion of the quotient space $X/\Ker\upsilon$ (with the usual quotient topology):
    $$
    Y\cong (X/\Ker\upsilon)^\triangledown.
    $$

\item[$4^\circ$.] If the representing morphism of an immediate quotient space $\upsilon:Y\oquarr X$ is a monomorphism, then $Y=X$.

\item[$5^\circ$.] The representing morphism of an immediate quotient space $\upsilon:Y\oquarr X$ is always relatively open, i.e. for each neoighbourhood of zero $U$ in $X$ the condition
    \bit{
    \item[(a)] each functional $f\in X^\star$, bounded on $U$,
    $$
    \sup_{x\in U}|f(x)|<\infty,
    $$
    can be extended along the map $\upsilon:X\to Y$ to some functional $g\in Y^\star$:
    $$
    f=g\circ\upsilon,
    $$
    }\eit
    implies the condition
    \bit{
    \item[(b)] there exists a neighbourhood of zero $V$ in $Y$ such that
    $$
 \upsilon(U)\supseteq V\cap\upsilon(X).
    $$
    }\eit

}\eit

The following example is dual to Example \ref{EX:zamk-neposr-podpr-so-strogo-mazhor-topol}:

\bex\label{EX:otkr-neposr-faktor-pr-ne-X/E} There exists a stereotype space $P$ with an immediate quotient space of the form $Y=(P/E)^\triangledown$, which cannot be represented in the form $Y=P/F$ for a subspace $F\subseteq P$.
\eex

Example \ref{EX:nezamk-neposr-podpr} implies

\bex\label{EX:nezamk-neposr-faktor-pr} There exists a stereotype space $P$ with an immediate quotient space $Y$, which cannot be represented in the form $Y=(P/E)^\triangledown$ for some subspace $E\subseteq P$.
\eex

\bit{

\item[$\bullet$] The {\it refinement} of a set of functionals $F\subseteq X^\star$ on a stereotype space $X$ is the minimal (in the sense of the pre-order defined in Diagram \eqref{DEF:le-dlya-faktor-prostranstv}) immediate quotient space of $X$, to which all functionals $f\in F$ can be linearly and continuously extended. It is denoted by $\Rf^X F$, or by $\Rf F$:
    \beq\label{DEF:Ref^XF}
    (\rho:\Rf^XF\oquarr X\quad\&\quad F\subseteq(\Rf^XF)^\star\circ\rho)\quad \&\quad \forall \upsilon: Y\oquarr X\quad (F\subseteq Y^\star\circ\upsilon\quad \Longrightarrow\quad \Rf^XF\le Y)
    \eeq

    }\eit

\medskip
\centerline{\bf Properties of $\Rf^XF$:\footnote{See proof in \cite[Section 4.6]{Akbarov-env}.}}

\bit{\it

\item[$1^\circ$.] The refinement $\Rf^XF$ always exist.

\item[$2^\circ$.] If $F\subseteq Y^\star$ and $\upsilon:Y\oquarr X$, then the refinement of a set of functionals $F$ on $Y$ is isomorphic to the refinement of the set of functionals $F\circ\upsilon$ on $X$:
$$
\Rf^YF\cong \Rf^X(F\circ\upsilon)\oquarr X
$$

\item[$3^\circ$.] The refinement $\Rf^XF$ of each set $F\subseteq X^\star$ is an immediate quotient space of the space  $X$, such that the dual space $(\Rf^XF)^\star$ contains $F$ as a total subset:
\beq\label{overline(Sp)^((Rf^XF)^star)F=(Rf^XF)^star}
\Rf^XF\oquarr X,
\qquad
\overline{\Sp F}^{(\Rf^X F)^\star}=(\Rf^X F)^\star.
\eeq

\item[$4^\circ$.]
Each quotient space $Y$ of a stereotype space $X$ is subordinated to the refinement $\Rf^X (Y^\star\circ\upsilon)$ of the system of functionals $Y^\star\circ\upsilon=\{g\circ\upsilon;\ g\in Y^\star\}$ on the space $X$, where  $\upsilon:Y\quarr X$ is the representing epimorphism of the space $Y$:
\beq\label{Y-quarr-X=>Y-quarr-Rf^XY^star}
\upsilon:Y\quarr X\quad\Longrightarrow\quad Y\le \Rf^X (Y^\star\circ\upsilon),
\eeq
and $Y$ is an immediate quotient space of $X$ if and only if it coincides (with the topology) with this refinement:
\beq\label{Y-oquarr-X<=>Y=Rf^XY^star}
\upsilon:Y\oquarr X\quad\Longleftrightarrow\quad Y=\Rf^X (Y^\star\circ\upsilon).
\eeq

\item[$5^\circ$.] If $\ph:Y\gets X$ is a morphism of stereotype spaces, such that the dual morphism  $\ph^\star:Y^\star\to X^\star$ maps the set of functionals $G\subseteq Y^\star$ into the set of functionals $F\subseteq X^\star$,
$$
\ph^\star(G)=G\circ\ph\subseteq F,
$$
then there exists a unique morphism $\ph'$ such that the following diagram is commutative:
 $$
 \xymatrix{
 Y\ar[d] & X\ar[l]_{\ph}\ar[d]\\
 \Rf^YG & \Rf^XF\ar@{-->}[l]_{\ph'}
 }
 $$
In the special cases:
 \begin{align}
& \begin{Bmatrix}\ph:Y  \quarr  X \\
G\circ\ph  \subseteq  F \end{Bmatrix}\quad\Longrightarrow\quad \text{$\ph'$ is an epimorphism}
\label{G-ph-subseteq-F,Y-quarr-X}
\\
& \begin{Bmatrix}\ph:Y \oquarr  X \\
G\circ\ph  \subseteq  F \end{Bmatrix}\quad\Longrightarrow\quad \text{$\ph'$ is an immediate epimorphism},
\label{G-ph-subseteq-F,Y-oquarr-X}
 \end{align}

\item[$6^\circ$.] If $F$ is closed in $X^\star$ (equivalently: $F$ coincides with its second annihilator, $F=({^\perp F})^\perp$), then the refinement $\Rf^XF$ is a pseudocompletion of the locally convex quotient space of $X$ by the annihilator ${^\perp F}=\bigcap_{f\in F}\Ker f$:
    $$
    \Rf^XF=(X/({^\perp F}))^\triangledown.
    $$

\item[$7^\circ$.] If ${^\perp F}=\bigcap_{f\in F}\Ker f=0$, then $\Rf^XF=X$.

}\eit

The following result follows immediately from the definitions of $\Env^XM$ and $\Rf^XF$.

\btm\label{TH:(Env^XM)^star=Rf^(X^star)M}
For each stereotype space $X$ we have
\beq\label{(Env^XM)^star=Rf^(X^star)M}
(\Env^XM)^\star=\Rf^{X^\star}M,\qquad M\subseteq X,
\eeq
\beq\label{(Rf^XF)^star=Env^(X^star)F}
(\Rf^XF)^\star=\Env^{X^\star}F,\qquad F\subseteq X^\star.
\eeq
\etm

\paragraph{Nodal decomposition in ${\tt Ste}$.}

\btm\cite[Theorem 4.100]{Akbarov-env} \label{TH:uzlovoe-razlozhenie-v-Ste} In the category ${\tt Ste}$ of stereotype spaces each morphism $\ph:X\to Y$ has nodal decomposition \eqref{DEF:oboznacheniya-dlya-uzlov-razlozh}. The operation $\ph\mapsto \ph^\star$ of taking the dual map establishes the following identities:
\begin{align}
&  (\im_\infty\ph)^\star=\coim_\infty \ph^\star
&&  (\coim_\infty\ph)^\star=\im_\infty \ph^\star
\label{eq4.10-infty}
\\
&  (\Im_\infty\ph)^\star=\Coim_\infty \ph^\star
&&  (\Coim_\infty\ph)^\star=\Im_\infty \ph^\star
\label{eq4.10-infty-BIG}
\end{align}
\etm

\btm\cite[Theorem 4.102]{Akbarov-env}\label{TH:Im_infty-ph=Sp_infty-ph(X)} For any morphism of stereotype spaces $\ph:X\to Y$
\bit{
\item[---] its nodal image $\Im_\infty\ph$ coincides with the envelope in $Y$ of its set of values $\ph(X)$:
    \beq\label{Im_infty-ph=Sp_infty-ph(X)}
    \Im_\infty\ph=\Env^Y\ph(X)
    \eeq

\item[---] its nodal coimage $\Coim_\infty\ph$ coincides with the refinement on $X$ of a set of functionals $\ph^\star(Y^\star)$:
    \beq\label{Coim_infty-ph=Cosp_infty-ph*(Y*)}
    \Coim_\infty\ph=\Rf^X\ph^\star(Y^\star)
    \eeq
    }\eit
\etm

\btm[characterization of strong monomorphisms]\cite[Theorem 4.104]{Akbarov-env}\label{TH:strogie-monomorfizmy-v-Ste} In the category ${\tt Ste}$ for a morphism $\mu:Z\to X$ the following conditions are equivalent:
  \bit{
 \item[$(i)$] $\mu$ is an immediate monomorphism,

 \item[$(i)'$] in diagram \eqref{predstavlenie-monomorfizma} the space $Y$ is an immediate subspace in $X$,

 \item[$(ii)$] $\mu$ is a strong monomorphism,

 \item[$(ii)'$] in diagram \eqref{predstavlenie-monomorfizma} the morphism $\sigma$ is a strong monomorphism,

 \item[$(iii)$] $\mu\cong\im_\infty\mu$,

 \item[$(iv)$] $\coim_\infty\mu$ and $\red_\infty\mu$ are isomorphisms.

 }\eit
\etm

\btm[characterization of strong epimorphisms]\cite[Theorem 4.105]{Akbarov-env}\label{TH:strogie-epimorfizmy-v-Ste} In the category ${\tt Ste}$ for a morphism $\e:Z\to X$ the following conditions are equivalent:
 \bit{

 \item[$(i)$] $\e$ is an immediate epimorphism,

 \item[$(i)'$] in diagram \eqref{predstavlenie-epimorfizma} the space $Y$ is an immediate quotient space for $X$,

 \item[$(ii)$] $\e$ is a strong epimorphism,

 \item[$(ii)'$] in diagram \eqref{predstavlenie-epimorfizma} the morphism $\pi$ is a strong epimorphism,

 \item[$(iii)$] $\e\cong\coim_\infty\e$,

 \item[$(iv)$] $\im_\infty\mu$ and $\red_\infty\mu$ are isomorphisms.
 }\eit
\etm

\paragraph{${\tt Ste}$ as a pre-abelian category and basic decomposition.}\label{SUBSEC:predabelevost}

Since any two parallel morphisms $\xymatrix{X\ar@/^1ex/[r]^{\ph}\ar@/_1ex/[r]_{\psi} & Y}$ in the category ${\tt Ste}$ of stereotype spaces can be added and subtracted one from another, it is clear that ${\tt Ste}$ is an additive category. In \cite{Akbarov} it was noticed that this category is pre-abelian:

\btm\cite[Theorem 4.17]{Akbarov}\label{th-4.17} In the category ${\tt Ste}$ of stereotype spaces for each morphism $\ph:X\to Y$ the formulas
\begin{align}\label{DEF:Ker,Coker,...-v-Ste}
& \Ker\ph=\Big(\ph^{-1}(0)\Big)^\vartriangle,
&& \Coker\ph=\Big(Y/ \overline{\ph(X)}\Big)^\triangledown,
&& \Coim\ph=\Big(X/\ph^{-1}(0)\Big)^\triangledown,
&& \Im\ph=\Big(\overline{\ph(X)}\Big)^\vartriangle
\end{align}
define respectively kernel, cokernel, coimage and image. The operation $\ph\mapsto \ph^\star$ of taking dual map establishes the following connections between these objects:
\begin{align}
&  (\ker\ph)^\star=\coker \ph^\star
&&  (\coker\ph)^\star=\ker \ph^\star
&&  (\im\ph)^\star=\coim \ph^\star
&&  (\coim\ph)^\star=\im \ph^\star
\label{eq4.10}
\\
& (\Ker \ph)^\perp{}^\vartriangle=\Im \ph^\star
&&
(\Im\ph)^\perp{}^\vartriangle=\Ker \ph^\star
&&
\Ker \ph=(\Im\ph^\star)^\perp{}^\vartriangle
&&
\Im \ph=(\Ker\ph^\star)^\perp{}^\vartriangle
\label{eq4.11}
\end{align}
\etm

The pre-abelian property of ${\tt Ste}$ implies

\btm Each morphism $\ph:X\to Y$ in ${\tt Ste}$ has basic decomposition \eqref{EX:bazis-razlozh}. The operation
$\ph\mapsto \ph^\star$ of taking dual map establishes the following identities:
\begin{align}
&  (\im\ph)^\star=\coim \ph^\star
&&  (\coim\ph)^\star=\im\ph^\star
\label{svyaz-im-i-coim}
\\
&  (\Im\ph)^\star=\Coim\ph^\star
&&  (\Coim\ph)^\star=\Im\ph^\star
\label{svyaz-Im-i-Coim}
\end{align}
\etm

Formulas \eqref{DEF:Ker,Coker,...-v-Ste} imply

\btm For any morphism of stereotype spaces $\ph:X\to Y$
 \bit{
\item[---] its kernel $\Ker\ph$ and image $\Im\ph$ are closed immediate subspaces (in $X$ and $Y$ respectively),

\item[---] its coimage $\Coim\ph$ and cokernel $\Coker\ph$ are open immediate quotient spaces (of $X$ and $Y$ respectively).
}\eit
\etm

By Theorem \ref{TH:svayz-bazisnogo-i-uzlovogo-razlozheniya}, {\it each morphism $\ph:X\to Y$ in the category $\tt Ste$ is decomposed into the diagram \eqref{svayz-bazisnogo-i-uzlovogo-razlozheniya}},
\beq\label{uzlovoe-razlozhenie-v-Ste}
\xymatrix  @R=2.5pc @C=4.0pc
{
X\ar[rrr]^{\ph}\ar[rd]^{\coim_\infty\ph}\ar[dd]_{\coim\ph} & & & Y \\
& \Coim_\infty\ph\ar[r]^{\red_\infty\ph} & \Im_\infty\ph\ar[ru]^{\im_\infty\ph}\ar@{-->}[rd]_{\tau} & \\
\Coim\ph\ar[rrr]^{\red\ph}\ar@{-->}[ru]_{\sigma} & & & \Im\ph \ar[uu]_{\im\ph}
}
\eeq
where the morphisms $\sigma$ and $\tau$ are defined uniquely (by $\ph$). At the same time, {\it in the category $\tt Ste$ the morphisms $\sigma$ and $\tau$ are not necessarily isomorphisms} \cite[Example 4.101]{Akbarov-env}.

\paragraph{Envelopes and refinements in ${\tt Ste}$.}

Since $\Ste$ has nodal decomposition, is complete, well-powered and co-well-powered, it has some envelopes and refinements.

\btm\cite[Theorem 4.106]{Akbarov-env}\label{TH:obolochki-i-nachinki-otn-klassa-morfizmov-v-Ste} In the category ${\tt Ste}$ of stereotype spaces

 \bit{
\item[(a)] each space $X$ has envelopes in the classes $\Epi$ of all epimorphisms and $\SEpi$ of all strong epimorphisms with respect to arbitrary class of morphisms $\varPhi$, among which there is at least one going from $X$; in addition,
     \bit{
\item[(i)] if $\varPhi$ differs morphisms on the outside in ${\tt Ste}$, then the envelope in $\Epi$ is also an envelope in the class $\Bim$ of all bimorphisms:
    $$
    \env_\varPhi^{\Epi} X=\env_\varPhi^{\Bim} X,
    $$

\item[(ii)] if $\varPhi$ differs morphisms on the outside and is an ideal in ${\tt Ste}$, then the envelope in $\Epi$ is also an envelope in the class $\Bim$ of all bimorphisms, and in any other class $\varOmega$ which contains $\Bim$ (for example, in the class $\Mor$ of all morphisms):
    $$
    \env_\varPhi^{\Epi} X=\env_\varPhi^{\Bim} X=\env_\varPhi^\varOmega X=\env_\varPhi X,\qquad \varOmega\supseteq\Bim.
    $$
}\eit

\item[(b)] each space $X$ has refinements in the classes $\Mono$ of all monomorphisms and $\SMono$ of all strong monomorphisms by means of arbitrary class of morphisms $\varPhi$, among which there is at least one coming to $X$; in addition,
      \bit{
\item[(i)] if $\varPhi$ differs morphisms on the inside in ${\tt Ste}$, then the refinement in $\Mono$ is also a refinement in the class $\Bim$ of all bimorphisms:
    $$
    \rf_\varPhi^{\Mono} X=\rf_\varPhi^{\Bim} X.
    $$

\item[(ii)] if $\varPhi$ differs morphisms on the inside and is a left ideal in ${\tt Ste}$, then the refinement in $\Mono$ is a refinement in the class  $\Bim$ of all bimorphisms, and of any other class $\varOmega$ which contains $\Bim$ (for example, in the class $\Mor$ of all morphisms):
    $$
    \rf_\varPhi^{\Mono} X=\rf_\varPhi^{\Bim} X=\rf_\varPhi^\varOmega X=\rf_\varPhi X,\qquad \varOmega\supseteq\Bim.
    $$

 }\eit
 }\eit
\etm

\btm\cite[Theorem 4.107]{Akbarov-env}\label{TH:Env^X(M)-obolochka-v-kateg-smysle} The envelope $\Env^X M$ of a set $M$ in a stereotype space $X$ coincides with the envelope of the space\footnote{We use here the notations of \cite[p.478]{Akbarov-stein-groups}.} $\C_M$ in the class $\Epi$ of all epimorphisms of the category ${\tt Ste}$ with respect to the morphism $\ph:\C_M\to X$, $\ph(\alpha)=\sum_{x\in M}\alpha_x\cdot x$,
$$
\Env^X M=\Env_{\ph}^{\Epi}\C_M.
$$
\etm

\btm\cite[Theorem 4.108]{Akbarov-env}\label{TH:Imp^X(F)-otpechatok-v-kateg-smysle} The refinement $\Rf^X F$ of a set $F$ of functionals on a stereotype space $X$ coinsides with the refinement of the space\footnote{The notations of \cite[p.477]{Akbarov-stein-groups} are used here.} $\C^F$ in the class $\Mono$ of all monomorphisms in the category ${\tt Ste}$ by means of the morphism $\ph:X\to\C^F$, $\ph(x)^f=f(x)$, $f\in F$
$$
\Rf^X F=\Rf_{\ph}^{\Mono}\C^F.
$$
\etm

\subsection{Space of operators and tensor products}

\paragraph{Space of operators and bilinear forms.}

\bit{

\item Let $X$ and $Y$ be stereotype spaces.
 \bit{
\item[---] For any sets $A\subseteq X$, $B\subseteq Y$ we denote by $B:A$ the system of morphisms $\ph:X\to Y$ which map $A$ into $B$:
\beq\label{DEF:B:A}
\ph\in B:A\quad\Longleftrightarrow\quad \ph :X\to Y\quad \&\quad \ph (A)\subseteq B.
\eeq
The following identities are justification of this notation:
\beq\label{eq5.4}
  (\lambda\cdot B):A=\lambda\cdot (B:A), \qquad
  B:(\lambda\cdot A)=\frac{1}{\lambda}\cdot (B:A),\qquad \lambda\in \C.
\eeq

\item[---] By $Y:X$\label{DEF:Y:X} we denote the space of linear continuous maps $\ph:X\to Y$ endowed with the topology of uniform convergence on compact sets in $X$.

\item[---] By $Y\oslash X$ we denote the pseudosaturation of the space $Y:X$,
 \beq\label{DEF:Y-oslash-X}
  Y\oslash X = (Y:X)^\vartriangle
 \eeq
The space $Y\oslash X$ is stereotype and is called the {\it inner space of operators} from $X$ into $Y$.

 }\eit
}\eit

Recall the the dual map $\ph^\star$ was defined in \eqref{ph^star(f)}.

\btm\cite[Theorem 6.2]{Akbarov} \label{th-6.2} For any two stereotype spaces $X$ and $Y$ the map $\ph\in Y\oslash X \mapsto \ph^\star\in X^\star\oslash Y^\star$ is an isomorphism of stereotype spaces
\beq\label{Y-oslash-X-cong-X^star-oslash-Y^star}
Y\oslash X\cong X^\star\oslash Y^\star
\eeq
\etm

\btm\label{TH:functorialnost-star}
The map $\ph\mapsto\ph^\star$ is a contravariant functor from ${\tt Ste}$ into ${\tt Ste}$:
\beq\label{functorialnost-star}
1_X^\star=1_{X^\star},\qquad \ph^\star\circ\psi^\star=(\psi\circ\ph)^\star.
\eeq
\etm

\begin{ex}\label{ex-6.3} If $X$ is a Smith space, and $Y$ a Banach space, then $Y\oslash X =
Y:X$ is a Banach space.
\end{ex}

\begin{ex}\label{ex-6.4}
If $X$ is a Banach space, and $Y$ a Smith space, then $Y\oslash X =Y:X$ is a Smith space.
\end{ex}

\begin{ex}\label{ex-6.5}
If $X$ is a Brauner space, and $Y$ a Fr\'{e}chet space, then $Y\oslash X= Y:X$ is a Fr\'{e}chet space.
\end{ex}

\begin{ex}\label{ex-6.6}
If $X$ is a Fr\'{e}chet space, and $Y$ a Brauner space, then $Y\oslash X= Y:X$ is a Brauner space.
\end{ex}

\bit{
\item Let $X,Y,Z$ be stereotype spaces. Then
 \bit{
\item[---]\label{DEF:bilin-otobrazhenie} we say that  a bilinear map $\beta:X\times Y\to Z$ is {\it continuous}\footnote{This type of continuity is called sometimes  $(\mathcal{K}(X),\mathcal{K}(Y))$-hypocontinuity (cf.\cite{Schaefer}), where $\mathcal{K}(X)$ and $\mathcal{K}(Y)$ are systems of compact sets in $X$ and $Y$ respectively.}, if
 \bit{
\item[(1)] for each compact set $K$ in $X$ and for each neighbourhood of zero $W$ in $Z$ there is a neighbourhood of zero $V$ in $Y$ such that
$$
\beta(K,V)\subseteq W,
$$
\item[(2)] for each compact set $L$ in $Y$ and for each neighbourhood of zero $W$ in $Z$ there is a neighbourhood of zero $U$ in $X$ such that
$$
\beta(U,L)\subseteq W,
$$
 }\eit

\item[---] we denote by $Z:(X,Y)$ the space of continuous bilinear maps $\beta:X\times Y\to Z$ endowed with the topology of uniform convergence on compact sets in $X\times Y$,

\item[---] we denote by $Z\oslash(X,Y)$ the pseudosaturation of the space $Z:(X,Y)$,
 \beq\label{DEF:Z-oslash-(X,Y)}
  Z\oslash(X,Y) = (Z:(X,Y))^\vartriangle
 \eeq
 }\eit
The space $Z\oslash(X,Y)$ is stereotype, and we call it the {\it inner space of bilinear maps} from $X\times Y$ into $Z$. Like $Z:(X,Y)$, it consists of continuous bilinear maps $\beta:X\times Y\to Z$, but the topologies of $Z:(X,Y)$ and $Z\oslash(X,Y)$ may be different.\footnote{Cf. footnote \ref{FOOTNOTE:Y-oslash-X=Y:X?}, the situation with $Z:(X,Y)$ and $Z\oslash(X,Y)$ is the same.}
 }\eit

\begin{ex}\label{ex-6.8} If $X$ and $Y$ are Smith spaces, and $Z$ a Banach space, then $Z\oslash (X,Y)=Z:(X,Y)$ is a Banach space.
\end{ex}

\begin{ex}\label{ex-6.9} If $X$ and $Y$ are Banach spaces, and $Z$ is a Smith space, then $Z\oslash (X,Y)=Z:(X,Y)$ is a Smith space.
\end{ex}

\begin{ex}\label{ex-6.10} If $X$ and $Y$ are Brauner spaces, and $Z$ a Fr\'{e}chet space, then $Z\oslash (X,Y)=Z:(X,Y)$ is a Fr\'{e}chet space.
\end{ex}

\begin{ex}\label{ex-6.11} If $X$ and $Y$ are Fr\'{e}chet spaces, and $Z$ a Brauner space, then $Z\oslash (X,Y)=Z:(X,Y)$ is a Brauner space.
\end{ex}

\btm \label{th-6.12} If $X,Y,Z$ are stereotype spaces, then the formula
\beq
\beta (x,y)= \ph(y)(x) \label{eq6.1}
\eeq
defines an isomorphism of stereotype spaces
\beq
Z\oslash (X,Y) = (Z\oslash X)\oslash Y \label{eq6.2}
\eeq
\etm
\begin{rem} In the special case when $Z=\C$ we have
\beq
\C\oslash (X,Y)= X^\star \oslash Y \label{eq6.3}
\eeq
\beq
Y\oslash X = \C\oslash (Y^\star,X) \label{eq6.4}
\eeq
\end{rem}

\btm \label{th-6.13} For all stereotype spaces $X,Y,Z$ the composition map
$$
(\beta,\alpha) \in (Z\oslash Y)\times (Y\oslash X) \mapsto  \beta \circ \alpha
\in (Z\oslash X)
$$
is a continuous bilinear form.
\etm

 \bit{
\item Let $\alpha :X\gets X'$ and $\beta :Y\to Y'$ be linear continuous maps of stereotype spaces. Denote by $\beta
\oslash \alpha$ the map
$$
(\beta \oslash \alpha): (Y\oslash X)\to (Y'\oslash X')
$$
acting by formula
 \beq\label{DEF:beta-oslash-alpha}
(\beta \oslash \alpha)(\ph)= \beta \circ \ph \circ \alpha,\qquad \ph\in Y\oslash X.
 \eeq
 }\eit

\btm \label{th-6.14} For all stereotype spaces $X,Y,X',Y'$ the bilinear map
\beq\label{eq6.5}
(\beta,\alpha)\in (Y'\oslash Y)\times (X\oslash X') \mapsto \beta \oslash \alpha
\in (Y'\oslash X')\oslash (Y\oslash X)
\eeq
is continuous.
\etm

\btm\label{TH:functorialnost-oslash}
The map $(\alpha,\beta)\mapsto\beta\odot\alpha$ is a covariant functor from ${\tt Ste}\times {\tt Ste}^{\text{op}}$ into ${\tt Ste}$:
\beq\label{beta'-oslash-alpha'-circ-beta-oslash-alpha}
1_Y\oslash 1_X=1_{Y\oslash X},\qquad
\beta'\oslash\alpha'\circ\beta\oslash\alpha=
(\beta'\circ\beta)\oslash(\alpha\circ\alpha').
\eeq
\etm

\paragraph{Tensor products.}

{\it A projective (stereotype) tensor product} $X\circledast Y$ of stereotype spaces $X$ and $Y$ is defined by the equality
\beq\label{eq7.1}
X\circledast Y=(X^\star \oslash Y)^\star
\eeq
or, equivalently, due to \eqref{eq6.3},
\beq
X\circledast Y=(\C\oslash (X,Y))^\star \label{eq7.2}
\eeq
For $x\in X$ and $y\in Y$ the elementary tensor $x\circledast y\in X\circledast Y$ is defined by the equality
\beq
(x\circledast y)(\ph)= \ph(y)(x) \label{eq7.3}
\eeq
(where $\ph \in X^\star \oslash Y$, and $x\circledast y$ is considered as the element of $(X^\star \oslash Y)^\star$), or, equivalently,
\beq
(x\circledast y)(\beta)= \beta(x,y) \label{eq7.4}
\eeq
(where $\beta \in \C\oslash(X,Y)$, and $x\circledast y$ is considered as an element of $\C\oslash (X,Y)^\star$).

\begin{prop}\label{prop-7.1} The map
$\iota :(x,y)\in X\times Y \mapsto x\circledast y \in X\circledast Y$ is a continuous bilinear form.
\end{prop}

\begin{prop}\label{prop-7.2} The algebraic tensor product $X\otimes Y$ is injectively and denseley embedded into the projective tensor product$X\circledast Y$ by the formula
$$
x\otimes y \mapsto x\circledast y
$$
\end{prop}

\btm [universality of projective tensor product]\label{th-7.3} For any stereotype spaces $X,Y,Z$ and for any continuous bilinear form $\beta :X\times Y\to Z$ there is a unique linear continuous map of stereotype spaces $\widetilde{\beta}:X\circledast Y\to Z$ such that the following diagram is commutative:
$$
\begin{diagram}
\node{X\times Y} \arrow[2]{e,t}{\iota} \arrow{se,b}{\beta}
\node[2]{X\circledast Y} \arrow{sw,b}{\widetilde\beta}
\\
\node[2]{Z}
\end{diagram},
$$
where $\iota$ is defined in Proposition \ref{prop-7.1}. Moreover, the maps $\beta \mapsto \widetilde{\beta}$ is an isomorphism of stereotype spaces
\beq
Z\oslash (X,Y)=Z\oslash (X\circledast Y) \label{eq7.7}
\eeq
\etm

For any morphisms $\alpha:X\to X'$ and $\beta:Y\to Y'$ their {\it projective stereotype tensor product} $\alpha\circledast\beta:X\circledast Y\to X'\circledast Y'$ is defined by the equality
\beq\label{alpha-circledast-beta}
\alpha\circledast\beta=(\alpha^\star \oslash \beta)^\star
\eeq
(where $\alpha^\star$ and $\alpha\oslash\beta$ are defined in \eqref{ph^star(f)} and \eqref{alpha-circledast-beta}).

\btm\label{TH:functorialnost-circledast}
The mapping $(\alpha,\beta)\mapsto\alpha\circledast\beta$ is a covariant functor from ${\tt Ste}\times {\tt Ste}$ into ${\tt Ste}$:
\beq\label{functorialnost-circledast}
1_X\circledast 1_Y=1_{X\circledast Y},\qquad
\alpha'\circledast\beta'\circ\alpha\circledast\beta=
(\alpha'\circ\alpha)\circledast(\beta'\circ\beta).
\eeq
\etm

An {\it injective (stereotype) tensor product $X\odot Y$} of stereotype spaces $X$ and $Y$ is defined by the formula
\beq
X\odot Y= Y\oslash X^\star \label{eq7.8}
\eeq
or, equivalently, due to \eqref{eq6.4}, by the formula
\beq
X\odot Y= \C \oslash (X^\star,Y^\star) \label{eq7.9}
\eeq
For $x\in X$ and $y\in Y$ the elementary operator $x\odot y\in X\odot Y$ is defined by
\beq
(x\odot y)(f)=f(x)y, \quad f\in X^\star \label{eq7.10}
\eeq
(if $x\odot y$ is considered as an element of $Y\oslash X^\star$), or by
\beq
(x\odot y)(f,g)=f(x)g(y), \quad f\in X^\star, g\in Y^\star \label{eq7.11}
\eeq
(if $x\odot y$ is considered as an element of $\C\oslash (X^\star,Y^\star)$).

\begin{prop}\label{prop-7.4} The map $\iota :(x,y)\in X\times Y \mapsto x\odot y \in X\odot Y$ is a continuous bilinear form.
\end{prop}

\begin{prop}\label{prop-7.5} The algebraic tensor product $X\otimes Y$ is injectively (but not necessarily dense) embedded into the injective tensor product $X\odot Y$ by the formula
$$
x\otimes y \mapsto x\odot y
$$
\end{prop}

\bex \label{th-7.24} If $X$ and $Y$ are Banach spaces, then $X\circledast Y$ are $X\odot Y$ Banach spaces.
\eex

\bex \label{th-7.25} If $X$ and $Y$ are Smith spaces, then $X\circledast Y$ and $X\odot Y$ are Smith spaces.
\eex

\bex \label{th-7.22} If $X$ and $Y$ are Fr\'{e}chet spaces, then $X\circledast Y$ are $X\odot Y$ Fr\'{e}chet spaces.
\eex

\bex \label{th-7.23} If $X$ and $Y$ are Brauner spaces, then $X\circledast Y$ and $X\odot Y$ are Brauner spaces.
\eex

For any morphisms $\alpha:X\to X'$ and $\beta:Y\to Y'$ their {\it injective tensor product} $\alpha\odot\beta:X\odot Y\to X'\odot Y'$ is defined by the equality
\beq\label{alpha-odot-beta}
\alpha\odot\beta=\beta \oslash \alpha^\star
\eeq
(where $\alpha^\star$ and $\alpha\oslash\beta$ are defined in \eqref{ph^star(f)} and \eqref{alpha-circledast-beta}).

\btm\label{TH:functorialnost-odot}
The mapping $(\alpha,\beta)\mapsto\alpha\odot\beta$ is a covariant functor from ${\tt Ste}\times {\tt Ste}$ into ${\tt Ste}$:
\beq\label{functorialnost-odot}
1_X\circledast 1_Y=1_{X\circledast Y},\qquad
\alpha'\odot\beta'\circ\alpha\odot\beta=
(\alpha'\circ\alpha)\odot(\beta'\circ\beta).
\eeq
\etm

\btm\label{TH:@} \cite[Theorem 7.8]{Akbarov} The identity
 \beq\label{DEF:@}
@_{X,Y}(x\circledast y)=x\odot y,\qquad x\in X,\quad y\in Y
 \eeq
defines a natural trnsformation $@_{X,Y}:X\circledast Y\to X\odot Y$ of the bifunctor $\circledast$ into the bifunctor $\odot$: for each morphisms $\alpha:X\to X'$ and $\beta:Y\to Y'$ the diagram
 \beq\label{@-preobr-functorov}
\xymatrix @R=4pc @C=4pc
{
X\circledast Y\ar[r]^{@_{X,Y}}\ar[d]_{\alpha\circledast\beta} & X\odot Y\ar[d]^{\alpha\odot\beta} \\
X'\circledast Y'\ar[r]^{@_{X',Y'}} & X'\odot Y'.
}
\eeq
is commutative.
\etm

\bit{

\item The mapping $@_{X,Y}:X\circledast Y\to X\odot Y$ defined in \eqref{DEF:@} is called the {\it Grothendieck transformation} for the spaces $X$ and $Y$.\label{DEF:preobr-Groth}

}\eit

\btm\label{TH:svoistva-tenz-proizvedenij} The tensor products $\circledast$, $\odot$, the fraction $\oslash$ and the star $\star$ are connected with each other in the category ${\tt Ste}$ throght the following asomorphisms of functors:
\bit{

\item[(i)] the connection between $\circledast$, $\odot$ and $\oslash$:
    \begin{align}
& \label{eq7.16}
(X\circledast Y)^\star  \cong   Y^\star \odot X^\star  &&
(X\odot Y)^\star
\cong   Y^\star \circledast X^\star  \\
& \label{eq7.17}
 Z\oslash (X\circledast Y)  \cong  (Z\oslash X)\oslash Y &&
(X\odot Y)\oslash Z \cong X\odot (Y\oslash Z)
 \end{align}

\item[(ii)] symmetry and monoidal property of both $\circledast$ and $\odot$:
    \begin{align}
& \label{eq7.18}
\C\circledast X \cong  X  \cong  X\circledast \C &&
\C\odot X \cong  X  \cong  X\odot \C \\
& \label{eq7.19} X\circledast Y \cong  Y\circledast X &&  X\odot Y \cong Y\odot
X \\
& \label{eq7.20} (X\circledast Y)\circledast Z  \cong  X\circledast
(Y\circledast Z) && (X\odot Y)\odot Z  \cong  X\odot (Y\odot Z)
 \end{align}

\item[(iii)] connection with the direct sums ($\bigoplus$) and direct products ($\prod$):
\begin{align}
& \label{(oplus-X_i)^star=prod-X_i^star}
\Big(\bigoplus_{i\in I} X_i\Big)^\star\cong \prod_{i\in I} X_i^\star &&
\Big(\prod_{i\in I} X_i\Big)^\star\cong \bigoplus_{i\in I} X_i^\star
\\
&
\label{perestanovochnost-oslash-s-summami}
Y\oslash\Big(\bigoplus_{i\in I} X_i\Big)\cong\prod_{i\in I}\Big(Y\oslash
X_i\Big), &&
\Big(\prod_{j\in J} Y_j\Big)\oslash X\cong\prod_{j\in J}\Big(Y_j\oslash X\Big),\\
& \label{perestanovochnost-circledast-s-summami}
\left(\bigoplus_{i\in I} X_i \right)\circledast \l\bigoplus_{j\in J}Y_j\r \cong
\bigoplus_{i\in I, j\in J} (X_i \circledast Y_j) && \left(\prod_{i\in I} X_i \right)\odot
\l\prod_{j\in J}Y_j\r \cong \prod_{i\in I, j\in J} (X_i \odot Y_j)
 \end{align}

\item[(iv)] connection with the projective ($\lim_{\infty\gets i}$) and injective ($\lim_{i\to\infty}$) limits:
\begin{align}
& \label{(lim_infty-gets-i-X_i)^star=lim_i-to-infty-X_i^star}
\Big(\lim_{i\to\infty} X_i\Big)^\star\cong
\lim_{\infty\gets i} X_i^\star &&
\Big(\lim_{\infty\gets i} X_i\Big)^\star\cong
\lim_{i\to\infty} X_i^\star
\\
&
\label{perestanovochnost-oslash-s-pred}
Y\oslash\Big(\lim_{i\to\infty} X_i\Big)\cong\lim_{\infty\gets i}\Big(Y\oslash
X_i\Big), &&
\Big(\lim_{\infty\gets j} Y_j\Big)\oslash X\cong\lim_{\infty\gets
j}\Big(Y_j\oslash X\Big),\\
& \label{perestanovochnost-circledast-i-oplus-s-predelami-2}
\Big(\lim_{i\to \infty}X_i\Big)\circledast \Big(\lim_{j\to
\infty}Y_j\Big)\cong \lim_{i,j\to\infty}\Big(X_i\circledast Y_j\Big), &&
\Big(\lim_{i\to\infty}X_i\Big)\odot \Big(\lim_{j\to\infty}Y_j\Big)\cong
\lim_{i,j\to\infty}\Big(X_i\odot Y_j\Big),
 \end{align}
($\circledast$  commutes with injective limits, and $\odot$ commutes with projective limits).
}\eit
\etm
 \bpr
1. The essential part of these identities are proved in \cite{Akbarov}. In particular, \eqref{eq7.16} is noticed in \cite[(7.16)]{Akbarov}. In the identities \eqref{eq7.17} the left one follows from \cite[Theorem 7.3]{Akbarov} and \cite[Theorem 6.12]{Akbarov}, and the right one is its corollary:
\begin{multline*}
(X\odot Y)\oslash Z=\eqref{Y-oslash-X-cong-X^star-oslash-Y^star}=Z^\star\oslash (X\odot Y)^\star=\eqref{eq7.16}=
Z^\star\oslash (Y^\star\circledast X^\star)=\\=\left(\text{left identity in \eqref{eq7.17}}\right)=
(Z^\star\oslash Y^\star)\oslash X^\star=\eqref{Y-oslash-X-cong-X^star-oslash-Y^star}=
(Y\oslash Z)\oslash X^\star=\eqref{eq7.8}=X\odot (Y\oslash Z).
\end{multline*}
The identities \eqref{eq7.18}, \eqref{eq7.19} and \eqref{eq7.20} are proved in \cite{Akbarov} as (7.18), (7.19) and (7.20). It remains to prove the identities in (iii) and in (iv). And (iii) are special cases of (iv), so it is sufficient to prove (iv).

2. In (iv) the identities \eqref{(lim_infty-gets-i-X_i)^star=lim_i-to-infty-X_i^star} are evident, since they follow immediately from the definitions of injective and projective limits, and from the autoduality of the category $\tt Ste$. In the rest the key identity is the right one in \eqref{perestanovochnost-oslash-s-pred}, and we prove it first. Let $Y_i\overset{\pi_i^j}{\gets} Y_j$ be morphisms of the projective system $\{Y_i\}$. Then the morphisms $Y_i\oslash X\overset{\pi_i^j\oslash \id_X}{\gets}Y_j\oslash X$ turn the family $Y_i\oslash X$ into a projective system. We have to prove that the space $\Big(\lim_{\infty\gets i} Y_i\Big)\oslash X$ is the projective limit
of the system $\{Y_i\oslash X; \pi_i^j\oslash \id_X\}$.

Put
$$
Y=\lim_{\infty\gets i} Y_i,
$$
and let
$$
Y_j\overset{\rho_j}{\gets} Y
$$
be natural projections of the projective limit $Y$ into the system $\{Y_i;\pi_i^j\}$.
I.e., first, $\{Y;\rho_j\}$ is a cone for $\{Y_i;\pi_i^j\}$,
 \beq\label{Y-konus-dlya-Y_i}
\begin{diagram}
\node[2]{Y}\arrow{sw,t}{\rho_j}\arrow{se,t}{\rho_k}  \\
\node{Y_j} \node[2]{Y_k}\arrow[2]{w,b}{\pi_j^k}
\end{diagram}\qquad j<k
 \eeq
and, second, if $\{Z;\sigma_j\}$ is another cone,
$$
\begin{diagram}
\node[2]{Z}\arrow{sw,t}{\sigma_j}\arrow{se,t}{\sigma_k}  \\
\node{Y_j} \node[2]{Y_k}\arrow[2]{w,b}{\pi_j^k}
\end{diagram}\qquad j<k
$$
then there exists a unique morphism $Y\overset{\tau}{\gets}Z$ such that
$$
\begin{diagram}
\node{Y}\arrow{se,b}{\rho_i}\node[2]{Z}\arrow[2]{w,t,--}{\tau}\arrow{sw,b}{\sigma_i}
\\
\node[2]{Y_i}
\end{diagram}\qquad i\in I
$$

a. Let us show that the space $Y\oslash X$ is a cone for the system
$\{Y_i\oslash X; \pi_i^j\oslash \id_X\}$ with respect to the family of morphisms
$$
Y_j\oslash X\overset{\rho_j\oslash \id_X}{\longleftarrow} Y\oslash X,
$$
i.e. for any $j<k$ the following diagram is commutative:
$$
\begin{diagram}
\node[2]{Y\oslash X}\arrow{sw,t}{\rho_j\oslash\id_X}\arrow{se,t}{\rho_k\oslash\id_X}  \\
\node{Y_j\oslash X} \node[2]{Y_k\oslash X}\arrow[2]{w,b}{\pi_j^k\oslash\id_X}
\end{diagram}
$$
Indeed, for each $\ph\in Y\oslash X$ we have
$$
(\pi_j^k\oslash\id_X)\big((\rho_k\oslash\id_X)(\ph)\big)=
(\pi_j^k\oslash\id_X)\big(\rho_k\circ\ph\big)=\pi_j^k\circ\rho_k\circ\ph=\eqref{Y-konus-dlya-Y_i}=\rho_j\circ\ph=
(\rho_j\oslash\id_X)(\ph).
$$
I.e.
$$
(\pi_j^k\oslash\id_X)\circ(\rho_k\oslash\id_X)=\rho_j\oslash\id_X.
$$

b. Let $\{Z;\sigma_j\}$ be another cone for the system
$\{Y_i\oslash X; \pi_i^j\oslash \id_X\}$:
$$
\begin{diagram}
\node[2]{Z}\arrow{sw,t}{\sigma_j}\arrow{se,t}{\sigma_k}  \\
\node{Y_j\oslash X} \node[2]{Y_k\oslash X}\arrow[2]{w,b}{\pi_j^k\oslash\id_X}
\end{diagram}\qquad j<k.
$$
For each $z\in Z$ we have:
$$
(\pi_j^k\oslash\id_X)\circ\sigma_k=\sigma_j
$$
$$
\Downarrow
$$
$$
\pi_j^k\Big(\sigma_k(z)\Big)=\Big((\pi_j^k\oslash\id_X)\circ\sigma_k\Big)(z)=\sigma_j(z)
$$
$$
\Downarrow
$$
$$
\begin{diagram}
\node[2]{X}\arrow{sw,t}{\sigma_j(z)}\arrow{se,t}{\sigma_k(z)}  \\
\node{Y_j} \node[2]{Y_k}\arrow[2]{w,b}{\pi_j^k}
\end{diagram}\qquad j<k
$$
Thus, $Y_j\overset{\sigma_j(z)}{\gets}X$ is a cone for the system $\{Y_i;\pi_i^j\}$. Hence, there is a unique morphism $Y\overset{\tau(z)}{\gets}X$ such that the following diagrams are commutative:
 \beq\label{Y<-tau(z)<-X}
\begin{diagram}
\node{Y}\arrow{se,b}{\rho_i}\node[2]{X}\arrow[2]{w,t,--}{\tau(z)}\arrow{sw,b}{\sigma_i(z)}
\\
\node[2]{Y_i}
\end{diagram}
 \eeq

Since $\tau(z)\in Y\oslash X$, we have a map
$$
Y\oslash X\overset{\tau}{\longleftarrow}Z.
$$
We have to show that this is the very same unique morphism for which the following diagram is commutative:
$$
\begin{diagram}
\node{Y\oslash X} \arrow{se,b}{\rho_i\oslash\id_X}\node[2]{Z}
\arrow[2]{w,t,--}{\tau}\arrow{sw,b}{\sigma_i}
\\
\node[2]{Y_i\oslash X}
\end{diagram}
$$

c. First, we have to verify that $\tau$ is a linear map. Take $z_1,z_2\in Z$ and $\lambda_1,\lambda_2\in\C$. Then for each $i$
 \begin{multline*}
\rho_i\circ\Big(\lambda_1\tau(z_1)+\lambda_2\tau(z_2)\Big)=
\lambda_1\rho_i\circ\tau(z_1)+\lambda_2\rho_i\circ\tau(z_2)=\eqref{Y<-tau(z)<-X}=\\=
\lambda_1\sigma_i(z_1)+\lambda_2\sigma_i(z_2)=\sigma_i\Big(\lambda_1z_1+\lambda_2
z_2\Big)
 \end{multline*}
So for each $i$ the following diagram is commutative:
$$
\dgARROWLENGTH=6em
\begin{diagram}
\node{Y}\arrow{se,b}{\rho_i}\node[2]{X}\arrow[2]{w,t,--}{\lambda_1\tau(z_1)+\lambda_2\tau(z_2)}
\arrow{sw,b}{\sigma_i(\lambda_1z_1+\lambda_2 z_2)}
\\
\node[2]{Y_i}
\end{diagram}
$$
But on the other hand, if we put $z=\lambda_1z_1+\lambda_2 z_2$ into
\eqref{Y<-tau(z)<-X}, we see that for all $i$ the following diagram is commutative:
$$
\dgARROWLENGTH=6em
\begin{diagram}
\node{Y}\arrow{se,b}{\rho_i}\node[2]{X}\arrow[2]{w,t,--}{\tau(\lambda_1z_1+\lambda_2z_2)}
\arrow{sw,b}{\sigma_i(\lambda_1z_1+\lambda_2 z_2)}
\\
\node[2]{Y_i}
\end{diagram}
$$
Recall that $Y=\lim_{\infty\gets i}Y_i$ is a projective limit, hence the morphism $Y\gets X$ for which these diagrams are commutative, is unique. We can conclude from this that
$$
\lambda_1\tau(z_1)+\lambda_2\tau(z_2)=\tau(\lambda_1z_1+\lambda_2z_2).
$$
Since this is true for any $z_1,z_2\in Z$ and $\lambda_1,\lambda_2\in\C$, the map $\tau$ is indeed linear.

d. Now let us prove that $\tau$ continuously maps $Z$ into\footnote{We use the notation $Y:X$ defined on page \pageref{DEF:Y:X}.} $Y:X$
$$
(Y:X)\overset{\tau}{\longleftarrow}Z.
$$
Let $U$ be a base neighbourhood of zero in $Y:X$, i.e.\footnote{We use the notation $B:A$ from \eqref{DEF:B:A}.}
$$
U=\Big(\bigcap_{s=1}^n \rho_{i_s}^{-1}(V_s)\Big):K,
$$
where $K$ is a compact set in $X$, and $V_1,...,V_n$ is a finite system of neighbourhoods of zeroes in the spaces $Y_{i_1},...,Y_{i_n}$ respectively. Note that for each $s$ the set
$$
W_s=\sigma_{i_s}^{-1}(V_s:K)
$$
is a neighbourhood of zero in $Z$, since the map $(Y_{i_s}:X)\overset{\sigma_{i_s}}{\longleftarrow}Z$ is continuous. Hence, the set
$$
W=\bigcap_{s=1}^n W_s=\bigcap_{s=1}^n\sigma_{i_s}^{-1}(V_s:K)
$$
is also a neighbourhood of zero in $Z$. Now for any $z\in Z$ we have:
 \begin{multline*}
\tau(z)\in U=\Big(\bigcap_{s=1}^n
\rho_{i_s}^{-1}(V_s)\Big):K\quad\Longleftrightarrow\quad \tau(z)(K)\subseteq
\bigcap_{s=1}^n \rho_{i_s}^{-1}(V_s)\quad\Longleftrightarrow\\
\Longleftrightarrow\quad \forall s\quad \tau(z)(K)\subseteq
\rho_{i_s}^{-1}(V_s)\quad \Longleftrightarrow\quad \forall s\quad
\underbrace{\rho_{i_s}\big(\tau(z)(K)\big)}_{\scriptsize\begin{matrix}\text{\rotatebox{90}{$=$}}
\\
\big(\rho_{i_s}\circ\tau(z)\big)(K)\\
\phantom{\tiny\eqref{Y<-tau(z)<-X}} \ \text{\rotatebox{90}{$=$}}\ {\tiny\eqref{Y<-tau(z)<-X}}\\
 \sigma_{i_s}(z)(K)\end{matrix}}\subseteq V_s
 \quad \Longleftrightarrow\quad \forall s\quad \sigma_{i_s}(z)(K) \subseteq V_s
 \quad\Longleftrightarrow\\
\Longleftrightarrow\quad \forall s\quad \sigma_{i_s}(z)\in V_s:K
 \quad
\Longleftrightarrow\quad \forall s\quad z\in \sigma_{i_s}^{-1}(V_s:K)=W_s
 \quad
\Longleftrightarrow\quad z\in \bigcap_{s=1}^nW_s=W
 \end{multline*}
I.e. $\tau^{-1}(U)=W$. The pre-image of the base neighbourhood of zero $U\subseteq Y:X$ is a neighbourhood of zero  $W\subseteq Z$, thus the map $(Y:X)\overset{\tau}{\longleftarrow}Z$ is continuous.

e. The map $(Y:X)\overset{\tau}{\longleftarrow}Z$ is continuous and the space $Z$ is pseudosaturated, hence the map
$$
Y\oslash X=(Y:X)^\vartriangle \overset{\tau}{\longleftarrow}Z
$$
is also continuous.

f. For each $i\in I$ and for each $z\in Z$ we have
$$
(\rho_i\oslash\id_X)(\tau(z))=\rho_i\circ\tau(z)=\eqref{Y<-tau(z)<-X}=\sigma_i(z),
$$
and this means that all the following diagram are commutative:
 $$
 \begin{diagram}
\node{Y\oslash X} \arrow{se,b}{\rho_i\oslash\id_X}
\node[2]{Z}\arrow[2]{w,t,--}{\tau}\arrow{sw,b}{\sigma_i}
\\
\node[2]{Y_i\oslash X}
 \end{diagram}
 $$

g. It remains to prove that the morphism $\tau$ with this property is unique. Suppose that $\tau'$ is another morphism with the same property, i.e. such that all the diagrams
 $$
 \begin{diagram}
\node{Y\oslash X} \arrow{se,b}{\rho_i\oslash\id_X}
\node[2]{Z}\arrow[2]{w,t,--}{\tau'}\arrow{sw,b}{\sigma_i}
\\
\node[2]{Y_i\oslash X}
 \end{diagram}
 $$
 are commutative. Then for any $z\in Z$ we have
$$
\rho_i\circ\tau'(z)=(\rho_i\oslash\id_X)(\tau'(z))=\sigma_i(z),
$$
i.e. the diagram
$$
\begin{diagram}
\node{Y}\arrow{se,b}{\rho_i}\node[2]{X}\arrow[2]{w,t,--}{\tau'(z)}\arrow{sw,b}{\sigma_i(z)}
\\
\node[2]{Y_i}
\end{diagram}
$$
is commutative. And this is true for any $i\in I$. But the morphism $Y\gets X$ with this property is unique and must coincide with the morphism $\tau(z)$ from diagrams \eqref{Y<-tau(z)<-X}. Thus,
$$
\tau'(z)=\tau(z).
$$
This is true for all $z\in Z$, therefore $\tau'=\tau$.

3. When the right identity in \eqref{perestanovochnost-oslash-s-pred} is proved, the left one becomes its corollary:
$$
Y\oslash\Big(\lim_{i\to\infty} X_i\Big)\cong\eqref{Y-oslash-X-cong-X^star-oslash-Y^star}
\cong\Big(\lim_{i\to\infty}
X_i\Big)^\star\oslash Y^\star\cong\eqref{(lim_infty-gets-i-X_i)^star=lim_i-to-infty-X_i^star}\cong\Big(\lim_{\infty\gets i}
X_i^\star\Big)\oslash Y^\star.
$$
Now it becomes possible to prove \eqref{perestanovochnost-circledast-i-oplus-s-predelami-2}:
\begin{multline*}
\Big(\lim_{i\to \infty}X_i\Big)\circledast \Big(\lim_{j\to
\infty}Y_j\Big)=\eqref{eq7.1}=\bigg( \Big(\lim_{i\to \infty}X_i\Big)^\star\oslash\Big(\lim_{j\to
\infty}Y_j\Big)\bigg)^\star=\eqref{(lim_infty-gets-i-X_i)^star=lim_i-to-infty-X_i^star}=
\bigg( \Big(\lim_{\infty\gets i}X_i^\star\Big)\oslash\Big(\lim_{j\to
\infty}Y_j\Big)\bigg)^\star
 =\\=\eqref{perestanovochnost-oslash-s-pred}=
\bigg( \lim_{\infty\gets i}\Big(X_i^\star\oslash\Big(\lim_{j\to
\infty}Y_j\Big)\Big)\bigg)^\star
=\eqref{perestanovochnost-oslash-s-pred}=
\bigg( \lim_{\infty\gets i,j}\Big(X_i^\star\oslash Y_j\Big)\bigg)^\star
=\eqref{(lim_infty-gets-i-X_i)^star=lim_i-to-infty-X_i^star}=
\lim_{i,j\to\infty}\Big(X_i^\star\oslash Y_j\Big)^\star=\\=
\eqref{eq7.1}=\lim_{i,j\to\infty}X_i\circledast Y_j
\end{multline*}
and
\begin{multline*}
\Big(\lim_{\infty\gets i}X_i\Big)\odot\Big(\lim_{\infty\gets j}Y_j\Big)=\eqref{eq7.8}=
\Big(\lim_{\infty\gets j}Y_j\Big)\oslash \Big(\lim_{\infty\gets i}X_i\Big)^\star =\eqref{(lim_infty-gets-i-X_i)^star=lim_i-to-infty-X_i^star}=
\Big(\lim_{\infty\gets j}Y_j\Big)\oslash \Big(\lim_{i\to\infty}X_i^\star\Big)=\\=
\eqref{perestanovochnost-oslash-s-pred}=
\lim_{\infty\gets j}\bigg( Y_j\oslash \Big(\lim_{i\to\infty}X_i^\star\Big)\bigg)=
\eqref{perestanovochnost-oslash-s-pred}=
\lim_{\infty\gets i,j}\Big( Y_j\oslash X_i^\star\Big)=\eqref{eq7.8}=
\lim_{\infty\gets i,j}\Big( X_i\odot Y_j\Big)
\end{multline*}
 \epr

\bcor\label{COR:Ste-monoidalnaya-kategoriya}
The category $\Ste$ is a symmetric monoidal category with respect to each of the tensor products $\circledast$ and $\odot$, and moreover, it is a closed monoidal category with respect to $\circledast$.
\ecor

The following proposition will be useful in Lemma \ref{LM:J^0_C(M)C(M)-circledast-A-cong-J^0_C(M)C(M,A)}:

\blm\label{PROP:[(X-circledast-Z)/(Y-circledast-Z)]^triangledown-cong-[(X-odot-Z)/(Y-odot-Z)]^triangledown}
Suppose $X$, $Y$ and $Z$ are stereotype spaces, and $Y$ is a complementable closed subspace in $X$. Then
\beq\label{[X/Y]^triangledown-circledast-Z}
[(X\circledast Z)/(Y\circledast Z)]^\triangledown\cong[X/Y]^\triangledown\circledast Z
\eeq
and
\beq\label{[X/Y]^triangledown-odot-Z}
[(X\odot Z)/(Y\odot Z)]^\triangledown\cong[X/Y]^\triangledown\odot Z
\eeq
In particular, if $Y$ is a closed subspace of finite codimension in $X$, then
\beq\label{[(X-circledast-Z)/(Y-circledast-Z)]^triangledown-cong-[(X-odot-Z)/(Y-odot-Z)]^triangledown}
[(X\circledast Z)/(Y\circledast Z)]^\triangledown\cong [X/Y]\circledast Z\cong [X/Y]\odot Z\cong
[(X\odot Z)/(Y\odot Z)]^\triangledown
\eeq
\elm
\bpr
Let $P$ be the complement of $Y$ in $X$:
$$
X=Y\oplus P.
$$
Then $X/Y\cong P$, and the quotient map $\ph:X\to P$ can be represented as a projection to the second component. The tensor product has the form
$$
X\circledast Z=(Y\oplus P)\circledast Z=\cite[(2.52)]{Akbarov-env}=(Y\circledast Z)\oplus (P\circledast Z),
$$
and $\ph\circledast\id_Z:X\circledast Z\to P\circledast Z$ also becomes a projection to the second component. We have
\eqref{[X/Y]^triangledown-circledast-Z}:
$$
[(X\circledast Z)/(Y\circledast Z)]^\triangledown=\Coim(\ph\circledast\id_Z)=P\circledast Z=[X/Y]^\triangledown\circledast Z.
$$
Similarly one can prove \eqref{[X/Y]^triangledown-odot-Z}. After that  \eqref{[(X-circledast-Z)/(Y-circledast-Z)]^triangledown-cong-[(X-odot-Z)/(Y-odot-Z)]^triangledown} becomes a corollary of the fact that the tensor products $\circledast$ and $\odot$ coincide when one of the multipliers is finite dimensional:
$$
P\circledast Z\cong P\odot Z,\qquad \dim P<\infty.
$$
\epr

\subsection{Involution on stereotype spaces}

\paragraph{Involution on a vector space.} We denote involution by the symbol $\bullet$ with the aim of not confusing it with the convolution $*$ in a group algebra (defined below in \eqref{DEF:svertka}).

For uniformity we sometimes denote the standard involution on $\C$ by the symbol $\bullet_{\C}$:
\beq\label{bullet_C}
\bullet_{\C}(x+yi)=\overline{x+yi}=x-yi,\qquad x,y\in\R.
\eeq

Let $X$ be a vector space over $\C$. An {\it involution} in $X$ is an arbitrary mapping $x\in X\mapsto x^\bullet\in X$ satisfying the following identities:
 \beq\label{DEF:involutsiya}
(x^\bullet)^\bullet=x,\qquad (x+y)^\bullet=x^\bullet+y^\bullet,\qquad (\lambda\cdot
x)^\bullet=\overline{\lambda}\cdot x^\bullet,\qquad x,y\in X,\quad \lambda\in\C.
 \eeq
A vector $x\in X$ is said to be {\it real}, if it coinsides with its involution:
 \beq\label{DEF:veshestv-vektor}
 x^\bullet=x.
 \eeq
Each involution $x\mapsto x^\bullet$ defines {\it real} and {\it imaginary part} of every vector $x\in X$:
 \beq\label{DEF:Re-i-Im}
\Real x=\frac{x+x^\bullet}{2},\qquad \Im x=\frac{x-x^\bullet}{2i}.
 \eeq

 \vglue10pt \centerline{\bf Properties of $\Real x$ and $\Im x$:} \vglue10pt
 \bit{\it
\item[$1^\circ.$] A vector $x\in X$ is real if and only if it coincides with its real part, or, equivalently, its imaginary part vanishes:
 \beq\label{kriterij-veshestvennosti}
x^\bullet=x\quad\Longleftrightarrow\quad \Real x=x \quad\Longleftrightarrow\quad \Im
x=0.
 \eeq

\item[$2^\circ.$] For any $x\in X$ the vectors $\Real x$ and $\Im x$ are real, and
 \beq\label{x=Re-x+i-Im-x}
x=\Real x+i\cdot\Im x.
 \eeq
Conversely, if for some real vectors $a,b\in X$ we have
 \beq\label{x=a+i-cdot-b}
x=a+i\cdot b,
 \eeq
then $a=\Real x$ and $b=\Im x$.

\item[$3^\circ.$] For any $x\in X$
 \begin{align}
& x^\bullet=\Real x-i\cdot\Im x, \label{x^*=Re-x-i-Im-x} \\
& \Real(i\cdot x)=-\Im x \label{Re-(ix)=-Im-x} \\
& \Im(i\cdot x)=\Real x. \label{Im-(ix)=Re-x}
 \end{align}

\item[$4^\circ.$] The mappings $\Real:X\to X$ and $\Im:X\to X$ are linear over $\R$, satisfy the equalities
$$
\Real^2=\Real,\quad \Real\circ\Im=\Im,\quad \Im\circ\Real=0,\quad  \Im^2=0,
$$
and have the $\R$-subspace of real vectors in $X$ as image:
 \beq\label{Real-X=Im-X}
\Real X=\Im X=\{x\in X:\ x^\bullet=x\}=\{x\in X:\ \Real x=x\}=\{x\in X:\ \Im x=x\}.
 \eeq
 }\eit

 \bit{

\item[$\bullet$] The $\R$-subspace \eqref{Real-X=Im-X} of real vectors in $X$ is called {\it real part} of $X$.

 }\eit

\bpr 1. The first property is proved in two chains of equivalences:
$$
x^\bullet=x\quad\Longleftrightarrow\quad \frac{x+x^\bullet}{2}=x
\quad\Longleftrightarrow\quad \Real x=x.
$$
$$
x^\bullet=x\quad\Longleftrightarrow\quad \frac{x-x^\bullet}{2i}=0
\quad\Longleftrightarrow\quad \Im x=0.
$$

2. The fact that  $\Real x$ and $\Im x$ are real is verified by computation, and the formula
\eqref{x=Re-x+i-Im-x} follows from definitions. Suppose $a$ and $b$ are real and
\eqref{x=a+i-cdot-b} holds. Then
$$
x^\bullet=a^\bullet+\overline{i}\cdot b^\bullet=a-i\cdot b
$$
Hence
$$
\Real x=\frac{x+x^\bullet}{2}=\frac{(a+i\cdot b)+(a-i\cdot b)^\bullet}{2}=a,\qquad \Im
x=\frac{x-x^\bullet}{2i}=\frac{(a+i\cdot b)-(a-i\cdot b)^\bullet}{2i}=b.
$$

3. The formula \eqref{x^*=Re-x-i-Im-x} follows from the fact that $\Real x$ and $\Im
x$ are real:
$$
x^\bullet=(\Real x+i\cdot \Im x)^\bullet=(\Real x)^\bullet+\overline{i}\cdot (\Im x)^\bullet=\Real x-
i\cdot \Im x,
$$
and \eqref{Re-(ix)=-Im-x} and \eqref{Im-(ix)=Re-x} from $2^\circ$:
$$
\Real(i\cdot x)=\Real(i\cdot \Real x+i^2\cdot\Im x)=-\Im x,
$$
$$
\Im(i\cdot x)=\Im(i\cdot \Real x+i^2\cdot\Im x)=\Real x.
$$

4. Let us prove the identity $\Real X=\Im X$. First,
$$
x\in\Real X\quad\Longleftrightarrow\quad \exists y\in X: \ x=\Real y
\quad\Longleftrightarrow\quad \Real x=\Real(\Real y)=\Real y=x
\quad\Longleftrightarrow\quad x^\bullet=x.
$$
And, second,
$$
x\in\Im X\quad\Longleftrightarrow\quad \exists y\in X: \ x=\Im y
\quad\Longleftrightarrow\quad \Real x=\Real(\Im y)=\Im y=x
\quad\Longleftrightarrow\quad x^\bullet=x.
$$
 \epr

\paragraph{Involution on a stereotype space.}

{\it Involution} in a stereotype space $X$ is an arbitrary involution $x\mapsto x^\bullet$ in $X$ as a vector space, which is in addition a continuous mapping.

If $X$ and $Y$ are stereotype spaces with involutions $\bullet_X$ and $\bullet_Y$, then each operator $\ph:X\to Y$ satisfying the identity
$$
\ph\circ\bullet_X=\bullet_Y\circ\ph.
$$
is said to be {\it involutive}\label{DEF:involutivnyi-operator}.

The formulas \eqref{beta-oslash-alpha} and \eqref{beta'-oslash-alpha'-circ-beta-oslash-alpha} above have sense when $\alpha$ and $\beta$ are antilinear mappings, since in this case the composition $\beta\circ\ph\circ\alpha$ is still linear over $\C$. This is used in the definition of involutions on the spaces of operators and in tensor products.

\blm For any operator $\ph:X\to Y$
\beq\label{bullet_C-oslash-(bullet_C-oslash-ph)}
\bullet_{\C}\oslash(\bullet_{\C}\oslash\ph)=i_Y\circ\ph\circ i_X^{-1}
\eeq
$$
 \xymatrix @R=2pc @C=.5pc
 {
X^{\star\star}\ar@{=}[r]&\C\oslash(\C\oslash X)\ar[rrrrr]^{\bullet_{\C}\oslash(\bullet_{\C}\oslash\ph)}& & & & & \C\oslash(\C\oslash Y)\ar@{=}[r]& Y^{\star\star}\\
 X\ar[rrrrrrr]^{\ph}\ar[u]^{i_X} & & & & & & & Y\ar[u]^{i_Y}
 }
$$
\elm
\bpr
For any $x\in X$ and $f\in Y^\star$ we have:
\begin{multline*}
(\bullet_{\C}\oslash(\bullet_{\C}\oslash\ph))(i_X(x))(f)=(\bullet_{\C}\circ i_X(x)\circ (\bullet_{\C}\oslash\ph))(f)=
(\bullet_{\C}\circ i_X(x))((\bullet_{\C}\circ f\circ \ph))=\\=
\bullet_{\C}(i_X(x)((\bullet_{\C}\circ f\circ \ph)))=
\bullet_{\C}(\bullet_{\C}(f(\ph(x))))=f(\ph(x)=i_Y(\ph(x))(f)
\end{multline*}
\epr

\btm\label{PROP:f^bullet(a)=overline(f(sigma(a)^bullet))}
Let $\bullet$ be an involution on a stereotype space $X$. Then the formula
\beq\label{DEF:inv-v-X^star}
(\bullet_{\C}\oslash\bullet)(f)=\bullet_{\C}\circ f\circ \bullet
\eeq
or, equivalently, formula
\beq\label{f^*(x)=overline(f(x^*))}
f^\bullet(x)=\overline{f(x^\bullet)},\qquad f\in X^\star.
\eeq
define an involution $\bullet_{\C}\oslash\bullet$ in the dual space $X^\star$.
The operators
 \beq\label{f-mapsto-f|_Real X}
P(f)=f|_{\Real X},\qquad Q(g)(x)=g(\Real x)+i\cdot g(\Im x),\qquad f\in\Real
(X^\star_{\C})\quad g\in(\Real X)^\star_{\R},
 \eeq
are mutually inverse isomorphisms between the real part of $X^\star$ and the dual space over the field $\R$ to the real part of $X$:
 \beq\label{Real(X^star_C)=(Real-X)^star_R}
\Real (X^\star_{\C})\cong (\Real X)^\star_{\R}
 \eeq
\etm
\bpr 1. The formula \eqref{DEF:inv-v-X^star} defines an involution on $X^\star$, since
$$
(\bullet_{\C}\oslash\bullet)\circ(\bullet_{\C}\oslash\bullet)=\eqref{beta'-oslash-alpha'-circ-beta-oslash-alpha}= (\bullet_{\C}\circ\bullet_{\C})\oslash(\bullet\circ\bullet) =1_{\C}\oslash 1_X=1_{\C\oslash X}.
$$
2. Let us verify that for any $g\in\Real (X^\star_{\C})$ the mapping $Q(g):X\to\C$ is a linear functional over $\C$. It is additive since $g$, $\Real$ and $\Im$ are additive mappings:
 \beq\label{Q(g)(x+y)=Q(g)(x)+Q(g)(y)}
Q(g)(x+y)=Q(g)(x)+Q(g)(y),\qquad x,y\in X^\star_{\C}.
 \eeq
By the same reasons it is linear over $\R$:
 \beq\label{Q(g)(lambda-cdot-x)=lambda-cdot-Q(g)(x)}
Q(g)(\lambda\cdot x)=\lambda\cdot Q(g)(x),\qquad x\in
X^\star_{\C},\quad\lambda\in\R.
 \eeq
Note also that
 \begin{multline*}
Q(g)(i\cdot x)=g(\Real (i\cdot x))+i\cdot g(\Im(i\cdot x))=
\eqref{Re-(ix)=-Im-x},\eqref{Im-(ix)=Re-x}=\\=g(-\Im x)+i\cdot g(\Real x)=
i\cdot\big(i\cdot g(\Im x)+g(\Real x)\big)=i\cdot Q(g)(x)
 \end{multline*}
Together with \eqref{Q(g)(x+y)=Q(g)(x)+Q(g)(y)} and \eqref{Q(g)(lambda-cdot-x)=lambda-cdot-Q(g)(x)} this gives linearity of $Q(g)$ over $\C$:
 \begin{multline*}
Q(g)(\lambda\cdot x)=Q(g)(\Real\lambda\cdot x+i\cdot \Im\lambda\cdot x)=
Q(g)(\Real\lambda\cdot x)+Q(g)(i\cdot \Im\lambda\cdot x)=\\=
Q(g)(\Real\lambda\cdot x)+i\cdot Q(g)(\Im\lambda\cdot x)=\Real\lambda\cdot
Q(g)(x)+i\cdot\Im\lambda\cdot Q(g)(x)=\\= (\Real\lambda+i\cdot \Im\lambda)\cdot
Q(g)(x)=\lambda\cdot Q(g)(x).
 \end{multline*}

3. Let us show that for any $g\in(\Real X)^\star_{\R}$ the mapping $Q(g):X\to\C$ is a real vector in $X^\star_{\C}$:
 \begin{multline*}
Q(g)^\bullet(x)=\overline{Q(g)(x^\bullet)}=\overline{g(\Real x^\bullet)+i\cdot g(\Im x^\bullet)}=
\overline{g(\Real x)+i\cdot g(-\Im x)}=\\=\overline{g(\Real x)-i\cdot g(\Im
x)}=g(\Real x)+i\cdot g(\Im x)=Q(g)(x).
 \end{multline*}

4. We see that $Q$ maps $(\Real X)^\star_{\R}$ into $\Real(X^\star_{\C})$.
Let us show now that this map is linear over $\R$. For any $g,h\in(\Real
X)^\star_{\R}$
$$
Q(g+h)(x)=(g+h)(\Real x)+i\cdot (g+h)(\Im x)= \big(g(\Real x)+i\cdot g(\Im
x)\big)+\big(h(\Real x)+i\cdot h(\Im x)\big)=Q(g)(x)+Q(h)(x).
$$
and for any $\lambda\in\R$
$$
Q(\lambda\cdot g)(x)=(\lambda\cdot g)(\Real x)+i\cdot (\lambda\cdot g)(\Im x)=
\lambda\cdot \big(g(\Real x)+i\cdot g(\Im x)\big)=\lambda\cdot Q(g)(x).
$$

5. Thus, $Q$ is a linear (over $\R$) operator from $(\Real X)^\star_{\R}$ into
$\Real(X^\star_{\C})$. As to $P$, it is obviously a linear operator from $\Real(X^\star_{\C})$ into $(\Real X)^\star_{\R}$. The continuity of these operators follow from the continuity of operations $\bullet$, $\Real$ and $\Im$. We now have to verify that they are mutually inverse. First, for each $f\in\Real
(X^\star_{\C})$ and $x\in X$ we have
$$
Q(P(f))(x)=P(f)(\Real x)+i\cdot P(f)(\Im x)=f(\Real x)+i\cdot f(\Im x)= f(\Real
x+i\cdot\Im x)=f(x).
$$
And, second, for each $g\in(\Real X)^\star_{\R}$ and $x\in\Real X$
$$
P(Q(g))(x)=Q(g)(x)=g(\Real x)+i\cdot g(\Im x)=g(x)+i\cdot g(0)=g(x).
$$
 \epr

\btm\label{TH:inv-v-Y-oslash-X}
Let $X$ and $Y$ be stereotype spaces with involutions $\bullet_X$ and $\bullet_Y$. Then
\bit{
\item[(i)] the formula
\beq\label{bullet_Y-oslash-bullet_X}
(\bullet_Y\oslash\bullet_X)(\ph)=\bullet_Y\circ\ph\circ\bullet_X,\qquad \ph\in Y\oslash X,
\eeq
defines an involution on the space of operators $Y\oslash X$;

\item[(ii)] the formula
\beq\label{bullet_X-odot-bullet_Y}
(\bullet_X\odot\bullet_Y)(\ph)=\bullet_X\oslash(\bullet_{\C}\oslash\bullet_Y)
\eeq
defines an involution on the injective tensor product $X\odot Y=X\oslash Y^\star=X\oslash (\C\oslash Y)$, and
\beq\label{bullet_X-odot-bullet_Y(x,y)}
(\bullet_X\odot\bullet_Y)(x\odot y)=\bullet_X(x)\odot\bullet_Y(y),\qquad x\in X,\ y\in Y
\eeq

\item[(iii)] the formula
\beq\label{bullet_X-circledast-bullet_Y}
(\bullet_X\circledast\bullet_Y)(\alpha)=\bullet_{\C}\oslash((\bullet_{\C}\oslash \bullet_X)\oslash\bullet_Y)
\eeq
defines an involution on the projective tensor product
$X\circledast Y=(X^\star\oslash Y)^\star=\C\oslash((\C\oslash X)\oslash Y)$, and
\beq\label{bullet_X-circledast-bullet_Y(x,y)}
(\bullet_X\circledast\bullet_Y)(x\circledast y)=\bullet_X(x)\circledast\bullet_Y(y),\qquad x\in X,\ y\in Y
\eeq

\item[(iv)] the involutions \eqref{bullet_X-odot-bullet_Y} and \eqref{bullet_X-circledast-bullet_Y} turn onto each other under the passing to the dual space \eqref{DEF:inv-v-X^star}:
\beq\label{bullet_X-circledast-bullet_Y->bullet_X-odot-bullet_Y}
\bullet_{\C}\oslash(\bullet_X\circledast\bullet_Y)\cong (\bullet_{\C}\oslash\bullet_X)\odot(\bullet_{\C}\oslash\bullet_Y)
\eeq
\beq\label{bullet_X-odot-bullet_Y->bullet_X-circledast-bullet_Y}
\bullet_{\C}\oslash(\bullet_X\odot\bullet_Y)\cong (\bullet_{\C}\oslash\bullet_X)\circledast (\bullet_{\C}\oslash\bullet_Y)
\eeq
}\eit
\etm
\bpr
The first proposition is proved in the following chain:
$$
(\bullet_Y\oslash\bullet_X)\circ(\bullet_Y\oslash\bullet_X)=\eqref{beta'-oslash-alpha'-circ-beta-oslash-alpha}=
(\bullet_Y\circ\bullet_Y)\oslash(\bullet_X\circ\bullet_X)=1_Y\oslash 1_X=1_{Y\oslash X}.
$$
The second one:
\begin{multline*}
(\bullet_X\odot\bullet_Y)\circ(\bullet_X\odot\bullet_Y)= \big(\bullet_X\oslash(\bullet_{\C}\oslash\bullet_Y)\big)
\circ
\big(\bullet_X\oslash(\bullet_{\C}\oslash\bullet_Y)\big)=
\eqref{beta'-oslash-alpha'-circ-beta-oslash-alpha}=\\=
(\bullet_X\circ\bullet_X)\oslash\big((\bullet_{\C}\oslash\bullet_Y)\circ
(\bullet_{\C}\oslash\bullet_Y)\big)
=1_X\oslash\big((\bullet_{\C}\circ\bullet_{\C})\oslash(\bullet_Y\circ\bullet_Y)\big)=\\
=1_X\oslash(1_{\C}\oslash 1_Y)
=1_X\oslash 1_{\C\oslash_Y}=1_{X\oslash(\C\oslash_Y)}=1_{X\odot Y}
\end{multline*}
The identity \eqref{bullet_X-odot-bullet_Y(x,y)}:
\begin{multline*}
(\bullet_X\odot\bullet_Y)(x\odot y)(f)=
\big(\bullet_X\oslash(\bullet_{\C}\oslash\bullet_Y)\big)(x\odot y)(f)=\big(\bullet_X\circ(x\odot y)\circ(\bullet_{\C}\oslash\bullet_Y)\big)(f)=\\=
\bullet_X\Big((x\odot y)\big((\bullet_{\C}\oslash\bullet_Y)(f)\big)\Big)=
\bullet_X\big((x\odot y)(\bullet_{\C}\circ f\circ\bullet_Y)\big)=
\bullet_X\big((\bullet_{\C}\circ f\circ\bullet_Y)(y)\cdot x\big)=\\=
\bullet_X\big(\overline{f(\bullet_Y(y))}\cdot x\big)=
\overline{\overline{f(\bullet_Y(y))}}\cdot\bullet_X(x)=
f(\bullet_Y(y))\cdot\bullet_X(x)=(\bullet_X(x)\odot\bullet_Y(y))(f).
\end{multline*}
The third one:
\begin{multline*}
(\bullet_X\circledast\bullet_Y)\circ(\bullet_X\circledast\bullet_Y)=\big[\bullet_{\C}\oslash((\bullet_{\C}\oslash \bullet_X)\oslash\bullet_Y)\big]
\circ
\big[\bullet_{\C}\oslash((\bullet_{\C}\oslash \bullet_X)\oslash\bullet_Y)\big]=
\eqref{beta'-oslash-alpha'-circ-beta-oslash-alpha}=\\=
(\bullet_{\C}\circ\bullet_{\C})\oslash\big[((\bullet_{\C}\oslash \bullet_X)\oslash\bullet_Y)\circ
((\bullet_{\C}\oslash \bullet_X)\oslash\bullet_Y)\big]
=1_{\C}\oslash\big[(\bullet_{\C}\oslash \bullet_X)\circ(\bullet_{\C}\oslash \bullet_X))\oslash(\bullet_Y\circ\bullet_Y)\big]=\\=
1_{\C}\oslash\big[\big((\bullet_{\C}\circ\bullet_{\C})\oslash(\bullet_X\circ
\bullet_X)\big)\oslash 1_Y\big]=
1_{\C}\oslash\big[(1_{\C}\oslash 1_X)\oslash 1_Y\big]=\\=
1_{\C}\oslash(1_{{\C}\oslash X}\oslash 1_Y)=
1_{\C}\oslash 1_{({\C}\oslash X)\oslash Y}=
1_{{\C}\oslash(({\C}\oslash X)\oslash Y)}=1_{X\circledast Y}
\end{multline*}
The identity \eqref{bullet_X-circledast-bullet_Y(x,y)}:
\begin{multline*}
(\bullet_X\circledast\bullet_Y)(x\circledast y)(\ph)=
\big[\bullet_{\C}\oslash((\bullet_{\C}\oslash \bullet_X)\oslash\bullet_Y)\big](x\circledast y)(\ph)=
\big[\bullet_{\C}\circ(x\circledast y)\circ((\bullet_{\C}\oslash \bullet_X)\oslash\bullet_Y)\big](\ph)=\\=
\overline{(x\circledast y)\Big(\big((\bullet_{\C}\oslash \bullet_X)\oslash\bullet_Y\big)(\ph)\Big)}=
\overline{(x\circledast y)\Big((\bullet_{\C}\oslash \bullet_X)\circ\ph\circ\bullet_Y\Big)}=
\overline{\Big((\bullet_{\C}\oslash \bullet_X)\circ\ph\circ\bullet_Y\Big)(y)(x)}=\\=
\overline{(\bullet_{\C}\oslash \bullet_X)(\ph(\bullet_Y(y)))(x)}=
\overline{(\bullet_{\C}\circ\ph(\bullet_Y(y))\circ\bullet_X)(x)}=
\overline{\overline{\ph(\bullet_Y(y))(\bullet_X(x))}}=\\=
\ph(\bullet_Y(y))(\bullet_X(x))=(\bullet_X(x)\circledast\bullet_Y(y))(\ph)
\end{multline*}
The identity \eqref{bullet_X-circledast-bullet_Y->bullet_X-odot-bullet_Y}:
\begin{multline*}
\bullet_{\C}\oslash(\bullet_X\circledast\bullet_Y)=
\bullet_{\C}\oslash\big[\bullet_{\C}\oslash\big((\bullet_{\C}\oslash\bullet_X)\oslash\bullet_Y\big)\big]
\cong\eqref{bullet_C-oslash-(bullet_C-oslash-ph)}\cong(\bullet_{\C}\oslash\bullet_X)\oslash\bullet_Y
\cong\eqref{bullet_C-oslash-(bullet_C-oslash-ph)}\cong\\ \cong
(\bullet_{\C}\oslash\bullet_X)\oslash \big[\bullet_{\C}\oslash(\bullet_{\C}\oslash\bullet_Y)\big]=
(\bullet_{\C}\oslash\bullet_X) \odot (\bullet_{\C}\oslash\bullet_Y).
\end{multline*}
And the identity \eqref{bullet_X-odot-bullet_Y->bullet_X-circledast-bullet_Y}:
\begin{multline*}
\bullet_{\C}\oslash(\bullet_X\odot\bullet_Y)=
\bullet_{\C}\oslash\big[\bullet_X\oslash(\bullet_{\C}\oslash\bullet_Y)\big]
\cong\eqref{bullet_C-oslash-(bullet_C-oslash-ph)}\cong
\bullet_{\C}\oslash\Big\{\big[\bullet_{\C}\oslash(\bullet_{\C}\oslash\bullet_X)\big]\oslash (\bullet_{\C}\oslash\bullet_Y)\Big\}=\\=
(\bullet_{\C}\oslash\bullet_X)\circledast(\bullet_{\C}\oslash\bullet_Y)
 \end{multline*}
\epr

\paragraph{Spin.}
Let $\bullet$ be an involution on a stereotype space $X$ over $\C$. Let us call a linear continuous mapping $\sigma:X\to X$ a {\it spin} for the involution $\bullet$, if
    \beq\label{bullet-circ-sigma}
\bullet\circ\ \sigma\circ\bullet\circ\sigma=\id_X
    \eeq
or, equivalently,
    \beq\label{sigma-circ-bullet}
\sigma\circ\bullet\circ\sigma\circ\bullet=\id_X
    \eeq

\bprop\label{PROP:podkrutka}
Let $\sigma$ be a spin for $\bullet$ in $X$. Then
 \bit{
\item[(i)] the operator $\sigma$ is invertible, and
    \beq\label{sigma^(-1)}
\sigma^{-1}=\bullet\circ\ \sigma\circ\bullet,
    \eeq
\item[(ii)] the following equalities and identities hold:
    \beq\label{bullet-circ-sigma=sigma^(-1)-circ-bullet}
\bullet\circ\ \sigma=\sigma^{-1}\circ\bullet,\qquad \sigma\circ\bullet=\bullet\circ\sigma^{-1},
    \eeq
    \beq\label{bullet-circ-sigma^n}
\bullet\circ\ \sigma^n\circ\bullet\circ\sigma^n=\id_X=\sigma^n\circ\bullet\circ\sigma^n\circ\bullet,
    \eeq
    \beq\label{sigma^k-circ-bullet-circ-sigma^(k+l)-circ-bullet-circ-sigma^l}
\sigma^k\circ\bullet\circ\ \sigma^{k+l}\circ\bullet\circ\sigma^l=\id_X,
    \eeq

\item[(iii)] for any $k\in\Z$ the operator $\sigma^k$ is also a spin for $\bullet$,

\item[(iv)]  $\bullet\circ\sigma$ is an involution on $X$, and $\sigma$ is a spin for it,

\item[(v)]  $1_{\C}\oslash\sigma$ is a spin for the involution $\bullet_{\C}\oslash\bullet$ in $X^\star$.
}\eit
\eprop

\bex
In the field of complex numbers $\C$ the spins for the standard involution $x+iy\mapsto x-iy$ are exactly multiplications by a number $\theta\in\C$ which module is equal to 1:
$$
\sigma(\lambda)=\theta\cdot\lambda,\qquad |\theta|=1.
$$
Indeed, since a spin is a linear operator over $\C$, it must be a multiplication by some number $\theta$. And the condition \eqref{bullet-circ-sigma} means that $|\theta|=1$.
\eex

In the theory of Hopf algebras the passage to the dual involution is usually complemented by a spin -- in this case  the dual involution becomes an involution of Hopf algebras.

\btm\label{PROP:podkrutka-involutsii} Let $\sigma$ be a spin for an involution $\bullet$ in a stereotype space $X$. Then the formula
\beq\label{DEF:inv-v-X^star-sigma}
\big(\bullet_{\C}\oslash(\bullet\circ\sigma)\big)(f)=\bullet_{\C}\circ f\circ \bullet\circ\sigma,
\qquad f\in X^\star
\eeq
or, equivalently, the formula
    \beq\label{f^bullet(a)=overline(f(sigma(a)^bullet))}
     f^{\bullet}(x)=\overline{f(\sigma(x)^\bullet)},\qquad x\in X,\ f\in X^\star
    \eeq
defines an involution $\bullet_{\C}\oslash(\bullet\circ\sigma)$ in the dual space $X^\star$, and the dual mapping  $\sigma^\star$ is a spin for this involution:
    \beq\label{f^bullet(a)=overline(f(sigma(a)^bullet))-1}
\sigma^\star\circ\big(\bullet_{\C}\oslash(\bullet\circ\sigma)\big) \circ\sigma^\star\circ\big(\bullet_{\C}\oslash(\bullet\circ\sigma)\big)=\id_{X^\star}
\eeq
\etm

\btm
Let $X$ and $Y$ be stereotype spaces with involutions $\bullet_X$ and $\bullet_Y$ and spins $\sigma_X$ and $\sigma_Y$ respectively. Then
\bit{
\item[(i)] the operator $\sigma_Y\oslash \sigma_X$ is a spin for the involution $\bullet_Y\oslash\bullet_X$ in the space of operators $Y\oslash X$,

\item[(ii)] the operator $\sigma_X\odot\sigma_Y=\sigma_X\oslash\sigma_Y^\star$ is a spin for the involution $\bullet_X\odot\bullet_Y=\bullet_X\oslash(\bullet_{\C}\oslash\bullet_Y)$, and this involution is not changed if  the passage to the fractions $\oslash$ is complemented by the spin:
\beq\label{bullet_X-odot-bullet_Y-3}
\bullet_X\odot\bullet_Y=\bullet_X\oslash[\bullet_{\C}\oslash\bullet_Y]=
\bullet_X\oslash\big\{[\bullet_{\C}\oslash(\bullet_Y\circ\sigma_Y)]\circ\sigma_Y^\star\big\}
\eeq

\item[(iii)] the operator $\sigma_X\circledast\sigma_Y=(\sigma_X^\star\oslash \sigma_Y)^\star$ is a spin for the involution $\bullet_X\circledast\bullet_Y=\bullet_{\C}\oslash\big\{[\bullet_{\C}\oslash \bullet_X]\oslash\bullet_Y\big\}$, and this involution is not changed if the passage to the fractions $\oslash$ is complemented by the spin:
\beq\label{bullet_X-circledast-bullet_Y-3}
\bullet_X\circledast\bullet_Y=\bullet_{\C}\oslash\big\{[\bullet_{\C}\oslash \bullet_X]\oslash\bullet_Y\big\}=
\bullet_{\C}\oslash\Big\{\big[(\bullet_{\C}\oslash(\bullet_X\circ\sigma_X))\oslash (\bullet_Y\circ\sigma_X)\big]\circ(\sigma_X^\star\oslash\sigma_Y)\Big\}
\eeq

\item[(iv)] the equalities \eqref{bullet_X-circledast-bullet_Y->bullet_X-odot-bullet_Y} and \eqref{bullet_X-odot-bullet_Y->bullet_X-circledast-bullet_Y} hold if the passage to the fractions is complemented by the spin:
\beq\label{bullet_X-circledast-bullet_Y->bullet_X-odot-bullet_Y-sigma}
\bullet_{\C}\oslash\big\{(\bullet_X\circledast\bullet_Y)\circ(\sigma_X\circledast\sigma_Y)\big\}\cong \big[\bullet_{\C}\oslash(\bullet_X\circ\sigma_X)\big]\odot\big[\bullet_{\C}\oslash(\bullet_Y\circ\sigma_Y)\big]
\eeq
\beq\label{bullet_X-odot-bullet_Y->bullet_X-circledast-bullet_Y-sigma}
\bullet_{\C}\oslash\big\{(\bullet_X\odot\bullet_Y)\circ(\sigma_X\odot\sigma_Y)\big\}\cong \big[\bullet_{\C}\oslash(\bullet_X\circ\sigma_X)\big]\circledast \big[\bullet_{\C}\oslash(\bullet_Y\circ\sigma_Y)\big]
\eeq

}\eit
\etm
\bpr
The proposition (i) is proved in the following chain:
$$
(\bullet_Y\oslash\bullet_X)\circ(\sigma_Y\oslash\sigma_X)\circ(\bullet_Y\oslash\bullet_X)\circ(\sigma_Y\oslash\sigma_X)=
(\bullet_Y\circ\sigma_Y\circ\bullet_Y\circ\sigma_Y)\oslash(\bullet_X\circ\sigma_X\circ\bullet_X\circ\sigma_X)=1_Y\oslash 1_X=1_{Y\oslash X}.
$$
Then (i) implies that $\sigma_X\odot\sigma_Y=\sigma_X\oslash\sigma_Y^\star$ is a spin for $\bullet_X\odot\bullet_Y=\bullet_X\oslash(\bullet_{\C}\oslash\bullet_Y)$, and \eqref{bullet_X-odot-bullet_Y-3} is proved by the chain
\begin{multline*}
\bullet_X\oslash\big[(\bullet_{\C}\oslash(\bullet_Y\circ\sigma_Y))\circ\sigma_Y^\star\big]=
\bullet_X\oslash\big[(\bullet_{\C}\oslash(\bullet_Y\circ\sigma_Y))\circ(1_{\C}\oslash\sigma_Y)\big]= \eqref{beta'-oslash-alpha'-circ-beta-oslash-alpha}=\\=
\bullet_X\oslash\big[(\bullet_{\C}\circ 1_{\C})\oslash(\sigma_Y\circ\bullet_Y\circ\sigma_Y))\big]=\eqref{bullet-circ-sigma}=
\bullet_X\oslash(\bullet_{\C}\circ \bullet_Y)
\end{multline*}
Similarly, (ii) implies that $\sigma_X\circledast\sigma_Y=(\sigma_X^\star\oslash \sigma_Y)^\star$ is a spin for $\bullet_X\circledast\bullet_Y=\bullet_{\C}\oslash\big\{[\bullet_{\C}\oslash \bullet_X]\oslash\bullet_Y\big\}$, and \eqref{bullet_X-circledast-bullet_Y-3} is proved as follows:
\begin{multline*}
\bullet_{\C}\oslash\Big\{\big[(\bullet_{\C}\oslash(\bullet_X\circ\sigma_X))\oslash (\bullet_Y\circ\sigma_Y)\big]\circ(\sigma_X^\star\oslash\sigma_Y)\Big\}=
\bullet_{\C}\oslash\Big\{\big[(\bullet_{\C}\oslash(\bullet_X\circ\sigma_X))\oslash (\bullet_Y\circ\sigma_Y)\big]\circ((1_{\C}\oslash\sigma_X)\oslash\sigma_Y)\Big\}=\\=
\eqref{beta'-oslash-alpha'-circ-beta-oslash-alpha}=
\bullet_{\C}\oslash\Big\{\big((\bullet_{\C}\oslash(\bullet_X\circ\sigma_X))\circ(1_{\C}\oslash\sigma_X)\big)\oslash (\sigma_Y\circ\bullet_Y\circ\sigma_Y)\Big\}=\eqref{beta'-oslash-alpha'-circ-beta-oslash-alpha},\eqref{bullet-circ-sigma} =\\=
\bullet_{\C}\oslash\Big\{\big((\bullet_{\C}\circ 1_{\C})\oslash(\sigma_X\circ\bullet_X\circ\sigma_X)\big)\oslash \bullet_Y\Big\}=
\bullet_{\C}\oslash\big\{(\bullet_{\C}\oslash\bullet_X)\oslash \bullet_Y\big\}
\end{multline*}
Then \eqref{bullet_X-circledast-bullet_Y->bullet_X-odot-bullet_Y-sigma}:
\begin{multline*}
\bullet_{\C}\oslash\big\{(\bullet_X\circledast\bullet_Y)\circ(\sigma_X\circledast\sigma_Y)\big\}=
\bullet_{\C}\oslash\Big\{\big[\bullet_{\C}\oslash\big((\bullet_{\C}\oslash\bullet_X)\oslash\bullet_Y\big)\big] \circ\big[\id_{\C}\oslash\big((\id_{\C}\oslash\sigma_X)\oslash\sigma_Y\big)\big]\Big\}
=\eqref{beta'-oslash-alpha'-circ-beta-oslash-alpha}=\\=
\bullet_{\C}\oslash\Big\{(\bullet_{\C}\circ\id_{\C})\oslash\big[\big((\id_{\C}\oslash\sigma_X)\oslash\sigma_Y\big)
\circ \big((\bullet_{\C}\oslash\bullet_X)\oslash\bullet_Y\big)\big]\Big\}=\eqref{beta'-oslash-alpha'-circ-beta-oslash-alpha}=\\=
\bullet_{\C}\oslash\Big\{\bullet_{\C}\oslash\big[\big((\id_{\C}\oslash\sigma_X)\circ(\bullet_{\C}\oslash\bullet_X)\big) \oslash(\bullet_Y\circ\sigma_Y)\big]\Big\}=\eqref{beta'-oslash-alpha'-circ-beta-oslash-alpha}=\\=
\bullet_{\C}\oslash\Big\{\bullet_{\C}\oslash\big[\big((\id_{\C}\circ\bullet_{\C})\oslash(\bullet_X\circ\sigma_X)\big) \oslash(\bullet_Y\circ\sigma_Y)\big]\Big\}=
\bullet_{\C}\oslash\Big\{\bullet_{\C}\oslash\big[\big(\bullet_{\C}\oslash(\bullet_X\circ\sigma_X)\big) \oslash(\bullet_Y\circ\sigma_Y)\big]\Big\}\cong \\ \cong\eqref{bullet_C-oslash-(bullet_C-oslash-ph)}\cong
\big[\bullet_{\C}\oslash(\bullet_X\circ\sigma_X)\big]\oslash(\bullet_Y\circ\sigma_Y) \cong\eqref{bullet_C-oslash-(bullet_C-oslash-ph)}\cong
\big[\bullet_{\C}\oslash(\bullet_X\circ\sigma_X)\big]\oslash \big[\bullet_{\C}\oslash(\bullet_{\C}\oslash(\bullet_Y\circ\sigma_Y))\big]=\\=
\big[\bullet_{\C}\oslash(\bullet_X\circ\sigma_X)\big] \odot \big[\bullet_{\C}\oslash(\bullet_Y\circ\sigma_Y)\big]
\end{multline*}
And \eqref{bullet_X-odot-bullet_Y->bullet_X-circledast-bullet_Y-sigma}:
\begin{multline*}
\bullet_{\C}\oslash\big\{(\bullet_X\odot\bullet_Y)\circ(\sigma_X\odot\sigma_Y)\big\}=
\bullet_{\C}\oslash\Big\{\big[\bullet_X\oslash(\bullet_{\C}\oslash\bullet_Y)\big]\circ \big[\sigma_X\oslash (\id_{\C}\oslash\sigma_Y)\big]\Big\}=\eqref{beta'-oslash-alpha'-circ-beta-oslash-alpha}=\\=
\bullet_{\C}\oslash\Big\{(\bullet_X\circ\sigma_X)\oslash\big[(\id_{\C}\oslash\sigma_Y)\circ(\bullet_{\C}\oslash\bullet_Y)\big] \Big\}=\eqref{beta'-oslash-alpha'-circ-beta-oslash-alpha}=
\bullet_{\C}\oslash\Big\{(\bullet_X\circ\sigma_X)\oslash \big[(\id_{\C}\circ\bullet_{\C})\oslash(\bullet_Y\circ\sigma_Y)\big] \Big\}=\\=
\bullet_{\C}\oslash\Big\{(\bullet_X\circ\sigma_X)\oslash \big[\bullet_{\C}\oslash(\bullet_Y\circ\sigma_Y)\big] \Big\}\cong\eqref{bullet_C-oslash-(bullet_C-oslash-ph)}\cong
\bullet_{\C}\oslash\Big\{\Big[\bullet_{\C}\oslash\big[\bullet_{\C}\oslash(\bullet_X\circ\sigma_X)\big]\Big]\oslash \big[\bullet_{\C}\oslash(\bullet_Y\circ\sigma_Y)\big] \Big\}=\\=
\big[\bullet_{\C}\oslash(\bullet_X\circ\sigma_X)\big]\circledast \big[\bullet_{\C}\oslash(\bullet_Y\circ\sigma_Y)\big] \end{multline*}
\epr

\section{Stereotype algebras}

\subsection{Stereotype algebras and stereotype modules}\label{SUBSEC:Ste^circledast}

Recall \cite{MacLane}, that {\it algebra}, or {\it monoid} in a symmetric mnoidal category $({\tt K},\otimes)$ is a triple $(A,\mu,\iota)$, where $A$ is an object in ${\tt K}$, and $\mu:A\otimes A\to A$ (multiplication) and $\iota:I\to A$ (unit) are morphisms such that the following diagrams (axioms od associalitive and of unit) are commutative
  \beq\label{monoid}
\begin{diagram}
\node{(A\otimes A)\otimes A} \arrow{s,l}{\mu\otimes 1}
\arrow[2]{e,l}{a_{A,A,A}} \node[2]{A\otimes (A\otimes A)}
\arrow{s,r}{1\otimes\mu}
\\
\node{A\otimes A} \arrow{e,t}{\mu} \node{A} \node{A\otimes A} \arrow{w,t}{\mu}
\end{diagram}
\qquad
\begin{diagram}
\node{I\otimes A} \arrow{se,b}{\iota\otimes 1} \arrow{e,l}{l_A} \node{A}
 \node{A\otimes I} \arrow{sw,b}{1\otimes\iota}\arrow{w,l}{r_A}
\\
\node[2]{A\otimes A} \arrow{n,r}{\mu}
\end{diagram},
 \eeq
where $I$ is the unit object in ${\tt K}$, $a$ the associativity isomorphism, and $l$ and $r$ are isomorphisms of multiplication by the unit. Since by Theorem \ref{TH:svoistva-tenz-proizvedenij}(ii), $\Ste$ is a symmetric monoidal category with respect to both stereotype tensor products, $\circledast$ and $\odot$, there are two natural ways to define stereotype algebra.

 \bit{

\item A {\it projective stereotype algebra}\label{DEF:inj-algebra-circledast}, or just a {\it stereotype algebra} is a monoid in a symmetric monoidal category $(\Ste,\circledast)$, i.e. a triple $(A,\mu,\iota)$, where $A$ is a stereotype space, and $\mu:A\circledast A\to A$ (multip0lication) and $\iota:\C\to A$ (unit) are morphisms satisfying \eqref{monoid} (where $\otimes$ and $I$ are replaced by $\circledast$ and $\C$).

\item An {\it injective stereotype algebra}\label{DEF:proj-algebra-odot} is a monoid in the symmetric monoidal category $(\Ste,\odot)$, i.e. a triple $(A,\mu,\iota)$, where $A$ is a stereotype space, and $\mu:A\odot A\to A$ (multiplication) and $\iota:\C\to A$ (unit) are morphisms satisfying \eqref{monoid} (where $\otimes$ and $I$ are replaced by $\odot$ and $\C$).

 }\eit

The projective algebras form a very wide class, in contrast to injective algebras, that is why the projective algebras deserve primary interest.

\paragraph{Projective stereotype algebras.}
This class of algebras can be defined in more customary terms, without categorical constructions.

 \bit{

\item A stereotype space $A$ over $\C$ is called a {\it stereotype algebra}\label{DEF:ster-alg}, if $A$ is endowed with a structure of unital associative algebra over $\C$, and the multiplication is a continuous bilinear form in the  sense of the definition on p.\pageref{DEF:bilin-otobrazhenie}: for any compact set $K$ in $A$ and for any neighborhood of zero $U$ in $A$ there exists a neighborhood of zero $V$ in $A$ such that
$$
K\cdot V\subseteq U\quad \& \quad V\cdot K\subseteq U.
$$

 }\eit

Certainly, each stereotype algebra $A$ is a topological algebra with the separately continuous multiplication (but not vise versa). The class of all stereotype algebras is denoted by $\SteAlg$. It forms a category with unital multiplicative linear continuous maps $\ph :A\to B$ as morphisms.

Below in the text we give various examples of stereotype algebras. The following three are the first general facts.

\begin{ex}\label{ex-10.1} {\sl The Fr\'{e}chet algebras and the Banach algebras.} For a Fr\'{e}chet space $A$ the property of being stereotype algebra is equivalent to the joint continuity of the multilication operation. As a corollary, {\it each inutal  Fr\'{e}chet algebra is a stereotype algebra}. In particular, each Banach algebra is stereotype.
\end{ex}

\begin{ex}\label{ex-10.2} {\sl The algebra of operators ${\mathcal L}(X)$.}
Theorem \ref{th-6.13} implies that for each stereotype space $X$ the corresponding space ${\mathcal L}(X)=X\oslash X$ of linear continuous operators $\ph :X\to X$ is a stereotype algebra with respect to the composition $\ph \circ \psi$.
\end{ex}

\bex\label{EX:siln-loc-vyp-topol}
Each unital associative algebra $A$ over $\C$ turns into a stereotype algebra, being endowed with the strongest locally convex topology.
\eex
\bpr
One can define this topology as the locally convex topology on $A$, where each seminorm $p:A\to\R_+$ is continuous. With this topology $A$ becomes isomorphic to a locally convex direct sum of some cardinal number of copies of the space $\C$, $A=\C_{\frak m}$, that is why $A$ will be a stereotype space. We have to verify only that the multiplication will be continuous in the described above sense. Suppose $K\subset A$ is compact with respect to this topology (this means automatically that $K$ is finite dimensional), and let $x_i$ be a net tending to zero in $A$ in this topology: $x_i\to 0$. This means that for any seminorm $p:A\to\R_+$ we have $p(x_i)\to 0$. We have to verify that
$$
\sup_{a\in K}p(x_i\cdot a)\to 0,\qquad \sup_{a\in K}p(a\cdot x_i)\to 0.
$$
The set $K$ is finite-dimensional, i.e. it is contained in a convex hull of a finite set $a_1,...,a_n\in A$. Each vector  $a\in K$ can be represented as a finite combination of $a_i$, hence
$$
 p(x\cdot a)=p\left(x\cdot\sum_{i=1,...,n}\lambda_i\cdot a_i\right)\le\sum_{i=1,...,n}\lambda_i\cdot p\left(x\cdot a_i\right)=\underbrace{\left(\sum_{i=1,...,n}\lambda_i\right)}_{\scriptsize\begin{matrix}\|\\ 1\end{matrix}}\cdot\max_{i=1,...,n}p(x\cdot a_i)=
\max_{i=1,...,n}p(x\cdot a_i).
$$
We can deduce that for any $x\in A$
$$
q(x)=\sup_{a\in K}p(x\cdot a)\le \max_{i=1,...,n}p(x\cdot a_i)<\infty,
$$
therefore $q$ is a semonorm on $A$. As a corollary, it tends to zero on the net $x_i\to 0$:
$$
q(x_i)=\sup_{a\in K}p(x_i\cdot a)\to 0.
$$
Similarly we have $\sup_{a\in K}p(a\cdot x_i)\to 0$.
\epr

Of course, the category $\SteAlg$ of stereotype algebras is quite unsimilar to the category $\Ste$ of stereotype spaces, for instance, $\SteAlg$ is not additive. But this is a category with nodal decomposition:

\btm\cite[Theorem 5.30]{Akbarov-env} \label{TH:uzlovoe-razlozhenie-v-Ste^circledast} In the category $\SteAlg$ of stereotype algebras each morphism $\ph:A\to B$ has nodal decomposition \eqref{DEF:oboznacheniya-dlya-uzlov-razlozh}.
\etm

\paragraph{Tensor product of stereotype algebras.}

\btm\cite[Theorem 10.15]{Akbarov}
The projective stereotype tensor product $A\circledast B$ of (projective) stereotype algebras $A$ and $B$ has a unique structure of (projective) stereotype algebra such that the multiplication satisfies the following identity:
$$
(a_1\circledast b_1)\cdot (a_2\circledast b_2)=(a_1\cdot a_2)\circledast(b_1\cdot b_2),\qquad a_1,a_2\in A,\ b_1,b_2\in B.
$$
\etm

\blm\label{LM:ph:A-circledast-B->C}
Let $\alpha:A\to C$ and $\beta:B\to C$ be morphisms of stereotype algebras with commuting images:
\beq\label{alpha-betta-kommut}
\alpha(a)\cdot\beta(b)=\beta(b)\cdot\alpha(a),\qquad a\in A,\ b\in B.
\eeq
Then the map $\ph:A\circledast B\to C$, defined by the identity
\beq\label{ph:A-circledast-B->C}
\ph(a\circledast b)=\alpha(a)\cdot\beta(b),\qquad a\in A,\ b\in B,
\eeq
is a morphism of stereotype algebras. Conversely, each morphism of stereotype algebras $\ph:A\circledast B\to C$ can be uniquely represented in the form \eqref{ph:A-circledast-B->C}, where $\alpha:A\to C$ and $\beta:B\to C$ are morphisms of stereotype algebras satisfying \eqref{alpha-betta-kommut}. If $A$, $B$, $C$ are involutive algebras and  $\ph$ an involutive morphism, then $\alpha$ and $\beta$ are also involutive.
\elm
\bpr Here only the second part of the proposition is not quite evident. Suppose $\ph:A\circledast B\to C$ is a morphism of algebras. Put
$$
\alpha(a)=\ph(a\circledast 1_B),\qquad \beta(b)=\ph(1_A\circledast b),\qquad a\in A,\ b\in B.
$$
Then, first,
$$
\ph(a\circledast b)=\ph\Big((a\circledast 1_B)\cdot(1_A\circledast b)\Big)=\ph(a\circledast 1_B)\cdot\ph(1_A\circledast b)=\alpha(a)\cdot\beta(b),
$$
and, second,
$$
\alpha(a)\cdot\beta(b)=\ph(a\circledast b)=\ph\Big((1_A\circledast b)\cdot(a\circledast 1_B)\Big)=\ph(1_A\circledast b)\cdot\ph(a\circledast 1_B)=\beta(b)\cdot\alpha(a).
$$
\epr

\paragraph{Stereotype modules.}

A stereotype space $X$ over $\C$ with a given structure of left (or right) $A$-module is called a {\it stereotype $A$-module}, if the multiplication by elements of $A$ is a continuous bilinear form in the sense of definition on page \pageref{DEF:bilin-otobrazhenie}. Theorem \ref{th-7.3} implies that $X$ is a stereotype (left) module over $A$ if and only if the multiplication $\mu$ by elements of $A$ can be continuously factored through the projective stereotype tensor product
$$
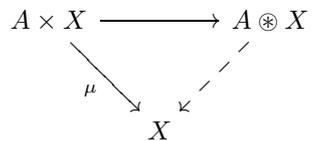
\begin{diagram}
\node{A\times X} \arrow[2]{e,t}{} \arrow{se,b}{\mu} \node[2]{A\circledast X}
\arrow{sw,b,--}{}
\\
\node[2]{X}
\end{diagram}.
$$

\begin{ex}\label{ex-11.1} Each stereotype space $X$ is a stereotype left module over the stereotype algebra $\mathcal  L(X)$ (from Example \ref{ex-10.2}).
\end{ex}

\btm [on representation]\label{th-11.2} \cite[Theorem 11.2]{Akbarov}
Let $A$ be a stereotype algebra. A stereotype space $X$ with the structure of left (right) $A$-module is a stereotype $A$-module if and only if the operation of multiplication by elements of $A$ defines a continuous homomorphism (respectively, antihomomorphism) of $A$ into $\mathcal L(X)$.
\etm

The classes $_A{\tt Ste}$ and ${\tt Ste}_A$ of left and right stereotype modules over a stereotype algebra $A$ form categories with continuous $A$-linear maps as morphisms.

\medskip
\centerline{\bf Properties of the categories $_A{\tt Ste}$ and ${\tt Ste}_A$ of stereotype modules:\footnote{Theorems 11.11, 11.17, 12.3 in \cite{Akbarov}.}}

\bit{\it

\item[$1^\circ$.] $_A{\tt Ste}$ and ${\tt Ste}_A$ are pre-Abelian categories.

\item[$2^\circ$.]\label{TH:STE-polnaya-kategoriya} $_A{\tt Ste}$ and ${\tt Ste}_A$ are complete:
 each covariant (and each contravariant) system has an injective and a projective limit.

\item[$3^\circ$.] $_A{\tt Ste}$ and ${\tt Ste}_A$ are enriched categories over the monoidal category ${\tt Ste}$.

}\eit

We will need the following notation. Let $X$ be a left stereotype module over a stereotype algebra $A$. If $M$ is a subspace in $A$ and $N$ a subspace in $X$, then we put
\beq\label{DEF:M-cdot-N} M\cdot N=\left\{ \sum_{i=1}^k m_i\cdot n_i ;\ m_i\in
M, n_i\in N, k\in \Bbb N\right\} \eeq

Obviously, for any $L,M,N\subseteq A$
\beq\label{(LM)N=L(MN)}
(L\cdot M)\cdot N=L\cdot (M\cdot N)
\eeq
and
\beq\label{(M^N)^=(MN)^}
\overline{M\cdot N}=\overline{\overline{M} \cdot N}= \overline{M\cdot
\overline{N}}=\overline{\overline{M}\cdot \overline{N}}
\eeq

\blm\label{LM:I-ideal=>ph(I)-ideal}
Let $\ph :A\to B$ be a morphism of stereotype algebras. Then
\bit{
\item[(i)] for any subsets $M,N\subseteq A$ we have
$$
\ph(\overline{M\cdot N})\subseteq \overline{\ph (M)\cdot \ph (N)}
$$
\item[(ii)] for any left (right) ideal $I$ in $A$ the set $\overline{\ph(I)}$ is a left (right) ideal in $\overline{\ph(A)}$.
}\eit
\elm
\bpr The part (i) was noticed in  \cite[Lemma 13.10]{Akbarov}, so it remains to prove (ii).
Take $b\in \overline{\ph(A)}$ and $y\in\overline{\ph(I)}$, then
$$
b\underset{\infty\gets i}{\longleftarrow} \ph(a_i),\qquad
y\underset{\infty\gets j}{\longleftarrow} \ph(x_j),
$$
for some nets $a_i\in A$ and $x_j\in I$. Hence
$$
b\cdot y\underset{\infty\gets i}{\longleftarrow} \ph(a_i)\cdot y\underset{\infty\gets j}{\longleftarrow} \ph(a_i)\cdot \ph(x_j)=\ph(a_i\cdot x_j).
$$
and therefore $b\cdot y\in \overline{\ph(I)}$.
\epr

\subsection{Involutive algebras}

\paragraph{Involutions on stereotype algebras and coalgebras.}

Let us denote by $\br$ the braiding morphism
\beq\label{DEF:zauzlivanie}
\br:a\otimes b\mapsto b\otimes a
\eeq
in an arbitrary braided category, in praticular in categories $(\tt{Ste},\circledast)$ and $(\tt{Ste},\odot)$. Note that in $(\tt{Ste},\circledast)$ and in $(\tt{Ste},\odot)$ the operation dual to braiding is also a braiding:
\beq\label{br^star=br}
\br^\star=\br
\eeq
(in $(\tt{Ste},\circledast)$ the braiding is defined by the formula \eqref{DEF:zauzlivanie}, and in $(\tt{Ste},\odot)$ it is defined exactly as the operator dual to the braiding in $(\tt{Ste},\circledast)$).

Each stereotype algebra $A$ in the sense of definition on page \pageref{DEF:ster-alg} is a monoid $(A,\mu,\iota))$ in the category $(\tt{Ste},\circledast)$ (to emphasize this we will say that $A$ is an algebra over $\circledast$). Let us call
\bit{

\item[---] an {\it antipode} in $A$ an arbitrary morphism $\sigma:A\to A$ such that
\beq\label{sigma-circ-mu=mu-circ-sigma-circledast-sigma-circ-br}
\sigma\circ\ \mu=\mu\circ \sigma\circledast \sigma\circ\br;
\eeq
in usual terms this means that $\sigma$ is an antiautomorphism of $A$:
$$
\sigma(x\cdot y)=\sigma(y)\cdot\sigma(x),\qquad x,y\in A.
$$

\item[---] an {\it involution} in $A$ an arbitrary involution $\bullet:A\to A$ in $A$ as a stereotype space such that
\beq\label{*-circ-mu=mu-circ-*-circledast-*-circ-op}
\bullet\circ\ \mu=\mu\circ \bullet\circledast \bullet\circ\br;
\eeq
this means that
\beq\label{(a-cdot-b)^*=b^*-cdot-a^*}
(x\cdot y)^\bullet=y^\bullet\cdot x^\bullet,\qquad x,y\in A.
\eeq
}\eit

By analogy we define the antipode and the involution on a stereotype algebra over $\odot$, i.e. in a monoid $(A,\mu,\iota)$ in the category $(\tt{Ste},\odot)$.

\brem If $\bullet$ is an involution on a stereotype algebra $(A,\mu,\iota)$ (no matter, in $(\tt{Ste},\circledast)$ or in $(\tt{Ste},\odot)$), then
    \beq\label{bullet-circ-iota-circ-bullet=iota}
    \bullet\circ\ \iota\circ\bullet_{\C}=\iota
    \eeq
or, in usual notations,
\beq\label{1^bullet=1}
1^\bullet=1.
\eeq
Indeed, independently on tensor product  -- $\circledast$ or $\odot$ -- we have:
$$
1^\bullet=1^\bullet\cdot 1=(1^\bullet\cdot 1^{\bullet\bullet})^\bullet=(1^\bullet\cdot 1)^\bullet=(1^\bullet)^\bullet=1.
$$
\erem

Let $K$ be a stereotype coalgebra over $\circledast$, i.e. a comonoid $(K,\varkappa,\e)$ in the category $(\tt{Ste},\circledast)$. Let us call
\bit{

\item[---] an {\it antipode} in $K$ an arbitrary morphism $\sigma:K\to K$ such that
\beq\label{varkappa-circ-sigma=br-circ-sigma-circledast-sigma-circ-varkappa}
\varkappa\circ\ \sigma=\br\circ\ \sigma\circledast \sigma\circ\varkappa,
\eeq

\item[---] an {\it involution} in $K$ an arbitrary involution $\bullet:K\to K$ in $K$ as a stereotype space, such that
\beq\label{varkappa-circ-*=*-circledast-*-circ-varkappa}
\varkappa\circ\ \bullet=\bullet\circledast \bullet\circ\varkappa,
\eeq
(without braiding $\br$).
}\eit

By analogy we define antipode and involution on a stereotype coalgebra over $\odot$, i.e. in a comonoid  $(K,\varkappa,\e)$ in the category $(\tt{Ste},\odot)$.

\brem If $\bullet$ is an involution on a stereotype coalgebra $(K,\varkappa,\e)$ (no matter, in  $(\tt{Ste},\circledast)$ or $(\tt{Ste},\odot)$), then
    \beq\label{bullet-circ-e-circ-bullet=e}
    \bullet_{\C}\circ\ \e\circ\bullet=\e
    \eeq
This becomes clear after the passage to the dual algebra $(K^\star,\varkappa^\star,\e^\star)$: then $\e^\star$ turns into the unit, and \eqref{bullet-circ-iota-circ-bullet=iota} will be equivalent to \eqref{bullet-circ-e-circ-bullet=e}.
\erem

Let us say that an antipode $\sigma$ and an involution $\bullet$ are {\it interconsistent} in an algebra (or in a coalgebra) $A$, if $\sigma$ is a spin for $\bullet$ (i.e. if \eqref{bullet-circ-sigma} and \eqref{sigma-circ-bullet} hold). In contrast to Proposition \ref{PROP:podkrutka}(iv), however, the composition $\bullet\circ\sigma$ is not necessary an involution on the algebra (or coalgebra) $A$, since the equality \eqref{*-circ-mu=mu-circ-*-circledast-*-circ-op} (respectively, \eqref{varkappa-circ-*=*-circledast-*-circ-varkappa}) may be not true.

 \vglue10pt
\centerline{\bf Properties of antipodes and involutions:}
 \vglue10pt

\bit{\it

\item[$1^\circ.$] If $\sigma:A\to A$ is an antipode in the algebra $(A,\mu,\iota)$ in the category  $(\tt{Ste},\circledast)$ (or in $(\tt{Ste},\odot)$), then the dual mapping $\sigma^\star:A^\star\to A^\star$ is an antipode in the coalgebra $(A^\star,\mu^\star,\iota^\star)$.

\item[$2^\circ.$] If $\sigma:K\to K$ is an antipode in the coalgebra $(K,\varkappa,\e)$ in the category  $(\tt{Ste},\circledast)$ (or in $(\tt{Ste},\odot)$), then the dual mapping $\sigma^\star:K^\star\to K^\star$ is an antipode in the algebra $(K^\star,\varkappa^\star,\e^\star)$.

\item[$3^\circ.$] Let $\sigma$ and $\bullet$ be interconsistent antipode and involution on the algebra $(A,\mu,\iota)$ in $(\tt{Ste},\circledast)$ (or in $(\tt{Ste},\odot)$). Then the formula \eqref{f^bullet(a)=overline(f(sigma(a)^bullet))} defines an involution on the dual coalgebra  $(A^\star,\mu^\star,\iota^\star)$, and it will be interconsistent with the dual antipode $\sigma^\star$.

\item[$4^\circ.$]\label{4^0:inv-v-sopryazh-alg} Let $\sigma$ and $\bullet$ be interconsistent antipode and involution on a coalgebra $(K,\varkappa,\e)$ in $(\tt{Ste},\circledast)$ (or in $(\tt{Ste},\odot)$).
    Then the formula \eqref{f^bullet(a)=overline(f(sigma(a)^bullet))} defines an involution on the dual algebra $(K^\star,\varkappa^\star,\e^\star)$, and it will be interconsistent with the dual antipode $\sigma^\star$.
}\eit

\bpr
The properties $1^\circ$ and $2^\circ$ are evident. Let us prove $3^\circ$. Suppose $\sigma$ and $\bullet$ are interconsistent antipode and involution on the algebra $(A,\mu,\iota)$ in $(\tt{Ste},\circledast)$ (the case of  $(\tt{Ste},\odot)$ is considered by analogy). By Theorem \ref{PROP:f^bullet(a)=overline(f(sigma(a)^bullet))}, the formula \eqref{f^bullet(a)=overline(f(sigma(a)^bullet))} defines an involution  $\bullet^\star=\bullet_{\C}\oslash(\bullet\circ\sigma)$ in the dual space $A^\star$. We have:
\begin{multline*}
(\bullet^\star\odot\bullet^\star)\circ \mu^\star=
\Big\{\big[\bullet_{\C}\oslash(\bullet\circ\sigma)\big]\odot\big[\bullet_{\C}\oslash(\bullet\circ\sigma)\big]\Big\} \circ(\id_{\C}\circ\mu)=\eqref{bullet_X-circledast-bullet_Y->bullet_X-odot-bullet_Y}=
\Big\{\bullet_{\C}\oslash\big[(\bullet\circledast\bullet)\circ (\sigma\circ\sigma)\big]\Big\} \circ(\id_{\C}\circ\mu)=\\=
(\bullet_{\C}\circ\id_{\C})\oslash \big[\mu\circ(\bullet\circledast\bullet)\circ (\sigma\circ\sigma)\big]=
\bullet_{\C}\oslash \big[\mu\circ(\bullet\circledast\bullet)\circ (\sigma\circ\sigma)\big]=\eqref{*-circ-mu=mu-circ-*-circledast-*-circ-op}=
\bullet_{\C}\oslash \big[\bullet\circ\mu\circ\br\circ (\sigma\circ\sigma)\big]=\\=
\bullet_{\C}\oslash \big[\bullet\circ\mu\circ (\sigma\circ\sigma)\circ\br\big]=
\eqref{sigma-circ-mu=mu-circ-sigma-circledast-sigma-circ-br}=
\bullet_{\C}\oslash (\bullet\circ\sigma\circ\mu)=
(\id_{\C}\oslash\mu)\circ\big[\bullet_{\C}\oslash (\bullet\circ\sigma)\big]=
\mu^\star\circ \bullet^\star.
\end{multline*}

Let us prove $4^\circ$. Suppose $\sigma$ and $\bullet$ are interconsistent antipode and involution on a coalgebra $(K,\varkappa,\e)$ in $(\tt{Ste},\circledast)$. Again, the formula \eqref{f^bullet(a)=overline(f(sigma(a)^bullet))} defines an involution $\bullet^\star=\bullet_{\C}\oslash(\bullet\circ\sigma)$ in the dual space $K^\star$. We have:
\begin{multline*}
\varkappa^\star\circ \bullet^\star\odot\bullet^\star\circ\br=\eqref{br^star=br}=
(\id_{\C}\oslash\varkappa)\circ \Big\{\big[\bullet_{\C}\oslash(\bullet\circ\sigma)\big]\odot\big[\bullet_{\C}\oslash(\bullet\circ\sigma)\big]\Big\}
\circ(\id_{\C}\oslash \br)=\eqref{bullet_X-circledast-bullet_Y->bullet_X-odot-bullet_Y}=\\=
(\id_{\C}\oslash\varkappa)\circ \Big\{\bullet_{\C}\oslash\big[(\bullet\circledast\bullet)\circ (\sigma\circ\sigma)\big]\Big\} \circ(\id_{\C}\oslash \br)=
(\id_{\C}\circ\bullet_{\C}\circ\id_{\C})\oslash\big[\br\circ(\bullet\circledast\bullet)\circ (\sigma\circ\sigma)\circ\varkappa\big]=
\\=
\bullet_{\C}\oslash\big[(\bullet\circledast\bullet)\circ \br\circ (\sigma\circ\sigma)\circ\varkappa\big]=
\eqref{varkappa-circ-sigma=br-circ-sigma-circledast-sigma-circ-varkappa}=
\bullet_{\C}\oslash\big[(\bullet\circledast\bullet)\circ \varkappa\circ \sigma\big]=
\eqref{varkappa-circ-*=*-circledast-*-circ-varkappa}=
\bullet_{\C}\oslash\big[\varkappa\circ\bullet\circ \sigma\big]=\\=
\big[\bullet_{\C}\oslash(\bullet\circ\sigma)\big]\circ(\id_{\C}\oslash\varkappa)=\bullet^\star\circ\varkappa^\star.
\end{multline*}
\epr

\btm Suppose $A$ and $B$ are two stereotype algebras with involutions $\bullet_A$ and $\bullet_B$. Then the formula
\beq\label{(a-circledast-b)^bullet=a^bullet-circledast-b^bullet}
\bullet_A\circledast\bullet_B=
\bullet_{\C}\oslash((\bullet_{\C}\oslash \bullet_A)\oslash\bullet_B)
\eeq
defines an involution on the stereotype algebra $A\circledast B$.
\etm
\bpr By Theorem \ref{TH:inv-v-Y-oslash-X}(iii) this formula defines involution on $A\circledast B$ as a stereotype space. The identity \eqref{(a-cdot-b)^*=b^*-cdot-a^*} can be verified on elementary tensors $a\circledast b, a'\circledast b'\in A\circledast B$:
\begin{multline*}
(\bullet_A\circledast\bullet_B)\big((a\circledast b)\cdot(a'\circledast b')\big)=
(\bullet_A\circledast\bullet_B)\big((a\cdot a')\circledast (b\cdot b')\big) =\eqref{bullet_X-circledast-bullet_Y(x,y)}=\bullet_A(a\cdot a')\circledast\bullet_B(b\cdot b')=\\=
\bullet_A(a')\cdot\bullet_A(a)\circledast\bullet_B(b')\cdot\bullet_B(b)=
(\bullet_A(a')\circledast\bullet_B(b'))\cdot(\bullet_A(a)\circledast\bullet_B(b))=
\eqref{bullet_X-circledast-bullet_Y(x,y)}=\\=
(\bullet_A\circledast\bullet_B)(a'\circledast b')\cdot(\bullet_A\circledast\bullet_B)(a\circledast b)
\end{multline*}
\epr

\paragraph{Involutive projective stereotype algebras.}

 \bit{
\item A projective stereotype algebra $A$ with a given involution $\bullet$ (on $A$ as a stereotype space and such that \eqref{(a-cdot-b)^*=b^*-cdot-a^*} holds) is called {\it involutive projective stereotype algebra}.

\item The class of all involutive (projective) stereotype algebras is denoted by $\InvSteAlg$. It forms a category where morphisms are morphisms of stereotype algebras $\ph:A\to B$ which preserve involution:
\beq\label{ph(x^bullet)=ph(x)^bullet}
\ph(x^\bullet)=\ph(x)^\bullet,\qquad x\in A.
\eeq

}\eit

\bex\label{EX:C*-algebra-kak-inv-ster-algebra}
$C^*$-algebras are obvious examples of involutive stereotype algebras.
\eex

In the category $\InvSteAlg$ we are interested in the envelopes in the class of the so-called dense epimorphisms.

 \bit{
\item[$\bullet$]\label{DEF:DEpi} A homomorphism of stereotype algebras $\ph:A\to B$ is called a {\it dense epimorphism}, if $\ph(A)$ is dense in $B$ (not any epimorphism in ${\SteAlg}$ is dense). The class of all dense epimorphisms will be denoted by $\DEpi$.

Certainly, dense morphisms are epimorphisms, that is why we call them also dense epimorphisms. The class of all dense epimorphisms in the category $\InvSteAlg$ will be denoted by $\DEpi$. It is connected with the classes $\Epi$ of all epimorphisms and $\SEpi$ of strong epimorphisms through the following embeddings:
\beq\label{SEpi-subset-DEpi-subset-Epi}
\SEpi\subset\DEpi\subset\Epi.
\eeq
 }\eit

\btm\label{TH:SMono-circledcirc-DEpi=InvSteAlg}
The class $\DEpi$ of dense epimorphisms is monomorphically complementable in the category $\InvSteAlg$.
\etm
\bpr
The monomorphic complement for $\DEpi$ is the class $\SMono_{\tt Ste}$ of all monomorphisms in $\InvSteAlg$, which are strong monomorphisms in $\Ste$:
\beq\label{SMono-circledcirc-DEpi=InvSteAlg}
\SMono_{\tt Ste}\circledcirc\DEpi=\InvSteAlg
\eeq
\epr

\btm\cite[Theorem 5.37]{Akbarov-env} \label{TH:obolochki-i-nachinki-otn-plotnyh-epimorfizmov-v-Ste^circledast}
In the category ${\tt Ste}^\circledast$ every algebra $A$ has an envelope in the class $\DEpi$ of dense epimorphisms with respect to the arbitrary {\rm class} of morphisms $\varPhi$ going from $A$. If in addition $\varPhi$ differs morphisms on the outside\footnote{See definition on page \pageref{DEF:varPhi-razlich-morfizmy-snaruzhi}} in ${\tt Ste}^\circledast$, then the envelope in $\DEpi$ is also an envelope in the class $\DBim$ of all dense bimorphisms:
    $$
    \env_\varPhi^{\DEpi} A=\env_\varPhi^{\DBim} A.
    $$
\etm

\btm\cite[Theorem 5.38]{Akbarov-env}\label{TH:sushestvovanie-seti-v-Ste^circledast-dlya-DEpi}
Let $\varPhi$ be a class of morphisms in $\InvSteAlg$, that goes from\footnote{See definition on page \pageref{DEF:goes-from}.} $\SteAlg$, is a right ideal,
$$
\varPhi\circ\Mor(\SteAlg)\subseteq\varPhi,
$$
and such that $\DEpi$ pushes\footnote{See definition on page \pageref{Omega-podderzhivaet-Phi}.} $\varPhi$.
Then the classes of morphisms $\DEpi$ and $\varPhi$ define a regular envelope $\Env_\varPhi^{\DEpi}$ in $\InvSteAlg$. \etm

\paragraph{Involutive injective stereotype algebras.}

\bit{

\item Injective stereotype algebra $A$ with a given involution $\bullet$ (on $A$ as a stereotype space) such that  \eqref{*-circ-mu=mu-circ-*-circledast-*-circ-op} holds (with $\odot$ instead of $\circledast$) is called {\it involutive injective stereotype algebra}\footnote{It is impossible replace here \eqref{*-circ-mu=mu-circ-*-circledast-*-circ-op} by \eqref{(a-cdot-b)^*=b^*-cdot-a^*}, since for the tensor product $\odot$ these condition are apparently not equivelent.}.

}\eit

In contrast to Example \ref{EX:C*-algebra-kak-inv-ster-algebra}, not all $C^*$-algebras $A$ are injective stereotype algebras. We shall describe here the class of such $C^*$-algebras, since it will be useful below in the building of the continuous duality theory. Denote by $X\check{\otimes}Y$ the injective tensor product of locally convex spaces ($X\tilde{\otimes}_{\e}X$ in the notations of \cite{Pietsch} and \cite{Jarchow}).

\bit{

\item A Banach algebra $A$ is said to be {\it strict}, if its multiplication can be extended to an operator on the injective tensor product $\mu:A\check{\otimes}A\to A$.

}\eit

\btm[\cite{Aristov-1}, \cite{Aristov-2}]\label{TH:Aristov} For a $C^*$-algebra $A$ the following conditions are equivalent:
\bit{
\item[$(i)$] $A$ is a strict algebra,

\item[$(ii)$] all unitary irreducible representations of $A$ are finite dimensional and their dimensions are bounded,

\item[$(iii)$] the minimal tensor product of $A$ with any $C^*$-algebra $B$ coincides with their injective tensor product:
    $$
     A\underset{\min}{\otimes}B=A\check{\otimes}B.
    $$

}\eit
\etm

\bex\label{EX:C(K)-strogaya-algebra} The algebra ${\mathcal C}(K)$ of continuous functions on a compact set $K$ is a strict $C^*$-algebra.
\eex

Theorem \ref{TH:Aristov} implies that any strict $C^*$-algebra $A$ is a CCR-algebra, hence a GCR-algebra, and therefore it is nuclear (as a $C^*$-algebra). This means in its turn that $A$ has the approximation property \cite{Cooper} (as a locally convex space). As a corollary (see \cite[Theorem 7.21]{Akbarov}) we have

\btm\label{TH:A-min-B=A-odot-B}
For a $C^*$-algebra $A$ the conditions (i)-(iii) of Theorem \ref{TH:Aristov} are equivalent to the condition
\bit{
\item[$(iii)'$] the minimal tensor product of $A$ with any $C^*$-algebra $B$ coincides with their maximal tensor product (as $C^*$-algebras), with their injective locally convex tensor product, and with thir injective stereotype tensor product:
    \beq\label{A-min-B=A-odot-B}
A\underset{\min}{\otimes}B=A\underset{\max}{\otimes}B=A\check{\otimes}B=A\odot B.
    \eeq
}\eit
\etm

\bcor Each strict $C^*$-algebra $A$ is an injective stereotype algebra.
\ecor
\bpr
By definition of strict algebra, the multiplication on $A$ can be extended to an operator $\mu:A\check{\otimes}A\to A$, hence, by Theorem \ref{TH:A-min-B=A-odot-B}, to an operator $\mu:A\odot A\to A$. The diagrams \eqref{monoid} are true for tensor product $\check{\otimes}$, since the algebraic tensor product $\otimes$ is dense in it. As a corollary, they are true for $\odot$, since $A$ has the approximation property, and thus by \cite[Theorem 7.21]{Akbarov}, $A\check{\otimes}A=A\odot A$ and $A\check{\otimes}A\check{\otimes}A=A\odot A\odot A$. The equality \eqref{*-circ-mu=mu-circ-*-circledast-*-circ-op} (with $\odot$ instead of $\circledast$) is true because of the density of $A\otimes A$ in $A\check{\otimes}A=A\odot A$.
\epr

\paragraph{Involution on stereotype Hopf algebras.}
Let us remind (see e.g. \cite{Akbarov-stein-groups}), that the notion of Hopf algebra can be defined in any symmetric monoidal category.

\bprop
If $(H,\mu,\iota,\varkappa,\e,\sigma)$ is a Hopf algebra in a symmetric monoidal category $\tt{K}$ with tensor product $\otimes$, then the multiplication and the comultiplication are connected with the antipode trough the identities \eqref{sigma-circ-mu=mu-circ-sigma-circledast-sigma-circ-br} and \eqref{varkappa-circ-sigma=br-circ-sigma-circledast-sigma-circ-varkappa}:
\beq\label{sigma-circ-mu=mu-circ-sigma-otimes-sigma-circ-op}
\sigma\circ\mu=\mu\circ\sigma\otimes\sigma\circ\br,\qquad
\varkappa\circ\sigma=\br\circ\ \sigma\otimes\sigma\circ\varkappa
\eeq
\eprop
\bpr
This is proved for example in \cite[Figure 9.14]{Majid}.
\epr

An {\it involution on a stereotype Hopf algebra} $(H,\mu,\iota,\varkappa,\e,\sigma)$ (over the tensor product $\circledast$) is an arbitrary mapping $a\in H\mapsto a^\bullet\in H$, which is an involution on $H$ as a stereotype space and satisfies the identities \eqref{*-circ-mu=mu-circ-*-circledast-*-circ-op} and \eqref{varkappa-circ-*=*-circledast-*-circ-varkappa} (and, as a corollary, the identities \eqref{bullet-circ-iota-circ-bullet=iota} and \eqref{bullet-circ-e-circ-bullet=e}).

By analogy we an define involution on a stereotype Hopf algebra over $\odot$.

\bprop
Let $\bullet$ be an involution on a stereotype Hopf algebra $H$ over $\circledast$ (respectively,  $\odot$). Then
\bit{
\item[(i)] the involution on $H$ is interconsistent with the antipode (i.e. \eqref{bullet-circ-sigma} holds) if and only if the antipode $\sigma$ is invertible in $H\oslash H$,

\item[(ii)] the formula \eqref{f^bullet(a)=overline(f(sigma(a)^bullet))} defines an involution on the dual Hopf algebra $H^\star$ over $\odot$ (respectively, over $\circledast$).
}\eit
\eprop
\bpr
If $\sigma$ is invertible, then we can notice that $\bullet\circ\sigma^{-1}\circ\bullet$ is also an antipode in $H$, and since the antipode in a Hopf algebra is unique, \eqref{sigma^(-1)} must hold. This proves (i), and (ii) follows from $3^\circ$ and $4^\circ$ on page \pageref{f^bullet(a)=overline(f(sigma(a)^bullet))}.
\epr

\subsection{Spectrum and tangent space}

\paragraph{Spectrum.}

 \bit{
 \item[$\bullet$] The {\it complex spectrum}\label{DEF:Spec_C(G)} ${\C}\Spec(A)$ of a stereotype algebra $A$ over $\C$ is a set of its {\it complex characters}, i.e. (continuous and unital) homomorphisms $s:A\to \C$. We endow this set with the topology of uniform convergence on totally bounded sets in $A$.

 \item[$\bullet$] Similarly, the {\it real spectrum}\label{DEF:RSpec(G)} ${\R}\Spec(A)$ of a stereotype algebra $A$ over $\R$ is the set of its {\it real characters}, i.e. (continuous and unital) homomorphisms $s:A\to \R$. We endow this set with the topology of uniform convergence on totally bounded sets in $A$.

 }\eit

\btm\label{TH:Real(CSpec[A])=RSpec[Re-A]} Let $\bullet$ be an involution on a stereotype algebra $A$ over $\C$, and suppose on the dual space $A^\star$ the involution is defined by formula \eqref{DEF:inv-v-X^star}. Then formulas \eqref{f-mapsto-f|_Real X} define isomorphism of spectra:
\beq\label{Real(CSpec[A])=RSpec[Re-A]}
\Real({\C}\Spec[A])\cong {\R}\Spec[\Real A].
\eeq
\etm

\bit{
 \item[$\bullet$] The space in \eqref{Real(CSpec[A])=RSpec[Re-A]} is called the {\it involutive spectrum}\label{DEF:Spec_R(G)} of the involutive stereotype algebra $A$ over $\C$, and is denoted by  $\Spec(A)$:
\beq\label{DEF:Spec[A]}
\Spec(A)=\Real({\C}\Spec[A])\cong {\R}\Spec[\Real A].
\eeq
     It consists of {\it involutive characters}, i.e. (continuous, unital and) involutive\footnote{Involutive operators were defines on page \pageref{DEF:involutivnyi-operator}.} homomorphisms $s:A\to \C$, and is endowed with the topology of uniform convergence on totally bounded sets in $A$.
 }\eit

\bpr
Suppose $t\in\Real({\C}\Spec[A])$ and
$$
s=P(t)=t|_{\Real A}.
$$
Then by Theorem \ref{PROP:f^bullet(a)=overline(f(sigma(a)^bullet))}, $s\in (\Real A)^\star_{\R}$, i.e. $s:\Real A\to\R$. Moreover, $s\in {\R}\Spec[\Real A]$, since $u,v\in \Real A$ imply $u\cdot v\in \Real A$, and therefore
$$
s(u\cdot v)=t(u\cdot v)=t(u)\cdot t(v)=s(u)\cdot s(v).
$$
Conversely, if $s\in {\R}\Spec[\Real A]$ and
$$
t(x)=Q(s)(x)=s(\Real x)+i\cdot s(\Im x),
$$
then by Theorem \ref{PROP:f^bullet(a)=overline(f(sigma(a)^bullet))}, $t\in\Real(A^\star_{\C})$. Moreover, $t\in\Real({\C}\Spec[A])$, since for any $x,y\in A$
\begin{multline*}
t(x\cdot y)=s(\Real(x\cdot y))+i\cdot s(\Im(x\cdot y))=
s(\Real x\cdot\Real y-\Im x\cdot\Im y)+i\cdot s(\Im x\cdot\Real y+\Real x\cdot\Im y)=\\=
s(\Real x\cdot\Real y)-s(\Im x\cdot\Im y)+i\cdot s(\Im x\cdot\Real y)+i\cdot s(\Real x\cdot\Im y)=\\=
s(\Real x)\cdot s(\Real y)-s(\Im x)\cdot s(\Im y)+i\cdot s(\Im x)\cdot s(\Real y)+i\cdot s(\Real x)\cdot s(\Im y)=\\=
\Big(s(\Real x)+i\cdot s(\Im x)\Big)\cdot\Big(s(\Real y)+i\cdot s(\Im y)\Big)=t(x)\cdot t(y).
\end{multline*}
\epr

\blm\label{LM:overline(ph(I_t)-cdot-B)=B} Let $\ph :A\to B$ be an involutive dense\footnote{I.e. $\ph(A)$ is dense in $B$.} epimorphism of an involutive commutative stereotype algebra $A$ into a $C^*$-algebra $B$. Then
\bit{

\item[(i)] for any point $t\in\Spec(B)$
\beq\label{overline(ph(Ker(t-circ-ph)))=Ker_t}
\overline{\ph\big(\Ker (t\circ\ph)\big)}=\Ker t
\eeq

\item[(ii)] for any point $s\in\Spec(A)\setminus\big(\Spec(B)\circ\ph \big)$
\beq\label{overline(ph(Ker_s))=B}
\overline{\ph(\Ker s)}=B
\eeq
}\eit
\elm
\bpr
1. Take $t\in\Spec(B)$. First of all,
$$
a\in \Ker (t\circ\ph)\quad\Longrightarrow\quad t(\ph(a))=(t\circ\ph)(a)=0\quad\Longrightarrow\quad \ph(a)\in\Ker t.
$$
Hence $\ph\big(\Ker (t\circ\ph)\big)\subseteq\Ker t$, and therefore $\overline{\ph\big(\Ker (t\circ\ph)\big)}\subseteq\Ker t$.

Conversely, suppose $b\in \Ker t$. Since $\ph(A)$ is dense in $B$, for each $\e>0$ one can find $a_{\e}\in A$ such that
$$
\norm{b-\ph(a_{\e})}<\e.
$$
Put
$$
a'_{\e}=a_{\e}-t\big(\ph(a_{\e})\big)\cdot 1_A.
$$
Then, first,
$$
t\big(\ph(a'_{\e})\big)=t\big(\ph(a_{\e})\big)-t\big(\ph(a_{\e})\big)\cdot t\big(\ph(1_A)\big)=0,
$$
i.e. $a'_{\e}\in\Ker (t\circ\ph)$. And, second,
$$
\norm{\ph(a_{\e})-\ph(a'_{\e})}=\norm{\ph(a_{\e})-\ph(a_{\e})+t\big(\ph(a_{\e})\big)\cdot\ph(1_A)}=\abs{t\big(\ph(a_{\e})\big)}=
\Big|\underbrace{t(b)}_{\scriptsize\begin{matrix}\text{\rotatebox{90}{$=$}}\\ 0\end{matrix}}-t\big(\ph(a_{\e})\big)\Big|\le \norm{t}\cdot \norm{b-\ph(a_{\e})}<\e.
$$
Thus
$$
\norm{b-\ph(a'_{\e})}\le \norm{b-\ph(a_{\e})}+\norm{\ph(a_{\e})-\ph(a'_{\e})}<2\e.
$$
We see that for each $b\in \Ker t$ and $\e>0$ there is an element $a'_{\e}\in\Ker(t\circ\ph)$ such that $\norm{b-\ph(a'_{\e})}<2\e$. This proves the inclusion $\overline{\ph\big(\Ker (t\circ\ph)\big)}\supseteq\Ker t$.

2. Take $s\in\Spec(A)\setminus\big(\Spec(B)\circ\ph \big)$. For any point $t\in \Spec(B)$ there is an element $x_t\in A$ that separates $s$ from $t\circ\ph$:
$$
s(x_t)\ne t(\ph(x_t)).
$$
As a corollary the element $a_t=x_t-s(x_t)\cdot 1_A$ must possess the following properties:
\beq\label{s(a_t)=0}
s(a_t)=s(x_t)-s(x_t)\cdot s(1_A)=s(x_t)-s(x_t)=0
\eeq
and
\beq\label{t(ph(a_t))-ne-0}
t(\ph(a_t))=t(\ph(x_t))-s(x_t)\cdot t(\ph(1_A))=t(\ph(x_t))-s(x_t)\ne 0.
\eeq
From \eqref{t(ph(a_t))-ne-0} we see that the sets
$$
U_t=\{r\in\Spec(A):\ r(a_t)\ne 0 \}
$$
cover the compact set $\Spec(B)\circ\ph$. Hence, they have finite subcovering $U_{t_1},...,U_{t_n}$:
$$
\bigcup_{i=1}^n U_{t_i}\supseteq \Spec(B)\circ\ph.
$$
Put
$$
a=\sum_{i=1}^n a_{t_i}\cdot a_{t_i}^\bullet.
$$
This element vanishes in the point $s$, since due to \eqref{s(a_t)=0} all $a_{t_i}$ vanish in $s$:
$$
s(a)=\sum_{i=1}^n s(a_{t_i})\cdot \overline{s(a_{t_i})}=0.
$$
Therefore, $a\in \Ker s$, and
$$
\ph(a)\in \overline{\ph(\Ker s)}.
$$
On the other hand, the element $\ph(a)$ in non-zero everywhere on $\Spec(B)$, since in each point $t\in\Spec(B)$ some element $\ph(a_{t_i})$ is non-zero:
$$
t(\ph(a))=\sum_{i=1}^n t(\ph(a_{t_i}))\cdot \overline{t(\ph(a_{t_i}))}=\sum_{i=1}^n \abs{t(\ph(a_{t_i}))}^2>0.
$$
In other words, $\ph(a)$ is non-zero everywhere on the spectrum $\Spec(B)$ of the commutative $C^*$-algebra $B\cong C(\Spec(B))$. As a corollary, it is invertible in $B$:
$$
b\cdot\ph(a)=1_B
$$
for some $b\in B$. Since $\ph(a)$ belongs to a left ideal $\overline{\ph(\Ker s)}$ of the algebra $B=\overline{\ph(A)}$ (by Lemma \ref{LM:I-ideal=>ph(I)-ideal}(ii)), we conclude that $1_B$ also belongs to $\overline{\ph(\Ker s)}$. Thus,
$$
\overline{\ph(\Ker s)}\supseteq B\cdot\overline{\ph(\Ker s)}\supseteq B\cdot 1_B=B.
$$
\epr

\bcor\label{COR:overline(ph(I_t)-cdot-B)=B} Let $\ph :A\to B$ be an involutive homomorphism of involutive stereotype algebras, where $A$ is commutative, and $B$ is a $C^*$-algebra. Then for each point $s\in\Spec(A)\setminus\Big(\Spec\big(\overline{\ph(A)}\big)\circ\ph\Big)$
\beq
\overline{\ph(\Ker s)\cdot B}=B
\eeq
\ecor
\bpr By Lemma \ref{LM:overline(ph(I_t)-cdot-B)=B}(ii), the set $\overline{\ph(\Ker s)}$ contains the unit of the algebra $\overline{\ph(A)}$, which coincides with the unit of $B$:
$$
\overline{\ph(\Ker s)}\owns 1_B.
$$
Hence,
$$
\overline{\ph(\Ker s)\cdot B}\supseteq \overline{\ph(\Ker s)}\cdot B\supseteq 1_B\cdot B=B.
$$
\epr

\paragraph{Tangent and cotangent spaces.}

 \bit{
 \item[$\bullet$]
A {\it complex} (respectively, a {\it real}) {tangent vector} to a stereotype algebra $A$ over $\C$ (respectively, over $\R$) in a point $s\in\C\Spec(A)$ (respectively, $s\in\R\Spec(A)$) is a linear continuous functional  $\tau:A\to\C$ ($\tau:A\to\R$), satisfying the identity:
 \beq\label{DEF:kasatelnyj-vektor}
 \tau(a\cdot b)=s(a)\cdot\tau(b)+\tau(a)\cdot s(b), \qquad
 a,b\in A
 \eeq
The set of all complex (respectively, real) tangent vectors to $A$ in the point $s$ is called the {\it complex} ({\it real}) {\it tangent space} to $A$ in the point
$s\in\C\Spec(A)$ and is denoted by ${\C}T_s[A]$ (respectively, ${\R}T_s[A]$). It is endowed with the topology, which is the pseudosaturation of the topology of uniform convergence on totally bounded sets in $A$. As a corollary, ${\C}T_s[A]$ (respectively, ${\R}T_s[A]$) is a stereotype space and an immediate subspace in $A^\star$.

 }\eit

\brem
Let us note that each tangent vector vanishes on the identity
 \beq\label{tau(1)=0}
 \tau(1_A)=0,\qquad
\tau\in T_s[A],
 \eeq
since
$$
\tau(1_A)=\tau(1_A\cdot 1_A)=
\underbrace{s(1_A)}_{\scriptsize\begin{matrix}\text{\rotatebox{90}{$=$}}\\
1\end{matrix}}\cdot\tau(1_A)+\tau(1_A)\cdot
\underbrace{s(1_A)}_{\scriptsize\begin{matrix}\text{\rotatebox{90}{$=$}}\\
1\end{matrix}}=2\tau(1_A).
$$
\erem

\brem
If $A$ is a stereotype algebra with an involution $\bullet$, and $\tau$ is its complex tangent vector in a point $s\in\C\Spec[A]$, then the formula
 \beq\label{DEF:invol-v-CT_s[A]}
\tau^\bullet(a)=\overline{\tau(a^\bullet)},\qquad a\in A.
 \eeq
defines a tangent vector to $A$ in the point $s^\bullet$. In the special case when $s$ is an involutive character, i.e. $s^\bullet=s$, then $\tau^\bullet\in\C T_s[A]$, i.e. \eqref{DEF:invol-v-CT_s[A]} defines an involution (in the sense of \eqref{DEF:involutsiya}) on the complex tangent space ${\C}T_s[A]$.
\erem

\btm\label{TH:Real(CT_s[A])=RT_s[Re-A]} Let $\bullet$ be an involution on a stereotype algebra $A$ over $\C$, and let the dual space $A^\star$ is endowed with the involution \eqref{DEF:inv-v-X^star}. Then for any point $s\in\Spec(A)$  formulas \eqref{f-mapsto-f|_Real X} define an isomorphism of tangent spaces:
\beq\label{Real(CT_s[A])=RT_s[Re-A]}
\Real({\C}T_s[A])\cong {\R}T_s[\Real A].
\eeq
\etm

\bit{
 \item[$\bullet$] The space in  \eqref{Real(CT_s[A])=RT_s[Re-A]} will be called the {\it involutive tangent space}\label{DEF:T_s[A]} of an involutive stereotype algebra $A$ in the point $s\in\Spec[A]$, and is denoted by  $T_s(A)$:
\beq\label{DEF:T_s[A]}
T_s[A]=\Real({\C}T_s[A])\cong {\R}T_s[\Real A].
\eeq
     Certainly, it consists of {\it involutive tangent vectors}, i.e. functionals $\tau:A\to \C$, satifsying (apart from \eqref{DEF:kasatelnyj-vektor}) the condition
 \beq\label{DEF:kasatelnyj-vektor-inv}
 \tau(a^\bullet)=\overline{\tau(a)}, \qquad a\in A.
 \eeq
The space $T_s[A]$ is endowed with the topology, which is the pseudosaturation of the topology of uniform convergence on totally bounded sets in $A$.
 }\eit

\bpr
Take $\tau\in\Real({\C}T_s[A])$ and
$$
\sigma=P(\tau)=\tau|_{\Real A}.
$$
By Theorem \ref{PROP:f^bullet(a)=overline(f(sigma(a)^bullet))}, $\sigma\in (\Real A)^\star_{\R}$, i.e. $\sigma:\Real A\to\R$. Moreover, $\sigma\in {\R}T_s[\Real A]$, since $u,v\in \Real A$ implies $u\cdot v\in \Real A$, and therefore
$$
\sigma(u\cdot v)=\tau(u\cdot v)=u(s)\cdot\tau(v)+\tau(u)\cdot v(s)=u(s)\cdot\sigma(v)+\sigma(u)\cdot v(s).
$$
Conversely, if $\sigma\in {\R}T_s[\Real A]$ and
$$
\tau(x)=Q(\sigma)(x)=\sigma(\Real x)+i\cdot\sigma(\Im x)
$$
then by Theorem \ref{PROP:f^bullet(a)=overline(f(sigma(a)^bullet))}, $\tau\in\Real(A^\star_{\C})$. Moreover, $\tau\in\Real({\C}T_s[A])$, since for each $x,y\in A$
\begin{multline*}
\tau(x\cdot y)=\sigma(\Real(x\cdot y))+i\cdot\sigma(\Im(x\cdot y))=
\sigma(\Real x\cdot\Real y-\Im x\cdot\Im y)+i\cdot\sigma(\Im x\cdot\Real y+\Real x\cdot\Im y)=\\=
\sigma(\Real x\cdot\Real y)-\sigma(\Im x\cdot\Im y)+i\cdot\sigma(\Im x\cdot\Real y)+i\cdot\sigma(\Real x\cdot\Im y)=\\=
\Real x(s)\cdot\sigma(\Real y)+\sigma(\Real x)\cdot \Real y(s)
-\Im x(s)\cdot\sigma(\Im y)-\sigma(\Im x)\cdot \Im y(s)+\\
+i\cdot \Im x(s)\cdot \sigma(\Real y)+i\cdot \sigma(\Real x)\cdot \Im y(s)
+i\cdot \Real x(s)\cdot\sigma(\Im y)+i\cdot\sigma(\Real x)\cdot\Im y(s)=\\=
\Real x(s)\cdot\sigma(\Real y)-\Im x(s)\cdot\sigma(\Im y)
+i\cdot \Im x(s)\cdot \sigma(\Real y)+i\cdot \Real x(s)\cdot\sigma(\Im y)+\\+
\sigma(\Real x)\cdot \Real y(s)
-\sigma(\Im x)\cdot \Im y(s)+i\cdot \sigma(\Real x)\cdot \Im y(s)
+i\cdot\sigma(\Real x)\cdot\Im y(s)=\\=
\Big(\Real x(s)+i\cdot \Im x(s)\Big)\cdot\Big(\sigma(\Real y)+i\cdot\sigma(\Im y)\Big)
+
\Big(\sigma(\Real x)+i\cdot \sigma(\Im x)\Big)\cdot\Big(\Real y(s)+i\cdot\Im y(s)\Big)=\\=
x(s)\cdot\tau(y)+\tau(x)\cdot y(s).
\end{multline*}
\epr

\bit{
 \item[$\bullet$]
A {\it complex cotangent space} of an involutive stereotype algebra
$A$ at a point $s\in\Spec(A)$ is the stereotype quotient space of the ideal
$I_s[A]$ by the subideal $\overline{I_s^2}[A]$:
 \beq\label{DEF:T_s^(C-star)[A]}
{\C}T_s^\star[A]:=\Big(I_s[A]\Big/\overline{I_s^2}[A]\Big)^\triangledown.
 \eeq
The elements of this space are called {\it complex cotangent vectors} (of the algebra $A$ at the point $s\in\Spec(A)$).

 \item[$\bullet$]
Since both ideals $I_s[A]$ and $\overline{I_s^2}[A]$ are closed under the involution, the formula
$$
(a+\overline{I_s^2}[A])^\bullet=a^\bullet+\overline{I_s^2}[A],\qquad a\in I_s[A]
$$
defines an involution on the quotient space
${\C}T_s^\star[A]=\Big(I_s[A]/\overline{I_s^2}[A]\Big)^\triangledown$ (as a stereotype $A$-bimodule). Of course, this definition is chosen so that it guarantees that the quotient map $\pi:I_s\to (I_s/\overline{I_s^2})^\triangledown$ turns the involution on $I_s$ induced from $A$ into the involution on $(I_s/\overline{I_s^2})^\triangledown$:
$$
\pi(a^\bullet)=\pi(a)^\bullet.
$$

\item[$\bullet$] A {\it cotangent vector} (or {\it real cotangent vector}) of an involutive stereotype algebra $A$ at a point $s\in\Spec(A)$ is an arbitrary real vector (in the sense of definition \eqref{DEF:veshestv-vektor}) in the space ${\C}T_s^\star[A]$, i.e. an arbitrary complex cotangent vector $\xi\in {\C}T_s^\star[A]$, which is stable under the involution $\bullet$ in ${\C}T_s^\star[A]$:
 \beq\label{DEF:kokasat-vektor}
\xi^\bullet=\xi.
 \eeq
The set of all cotangent vectors of the algebra $A$ at a point $s\in\Spec(A)$ is called the {\it cotangent space} (of the algebra $A$ at the point $s\in\Spec(A)$) and is denoted by $T_s^\star[A]$. Certainly, $T_s^\star[A]$ is the real part of ${\C}T_s^\star[A]$:
$$
T_s^\star[A]=\Real{\C}T_s^\star[A].
$$
 }\eit

\btm\label{TH:svyaz-T(M)-i-T^*(M)} The formula
\beq\label{tau(x)=(f-circ-pi)(x-s(x)-cdot-1_A)}
\tau(a)=(f\circ\pi)(a-s(a)\cdot
1_A),\qquad a\in A,
\eeq
establishes
 \bit{ \item[---] a bijection between the complex tangent vectors $\tau\in T_s[A]$ to the algebra $A$ at the point $s$ and the $\C$-linear continuous functionals $f:{\C}T_s^\star[A]\to\C$ on the complex cotangent space ${\C}T_s^\star[A]$, and this bijection is an isomorphism of stereotype spaces
 \beq
{\C}T_s[A]\cong {\C}T_s^\star[A]^\star_{\C},
 \eeq

\item[---]  a bijection between the tangent vectors $\tau\in T_s[A]$ to the algebra $A$ at the point $s$ and the $\R$-linear continuous functionals $f:T_s^\star[A]\to\R$, and this bijection is an isomorphism of stereotype spaces
 \beq\label{T_s[A]-cong-Real(I_s/overline(I_s^2)^star)}
T_s[A]\cong T_s^\star[A]^\star_{\R}.
 \eeq
  }\eit
 \etm
\bpr For us the second part of this proposition is important, so we concentrate on it.

1. Let us first show that for each functional $f\in T_s^\star[A]^\star_{\R}$ the formula \eqref{tau(x)=(f-circ-pi)(x-s(x)-cdot-1_A)} defines a tangent vector $\tau\in T_s[A]$. This functional $\tau$ is obviously linear, continuous and real, and we only have to check
\eqref{DEF:kasatelnyj-vektor}. Indeed, for any $a,b\in A$ we have:
\begin{multline*}
\tau(a\cdot b)-s(a)\cdot\tau(b)-\tau(a)\cdot s(b)=(f\circ\pi)(a\cdot b-s(a\cdot b)\cdot 1_A)-s(a)\cdot(f\circ\pi)(b-s(b)\cdot 1_A)-(f\circ\pi)(a-s(a)\cdot 1_A)\cdot s(b)=\\
=(f\circ\pi)(a\cdot b-s(a)\cdot b-a\cdot s(b)+s(a)\cdot s(b)\cdot 1_A)=(f\circ\pi)\Big(\underbrace{(a-s(a)\cdot 1_A)}_{\scriptsize\begin{matrix}\text{\rotatebox{90}{$\owns$}}\\ I_s\end{matrix}}\cdot \underbrace{(b-s(b)\cdot 1_A)}_{\scriptsize\begin{matrix}\text{\rotatebox{90}{$\owns$}}\\ I_s\end{matrix}}\Big)=0
\end{multline*}
(since $\pi\big|_{I_s^2}=0$).

2. On the contrary, if $\tau\in T_s[A]$, then for any $a,b\in I_s$ we have
$$
 \tau(a\cdot b)=\underbrace{s(a)}_{\scriptsize\begin{matrix}\text{\rotatebox{90}{$=$}}\\ 0\end{matrix}}\cdot\tau(b)+\tau(a)\cdot \underbrace{s(b)}_{\scriptsize\begin{matrix}\text{\rotatebox{90}{$=$}}\\ 0\end{matrix}}=0,
$$
hence $\tau\big|_{I_s^2}=0$, and thus $\tau\big|_{I_s}$ factors through the quotient map $\pi:I_s\to (I_s/\overline{I_s^2})^\triangledown$:
$$
\tau\big|_{I_s}=f\circ\pi
$$
for each $f\in\big((I_s/\overline{I_s^2})^\triangledown\big)^\star_{\C}$. That is,
$$
\tau(x)=f(\pi(x)),\qquad x\in I_s.
$$
At that same time, since $\tau$ and $\pi$ preserve involution, and $\pi$ has dense image in the domain of $f$, the functional $f$ must also preserve involution, so
$$
f\in \Real\big((I_s/\overline{I_s^2})^\triangledown\big)^\star_{\C}=
\Real\big({\C}T_s^\star[A]\big)^\star_{\C}=\eqref{Real(X^star_C)=(Real-X)^star_R}=
\big(\Real{\C}T_s^\star[A]\big)^\star_{\R}=T_s^\star[A]^\star_{\R}.
$$
Finally, for each $a\in A$
$$
\tau(a)=\tau(a)-s(a)\cdot \underbrace{\tau(1_A)}_{\scriptsize\begin{matrix}\text{\rotatebox{90}{$=$}}\\ 0\end{matrix}}=\tau(\underbrace{a-s(a)\cdot 1_A}_{\scriptsize\begin{matrix}\text{\rotatebox{90}{$\owns$}}\\ I_s\end{matrix}})=f\Big(\pi\big(a-s(a)\cdot 1_A\big)\Big).
$$

3. Now we have to verify that the mapping $\tau\mapsto f$ is a homeomorphism. Let us start with the space ${\C}T_s[A]$ of complex tangent vectors. The topology of ${\C}T_s[A]$ is the pseudosaturation of the topology of uniform convergence on totally bounded sets in $A$, i.e. of the topology generated by seminorms
$$
\abs{\tau}_K=\sup_{a\in K}\abs{\tau(a)}
$$
where $K$ runs through the system of all totally bounded sets in $A$. For any such $K\subseteq A$ we can consider the set
$$
T_K=\{a-s(a)\cdot 1_A;\ a\in K\},
$$
which, in contrast to $K$, lies in the ideal $I_s$. But, similarly to $K$, the set $T_K$ is totally bounded in $A$, hence in $I_s$ as well (if we endow $I_s$ with the topology of an immediate stereotype space in $A$). We have:
$$
\abs{\tau}_{T_K}=\sup_{x\in T_K}\abs{\tau(x)}=\sup_{a\in K}\abs{\tau(a-s(a)\cdot 1_A)}=\sup_{a\in K}\abs{\tau(a)}=\abs{\tau}_K
$$
We can conclude that the topology of ${\C}T_s[A]$ coincides with the pseudosaturation of the topology of uniform convergence on totally bounded sets in $I_s$ (and not just in $A$). In addition, each functional $\tau\in
{\C}T_s[A]$ is uniquely  defined by its restriction on the ideal $I_s$. Hence we can think that
${\C}T_s[A]$ is an immediate subspace in the stereotype space $I_s^\star$, dual to $I_s$:
$$
{\C}T_s[A]\subseteq I_s^\star.
$$
Then ${\C}T_s[A]$ annihilates the subspace $I_s^2$, and thus the subspace $\overline{I_s^2}$, hence ${\C}T_s[A]$ is an immediate subspace in the annihilator $\Big(\overline{I_s^2}\Big)^\perp$ of the space
$\overline{I_s^2}$ with the topology of the immediate subspace in
$I_s^\star$:
$$
{\C}T_s[A]\subseteq \Big(\overline{I_s^2}\Big)^{\perp\vartriangle}.
$$
By \cite[(4.5)]{Akbarov}, the space in the right side is isomorphic to the dual space to the stereotype quotient space $\big(I_s/\overline{I_s^2}\big)^\triangledown$:
 \beq\label{T_s^C[A]subseteq(overline(I_s^2))^perp-vartriangle}
{\C}T_s[A]\subseteq \Big(\overline{I_s^2}\Big)^{\perp\vartriangle}\cong
\l\big(I_s/\overline{I_s^2}\big)^\triangledown\r^\star={\C}T_s^\star[A]^\star_{\C}.
 \eeq
We obtain an injection
${\C}T_s[A]\to\big(I_s/\overline{I_s^2}\big)^\triangledown$ which is exactly the mapping  $\tau\mapsto f$, defined by \eqref{tau(x)=(f-circ-pi)(x-s(x)-cdot-1_A)}. Since as we already understood, this formula defines a bijection, we see that the injection in the chain \eqref{T_s^C[A]subseteq(overline(I_s^2))^perp-vartriangle} must be an equality, moreover, an equality of stereotype spaces (since ${\C}T_s[A]$ is an immediate subspace in
$\Big(\overline{I_s^2}\Big)^{\perp\vartriangle}$):
$$
{\C}T_s[A]=\Big(\overline{I_s^2}\Big)^{\perp\vartriangle}\cong
\l\big(I_s/\overline{I_s^2}\big)^\triangledown\r^\star={\C}T_s^\star[A]^\star_{\C}.
$$
Now we pass to the real part, and we obtain
$$
T_s[A]=\Real{\C}T_s[A]=\Real\Big({\C}T_s^\star[A]^\star_{\C}\Big)=\eqref{Real(X^star_C)=(Real-X)^star_R}=
\Big(\Real{\C}T_s^\star[A]\Big)^\star_{\R}=\Big(T_s^\star[A]\Big)^\star_{\R}.
$$
  \epr

\bprop For each homomorphism of involutive stereotype algebras $\ph:A\to B$ and for any point $t\in\Spec(B)$ the formula
 \beq
\ph^\star_t(\tau)(a)=\tau(\ph(a)),\qquad \tau\in T_t(B)
 \eeq
defines a linear continuous map of tangent spaces
$\ph^\star_t:T_{t\circ\ph}[A]\gets T_t[B]$.
 \eprop

\brem
If $\ph:A\to B$ is a dense epimorphism, then $\ph^\star_t:T_{t\circ\ph}[A]\gets T_t[B]$ is a monomorphism:
 \begin{multline*}
(\tau\in T_t[B],\quad \tau\ne 0)\qquad\Longrightarrow\qquad \exists b\in B\quad \tau(b)\ne 0
\qquad\Longrightarrow\\ \Longrightarrow\qquad \exists a_i\in A\quad \ph(a_i)\overset{B}{\underset{i\to\infty}{\longrightarrow}} b\quad\&\quad \tau(b)\ne 0
\qquad\Longrightarrow\qquad \exists a_i\in A\quad \ph^T_t(a_i)=\tau(\ph(a_i))\ne 0.
 \end{multline*}
If $\ph$ is a monomorphism, then $\ph^T_t$ is not necessarily an epimorphism, as the following example shows.
\erem

\bex There is a monomorphism of involutive stereotype algebras $\ph:A\to B$ such that the mapping of tangent spaces  $\ph^\star_t:T_{t\circ\ph}[A]\gets T_t[B]$ is not an epimorphism for some $t\in\Spec(B)$.
 \eex
\bpr The natural embedding of the functional algebras
${\mathcal E}(\R)\subset{\mathcal C}(\R)$ has this property. For each point $t\in\R$
the corresponding mapping of the tangent spaces
$$
\ph^\star_t:T_t[{\mathcal E}(\R)]\cong\R\gets 0\cong T_t[{\mathcal
C}(\R)],
$$
cannot be an epimorphism (in this case, a surjection), of course.
 \epr

Suppose $A$ and $B$ are two involutive stereotype algebras, $s\in\Spec(A)$, $t\in\Spec(B)$, and $\sigma\in T_s[A]$. Then the formula
\beq\label{sigma-odot-t}
\sigma\odot t(a\circledast b)=\sigma(a)\cdot t(b),\qquad a\in A,\ b\in B,
\eeq
defines a tangent vector to the algebra $A\circledast B$ at the point $s\odot t\in\Spec(A\circledast B)$, sunce
\begin{multline*}
\sigma\odot t\big((a\circledast b)\cdot(a'\circledast b')\big)=
\sigma\odot t\big((a\cdot a')\circledast(b\cdot b')\big)=
\sigma(a\cdot a')\cdot t(b\cdot b')=
\big(\sigma(a)\cdot s(a')+s(a)\cdot\sigma(a')\big)\cdot t(b)\cdot t(b')=\\=
\sigma(a)\cdot s(a')\cdot t(b)\cdot t(b')+s(a)\cdot\sigma(a')\cdot t(b)\cdot t(b')=
\sigma\odot t(a\circledast b)\cdot s\odot t(a'\circledast b')+s\odot t(a\circledast b)\cdot \sigma\odot t(a'\circledast b')
\end{multline*}
Similarly, if $\tau\in T_t[B]$ is a tangent vector to $B$ at $t$, then the formula
\beq\label{s-odot-tau}
s\odot\tau(a\circledast b)=s(a)\cdot\tau(b),\qquad a\in A,\ b\in B,
\eeq
defines a tangent vector to the algebra $A\circledast B$ at the point $s\odot t\in\Spec(A\circledast B)$.

\bprop\label{PROP:T_s[A]-oplus-T_t[B]}
The mapping
\beq\label{sigma-oplus-tau=sigma-odot-t+s-odot-tau}
\sigma\oplus\tau\mapsto \sigma\odot t+s\odot\tau,\qquad \sigma\in T_s[A],\quad \tau\in T_t[B]
\eeq
is an isomorphism of stereotype spaces
\beq\label{T_s[A]-oplus-T_t[B]}
T_s[A]\oplus T_t[B]\cong T_{s\odot t}[A\circledast B]
\eeq
\eprop
\bpr
The injectivity of this mapping is obvious. Let us prove the surjectivity. For any tangent vector  $\upsilon\in T_{s\odot t}[A\circledast B]$ one can consider the functional
$$
\sigma(a)=\upsilon(a\circledast 1),\qquad a\in A,
$$
and this is a tangent vector to $A$ at $s$, since
\begin{multline*}
\sigma(a\cdot a')=\upsilon\big((a\cdot a')\circledast 1\big)=\upsilon\big((a\circledast 1)\cdot (a'\circledast 1)\big)=\upsilon(a\circledast 1)\cdot s\odot t(a'\circledast 1)+s\odot t(a\circledast 1)\cdot\upsilon(a'\circledast 1)=\\=\sigma(a)\cdot s(a)+s(a)\cdot\sigma(a').
\end{multline*}
Similarly, the functional
$$
\tau=\upsilon(1\circledast b),\qquad b\in B,
$$
is a tangent vector to $B$ at $t$. We have
\begin{multline*}
\upsilon(a\circledast b)=\upsilon\big((a\circledast 1)\cdot (1\circledast b)\big)=
\upsilon(a\circledast 1)\cdot s\odot t(1\circledast b)+s\odot t(a\circledast 1)\cdot\upsilon(1\circledast b)=\\=
\sigma(a)\cdot t(b)+s(a)\cdot\tau(b)=(\sigma\odot t+s\odot\tau)(a\circledast b).
\end{multline*}
The continuity of this map in both directions is obvious.
\epr

\paragraph{Functional algebras.}

The key examples that illustrate the notion of spectrum and (co)tangemt space are the standard functional algebras.

\begin{ex}\label{ex-10.3} {\sl Algebra $\mathcal{C}(M)$ of continuous functions on a (Hausdorff) paracompact locally compact space $M$}, endowed with (the poinwise multiplication and) the topology of uniform convergence on compact sets $S\subseteq M$, is a stereotype algebra \cite[Section 8.1]{Akbarov}. The space $M$ is embedded into the involutive spectrum of the algebra $\mathcal{C}(M)$ via delta-functionals
$$
s\in M\mapsto \delta^s\in \Spec \mathcal{C}(M),
$$
and this mapping is a homeomorphism of topological spaces:
\beq\label{Spec-C(M)=M}
\Spec \mathcal{C}(M)=M.
\eeq
The involutive tangent and cotangent spaces vanish in any point $s\in M$:
\beq\label{T_sC(M)=0}
T_s\mathcal{C}(M)=0,\qquad T_s^\star\mathcal{C}(M)=0.
\eeq
\end{ex}
\bpr
1. Let us prove \eqref{Spec-C(M)=M}. First, the delta-functionals ar involutive and multiplicative,
$$
\overline{\delta_s(\overline{f})}=\overline{\overline{f(a)}}=f(a)=\delta_a(f),\qquad f\in {\mathcal C}(M),
$$
therefore they map $M$ into $\Spec \mathcal{C}(M)$. Further, the map $\delta:M\to\Spec{\mathcal C}(M)$ is injective, since on a paracompsct space the continuous functions separate points. It is not quite clear, thst it is surjective. Let $\chi:{\mathcal C}(M)\to\C$ be an element of spectrum, i.e. an involutive continuous and {\it unital} homomorphism into $\C$. Let us show that there is a point $s\in M$ such that
\beq\label{PROOF:Spec-C(M)=M-2}
\Ker\chi\subseteq\Ker\delta_s.
\eeq
Suppose this is not true:
\beq\label{PROOF:Spec-C(M)=M-1}
\forall s\in M\quad \exists f_s\in\Ker \chi:\quad f_s(s)\ne 0.
\eeq
Then the sets $U_s=\{t\in M:\ f_s(t)\ne 0\}$ form a covering of $M$. Hence there is a subordinated partition of unity:
\beq\label{PROOF:Spec-C(M)=M}
\eta_s\in{\mathcal C}(M)\quad \supp\eta_s\subseteq U_s\quad \sum_{s\in M}\eta_s=1.
\eeq
Consider the function
$$
g=\sum_{s\in M}f_s\cdot\overline{f_s}\cdot\eta_s.
$$
The series in the right is locally finite, so it converges in the space ${\mathcal C}(M)$. And its elements belong to the ideal $\Ker\chi$, since $f_s\in\Ker\chi$. Therefore,
$$
g\in \Ker\chi.
$$
Let us show that $g>0$ everywhere. Take $t\in M$. From \eqref{PROOF:Spec-C(M)=M} we have that there is a point $s_t\in M$ such that
$$
\eta_{s_t}(t)>0.
$$
Since $\supp\eta_s\subseteq U_s=\{t\in M:\ f_s(t)\ne 0\}$, we have
$$
f_{s_t}(t)\ne 0,
$$
hence,
$$
g(t)=\sum_{s\in M}f_s(t)\cdot\overline{f_s(t)}\cdot\eta_s(t)= \sum_{s\in M}\abs{f_s(t)}^2\cdot\eta_s(t)\ge
\underbrace{\abs{f_{s_t}(t)}^2}_{\scriptsize\begin{matrix}\text{\rotatebox{90}{$<$}}\\ 0\end{matrix}}\cdot\underbrace{\eta_{s_t}(t)}_{\scriptsize\begin{matrix}\text{\rotatebox{90}{$<$}}\\ 0\end{matrix}}
>0.
$$
This means that $g$ is an invertible element in the algebra ${\mathcal C}(M)$. At the same time it belongs to the ideal $\Ker\chi$. Thus, $\Ker\chi={\mathcal C}(M)$. And this contradicts the condition of unitality: $\chi(1)=1$.

So our supposition \eqref{PROOF:Spec-C(M)=M-1} is not true, and, as a corollary, there exists a point $s\in M$ such that \eqref{PROOF:Spec-C(M)=M-2} holds. Since $\Ker\chi$ is the kernel of a linear functional, this is a maximal ideal. This implies in its turn that the embedding \eqref{PROOF:Spec-C(M)=M-2} can be replaced by the equality:
$$
\Ker\chi=\Ker\delta_s.
$$
We see that the functionals $\chi$ and $\delta_s$ have the same kernels and they coincide at the element 1:
$$
\chi(1)=1=\delta_s(1).
$$
Thus, they coincide: $\chi=\delta_s$. This proves the surjectivity of the mapping $\delta:M\to\Spec{\mathcal C}(M)$. Its injectivity was noticed, and the continuity in both directions follows from the Ascoli theorem \cite[Theorem 8.2.10]{Engelking}.

2. Let us prove \eqref{T_sC(M)=0}. Take $s\in M$, $\tau\in T_s[{\mathcal C}(M)]$, and show that for each function  $f\in {\mathcal C}(M)$ we have $\tau(f)=0$. It is convenient to consider several cases.

a) Suppose first that $f(t)\ge 0$, $t\in M$, and $f(s)=0$. Then the function $g(t)=\sqrt{f(t)}$ also belongs to ${\mathcal C}(M)$, and we have
$$
\tau(f)=\tau(g^2)=2 g(s)\cdot\tau(g)=2\sqrt{f(s)}\cdot\tau(g)=0.
$$

b) Suppose then that $f$ is a real function (not necessarily non-negative), such that $f(s)=0$. Then it can be represented as a difference of two non-negative functions with the same property
$$
f=g-h,\quad g\ge 0,\quad h\ge 0,\quad g(s)=h(s)=0.
$$
By what we have already proved,
$$
\tau(f)=\tau(g)-\tau(h)=0.
$$

c) Suppose $f$ is an arbitrary (not necessarily real) function, such that $f(s)=0$.
Then by the previous,
$$
\tau(f)=\tau(\Real f+i\cdot \Im f)=\tau(\Real f)+i\cdot\tau(\Im f)=0.
$$

d) Suppose $f$ is an arbitrary function (not necessarily vanishing at $s$). Then
$$
\tau(f)=\tau(f-f(s)\cdot 1+f(s)\cdot 1)=\tau(f-f(s)\cdot 1)+\tau(f(s)\cdot 1)=
\kern-5pt \underbrace{\tau(f-f(s)\cdot 1)}_{\scriptsize\begin{matrix}\|\\ \phantom{,} 0,\\ \text{since}\\ (f-f(s)\cdot 1)(s)=0\end{matrix}}\kern-5pt +f(s)\cdot\kern-8pt\underbrace{\tau(1)}_{\scriptsize\begin{matrix}\|\\ \phantom{,}0,\\ \text{by \eqref{tau(1)=0}}\end{matrix}}\kern-8pt =0.
$$

We proved the first equality in \eqref{T_sC(M)=0}. The second one is its corollary.
\epr

\begin{ex}\label{ex-10.4} {\sl Algebra $\mathcal{E}(M)$ of smooth fnctions on a smooth manifold $M$} (with the usual pointwise multiplication and the topology of uniform convergence by each partial derivative on compact sets), being a Fr\'{e}chet algebra, is a stereotype algebra. Like in the previous case, the involutive spectrun of this algebra coincides with $M$ (as a topological space),
$$
\Spec \mathcal{E}(M)=M,
$$
but the tangent and the cotangent spaces in an arbitrary point $s\in M$ do not degenerate, and they coincide with the tangent and the cotangent spaces to this manifold $M$:
$$
T_s\mathcal{E}(M)=T_s(M),\qquad T_s^\star\mathcal{E}(M)=T_s^*(M).
$$
\end{ex}

\begin{ex}\label{ex-10.5} {\sl Algebra $\mathcal{O}(M)$ of holomorphic functions on a Stein manifold $M$} (with the pointwise multiplication and the topology of uniform convergence on compact sets in $M$), being a Fr\'{e}chet algebra, is a stereotype algebra. H.~Rossi in \cite{Rossi} proves the homeomorphism
$$
\Spec \mathcal{O}(M)=M.
$$
The equalities
$$
T_s\mathcal{O}(M)=T_s(M),\qquad T_s^\star\mathcal{O}(M)=T_s^*(M)
$$
(obvuois for open subsets $M$ in $\C^n$), apparently are not proved in the general case.
\end{ex}

\begin{ex}\label{EX:kasat-prostr-k-P(M)} {\sl Algebra of polynomials (regular functions) $\mathcal{P}(M)$
on an affine algebraic manifold $M$}. Let us recall \cite{Humphreys}, \cite{Vinberg-Onischik}, that an {\it affine algebraic variety} over a field $K$ is the common set $M$ of zeroes of an arbitrary given (not necessarily finite) set of polynomials $F$ on $K^n$:
$$
M=\{x\in\R^n:\ \forall u\in F\ u(x)=0\}
$$
If we denote by the symbol ${\mathcal P}(K^n)$ the algebra of polynomials on $K^n$, and by $I$ the ideal in it, generated by the set $F$, then we can perceive the quotient algebra ${\mathcal P}(K^n)/I$ as an algebra of functions on $M$. It is called the {\it algebra of polynomials on $M$}, and we denote it by ${\mathcal P}(M)$. For $K=\R$ the variety $M$ is said to be {\it real}, and for $K=\C$ {\it complex}. A point $a\in M$ of a real affine algebraic manifold $M\subseteq\R^n$ is said to be {\it simple} \cite[3.3.10]{Bochnak-Coste-Roy}, if there is a finite set of polynomials $u_1,...,u_k\in{\mathcal P}(\R^n)$ with the following properties:
 \bit{
 \item[(i)] the functions $u_1,...,u_k$ are independent in the point $a$, i.e. their differentials in this point $\d u_1(a),...,\d u_k(a)$ are linearly independent (as functionals on $\R^n$),

 \item[(ii)] the functions $u_1,...,u_k$ vanish on $M$:
 $$
 u_i\Big|_M=0,\qquad i=1,...,k,
 $$

 \item[(iii)] there is a neighbourhood $V$ (in the Zarisski topology) of the point $a$ in $\R^n$, where the set of common zeroes of the functions $u_1,...,u_k$ coincide with the intersection $V\cap M$:
     $$
     \{x\in V:\quad \forall i=1,...,k\quad u_i(x)=0\}=V\cap M.
     $$

 }\eit

For the case $K=\R$ we need the following two propositions:
 \bit{
 \item[(a)] the real spectrum of ${\mathcal P}(M)$, and the involutive spectrum of its complexification $\C{\mathcal P}(M)$ coincide with $M$:
 \beq\label{Spec[CP(M)]=RSpec[P(M)]=M}
\Spec[\C{\mathcal P}(M)]=\R\Spec[{\mathcal P}(M)]=M
 \eeq

 \item[(b)] in each simple point $a\in M$ every real tangent vector $\tau\in T_a[{\mathcal P}(M)]$ to the algebra  ${\mathcal P}(M)$, and every involutive tangent vector to its complexification $\C{\mathcal P}(M)$, are derivatives along some smooth curve $\gamma:\R\to M$, going from $a$, $\gamma(0)=a$:
\beq\label{kasat-prostr-k-P(G)}
\tau(u)=\lim_{t\to 0}\frac{u(\gamma(t))-u(a)}{t},\qquad u\in {\mathcal P}(M),
\eeq
thus,
 \beq\label{T_a[CP(M)]=RT_a[P(M)]=T_a(M)}
T_a[\C{\mathcal P}(M)]=\R T_a[{\mathcal P}(M)]=T_a(M)
 \eeq
}\eit
\eex
\bpr
1. Consider first the case when $M=\C^n$.
If $s:{\mathcal P}(\R^n)\to\R$ is a (unital) homomorphism of algebras, then it is a linear functional on the space $(\R^n)_{\R}^\star$ of linear functionals on $\R^n$. I.e. $s$ is an element of the second dual space  $((\R^n)_{\R}^\star)_{\R}^\star\cong\R^n$. Hence, there is a point $x\in\R^n$ such that
$$
s(u)=u(x),
$$
for any linear functional $u:\R^n\to\R$. This identity is extended by multiplicativity to the polynomials  $u\in{\mathcal P}(\R^n)$, and we get the second identity in \eqref{Spec[CP(M)]=RSpec[P(M)]=M}. The first one follows from \eqref{DEF:Spec[A]}.

Further, take $a\in\R^n$ and $\tau\in\R T_a[{\mathcal P}(\R^n)]$. Note that each polynomial $u\in{\mathcal P}(M)={\mathcal P}(\C^n)$ has an expansion
$$
u(z)=u(a)+\sum_{k=1}^n\partial_k u(a)\cdot (z_k-a_k)+\sum_{k=1}^n\partial_k v_k(z)\cdot w_k(z),
$$
where $v_k$ and $w_k$ are polynomials vanishing in $a$. As a corollary, the action of $\tau$ on $u$ is uniquely defined by its action on monomials $z_k-a_k$:
$$
\tau(u)=0+\sum_{k=1}^n\partial_k u(a)\cdot \tau(z_k-a_k)+0.
$$
If we take the vector $b\in\C^n$ with coordinates $\tau(z_k-a_k)$, then the differentiation along the curve
$$
\gamma(t)=a+t\cdot b,
$$
coincides with the action of the functional $\tau$, i.e. \eqref{kasat-prostr-k-P(G)} holds.

2. Now the general case. By definition, the algebra ${\mathcal P}(M)$ is a quotient algebra of the algebra of polynomials ${\mathcal P}(\C^n)$ on some $\C^n$:
$$
{\mathcal P}(M)\cong {\mathcal P}(\C^n)/I
$$
where $I$ is an ideal in ${\mathcal P}(\C^n)$. Consider the quotient mapping $\pi:{\mathcal P}(\C^n)\to{\mathcal P}(\C^n)/I\cong {\mathcal P}(M)$. The composition
$$
\tau\circ\pi:{\mathcal P}(\C^n)\to\C
$$
is a tangent vector to the algebra ${\mathcal P}(\C^n)$ in the point $a\in M\subseteq\C^n$, and we already proved that it coincides with the differentiation along a smooth curve $\gamma:\R\to \C^n$, going from $a$, $\gamma(0)=a$:
$$
(\tau\circ\pi)(v)=\lim_{t\to 0}\frac{v(\gamma(t))-v(a)}{t},\qquad v\in {\mathcal P}(\C^n).
$$
On each polynomial $v\in I$ this vector vanishes, since $\pi(v)=0$,
$$
0=\tau(\pi(v))=(\tau\circ\pi)(v)=\lim_{t\to 0}\frac{v(\gamma(t))-v(a)}{t},\qquad v\in I.
$$
This implies, first, \eqref{kasat-prostr-k-P(G)}, and, second, that if $a$ is a simple point of $M$, then the ideal  $I$ has a finite set of functions $v_1,...,v_l$ such that $M$ is a non-degenerate level set in a neighbourhood of $a$, and since the vector $\gamma'(0)$ vanishes on $v_1,...,v_l$, it must be a tangent vector to the manifold $M$. Hence, one can change $\gamma$ in such a way that it will belong to $M$. This proves \eqref{kasat-prostr-k-P(G)}. We obtain the second equality in \eqref{T_a[CP(M)]=RT_a[P(M)]=T_a(M)}, and the first one follows from \eqref{DEF:T_s[A]}.
\epr

The following identities hold \cite[Theorems 8.4, 8.10, 8.13, 8.16]{Akbarov}:
\begin{align}
& {\mathcal C}(M)\odot{\mathcal C}(N)\cong{\mathcal C}(M\times N), \label{C(M-times-N)} \\
& {\mathcal E}(M)\circledast{\mathcal E}(N)\cong{\mathcal E}(M)\odot{\mathcal E}(N)\cong{\mathcal E}(M\times N),
\label{E(M-times-N)} \\
& {\mathcal O}(M)\circledast{\mathcal O}(N)\cong{\mathcal O}(M)\odot{\mathcal O}(N)\cong{\mathcal O}(M\times N),
\label{O(M-times-N)} \\
& {\mathcal P}(M)\circledast{\mathcal P}(N)\cong{\mathcal P}(M)\odot{\mathcal P}(N)\cong{\mathcal P}(M\times N),
\label{P(M-times-N)}
\end{align}
and they imply

\btm\label{TH:C(M),E(M),O(M),P(M)-inj-alg}
Functional algebras ${\mathcal C}(M)$, ${\mathcal E}(M)$, ${\mathcal O}(M)$, ${\mathcal P}(M)$ are injective stereotype algebras\footnote{In the sense of definition on page \pageref{DEF:proj-algebra-odot}.}.
\etm

The following result will be useful in Lemma \ref{LM-3-0-dlya-TH:C^infty-obolochka-podalgebry-v-C^infty}.

\blm\label{LM:o-plotnom-ideale}
Let $A$ be a dense unital subalgebra in the functional algebra $B={\mathcal C}(M)$ or in $B={\mathcal E}(M)$ (in particular, $A$ contains constants). Then for any point $t\in M$ the ideal
$$
I_t^A=\{a\in A:\quad a(t)=0\}
$$
is dense in the ideal
$$
I_t^B=\{u\in B:\quad u(t)=0\}.
$$
\elm
\bpr
Take $u\in I_t^B$, i.e. $u(t)=0$. Since $A$ is dense in ${\mathcal C}(K)$, there is a net $a_i\in A$ such that
$$
a_i\overset{B}{\underset{i\to\infty}{\longrightarrow}}u.
$$
Put $b_i=a_i-a_i(t)$. Then $b_i(t)=0$, hence $b_i\in I_t^A$. On the other hand,
$$
b_i=a_i-a_i(t)\overset{B}{\underset{i\to\infty}{\longrightarrow}}u-\underbrace{u(t)}_{\scriptsize\begin{matrix}\|\\ 0\end{matrix}}=u.
$$
\epr

\subsection{Group algebras ${\mathcal C}^\star(G)$ and ${\mathcal E}^\star(G)$}

If a maniufold $M$ is endowed with a group structure, so that it becomes a group $G$ in the given geometry (i.e. a locally compact group, or a Lie group, real or complex, or an affine algebraic group), then the spaces ${\mathcal C}^\star(G)$, ${\mathcal E}^\star(G)$, ${\mathcal O}^\star(G)$, ${\mathcal R}^\star(G)$, dual to the functional algebras ${\mathcal C}(M)$, ${\mathcal E}(M)$, ${\mathcal O}(M)$, ${\mathcal R}(M)$, become stereotype algebras with trespect to the convolution. We are interested in two of these algebras, ${\mathcal C}^\star(G)$ and ${\mathcal E}^\star(G)$, and here we notice some general facts about them.

\paragraph{Convolution and involution on ${\mathcal C}^\star(G)$.}\label{paragraph:C^star(G)}
Let $G$ be a locally compact group, and ${\mathcal C}(G)$ the algebra of continuous functions on $G$. The dual stereotype space ${\mathcal C}^\star(G)$ consists of measures with compact support in $G$. It is an algebra with respect to the convolution $*$. We can define this operation as follows.

First, for functions $u\in \mathcal{C}(G)$ and measures $\alpha\in
\mathcal{C}^\star (G)$ we denote by $\widetilde{u}$ and $\widetilde{\alpha}$
their {\it anitpodes}
\beq
\widetilde{u}(t)=u(t^{-1}), \qquad \widetilde{\alpha}(u)=\alpha(\widetilde{u})
\label{eq10.4}
\eeq
and by $a\cdot u, u\cdot a$ and $a\cdot \alpha, \alpha \cdot a$ their {\it shifts}:
\beq\label{eq10.5}
(a\cdot u)(t)=u(t\cdot a), \qquad (u\cdot a)(t)=u(a\cdot t)
\eeq
\beq
(a\cdot \alpha)(u)=\alpha(u\cdot a), \qquad (\alpha\cdot a)(u)=\alpha(a\cdot u)
\label{eq10.6}
\eeq
Obviously,
\beq
\widetilde{a\cdot u}=\widetilde{u}\cdot a^{-1}, \qquad \widetilde{u\cdot
a}=a^{-1}\cdot \widetilde{u} \label{eq10.7}
\eeq
\beq
\widetilde{a\cdot \alpha}=\widetilde{\alpha }\cdot a^{-1}, \qquad
\widetilde{\alpha \cdot a}=a^{-1}\cdot \widetilde{\alpha} \label{eq10.8}
\eeq
If we denote by $\delta^a$ the delta-functional
\beq\label{delta^a(u)=u(a)}
\delta^a(u)=u(a),
\eeq
then
\beq
\widetilde{\delta^a}=\delta^{a^{-1}}, \qquad a\cdot \delta^b=\delta^{a\cdot b},
\quad \delta^b\cdot a=\delta^{b\cdot a} \label{eq10.9}
\eeq
The convolution of a function with a measure is defined by the formula
\beq
\alpha * u (t)=\alpha\l \widetilde{t\cdot u}\r, \qquad u * \alpha (t)=\alpha\l
\widetilde{u\cdot t}\r \label{eq10.10}
\eeq
and the following identities are proved subsequently:
\beq
  \delta^a * u=u\cdot a^{-1}, \qquad u * \delta^a=a^{-1}\cdot u
\label{eq10.11}
\eeq
\beq
\widetilde{\alpha * u}=\widetilde{u} * \widetilde{\alpha},
 \qquad
\widetilde{u * \alpha}=\widetilde{\alpha} * \widetilde{u} \label{eq10.12}
\eeq
\beq
\alpha (\beta * u)=\beta(\alpha * \widetilde{u}),
 \qquad
\alpha (u * \beta)=\beta(\widetilde{u} * \alpha) \label{eq10.13}
\eeq
They imply the following chain of identities,
\beq\label{eq10.14}
\alpha \l u * \widetilde{\beta} \r= \alpha \l \widetilde{\beta * \widetilde{u}}
\r= \widetilde{\alpha} \l \beta * \widetilde{u} \r= \beta \l \widetilde{\alpha}
* u \r
\eeq
and this chain is called a {\it convolution of measures}:
\beq\label{DEF:svertka}
\alpha * \beta (u)=\alpha \l u * \widetilde{\beta} \r= \alpha \l
\widetilde{\beta * \widetilde{u}} \r= \widetilde{\alpha} \l \beta *
\widetilde{u} \r= \beta \l \widetilde{\alpha}
* u \r
\eeq
It is easily verified that $(\alpha, \beta) \mapsto \alpha * \beta$ is a continuous multiplication on  $\mathcal{C}(G)$ in the sense of definition at page \pageref{DEF:ster-alg}. Besides this, the following identities hold:
\beq\label{DEF:svertka-s-delta^a}
\delta^a * \beta=a\cdot \beta, \qquad \beta * \delta^a=\beta \cdot a,
\eeq
\beq\label{delta^a*delta^b=delta^(a-cdot-b)}
\delta^a * \delta^b=\delta^{a\cdot b}
\eeq
\beq\label{eq10.18}
\widetilde{\alpha * \beta}=\widetilde{\beta} * \widetilde{\alpha}
\eeq
\beq\label{sevrtka-kak-dvoinoi-integral}
\alpha * \beta (u)=(\alpha\otimes \beta)(w)\Big|_{w(s,t)=u(s\cdot t)}=
\int_G \l \int_G u(s\cdot t) \d \alpha(s) \r \d \beta(t)= \int_G \l \int_G
u(s\cdot t) \d \beta(t)\r \d \alpha(s)
\eeq
The first three of them are verified by direct computation, and the last one is proved for delta-functionals (this is sufficient, since they are total in ${\mathcal C}(G)$):
$$
\delta^a*\delta^b(u)=\eqref{delta^a*delta^b=delta^(a-cdot-b)}=\delta^{a\cdot b}(u)=u(a\cdot b)=\int_G u(a\cdot t)\cdot\delta^b(\d t)=\int_G \left(\int_G u(s\cdot t)\cdot\delta^a(\d s)\right)\delta^b(\d t).
$$

The involution of a function $u\in\mathcal{C}(G)$ is pointwise:
 \beq
\overline{u}(t)=\overline{u(t)},\qquad t\in G.
 \eeq
while for a measure $\alpha\in\mathcal{C}^\star (G)$ it is defined by the formula
 \beq
\alpha^\bullet(u)=\overline{\alpha(\overline{\widetilde{u}})}.
 \eeq
In particular,
$$
(\delta^a)^\bullet=\delta^{a^{-1}}.
$$
This meets the definition of involution on $L_1(G)$, since for the left invariant Haar measure $\mu$ on $G$ we have
$$
\mu(\overline{u})=\overline{\mu(u)},\qquad \mu(\widetilde{u})=\mu\left(\frac{u}{\Delta}\right)
$$
for finite functions $u$ on $G$, \cite[15.5 and 15.15]{Hewitt-Ross}, here $\Delta$ is the modular function:
$$
\mu(x^{-1}\cdot u)=\Delta(x)\cdot\mu(u).
$$
If for a measure $\mu\ge 0$ and a function $f\in L_1(\mu)$ we denote by $f\cdot\mu$ the measure
\beq\label{DEF:f-cdot-alpha(u)}
f\cdot\mu(u)=\int u(s)\cdot f(s)\cdot\mu(\d s),
\eeq
then for any function $f\in L_1(\mu)$ we have
$$
(f\cdot\mu)^\bullet(u)=\overline{(f\cdot\mu)(\overline{\widetilde{u}})}=\overline{\mu(f\cdot\overline{\widetilde{u}})}=
\overline{\mu\left(\overline{\widetilde{\overline{\widetilde{f}}\cdot u}}\right)}=
\mu\left(\widetilde{\overline{\widetilde{f}}\cdot u}\right)=\mu\left(\frac{\overline{\widetilde{f}}\cdot u}{\Delta}\right)=\left(\frac{\overline{\widetilde{f}}}{\Delta}\cdot\mu\right)(u)
$$
i.e. the usual formula for involution on $L_1(G)$ holds:
$$
f^\bullet=\frac{\overline{\widetilde{f}}}{\Delta}.
$$

\paragraph{The case of compact group.}

Suppose $K$ is a compact group, and $\mu_K$ is its normed Haar measure:
$$
\mu_K(K)=1.
$$
Then one can also define convolution on the space ${\mathcal C}(K)$:
\beq\label{f*g}
f*g(s)=\int_G f(t)\cdot g(t^{-1}\cdot s)\cdot \mu_K(\d t),\qquad f,g\in{\mathcal C}(K).
\eeq
The following identities hold\footnote{Here $f\cdot\mu_K$ is defined in \eqref{DEF:f-cdot-alpha(u)}.}
\beq\label{widetilde(f*g)}
\widetilde{f*g}=\widetilde{g}*\widetilde{f},\qquad f,g\in{\mathcal C}(K).
\eeq
\beq\label{(f-mu_K)*(g-mu_K)=(f*g)-mu_K}
(f\cdot\mu_K)*(g\cdot\mu_K)=(f*g)\cdot\mu_K,\qquad f,g\in {\mathcal C}(K).
\eeq
\bpr
First,
\begin{multline*}
\widetilde{f*g}(s)=\eqref{widetilde(u)(t)=u(t^(-1))}=f*g(s^{-1})=\eqref{widetilde(f*g)}=
\int_K f(t)\cdot g(t^{-1}\cdot s^{-1})\cdot \mu_K(\d t)=\int_K f(t)\cdot g((s\cdot t)^{-1})\cdot \mu_K(\d t)=\\=
\int_K \widetilde{f}(t^{-1})\cdot \widetilde{g}(s\cdot t)\cdot \mu_K(\d t)=
\begin{vmatrix}s\cdot t=r\\ t^{-1}=r^{-1}\cdot s\\ \d t=\d r\end{vmatrix}=
\int_K \widetilde{f}(r^{-1}\cdot s)\cdot \widetilde{g}(r)\cdot \mu_K(\d r)=\eqref{widetilde(f*g)}=
(\widetilde{g}*\widetilde{f})(s).
\end{multline*}
And, second,
\begin{multline*}
(f\cdot\mu_K)*(g\cdot\mu_K)(u)=\eqref{sevrtka-kak-dvoinoi-integral}=\int_K\left(\int_Ku(s\cdot t)\cdot g(t)\cdot\mu_K(\d t)\right)\cdot f(s)\cdot\mu_K(\d s)=\begin{vmatrix}s\cdot t=r\\ t=s^{-1}\cdot r\\ \d t=\d r\end{vmatrix}=\\=
\int_K\left(\int_Ku(r)\cdot g(s^{-1}\cdot r)\cdot\mu_K(\d r)\right)\cdot f(s)\cdot\mu_K(\d s)=
\int_K u(r)\cdot\left(\int_K f(s)\cdot g(s^{-1}\cdot r)\cdot\mu_K(\d s)\right)\cdot \mu_K(\d r)=\\=
\int_K u(r)\cdot(f*g)(r)\cdot \mu_K(\d r)=(f*g)\cdot\mu_K(u).
\end{multline*}
\epr

We need the following proposition that belongs to folklore:

\btm\label{TH:mu(ph(X))=mu(X)}
Each (continuous) automorphism $\ph:K\to K$ of a compact group $K$ preserves the Haar measure on $K$:
\beq\label{mu(ph(X))=mu(X)}
\mu(\ph(X))=\mu(X),\qquad X\subseteq K.
\eeq
\etm
\bpr
Since $\ph$ is a homeomorphism, it preserves the Borel sets, so we can consider the measure
$$
m(X)=\mu(\ph(X))
$$
on Borel sets $X\subseteq K$. For each point $t\in K$ we have
$$
m(t\cdot X)=\mu(\ph(t\cdot X))=\mu(\ph(t)\cdot\ph(X))=\mu(\ph(X))=m(X).
$$
I.e., $m$ is invariant on $K$, like $\mu$. Hence, these measures differ in a constant multiplier:
$$
m(X)=C\cdot\mu(X).
$$
On the other hand, $\ph(K)=K$, therefore
$$
m(K)=\mu(\ph(K))=\mu(K),
$$
and $C=1$.
\epr

\paragraph{Representations of locally compact groups.}
Let $G$ be a locally compact group and $A$ a stereotype algebra. A mapping $\pi:G\to A$ is called a {\it representation of $G$ in $A$}, if it is continuous and multiplicative:
$$
\pi(x\cdot y)=\pi(x)\cdot\pi(y),\qquad \pi(1)=1.
$$
In this case (see \cite[Theorem 10.12]{Akbarov}) there is a unique homomorphism of stereotype algebras  $\dot{\pi}:{\mathcal C}^\star(G)\to A$ such that the following diagram is commutative:
\begin{gather}\label{C*(G)=grup-alg-A}
\begin{diagram}
\node{G} \arrow[2]{e,t}{\delta} \arrow{se,b}{\pi}
\node[2]{\mathcal{C}^\star(G)} \arrow{sw,b}{\dot{\pi}}
\\
\node[2]{A}
\end{diagram}
\end{gather}
(here $\delta$ is the embedding as delta-functions).

\btm\label{TH:nepr-predstavleniya} Suppose a locally compact group $G$ acts by unitary operators in a Hilbert space  $X$, i.e. a mapping $\pi:G\to{\mathcal L}(X)$ is defined such that
\beq\label{DEF:representation-of-G-in-X }
\pi(1_G)=1_X,\qquad \pi(s\cdot t)=\pi(s)\cdot\pi(t),\qquad \pi(s^{-1})=\pi(s)^\bullet,\qquad s,t\in G.
\eeq
Then the following conditions are equivalent:
 \bit{
\item[(i)] the mapping $\pi:G\to{\mathcal L}(X)$ is continuous with respect to the strong operator topology in ${\mathcal L}(X)$,

\item[(ii)] the mapping $(s,x)\in G\times X\mapsto \pi(s)x\in X$ is continuous,

\item[(iii)] the mapping $\pi:G\to{\mathcal L}(X)$ is continuous with respect to the stereotype topology of the space ${\mathcal L}(X)$.

\item[(iv)] there is a (necessarily unique) involutive continuous homomorphism of stereotype algebras $\dot{\pi}:\mathcal{C}^\star(G)\to {\mathcal L}(X)$ such that the following diagram is commutative:
\begin{gather}\label{C*(G)=grup-alg}
\begin{diagram}
\node{G} \arrow[2]{e,t}{\delta} \arrow{se,b}{\pi}
\node[2]{\mathcal{C}^\star(G)} \arrow{sw,b}{\dot{\pi}}
\\
\node[2]{{\mathcal L}(X)}
\end{diagram},
\end{gather}
 }\eit
\etm
\bit{
\item[$\bullet$] Further by {\it unitary representation of a locally compact group $G$} (in a Hilbert space $X$) we mean a mapping $\pi:G\to{\mathcal L}(X)$ that satisfies \eqref{DEF:representation-of-G-in-X } and the conditions (i)-(iv) of this Theorem.
}\eit

\bpr
The equivalence (i)$\Leftrightarrow$(ii) is a standard fact (see e.g. \cite[Chapter 5, \S 1]{Barut-Raczka}), and
(iii)$\Leftrightarrow$(iv) is proved in \cite[Theorem 10.12]{Akbarov}. Let us prove (ii)$\Rightarrow$(iii). Suppose the mapping $(s,x)\in G\times X\mapsto \pi(s)x\in X$ is continuous. Then for any compact set $K\subseteq X$
$$
\pi(g)x\underset{g\to g_0}{\longrightarrow}\pi(g_0)x
$$
uniformly by $x\in K$. This means that the mapping $\pi:G\to X:X$ is continuous (here $X:X$ means the space of linear continuous mappings from $X$ into $X$, endowed with the topology of uniform convergence on compact sets). For any compact set $T\subseteq G$ its image $\pi(T)$ is compact in $X:X$, so the topology on $\pi(T)$ deos not change under the pseudosaturation of the space $X:X$ (i.e. ubder the passage from $X:X$ to ${\mathcal L}(X)$). This means that the restriction $\pi|_T:T\to {\mathcal L}(X)$ is continuous. This is true for any compact set $T\subseteq G$, hence the map $\pi:G\to {\mathcal L}(X)$ is continuous as well. We proved that (ii)$\Rightarrow$(iii). The reverse implication  $(ii)\Leftarrow(iii)$ is proved by the reverse reasoning.
\epr

\paragraph{The algebra  ${\mathcal E}^\star(G)$.}
If $G$ is a real Lie  group, then one can consider also the group algebra ${\mathcal E}^\star(G)$ of distributions with compact support on $G$. It is defined by analogy: first we consider the algebra  ${\mathcal E}(G)$ of smooth functions on $G$ with the usual topology of uniform convergence on compact sets of each partial derivative in each local chart. Then ${\mathcal E}^\star(G)$ is defined as the stereotype dual of ${\mathcal E}(G)$. The same formulas define multiplication and involution on ${\mathcal E}^\star(G)$, and similarly there is a connection in an analog of diagram \eqref{C*(G)=grup-alg-A} between the representations of $G$ and of ${\mathcal E}^\star(G)$:
\begin{gather}\label{E*(G)=grup-alg-A}
\begin{diagram}
\node{G} \arrow[2]{e,t}{\delta} \arrow{se,b}{\pi}
\node[2]{\mathcal{E}^\star(G)} \arrow{sw,b}{\ddot{\pi}}
\\
\node[2]{A}
\end{diagram}
\end{gather}
but the difference is that in this case $\pi$ must be a smooth (in the usual sense) mapping \cite[Theorem 10.12]{Akbarov}.

\paragraph{Decomposition by the characters of normal compact subgroup.}

Let us remind some notions in the representation theory \cite{Kowalski}. Suppose $\{X_i;\ i\in I\}$ is a family of Hilbert spaces. Consider the set $X=\dot{\bigoplus}_{i\in I}X_i$ of families $x=\{x_i;\ x_i\in X_i\}$ such that
$$
x=\{x_i;\ i\in I\}\in \dot{\bigoplus_{i\in I}}X_i\qquad\Longleftrightarrow\qquad  \forall i\in I\quad x_i\in X_i\quad
\&
\quad \sum_{i\in I} \norm{x_i}_{X_i}^2<\infty.
$$
This is a Hilbert space with the coordinate-wise algebraic operations and the scalar product
$$
\langle x,y \rangle=\sum_{i\in I}\langle x_i,y_i\rangle_i
$$
(here $\langle x_i,y_i\rangle_i$ is the scalar product in $X_i$). This space is called a {\it Hilbert direct sum of the spaces $X_i$}.

If $G$ is a locally compact group and $\{\pi_i:G\to {\mathcal L}(X_i)\}$ a family of unitary representations in the spaces $X_i$, then the mapping
$$
\pi:G\to {\mathcal L}\left(\dot{\bigoplus_{i\in I}}X_i\right)\qquad\Big|\qquad \pi(g)\{x_i;\ x_i\in X_i\}=\{\pi_i(g)x_i;\ x_i\in X_i\}
$$
is a unitary representation of $G$ and is called {\it Hilbert direct sum of the representations $\pi_i$}.

A linear continuous mapping of Hilbert spaces $\alpha:X\to Y$ is called
 \bit{
\item[---] a {\it morphism of a representation $\pi:G\to{\mathcal L}(X)$ into a representation $\rho:G\to{\mathcal L}(Y)$}, if
$$
\rho(t)\alpha(x)=\alpha(\pi(t)x),\qquad t\in G,\ x\in X;
$$

\item[---] a {\it subrepresentation} of a representation $\rho:G\to{\mathcal L}(Y)$, if $\alpha$ is a morphism, it is injective, and $\alpha(X)$ is closed in $Y$;

\item[---] an {\it isomorphism of representations} $\pi:G\to{\mathcal L}(X)$ and $\rho:G\to{\mathcal L}(Y)$, if  $\alpha$ is a morphism of these representations and an isomorphism of Hilbert spaces $X$ and $Y$.

}\eit

A representation $\pi:G\to{\mathcal L}(X)$ is said to be
 \bit{
\item[---] {\it irreducible}, if it has no non-trivial (i.e. different from $X$ and 0) subrepresentations,

\item[---] {\it semisimple}, if it is isomorphic to a Hilbert direct sum of a family of irreducible representations:
 \beq\label{DEF:semisimple}
\pi=\dot{\bigoplus_{i\in I}}\sigma_i,
 \eeq

\item[---] {\it isotypic}, if it is semisimple, and in the decomposition \eqref{DEF:semisimple} all the representations $\sigma_i$ are isomorphic.

 }\eit

For a locally compact group $G$ we denote by $\widehat{G}$ its {\it dual object}, i.e. a set of irreducible unitary representations of $G$ such that
 \bit{
\item[---] any two representations $\rho,\sigma\in\widehat{G}$ are isomorphic if and only if they coincide:
$$
\rho\cong\sigma\qquad\Longleftrightarrow\qquad \rho=\sigma,
$$

\item[---] each irreducible unitary representation $\pi:G\to{\mathcal L}(X)$ is isomorphic to some $\sigma\in\widehat{G}$.
}\eit
The dual object $\widehat{G}$ exists for all locally compact groups $G$.

Suppose $K$ is a compact normal subgroup in a locally compact group $G$, and $\mu_K$ the normed Haar measure on $K$. For each representation $\sigma\in\widehat{K}$, $\sigma:K\to {\mathcal B}(X_\sigma)$ let us define a measure $\nu_\sigma\in{\mathcal C}^\star(G)$ by the formula
\beq\label{DEF:nu_sigma-1}
\nu_\sigma(u)=\dim X_\sigma\cdot \int_K \overline{\tr\sigma(s)}\cdot u(s)\cdot \mu_K(\d s)=\dim X_\sigma\cdot \int_K \tr\sigma(s^{-1})\cdot u(s)\cdot \mu_K(\d s),\qquad u\in{\mathcal C}(G).
\eeq
or, equivalently, by the formula
\beq\label{DEF:nu_sigma}
\nu_\sigma=\dim X_\sigma\cdot \int_K \overline{\tr\sigma(s)}\cdot\delta_s\cdot \mu_K(\d s)=\dim X_\sigma\cdot \int_K \tr\sigma(s^{-1})\cdot\delta_s\cdot \mu_K(\d s)
\eeq

 \vglue10pt
\centerline{\bf Properties of $\nu_\sigma$:}
 \vglue10pt

\bit{\it

\item[$1^\circ$.] For each $\sigma\in\widehat{K}$ the measure $\nu_\sigma$ is central in ${\mathcal C}^\star(G)$:
 \beq\label{chi_pi-centr}
 \nu_\sigma*\alpha=\alpha*\nu_\sigma,\qquad \alpha\in {\mathcal C}^\star(G),
 \eeq

\item[$2^\circ$.] The measures $\nu_\sigma$ for a system of orthogonal projectors in the algebra ${\mathcal C}^\star(G)$:
 \beq\label{chi_pi-ortog-proj}
 \nu_\sigma*\nu_\sigma=\nu_\sigma,\qquad \nu_\sigma*\nu_\rho=0,\qquad \sigma\ne\rho\in \widehat{K}.
 \eeq

}\eit

\bpr
1. Since delta-functionals are total in ${\mathcal C}^\star(G)$, \eqref{chi_pi-centr} is equivalent to the identity
$$
\delta^t*\nu_\sigma*\delta^{t^{-1}}=\nu_\sigma,\qquad t\in G.
$$
or, by \eqref{DEF:svertka-s-delta^a}, to the identity
$$
t\cdot\nu_\sigma\cdot t^{-1}=\nu_\sigma,\qquad t\in G.
$$
We use Theorem \ref{TH:mu(ph(X))=mu(X)} for its proof: since the map $x\mapsto t\cdot x\cdot t^{-1}$ is an automorphism of the group $K$, it preserves the Haar measure $\mu_K$. As a corollary, in the following chain the change of variable is valid:
\begin{multline*}
t\cdot\nu_\sigma\cdot t^{-1}(u)=\nu_\sigma(t^{-1}\cdot u\cdot t)=\dim X_\sigma\cdot \int_K\tr\sigma(s^{-1})
(t^{-1}\cdot u\cdot t)(s)\cdot \mu_K(\d s)=\\=
\dim X_\sigma\cdot \int_K\tr\sigma(s^{-1})
u(t\cdot s\cdot t^{-1})\cdot \mu_K(\d s)=\begin{vmatrix}t\cdot s\cdot t^{-1}=r \\ s=t^{-1}\cdot r\cdot t\\ \mu_K(\d s)=\mu_K(\d r)\end{vmatrix}=
\dim X_\sigma\cdot \int_K\tr\sigma(t^{-1}\cdot r^{-1}\cdot t) u(r)\cdot \mu_K(\d r)=\\=
\dim X_\sigma\cdot \int_K\tr\sigma(r^{-1}) u(r)\cdot \mu_K(\d r)=\nu_\sigma(u).
\end{multline*}

2. Denote by $\chi_\sigma$ the characters of the representations $\sigma\in\widehat{K}$:
$$
\chi_\sigma(s)=\tr\sigma(s).
$$
As is known \cite[(27.24)]{Hewitt-Ross-2},
\beq\label{e_sigma*e_tau}
\chi_\sigma*\chi_\tau=\begin{cases}0, &\sigma\ne\tau\\ \frac{1}{\dim X_\sigma}\cdot\chi_\sigma,& \sigma=\tau\end{cases}.
\eeq
On the other hand, obviously,
\beq\label{chi_sigma=dim X_sigma-widetilde(e_sigma)-mu_K}
\nu_\sigma=\dim X_\sigma\cdot\widetilde{\chi_\sigma}\cdot\mu_K.
\eeq
Thus,
\begin{multline*}
\nu_\sigma*\nu_\tau=(\dim X_\sigma\cdot\widetilde{\chi_\sigma}\cdot\mu_K)*(\dim X_\tau\cdot\widetilde{\chi_\tau}\cdot\mu_K)=\eqref{(f-mu_K)*(g-mu_K)=(f*g)-mu_K}=
\dim X_\sigma\cdot\dim X_\tau\cdot(\widetilde{\chi_\sigma}*\widetilde{\chi_\tau})\cdot\mu_K=\eqref{widetilde(f*g)}=\\=
\dim X_\sigma\cdot\dim X_\tau\cdot(\widetilde{\chi_\tau*\chi_\sigma})\cdot\mu_K=\eqref{e_sigma*e_tau}=
\begin{Bmatrix}
0, &\sigma\ne\tau\\ \dim X_\sigma\cdot\dim X_\sigma\cdot\frac{1}{\dim X_\sigma}\cdot\widetilde{\chi_\sigma}
\cdot\mu_K,& \sigma=\tau
\end{Bmatrix}=\\=
\begin{Bmatrix}
0, &\sigma\ne\tau\\ \dim X_\sigma\cdot\widetilde{\chi_\sigma}
\cdot\mu_K,& \sigma=\tau
\end{Bmatrix}=\eqref{chi_sigma=dim X_sigma-widetilde(e_sigma)-mu_K}=
\begin{Bmatrix}
0, &\sigma\ne\tau\\ \nu_\sigma,& \sigma=\tau
\end{Bmatrix}.
\end{multline*}
\epr

Let $\pi:G\to {\mathcal L}(X)$ be a unitary representation. For any $\sigma\in\widehat{K}$ we put
 \beq\label{DEF:T_pi}
 \pi_\sigma(t)=\dot{\pi}(\nu_\sigma*\delta_t),\qquad  \dot{\pi}_\sigma(\alpha)=\dot{\pi}(\nu_\sigma*\alpha),\qquad t\in G,\quad \alpha\in {\mathcal C}^\star(G),\qquad \sigma\in \widehat{K}.
 \eeq

 \vglue10pt
\centerline{\bf Properties of $\pi_\sigma$:}
 \vglue10pt

\bit{\it

\item[$1^\circ$.] For any $\sigma\in\widehat{K}$ the mapping $\pi_\sigma:G\to {\mathcal L}(X)$ is a representation of $G$, and $\dot{\pi}_\sigma:{\mathcal C}^\star(G)\to {\mathcal L}(X)$ is its extension to the group algebra:
\beq
\dot{\pi}_\sigma=(\pi_\sigma)^\cdot
\eeq

\item[$2^\circ$.] For any $\sigma\in\widehat{K}$ the space
\beq
X_\sigma=\overline{\dot{\pi}_\sigma({\mathcal C}^\star(G))X}
\eeq
is invariant with respect to $\pi$ (and therefore defines a subrepresentation of $\pi$).

\item[$3^\circ$.] The space $X$ can be decomposed into a Hilbert direct sum of the spaces $X_\sigma$:
\beq\label{X=oplus_pi-X_pi}
X=\dot{\bigoplus_{\sigma\in\widehat{K}}}X_\sigma.
\eeq
and the representation $\pi$ into the Hilbert direct sum of representations $\pi_\sigma$:
\beq\label{T=oplus_pi-T_pi}
\pi=\dot{\bigoplus_{\sigma\in\widehat{K}}}\pi_\sigma.
\eeq

\item[$4^\circ$.] Being restricted to the subgroup $K$ the decomposition \eqref{T=oplus_pi-T_pi} turns into the decomposition $\pi|_K$ by isotypic components multiple to $\sigma$:
\beq\label{T=oplus_pi-T_pi|_K}
\pi|_K^\cdot=\dot{\bigoplus_{\sigma\in\widehat{K}}}\pi_\sigma|_K^\cdot, \qquad \pi_\sigma|_K^\cdot=\dot{\bigoplus_{i\in I_\sigma}}\sigma_i^\cdot,\qquad \sigma_i\cong\sigma\qquad (i\in I_\sigma).
\eeq
}\eit
\bpr
1. First,
\begin{multline*}
\pi_\sigma(s\cdot t)=\dot{\pi}(\nu_\sigma*\delta_{s\cdot t})=\dot{\pi}(\nu_\sigma*\delta_s*\delta_t)=\eqref{chi_pi-ortog-proj}=
\dot{\pi}(\nu_\sigma*\nu_\sigma*\delta_s*\delta_t)=\\=\eqref{chi_pi-centr}=\dot{\pi}(\nu_\sigma*\delta_s *\nu_\sigma*\delta_t)=
\dot{\pi}(\nu_\sigma*\delta_s)\cdot \dot{\pi}(\nu_\sigma*\delta_t)=\pi_\sigma(s)\cdot \pi_\sigma(t)
\end{multline*}

2. Take $x\in \dot{\pi}_\sigma({\mathcal C}^\star(G))X$, i.e.
$$
x=\sum_{i=1}^n \dot{\pi}_\sigma(\beta_i)x_i,\qquad \beta_i\in{\mathcal C}^\star(G),\ x_i\in X.
$$
Then for any $\alpha\in{\mathcal C}(G)$
\begin{multline*}
\dot{\pi}(\alpha)x=\sum_{i=1}^n \dot{\pi}(\alpha)\dot{\pi}_\sigma(\beta_i)x_i=\sum_{i=1}^n \dot{\pi}(\alpha)\dot{\pi}(\nu_\sigma*\beta_i)x_i=\sum_{i=1}^n \dot{\pi}(\alpha*\nu_\sigma*\beta_i)x_i=\\=
\eqref{chi_pi-centr}=\sum_{i=1}^n \dot{\pi}(\nu_\sigma*\alpha*\beta_i)x_i=\sum_{i=1}^n \dot{\pi}_\sigma(\alpha*\beta_i)x_i
\in \dot{\pi}_\sigma(G)X.
\end{multline*}
Hence
$$
\dot{\pi}({\mathcal C}^\star(G))\Big(\dot{\pi}_\sigma({\mathcal C}^\star(G))X\Big)\subseteq \dot{\pi}_\sigma({\mathcal C}^\star(G))X,
$$
and
$$
\dot{\pi}({\mathcal C}^\star(G))\Big(\overline{\dot{\pi}_\sigma({\mathcal C}^\star(G))X}\Big)\subseteq \overline{\dot{\pi}_\sigma({\mathcal C}^\star(G))X}.
$$

3. We prove $3^\circ$ and $4^\circ$ together. Consider the restriction $\pi|_K:K\to{\mathcal L}(X)$. As is known, it is decomposed into the Hilbert direct sum of irreducible representations:
$$
\pi|_K=\dot{\bigoplus_{\sigma\in\widehat{K}}}\dot{\bigoplus_{i\in I_\sigma}}\sigma_i,\qquad \sigma_i\cong\sigma\qquad (i\in I_\sigma).
$$
The projection on the isotypic component
$$
\varPhi_\sigma:X=\dot{\bigoplus_{\sigma\in\widehat{K}}}\dot{\bigoplus_{i\in I_\sigma}}X_i\to \dot{\bigoplus_{i\in I_\sigma}}X_i
$$
and at the same time the intertwinner between the representations $\pi$ and $\pi_\sigma$ (where $\pi_\sigma$ is defined in \eqref{DEF:T_pi}),
$$
\varPhi_\sigma\circ\pi=\pi_\sigma\circ\varPhi_\sigma,
$$
is described by the formula \cite[Theorem 5.5.1(2)]{Kowalski}
\beq\label{DEF:varPhi_sigma}
\varPhi_\sigma=\dim X_\sigma\cdot\int_K\overline{\tr\sigma(s)}\cdot \pi(s)\cdot\mu_K(\d s)=\dim X_\sigma\cdot\int_K \tr\sigma(s^{-1})\cdot \pi(s)\cdot\mu_K(\d s)
\eeq
Note that
$$
\dot{\pi}_\sigma(\alpha)=\dot{\pi}(\nu_\sigma*\alpha)=
\dot{\pi}(\nu_\sigma)\cdot\dot{\pi}(\alpha)=\varPhi_\sigma\cdot \dot{\pi}(\alpha),\qquad \alpha\in{\mathcal C}^\star(G).
$$
As a corollary,
$$
X_\sigma=\overline{\dot{\pi}_\sigma({\mathcal C}^\star(G))X}=\overline{\varPhi_\sigma \dot{\pi}({\mathcal C}^\star(G))X}=
\overline{\varPhi_\sigma X}=\overline{\dot{\bigoplus_{i\in I_\sigma}}X_i}=\dot{\bigoplus_{i\in I_\sigma}}X_i.
$$
We have a decomposition into a Hilbert sum
$$
X=\dot{\bigoplus_{\sigma\in\widehat{K}}}\dot{\bigoplus_{i\in I_\sigma}}X_i=\dot{\bigoplus_{\sigma\in\widehat{K}}}X_\sigma.
$$
This proves \eqref{T=oplus_pi-T_pi|_K} and \eqref{T=oplus_pi-T_pi}.
\epr

\subsection{Norm continuous representations}

\paragraph{Central groups, SIN-groups and Moore groups.}
Let us remind that the center $Z(G)$ of a group $G$ is the set of all elements $a\in G$ which commute with the all other elements:
$$
a\cdot x=x\cdot a,\qquad x\in G.
$$
A locally compact group $G$ is called {\it central}\label{DEF:zentralnaya-gruppa}, if its quotient group by the center $G/Z(G)$ is compact.

\btm[H.Freudenthal \cite{Freudenthal}, see also \cite{Palmer}]\label{TH:stroenie-zentralnoi-gruppy}
A connected locally compact group $G$ is central if and only if it is a direct product of a vector group $\R^n$ and a compact group $K$:
 \beq\label{stroenie-zentralnoi-gruppy}
 G=\R^n\times K.
 \eeq
\etm

Let us call a locally compact group $G$ a {\it compact buildup of a commutative group}\label{DEF:komp-nadstr-abelevoi-gruppy}, if there are closed subgroups $C$ and $K$ in $G$ such that:
\bit{

\item[1)] $C$ is an commutative group,

\item[2)] $K$ is a compact group,

\item[3)] $C$ and $K$ commute:
$$
\forall a\in C,\quad \forall y\in K\qquad a\cdot y=y\cdot a,
$$

\item[4)] the product of $C$ and $K$ is $G$:
$$
\forall x\in G\qquad \exists a\in C\quad\exists y\in K\qquad x=a\cdot y.
$$
}\eit
Certainly, every such a group is central.

\btm[Yu.~Kuznetsova \cite{Kuznetsova}]\label{TH:zentral-Lie}
Each central Lie group $G$ is a finite extension of some compact buildup of a commutative group.
\etm
\bpr
By \cite[Theorem 4.4]{Grosser-Moskowitz-1}, $G$ can be represented in the form $G=\R^n\times G_1$, where $G_1$ contains an open compact normal subgroup $K$. Let $C$ be the center of $G$. Let us identify $K$ with $\{0\}\times K$ and consider the group $H=C\cdot K$. Of course, this is a compact buildup of a commutative group. Besides this, $H$ contains an open subgroup $\R^n\times K$, hence $H$ is open. Since $K$ is normal in $G_1$ the group $K=\{0\}\times K$ is normal in $G$, therefore $H$ is normal in $G$:
$$
x\cdot H=x\cdot C\cdot K=C\cdot x\cdot K=C \cdot K\cdot x=H\cdot x.
$$
Since $G$ is central, the quotient group $G/C$ is compact. As a corollary, $G/H$ is also compact (since it is a continuous image of a compact group $G/C$). On the other hand, since $H=C\cdot K$ is open, $G/H$ must be discrete. Hence, $G/H$ is finite.
\epr

A set $U$ in a group $G$ is said to be {\it normal}, if it is invariant with respect to conjugations:
$$
a\cdot U\cdot a^{-1}\subseteq U,\qquad a\in G.
$$
A locally compact group $G$ is called a {\it SIN-group}\label{DEF:SIN-gruppa}, if normal neighbourhoods of identity in $G$ form a local base. An equivalent condition: the left and the right uniform structures in $G$ are equivalent. In particular, all such groups are unimodular.

The class of SIN-groups includes Abelian, compact and discrete groups. The following result belongs to S.~Grosser and M.~Moskowitz \cite[2.13]{Grosser-Moskowitz-2}:

\btm\label{TH:stroenie-SIN}
Each SIN-group $G$ is a discrete extension of the group $\R^n\times K$, where $n\in\Z_+$, and $K$ is a compact group:
\beq\label{SIN-kak-rasshirenie}
1\to \R^n\times K=N\to G\to D\to 1
\eeq
($D$ is a discrete group).
\etm

A locally compact group $G$ is called a {\it Moore group}\label{DEF:Moore-gruppa}, if all its continuous (in the usual sense, i.e. as it is described in Theorem \ref{TH:nepr-predstavleniya}) unitary continuous representations are finite-dimensional.

\btm\label{TH:Moore->SIN}
Each Moore group is a SIN-group.\footnote{See \cite[p.1452]{Palmer}.}
\etm

\btm\label{TH:Moore->AM}
Each Moore group is amenable.\footnote{See \cite[p.1486]{Palmer}.}
\etm

\btm\label{TH:Moore->quotient}
Each (Hausdorff) quotient group $G/H$ of a Moore group $G$ is a Moore group.\footnote{This proposition is obvious.}
\etm

\bcor\label{TH:G=Moore->D=Moore}
If $G$ is a Moore group, then in its representation \eqref{SIN-kak-rasshirenie} the group $D$ is also a Moore group (and in particular it is amenable).
\ecor

\btm\label{TH:D=Moore->konech-rassh-Abelevoi}
Each discrete Moore group is a finite extension of a commutative group.\footnote{See \cite[Theorem 12.4.26 and p.1397]{Palmer}.}
\etm

\btm[Yu.~N.~Kuznetsova \cite{Kuznetsova}]\label{TH:Lie-Moore}
Each Lie-Moore group $G$ is a finite extension of some compact buildup of a commutative group.
\etm
\bpr
By \cite[Theorem 12.4.27]{Palmer}, $G$ is a finite extension of some central group $G_1$. By Theorem  \ref{TH:zentral-Lie}, $G_1$ is a finite extension of some compact buildup of a commutative group.
\epr

\paragraph{Norm continuous representations.}

Let us return back to the situation when we have a representation $\pi:G\to A$ of a group $G$ in a stereotype algebra $A$ (see diagram \eqref{C*(G)=grup-alg-A}). Suppose that $A$ is a Banach algebra. Then the continuity of $\pi$ will be the norm-continuity:
$$
x_i\to x\quad\Longrightarrow\quad \norm{\pi(x_i)-\pi(x)}\to 0.
$$

\bex\label{EX:gladkost-nepr-po-norme-predstavleniya}
If $G$ is a real Lie group, then every its norm-continuous representation $\pi:G\to A$ in an arbitrary Banach algebra $A$ is a smooth mapping, and, as a corollary, $\pi$ has a (unique) extension (in the sense of diagram \eqref{E*(G)=grup-alg-A}) to a continuous homomorphism $\ddot{\pi}:{\mathcal E}^\star(G)\to A$ of the group algebra ${\mathcal E}^\star(G)$ of distributions with compact support (described at page \pageref{E*(G)=grup-alg-A}).
\eex
\bpr
The smoothness of $\pi$ is proved by analogy with the classical result, when $A$ is a finite-dimensional algebra, see \cite[4.21]{Michor}. After that the existence of $\ddot{\pi}$ follows from \cite[10.12]{Akbarov}.
\epr

Let us denote by ${\mathcal B}(X)$ the usual Banach space of linear continuous operators in a Hilbert space $X$ (i.e. the same space ${\mathcal L}(X)$, but with the topology, generated by the norm of operator). Let $\iota:{\mathcal B}(X)\to{\mathcal L}(X)$ be the usual embedding (certainly, it is continuous).

The following proposition is obvious.

\btm\label{TH:nepr-po-norme-predstavleniya} For a representation $\pi:G\to{\mathcal L}(X)$ of a locally compact group $G$ in a Hilbert space $X$ the following conditions are equivalent:
 \bit{
\item[(i)] the mapping $\pi:G\to{\mathcal B}(X)$ is continuous,

\item[(ii)] the morphism of stereotype algebras $\dot{\pi}:\mathcal{C}^\star(G)\to {\mathcal L}(X)$ from \eqref{C*(G)=grup-alg} has a lifting to a morphism of stereotype algebras $\ph:\mathcal{C}^\star(G)\to {\mathcal B}(X)$:
\begin{gather}\label{C*(G)=grup-alg-B(X)}
\begin{diagram}
\node[2]{\mathcal{C}^\star(G)} \arrow{sw,t}{\ph}\arrow{se,t}{\dot{\pi}}
\\
\node{{\mathcal B}(X)}\arrow[2]{e,b}{\iota} \node[2]{{\mathcal L}(X)}
\end{diagram}.
\end{gather}
 }\eit
\etm

\bit{
\item[$\bullet$] The representation $\pi:G\to{\mathcal L}(X)$ satisfying these conditions is said to be {\it norm-continuous}.\label{DEF:nepr-po-norme-predst}

\item[$\bullet$] In the special case when the operators $\pi(t)$, $t\in G$, are unitary, the representation  $\pi:G\to{\mathcal L}(X)$ that satisfies the conditions of this theorem is called a {\it norm-continuous unitary representation}.\label{DEF:nepr-po-norme--unit-predst}

}\eit

\btm \cite[Corollary 2]{Shtern}\label{TH:nepr-po-norme-perdst-K} A unitary representation $\pi:K\to{\mathcal L}(X)$ of a compact group $K$ is norm-continuous if and only if in its decomposition \eqref{T=oplus_pi-T_pi|_K}
$$
\dot{\pi}=\dot{\bigoplus_{\sigma\in\widehat{K}}}\dot{\pi}_\sigma,
$$
only finite number of isotypic components $\dot{\pi}_\sigma$ do not vanish.
\etm

\btm \cite[Theorem 4.8]{Kuznetsova}\label{TH:Kuz-2}
If $K$ is a normal compact subgroup in $G$ and a unitary representation $\pi:G\to {\mathcal L}(X)$ is norm-continuous, then in the decomposition \eqref{X=oplus_pi-X_pi} $X_\pi\ne 0$ only for finite number of indices $\sigma\in\widehat{G}$. As a corollary, the sum \eqref{T=oplus_pi-T_pi} is finite in this case.
\etm

\paragraph{Induced representations.}

Let us remind the construction of induced representation. Suppose $N$ is an open normal subgroup in a locally compact group $G$, and $\pi:N\to {\mathcal L}(X)$ -a unitary representation of $N$. Let us choose a mapping $\sigma:D\to G$ which is a coretraction for the quotient mapping $\ph:G\to G/N=D$
\beq\label{ph(sigma(t))=t}
\ph(\sigma(t))=t,\qquad t\in D,
\eeq
and preserves the unit:
\beq\label{sigma(1_D)=1_G}
\sigma(1_D)=1_G.
\eeq
Then for any $g\in G$ the elements $g$ and $\sigma(\ph(g))$ belong to the same coset with respect to $N$,
$$
g\in \sigma(\ph(g))\cdot N
$$
i.e.
\beq\label{g-sigma(ph(g))^(-1)-in-N}
g\cdot \sigma(\ph(g))^{-1}\in N,\qquad g\in G.
\eeq

Consider the space $L_2(D,X)$ of square-summable mappings $\xi:D\to X$ (with respect to the counting measure $\card$ on $D$). Then the induced representation $\pi':G\to {\mathcal L}(L_2(D,X))$ is defined by the formula
\beq\label{ind-representation}
\pi'(g)(\xi)(t)=\pi\Big(\underbrace{\sigma(t)\cdot g\cdot \sigma\big(\ph(\sigma(t)\cdot g)\big)^{-1}}_{\scriptsize\begin{matrix}\phantom{\quad\tiny \eqref{g-sigma(ph(g))^(-1)-in-N}}
\text{\rotatebox{90}{$\owns$}}{\quad\tiny \eqref{g-sigma(ph(g))^(-1)-in-N}}\\ N\end{matrix}}\Big)\Big(\xi\big(\underbrace{\ph(\sigma(t)\cdot g)}_{\scriptsize\begin{matrix}
\text{\rotatebox{90}{$\owns$}}\\ D\end{matrix}}\big)\Big),\qquad \xi\in L_2(D,X),\quad t\in D,\quad g\in G.
\eeq

\btm\cite[Lemma 3.6]{Kuznetsova}\label{TH:Kuz-1}
Let $G$ be a SIN-group, $N$ its open subgroup in the chain \eqref{SIN-kak-rasshirenie}, and $\pi:N\to{\mathcal L}(X)$ a norm-continuous unitary representation. Then the induced representation $\pi':G\to{\mathcal L}(L_2(D,X))$ is also norm-continuous.
\etm
\bpr
Take $\e>0$. Since $\pi$ is norm-continuous, there is a neighbourhood of identity $U\subseteq N$ such that
$$
\norm{\pi(h)-1_X}<\e,\qquad h\in U.
$$
At the same time, $N$ is open in $G$, therefore $U$ is also a neighbourhood of identity in $G$. Since $G$ is a SIN-group, there is a neighbourhood of identity $V\subseteq U$ invariant with respect to conjugations: $g\cdot V\cdot g^{-1}\subseteq V$ for all $g\in G$. Now for any $h\in V\subseteq U\subseteq N$ and for any $t\in D$ we have
$$
\ph(\sigma(t)\cdot h)=\ph(\sigma(t))\cdot\ph(h)=t\cdot 1=t,
$$
and thus, for $h\in V$ and $\xi\in L_2(D,X)$
\begin{multline*}
\norm{\pi'(h)(\xi)-\xi}^2=\sum_{t\in D}\norm{\pi'(h)(\xi)(t)-\xi(t)}^2=
\sum_{t\in D}\norm{\pi\Big(\sigma(t)\cdot h\cdot\sigma\big(\ph(\sigma(t)\cdot h)\big)^{-1}\Big)\Big(\xi\big(\ph(\sigma(t)\cdot h)\big)\Big)-\xi(t)}^2=\\=
\sum_{t\in D}\norm{\pi\big(\sigma(t)\cdot h\cdot\sigma(t)^{-1}\big)\big(\xi(t)\big)-\xi(t)}^2\le
\sum_{t\in D}\|\pi\big(\underbrace{\sigma(t)\cdot h\cdot\sigma(t)^{-1}}_{\scriptsize\begin{matrix}\text{\rotatebox{90}{$\owns$}}\\ V\end{matrix}}\big)-1_X\|^2\cdot\norm{\xi(t)}^2<
\e^2\cdot \sum_{t\in D}\norm{\xi(t)}^2=\e^2\cdot\norm{\xi}^2
\end{multline*}
\epr

\paragraph{The space $\Trig(G)$ of norm-continuous trigonometric polynomials.}
Let us call a {\it norm-continuous trigonometric polynomial} on $G$ an arbitrary linear combination of matrix elements of norm-continuous unitary irreducible representations of $G$
\beq\label{DEF:Trig(G)}
 u(t)=\sum_{i=1}^n\lambda_i\cdot\langle \sigma_i(t)x_i,y_i\rangle,\quad t\in G,
\eeq
($\sigma_i:G\to{\mathcal B}(X_i)$ are norm continuous unitary irreducible representations, $x_i,y_i\in X_i$, $\lambda_i\in\C$).

\bex\label{EX:trig-monom}
Let $G$ be a locally compact group, $U$ an arbitrary $C^*$-neighbourhood of zero in ${\mathcal C}^\star(G)$ and $v$ a pure state of the $C^*$-algebra ${\mathcal C}^\star(G)/U$. Then the function
$$
G\overset{\delta}{\longrightarrow}{\mathcal C}^\star(G)\overset{\pi_U}{\longrightarrow}{\mathcal C}^\star(G)/U\overset{v}{\longrightarrow}\C
$$
is a norm-continuous trigonometric polynomial on $G$.
\eex
\bpr
Consider the GNS-representation $\sigma:{\mathcal C}^\star(G)/U\to {\mathcal B}(X)$, generated by the functional $v$:
$$
v(a)=\langle \sigma(a)x,x\rangle,\qquad a\in {\mathcal C}^\star(G)/U,
$$
($x\in X$). Since $v$ is a pure state, the representation $\sigma$ is irreducible \cite[10.2.3]{Kad-Ring-2}. Since $G$ is totally mapped into the algebra ${\mathcal C}^\star(G)$ (i.e. the linear combinations of the delta-functions are dense in ${\mathcal C}^\star(G)$) and therefore in the quotient algebra ${\mathcal C}^\star(G)/U$, we obtain that the composition
$$
G\overset{\delta}{\longrightarrow}{\mathcal C}^\star(G)\overset{\pi_U}{\longrightarrow}{\mathcal C}^\star(G)/U\overset{\sigma}{\longrightarrow}{\mathcal B}(X)
$$
is an irreducible representation of the group $G$. Hence the function
$$
u(t)=v(\pi_U(\delta(t)))=\langle \sigma(\pi_U(\delta(t)))x,x\rangle
$$
is a norm-continuous trigonometric polynomial.
\epr

\brem
If the group $G$ is compact or Abelian, then the space $\Trig(G)$ is an algebra with respect to the operation of pointwise multiplication. However, apparently, there exist (non-compact and non-Abelian) groups where this is not so. If we don't demand that the trigonometric polynomial is continuous with respect to a $C^*$-seminorm on ${\mathcal C}^\star(G)$, then the counterexample is the group $G=H_3(\R)$ of upper triangle matices with real coefficients and units on the diagonal: the space $\Trig(G)$ (without the claim of ``$C^*$-seminorm-continuity'') is not closed with respect to the pointwise multiplication -- this example belongs to Yemon Choi.
\erem

\paragraph{The algebra $k(G)$ of norm-continuous matrix elements.}
 A {\it norm-continuous matrix element} on a group $G$ is a function of the form
\beq\label{DEF:k(G)}
 u(t)=\langle \pi(t)x,y\rangle,\quad t\in G,
\eeq
where $\pi:G\to{\mathcal B}(X)$ is a norm-continuous unitary representation in a Hilbert space $X$, and $x,y\in X$. The set of all such functions will be denoted by $k(G)$.

\btm
For each locally compact group $G$ the space $k(G)$ forms an involutive algebra over $\C$ with the pointwise algebraic operations.
\etm
\bit{
\item[$\bullet$] We endow the algebra $k(G)$ with the strongest locally convex topology over $\C$. This turns $k(G)$ into the involutive stereotype algebra (due to Example \ref{EX:siln-loc-vyp-topol}).
}\eit
\bpr
The multiplication by a scalar in $k(G)$ is equivalent to the multiplication of the vector $x$ by this scalar in the formula \eqref{DEF:k(G)}. Let us take two functions of the form \eqref{DEF:k(G)}
$$
u(t)=\langle \pi(t)x,y\rangle,\quad v(t)=\langle \pi'(t)x',y'\rangle,\qquad t\in G,
$$
where $\pi:G\to{\mathcal B}(X)$ and $\pi':G\to{\mathcal B}(X')$ are norm continuous unitary representations. If we consider the Hilbert direct sum $X\oplus X'$, i.e. the space of pairs $(\xi,\xi')$, $\xi\in X$, $\xi'\in X'$, with the scalar product
$$
\langle (\xi,\xi'),(\upsilon,\upsilon')\rangle=\langle \xi,\upsilon\rangle+\langle \xi',\upsilon'\rangle,
$$
then the formula
$$
(\pi\oplus\pi')(t)(\xi,\xi')=(\pi(t)\xi,\pi'(t)\xi'),\qquad t\in G,
$$
defines a norm-continuous representation $\pi\oplus\pi':G\to{\mathcal B}(X\oplus X')$ with the sum $u+v$ as a matrix element:
$$
\langle (\pi\oplus\pi')(t)(x,x'),(y,y')\rangle=\langle (\pi(t)x,\pi'(t)x'),(y,y')\rangle=
\langle \pi(t)x,y\rangle+\langle \pi'(t)x',y'\rangle=u(t)+v(t).
$$
On the other hand, if we consider the Hilbert tensor product $X\dot{\otimes} X'$ \cite{Kad-Ring}, i.e. the completion of the algebraic tensor product $X\otimes X'$ with respect to the scalar product
$$
\langle \xi\otimes\xi',\upsilon\otimes\upsilon'\rangle=\langle \xi,\upsilon\rangle\cdot \langle \xi',\upsilon'\rangle, \qquad \xi,\upsilon\in X,\quad \xi',\upsilon'\in X',
$$
then the formula
$$
(\pi\otimes\pi')(t)\xi\otimes\xi'=\pi(t)\xi\otimes\pi'(t)\xi',\qquad t\in G,
$$
defines a norm-continuous representation $\pi\otimes\pi':G\to{\mathcal B}(X\dot{\otimes} X')$ with the product $u\cdot v$ as a matrix element:
$$
\langle (\pi\otimes\pi')(t)x\otimes x',y\otimes y'\rangle=\langle \pi(t)x\otimes\pi'(t)x',y\otimes y'\rangle=
\langle \pi(t)x,y\rangle\cdot\langle \pi'(t)x',y'\rangle=u(t)\cdot v(t).
$$
From the existence of the conjugate linear finction on $X$ \cite[27.25,27.26]{Hewitt-Ross-2} it follows that $k(G)$ is closed with respect to the pointwise involution. The unit belongs to $k(G)$ since it can be represented as  \eqref{DEF:k(G)}, if we put $\pi(t)=1$ and take $x$ and $y$ such that $\langle x,y\rangle=1$.
\epr

Let us say that a function $u\in{\mathcal C}(G)$ {\it is subordinated to a  $C^*$-seminorm},\label{DEF:positive-def-func-sub-C*} if for some (continuous) $C^*$-seminorm $p$ on ${\mathcal C}^\star(G)$ we have
$$
|\alpha(u)|\le p(\alpha),\qquad \alpha\in{\mathcal C}^\star(G).
$$
Certainly, all the functions from $k(G)$ are subordinated to $C^*$-semonorms.

Recall also that a function $u\in{\mathcal C}(G)$ is said to be {\it positive definite} \cite[\S 32]{Hewitt-Ross-2}, if
\beq\label{(alpha^bullet*alpha)(u)-ge-0}
(\alpha^\bullet*\alpha)(u)\ge 0,\qquad \alpha\in {\mathcal C}^\star(G).
\eeq

\btm\label{TH:u(t)=langle-pi(t)x,x-rangle} For a function $u\in{\mathcal C}(G)$ the following conditions are equivalent:
\bit{
\item[(i)] $u$ can be represented in the form
\beq\label{u(t)=langle-pi(t)x,x-rangle}
u(t)=\langle\pi(t)x,x\rangle,\qquad t\in G,
\eeq
for some norm-continuous unitary representation $\pi:G\to{\mathcal B}(X)$ and some $x\in X$ (and, as a corollary,  $u\in k(G)$);

\item[(ii)] $u$ is positive definite and subordinated to some $C^*$-seminorm.
}\eit
\etm
\bpr We only have to verify the implication (i)$\Longleftarrow$(ii). Suppose (ii) holds.
Since $u$ is subordinated to a $C^*$-seminorm $p$, this function can be extended as a functional on the  $C^*$-quotient algebra ${\mathcal C}^\star(G)/p$:
$$
u=v\circ\rho,
$$
where $\rho:{\mathcal C}^\star(G)\to {\mathcal C}^\star(G)/p$ is the quotient mapping, and $v:{\mathcal C}^\star(G)/p\to\C$ is some continuous functional on the $C^*$-algebra ${\mathcal C}^\star(G)/p$. Therewith, since $u$ is positive definite, and the involution and the multiplication in ${\mathcal C}^\star(G)/p$ are inherited from ${\mathcal C}^\star(G)$, $v$ is a positive functional. Consider the GNS-representation $\sigma:{\mathcal C}^\star(G)/p\to {\mathcal B}(X)$, generated by $v$. Then
$$
v(a)=\langle \sigma(a)x,x\rangle,\qquad a\in {\mathcal C}^\star(G)/p,
$$
for some $x\in X$, and if we pout $\pi=\sigma\circ\rho$, then
$$
\alpha(u)=v(\rho(\alpha))=\langle \sigma(\rho(\alpha))x,x\rangle=\langle \pi(\alpha)x,x\rangle,
\qquad \alpha\in {\mathcal C}^\star(G).
$$
\epr

Obviously,
\beq\label{Trig(G)-subseteq-k(G)}
\Trig(G)\subseteq k(G).
\eeq

\btm\label{TH:Trig(G)=k(G)-na-kompaktah}
Suppose $G$ is a compact group. Then
\bit{

\item[(a)] the embedding \eqref{Trig(G)-subseteq-k(G)} turns into an equality
$$
\Trig(G)=k(G),
$$

\item[(b)] the algebra $k(G)$ coincides with the set of matrix elements \eqref{DEF:k(G)} of various {\rm finite-dimensional} continuous unitary representations $\pi:G\to{\mathcal B}(X)$;

\item[(c)] the algebra $k(G)$ coincides with the set of functions of the form
 \beq\label{u(t)=f(pi(t))}
u(t)=f(\pi(t)),\qquad t\in G,
\eeq
where $\pi:G\to{\mathcal B}(X)$ is an arbitrary norm-continuous unitary representation, and $f:{\mathcal B}(X)\to\C$ an arbitrary linear continuous functional.
}\eit
\etm
\bpr
It is sufficient here to verify that each function of the form \eqref{u(t)=f(pi(t))} belongs to $\Trig(G)$. Let  $\pi:G\to{\mathcal B}(X)$ be a norm continuous unitary representation, and $f:{\mathcal B}(X)\to\C$ a linear continuous functional.

1. Consider the case when the Hilbert space $X$ is finite dimensional. Then the representation $\pi:G\to{\mathcal B}(X)$ can be decomposed into the sum of irreducible ones,
$$
\pi=\bigoplus_{i=1}^n\pi_i:G\to{\mathcal B}\left(\bigoplus_{i=1}^nX_i\right).
$$
Take in each $X_i$ an orthogonal normed basis $e_{ij}$. Then the functionals
$$
A\in{\mathcal B}(X)\mapsto \lambda_{kl}^{ij}\cdot\langle Ae_{ij},e_{kl}\rangle\in C
$$
form a basis in the space ${\mathcal B}^\star(X)$. Hence $f\in{\mathcal B}^\star(X)$ is decomposed into the sum
$$
u(t)=f(\pi(t))=\sum_{i,j,k,l}\lambda_{kl}^{ij}\cdot\langle\pi(t)e_{ij},e_{kl}\rangle= \sum_{i,j,k,l}\lambda_{kl}^{ij}\cdot\langle\underbrace{\pi_i(t)e_{ij}}_{ \scriptsize\begin{matrix}\text{\rotatebox{90}{$\owns$}}\\ X_i\end{matrix}}, \underbrace{e_{kl}}_{ \scriptsize\begin{matrix}\text{\rotatebox{90}{$\owns$}}\\ X_k\end{matrix}}\rangle=
\sum_{i,j,l}\lambda_{kl}^{ij}\cdot\langle\underbrace{\pi_i(t)e_{ij}}_{ \scriptsize\begin{matrix}\text{\rotatebox{90}{$\owns$}}\\ X_i\end{matrix}}, \underbrace{e_{il}}_{ \scriptsize\begin{matrix}\text{\rotatebox{90}{$\owns$}}\\ X_i\end{matrix}}\rangle
$$
The last sum lies in $\Trig(G)$, since each $\pi_i$ is irreducible.

2. If $X$ is infinite-dimensional, then by the A.~I.~Shtern theorem \ref{TH:nepr-po-norme-perdst-K} each norm-continuous representation $\pi:G\to{\mathcal B}(X)$ is factored through a finite dimensional unitary representation $\sigma$ (see diagram below, where $\ph$ is a homomorphism of involutive algebras). Put $g=f\circ\ph$. This is a functional on a finite dimensional algebra ${\mathcal B}(Y)$, that generates the same function $u$
$$
\xymatrix @R=2.pc @C=3.0pc % @M=14pt
{
G\ar@/_4ex/[ddr]_{u}\ar[rr]^{\sigma}\ar[dr]_{\pi} &  & {\mathcal B}(Y)\ar@{-->}[dl]^{\ph}\ar@/^4ex/@{-->}[ddl]^{g} \\
& {\mathcal B}(X)\ar[d]^{f} & \\
& \C &
}
$$
and we have already proved that $u\in\Trig(G)$.
\epr

\brem
If the group $G$ is not compact, then \eqref{Trig(G)-subseteq-k(G)} is not necessarily an equality. For example, we can take $G=\Z$. This was told to the author by A.~I.~Degtyarev: we can take the regular representation $\pi:\Z\to L_2(\Z)$, $\pi(k)u(l)=u(l+k)$, and put
$$
\chi_0(k)=\begin{cases}0,& k\ne 0\\ 1,& k=1 \end{cases}
$$
Then we obtain the matrix element
$$
f(k)=\langle \pi(k)\chi_0,\chi_0\rangle=\chi_0(k),\qquad k\in\Z,
$$
which is not representable as a trigonometric polynomial. Indeed, suppose $f$ is a trigonometric polynomial on $\Z$, i.e. a function of the form
$$
f(k)=\sum_{n=1}^N\lambda_n\cdot\e_n^k,
$$
where $\e_n\in\C$, $|\e_n|=1$, $\lambda_n\in\C$. Then for computing the coefficients $\lambda_k$ we can consider the system
$$
\sum_{n=1}^N\lambda_n\cdot\e_n^k=f(k)=\chi_0(k)=0,\qquad k=1,...,N
$$
which is linearly independent, since its determinant is the Vandermond determinant, and it does not vanish. The column of free terms vanish, hence the coefficients $\lambda_n$ also vanish. I.e. $f=0$, and this is impossible since $f(0)=\chi_0(0)=1$.
\erem

Recall that in the theory of algebraic groups (see \cite[Chapter VI]{Chevalley} of \cite[45.6]{Zhelobenko-2}) to each compact Lie group $G$ one assign the algebra of functions on it, called {\it representation ring}. It consists of various matrix elements of continuous finite dimensional representations of $G$. This is exactly the set of functions mentioned in proposition (b) of Theorem \ref{TH:Trig(G)=k(G)-na-kompaktah}. As a corollary, we have

\btm
If $G$ is a compact Lie group, then the algebra $k(G)=\Trig(G)$ coincides with the representation ring of the group $G$.
\etm

Let us recall that a {\it polynomial} on a general linear group $GL_n(\R)$ (or $GL_n(\C)$) is a function of the form
\beq\label{DEF:mnogochlen-na-GL_n(R)}
u(g)=\frac{f(g)}{(\det g)^k},\qquad g\in GL_n(\R),
\eeq
where $f$ is a polynomial of the coefficients of the matrix $g$, and $k\in\Z_+$. The set of all polynomials on $GL_n(\C)$ will be denoted by ${\mathcal P}(GL_n(\C))$. It forms an involutive algebra with respect to the pointwise operatios.

The following result is classical (see \cite[Cahpter VI]{Chevalley} or \cite[45.6]{Zhelobenko-2}):

\btm
For each compact real Lie group $G$

\bit{

\item[--] there is an integer $n\in\N$ and a continuous injective homomorphism (an embedding) of groups $\sigma:G\to GL_n(\C)$;

\item[--] for any such an embedding $\sigma:G\to GL_n(\C)$ the algebra $k(G)$ is isomorphic to the quotient algebra of the algebra ${\mathcal P}(GL_n(\C))$ by the ideal $I$ of polynomails vanishing on $G$,
    \beq\label{k(G)=P(GL_n(C))/I}
    k(G)={\mathcal P}(GL_n(\C))/I,
    \eeq
    and the same isomorphism holds for real parts:
    \beq\label{Re-k(G)=Re-P(GL_n(C))/Re-I}
    \Real k(G)=\Real{\mathcal P}(GL_n(\C))/\Real I;
    \eeq

\item[--] as a corollary, $G$ has a natural structure of a real algebraic group, for which $\Real k(G)$ is the algebra of (real) polynomials,
    $$
    {\mathcal P}(G)=\real k(G)
    $$
    and the set $\C G$ of common zeroes of $I$
    \beq\label{DEF:CG-dlya-komp-grupp}
    \C G=\{g\in GL_n(\C):\ \forall u\in I\ u(g)=0 \}
    \eeq
    (or, what is the same the complex spectrum $\C\Spec k(G)$ of $k(G)$) has a natural structure of a complex algebraic group with $k(G)$ as the algebra of (complex) manifolds:
    $$
    {\mathcal P}(\C G)=k(G).
    $$

\item[--] in the natural embedding $G\subseteq\C G$ the group $G$ is a real form of the group $\C G$, i.e.
\bit{
\item[(a)] $G$ has non-empty intersection with each connected component of $\C G$, and

\item[(b)] the tangent algebra $L(G)$ is the real part of the tangent algebra $L(\C G)$:
\beq\label{L(G)=Real-L(CG)}
L(G)=\Real L(\C G)
\eeq

}\eit
}\eit
\etm

\bcor\label{COR:kasat-prostr-k-Trig(G)}
Let $G$ be a compact real Lie group. Then
\bit{

\item[(i)] each involutive character $s\in\Spec(k(G))$ is a value in a point $a\in G$
\beq\label{Spec(k(G))}
s(u)=u(a),\qquad u\in k(G),
\eeq
hence the spectrum of the algebra $k(G)$ coincides with $G$:
\beq\label{Spec-k(G)=R-Spec-Real-k(G)=G}
\Spec k(G)=\R\Spec\Real k(G)=G
\eeq

\item[(ii)] each involutive tangent vector $\tau\in T_a[k(G)]$ to the algebra $k(G)$ in an arbitrary point $a\in G$ is a derivation along a one-parameter subgroup $x:\R\to G$:
\beq\label{kasat-prostr-k-Trig(G)}
\tau(u)=\lim_{t\to 0}\frac{u(a\cdot x(t))-u(a)}{t},\qquad u\in k(G),
\eeq
and is uniquely extended to a tangent vector of the algebra ${\mathcal E}(G)$ in the point $a\in G$; hence the tangent spaces to $k(G)$ and ${\mathcal E}(G)$ coincide:
\beq\label{T_a[k(G)]=T_a(G)}
T_a[k(G)]=T_a[{\mathcal E}(G)]=T_a(G).
\eeq
}\eit
\ecor
\bpr
Both propositions follow from Example \ref{EX:kasat-prostr-k-P(M)}: formula \eqref{Spec-k(G)=R-Spec-Real-k(G)=G} follows from \eqref{Spec[CP(M)]=RSpec[P(M)]=M}\footnote{Formula \eqref{Spec-k(G)=R-Spec-Real-k(G)=G} holds for a wider class of all compact groups, not necessarily Lie, see \cite[(30.30)]{Hewitt-Ross-2}.}, since $G$ is a real algebraic manifold, and $k(G)$ a complexification of $\Real k(G)$. On the other hand, \eqref{kasat-prostr-k-P(G)} imply that in \eqref{kasat-prostr-k-Trig(G)} one can treat $\tau$ as a derivation along a smooth curve in $G$, and this is equivalent to a derivation along a one-parameter subgroup.
\epr

\section{Locally convex bundles and constructions of differential geometry}

We need some facts from the theory of bundles of topological vector spaces. In this exposition we follow the ideology of M.~J.~Dupr\'e and R.~M.~Gillette \cite{Dupre-Gillette}. First, it will be convenient to introduce the following definition.

\bit{

\item[$\bullet$] A {\it vector bundle} over the field $\C$ is a set of seven
$(\Xi,M,\pi,\cdot,+,\{0_t;\ t\in M\},-)$, where
\bit{

\item[1)] $\Xi$ is a set, called the {\it space of the bundle},

\item[2)] $M$ is a set, called the {\it base of the bundle},

\item[3)] $\pi:\Xi\to M$ is a surjective mapping, called the {\it projection of the bundle},

\item[4)] $\cdot:\C\times\Xi\to \Xi$ is a mapping, called the {\it fiberwise multiplication by scalars},

\item[5)] $+:\Xi\underset{M}{\sqcap} \Xi\to \Xi$ is a mapping, called the {\it fiberwise summing}\footnote{\label{FOOTNOTE:Xi-_M-Xi} Here  $\Xi\underset{M}{\sqcap} \Xi$ is the fiber square of $\Xi$ over $M$, i.e. the subset of the Cartesian square $\Xi\times \Xi$ consisting of pairs $(\xi,\zeta)$ such that $\pi(\xi)=\pi(\zeta)$.}

\item[6)] $\{0_t;\ t\in M\}$ is a family of elements of $\Xi$, called {\it zeroes}.

\item[7)] $-:\Xi\to \Xi$ is a mapping, called the {\it fiberwise minus},

}\eit\noindent
and for each point $t\in M$ the operation $\cdot,+,-$ with the element $0_t$ define on the inverse image  $\Xi_t:=\pi^{-1}(t)\subseteq \Xi$ (called a {\it fiber over the point $t$}) the structure of vector space over $\C$ with  $0_t$ as zero.

\item[$\bullet$] A {\it section} of a vector bundle $(\Xi,M,\pi,\cdot,+,\{0_t;\ t\in M\},-)$ is an arbitrary mapping $x:M\to \Xi$ such that
$$
\pi\circ x=\id_M.
$$

\item[$\bullet$] A {\it subbundle} of a vector bundle $(\Xi,M,\pi,\cdot,+,\{0_t;\ t\in M\},-)$ is a subset $\Psi$ in $\Xi$ whose intersection with each fiber $\pi^{-1}(t)$ is a (non-zero) vector subspace in $\pi^{-1}(t)$. Certainly, $\Psi$ is a vector bundle over the same base with the same structure elements.

}\eit

\subsection{Locally convex bundles}

\bit{

\item[$\bullet$] Let a vector bundle $(\Xi,M,\pi,\cdot,+,\{0_t;\ t\in M\},-)$ over $\C$ be endowed with the following supplementary structure:
\bit{

\item[1)] the space of the bundle $\Xi$ and the base of the bundle $M$ are endowed with topologies in such a way that the projection $\pi:\Xi\to M$ is (not only surjective, but also) a continuous and an open mapping,

\item[2)] a set $\mathcal P$ of functions $p:\Xi\to\R_+$, called {\it seminorms} is defined,

}\eit\noindent
and the following conditions hold:
\bit{

\item[(a)] on each bundle $\pi^{-1}(t)$ the restrictions $p|_{\pi^{-1}(t)}:\Xi_t\to\R_+$ are the system of seminorms that defines a structure of a (Hausdorff) locally convex space on $\pi^{-1}(t)$, and the topology of this space coincides with the topology induced from $\varXi$,

\item[(b)] the fiberwise multiplication by scalars $\C\times \Xi\to \Xi$ is continuous:
$$
\Big(\lambda_i\overset{\C}{\underset{i\to\infty}{\longrightarrow}}\lambda,\quad
\xi_i\overset{\varXi}{\underset{i\to\infty}{\longrightarrow}} \xi\Big)
\quad\Longrightarrow\quad
\lambda_i\cdot\xi_i\overset{\varXi}{\underset{i\to\infty}{\longrightarrow}}
\lambda\cdot\xi.
$$

\item[(c)] the fiberwise summing\footnote{See footnote \ref{FOOTNOTE:Xi-_M-Xi}.} $\Xi\underset{M}{\sqcap} \Xi\to \Xi$ is continuous:
$$
\Big(\xi_i\overset{\varXi}{\underset{i\to\infty}{\longrightarrow}} \xi,\quad \zeta_i\overset{\varXi}{\underset{i\to\infty}{\longrightarrow}} \zeta,\quad \pi(\xi_i)=\pi(\zeta_i),\quad \pi(\xi)=\pi(\zeta)\Big)
\quad\Longrightarrow\quad
\xi_i+\zeta_i\overset{\varXi}{\underset{i\to\infty}{\longrightarrow}} \xi+\zeta.
$$

\item[(d)] each semonorm $p\in{\mathcal P}$ is an upper semicontinuous mapping $p:\Xi\to\R_+$, i.e. for any $\e>0$ the set $\{\xi\in \Xi: \ p(\xi)<\e\}$ is open; equivalently,
$$
\Big(\xi_i\overset{\varXi}{\underset{i\to\infty}{\longrightarrow}} \xi,\quad p(\xi)<\e\Big)
\quad\Longrightarrow\quad
\underbrace{p(\xi_i)<\e}_{\text{for almost all $i$}},
$$

\item[(e)]\label{uslovie-e-DEF-rassloeniya} for any point $t\in M$ and for any neighbourhood $V$ of $0_t$ in $\Xi$ there is a seminorm $p\in{\mathcal P}$, a number $\e>0$, and an open set $U$ in $M$, containing $t$, such that
$$
\{ \xi\in \pi^{-1}(U): \ p(\xi)<\e\}\subseteq V;
$$
in other words an implication holds:
$$
\Big(\pi(\xi_i)\overset{M}{\underset{i\to\infty}{\longrightarrow}} t\quad\&\quad \forall p\in{\mathcal P}\quad p(\xi_i)\underset{i\to\infty}{\longrightarrow} 0\Big)
\quad\Longrightarrow\quad
\xi_i\overset{\varXi}{\underset{i\to\infty}{\longrightarrow}} 0_t.
$$
}\eit
Then the system $(\Xi,M,\pi,\cdot,+,\{0_t;\ t\in M\},-)$ with the described topologies on $M$ and $\Xi$ and with the system of seminorms $\mathcal P$ is called a {\it locally convex bundle}.
}\eit

\brem One can relpace the condition (b) in this list by a formally weaker condition: for any $\lambda\in\C$ the multiplication $\xi\mapsto\lambda\cdot\xi$ is continuous from $\varXi$ into $\varXi$:
$$
\xi_i\overset{\varXi}{\underset{i\to\infty}{\longrightarrow}} \xi
\quad\Longrightarrow\quad \forall\lambda\in\C\quad
\lambda\cdot\xi_i\overset{\varXi}{\underset{i\to\infty}{\longrightarrow}} \lambda\cdot\xi.
$$
Indeed, if this holds, then $\lambda_i\to \lambda$ and $\xi_i\to\xi$ imply, on the one hand,
$$
\forall p\in{\mathcal P}\qquad p(\lambda_i\cdot\xi_i-\lambda\cdot\xi_i)\le \abs{\lambda_i-\lambda}\cdot p(\xi_i)\overset{\R}{\underset{i\to\infty}{\longrightarrow}} 0
$$
and, on the other,
$$
\pi(\xi_i)\to\pi(\xi)\quad\Longrightarrow\quad \pi(\lambda_i\cdot\xi_i-\lambda\cdot\xi_i)=\pi(\xi_i)\overset{M}{\underset{i\to\infty}{\longrightarrow}}\pi(\xi)
$$
This, due to (e), gives
$$
\lambda_i\cdot\xi_i-\lambda\cdot\xi_i\overset{\varXi}{\underset{i\to\infty}{\longrightarrow}} 0_{\pi(\xi)},
$$
and this, due to (c), gives
$$
\lambda_i\cdot\xi_i=\underbrace{\lambda_i\cdot\xi_i-\lambda\cdot\xi_i}_{\scriptsize\begin{matrix}\downarrow\\ 0_{\pi(\xi)} \end{matrix}}+\underbrace{\lambda\cdot\xi_i}_{\scriptsize\begin{matrix}\downarrow\\ \lambda\cdot\xi \end{matrix}}\overset{\varXi}{\underset{i\to\infty}{\longrightarrow}} 0_{\pi(\xi)}+\lambda\cdot\xi=\lambda\cdot\xi.
$$
\erem

\brem From (b) and (c) it follows that in (c) one can replace the fiberwise summing by the fiberwise subtraction:
\beq\label{xi_i->xi&zeta_i->zeta=>xi_i-zeta_i->xi-zeta}
\Big(\xi_i\overset{\varXi}{\underset{i\to\infty}{\longrightarrow}} \xi,\quad \zeta_i\overset{\varXi}{\underset{i\to\infty}{\longrightarrow}} \zeta,\quad \pi(\xi_i)=\pi(\zeta_i),\quad \pi(\xi)=\pi(\zeta)\Big)
\quad\Longrightarrow\quad
\xi_i-\zeta_i\overset{\varXi}{\underset{i\to\infty}{\longrightarrow}} \xi-\zeta.
\eeq
\erem

\bprop\label{PROP:harakt-shod-v-Xi}
In a locally convex bundle $(\Xi,M,\pi)$ the relation
$$
\xi_i\overset{\varXi}{\underset{i\to\infty}{\longrightarrow}} \xi
$$
is equivalent to the following two conditions:
\bit{

\item[(i)] $\pi(\xi_i)\overset{M}{\underset{i\to\infty}{\longrightarrow}} \pi(\xi)$,

\item[(ii)] for any seminorm $p\in{\mathcal P}$ and for any $\e>0$ there is a net $\zeta_i\in \varXi$ and an element  $\zeta\in\varXi$ such that
\beq\label{trebovanie-dlya-shod-napravl}
\zeta_i\overset{\varXi}{\underset{i\to\infty}{\longrightarrow}} \zeta
\quad
\underbrace{\pi(\zeta_i)=\pi(\xi_i)}_{\text{for almost all $i$}},
\quad
\pi(\zeta)=\pi(\xi),
\quad
p(\zeta-\xi)<\e,
\quad
\underbrace{p(\zeta_i-\xi_i)<\e}_{\text{for almost all $i$}}.
\eeq
}\eit
\eprop
 \bpr
We have to prove sufficiency. Suppose (i) and (ii) hold. Take $p\in{\mathcal P}$ and $\e>0$, and find a net $\zeta_i$ described in (ii). Let $V$ be an arbitrary neighbourhood of the point $\xi$ in $\varXi$. Since the map $\pi$ is open, the image $\pi(V)$ of the set $V$ must be a neighbourhood of the point $\pi(\xi)$ in $M$. Hence from the relation  $\pi(\xi_i)\overset{M}{\underset{i\to\infty}{\longrightarrow}} \pi(\xi)$ we have the almost all $\pi(\xi_i)$ belong to $\pi(V)$:
$$
\exists i_V:\quad \forall i\ge i_V \quad \pi(\xi_i)\in \pi(V).
$$
Hence, there are $\{\zeta^V_i;\ i\ge i_V\}$ such that
$$
\zeta^V_i\in V,\qquad \pi(\zeta^V_i)=\pi(\xi_i).
$$
We get a double net $\{\zeta^V_i;\ V\in{\mathcal U}(\xi),\ i\ge i_V\}$, where the upper index $V$ runs through the system ${\mathcal U}(\xi)$ of all neighbourhoods of the point $\xi$ in $\varXi$ (ordered by inclusion and directed by narrowing) with the following properties:
$$
\pi(\xi_i)=\pi(\zeta^V_i),\qquad \zeta^V_i\overset{\varXi}{\underset{\scriptsize\begin{matrix}i\to\infty\\ V\to \{\xi\}\end{matrix}}{\longrightarrow}}\xi
$$
Together with the conditions $\pi(\zeta)=\pi(\xi)$, $\pi(\zeta_i)=\pi(\xi_i)$ (for almost all $i$) and $\zeta_i\overset{\varXi}{\underset{i\to\infty}{\longrightarrow}} \zeta$ from \eqref{trebovanie-dlya-shod-napravl} this gives by \eqref{xi_i->xi&zeta_i->zeta=>xi_i-zeta_i->xi-zeta} the relation
$$
\zeta_i-\zeta^V_i\overset{\varXi}{\underset{\scriptsize\begin{matrix}i\to\infty\\ V\to \{\xi\}\end{matrix}}{\longrightarrow}}\zeta-\xi
$$
This, together with the inequality $p(\zeta-\xi)<\e$ from \eqref{trebovanie-dlya-shod-napravl}, gives due to (d) an inequality
$$
p(\zeta_i-\zeta^V_i)<\e
$$
(true for almost all $i$). It implies
$$
p(\xi_i-\zeta^V_i)=p(\xi_i-\zeta_i+\zeta_i-\zeta^V_i)\le \overbrace{p(\xi_i-\zeta_i)}^{\scriptsize\begin{matrix}\e \\ \phantom{\tiny{\eqref{trebovanie-dlya-shod-napravl}}}\ \text{\rotatebox{90}{$<$}}\ \tiny{\eqref{trebovanie-dlya-shod-napravl}}\end{matrix}}+p(\zeta_i-\zeta^V_i)<\e+\e=2\e.
$$
(also true for almost all $i$). Now we add an obvious relation
$$
\pi(\xi_i-\zeta^V_i)=\pi(\xi_i)\overset{M}{\underset{\scriptsize\begin{matrix}i\to\infty\\ V\to \{\xi\}\end{matrix}}{\longrightarrow}}\pi(\xi)
$$
and, by (e), we get:
$$
\xi_i-\zeta^V_i\overset{\varXi}{\underset{\scriptsize\begin{matrix}i\to\infty\\ V\to \{\xi\}\end{matrix}}{\longrightarrow}}0_{\pi(\xi)}.
$$
Applying (c), we have:
$$
\xi_i=\underbrace{\xi_i-\zeta^V_i}_{\scriptsize\begin{matrix}\downarrow\\ 0_{\pi(\xi)} \end{matrix}}+\underbrace{\zeta^V_i}_{\scriptsize\begin{matrix}\downarrow\\ \xi \end{matrix}}\overset{\varXi}{\underset{\scriptsize\begin{matrix}i\to\infty\\ V\to \{\xi\}\end{matrix}}{\longrightarrow}} 0_{\pi(\xi)}+\xi=\xi.
$$
 \epr

\paragraph{Continuous sections of a locally convex bundle.}

We consider here {\it continuous sections of locally convex bundles} $\pi:\Xi\to M$, i.e. continuous mappings
$x:M\to \Xi$ such that
$$
\pi\circ x=\id_M.
$$
The set of all continuous sections is denoted by $\Sec(\pi)$. It is endowed with a structure of a left
${\mathcal C}(M)$-module and a topology of uniform convergence on compact sets in $M$.

\vglue10pt \centerline{\bf Properties of continuous sections:} \vglue10pt
{\it
\begin{itemize}
\item[$1^\circ$.] The space $\Sec(\pi)$ of continuous sections of any locally convex bundle $\pi:\varXi\to M$ over an arbitrary paracompact locally compact space $M$ is a locally convex ${\mathcal C}(M)$-module  with respect to the fiberwise multiplication (which is jointly continuous)
    $$
    (a\cdot x)(t)=a(t)\cdot x(t),\qquad a\in {\mathcal C}(M),\ x\in \Sec(\pi)
    $$
    and seminorms
$$
p_T(x)=\sup_{t\in T}p(x(t)),\qquad p\in{\mathcal P}.
$$
where $T$ are compact sets in $M$.

\item[$2^\circ$.] \label{LM:poltnost-v-sloyah=>plotnost-v-secheniyah} Let $M$ be a paracompact locally compact space, $\pi:\varXi\to M$ a locally convex bundle, and $X$ a ${\mathcal C}(M)$-submodule in the ${\mathcal C}(M)$-module of continuous sections of $\pi$,
$$
X\subseteq\Sec(\pi)
$$
Suppose in addition that $X$ is dense in each fiber:
\beq\label{plotnost-v-sloyah}
\forall t\in M\qquad \overline{\{x(t),\ x\in X\}}=\pi^{-1}(t).
\eeq
Then $X$ is dense in $\Sec(\pi)$:
$$
\overline{X}=\Sec(\pi).
$$
\end{itemize}
} \vglue10pt

\bpr
1. If $x\in\Sec(\pi)$ and $a\in {\mathcal C}(M)$, then the fiberwise product $a\cdot x$ is a continuous section due to property (b):
$$
t_i\overset{M}{\underset{i\to\infty}{\longrightarrow}}t \quad\Longrightarrow\quad
\Big(
a(t_i)\overset{\C}{\underset{i\to\infty}{\longrightarrow}}a(t),\quad x(t_i)\overset{\varXi}{\underset{i\to\infty}{\longrightarrow}} x(t)
\Big)
\quad\Longrightarrow\quad
a(t_i)\cdot x(t_i)\overset{\varXi}{\underset{i\to\infty}{\longrightarrow}} a(t)\cdot x(t).
$$
If $T$ is a compact set in $M$, then for each continuous section $x\in\Sec(\pi)$ the image $x(T)$
is a compact space in $\Xi$. As a corollary, for every $p\in{\mathcal P}$ the image $p(x(T))$ also must be compact, if we endow $\R$ with the topology generated by the base $(-\infty,\e)$, $\e\in\R$ (then the map $p:\Xi\to\R$ becomes continuous). Thus, each covering of  $p(x(T))$ by the sets of the form $(-\infty,\e)$ contains a finite subcovering, and this means that $p(x(T))$ is bounded in $\R$ in the usual sense. In other words,
$$
p_T(x)=\sup_{t\in T}p(x(t))<\infty.
$$
Obviously, this is a seminorm on $\Sec(\pi)$, and from the chain
$$
p_T(a\cdot x)=\sup_{t\in T}p(a(t)\cdot x(t))\le \sup_{t\in
T}\Big(\abs{a(t)}\cdot p(x(t))\Big)\le \sup_{t\in T}\abs{a(t)}\cdot\sup_{t\in
T}p(x(t))
$$
it follows that these seminorms turn $\Sec(\pi)$ into a locally convex ${\mathcal C}(M)$-module (with a jointly continuous multiplication).

2. Suppose $y\in \Sec(\pi)$ and $\e>0$. From \eqref{plotnost-v-sloyah} it follows that for each seminorm $p\in{\mathcal P}$ and for each point $t\in M$ there is a continuous section $x_t\in\Sec(\pi)$ such that
$$
p(x_t(t)-y(t))<\e.
$$
From the fact that $p$ is upper semicontinuous (and the mappings $x_t$ and $y$ are continuous), it follows that the set
$$
U_t=\{s\in M: \ p(x_t(s)-y(s))<\e\}
$$
is open, hence it is a neighbourhood of $t$. We can conclude that the family $\{U_t;\ t\in M\}$ is an open covering of $M$.

Let us consider now a compact set $T\subseteq M$. The covering $\{U_t;\ t\in T\}$ contains a finite subcovering  $\{U_{t_1},...,U_{t_n}\}$ of $T$. Let us take a subordinate partition of unity
$$
0\le a_i\le 1,\qquad \supp a_i\subseteq U_{t_i},\qquad \sum_{i=1}^n a_i(t)=1\quad (t\in T)
$$
and put
$$
x=\sum_{i=1}^n a_i\cdot x_{t_i}
$$
(since $\Sec(\pi)$ is a ${\mathcal C}(M)$-module, $x\in \Sec(\pi)$). Then
\begin{multline*}
p_T(x-y)=\sup_{t\in T}p(x(t)-y(t))=\sup_{t\in T}p\l \sum_{i=1}^n a_i(t)\cdot
x_{t_i}(t)-\sum_{i=1}^n a_i(t)\cdot y(t)\r=\\= \sup_{t\in T}p\l \sum_{i=1}^n
a_i(t)\cdot \big(x_{t_i}(t)-y(t)\big)\r\le \sup_{t\in T}\sum_{i=1}^n a_i(t)\cdot
p\l \big(x_{t_i}(t)-y(t)\big)\r<\sup_{t\in T}\sum_{i=1}^n a_i(t)\cdot \e=\e.
\end{multline*}
\epr

\paragraph{Locally convex bundles generated by systems of sections and seminorms.}

\bprop\label{PROP:sushestv-topologii-v-Xi} Suppose we have
\bit{

\item[1)] a vector bundle $(\Xi,M,\pi,\cdot,+,\{0_t;\ t\in M\},-)$ over $\C$,

\item[2)] a vector space $X$ of its sections,

\item[3)] a system ${\mathcal P}$ of functions on $\varXi$,

\item[4)] a topology on the base $M$,
}\eit\noindent
and the following conditions hold:
\bit{

\item[(i)] on each fiber the restrictions $p|_{\pi^{-1}(t)}$ of functions $p\in{\mathcal P}$ form a system of seminorms, which turn $\pi^{-1}(t)$ into a (Hausdorff) locally convex space;

\item[(ii)] the system ${\mathcal P}$ is directed in ascending order: for any two functions $p,q\in{\mathcal P}$ there is a function $r\in{\mathcal P}$ that majorates $p$ and $q$:
 \beq\label{p<r&q<r}
p(\upsilon)\le r(\upsilon),\qquad q(\upsilon)\le r(\upsilon),\qquad \upsilon\in \varXi
 \eeq

\item[(iii)] for any section $x\in X$ and for any seminorm $p\in{\mathcal P}$ the function $t\in M\mapsto p(x(t))$ is upper semicontinuous on $M$,

\item[(iv)] for any point $t\in M$ the set $\{x(t);\ x\in X\}$ is dense in the locally convex space $\pi^{-1}(t)$.

}\eit
Then there us a unique topology on $\varXi$, such that the system $(\Xi,M,\pi,\cdot,+,\{0_t;\ t\in M\},-)$ with the given topology on $M$ and the system of seminorms $\mathcal P$ turns into a locally convex bundle, whose set of continuous section contains $X$:
$$
X\subseteq\Sec(\pi).
$$
Moreover, the sets
 \beq\label{baza-topologii-v-varXi}
W(x,U,p,\e)=\{\xi\in\varXi: \ \pi(\xi)\in U\ \& \ p(\xi-x(\pi(\xi)))<\e\},
 \eeq
where $x\in X$, $p\in{\mathcal P}$, $\e>0$ and $U$ is an open set in $M$, form a base of this topology in $\varXi$.
\eprop

\bpr 1. Let us show at the beginning, that the sets \eqref{baza-topologii-v-varXi} indeed form a base of some topology in $\varXi$. First, they cover $\varXi$, since if $\xi\in\varXi$, then (iii) implies that for any $\e>0$ and $p\in{\mathcal P}$ there is $x\in X$ such that
$$
p(\xi-x(\pi(\xi)))<\e,
$$
and if we choose now an open neighbourhood $U$ of $\pi(\xi)$, then $\xi$ lies in the set $W(x,U,p,\e)$.

Let us check the second axiom of base: consider a point $\xi$, and arbitrary base neighbourhoods $W(x,U,p,\e)$ and $W(y,V,q,\delta)$ of $\xi$,
 \beq\label{sushestv-topologii-v-Xi-8}
\xi\in W(x,U,p,\e)\cap W(y,V,q,\delta).
 \eeq
We have to show that there is a base neighbourhood $W(z,O,r,\sigma)$ of $\xi$ such that
 \beq\label{sushestv-topologii-v-Xi-7}
\xi\in W(z,O,r,\sigma)\subseteq W(x,U,p,\e)\cap W(y,V,q,\delta).
 \eeq
The inclusion \eqref{sushestv-topologii-v-Xi-8} means that
 \beq\label{sushestv-topologii-v-Xi-6}
\pi(\xi)\in U,\qquad p(\xi-x(\pi(\xi)))<\e,\qquad  \pi(\xi)\in V,\qquad
q(\xi-y(\pi(\xi)))<\delta.
 \eeq
Consider the fiber $\pi^{-1}(\pi(\xi))$. The conditions
 \beq\label{sushestv-topologii-v-Xi-0}
p(\xi-x(\pi(\xi)))<\e,\qquad  q(\xi-y(\pi(\xi)))<\delta
 \eeq
can be understood so that the point $\xi$ lies in the intersections of the neighbourhoods of the points $x(\pi(\xi))$ and $y(\pi(\xi))$, defined by seminorms $p$ and $q$ with the radii $\e$ and $\delta$. Hence (by (ii)) there exists a seminorm $r\in{\mathcal P}$ and a number $\sigma>0$ such that the $r$-neighbourhood of $\xi$ of radius $2\sigma$ is contained in those $p$- and $q$-neighbourhoods:
 \beq\label{sushestv-topologii-v-Xi-1}
\forall \zeta\in\pi^{-1}(\pi(\xi))\qquad r(\zeta-\xi)<2\sigma\quad\Longrightarrow\quad
\Big(p(\zeta-x(\pi(\xi)))<\e\quad\&\quad  q(\zeta-y(\pi(\xi)))<\delta\Big).
 \eeq
Let us reduce $\sigma$, if necessary, so that \eqref{sushestv-topologii-v-Xi-0} can be replaced by
 \beq\label{sushestv-topologii-v-Xi-4}
p(\xi-x(\pi(\xi)))<\e-2\sigma,\qquad  q(\xi-y(\pi(\xi)))<\delta-2\sigma.
 \eeq
Then let us use (iv) and take $z\in X$ such that
 \beq\label{sushestv-topologii-v-Xi-5}
r(z(\pi(\xi))-\xi)<\sigma.
 \eeq
We have a chain of implications:
 \begin{multline*}
r(\zeta-z(\pi(\xi)))<\sigma\quad\Longrightarrow\quad
r(\zeta-\xi)\le r(\zeta-z(\pi(\xi)))+r(z(\pi(\xi))-\xi)<\sigma+\sigma=2\sigma
\quad\Longrightarrow
\\ \overset{\eqref{sushestv-topologii-v-Xi-1}}{\Longrightarrow}\quad
\Big(p(\zeta-x(\pi(\xi)))<\e\quad\&\quad  q(\zeta-y(\pi(\xi)))<\delta\Big).
 \end{multline*}
In other words, in the fiber $\pi^{-1}(\pi(\xi))$ the $r$-neighbourhood of radius $\sigma$ of the point $z(\pi(\xi))$ is also contained in those $p$- and $q$-neighbourhoods:
 \beq\label{sushestv-topologii-v-Xi-2}
\forall \zeta\in\pi^{-1}(\pi(\xi))\qquad r(\zeta-z(\pi(\xi)))<\sigma\quad\Longrightarrow\quad
\Big(p(\zeta-x(\pi(\xi)))<\e\quad\&\quad  q(\zeta-y(\pi(\xi)))<\delta\Big).
 \eeq
By (ii), we can think that the seminorm $r$ majorates $p$ and $q$, i.e. \eqref{p<r&q<r} holds.
Take this and consider the set
 \beq\label{sushestv-topologii-v-Xi-3}
O=\left\{s\in U\cap V:\quad p(z(s)-x(s))<\e-\sigma \quad \&\quad q(z(s)-y(s))<\delta-\sigma \right\}.
 \eeq
By (iii), it is open in $M$. At the same time it contains $\pi(\xi)$, since, first, $\pi(\xi)\in U\cap V$ by  \eqref{sushestv-topologii-v-Xi-6}, second,
$$
p(z(\pi(\xi))-x(\pi(\xi)))\le \underbrace{p(z(\pi(\xi))-\xi)}_{\scriptsize\begin{matrix}
\phantom{\tiny{\eqref{p<r&q<r}}}\quad\text{\rotatebox{90}{$\ge$}}\quad\tiny{\eqref{p<r&q<r}}  \\
r(z(\pi(\xi))-\xi) \\
\phantom{\tiny{\eqref{sushestv-topologii-v-Xi-5}}}\quad\text{\rotatebox{90}{$>$}}\quad\tiny{\eqref{sushestv-topologii-v-Xi-5}} \\
\sigma
\end{matrix}}+\underbrace{p(\xi-x(\pi(\xi)))}_{\scriptsize\begin{matrix}
\phantom{\tiny{\eqref{sushestv-topologii-v-Xi-5}}}\quad\text{\rotatebox{90}{$>$}}\quad\tiny{\eqref{sushestv-topologii-v-Xi-4}} \\
\e-2\sigma
\end{matrix}}<\e-\sigma,
$$
and, third,
$$
p(z(\pi(\xi))-y(\pi(\xi)))\le \underbrace{p(z(\pi(\xi))-\xi)}_{\scriptsize\begin{matrix}
\phantom{\tiny{\eqref{p<r&q<r}}}\quad\text{\rotatebox{90}{$\ge$}}\quad\tiny{\eqref{p<r&q<r}}  \\
r(z(\pi(\xi))-\xi) \\
\phantom{\tiny{\eqref{sushestv-topologii-v-Xi-5}}}\quad\text{\rotatebox{90}{$>$}}\quad\tiny{\eqref{sushestv-topologii-v-Xi-5}} \\
\sigma
\end{matrix}}+\underbrace{p(\xi-y(\pi(\xi)))}_{\scriptsize\begin{matrix}
\phantom{\tiny{\eqref{sushestv-topologii-v-Xi-5}}}\quad\text{\rotatebox{90}{$>$}}\quad\tiny{\eqref{sushestv-topologii-v-Xi-4}} \\
\delta-2\sigma
\end{matrix}}<\delta-\sigma.
$$
The inclusion $\pi(\xi)\in O$ together with \eqref{sushestv-topologii-v-Xi-5} imply
 \beq\label{sushestv-topologii-v-Xi-9}
\xi\in W(z,O,r,\sigma).
 \eeq
Further, for any $\zeta\in\varXi$ we have
$$
\pi(\zeta)\in O\quad\&\quad r(\zeta-z(\pi(\zeta)))<\sigma
\quad\Longrightarrow\quad
p(\zeta-x(\pi(\zeta)))\le \underbrace{p(\zeta-z(\pi(\zeta)))}_{\scriptsize\begin{matrix}
\phantom{\tiny{\eqref{p<r&q<r}}}\quad\text{\rotatebox{90}{$\ge$}}\quad\tiny{\eqref{p<r&q<r}}  \\
r(\zeta-z(\pi(\zeta))) \\
\text{\rotatebox{90}{$>$}} \\
\sigma
\end{matrix}}
+\underbrace{p(z(\pi(\zeta))-x(\pi(\zeta)))}_{\scriptsize\begin{matrix}
\phantom{\tiny{\eqref{sushestv-topologii-v-Xi-3}}}\quad\text{\rotatebox{90}{$>$}}\quad\tiny{\eqref{sushestv-topologii-v-Xi-3}} \\
\e-\sigma
\end{matrix}}<\e
$$
hence
$$
W(z,O,r,\sigma)\subseteq W(x,U,p,\e).
$$
By analogy,
$$
W(z,O,r,\sigma)\subseteq W(y,V,q,\delta).
$$
Together with \eqref{sushestv-topologii-v-Xi-9} this gives \eqref{sushestv-topologii-v-Xi-7}.

2. Thus, we understood that the sets \eqref{baza-topologii-v-varXi} form a local base of some topology in $\Xi$. Note then that with respect to this topology the mapping $\pi:\Xi\to M$ is continuous and open. Suppose
$$
\xi_i\overset{\Xi}{\underset{i\to\infty}{\longrightarrow}}\xi.
$$
Consider an arbitrary neighbourhood $U$ of the point $\pi(\xi)$. Take a seminorm $p\in{\mathcal P}$, and, using (iv), find a section $x\in X$ such that $p(\xi-x(\pi(\xi)))<1$. Then the set $W(x,U,p,1)$ is a neighbouhood of $\xi$, hence $\xi\in W(x,U,p,1)$ for almost all $i$, and this implies $\pi(\xi_i)\in U$ for almost all $i$. This proves the relation
$$
\pi(\xi_i)\overset{M}{\underset{i\to\infty}{\longrightarrow}}\pi(\xi),
$$
i.e. the continuity of $\pi$.

Now let us take an arbitrary base neighbourhood $W(x,U,p,\e)$. The projection $\pi$ maps it onto an open set $U$, so each point $t\in U$ has an inverse image $x(t)\in W(x,U,p,\e)$. This means that $\pi$ is an open mapping.

3. Let us show then that with respect to this topology in $\Xi$ each section $x\in X$ is a continuous mapping. Suppose
$$
t_i\overset{M}{\underset{i\to\infty}{\longrightarrow}}t.
$$
Consider a base neighbourhood $W(y,U,p,\e)$ of $x(t)$. The condition $x(t)\in W(y,U,p,\e)$ means that
 \beq\label{sushestv-topologii-v-Xi-10}
t\in U,\qquad p(x(t)-y(t))<\e.
 \eeq
Put
$$
V=\{s\in U:\quad p(x(s)-y(s))<\e\}.
$$
By (iii), this set is open in $M$. By \eqref{sushestv-topologii-v-Xi-10} it contains $t$, and thus it is a neighbourhood of $t$. Therefore for almost all indices $i$ have
$$
t_i\in V\quad\Longrightarrow\quad p(x(t_i)-y(t_i))<\e\quad\Longrightarrow\quad x(t_i)\in W(y,V,p,\e)\subseteq W(y,U,p,\e).
$$
This proves the relation
 \beq\label{sushestv-topologii-v-Xi-12}
x(t_i)\overset{\Xi}{\underset{i\to\infty}{\longrightarrow}}x(t),
 \eeq
i.e. the continuity of $x$.

4. Now we start to check that with the described topology on $\Xi$ the triple $(\Xi,M,\pi)$ is a locally convex bundle. Let us first verify the continuity of the fiberwise multiplication by scalars. Suppose
$$
\lambda_i\overset{\C}{\underset{i\to\infty}{\longrightarrow}}\lambda,\quad
\xi_i\overset{\varXi}{\underset{i\to\infty}{\longrightarrow}} \xi.
$$
Consider first the case of $\lambda\ne 0$. Take any neighbourhood $W(x,U,p,\e)$ of the point $\lambda\cdot\xi$:
$$
\pi(\lambda\cdot\xi)\in U,\qquad p(\lambda\cdot\xi-x(\pi(\xi)))<\e.
$$
Then
$$
\pi(\xi)\in U,\qquad p\l \xi-\frac{1}{\lambda}\cdot x(\pi(\xi))\r <\frac{\e}{\abs{\lambda}}.
$$
i.e. the set $W\l \frac{1}{\lambda}\cdot x,U,p,\frac{\e}{\abs{\lambda}}\r$ is a neighbourhood for $\xi$, and thus  $\xi_i\in W\l \frac{1}{\lambda}\cdot x,U,p,\frac{\e}{\abs{\lambda}}\r$ for almost all $i$:
 \beq\label{sushestv-topologii-v-Xi-11}
\exists i_0\qquad\forall i\ge i_0\qquad
\pi(\xi_i)\in U,\qquad p\l \xi_i-\frac{1}{\lambda}\cdot x(\pi(\xi_i))\r <\frac{\e}{\abs{\lambda}}.
 \eeq
Now for almost all $i$ we have
\begin{multline*}
p(\lambda_i\cdot\xi_i-x(\pi(\lambda_i\cdot\xi_i)))=p(\lambda_i\cdot\xi_i-x(\pi(\xi_i)))=
\abs{\lambda_i}\cdot p\l \xi_i-\frac{1}{\lambda_i}\cdot x(\pi(\xi_i))\r\le \\ \le
\abs{\lambda_i}\cdot p\l \xi_i-\frac{1}{\lambda}\cdot x(\pi(\xi))\r+\abs{\lambda_i}\cdot p\l \frac{1}{\lambda}\cdot x(\pi(\xi))-\frac{1}{\lambda_i}\cdot x(\pi(\xi_i))\r=\\=
\abs{\lambda_i}\cdot p\l \xi_i-\frac{1}{\lambda}\cdot x(\pi(\xi))\r+p\l \frac{\lambda_i}{\lambda}\cdot x(\pi(\xi))-x(\pi(\xi_i))\r=
\\=
\abs{\lambda_i}\cdot p\l \xi_i-\frac{1}{\lambda}\cdot x(\pi(\xi))\r+p\l \frac{\lambda_i}{\lambda}\cdot x(\pi(\xi))-x(\pi(\xi))\r
+p\Big( x(\pi(\xi))-x(\pi(\xi_i))\Big)=
\\=
\underbrace{\abs{\lambda_i}}_{\scriptsize\begin{matrix}
\text{\rotatebox{90}{$>$}} \\
2\abs{\lambda}
\end{matrix}}\cdot \underbrace{p\l \xi_i-\frac{1}{\lambda}\cdot x(\pi(\xi))\r}_{\scriptsize\begin{matrix}
\phantom{\tiny{\eqref{sushestv-topologii-v-Xi-11}}}\quad\text{\rotatebox{90}{$>$}}\quad\tiny{\eqref{sushestv-topologii-v-Xi-11}} \\
\frac{\e}{\abs{\lambda}}
\end{matrix}}
+\underbrace{\abs{\frac{\lambda_i}{\lambda}-1}}_{\scriptsize\begin{matrix}
\text{\rotatebox{90}{$>$}} \\
\e
\end{matrix}}\cdot p\Big(x(\pi(\xi))\Big)
+\underbrace{p\Big( x(\pi(\xi))-x(\pi(\xi_i))\Big)}_{\scriptsize\begin{matrix}
\phantom{\tiny{\eqref{sushestv-topologii-v-Xi-12}}}\quad\text{\rotatebox{90}{$>$}}\quad\tiny{\eqref{sushestv-topologii-v-Xi-12}} \\
\e
\end{matrix}}<3\e+\e\cdot p\big(x(\pi(\xi))\big)
\end{multline*}
If we add the condition $\pi(\xi_i)\in U$ from \eqref{sushestv-topologii-v-Xi-11}, we obtain
$$
\lambda_i\cdot\xi_i\in W(x,U,p,3\e+\e\cdot p(x(\pi(\xi))))
$$
for almost all $i$. Since initially $\e$ was arbitrary,
 \beq\label{lambda_i-xi_i->lambda-xi}
\lambda_i\cdot\xi_i\overset{\varXi}{\underset{i\to\infty}{\longrightarrow}}
\lambda\cdot\xi.
 \eeq

Now consider the case of $\lambda=0$. Take a neighbourhood $W(x,U,p,\e)$ of $0_{\pi(\xi)}$:
$$
\pi(\xi)\in U,\qquad p(x(\pi(\xi)))<\e.
$$
Find $\delta>0$ such that
$$
p(x(\pi(\xi)))<\e-\delta.
$$
Since $\xi_i\to\xi$, for almost all $i$ we have
$$
\pi(\xi_i)\in U,\qquad p(x(\pi(\xi_i)))<\e-\delta.
$$
Now if $p(\xi)\ne 0$, then for almost all $i$
$$
p(\lambda_i\cdot\xi_i-x(\pi(\xi_i)))\le \underbrace{\abs{\lambda_i}}_{\scriptsize\begin{matrix}
\text{\rotatebox{90}{$>$}}\\ \frac{\delta}{2p(\xi)} \end{matrix}}\cdot \underbrace{p(\xi_i)}_{\scriptsize\begin{matrix}
\text{\rotatebox{90}{$>$}} \\ 2p(\xi)\end{matrix}} +\underbrace{p(x(\pi(\xi_i)))}_{\scriptsize\begin{matrix}
\text{\rotatebox{90}{$>$}} \\ \e-\delta\end{matrix}}<\e.
$$
If $p(\xi)=0$, then $\xi_i\to\xi$ implies $p(\xi_i)<1$ for almost all $i$, hence
$$
p(\lambda_i\cdot\xi_i-x(\pi(\xi_i)))\le \underbrace{\abs{\lambda_i}}_{\scriptsize\begin{matrix}
\text{\rotatebox{90}{$>$}}\\ \delta \end{matrix}}\cdot \underbrace{p(\xi_i)}_{\scriptsize\begin{matrix}
\text{\rotatebox{90}{$>$}} \\ 1\end{matrix}} +\underbrace{p(x(\pi(\xi_i)))}_{\scriptsize\begin{matrix}
\text{\rotatebox{90}{$>$}} \\ \e-\delta\end{matrix}}<\e.
$$
In any case for almost all $i$ we have
$$
\pi(\xi_i)\in U,\qquad p(\lambda_i\cdot\xi_i-x(\pi(\xi_i)))<\e.
$$
i.e. $\lambda_i\cdot\xi_i\in W(x,U,p,\e)$. This also means that \eqref{lambda_i-xi_i->lambda-xi}.

5. Let us prove the continuity of the fiberwise summing. Suppose
$$
\xi_i\overset{\varXi}{\underset{i\to\infty}{\longrightarrow}} \xi,\quad \upsilon_i\overset{\varXi}{\underset{i\to\infty}{\longrightarrow}} \upsilon,\quad \pi(\xi_i)=\pi(\upsilon_i),\quad \pi(\xi)=\pi(\upsilon)
$$
Take a base neighbourhood $W(z,U,p,\e)$ of $\xi+\upsilon$:
$$
\pi(\xi+\upsilon)\in U,\qquad p(\xi+\upsilon-z(\pi(\xi+\upsilon)))<\e.
$$
Find a number $\sigma>0$ such that
\beq\label{sushestv-topologii-v-Xi-18}
p(\xi+\upsilon-z(\pi(\xi+\upsilon)))<\e-2\sigma.
\eeq
In the fiber $\pi^{-1}(\pi(\xi+\upsilon))=\pi^{-1}(\pi(\xi))=\pi^{-1}(\pi(\upsilon))$ the operation of continuity is continuous, hence there are base neighbourhoods $W(x,V_x,q,\delta)$ and $W(y,V_y,r,\delta)$ of $\xi$ and $\upsilon$ respectively, such that
$$
\forall \xi'\in W(x,V_x,q,\delta)\cap \pi^{-1}(\pi(\xi))\quad \forall \upsilon'\in W(y,V_y,r,\delta)\cap \pi^{-1}(\pi(\upsilon))\quad
\xi'+\upsilon'\in W(z,U,p,\e-2\sigma).
$$
Thus,
\begin{multline} \label{sushestv-topologii-v-Xi-13}
\forall \xi',\upsilon'\in \pi^{-1}(\pi(\xi+\upsilon))\quad
q(\xi'-x(\pi(\xi)))<\delta\quad\&\quad r(\upsilon'-y(\pi(\upsilon)))<\delta \quad\Longrightarrow \\ \Longrightarrow\quad
p(\xi'+\upsilon'-z(\pi(\xi+\upsilon)))<\e-2\sigma.
\end{multline}
At the same time from the inclusions $\xi\in W(x,V_x,q,\delta)$ and $\upsilon\in W(y,V_y,r,\delta)$ it follows, first, that
 \beq\label{sushestv-topologii-v-Xi-17}
q(\xi-x(\pi(\xi)))<\delta\quad\&\quad r(\upsilon-y(\pi(\upsilon)))<\delta,
 \eeq
and, second, that
$$
\pi(\xi)\in V_x,\qquad \pi(\upsilon)\in V_y,
$$
That is $V_x$ and $V_y$ are neighbourhoods for $\pi(\xi)$ and $\pi(\upsilon)$ respectively. We can narrow the neighbourhoods $V_x$ and $V_y$ in such a way that they will coincide and lie in $U$:
$$
V_x=V_y=V\subseteq U.
$$

At the same time one can obviously choose the seminorms $q$ and $r$ such that they majorate $p$ (everywhere on $\varXi$)
\beq\label{sushestv-topologii-v-Xi-14}
p\le q, \qquad p\le r,\qquad
\eeq
and the number $\delta$ such that
\beq\label{sushestv-topologii-v-Xi-15}
\delta<\frac{\sigma}{2}.
\eeq

Let us now consider the set
 \beq\label{sushestv-topologii-v-Xi-16}
O=\left\{s\in V:\quad p(x(s)+y(s)-z(s))<\e-\sigma \right\}.
 \eeq
(it is open in $M$ by (iii)).

The open set $W(x,O,q,\delta)$ is a neighbourhood of $\xi$, since, first,
 $$
p(x(\pi(\xi))+y(\underbrace{\pi(\xi)}_{\scriptsize\begin{matrix}
\|
 \\
\pi(\upsilon)
 \end{matrix}})-z(\underbrace{\pi(\xi)}_{\scriptsize\begin{matrix}
\|
 \\
\pi(\xi+\upsilon)
 \end{matrix}})\le \underbrace{p(x(\pi(\xi))-\xi)}_{\scriptsize\begin{matrix}
\phantom{\tiny{\eqref{sushestv-topologii-v-Xi-14}}}\quad\text{\rotatebox{90}{$\ge$}}\quad\tiny{\eqref{sushestv-topologii-v-Xi-14}}
 \\
q(x(\pi(\xi))-\xi)
 \\
\phantom{\tiny{\eqref{sushestv-topologii-v-Xi-17}}}\quad\text{\rotatebox{90}{$>$}}\quad\tiny{\eqref{sushestv-topologii-v-Xi-17}}
 \\
 \delta
 \\
\phantom{\tiny{\eqref{sushestv-topologii-v-Xi-15}}}\quad\text{\rotatebox{90}{$>$}}\quad\tiny{\eqref{sushestv-topologii-v-Xi-15}}
 \\
\frac{\sigma}{2}
 \end{matrix}}
+\underbrace{p(y(\pi(\upsilon))-\upsilon)}_{\scriptsize\begin{matrix}
\phantom{\tiny{\eqref{sushestv-topologii-v-Xi-14}}}\quad\text{\rotatebox{90}{$\ge$}}\quad\tiny{\eqref{sushestv-topologii-v-Xi-14}}
 \\
r(x(\pi(\upsilon))-\upsilon)
 \\
\phantom{\tiny{\eqref{sushestv-topologii-v-Xi-17}}}\quad\text{\rotatebox{90}{$>$}}\quad\tiny{\eqref{sushestv-topologii-v-Xi-17}}
 \\
\delta \\
\phantom{\tiny{\eqref{sushestv-topologii-v-Xi-15}}}\quad\text{\rotatebox{90}{$>$}}\quad\tiny{\eqref{sushestv-topologii-v-Xi-15}}
 \\
\frac{\sigma}{2}
\end{matrix}}
+\underbrace{p(\xi+\upsilon-z(\pi(\xi+\upsilon)))}_{\scriptsize\begin{matrix}
\phantom{\tiny{\eqref{sushestv-topologii-v-Xi-18}}}\quad\text{\rotatebox{90}{$>$}}\quad\tiny{\eqref{sushestv-topologii-v-Xi-18}}
 \\
\e-2\sigma
 \end{matrix}} <\e-\sigma
 $$
and therefore, $\pi(\xi)\in O$. And, second, by \eqref{sushestv-topologii-v-Xi-17}, $q(\xi-x(\pi(\xi)))<\delta$.

Similarly, $\upsilon\in W(y,O,r,\delta)$.

Let us show now that the neighbourhoods $W(x,O,q,\delta)$ and $W(y,O,r,\delta)$ of $\xi$ and $\upsilon$ satisfy the condition
 \beq\label{sushestv-topologii-v-Xi-19}
\forall \xi'\in W(x,O,q,\delta)\quad \forall \upsilon'\in W(y,O,r,\delta)\quad \Big( \pi(\xi')=\pi(\upsilon')\quad\Longrightarrow\quad
\xi'+\upsilon'\in W(z,U,p,\e)\Big).
 \eeq
Indeed, $\xi'\in W(x,O,q,\delta)$, $\upsilon'\in W(y,O,r,\delta)$ and $\pi(\xi')=\pi(\upsilon')$ imply, first, that
$$
\pi(\xi')=\pi(\upsilon')\in O\subseteq V\subseteq U.
$$
And, second, that
$$
p(\xi'+\upsilon'-z(\pi(\xi'+\upsilon')))\le
\underbrace{p(\xi'-x(\pi(\xi')))}_{\scriptsize\begin{matrix}
\phantom{\tiny{\eqref{sushestv-topologii-v-Xi-14}}}\quad\text{\rotatebox{90}{$\ge$}}\quad\tiny{\eqref{sushestv-topologii-v-Xi-14}}
 \\
q(\xi'-x(\pi(\xi')))
 \\
\text{\rotatebox{90}{$>$}}
 \\
 \delta
 \\
\phantom{\tiny{\eqref{sushestv-topologii-v-Xi-15}}}\quad\text{\rotatebox{90}{$>$}}\quad\tiny{\eqref{sushestv-topologii-v-Xi-15}}
 \\
\frac{\sigma}{2}
 \end{matrix}}
+\underbrace{p(\upsilon'-y(\pi(\upsilon')))}_{\scriptsize\begin{matrix}
\phantom{\tiny{\eqref{sushestv-topologii-v-Xi-14}}}\quad\text{\rotatebox{90}{$\ge$}}\quad\tiny{\eqref{sushestv-topologii-v-Xi-14}}
 \\
r(\upsilon'-y(\pi(\upsilon')))
 \\
\text{\rotatebox{90}{$>$}}
 \\
\delta \\
\phantom{\tiny{\eqref{sushestv-topologii-v-Xi-15}}}\quad\text{\rotatebox{90}{$>$}}\quad\tiny{\eqref{sushestv-topologii-v-Xi-15}}
 \\
\frac{\sigma}{2}
\end{matrix}}
+\underbrace{p(x(\xi')+y(\upsilon')-z(\pi(\xi'+\upsilon')))}_{\scriptsize\begin{matrix}
\phantom{\tiny{\eqref{sushestv-topologii-v-Xi-16}}}\quad\text{\rotatebox{90}{$>$}}\quad\tiny{\eqref{sushestv-topologii-v-Xi-16}}
 \\
\e-\sigma
 \end{matrix}} <\e
$$
We now obtain a chain
$$
\xi_i\overset{\varXi}{\underset{i\to\infty}{\longrightarrow}} \xi,\quad \upsilon_i\overset{\varXi}{\underset{i\to\infty}{\longrightarrow}} \upsilon,\quad \pi(\xi_i)=\pi(\upsilon_i),\quad \pi(\xi)=\pi(\upsilon)
$$
$$
\Downarrow
$$
$$
\xi_i\in W(x,O,q,\delta),\quad \upsilon_i\in W(y,O,r,\delta)\qquad \text{for almost all $i$}
$$
$$
\phantom{\scriptsize\text{\eqref{sushestv-topologii-v-Xi-19}}}\quad\Downarrow\quad{\scriptsize\text{\eqref{sushestv-topologii-v-Xi-19}}}
$$
$$
\xi_i+\upsilon_i\in W(z,U,p,\e)\qquad \text{for almost all $i$}
$$
Here $W(z,U,p,\e)$ was an arbitrary base neighbourhood of the point $\xi+\upsilon$. Hence
$$
\xi_i+\upsilon_i\overset{\varXi}{\underset{i\to\infty}{\longrightarrow}} \xi+\upsilon.
$$

6. Let $p$ be an arbitrary neighbourhood from $\mathcal P$ and $\e>0$. Let us show that the set $W=\{\xi\in\Xi: \ p(\xi)<\e\}$ is open in $\Xi$. Take a point $\xi\in W$. The condition $p(\xi)<\e$ imply that there is a $\sigma>0$ such that $p(\xi)<\e-2\sigma$. By (iv), there is $x\in X$ such that $p(\xi-x(\pi(\xi)))<\sigma$. Put
$$
O=\{t\in M:\ p(x(t))<\e-\sigma\}.
$$
Then the base neighboruhood $W(x,O,p,\sigma)$ contains the point $\xi$, since, first,
$$
p(x(\pi(\xi)))\le \underbrace{p(x(\pi(\xi))-\xi)}_{\scriptsize\begin{matrix}
\text{\rotatebox{90}{$>$}} \\
\sigma
 \end{matrix}}+\underbrace{p(\xi)}_{\scriptsize\begin{matrix}
\text{\rotatebox{90}{$>$}}
 \\
\e-2\sigma
 \end{matrix}}<\e-\sigma
 \quad\Longrightarrow\quad \pi(\xi)\in O,
$$
and, second, by the choice of $x$, we have $p(\xi-x(\pi(\xi)))<\sigma$. That is, $W(x,O,p,\sigma)$ is a neighbourhood of the point $\xi$.

On the other hand, $W(x,O,p,\sigma)$ is contained in the set $W$, since if $\upsilon\in W(x,O,p,\sigma)$, then
$$
p(\upsilon)\le \underbrace{p(\upsilon-x(\pi(\upsilon)))}_{\scriptsize\begin{matrix}
\text{\rotatebox{90}{$>$}} \\
\sigma
 \end{matrix}}+\underbrace{x(\pi(\upsilon))}_{\scriptsize\begin{matrix}
\text{\rotatebox{90}{$>$}}
 \\
\e-\sigma
 \end{matrix}}<\e.
$$

7. let us show that the condition (e) on page \pageref{uslovie-e-DEF-rassloeniya} holds:
for any point $t\in M$ and for any neighbourhood $V$ of the point $0_t$ in $\Xi$ there is a seminorm $p\in{\mathcal P}$, a number $\sigma>0$ and an open set $O$ in $M$, such that $t\in O$, and
$$
\{ \xi\in \pi^{-1}(O): \ p(\xi)<\sigma\}\subseteq V.
$$
Since the topology in $\Xi$ is generated by the neighbourhoods \eqref{baza-topologii-v-varXi}, there is a base neighbourhood $W(x,U,p,\e)$ of $0_t$, such that $W(x,U,p,\e)\subseteq V$, and
$$
0_t\in W(x,U,p,\e)\subseteq V.
$$
This means the following two conditions:
$$
t\in U\quad\&\quad p(0_t-x(t))<\e.
$$
Take $\sigma>0$ such that
$$
p(0_t-x(t))<\e-\sigma
$$
and put
$$
O=\{s\in U:\ p(x(s))<\e-\sigma\}.
$$
By (iii), this is an open set in $M$. It contains $t$, since
$$
p(x(t))\le \underbrace{p(x(t)-0_t)}_{\scriptsize\begin{matrix}
\text{\rotatebox{90}{$>$}}
 \\
\e-\sigma
 \end{matrix}}+\underbrace{p(0_t)}_{\scriptsize\begin{matrix}
\|
 \\
0
 \end{matrix}}<\e-\sigma
$$
Note that if $\xi\in \pi^{-1}(O)$ and $p(\xi)<\sigma$, then
$$
p(\xi-x(\pi(\xi)))\le \underbrace{p(\xi)}_{\scriptsize\begin{matrix}
\text{\rotatebox{90}{$>$}}
 \\
\sigma
 \end{matrix}}+\kern-25pt\underbrace{p(x(\pi(\xi)))}_{\scriptsize\begin{matrix}
\phantom{\tiny{(\pi(\xi)\in O)}}\quad\text{\rotatebox{90}{$>$}}\quad\tiny{(\pi(\xi)\in O)}
 \\
\e-\sigma
 \end{matrix}}\kern-25pt<\e
$$
Now we have the following chain:
$$
\{ \xi\in \pi^{-1}(O): \ p(\xi)<\sigma\}\subseteq W(x,O,p,\e)\subseteq W(x,U,p,\e)\subseteq V.
$$

8. Let us prove the uniqueness of this topology. It follows from the fact that the convergence of a net in it
\beq\label{PROOF:edinstvennost-topologii-na-rassloenii}
\xi_i\overset{\varXi}{\underset{i\to\infty}{\longrightarrow}} \xi
\eeq
is uniquely defined by the behaviour of $\pi$, ${\mathcal P}$ and $X$ in the points
$\xi_i$ and $\xi$, namely, by the following two conditions:

\bit{

\item[(a)] $\pi(\xi_i)\overset{M}{\underset{i\to\infty}{\longrightarrow}} \pi(\xi)$,

\item[(b)] for any section $x\in X$, any seminorm $p\in{\mathcal P}$ and any $\e>0$ the condition
$$
p(\xi_i-x(\pi(\xi_i)))<p(\xi-x(\pi(\xi)))+\e
$$
holds for almost all indices $i$.
}\eit
Indeed, if \eqref{PROOF:edinstvennost-topologii-na-rassloenii} is true, then the condition (a) will follow from the continuity of the mapping $\pi:\varXi\to M$. And the condition (b) is proved as follows: from  $\xi_i\overset{\varXi}{\underset{i\to\infty}{\longrightarrow}} \xi$ it follows that  $x(\pi(\xi_i))\overset{\varXi}{\underset{i\to\infty}{\longrightarrow}} x(\pi(\xi))$, and together, due to the upper semicontinuity of $p$, this gives a chain of inequalities, which are true for almost all $i$:
$$
p(\xi_i-x(\pi(\xi_i)))\le \kern-10pt\underbrace{p(\xi_i-\xi)}_{\scriptsize\begin{matrix}\text{\rotatebox{90}{$>$}}\\ \frac{\e}{2}\\ \text{for almost all $i$}\end{matrix}} \kern-10pt+p(\xi-x(\pi(\xi)))+ \underbrace{p(x(\pi(\xi))-x(\pi(\xi_i)))}_{\scriptsize\begin{matrix}\text{\rotatebox{90}{$>$}}\\ \frac{\e}{2}\\ \text{for almost all $i$}\end{matrix}}<\frac{\e}{2}+p(\xi-x(\pi(\xi)))+\frac{\e}{2}.
$$
On the contrary, suppose (a) and (b) are true. Then by (iii), for each seminorm $p\in{\mathcal P}$ and for any $\e>0$ one can find a section $x\in X$ such that
\beq\label{PROOF:edinstvennost-topologii-na-rassloenii-1}
p(\xi-x(\pi(\xi)))<\e.
\eeq
After that we have, first,
$$
x(\pi(\xi_i))\overset{\varXi}{\underset{i\to\infty}{\longrightarrow}} x(\pi(\xi))
$$
(by the condition (a) and the continuity of the map $x$). Second, for all $i$
$$
\pi(x(\pi(\xi)))=\pi(\xi),\qquad \pi(x(\pi(\xi_i)))=\pi(\xi_i)
$$
(since $x$ is a section of $\pi$). And, third, for almost all $i$
$$
p(\xi_i-x(\pi(\xi_i)))\le \underbrace{p(\xi-x(\pi(\xi)))}_{\scriptsize\begin{matrix}\phantom{\tiny{\eqref{PROOF:edinstvennost-topologii-na-rassloenii-1}}}\ \text{\rotatebox{90}{$>$}}\ \tiny{\eqref{PROOF:edinstvennost-topologii-na-rassloenii-1}}\\ \e\end{matrix}}+\e<2\e.
$$
(by (b)). Together this means that the points $\zeta_i=x(\pi(\xi_i))$, $\zeta=x(\pi(\xi))$ satisfy the condition  \eqref{trebovanie-dlya-shod-napravl} (where $\e$ is replced by $2\e$), and since $p\in{\mathcal P}$ and $\e>0$ were arbitrary, by Proposition \ref{PROP:harakt-shod-v-Xi} this implies \eqref{PROOF:edinstvennost-topologii-na-rassloenii}.
\epr

\paragraph{Morphisms of bundles.}

Suppose we have two locally convex bundles $\pi:\varXi\to M$, $\rho:\varOmega\to N$, and two continuous maps  $\mu:\varXi\to\varOmega$ and $\sigma:M\to N$ such that the following diagram is commutative:
\beq\label{predmorfizm-rassloenij}
 \xymatrix @R=2pc @C=3pc
 {
 \varXi\ar[r]^{\mu}\ar[d]_{\pi} & \varOmega\ar[d]^{\rho}\\
   M\ar[r]^{\sigma} & N
 }
\eeq
Then for any $t\in M$ and any $\xi\in\pi^{-1}(t)$ we have
$$
\rho\big(\mu(\xi)\big)=\sigma(\pi(\xi))=\sigma(t),
$$
and this means that $\mu$ maps the fiber $\pi^{-1}(t)$ into the fiber $\rho^{-1}(\sigma(t))$:
$$
\mu\Big(\pi^{-1}(t)\Big)\subseteq \rho^{-1}(\sigma(t)).
$$

Let us call a {\it morphism} of a locally convex bundle $\pi:\varXi\to M$ into a locally convex bundle  $\rho:\varOmega\to N$ each pair of continuous mappings $\mu:\varXi\to\varOmega$ and $\sigma:M\to N$, such that the diagram \eqref{predmorfizm-rassloenij} is commutative and which is linear on each fiber, i.e. for each $t\in M$ the mappings of fibers
$$
\mu:\pi^{-1}(t)\to\rho^{-1}(\sigma(t))
$$
is linear (and continuous due to the continuity of $\mu$).

\paragraph{Dual bundle.}

To each locally convex bundle $\pi:\varXi\to M$ over $\C$ one can assign a vector bundle over $\C$
$$
\pi^\star:\varXi^\star\to M\qquad\Big|\qquad \varXi^\star=\bigsqcup_{t\in M}\pi^{-1}(t)^\star,\qquad \forall u\in \pi^{-1}(t)^\star\quad \pi^\star(u)=t.
$$
We shall call it a {\it dual vector bundle} to the bundle $\pi:\varXi\to M$.
For any compact $K\subset \varXi$ and any point $u\in \varXi^\star$ we put
$$
p_K(u)=\sup_{\xi\in K\cap\pi^{-1}\big(\pi^\star(u)\big)}\abs{u(\xi)}.
$$

Consider the locally convex bundle called {\it trivial} with the fiber $\C$:
$$
\pi_M:\C\times M\to M\qquad\Big|\qquad \pi_M(\lambda,t)=t,\qquad \lambda\in\C,\ t\in M.
$$
Denote by $\pi_{\C}$ the projection to the first component:
$$
\pi_{\C}:\C\times M\to\C\qquad\Big|\qquad \pi_{\C}(\lambda,t)=\lambda,\qquad  \lambda\in\C,\ t\in M.
$$
The fiber $\pi_M^{-1}(t)$ of any point $t\in M$ in this bundle is the set $\C\times\{t\}$, which the mapping $\pi_{\C}$ idenitfies with the field $\C$. As a corollary, if $\mu:\varXi\to \C\times M$ is a morphism of bundles, then on each bundle the following composition is defined:
$$
\pi^{-1}(t)\overset{\mu}{\longrightarrow}\pi_M^{-1}(t)=\C\times\{t\}\overset{\pi_{\C}}{\longrightarrow}\C,
$$
and it is a linear continuous functional on $\pi^{-1}(t)$. We can conclude that the formula
\beq\label{x(t)=pi_C-circ-mu}
x(t)=\pi_{\C}\circ\mu\Big|_{\pi^{-1}(t)},\qquad t\in M,
\eeq
defines a section $x:M\to \varXi^\star$ of a dual bundle of vector spaces $\pi^\star:\varXi^\star\to M$. Denote by $X$ the set of all such sections.

Consider now in each fiber $\pi^{-1}(t)^\star$ the set $X_t$ of functionals $u$, which can be represented in the form $u=\pi_{\C}\circ\mu\Big|_{\pi^{-1}(t)}$ for some morphism of bundles $\mu:\varXi\to\C\times M$. And let  $\overline{X_t}$ denote the closure of this space in the space $\pi^{-1}(t)^\star$ (with respect to the topology, defined by seminorms ${\mathcal P}=\{p_K;\ K\subseteq \varXi\}$).

\bprop\label{PROP:o-dvoistvennom-rasloenii}
If the base $M$ is a Hausdorff space, then on the subbundle
$$
\pi_\star:\varXi_\star=\bigsqcup_{t\in M}\overline{X_t}\to M
$$
the set of sections $X$ and the set of seminorms $\mathcal P$ define a topology, that turns $\varXi_\star$ into a locally convex bundle according to Proposition \ref{PROP:sushestv-topologii-v-Xi}.
\eprop

\bit{

\item[$\bullet$] The bundle $\pi_\star:\varXi_\star\to M$ will be called the {\it dual (locally convex) bundle} to the (locally convex) bundle $\pi:\varXi\to M$.
}\eit

\bpr
We have to check the conditions (i)-(iv) of Proposition \ref{PROP:sushestv-topologii-v-Xi}.

1. The seminorms $p_K$ define on each fiber $\pi^{-1}(t)^\star$ a locally convex topology, which is Hausdorff, since for example it is stronger than the topology of pointwise convergence, and which is defined by the seminorms of the form $p_{\xi}$, where $\xi$ runs through the fiber $\pi^{-1}(t)$.

2. The system $\mathcal P$ is directed in ascending order: each two seminorms $p_K$ and $p_L$ are majorated by a seminorm $p_{K\cup L}$.

3. Let us show that for any section $x\in X$ and for any seminorm $p_K$ the function $t\in T\mapsto p_K(x(t))$ is upper semicontinuous. Take $\e>0$ and consider the set $O=\{t\in M:\ p_K(x(t))<\e\}$. If it is empty then it is open automatically, so we assume that it is non-empty. Take $t_0\in O$, i.e.
$$
p_K(x(t_0))=\sup_{\xi\in K\cap\pi^{-1}(t_0)}\abs{\pi_{\C}\big(\mu(\xi)\big)}<\e
$$
Consider the set $U=\{\xi\in\varXi:\ \abs{\pi_{\C}\big(\mu(\xi)\big)}<\e\}$.
It is open and contains $K\cap\pi^{-1}(t_0)$:
\beq\label{K-cap-pi^(-1)(t_0)-subseteq-U}
K\cap\pi^{-1}(t_0)\subseteq U.
\eeq
We need to show that there is a neighbourhood $V$ of $t_0$ such that
$$
K\cap\pi^{-1}(V)\subseteq U.
$$
Let us assume that this is not so, then for any neighbourhood $V$ of $t_0$ there is a point
$$
\xi_V\in K\cap\pi^{-1}(V)\setminus U.
$$
Since the net $\{\xi_V;\ V\to\{t_0\}\}$ lies in the compact set $K$, it must have a limit point $\xi\in K$ \cite[Theorem 3.1.23]{Engelking}. The projection $\pi$ maps $\{\xi_V;\ V\to\{t_0\}\}$ into a net $\{\pi(\xi_V);\ V\to\{t_0\}\}$ with a limit point $\pi(\xi)$. Therewith $\{\pi(\xi_V);\ V\to\{t_0\}\}$ tends to $t_0$. Since the space $M$ is Hausdorff, the limit point $\pi(\xi)$ coincides with the limit $t_0$. We have:
$$
\xi\in K\cap\pi^{-1}(t_0)\setminus U.
$$
This contradicts to \eqref{K-cap-pi^(-1)(t_0)-subseteq-U}.

4. For any point $t\in M$ the set $\{x(t);\ x\in X\}$ is dense in the space $\overline{X_t}$ by definition of  $\overline{X_t}$.
\epr

\subsection{Value bundles and morphisms of modules}

\paragraph{Value bundle of a module over a commutative involutive algebra.}

Let $A$ be a commutative involutive stereotype algebra. For each point $t\in\Spec(A)$ we denote by $I_t$ the kernel of $t$:
$$
I_t=\{a\in A:\ t(a)=0\}.
$$
Let then $X$ be a left stereotype module over $A$. Following \eqref{DEF:M-cdot-N}, we denote by $I_t\cdot X$ the submodule in $X$, which consists of sums of elements of the form $a\cdot x$, where $a\in I_t$ and $x\in X$:
$$
I_t\cdot X=\left\{ \sum_{i=1}^k a_i\cdot x_i ;\ a_i\in I_t, \ x_i\in X, k\in \Bbb N\right\}
$$
and $\overline{I_t\cdot X}$ is the closure. Put
$$
\Jet^0_A X=\bigsqcup_{t\in M} (X/\overline{I_t\cdot X})^\triangledown,\qquad
\pi^0_{A,X}(x+\overline{I_t\cdot X})=t.
$$
This is a surjection $\pi^0_{A,X}: \Jet^0_A X\to\Spec(A)$. For each element
$x\in X$ we denote by $\jet^0(x)$ the section of $\pi^0_{A,X}$, acting bu the formula
$$
\jet^0(x)(t)=x+\overline{I_t\cdot X},\qquad t\in M.
$$
Denote by ${\mathcal P}(X)$ the set of continuous seminorms $p:X\to\R_+$ of the locally convex space $X$. Each seminorm $p\in{\mathcal P}(X)$ generates a seminorm $p^0$ on the stereotype quotient space
$(X/\overline{I_t\cdot X})^\triangledown$ by the formula
\beq\label{norma-v-X/K_t}
p^0(x+\overline{I_t\cdot X}):=\inf_{y\in \overline{I_t\cdot
X}}p(x+y)=\inf_{y\in I_t\cdot X}p(x+y),\qquad x\in X.
\eeq
We can consider $p^0$ as a function on $\Jet^0_AX$, whose action on each fiber is defined by the formula \eqref{norma-v-X/K_t}.

Denote by $\sigma:A\to {\mathcal C}(M)$ the natural mapping of an algebra $A$ into the algebra of continuous functions on its onvolutive spectrum $M=\Spec(M)$:
\beq\label{A->C(M)} \sigma(a)(t)=t(a),\qquad a\in A,\ t\in M
\eeq
and note the following identity:
 \beq\label{j^0(a-cdot-x)=sigma(a)-cdot-j^0(x)} \jet^0(a\cdot
x)=\sigma(a)\cdot \jet^0(x),\qquad a\in A,\ x\in X.
\eeq
It is proved by sending everything to the left side and substituting the argument $t\in M$:
$$
\jet^0(a\cdot x)(t)-(\sigma(a)\cdot \jet^0(x))(t)=\underbrace{a\cdot
x+\overline{I_t\cdot X}}_{\scriptsize\begin{matrix}\text{\rotatebox{90}{$=$}}
\\ \jet^0(a\cdot
x)(t)\end{matrix}}-\underbrace{t(a)}_{\scriptsize\begin{matrix}\text{\rotatebox{90}{$=$}}
\\ \sigma(a)(t)\end{matrix}}\cdot \underbrace{(x+\overline{I_t\cdot
X})}_{\scriptsize\begin{matrix}\text{\rotatebox{90}{$=$}} \\
\jet^0(x)(t)\end{matrix}}\subseteq
\underbrace{(a-t(a))}_{\scriptsize\begin{matrix}\text{\rotatebox{90}{$\owns$}}
\\ I_t\end{matrix}}\cdot x+\overline{I_t\cdot X}\subseteq \overline{I_t\cdot X}
$$

Now to apply Proposition \ref{PROP:sushestv-topologii-v-Xi} we need the following

\blm\label{LM:poluneprer-sverhu-p^0(j^0(x)(t))} For each element $x\in X$ and each continuous seminorm $p:X\to\R_+$ the mapping $t\in\Spec(A)\mapsto p^0(\jet^0(x)(t))\in\R_+$ is upper semicontinuous.
\elm
\bpr Take $t\in\Spec(A)$ and $\e>0$. The condition
$$
p^0(\jet^0(x)(t))=p^0(x+\overline{I_t\cdot X})=\inf_{y\in \overline{I_t\cdot
X}}p(x+y)=\inf_{y\in I_t\cdot X}p(x+y)<\e
$$
means that for some $y\in I_t\cdot X$ we have
$$
p(x+y)<\e.
$$
This in its turn means that there are $m\in\N$, $y_1,...,y_m\in X$ and $a^1,...,a^m\in I_t$ such that
\beq\label{p-l-x+sum_(k=1)^m-l-a^k-cdot-y_k-r<e}
p\l x+\sum_{k=1}^m a^k\cdot y_k\r<\e.
\eeq
For any sequence of numbers $\lambda=(\lambda^1,...,\lambda^m)$, $\lambda^k\in\C$, we put
$$
f(\lambda)=p\l x+\sum_{k=1}^m (a^k-\lambda^k\cdot 1_A)\cdot y_k\r.
$$
The function $\lambda\mapsto f(\lambda)$ in the point $\lambda=0$ coincides with the left side of \eqref{p-l-x+sum_(k=1)^m-l-a^k-cdot-y_k-r<e}, hence it satisfies the inequality
$$
f(0)<\e.
$$
On the other hand, it is continuous, as a composition of an affine mapping from $\C^m$ into the locally convex space  $X$ and a continuous function $p$ on $X$. Hence there is a number $\delta>0$ such that
\beq\label{max_(k)|lambda^k|<delta=>f(lambda)<e}
\forall \lambda \qquad \max_{1\le k\le m}\abs{\lambda^k}<\delta\quad\Longrightarrow\quad f(\lambda)<\e.
\eeq

Consider the set
$$
U=\{s\in\Spec(A):\ \forall k\in\{1,...,m\}\quad \abs{s(a^k)}<\delta\}.
$$
It is open and contains the point $t$ (since the inclusions $a^k\in I_t$ mean the system of equalities $t(a^k)=0$, $1\le k\le m$). On the other hand, for each point $s\in U$ we can consider the sequence
$$
\lambda^k=s(a^k),
$$
and then, first, $s(a^k-\lambda^k\cdot 1_A)=s(a^k)-\lambda^k\cdot s(1_A)=s(a^k)-s(a^k)=0$, i.e.
$$
a^k-\lambda^k\cdot 1_A\in I_s
$$
and, second, $\max_{k}\abs{\lambda^k}=\max_{k}\abs{s(a^k)}<\delta$, i.e., due to \eqref{max_(k)|lambda^k|<delta=>f(lambda)<e},
$$
p\bigg( x+\sum_{k=1}^m \underbrace{(a^k-\lambda^k\cdot 1_A)}_{\scriptsize\begin{matrix}\text{\rotatebox{90}{$\owns$}}\\ I_s\end{matrix}}\cdot y_k\bigg)=f(\lambda)<\e
$$
This means that for some $z\in I_s\cdot X$ we have
$$
p(x+z)<\e,
$$
and thus,
$$
p^0(\jet^0(x)(s))=p^0(x+\overline{I_s\cdot X})=\inf_{z\in \overline{I_s\cdot
X}}p(x+z)=\inf_{z\in I_s\cdot X}p(x+z)<\e
$$
This is true for each point $s$ from the neighbourhood $U$ of $t$, and this is what we had to prove.
\epr

\btm\label{TH:J^0_AX}
For each stereotype module $X$ over a commutative involutive algebra $A$ the direct sum of the stereotype quotient modules
$$
 \Jet^0_AX=\bigsqcup_{t\in \Spec(A)} (X/\overline{I_t\cdot X})^\triangledown
$$
has a unique topology such that the projection
$$
\pi^0_{A,X}: \Jet^0_A X\to\Spec(A),\qquad \pi^0_{A,X}(x+\overline{I_t\cdot
X})=t,\qquad t\in\Spec(A),\ x\in X
$$
becomes a locally convex bundle with the system of seminorms $\{p^0;\ p\in{\mathcal P}(X)\}$, and the mapping
$$
x\in X\mapsto \jet^0(x)\in\Sec(\pi^0_{A,X})\quad\Big|\quad
\jet^0(x)(t)=x+\overline{I_t\cdot X},\quad t\in\Spec(A),
$$
maps $X$ into a stereotype $A$-module $\Sec(\pi^0_{A,X})$ of continuous sections of $\pi^0_{A,X}$. Therewith,

 \bit{
\item[(i)] the sets
$$
W(x,U,p,\e)=\set{\xi\in  \Jet^0_AX: \ \pi^0_{A,X}(\xi)\in U\ \& \
p^0\Big(\xi-\jet^0(x)\big(\pi^0_{A,X}(\xi)\big)\Big)<\e }
$$
(where $x\in X$, $p\in{\mathcal P}(X)$, $\e>0$, $U$ is an open set in $\Spec(A)$) are the base of topology in $\Jet^0_AX$;

\item[(ii)] for each $x\in X$ and $s\in M$ the sets
$$
W(x,U,p,\e)=\set{\xi\in  \Jet^0_AX: \ \pi^0_{A,X}(\xi)\in U\ \& \
p^0\Big(\xi-\jet^0(x)\big(\pi^0_{A,X}(\xi)\big)\Big)<\e }
$$
(where $p\in{\mathcal P}(X)$, $\e>0$, $U$ is a neighbourhood of $s=\pi^0(\zeta)$ in $\Spec(A)$)
form a local base of the topology of $\Jet^0_AX$ in the point $\zeta=\jet^0(x)(s)$

\item[(iii)] if in addition the spectrum $\Spec(A)$ is a paracompact locally compact space, then the mapping  $\jet^0:X\to\Sec(\pi^0_{A,X})$ is a dense epimorphism.

    }\eit \etm

\bit{

\item[$\bullet$]\label{DEF:rassloenie-znachenij} The locally convex bundle $\pi^0_A: \Jet^0_AX\to \Spec(A)$ is called the {\it value bundle} of the module $X$ over the algebra $A$.

}\eit

\bpr
(i). Put $M=\Spec(A)$. From Lemma \ref{LM:poluneprer-sverhu-p^0(j^0(x)(t))} and Proposition
\ref{PROP:sushestv-topologii-v-Xi} it follows that there is a unique topology on $ \Jet^0_AX$ such that the projection $\pi^0: \Jet^0_A X\to\Spec(A)$ is a locally convex bundle with the seminorms $p^0$, and the sections $\jet^0(x)$, $x\in X$, are continuous. The continuity of the mapping $x\in
X\mapsto \jet^0(x)\in\Sec(\pi^0_{A,X})$ is proved by the implication
$$
p^0(\jet^0(x)(t))=\inf_{y\in I_t\cdot X}p(x+y)\le p(x)\quad \Longrightarrow \quad
p^0_T(\jet^0(x)(t))=\sup_{t\in T}p^0(\jet^0(x)(t))\le p(x)
$$
for any compect set $T\subseteq M$. The property (i) also follows from Proposition
\ref{PROP:sushestv-topologii-v-Xi}.

Let us prove (ii). Take $x\in X$, $s\in M$ and $\zeta=\jet^0(x)(s)$. Consider a base neighbourhood $W(y,V,p,\delta)$ of $\zeta$, i.e. $s\in V$ and $p^0(\zeta-\jet^0(y)(s))<\delta$. Find $\e>0$ such that
$$
p^0(\zeta-\jet^0(y)(s))<\delta-\e.
$$
By Proposition \ref{PROP:sushestv-topologii-v-Xi} (iii), the set
$$
U=\{t\in V: p^0(\jet^0(x)(t)-\jet^0(y)(t))<\delta-\e\}
$$
is a neighbourhood of $s$ in $M$. We obtain that $\zeta\in W(x,U,p,\e)\subseteq W(y,V,p,\delta)$. Indeed, if
$\xi\in W(x,U,p,\e)$, then, first, $\pi^0(\xi)=t\in U\subseteq V$, and, second, $p^0(\xi-\jet^0(x)(t))<\e$, hence
$$
p^0(\xi-\jet^0(y)(t))\le p^0(\xi-\jet^0(x)(t))+p^0(\jet^0(x)(t)-\jet^0(y)(t))<\e+\delta-\e=\delta.
$$

Let us prove (iii). Note first that the set $\jet^0(X)=\{\jet^0(x);\ x\in X\}$ is dense in the ${\mathcal C}(\Spec(A))$-module that it generates
$$
{\mathcal C}(M)\cdot \jet^0(X)=\set{\sum_{i=1}^k b_i\cdot \jet^0(x_i);\ b_i\in {\mathcal C}(M),\ x_i\in X,\ k\in\N}.
$$
Take $b\in {\mathcal C}(M)$ and $x\in X$ for this. Since $A$ is an involutive algebra, the mapping $\sigma:A\to {\mathcal C}(M)$ has a dense image in ${\mathcal C}(M)$. Hence there is a net $a_i\in A$ such that
$$
\sigma(a_i)\overset{{\mathcal C}(M)}{\underset{i\to\infty}{\longrightarrow}}b.
$$
As a corollary,
$$
\jet^0(a_i\cdot x)=\eqref{j^0(a-cdot-x)=sigma(a)-cdot-j^0(x)}=\sigma(a_i)\cdot
\jet^0(x)\overset{{\mathcal C}(M)}{\underset{i\to\infty}{\longrightarrow}}b\cdot \jet^0(x).
$$
Note the that since the image of $X$ under each projection $\pi^0_t$ is dense in the quotient space  $(X/\overline{I_t\cdot X})^\triangledown$, each set
$\jet^0(X)(\Spec(A))$ is dense in each fiber of the budnle $\Jet^0_AX$. Therefore its superset $({\mathcal C}(M)\cdot \jet^0(X))(\Spec(A))$ is also dense in each firber of the bundle $\Jet^0_AX$. Now we can apply property $2^\circ$ on page \pageref{LM:poltnost-v-sloyah=>plotnost-v-secheniyah}: if $M$ is paracompact and locally compact, then the  ${\mathcal C}(M)$-module ${\mathcal C}(M)\cdot \jet^0(X)$ is dense in $\Sec(\pi^0_{A,X})$ (and, as we already understood, $\jet^0(X)$ is dense in ${\mathcal C}(M)\cdot\jet^0(X)$). \epr

\btm\label{TH:Sec(J^0X)-kak-A-modul}
Suppose $X$ is a left stereotype module over a commutative involutive stereotype algebra $A$ with a paracompact locally compact involutive spectrum $\Spec(A)$. Then the formula
\beq\label{DEF:Sec(pi^0_(A,X))-kak-A-modul}
(a\cdot x)(t)=t(a)\cdot x(t),\qquad x\in \Sec(\pi^0_{A,X}),\quad a\in A,\quad t\in\Spec(A),
\eeq
endows the space $\Sec(\pi^0_{A,X})$ of continuous sections of the value bundle with the structure of left stereotype $A$-module, for which the mapping $\jet^0:X\to\Sec(\pi^0_{A,X})$ is a morphism of stereotype $A$-modules.
\etm
\bpr
1. Let us first show that the formula \eqref{DEF:Sec(pi^0_(A,X))-kak-A-modul} endows $\Sec(\pi^0_{A,X})$ with a structure of stereotype $A$-module. Suppose
$$
a_i\overset{A}{\underset{i\to\infty}{\longrightarrow}}0.
$$
Then for any compact set $T\subseteq\Spec(A)$ we have
$$
t(a_i)\overset{\C}{\underset{i\to\infty, t\in T}{\rightrightarrows}}0,
$$
(uniformly by $t\in T$), hence for any compact set $K\subseteq\Sec(\pi^0_{A,X})$
$$
(a_i\cdot x)(t)=t(a_i)\cdot x(t)\overset{X}{\underset{i\to\infty, t\in T, x\in K}{\rightrightarrows}}0.
$$
On the conrary, if
$$
x_i\overset{\Sec(\pi^0_{A,X})}{\underset{i\to\infty}{\longrightarrow}}0,
$$
then for any compact sets $K\subseteq A$ and $T\subseteq\Spec(A)$ the set $\{t(a);\ t\in T, a\in K\}$ is compact in $\C$, hence
$$
(a\cdot x_i)(t)=t(a)\cdot x_i(t)\overset{X}{\underset{i\to\infty, t\in T, x\in K}{\rightrightarrows}}0.
$$

2. Now let us verify that the mapping $\jet^0:X\to\Sec(\pi^0_{A,X})$ is a morphism of $A$-modules. For $a\in A$ and $x\in X$ we have:
\begin{multline*}
\jet^0(a\cdot x)(t)-\big(a\cdot \jet^0(x)\big)(t)=a\cdot x+\overline{I_t\cdot X}-t(a)\cdot \Big(x+\overline{I_t\cdot X}\Big)=a\cdot x-t(a)\cdot x+\overline{I_t\cdot X}-t(a)\cdot\overline{I_t\cdot X}=\\=
\underbrace{\big(a-t(a)\big)}_{\scriptsize\begin{matrix}\text{\rotatebox{90}{$\owns$}}\\ I_t \end{matrix}}\cdot x+\overline{I_t\cdot X}-t(a)\cdot\overline{I_t\cdot X}\subseteq \overline{I_t\cdot X}+
\overline{I_t\cdot X}-t(a)\cdot\overline{I_t\cdot X}\subseteq \overline{I_t\cdot X}
\end{multline*}
therefore
$$
\jet^0(a\cdot x)(t)=\big(a\cdot \jet^0(x)\big)(t).
$$
\epr

\paragraph{Morphisms of modules and their connection with the morphisms of value bundles.}

\btm\label{TH:morf-modul->morf-rassl-znach} Every morphism of stereotype modules
$D:X\to Y$ over a commutative involutive algebra $A$ defines a unique morphism of value bundles $\jet^0(D): \Jet^0_AX\to \Jet^0_A(Y)$,
$$
 \xymatrix @R=2pc @C=1.2pc
 {
  \Jet^0_AX\ar[rr]^{\jet^0(D)}\ar[dr]_{\pi^0_{A,X}} & &  \Jet^0_AY\ar[dl]^{\pi^0_{A,Y}}\\
  & \Spec(A) &
 }
$$
that satisfies the identity
 \beq\label{j^0(Dx)=j^0(D)-circ-j^0(x)}
\jet^0(Dx)=\jet^0(D)\circ \jet^0(x),\qquad x\in X.
 \eeq
$$
 \xymatrix @R=2pc @C=1.2pc
 {
  \Jet^0_AX\ar[rr]^{\jet^0(D)} & &  \Jet^0_AY\\
  & \Spec(A)\ar[ul]^{\jet^0(x)}\ar[ur]_{\jet^0(Dx)} &
 }
$$
\etm
\bpr
The obvious injection
\beq\label{D((I_t-X))-subseteq-(I_t-Y)}
D\l \overline{I_t\cdot X}\r\subseteq \overline{I_t\cdot Y}
\eeq
implies the existence of a natural mapping of quotient spaces:
$$
X/\ \overline{I_t\cdot X}\owns x+\overline{I_t\cdot X}\mapsto Dx+\overline{I_t\cdot Y}\in Y/\ \overline{I_t\cdot Y}
$$
It is continuous, since the initial mapping $D$ is continuous, hence, there exists a unique (and also, continuous) mapping of stereotype quotient spaces (i.e. a mapping of pseudocompletions of usual quotient spaces):
$$
\jet^0(D):\Big( X/\ \overline{I_t\cdot X}\Big)^\triangledown\to \Big( Y/\
\overline{I_t\cdot Y}\Big)^\triangledown
$$
This is true for any $t\in\Spec(A)$, so a mapping of direct sums appears:
$$
\jet^0(D):\bigsqcup_{t\in\Spec(A)}\Big( X/\ \overline{I_t\cdot
X}\Big)^\triangledown\to \bigsqcup_{t\in\Spec(A)}\Big( Y/\ \overline{I_t\cdot
Y}\Big)^\triangledown
$$

The identity \eqref{j^0(Dx)=j^0(D)-circ-j^0(x)} is verified by computation:
for any $t\in\Spec(A)$ and $x\in X$
$$
\big(\jet^0(D)\circ \jet^0(x)\big)(t)=\jet^0(D)\big(\jet^0(x)(t)\big)=\jet^0(D)\l
x+\overline{I_t\cdot X}\r=Dx+\overline{I_t\cdot Y}=\jet^0(Dx)(t).
$$

It remains to prove the continuity of the mapping $\jet^0(D)$. Take a point
$\zeta=\jet^0(x)(t)\in  \Jet^0_AX$, $x\in X$, $t\in T$, and consider a base neighbourhood of its image
$\jet^0(D)(\zeta)=\jet^0(D)\big(\jet^0(x)(t)\big)=\jet^0(Dx)(t)$ under the mapping $\jet^0(D)$:
$$
W(y,V,q,\e)=\set{\upsilon\in  \Jet^0_A(Y): \ \pi^0_{A,Y}(\upsilon)\in U\ \& \
q^0\Big(\upsilon-\jet^0(y)\big(\pi^0_{A,Y}(\upsilon)\big)\Big)<\e }.
$$
(here $q:Y\to\R_+$ is an arbitrary continuous seminorm on $Y$, $V$ a neighbourhood of $t$ in $\Spec(A)$, $y\in Y$, $\e>0$). Since $\jet^0(Dx)(t)\in W(y,U,q,\e)$, we have
$$
q^0\Big(\jet^0(Dx)(t)-\jet^0(y)(t)\Big)<\e.
$$
We can conclude that there exist $\delta>0$ and a neighbourhood $U$ of $t$, containing in $V$, such that
$$
\pi^0_{A,Y}(\upsilon)\in U\quad \& \quad
q^0\Big(\upsilon-\jet^0(Dx)\big(\pi^0_{A,Y}(\upsilon)\big)\Big)<\delta
\qquad\Longrightarrow\qquad
q^0\Big(\jet^0(Dx)\big(\pi^0_{A,Y}(\upsilon)\big)-\jet^0(y)\big(\pi^0_{A,Y}(\upsilon)\big)\Big)<\e
$$
This means that the neighbourhood
$$
W(Dx,U,q,\delta)=\set{\upsilon\in  \Jet^0_A(Y): \ \pi^0_{A,Y}(\upsilon)\in U\quad \& \quad
q^0\Big(\upsilon-\jet^0(Dx)\big(\pi^0_{A,Y}(\upsilon)\big)\Big)<\delta }.
$$
is contained in the neighbourhood $W(y,V,q,\e)$ of the point $\jet^0(Dx)(t)$:
$$
W(Dx,U,q,\delta)\subseteq W(y,V,q,\e).
$$

Further, let us recall that $D:X\to Y$ is a continuous linear mapping. As a corollary, there is a continuous seminorm  $p:X\to\R_+$, satisfying the condition
\beq\label{q(Dx)-le-p(x)} q(Dx)\le p(x),\qquad x\in X.
\eeq
Note that this inequality implies an inequality for the seminorms on the value bundles:
\beq\label{q^0(j^0(D)(xi))-le-p^0(xi)}
q^0\Big(\jet^0(D)(\xi)\Big)\le p^0(\xi),\qquad t\in\Spec(A),\ \xi\in  \Jet^0_AX.
\eeq
It is suffitient to prove this for the points $\xi=\jet^0(x)(t)$, $x\in X$, since they are dense in each fiber:
\begin{multline*}
q^0\Big(\jet^0(D)\big(\jet^0(x)(t)\big)\Big)=\eqref{j^0(Dx)=j^0(D)-circ-j^0(x)}=q^0\Big(\jet^0(Dx)(t)\Big)=q^0\Big(Dx+\overline{I_t\cdot
Y}\Big)= \inf_{v\in I_t\cdot Y}q(Dx+v)
\overset{\scriptsize\eqref{D((I_t-X))-subseteq-(I_t-Y)}}{\le}
 \\ \le \inf_{u\in I_t\cdot X}q(Dx+Du)=
\inf_{u\in I_t\cdot X}q\Big(D(x+u)\Big)\le\eqref{q(Dx)-le-p(x)}\le \inf_{u\in
I_t\cdot X}p(x+u)=p^0\Big(x+\overline{I_t\cdot X}\Big)=p^0\Big(\jet^0(x)(t)\Big)
\end{multline*}

Now we can show that the neighbourhood
$$
W(x,U,p,\delta)=\set{\xi\in  \Jet^0_AX: \ \pi^0_{A,X}(\xi)\in U\ \& \
p^0\Big(\xi-\jet^0(x)\big(\pi^0_{A,X}(\xi)\big)\Big)<\delta }.
$$
of the point $\zeta=\jet^0(x)(t)$ is turned by the mapping $\jet^0(D)$ into the neighbourhood
$W(Dx,U,q,\delta)$ of the point $\jet^0(D)(\zeta)=\jet^0(Dx)(t)$:
\beq\label{j^0(D)(W(x,U,p,delta))-subseteq-W(Dx,U,q,delta)}
\jet^0(D)\Big( W(x,U,p,\delta) \Big)\subseteq W(Dx,U,q,\delta)
\eeq
Indeed, for each point $\xi\in W(x,U,p,\delta)$ the condition $\pi^0_{A,X}(\xi)\in U$ is useful at the end of the chain
$$
\pi^0_{A,Y}(\jet^0(D)(\xi))=\pi^0_{A,X}(\xi)\in U
$$
and the condition $p^0\Big(\xi-\jet^0(x)\big(\pi^0_{A,X}(\xi)\big)\Big)<\delta$ at the end of the chain
\begin{multline*}
q_{\jet^0(D)(\xi)}\Big(\jet^0(D)(\xi)-\jet^0(Dx)\big(\pi^0_{A,Y}(\jet^0(D)(\xi))\big)\Big)=
q_{\pi^0_{A,X}(\xi)}\Big(\jet^0(D)(\xi)-\jet^0(Dx)\big(\pi^0_{A,X}(\xi)\big)\Big)=\eqref{j^0(Dx)=j^0(D)-circ-j^0(x)}=\\=
q_{\pi^0_{A,X}(\xi)}\bigg(\jet^0(D)(\xi)-\jet^0(D)\Big(\jet^0(x)\big(\pi^0_{A,X}(\xi)\big)\Big)\bigg)=
q_{\pi^0_{A,X}(\xi)}\bigg(\jet^0(D)\Big(\xi-\jet^0(x)\big(\pi^0_{A,X}(\xi)\big)\Big)\bigg)\le \eqref{q^0(j^0(D)(xi))-le-p^0(xi)}\le \\
\le p^0\Big(\xi-\jet^0(x)\big(\pi^0_{A,X}(\xi)\big)\Big)<\delta.
\end{multline*}
Together this means that $\jet^0(D)(\xi)\in W(Dx,U,q,\delta)$, and this is what we had to show. \epr

\paragraph{Morphisms with values in a $C^*$-algebra and the Dauns-Hoffman theorem.}

Let $B$ be an involutive subalgebra in the center of an involutive stereotype algebra $F$, and the involutive spectrum $\Spec(B)$ is a paracompact locally compact space.
Then $F$ can be considered as a (formally, left) module over $B$. Consider the value bundle\footnote{Value bundle was defined on page \pageref{DEF:rassloenie-znachenij}.} $\pi^0_B: \Jet^0_B F\to \Spec(B)$. Each point $t\in\Spec(B)$
$$
\overline{I_t\cdot F}=\overline{F\cdot I_t}.
$$
This means that the modules  $\overline{I_t\cdot F}$ are two-sided ideals in $F$. Hence each fiber
$$
\Big(F/\overline{I_t\cdot F}\Big)^\triangledown
$$
is a stereotype algebra, and the projection $F\to \Big(F/\overline{I_t\cdot F}\Big)^\triangledown$ is a homomorphism of stereotype algebras. We can conclude that the space of continuous sections $\Sec(\pi^0_{B,F})$ is also endowed with the structure of stereotype algebra, and the mapping $v:F\to\Sec(\pi^0_{B,F})$ is a homomorphism of stereotype algebras.

In the special case when $F$ is a $C^*$-algebra, the fibers $\Big(F/\overline{I_t\cdot F}\Big)^\triangledown$ and the algebra of section $\Sec(\pi^0_{B,F})$ are $C^*$-algebras as well.

The following variant of the Dauns-Hofmann theorem \cite{Dauns-Hofmann} was stated in the M.~J.~Dupr\'e and R.~M.~Gillette monograph \cite[Theorem 2.4]{Dupre-Gillette} (and for the case of $B=Z(F)$ in the T.~Becker work \cite{Becker}):

\btm[J.~Dauns, K.~H.~Hofmann]\label{TH:Dauns-Hofmann} Let $F$ be a $C^*$-algebra and $B$ its closed involutive subalgebra, that lies in the center of $F$:
$$
B\subseteq Z(F).
$$
Then the mapping $v:F\to\Sec(\pi^0_{B,F})$, that turns $F$ into the algebra of continuous sections of the value bundle $\pi^0_{B,F}: \Jet^0_BF\to\Spec(B)$ over the algebra $B$, is an isomorphism of $C^*$-algebras:
$$
F\cong\Sec(\pi^0_{B,F}).
$$
\etm
\bpr The algebra $B$ is a commutative $C^*$-algebra, hence its spectrum must be a compact space. This implies by Theorem \ref{TH:J^0_AX} that the mapping $v:F\to\Sec(\pi^0_{B,F})$ is not only continuous, but has a dense image in
$\Sec(\pi^0_{B,F})$.

If now $\pi:F\to {\mathcal B}(X)$ is a continuous representation of $F$, then it maps the center $Z(F)$ into the scalar multiples of the identity. In other words, $\pi$  maps $Z(F)$ (and $B$) into the subalgebra $\C\cdot 1_{{\mathcal B}(X)}$ of the algebra
${\mathcal B}(X)$. As a corollary there is a character $t\in\Spec(B)$ such that
$$
\pi(a)=t(a)\cdot 1_{{\mathcal B}(X)},\qquad a\in B.
$$
This implies in its turn that $\pi$ vanishes on $I_t\cdot F$, since
$$
\pi(a\cdot
x)=\pi(a)\cdot\pi(x)=\underbrace{t(a)}_{\scriptsize\begin{matrix}\text{\rotatebox{90}{$=$}}\\
0\end{matrix}}\cdot 1_{{\mathcal B}(X)}\cdot\pi(x)=0,\qquad a\in I_t,\ x\in F.
$$
Hence, $\pi$ vanishes on $\overline{I_t\cdot F}$ as well:
$$
\pi\big|_{\overline{I_t\cdot F}}=0.
$$
We can conclude that if $x$ is a non-zero vector from $F$, then (since there is an irreducible representation $\pi:F\to {\mathcal B}(X)$, that does not vanish on $x$), there is a point $t\in\Spec(B)$ such that
$x\notin\overline{I_t\cdot F}$. This means that $\jet^0(x)(t)\ne 0$.

We see now that the mapping $v:F\to\Sec(\pi^0_BF)$ is injective.
On the other hand, as we already noticed, it has a dense image. Since a continuous homomorphism of $C^*$-algebras always has a closed image \cite[Theorems 3.1.5 and 3.1.6]{Murphy}, the inclusion $F\to \Sec(\pi)$ is an isomorphism of $C^*$-algebras. \epr

The Dauns-Hofmann theorem \ref{TH:Dauns-Hofmann} allows to strengthen Theorem \ref{TH:Sec(J^0X)-kak-A-modul} in the important special case, when the $A$-module $X$ is a $C^*$-algebra, and $A$ is mapped into $X$ by a homomorphism.

\btm\label{TH:B-cong-Sec(val_AB)} Let $\ph :A\to F$ be a homomorphism of stereotype algebras, and $A$ is commutative, $F$ is a $C^*$-algebra, and $\ph(A)$ belongs to the center of $F$:
$$
\ph(A)\subseteq Z(F).
$$
Then the mapping $v:F\to\Sec(\pi^0_{A,F})$, that turns $F$ into the algebra of continuous section of the value bundle $\pi^0_{A,F}: \Jet^0_AF\to\Spec(A)$ over the algebra $A$, is an isomorphism of $C^*$-algebras:
\beq\label{B-cong-Sec(val_AB)}
 F\cong\Sec(\pi^0_{A,F})
\eeq
\etm

\bpr Put $B=\overline{\ph(A)}$. By the Dauns-Hofmann theorem \ref{TH:Dauns-Hofmann},
$$
F\cong \Sec(\pi^0_{B,F}).
$$
So it is sufficient to check the identity
$$
\Sec(\pi^0_{B,F})\cong\Sec(\pi^0_{A,F}).
$$
From the condition (i) of Lemma \ref{LM:overline(ph(I_t)-cdot-B)=B} it follows that in each point $t\in\Spec(B)$ the fibers $\Jet^0_BF(t)$ and $\Jet^0_AF(t\circ\ph)$ coincide:
$$
(\pi^0_{A,F})^{-1}(t\circ\ph)=F\Big/\
\overline{\Ker (t\circ\ph)\kern-15pt\underset{\scriptsize\begin{matrix}\uparrow\\ \text{action}\\ \text{of $A$ on $F$}\end{matrix}}{\cdot}\kern-15pt F}=F\Big/\ \overline{\ph(\Ker (t\circ\ph))\kern-15pt\underset{\scriptsize\begin{matrix}\uparrow\\ \text{action}\\ \text{of $F$ on $F$}\end{matrix}}{\cdot}\kern-15pt F}=F\Big/\
\overline{\underbrace{\overline{\ph(\Ker (t\circ\ph))}}_{\scriptsize\begin{matrix}\phantom{\tiny{\eqref{overline(ph(Ker(t-circ-ph)))=Ker_t}}}\ \|\ \tiny{\eqref{overline(ph(Ker(t-circ-ph)))=Ker_t}}\\ \Ker t\end{matrix}}\cdot F}=
F\Big/\ \overline{\Ker t\cdot F}=(\pi^0_{B,F})^{-1}(t).
$$
And from the condition (ii) of Lemma \ref{LM:overline(ph(I_t)-cdot-B)=B} it follows that in each point  $s\in\Spec(A)$, that lie outside of $\Spec(B)\circ\ph$, the fibers of $\pi^0_{A,F}$ vanish:
$$
\big(\pi^0_{A,F}\big)^{-1}(s)=F/\overline{\Ker
s\kern-15pt\underset{\scriptsize\begin{matrix}\uparrow\\ \text{action}\\ \text{of $A$ on $F$}\end{matrix}}{\cdot}\kern-15pt F}=F\Big/\ \overline{\ph(\Ker
s)\kern-15pt\underset{\scriptsize\begin{matrix}\uparrow\\ \text{action}\\ \text{of $F$ on $F$}\end{matrix}}{\cdot}\kern-15pt F}=F\Big/\ \overline{\underbrace{\overline{\ph(\Ker
s)}}_{\scriptsize\begin{matrix}\phantom{\tiny{\eqref{overline(ph(Ker_s))=B}}}\ \|\ \tiny{\eqref{overline(ph(Ker_s))=B}}\\ F\end{matrix}}\cdot F}=F/\overline{F\cdot F}=F/F=\{0\}.
$$
This implies, first, that the section $x\in\Sec(\pi^0_{A,F})$ is defined by its values on the compact set $\Spec(B)\circ\ph$, i.e. by its restriction on $\Spec(B)\circ\ph$. And, second, that the norm $x$ coincides with the norm of its restriction on $\Spec(B)\circ\ph$. We see that $\Sec(\pi^0_{B,F})$ and $\Sec(\pi^0_{A,F})$ have the same variety of elements and the same norm, hence they are isomorphic as Banach spaces. \epr

\subsection{Jet bundles and differential operators}

If $I$ is a left ideal in an algebra $A$, then, following notations on page \pageref{DEF:M-cdot-N}, for each $n\in\N$ we define the power $I^n$ as the linear space generated by various products of elements from $I$ of the length $n$:
$$
I^n=\sp\{a_1\cdot...\cdot a_n;\ a_1,...,a_n\in I\}.
$$
And the closed power $\overline{I^n}$ is the closure of $I^n$:
$$
\overline{I^n}=\csp\{a_1\cdot...\cdot a_n;\ a_1,...,a_n\in I\}.
$$
Certainly, this is a closed ideal in $A$.

\paragraph{Jet bundle.}

For each $t\in \Spec(A)$ we still denote $I_t=\{a\in A:\ t(a)=0\}$, the ideal in $A$ consisting of elements, vanishing in the point $t$. If now $X$ is a left module over $A$, then for each number $n\in\Z_+$ we consider the ideal  $I_t^{n+1}$, the corresponding submodule $\overline{I_t^{n+1}\cdot X}$ in $X$ and the quotient module
 \beq\label{DEF:prostr-struj}
 \Jet_t^n(X)=\big( X/\overline{I_t^{n+1}\cdot X}\big)^\triangledown.
 \eeq
It is called the {\it jet module} of the degree $n$ of the module $X$ in the point $t$.

Again, each continuous seminorm $p:X\to\R_+$ defines a seminorm on the quotient space $\Jet_t^n(X)=\big( X/\overline{I_t^{n+1}\cdot X}\big)^\triangledown$ by the formula
 \beq\label{polunorma-na-X/I_t^n-X}
p_t^n(x+\overline{I_t^{n+1}\cdot X}):=\inf_{y\in \overline{I_t^{n+1}\cdot
X}}p(x+y),\qquad x\in X.
 \eeq
Consider the direct sum of the sets
$$
\Jet_A^n X=\bigsqcup_{t\in\Spec(A)}  \Jet_t^n(X)=\bigsqcup_{t\in\Spec(A)} \Big( X/\
\overline{I_t^{n+1}\cdot X}\ \Big)^\triangledown
$$
and denote by $\pi^n_{A,X}$ the natural projection of $\Jet_A^n X$ into $\Spec(A)$:
$$
\pi^n_{A,X}:\Jet_A^n X\to \Spec(A),\qquad \pi^n_{A,X}\Big(
x+\overline{I_t^{n+1}\cdot X}\Big)=t,\qquad t\in \Spec(A),\ x\in X.
$$
Besides this, for any vector $x\in X$ we consider the mapping
$$
\jet^n_{A,X}(x):\Spec(A)\to \Jet_A^n X\quad\Big|\quad
\jet_{A,X}^n(x)(t)=x+\overline{I_t^{n+1}\cdot X}.
$$
Of course, for each $x\in X$
$$
\pi^n_{A,X}\circ \jet_A^n(x)=\id_{\Spec(A)}.
$$

\blm\label{LM:poluneprer-sverhu-p^n_t(j(x)(t))} For each element $x\in X$ and each continuous seminorm $p:X\to\R_+$ the mapping $t\in\Spec(A)\mapsto p^n_t(\jet_A^n(x)(t))\in\R_+$ is upper semicontinuous.
\elm
\bpr Take $t\in\Spec(A)$ and $\e>0$. The condition
$$
p^n_t(\jet_A^n(x)(t))=p^n_t(x+\overline{I_t^{n+1}\cdot X})=\inf_{y\in
\overline{I_t^{n+1}\cdot X}}p(x+y)=\inf_{y\in I_t^{n+1}\cdot X}p(x+y)<\e
$$
means that for some $y\in I_t^{n+1}\cdot X$
$$
p(x+y)<\e.
$$
This implies the existence of a number $m\in\N$, a system of vectors $y_1,...,y_m\in X$ and a matrix of vectors $\{a^k_i;\ 1\le i\le n+1,\ 1\le k\le m\}\subseteq I_t$ such that
\beq\label{p-l-x+sum_(k=1)^m-l-prod_(i=1)^(n+1) a^k_i-r-cdot-y_k-r<e}
p\l x+\sum_{k=1}^m \l \prod_{i=1}^{n+1} a^k_i\r\cdot y_k\r<\e.
\eeq
For any number matrix
$$
\lambda=\{\lambda^k_i;\ 1\le i\le n+1,\ 1\le k\le m\},\qquad \lambda^k_i\in\C
$$
we put
$$
f(\lambda)=p\l x+\sum_{k=1}^m\l\prod_{i=1}^{n+1} (a^k_i-\lambda^k_i\cdot 1_A)\r\cdot y_k\r.
$$
The function $\lambda\mapsto f(\lambda)$ in a point $\lambda=0$ coincides with the left side of \eqref{p-l-x+sum_(k=1)^m-l-prod_(i=1)^(n+1) a^k_i-r-cdot-y_k-r<e}, hence, satisfies the inequality
$$
f(0)<\e.
$$
On the other hand, it is continuous, as a composition of a polynomial of $m\cdot (n+1)$ complex variables with the values in the locally convex space $X$, and a continuous function $p$ on $X$. Thus, there exists a number $\delta>0$ such that
\beq\label{max_(i,k)|lambda^k_i|<delta=>f(lambda)<e}
\forall \lambda \qquad \max_{i,k}\abs{\lambda^k_i}<\delta\quad\Longrightarrow\quad f(\lambda)<\e.
\eeq

Consider the set
$$
U=\{s\in\Spec(A):\ \forall i,k\quad \abs{s(a^k_i)}<\delta\}.
$$
It is open and it contains the point $t$ (since the inclusion $a^k_i\in I_t$ means the system of equalities  $t(a^k_i)=0$, $1\le i\le n+1,\ 1\le k\le m$). On the other hand, for each point $s\in U$ we can consider the matrix
$$
\lambda^k_i=s(a^k_i),
$$
and we obtain, first, $s(a^k_i-\lambda^k_i\cdot 1_A)=s(a^k_i)-\lambda^k_i\cdot s(1_A)=s(a^k_i)-s(a^k_i)=0$, i.e.
$$
a^k_i-\lambda^k_i\cdot 1_A\in I_s
$$
and, second, $\max_{i,k}\abs{\lambda^k_i}=\max_{i,k}\abs{s(a^k_i)}<\delta$, i.e., by \eqref{max_(i,k)|lambda^k_i|<delta=>f(lambda)<e},
$$
p\bigg( x+\sum_{k=1}^m \bigg(\prod_{i=1}^{n+1} \underbrace{(a^k_i-\lambda^k_i\cdot 1_A)}_{\scriptsize\begin{matrix}\text{\rotatebox{90}{$\owns$}}\\ I_s\end{matrix}}\bigg)\cdot y_k\bigg)=f(\lambda)<\e
$$
This can be understood as if for some $z\in I_s^{n+1}\cdot X$ we had the inequality
$$
p(x+z)<\e,
$$
which implies
$$
p^n_s(\jet_A^n(x)(s))=p^n_s(x+\overline{I_s^{n+1}\cdot X})=\inf_{z\in
\overline{I_s\cdot X}}p(x+z)=\inf_{z\in I_s^{n+1}\cdot X}p(x+z)<\e
$$
This is true for each point $s$ in the neighbourhood $U$ of $t$, and this was what wee had to check.
\epr

\btm\label{TH:rassloenie-struj} For each stereotype module $X$ over the involutive stereotype algebra $A$ the direct sum of the stereotype modules
$$
\Jet_A^nX=\bigsqcup_{t\in\Spec(A)}  \Jet_t^n(X)=\bigsqcup_{t\in\Spec(A)}
(X/\overline{I_t^{n+1}\cdot X})^\triangledown
$$
has a unique topology, that turns the projection
$$
\pi^n_{A,X}:\Jet_A^n X\to\Spec(A),\qquad \pi^n_{A,X}(x+\overline{I_t^{n+1}\cdot
X})=t,\qquad t\in\Spec(A),\ x\in X
$$
into a locally convex bundle with the system of seminorms $\{p^n;\ p\in{\mathcal P}(X)\}$, for which the mapping
$$
x\in X\mapsto \jet^n(x)\in\Sec(\pi^n)\quad\Big|\quad
\jet^n(x)(t)=x+\overline{I_t^{n+1}\cdot X},\quad t\in\Spec(A),
$$
continuously maps $X$ into the stereotype $A$-module $\Sec(\pi^n_{A,X})$ of continuous sections of $\jet^n_AX$. The sets
$$
W(x,U,p,\e)=\set{\xi\in \Jet_A^n X: \ \pi^n(\xi)\in U\ \& \
p^n_{\pi^n(\xi)}\Big(\xi-\jet^n(x)\big(\pi^n(\xi)\big)\Big)<\e }
$$
(where  $x\in X$, $p\in{\mathcal P}(X)$, $\e>0$, $U$ is an open set in $M$) form a base of the topology in $\Jet_A^nX$.
 \etm

\bit{

\item[$\bullet$] The locally convex bundle $\pi^n_{A,X}:\Jet_A^n X\to \Spec(A)$ is called the {\it jet bundle of the order $n$} of the module $X$ over the algebra $A$.

}\eit

\bpr From Lemma \ref{LM:poluneprer-sverhu-p^n_t(j(x)(t))} and Proposition
\ref{PROP:sushestv-topologii-v-Xi} it follows that there exist a unique topology $\Jet_A^nX$ such that the projection $\jet^n_A X:\Jet_A^n X\to\Spec(A)$ is a locally convex bundle with the seminorms $p^n_t$, and the sections of the form $\jet^n(x)$, $x\in X$, are continuous. The continuity of the mapping
$x\in X\mapsto \jet^n(x)\in\Sec(\jet^n_AX)$ is proved by the implication
$$
p^n_t(\jet^n(x)(t))=\inf_{y\in I_t^{n+1}\cdot X}p(x+y)\le p(x)\quad
\Longrightarrow \quad \sup_{t\in T}p^n_t(\jet^n(x)(t))\le p(x)
$$
for any compact set $T\subseteq M$. The property (i) also follows from the Proposition
\ref{PROP:sushestv-topologii-v-Xi}.

Let us prove (ii). Note first that the set $\jet^n(X)=\{\jet^n(x);\ x\in X\}$ is dense in the ${\mathcal C}(M)$-module
$$
{\mathcal C}(M)\cdot \jet^n(X)=\set{\sum_{i=1}^k b_i\cdot \jet^n(x_i);\ b_i\in {\mathcal C}(M),\ x_i\in X}.
$$
Take $b\in {\mathcal C}(M)$ and $x\in X$. Since $A$ is an involutive algebra, the mapping $\sigma:A\to {\mathcal C}(M)$ has a dense image in ${\mathcal C}(M)$. Hence there is a net $a_i\in A$ such that
$$
\sigma(a_i)\overset{{\mathcal C}(M)}{\underset{i\to\infty}{\longrightarrow}}b.
$$
As a corollary,
$$
\jet^n(a_i\cdot x)=\eqref{j^0(a-cdot-x)=sigma(a)-cdot-j^0(x)}=\sigma(a_i)\cdot
\jet^n(x)\overset{{\mathcal C}(M)}{\underset{i\to\infty}{\longrightarrow}}b\cdot \jet^n(x).
$$
Note further, that since the image of $X$ under each projection $\pi^n_t$ is dense in the quotient space  $(X/\overline{I_t\cdot X})^\triangledown$, the set
$\jet^n(X)(\Spec(A))$ is dense in each fiber of the bundle $\Jet^n_AX$. Therefore its superset $({\mathcal C}(M)\cdot \jet^n(X))(\Spec(A))$ is also dense in each fiber of the bundle $\Jet^n_AX$. Now we can apply the property $2^\circ$ on page \pageref{LM:poltnost-v-sloyah=>plotnost-v-secheniyah}: if $M$ is paracompact and locally compact, then the ${\mathcal C}(M)$-module ${\mathcal C}(M)\cdot \jet^n(X)$ is dense in
$\Sec(\pi^n)$ (and, as we already understood, $\jet^n(X)$ is dense in ${\mathcal C}(M)\cdot \jet^n(X)$). \epr

\paragraph{Differential operators and their relations with morphisms of jet bundles.}

Let $X$ and $Y$ be two left stereotype modules over a stereotype algebra $A$. For each linear (over $\C$) mapping $D:X\to Y$ and each element $a\in A$ the mapping $[D,a]:X\to Y$ acting by formula
\beq\label{DEF:[P,a]}
[D,a](x)=D(a\cdot x)-a\cdot D(x),\qquad x\in X,
\eeq
is called a {\it commutator} of the mapping $D$ with the element $a$. If we have two elements $a,b\in
A$, then the commutator of $D$ with the pair of elements $(a,b)$ is defined as the commutator $[[D,a],b]:X\to Y$ of the mapping $[D,a]$ with the element $b$. Similarly, by induction the commutator with a sequence of elements $(a_0,...,a_n)$ is defined:
$$
[...[D,a_0],... a_n]
$$

A linear (over $\C$) continuous mapping $D:X\to Y$ is called a {\it differential operator} from $X$ into $Y$, if there is a number $n\in\Z_+$ such that for any $a_0,...,a_n\in A$ we have
\beq\label{DEF:diff-oper} [...[D,a_0],...a_n]=0.
\eeq
The least number $n\in\Z_+$ satisfying \eqref{DEF:diff-oper} is called the {\it order} of the differential operator  $D$ and is denoted by $\ord D$.

We denote by $\Diff^n(X,Y)$ the set of all differential operators from $X$ into $Y$ of order not higher than $n$. If we put
$$
\Diff^{-1}(X,Y)=\{D\in Y\oslash X: \ D=0\},
$$
then this sequence of spaces can be defined by the following inductive rule:
\begin{align}
& \Diff^{n+1}(X,Y)=\{D\in Y\oslash X:\ \forall a\in A\quad [D,a]\in \Diff^n(X,Y)\} \label{D^n(X,Y)}
\end{align}
Obviously,
\beq\label{[D^n,A]-subseteq-D^(n-1)}
[\Diff^n,A]\subseteq \Diff^{n-1}
\eeq

\blm
If $D:X\to Y$ is a differential operator of order $n\in\Z_+$, then it maps the submodule $\overline{I^{n+1}\cdot X}$ into the submodule $\overline{I\cdot Y}$:
\beq\label{D((I^(n+1)-X))-subseteq-(I-Y)}
D\in \Diff^n(X,Y)\qquad\Longrightarrow\quad D\l \overline{I^{n+1}\cdot X}\r\subseteq \overline{I\cdot Y}.
\eeq
\elm
\bpr
Since $D$ is continuous, it is sufficient to prove the inclusion
\beq\label{D(I^(n+1)-X)-subseteq-I-Y}
D\in \Diff^n(X,Y)\qquad\Longrightarrow\quad D(I^{n+1}\cdot X)\subseteq I\cdot Y.
\eeq
This is done by the induction by $n\in\Z_+$. For $n=0$ the proposition takes the form
$$
D\in \Diff^0(X,Y)\qquad\Longrightarrow\quad D(I\cdot X)\subseteq I\cdot Y,
$$
and this is evident, since a differential operator of order $n=0$ is just a homogenious mapping
$$
D(a\cdot x)=a\cdot D(x),\qquad a\in A,\ x\in X.
$$

Suppose we have already proved \eqref{D(I^(n+1)-X)-subseteq-I-Y} for $n=k$:
$$
D\in \Diff^k(X,Y)\qquad\Longrightarrow\quad D(I^{k+1}\cdot X)\subseteq I\cdot Y.
$$
For $n=k+1$ we have: if $D\in \Diff^{k+1}(X,Y)$, then for each $a\in A$ we have $[D,a]\in \Diff^k(X,Y)$, and by the assumption of the induction, $[D,a](I^{k+1}\cdot X)\subseteq I\cdot Y$. This means that for any vector $x\in X$ and for any sequence $a_0,...,a_k\in I$
$$
\underbrace{[D,a](a_0\cdot ...\cdot a_k\cdot x)}_{\scriptsize\begin{matrix}\text{\rotatebox{90}{$=$}}
\\ D(a\cdot a_0\cdot ...\cdot a_k\cdot x)-a\cdot D(a_0\cdot ...\cdot a_k\cdot x)\end{matrix}}\kern-30pt\in I\cdot Y
$$
and hence
$$
D(a\cdot a_0\cdot ...\cdot a_k\cdot x)\in I\cdot Y+\underbrace{a\cdot D(a_0\cdot ...\cdot a_k\cdot x)}_{\scriptsize\begin{matrix}\text{\rotatebox{90}{$\owns$}}
\\ I\cdot Y\end{matrix}}\subseteq I\cdot Y
$$
Since this is true for any $a,a_0,...,a_k\in I$, we obtain what we need: $D(I^{k+2}\cdot X)\subseteq I\cdot Y$.
\epr

\btm\label{TH:diff-oper->rassl-struuj} Each differential operator $D:X\to
Y$ of order $n$ defines a morphism of jet bundles $\jet_n[D]:\Jet_A^n X\to
\Jet_A^0(Y)$,
$$
 \xymatrix @R=2pc @C=1.2pc
 {
  \Jet^n_AX\ar[rr]^{\jet^n[D]}\ar[dr]_{\pi^n_{A,X}} & &  \Jet^0_AY\ar[dl]^{\pi^0_{A,Y}}\\
  & \Spec(A) &
 }
$$
such that
 \beq\label{j^0(Dx)=j_n[D]-circ-j^n(x)}
\jet^0(Dx)=\jet_n[D]\circ \jet^n(x),\qquad x\in X.
 \eeq
$$
 \xymatrix @R=2pc @C=1.2pc
 {
  \Jet^n_AX\ar[rr]^{\jet^n[D]} & &  \Jet^0_AY\\
  & \Spec(A)\ar[ul]^{\jet^n(x)}\ar[ur]_{\jet^0(Dx)} &
 }
$$
\etm \bpr Inclusion \eqref{D((I^(n+1)-X))-subseteq-(I-Y)} being applied to the ideal $I_t$, where
$t\in\Spec(A)$
\beq\label{D((I_t^(n+1)-X))-subseteq-(I_t-Y)} D\l
\overline{I_t^{n+1}\cdot X}\r\subseteq \overline{I_t\cdot Y}
\eeq
gives the existence of the mapping of quotient spaces:
$$
X/\ \overline{I_t^{n+1}\cdot X}\owns x+\overline{I_t^{n+1}\cdot X}\mapsto Dx+\overline{I_t\cdot Y}\in Y/\ \overline{I_t\cdot Y}
$$
It is continuous, since the initial mapping $D$ is continuous. Hence there is a natural (continuous) mapping of stereotype quotient spaces (i.e. pseudocompletion of usual quotient spaces):
$$
\jet_n[D]_t:\Big( X/\ \overline{I_t^{n+1}\cdot X}\Big)^\triangledown\to \Big( Y/\ \overline{I_t\cdot Y}\Big)^\triangledown
$$
This is true for any $t\in\Spec(A)$, hence a mapping of direct sums appears:
$$
\jet_n[D]:\bigsqcup_{t\in\Spec(A)}\Big( X/\ \overline{I_t^{n+1}\cdot X}\Big)^\triangledown\to \bigsqcup_{t\in\Spec(A)}\Big( Y/\ \overline{I_t\cdot Y}\Big)^\triangledown
$$

Identity \eqref{j^0(Dx)=j_n[D]-circ-j^n(x)} is verified by computation: for any $t\in\Spec(A)$ and $x\in X$
$$
\big(\jet_n[D]\circ \jet^n(x)\big)(t)=\jet_n[D]\big(\jet^n(x)(t)\big)=\jet_n[D]\l x+\overline{I_t^{n+1}\cdot X}\r=Dx+\overline{I_t\cdot Y}=\jet^0(Dx)(t).
$$

It remains to check the continuity of the mapping $\jet_n[D]$. Take a point $\zeta=\jet^n(x)(t)\in \Jet_A^n X$, $x\in X$, $t\in T$, and consider a base neighbourhood of its image $\jet_n[D](\zeta)=\jet_n[D]\big(\jet^n(x)(t)\big)=\jet^0(Dx)(t)$ under the mapping $\jet_n[D]$:
$$
W(y,V,q,\e)=\set{\upsilon\in \Jet_A^n(Y): \ \pi_Y(\upsilon)\in U\ \& \ q_{\pi_Y(\upsilon)}^n\Big(\upsilon-\jet^0(y)\big(\pi_Y(\upsilon)\big)\Big)<\e }.
$$
(here $q:Y\to\R_+$ is a continuous seminorm on $Y$, $V$ a neighbourhood of $t$ in $\Spec(A)$, $y\in Y$, $\e>0$). Since $\jet^0(Dx)(t)\in W(y,U,q,\e)$, we have
$$
q_t^0\Big(\jet^0(Dx)(t)-\jet^0(y)(t)\Big)<\e.
$$
Hence there is a number $\delta>0$ and a neighbourhood $U$ of $t$, lying in $V$, such that
$$
\pi_Y(\upsilon)\in U\& \ q_{\pi_Y(\upsilon)}^n\Big(\upsilon-\jet^0(Dx)\big(\pi_Y(\upsilon)\big)\Big)<\delta \qquad\Longrightarrow\qquad q_{\pi_Y(\upsilon)}^0\Big(\jet^0(Dx)\big(\pi_Y(\upsilon)\big)-\jet^0(y)\big(\pi_Y(\upsilon)\big)\Big)<\e
$$
Therefore the nieghbourhood
$$
W(Dx,U,q,\delta)=\set{\upsilon\in \Jet_A^n(Y): \ \pi_Y(\upsilon)\in U\ \& \ q_{\pi_Y(\upsilon)}^0\Big(\upsilon-\jet^0(Dx)\big(\pi_Y(\upsilon)\big)\Big)<\delta }.
$$
is contained in the neighbourhood $W(y,V,q,\e)$ of $\jet^0(Dx)(t)$:
$$
W(Dx,U,q,\delta)\subseteq W(y,V,q,\e).
$$

Further, let us recall that $D:X\to Y$ is a continuous linear mapping. As a corollary, there is a continuous seminorm $p:X\to\R_+$ such that
\beq\label{q(Dx)-le-p(x)}
q(Dx)\le p(x),\qquad x\in X.
\eeq
This implies an inequality for seminorms on the jet bundles:
\beq\label{q_t^0(j_n[D](xi))-le-p_t^n(xi)}
q_t^0\Big(\jet_n[D](\xi)\Big)\le p_t^n(\xi),\qquad t\in\Spec(A),\ \xi\in \Jet_A^n X.
\eeq
It is sufficient to prove this for the points $\xi=\jet^n(x)(t)$, $x\in X$, since they are dense in each fiber:
\begin{multline*}
q_t^0\Big(\jet_n[D]\big(\jet^n(x)(t)\big)\Big)=\eqref{j^0(Dx)=j_n[D]-circ-j^n(x)}=q_t^0\Big(\jet^0(Dx)(t)\Big)=q_t^0\Big(Dx+\overline{I_t\cdot Y}\Big)=
\inf_{v\in I_t\cdot Y}q(Dx+v)
\kern-30pt
\overset{\scriptsize\begin{matrix}\eqref{D(I^(n+1)-X)-subseteq-I-Y}
\\
\Downarrow\\
I_t\cdot Y\supseteq D(I_t^{n+1}\cdot Y)
\\
\Downarrow
\end{matrix}}{\le}
\kern-20pt
 \\ \le \inf_{u\in I_t^{n+1}\cdot X}q(Dx+Du)=
\inf_{u\in I_t^{n+1}\cdot X}q\Big(D(x+u)\Big)\le\eqref{q(Dx)-le-p(x)}\le \inf_{u\in I_t^{n+1}\cdot X}p(x+u)=p_t^n\Big(x+\overline{I_t^{n+1}\cdot X}\Big)=p_t^n\Big(\jet^n(x)(t)\Big)
\end{multline*}

Now we can show that the neighbourhood
$$
W(x,U,p,\delta)=\set{\xi\in \Jet_A^n X: \ \pi_X(\xi)\in U\ \& \ p_{\pi_X(\xi)}^n\Big(\xi-\jet^n(x)\big(\pi_X(\xi)\big)\Big)<\delta }.
$$
of the point $\zeta=\jet^n(x)(t)$ under the mapping $\jet_n[D]$ turns into the neighbourhood $W(Dx,U,q,\delta)$ of the point $\jet_n[D](\zeta)=\jet^0(Dx)(t)$:
\beq\label{j_n[D](W(x,U,p,delta))-subseteq-W(Dx,U,q,delta)}
\jet_n[D]\Big( W(x,U,p,\delta) \Big)\subseteq W(Dx,U,q,\delta)
\eeq
Indeed, for each point $\xi\in W(x,U,p,\delta)$ the condition $\pi_X(\xi)\in U$ is useful at the end of the following chain:
$$
\pi_Y(\jet_n[D](\xi))=\pi_X(\xi)\in U
$$
and the condition $p_{\pi_X(\xi)}^n\Big(\xi-\jet^n(x)\big(\pi_X(\xi)\big)\Big)<\delta$ at the end of the chain
\begin{multline*}
q_{\jet_n[D](\xi)}^0\Big(\jet_n[D](\xi)-\jet^0(Dx)\big(\pi_Y(\jet_n[D](\xi))\big)\Big)=
q_{\pi_X(\xi)}^0\Big(\jet_n[D](\xi)-\jet^0(Dx)\big(\pi_X(\xi)\big)\Big)=\eqref{j^0(Dx)=j_n[D]-circ-j^n(x)}=\\=
q_{\pi_X(\xi)}^0\bigg(\jet_n[D](\xi)-\jet_n[D]\Big(\jet^n(x)\big(\pi_X(\xi)\big)\Big)\bigg)=
q_{\pi_X(\xi)}^0\bigg(\jet_n[D]\Big(\xi-\jet^n(x)\big(\pi_X(\xi)\big)\Big)\bigg)\le \eqref{q_t^0(j_n[D](xi))-le-p_t^n(xi)}\le \\
\le
p_{\pi_X(\xi)}^n\Big(\xi-\jet^n(x)\big(\pi_X(\xi)\big)\Big)<\delta.
\end{multline*}
Together this means that $\jet_n[D](\xi)\in W(Dx,U,q,\delta)$, and this is what wee need.
\epr

\paragraph{Differential operators on algebras.}

Each homomorphism of algebras $\ph:A\to B$ defines on $B$ a steructure of left $A$-module by the formula:
$$
a\cdot y=\ph(a)\cdot y,\qquad a\in A,\ y\in Y=B.
$$
If $D:A\to B$ is an arbitrary linear (over $\C$) mapping, then formula \eqref{DEF:[P,a]} for its commutator with an element $a\in A$ turns into the formula
$$
[D,a](x)=D(a\cdot x)-\ph(a)\cdot D(x),\qquad a,x\in A,
$$
We call this operator the {\it commutator of the operator $D:A\to B$ and the element $a\in A$ with respect to the homomorphism} $\ph:A\to B$. Besides this for each element $b\in B$ we shall consider the linear (over $\C$) mapping  $b\cdot D:A\to B$, defined by formula
$$
(b\cdot D)(x)=b\cdot D(x),\qquad x\in A.
$$

\bprop The following identities hold:
\begin{align}
& [\ph,a]=0, && a\in A \label{[ph,a]=0} \\
& [b\cdot D, a]=b\cdot [D,a]+[b,\ph(a)]\cdot D, && a\in A,\ b\in B,\ D\in B\oslash A, \label{[b-cdot-P, a]} \\
& [b\cdot\ph,a]=[b,\ph(a)]\cdot\ph, && a\in A, \ b\in B. \label{[b-cdot-ph,a_0]}
\end{align}
\eprop
\bpr
The first and the second identities are verified by computation: for $x\in A$ we have
$$
[\ph,a](x)=\ph(a\cdot x)-\ph(a)\cdot\ph(x)=0,
$$
and
\begin{multline*}
[b\cdot D, a](x)=(b\cdot D)(a\cdot x)-\ph(x)\cdot(b\cdot D)(x)=b\cdot D(a\cdot x)-\ph(x)\cdot b\cdot D(x)=\\=
b\cdot D(a\cdot x)-b\cdot \ph(x)\cdot D(x)+b\cdot \ph(x)\cdot D(x)-\ph(x)\cdot b\cdot D(x)=
b\cdot \Big( D(a\cdot x)-\ph(x)\cdot D(x)\Big)+\Big(b\cdot \ph(x)-\ph(x)\cdot b\Big)\cdot D(x)=\\=
b\cdot [D,a](x)+[b,\ph(x)]\cdot D(x)=\Big(b\cdot [D,a]+[b,\ph(a)]\cdot D\Big)(x).
\end{multline*}
And the third one becomes a corollary of the first and the second:
$$
[b\cdot\ph,a]=\eqref{[b-cdot-P, a]}=b\cdot\kern-8pt\underbrace{[\ph,a_0]}_{\tiny\begin{matrix} \phantom{\eqref{[ph,a]=0}} \ \text{\rotatebox{90}{$=$}} \ \eqref{[ph,a]=0}
\\  0  \end{matrix}}\kern-8pt +[b,\ph(a)]\cdot \ph
$$
\epr

A given homomorphism of algebras $\ph:A\to B$ defines a system of differential operators from $A$ into $B$, which we shall denote by $\Diff^n(\ph)$, or just $\Diff^n$. It is defined by the following inductive rules:
\begin{align}
& \Diff^{-1}(\ph)=\{D\in B\oslash A: \ D=0\}, \label{D^(-1)(ph)} \\
& \Diff^{n+1}(\ph)=\{D\in B\oslash A:\ \forall a\in A\quad [D,a]\in \Diff^n(\ph)\} \label{D^n(ph)}
\end{align}
The spaces $\Diff^n(\ph)$ form an expanding sequence:
$$
0=\Diff^{-1}(\ph)\subseteq \Diff^0(\ph)\subseteq \Diff^1(\ph)\subseteq ...\subseteq \Diff^n(\ph)\subseteq \Diff^{n+1}(\ph)\subseteq...
$$
Besides this we shall need a sequence of spaces in $B$, denoted by $Z^n(\ph)$, or just $Z^n$, and defined by the following inductive rules:
\begin{align}
& Z^0(\ph)=0 \label{Z^0}\\
& Z^{n+1}(\ph)=\{b\in B:\ \forall a\in A\quad [b,\ph(a)]\in Z^n(\ph)\} \label{Z^n}
\end{align}
The spaces $Z^n(\ph)$ also form an expanding sequence:
\beq\label{Z^n(ph)}
0=Z^0(\ph)\subseteq Z^1(\ph)\subseteq ...\subseteq Z^n(\ph)\subseteq Z^{n+1}(\ph)\subseteq...
\eeq

\bprop The following inclusions hold:
\begin{align}
& [Z^q\cdot \Diff^p,A]\subseteq Z^q\cdot \Diff^{p-1}+Z^{q-1}\cdot \Diff^p, \label{[Z^q-cdot-D^p,A]}\\
& Z^0\cdot \Diff^p\subseteq \Diff^{-1}=0,\qquad p\ge 0 \label{Z.D^p=0} \\
& Z^q\cdot \Diff^0\subseteq \Diff^{q-1},\qquad q\ge 0 \label{Z^q.D^0->D^(q-1)} \\
& Z^q\cdot \Diff^p\subseteq \Diff^{q+p-1},\qquad q\ge 0 \label{Z^q.D^p->D^(q+p-1)}
\end{align}
\eprop
\bpr

1. For proving \eqref{[Z^q-cdot-D^p,A]} we take arbitrary $b\in Z^q$, $D\in \Diff^p$ and $a\in A$. Then
$$
[b\cdot D, a]=\overset{\tiny\begin{matrix}Z^q \\ \text{\rotatebox{90}{$\in$}} \end{matrix}}{b}\cdot \underbrace{[\overset{\tiny\begin{matrix}\Diff^p \\ \text{\rotatebox{90}{$\in$}} \end{matrix}}{D},a]}_{\tiny\begin{matrix}\text{\rotatebox{90}{$\owns$}} \\ \Diff^{p-1}\end{matrix}}+\underbrace{[\overset{\tiny\begin{matrix}Z^q \\ \text{\rotatebox{90}{$\in$}} \end{matrix}}{b},\ph(a)]}_{\tiny\begin{matrix}\text{\rotatebox{90}{$\owns$}} \\ Z^{q-1}\end{matrix}}\cdot \overset{\tiny\begin{matrix}\Diff^p \\ \text{\rotatebox{90}{$\in$}} \end{matrix}}{D}\in Z^q\cdot \Diff^{p-1}+Z^{q-1}\cdot \Diff^p
$$

2. Identity \eqref{Z.D^p=0} is obvious (the multiplication by zero always gives zero).

3. Identity \eqref{Z^q.D^0->D^(q-1)} is proved by induction. For $q=0$ it is true, since it is a special case of  \eqref{Z.D^p=0} with $p=0$:
$$
Z^0\cdot \Diff^0\subseteq \Diff^{-1}=0.
$$
Suppose we proved it for some $q=n$:
\beq\label{Z^n-cdot-D^0-subseteq-D^(n-1)}
Z^n\cdot \Diff^0\subseteq \Diff^{n-1}
\eeq
Then for $q=n+1$ we have: if $b\in Z^q=Z^{n+1}$ and $D\in \Diff^0$, then
$$
\Big(\forall a\in A\quad
[b\cdot D, a]=b\cdot \underbrace{[\overset{\tiny\begin{matrix}\Diff^0 \\ \text{\rotatebox{90}{$\in$}} \end{matrix}}{D},a]}_{\tiny\begin{matrix}\text{\rotatebox{90}{$=$}} \\ 0\end{matrix}}+\underbrace{[b,\ph(a)]}_{\tiny\begin{matrix}\text{\rotatebox{90}{$\owns$}} \\ Z^n\end{matrix}}\cdot \overset{\tiny\begin{matrix}\Diff^0 \\ \text{\rotatebox{90}{$\in$}} \end{matrix}}{D}\overset{\eqref{Z^n-cdot-D^0-subseteq-D^(n-1)}}{\in} \Diff^{n-1}\Big)
\quad\Longrightarrow\quad b\cdot D\in \Diff^n=\Diff^{q-1}
$$

4. Formula \eqref{Z^q.D^p->D^(q+p-1)} is also proved by induction by $p$. For $p=0$ it turns into \eqref{Z^q.D^0->D^(q-1)} which is already proved. Suppose we have proved it for some $p=n$:
\beq\label{Z^q.D^n->D^(q+n-1)}
Z^q\cdot \Diff^n\subseteq \Diff^{q+n-1},\qquad q\ge 0
\eeq
Then for $p=n+1$ and $a_1,...,a_q\in A$ we have
$$
[Z^q\cdot \Diff^{n+1},a_1]\overset{\eqref{[Z^q-cdot-D^p,A]}}{\subseteq}
\underbrace{Z^q\cdot \Diff^n}_{\tiny\begin{matrix}
\phantom{\eqref{Z^q.D^n->D^(q+n-1)}}
\ \text{\rotatebox{90}{$\supseteq$}}\ \eqref{Z^q.D^n->D^(q+n-1)} \\ \Diff^{q+n-1}\end{matrix}}
+Z^{q-1}\cdot \Diff^{n+1}\subseteq \Diff^{q+n-1}+Z^{q-1}\cdot \Diff^{n+1}
$$
$$
\Downarrow
$$
\begin{multline*}
[[Z^q\cdot \Diff^{n+1},a_1],a_2]\subseteq \underbrace{[\Diff^{q+n-1},a_2]}_{\tiny\begin{matrix}
\phantom{\eqref{[D^n,A]-subseteq-D^(n-1)}}
\ \text{\rotatebox{90}{$\supseteq$}}\ \eqref{[D^n,A]-subseteq-D^(n-1)} \\ \Diff^{q+n-2}\end{matrix}}+[Z^{q-1}\cdot \Diff^{n+1},a_2]\subseteq \\ \subseteq
\Diff^{q+n-2}+\underbrace{Z^{q-1}\cdot \Diff^n}_{\tiny\begin{matrix}
\phantom{\eqref{Z^q.D^n->D^(q+n-1)}}
\ \text{\rotatebox{90}{$\supseteq$}}\ \eqref{Z^q.D^n->D^(q+n-1)} \\ \Diff^{q+n-2}\end{matrix}}
+Z^{q-1}\cdot \Diff^{n+1}\subseteq \Diff^{q+n-2}+Z^{q-2}\cdot \Diff^{n+1}
\end{multline*}
$$
\Downarrow
$$
$$
\cdots
$$
$$
\Downarrow
$$
$$
[...[[Z^q\cdot \Diff^{n+1},a_1],a_2],... a_q]\subseteq \Diff^{q+n-q}+Z^{q-q}\cdot \Diff^{n+1}=
\Diff^n+Z^0\cdot \Diff^{n+1}=\Diff^n
$$
This is true for any $a_1,...,a_q\in A$, hence
$$
Z^q\cdot \Diff^{n+1}\subseteq \Diff^{q+n}=\Diff^{q+p-1}
$$
\epr

\bprop
For each $b\in B$
\beq\label{b.ph-in-D^n(ph)<=>b-in-Z^(n+1)(ph)}
b\cdot\ph\in\Diff^n(\ph)\quad\Longleftrightarrow\quad b\in Z^{n+1}(\ph)
\eeq
\eprop
\bpr
We use here formula \eqref{[b-cdot-ph,a_0]}, which we generalize to the identity
$$
[...[b\cdot\ph,a_0],...a_n]=[...[b,\ph(a_0)],...\ph(a_n)]\cdot\ph, \qquad a_0,...,a_n\in A.
$$
It is seen here that if $b\in Z^{n+1}(\ph)$, then the coefficient of $\ph$ on the right is zero, hence the right side is zero. This is true for any $a_0,...,a_n\in A$, therefore $b\cdot\ph\in\Diff^n(\ph)$. On the contrary, if  $b\cdot\ph\in\Diff^n(\ph)$, then the left side is zero, and, in particular, when we substitute unit as an argument we have
$$
0=[...[b\cdot\ph,a_0],...a_n](1_A)=[...[b,\ph(a_0)],...\ph(a_n)]\cdot\ph(1_A)=[...[b,\ph(a_0)],...\ph(a_n)]\cdot 1_B=[...[b,\ph(a_0)],...\ph(a_n)].
$$
Again this is true for any $a_0,...,a_n\in A$, hence $b\in Z^{n+1}(\ph)$.
\epr

\paragraph{Differential operators with values in $C^*$-algebras.}

Let $A$ be an involutive closed subalgebra in a $C^*$-algebra $B$. Consider the sequence $Z^n(A)$ of subspaces in  $B$, defined by the following inductive rules:
\begin{align}
& Z^0(A)=0 \label{Z^0(A)}\\
& Z^{n+1}(A)=\{b\in B:\ \forall a\in A\quad [b,a]\in Z^n(A)\} \label{Z^n(A)}
\end{align}

\btm\label{TH:stabilizatsiya-Z^n} For any $C^*$-algebra $B$ and any its closed involutive subalgebra $A$ the sequence of subspaces $Z^n(A)$ is stabilized starting from the number $n=1$:
\beq\label{Z^n(A)=Z^(n+1)(A)}
Z^n(A)=Z^{n+1}(A),\qquad n\ge 1.
\eeq
\etm

The idea of the proof belongs to Yu.~N.~Kuznetsova and is based on two lemmas.

Let us call a {\it projector} in a locally convex space $X$ an arbitrary linear continuous operator $P:X\to X$ with the idempotency property:
$$
P^2=P.
$$
Every such an operator acts as an identity operator on its image \cite[p.37]{Bourbaki-tvs}:
\beq\label{y-in-Im(P)=>P(y)=y}
\forall y\in \Im P\qquad P(y)=y.
\eeq

\blm\label{LM:TP=PT} For any locally convex space $X$, any projector $P$ in $X$ and any linear continuous operator  $T$ in $X$ the commutation condition
$$
[T,P]=0
$$
is equivalent to the fact that the image $\Im P$ and the kernel $\Ker P$ of $P$ are invariant with respect to $T$:
\beq\label{TP=PT-2}
T(\Im P)\subseteq \Im P\quad\&\quad T(\Ker P)\subseteq \Ker P.
\eeq
\elm
\bpr
1. Suppose $TP=PT$. Then for $x\in \Im P$, i.e. $x=Px$, we have $Tx=TPx=PTx\in \Im P$. On the other hand, if $x\in \Ker P$, i.e. $Px=0$, then $PTx=TPx=0$, and thus, $Tx\in \Ker P$.

2. On the conrary, suppose \eqref{TP=PT-2} is true. Then for any vector $x\in X$ we can put
$$
x_\parallel=Px,\qquad x_\perp=x-x_\parallel,
$$
and then
$$
\begin{cases}x_\parallel\in\Im P\quad\Longrightarrow\quad Tx_\parallel\in\Im P\\
x_\perp\in\Ker P\quad\Longrightarrow\quad Tx_\perp\in\Ker P\end{cases},
$$
hence
$$
PTx=PT(x_\parallel+x_\perp)=P\Big(\underbrace{Tx_\parallel}_{\tiny\begin{matrix}\text{\rotatebox{90}{$\owns$}}\\
\Im P\end{matrix}}+\underbrace{Tx_\perp}_{\tiny\begin{matrix}\text{\rotatebox{90}{$\owns$}}\\
\Ker P\end{matrix}}\Big)=P\underbrace{Tx_\parallel}_{\tiny\begin{matrix}\text{\rotatebox{90}{$\owns$}}\\
\Im P\end{matrix}}=\eqref{y-in-Im(P)=>P(y)=y}=Tx_\parallel=TPx.
$$
\epr

\blm\label{LM:TP=PT-(2)} For any locally convex space $X$, any projector $P$ in $X$ and any linear continuous operator $T$ in $X$ the condition
$$
[[T,P],P]=0
$$
implies the condition
$$
[T,P]=0.
$$
\elm
\bpr From Lemma \ref{LM:TP=PT} we have that the spaces $\Im P$ and $\Ker P$ are invariant for the operator $[T,P]$:
$$
[T,P](\Im P)\subseteq \Im P\quad\&\quad [T,P](\Ker P)\subseteq \Ker P.
$$
The first of these conditions means that for each $x\in\Im P$ we have
$$
TPx-\underbrace{PTx}_{\tiny\begin{matrix}\text{\rotatebox{90}{$\owns$}}\\
\Im P\end{matrix}}\in \Im P\quad\Longrightarrow\quad T\underbrace{Px}_{\tiny\begin{matrix}\text{\rotatebox{90}{$=$}}\\
x\end{matrix}}\in \Im P\quad\Longrightarrow\quad Tx\in \Im P.
$$
And the second that for each $x\in\Ker P$ we have
$$
T\underbrace{Px}_{\tiny\begin{matrix}\text{\rotatebox{90}{$=$}}\\ 0\end{matrix}}-PTx\in \Ker P\quad\Longrightarrow\quad PTx\in \Ker P\quad\Longrightarrow\quad \underbrace{P^2Tx}_{\tiny\begin{matrix}\text{\rotatebox{90}{$=$}}\\ PTx\end{matrix}}=0
\quad\Longrightarrow\quad Tx\in \Ker P.
$$
These two chains together mean \eqref{TP=PT-2}, and by Lemma \ref{LM:TP=PT} this is equivalent to $[T,P]=0$.
\epr

\bpr[Proof of Theorem \ref{TH:stabilizatsiya-Z^n}] 1. Suppose first that $B$ is the algebra of bounded operators in a Hilbert space $X$: $B={\mathcal B}(X)$. Note that if we replace $A$ by its bicommutatn $A^{!!}$ in ${\mathcal B}(X)$ the sequence of subspaces $Z^n(A)$ doesn't change:
$$
Z^n(A^{!!})=Z^n(A).
$$
So we can think from the very beginning that $A$ coincides with its bicommutant, i.e. is a von Neumenn algebra in ${\mathcal B}(X)$. Then $A$ is generated by the system of its orthogonal projectors $P\in A$: $A=\{P\}^{!!}$. For any such a projector $P$ and for any element $T\in Z^2(A)$ we have the equality $[[T,P],P]=0$, which by Lemma  \ref{LM:TP=PT-(2)} implies $[T,P]=0$. Since the projectors $P$ generate $A$, this implies the equality
$$
[T,a]=0,\qquad a\in A.
$$
I.e. $T\in Z^1(A)$, and we obtain $Z^1(A)\supseteq Z^2(A)$. This is equivalent to the equality $Z^1(A)=Z^2(A)$, which in its turn implies \eqref{Z^n(A)=Z^(n+1)(A)}.

2. Consider now the general case when $B$ is an arbitrary $C^*$-algebra. Then $B$ is included into some algebra ${\mathcal B}(X)$, and we have the following chain:
$$
A\subseteq B\subseteq {\mathcal B}(X).
$$
Again, suppose $Z^n(A)$ is a sequence of subspaces in $B$, defined by the equalities \eqref{Z^n(A)}, and $Z^n_X(A)$ is the same sequence in ${\mathcal B}(X)$ (i.e. $Z^n_X(A)$ is defined by \eqref{Z^n(A)}, where $B$ is replaced by  ${\mathcal B}(X)$). Let us show that
 \beq\label{Z^n(A)=B-cap-Z^n_H(A)}
Z^n(A)=B\cap Z^n_X(A),\qquad n\in \Z_+.
 \eeq
For $n=0$ this is true by definition. Suppose we already proved this for all indices not greater than some $n\in\Z_+$. Then for $n+1$ we have the following chain:
\begin{multline*}
c\in B\cap Z^{n+1}_X(A)\quad\Leftrightarrow\quad \begin{cases}c\in B \\ \forall a\in A\quad [c,a]\in Z^n_X(A) \end{cases}
\quad\Leftrightarrow\quad \begin{cases} c\in B \\ \forall a\in A\quad  [c,a]=c\cdot a-a\cdot c\in B \\ \forall a\in A\quad  [c,a]\in Z^n_X(A) \end{cases}
\quad\Leftrightarrow\\ \Leftrightarrow\quad \forall a\in A\quad [c,a]\in B\cap Z^n_X(A)\kern-13pt\underset{\tiny\begin{matrix}\uparrow\\ \text{inductive}\\ \text{assumption}\end{matrix}}{=}\kern-13pt Z^n(A) \quad\Leftrightarrow\quad c\in Z^{n+1}(A).
\end{multline*}

For the sequence $Z^n_X(A)$ we already proved that
$$
Z^1_X(A)=Z^n_X(A),\qquad n\ge 1.
$$
Together with \eqref{Z^n(A)=B-cap-Z^n_H(A)} this implies
$$
Z^n(A)=B\cap Z^n_X(A)=B\cap Z^1_X(A)=Z^1(A).
$$
\epr

Theorem \ref{TH:stabilizatsiya-Z^n} implies

\btm\label{TH:stabilizatsiya-Z^n(ph)} If $A$ is an involutive stereotype algebra, and $B$ a $C^*$-algebra, then for any involutive homomorphism of stereotype algebras $\ph:A\to B$ the sequence of subspaces $Z^n(\ph)$, defined by  \eqref{Z^0} and \eqref{Z^n}, is stabilized starting from 1:
\beq\label{Z^n(ph)=Z^(n+1)(ph)}
Z^n(\ph)=Z^{n+1}(\ph),\qquad n\ge 1.
\eeq
\etm
\bpr
The closed image $\overline{\ph(A)}$ of the algebra $A$ under the mapping $\ph$ is an involutive subalgebra in $B$, and
$$
Z^n(\ph)=Z^n\Big(\overline{\ph(A)}\Big).
$$
Hence \eqref{Z^n(ph)=Z^(n+1)(ph)} is just a corollary of \eqref{Z^n(A)=Z^(n+1)(A)}.
\epr

\subsection{Tangent and cotangent bundles}

Recall that tangent and cotangent spaces were defined on pages \pageref{DEF:kasatelnyj-vektor} and \pageref{DEF:T_s^(C-star)[A]}.

\paragraph{Cotangent bundle $T^\star[A]$.}

Let $A$ be an involutive stereotype algebra, and $t\in\Spec(A)$. Note that each continuous seminorm $p:A\to\R_+$ defines a seminorm
$p'_t$ on the quotient space $\C T^\star_t[A]=\big(
I_t/\overline{I_t^2}\big)^\triangledown$ by the formula
 \beq\label{polunorma-na-CT^*(X)}
p'_t(x+\overline{I_t^2}):=\inf_{y\in
\overline{I_t^2}}p(x+y),\qquad x\in X.
 \eeq
Consider the direct sum of the sets $\C T^\star_t(A)$
$$
\C T^\star(A)=\bigsqcup_{t\in\Spec(A)} \C T^\star_t(A)=\bigsqcup_{t\in\Spec(A)} \Big(
I_t/\ \overline{I_t^2}\ \Big)^\triangledown
$$
and denote by $\pi$ the natural projection of $\C T^\star[A]$ onto $\Spec(A)$:
$$
\pi:\C T^\star[A]\to \Spec(A),\qquad \pi\Big( x+\overline{I_t^2}\Big)=t,\qquad t\in \Spec(A),\ x\in I_t[A].
$$
Besides this for each vector $x\in A$ we consider the mapping
$$
T^\star(x):\Spec(A)\to\C T^\star(X)\quad\Big|\quad T^\star(x)(t)=x-t(x)\cdot
1_A+\overline{I_t^2}.
$$
Certainly, for any $x\in A$
$$
\pi\circ T^\star(x)=\id_{\Spec(A)},
$$
i.e. $T^\star(x)$ is a section of the vector bundle $\pi$.

\blm\label{LM:poluneprer-sverhu-p^*_t(T^*(x)(t))} For any element $x\in A$ and for any continuous seminorm  $p:A\to\R_+$ the mapping $t\in\Spec(A)\mapsto
p'_t(T^\star(x)(t))\in\R_+$ is upper semicontinuous.
\elm
\bpr Take $t\in\Spec(A)$ and $\e>0$. The condition
$$
p'_t(T^\star(x)(t))=\inf_{y\in
\overline{I_t^2}}p(x-t(x)\cdot 1_A+y)=\inf_{y\in I_t^2}p(x-t(x)\cdot 1_A+y)<\e
$$
means that for some $y\in I_t^2$ we have
$$
p(x-t(x)\cdot 1_A+y)<\e.
$$
This means in its turn that there exists a number $m\in\N$ and vectors
$a^1,...,a^m,b^1,...,b^m\in I_t$ such that
 \beq\label{p-l-x+sum_(k=1)^m-a^k_i-cdot-b^k-r<e}
p\l x-t(x)\cdot 1_A+\sum_{k=1}^m a^k\cdot b^k\r<\e.
 \eeq
For any numbers
$$
\lambda\in\C,\qquad \alpha=\{\alpha^k;\ 1\le k\le m\}\subseteq\C,
\beta=\{\beta^k;\ 1\le k\le m\}\subseteq\C
$$
let us put
$$
f(\lambda,\alpha,\beta)=p\l x-\lambda\cdot 1_A+\sum_{k=1}^m (a^k-\alpha^k\cdot
1_A)\cdot (b^k-\beta^k\cdot 1_A)\r.
$$
The function $(\lambda,\alpha,\beta)\mapsto f(\lambda,\alpha,\beta)$ in the point
$(\lambda,\alpha,\beta)=(t(x),0,0)$ concides with the left side of
\eqref{p-l-x+sum_(k=1)^m-a^k_i-cdot-b^k-r<e}, hence satisfies the inequality
$$
f(t(x),0,0)<\e.
$$
On the other hand it is continuous as a composition of a polynomial of $2m+1$ complex variables with the values in the locally convex space $X$ and a continuous function $p$ on $X$. Hence, there exists a number $\delta>0$ such that
 \beq\label{max_k|lambda^k|<delta=>f(lambda)<e} \forall \lambda,\alpha,\beta \qquad
 |\lambda-t(x)|<\delta\quad\&\quad\max_k\abs{\alpha^k}<\delta
 \quad\&\quad\max_k\abs{\beta^k}<\delta\quad\Longrightarrow\quad f(\lambda,\alpha,\beta)<\e.
 \eeq

Consider the set
$$
U=\{s\in\Spec(A):\ \abs{s(x)-t(x)}<\delta\quad \&\quad
\max_k\abs{s(a^k)}<\delta\quad \&\quad \max_k\abs{s(b^k)}<\delta\}.
$$
It is open and contains the point $t$ (since $\abs{t(x)-t(x)}=0<\delta$, and the inclusions $a^k,b^k\in I_t$ imply the system of inequalities $t(a^k)=t(b^k)=0$, $1\le
k\le m$). On the other hand, for any point $s\in U$ we can consider the numbers
$$
\lambda=s(x),\qquad \alpha^k=s(a^k),\qquad \beta^k=s(b^k),
$$
and we obtain, first, $s(a^k-\alpha^k\cdot 1_A)=s(a^k)-\alpha^k\cdot
s(1_A)=s(a^k)-s(a^k)=0$, i.e.
$$
a^k-\alpha^k\cdot 1_A\in I_s
$$
second, by the same reasons,
$$
b^k-\beta^k\cdot 1_A\in I_s
$$
and, third,
$$
\abs{\lambda-t(x)}<\delta,\qquad
\max_k\abs{\alpha^k}=\max_k\abs{s(a^k)}<\delta, \qquad
\max_k\abs{\beta^k}=\max_k\abs{s(b^k)}<\delta,
$$
i.e. by \eqref{max_k|lambda^k|<delta=>f(lambda)<e},
$$
p\bigg( x-\lambda\cdot 1_A+\sum_{k=1}^m  \underbrace{(a^k-\alpha^k\cdot
1_A)}_{\scriptsize\begin{matrix}\text{\rotatebox{90}{$\owns$}}\\
I_s\end{matrix}}\cdot \underbrace{(b^k-\beta^k\cdot
1_A)}_{\scriptsize\begin{matrix}\text{\rotatebox{90}{$\owns$}}\\
I_s\end{matrix}}\bigg)=f(\lambda,\alpha,\beta)<\e
$$
This can be understood as follows: for some point $z\in I_s^2$ we have the inequality
$$
p(x-s(x)\cdot 1_A+z)<\e,
$$
which in its turn imply the inequality
$$
p_s(T^\star(x)(s))=p_s(x-s(x)\cdot 1_A+\overline{I_s^2})=\inf_{z\in
\overline{I_s^2}}p(x-s(x)\cdot 1_A+z)=\inf_{z\in I_s^2}p(x-s(x)\cdot 1_A+z)<\e
$$
This is true for any point $s$ from  the neighbourhood $U$ of the point $t$, and this is what we had to prove. \epr

In the following theorem we describe the topology on the cotangent bundle
$T^\star(A)=\bigsqcup_{t\in\Spec(A)} T^\star_t(A)$, but exactly in the same way one defines the topology on the complex tangent bundle $\C T^\star(A)=\bigsqcup_{t\in\Spec(A)} \C T^\star_t(A)$.

\btm\label{TH:kokasatelnoe-rassloenie} For any involutive stereotype algebra $A$ the direct sum of the stereotype quotient modules
$$
T^\star(A)=\bigsqcup_{t\in\Spec(A)} T^\star_t(A)=\bigsqcup_{t\in\Spec(A)} \Real\Big(
I_t/\ \overline{I_t^2}\ \Big)^\triangledown
$$
has a unique topology that turns the projection
$$
\pi:T^\star(A)\to\Spec(A),\qquad \pi(x+\overline{I_t^2})=t,\qquad
t\in\Spec(A),\ x\in A
$$
into a locally convex bundle with the system of seminorms $\{p';\ p\in{\mathcal P}(A)\}$, for which the mapping $$
x\in A\mapsto T^\star(x)\in\Sec(\pi)\quad\Big|\quad
T^\star(x)(t)=x+\overline{I_t^2},\quad t\in\Spec(A),
$$
continuously maps $A$ into the stereotype space $\Sec(\pi)$ of continuous sections of $\pi$. Herewith the base of the topology in $T^\star(A)$ are the sets
$$
W(x,U,p,\e)=\set{\xi\in T^\star(A): \ \pi(\xi)\in U\ \& \
p_{\pi(\xi)}\Big(\xi-T^\star(x)\big(\pi(\xi)\big)\Big)<\e }
$$
where $x\in A$, $p\in{\mathcal P}(A)$, $\e>0$, and $U$ is an open set in $M$.
 \etm

\bit{

\item[$\bullet$] The locally convex bundle $\pi: T^\star(A)\to \Spec(A)$ is called the {\it cotangent bundle} of the algebra $A$.

}\eit

\bpr Lemma \ref{LM:poluneprer-sverhu-p^*_t(T^*(x)(t))} and Proposition
\ref{PROP:sushestv-topologii-v-Xi} imply the existence and uniqueness of a topology on $T^\star(A)$, for which the projection $\pi:T^\star(A)\to\Spec(A)$ is a locally convex bundle with seminorms $p'$, and the sections of the form $T^\star(x)$, $x\in X$, are continuous. The continuity of the mapping $x\in
A\mapsto T^\star(x)\in\Sec(\pi)$ is proved by the implication
$$
p'_t(\pi(x)(t))=\inf_{y\in I_t^2}p(x+y)\le p(x)\quad \Longrightarrow \quad
\sup_{t\in T}p'_t(\pi^n(x)(t))\le p(x)
$$
for each compact set $T\subseteq M$. The structure of the base of the topology in $T^\star(A)$ also follows from Proposition \ref{PROP:sushestv-topologii-v-Xi}. \epr

\paragraph{Tangent bundle $T[A]$.}

Let $\ph:A\to B$ be a homomorphism of involutive stereotype algebras. A $\C$-linear mapping $D:A\to B$ is called a  {\it derivative} from $A$ into $B$ with respect to $\ph$, if the following identity holds:
 \beq\label{DEF:differentsirovanie}
D(x^\bullet)=D(x)^\bullet,\qquad D(x\cdot y)=D(x)\cdot\ph(y)+\ph(x)\cdot D(y),\qquad x,y\in
A.
 \eeq
The set of all derivatives from $A$ into $B$ with respect to $\ph$ is denoted by
$\Der_{\ph}(A,B)$. In the special case when $A=B$ and $\ph=\id_A$, we speak about derivatives in $A$, and we use the notation
$$
\Der(A):=\Der_{\id_A}(A,A).
$$
The spaces $\Der_{\ph}(A,B)$ and $\Der(A)$ are endowed with the topologies of immediate subspaces in $\Mor(A,B)$ and $\Mor(A,A)$.

\btm\label{TH:differentsirovanie->rassl-struuj} If the set of values of a derivative $D:A\to B$ is contained in the center of the algebra $B$,
 \beq\label{D(A)-subseteq-Z(A)}
D(A)\subseteq Z(B),
 \eeq
then $D$ is a differential operator of order 1. As a corollary, in this case $D$ defines a morphism of the jet bundles $\jet_1[D]:\Jet_A^1(A)\to \Jet_A^0(B)$, that satisfies the identity
 \beq\label{j^0(Dx)=j_1[D]-circ-j^1(x)}
\jet^0(Dx)=\jet_1[D]\circ \jet^1(x),\qquad x\in A.
 \eeq
 \etm
\bpr From the identity \eqref{DEF:differentsirovanie} we have
$$
[D,a](x)=D(a\cdot x)-a\cdot D(x)=D(a)\cdot x+a\cdot D(x)-a\cdot D(x)=D(a)\cdot
x,
$$
and this implies
$$
[[D,a],b](x)=[D,a](b\cdot x)-b\cdot [D,a](x)=D(a)\cdot b\cdot x-b\cdot
D(a)\cdot x=\eqref{D(A)-subseteq-Z(A)}=0.
$$
Thus, $D$ is a differential operator of order 1. After that we apply Theorem \ref{TH:diff-oper->rassl-struuj}.
 \epr

\btm
Each derivative $D:A\to A$ of the algebra $A$ defines a tangent vector in each point $t\in\Spec(A)$ by the formula
$$
D^T(t)(x)=t(D(x)),\qquad x\in A.
$$
\etm
\bpr
Indeed,
$$
D^T(t)(a\cdot b)=t(D(a\cdot b))=t(Da\cdot b+a\cdot Db)=t(Da)\cdot
t(b)+t(a)\cdot t(Db)=D^T(t)(a)\cdot t(b)+t(a)\cdot D^T(t)(b).
$$
\epr

From Proposition \ref{PROP:o-dvoistvennom-rasloenii} we have

\btm\label{kasatelnoe-rassloenie} Suppose an involutive stereotype algebra $A$ has a Hausdorff spectrum $\Spec(A)$, and the set $\Der(A)$ of derivatives has a dense trace in each tangent space:
$$
\overline{\{D_t;\ D\in\Der(A)\}}=T_t[A].
$$
Then the direct sum of stereotype spaces
$$
T[A]=\bigsqcup_{t\in\Spec(A)}T_t[A]
$$
has a topology such that the projection
$$
T_A:T[A]\to\Spec(A),\qquad T_A(\tau)=t,\ \tau\in T_t[A]
$$
is turned into the dual bundle to the cotangent bundle $T^\star[A]$. At the same time the mapping
$$
D\in\Der(A)\mapsto D^T\in\Sec(T[A])\quad\Big|\quad D^T(s)(a)=s(Da),\qquad a\in
A,\quad s\in\Spec(A),
$$
continuously maps $\Der(A)$ into the stereotype space $\Sec(T[A])$ of continuous sections of the bundle $T[A]$.  \etm

\bit{

\item[$\bullet$] The locally convex bundle $T_A:T[A]\to\Spec(A)$ is called the {\it tangent bundle} of the involutive stereotype algebra $A$.

}\eit

\paragraph{The Nachbin theorem.}

Let ${\mathcal E}(M)$ be the algebra of complex valued smooth functions on a smooth manifold $M$, endowed with the usual topology of uniform convergence on compact sets with respect to each differential operator. The following variant of the Nachbin theorem on the subalgebras in ${\mathcal E}(M)$ (\cite{Nachbin}, see also the monohraph \cite{Llavona}) plays in differential geometry the same role as the Stone-Weierstrass theorem does in topology.

\btm\label{TH:Nachbin} Let $A$ be an involutive stereotype subalgebra in the algebra ${\mathcal E}(M)$ of smooth functions on a smooth manifold $M$. Then $A$ is dense in ${\mathcal E}(M)$ if and only if the following two conditions hold:
 \bit{
\item[(i)] $A$ separates the points of $M$: if the points $s,t\in M$ are different, then there is a function $a\in A$ that takes different values in them,
 $$
 \forall s\ne t\in M\qquad \exists a\in A\qquad a(s)\ne a(t),
 $$

\item[(ii)] $A$ separates the tangent vectors of $M$: for each point $s\in M$ and for each nonzero tangent vector
$\tau\in T_s(M)$ there is an element $a\in A$ that has nonzero value on $\tau$,
 $$
 \forall \tau\in T_s(M)\setminus\{0\}\qquad \exists a\in A\qquad \tau(a)\ne 0,
 $$

 }\eit
 \etm

\section{Continuous envelopes and continuous duality}

As we told in $\S 0$, the material of this section sufficiently intersects the J.~N.~Kuznetsova paper \cite{Kuznetsova}.

\subsection{$C^*$-seminorms and $C^*$-algebras.}

\paragraph{$C^*$-seminorms.} Let $A$ be an involutive stereotype algebra. A seminorm $p:A\to\R_+$ is called a {\it $C^*$-seminorm}, if
\beq\label{DEF:C^*-polunorma}
p(a^\bullet\cdot a)=p(a)^2,\qquad a\in A.
\eeq
By the famous Z.~Sebesty\'en theorem \cite{Sebestyen}, every such a seminorm is submultiplicative:
$$
p(a\cdot b)\le p(a)\cdot p(b),\qquad a,b\in A.
$$
The set of all continuous $C^*$-seminorms on $A$ will be denoted by ${\mathcal P}(A)$.

\btm\label{TH:p(a)=norm(pi(a))}
Each continuous $C^*$-seminorm $p$ on $A$ can be represented as a norm of some unitary (norm-continuous\footnote{See definition on page \pageref{DEF:nepr-po-norme--unit-predst}.}) representation $\pi:A\to{\mathcal B}(X)$:
\beq\label{p(a)=norm(pi(a))}
p(a)=\norm{\pi(a)},\qquad a\in A.
\eeq
\etm
\bpr
The space $\Ker p=\{x\in A:\ p(x)=0\}$ is a closed ideal in $A$. The seminorm $p$ can be factored through a $C^*$-norm $p'$ on the quotient algebra $A/\Ker p$. Let $B$ be a completion of $A/\Ker p$ with respect to $p'$. Certainly, $B$ is a $C^*$-algebra, so it can be isometrically embedded into some algebra of the form ${\mathcal B}(X)$, where $X$ is a Hilbert space \cite[Theorem 3.4.1]{Murphy}. The composition of the maps $A\to A/\Ker p\to B\to {\mathcal B}(X)$ is the representation $\pi$ we look for.
\epr

\paragraph{Definition of continuous envelope and functoriality.}

Let us denote by ${\tt C}^*$ the class of all $C^*$-algebras.

 \bit{
\item A {\it continuous envelope}\label{DEF:nepr-obolochka} $\env_{\mathcal C} A:A\to\Env_{\mathcal C} A$ of an involutive stereotype algebra $A$ is its envelope in the class $\DEpi$ of dense epimorphisms in the category $\InvSteAlg$ of involutive stereotype algebras with respect to the class of all homomorphisms into $C^*$-algebras:
$$
\Env_{\mathcal C} A=\Env_{{\tt C}^*}^{\DEpi}A
$$
 }\eit

In detail, a {\it continuous extension} of an involutive stereotype algebra $A$ is a dense epimorphism $\sigma:A\to A'$ of involutive stereotype algebras such that for each $C^*$-agebra $B$ and for any involutive continuous homomorphism $\ph:A\to B$ there is a (necessarily, unique) involutive continuous homomorphism $\ph':A'\to B$ such that the following diagram is commutative:
\beq\label{DEF:diagr-nepr-rasshirenie}
 \xymatrix @R=2pc @C=1.2pc
 {
  A\ar[rr]^{\sigma}\ar[dr]_{\ph} & & A'\ar@{-->}[dl]^{\ph'} \\
  & B &
 }
\eeq
A {\it continuous envelope} of an involutive stereotype algebra $A$ is defined as a continuous extension $\rho:A\to \Env_{\mathcal C} A$ such that for any continuous extension $\sigma:A\to A'$ there is a (necessarily, unique) homomorphism of involutive stereotype algebras $\upsilon:A'\to \Env_{\mathcal C} A$ such that the following diagram is commutative:
$$
 \xymatrix @R=2pc @C=1.2pc
 {
  & A\ar[ld]_{\sigma}\ar[rd]^{\rho} &   \\
  A'\ar@{-->}[rr]_{\upsilon} &  & \Env_{\mathcal C} A
 }
$$

\btm\label{TH:Env_C-reg-obolochka}
The continuous envelope $\Env_{\mathcal C}$ is regular and coherent with the projective tensor product\footnote{See definitions on pages \pageref{DEF:reg-obolochka} and \pageref{DEF:obolochka-soglasovana-s-tenz-proizv}.} $\circledast$ in $\InvSteAlg$.
\etm
\bpr
1. First, let us prove the regularity, i.e. the validity of conditions R.1 - R.5 of Theorem \ref{TH:reg-obolochka}. Denote by $\varPhi=\Mor(\InvSteAlg,{\tt C}^*)$ the class of homomorphisms into $C^*$-algebras.
\bit{
\item[R.1:] The category $\InvSteAlg$ of involutive stereotype algebras is complete, in particular, projectively complete (each functor from a small category into $\InvSteAlg$ has projective limit).
\item[R.2:] By Theorem \ref{TH:SMono-circledcirc-DEpi=InvSteAlg}, the class $\DEpi$ of dense epimorphisms (where we construct the envelope) is monomorphically complement in the category $\InvSteAlg$.
\item[R.3:] The category $\InvSteAlg$ is co-well-powered in the class $\DEpi$ (since ${\tt Ste}$ is co-well-powered in the class $\Epi$).
\item[R.4:] For each involutive stereotype algebra $A$, there always exists a morphism $\ph:A\to B$ into some $C^*$-algebra $B$ (for example, one can take $B=0$ and $\ph=0$). By definition on page \pageref{DEF:goes-from}, this means that the class $\varPhi$ of morphisms into $C^*$-algebras goes from the category $\InvSteAlg$. Apart from that this class is a right ideal in $\InvSteAlg$ (the composition $\ph\circ\psi$ of any morphism $\ph:A\to B$ into a $C^*$-algebra $B$ with any morphism $\psi:A'\to A$ is also a morphism into a $C^*$-algebra).

\item[R.5:] If $\psi\circ\sigma\in\varPhi$ and $\sigma\in\DEpi$, then the composition $\psi\circ\sigma:A\to C$ takes values in a $C^*$-algebra $C$, hence $\psi$ also takes values in a $C^*$-algebra, therefore, $\psi\in\varPhi$. This means that the class $\DEpi$ pushes the class $\varPhi$.
 }\eit

2. Let us verify the coherency with the tensor product $\circledast$, i.e. the conditions T.1 and T.2 on page \pageref{DEF:obolochka-soglasovana-s-tenz-proizv}.
\bit{
\item[T.1:] Let $\rho:A\to A'$ and $\sigma:B\to B'$ be continuous extensions. Take a homomorphism $\ph:A\circledast B\to C$ into a $C^*$-algebra $C$. By Lemma \ref{LM:ph:A-circledast-B->C} it can be represented in the form
$$
\ph(a\circledast b)=\alpha(a)\cdot\beta(b),\qquad a\in A,\ b\in B,
$$
where $\alpha:A\to C$ and $\beta:B\to C$ are morphisms of stereotype algebras with commuting values:
$$
\alpha(a)\cdot\beta(b)=\beta(b)\cdot\alpha(a),\qquad a\in A,\ b\in B.
$$
Let us extend $\alpha$ and $\beta$ to some morphisms $\alpha'$ and $\beta'$ in such a way that the following diagrams are commutative:
$$
\xymatrix @R=2.pc @C=2.0pc % @M=14pt
{
A\ar[dr]_{\alpha}\ar[rr]^{\rho} & & A'\ar@{-->}[dl]^{\alpha'}\\
& C &
}
\qquad
\xymatrix @R=2.pc @C=2.0pc % @M=14pt
{
B\ar[dr]_{\beta}\ar[rr]^{\sigma} & & B'\ar@{-->}[dl]^{\beta'}\\
& C &
}
$$
Since $\rho$ and $\sigma$ are dense epimorphisms, the values of $\alpha'$ and $\beta'$ also commute:
$$
\alpha'(a)\cdot\beta'(b)=\beta'(b)\cdot\alpha'(a),\qquad a\in A',\ b\in B'.
$$
Hence we again can use Lemma \ref{LM:ph:A-circledast-B->C} and define a morphism
$$
\ph'(a\circledast b)=\alpha'(a)\cdot\beta'(b),\qquad a\in A',\ b\in B'.
$$
Obviously, it is an extension of $\ph$ (and it is unique, since $\rho$ and $\sigma$ are dense epimorphisms).

\item[T.2:] The identity map $1_{\C}:\C\to\C$ is a continuous extension (since, for instance, it is an isomorphism). On the other hand, if $\rho:\C\to A'$ is another continuous extension, then, since it is dense, it must be surjective. At the same time, $A'$ cannot be zero, since otherwise we could build the diagram
$$
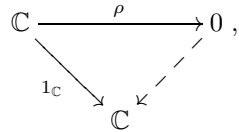
\xymatrix @R=2.pc @C=2.0pc % @M=14pt
{
\C\ar[rr]^{\rho}\ar[dr]_{1_{\C}} & & 0\ar@{-->}[dl] \\
& \C &
},
$$
what is impossible. We see that $\rho:\C\to A'$ must be an isomorphism of algebras, and we obtain the diagram
$$
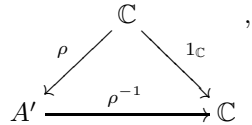
\xymatrix @R=2.pc @C=2.0pc % @M=14pt
{
& \C\ar[dl]_{\rho}\ar[dr]^{1_{\C}} &  \\
A'\ar[rr]^{\rho^{-1}} & & \C
},
$$
which means that the extension $\rho$ is embedded into the extension $1_{\C}$.
}\eit
\epr

\bcor\label{COR:Env_C-idemp-funktor}
The continuous envelope can be defined as an idempotent covariant functor from $\InvSteAlg$ into $\InvSteAlg$: there exist
\bit{
\item[1)] a map $A\mapsto (\Env_{\mathcal C}A,\env_{\mathcal C}A)$, which assigns to each involutive stereotype algebra $A$ an involutive stereotype algebra $\Env_{\mathcal C}A$ and a morphism of stereotype algebras $\env_{\mathcal C}A:A\to\Env_{\mathcal C}A$, which a continuous envelope of the algebra $A$, and

\item[2)] a map $\ph\mapsto \Env_{\mathcal C}(\ph)$, which assigns to each morphism of involutive stereotype algebras $\ph:A\to B$ a morphism of involutive stereotype algebras $\Env_{\mathcal C}(\ph):\Env_{\mathcal C}A\to \Env_{\mathcal C}B$, such that the following diagram is commutative
\beq\label{DIAGR:funktorialnost-env_C}
\xymatrix @R=2.pc @C=5.0pc % @M=14pt
{
A\ar[d]^{\ph}\ar[r]^{\env_{\mathcal C}A} & \Env_{\mathcal C}A\ar@{-->}[d]^{\Env_{\mathcal C}(\ph)} \\
B\ar[r]^{\env_{\mathcal C}B} & \Env_{\mathcal C}B \\
}
\eeq
}\eit
and the following adintities hold:
\beq\label{funktorialnost-env_C-v-Ste^circledast}
\Env_{\mathcal C}(1_A)=1_{\Env_{\mathcal C}A},\quad \Env_{\mathcal C}(\beta\circ\alpha)=\Env_{\mathcal C}(\beta)\circ \Env_{\mathcal C}(\alpha),
\eeq
\beq\label{funktorialnost-env_C-v-Ste^circledast-2}
\Env_{\mathcal C}(\Env_{\mathcal C}A)=\Env_{\mathcal C}A,\quad \env_{\mathcal C}\Env_{\mathcal C}A=1_{\Env_{\mathcal C}A},
\eeq
\beq\label{Env_C(C)=C}
\Env_{\mathcal C}\C=\C
\eeq
\ecor

\paragraph{The net of $C^*$-quotient mappings.}
The construction of continuous envelope can be described more visually in the following way. Let us call a neighbourhood of zero $U$ in an involutive stereotype algebra $A$ a {\it $C^*$-neighbourhood of zero}\label{DEF:C*-neighbourhood-of-zero}, if it is an inverse image of the unit ball under a (continuous) homomorphism $D:A\to B$ into a $C^*$-algebra $B$:
$$
U=\{x\in A:\ \norm{D(x)}\le 1\}
$$
Equivalently, $U$ must be the unit ball of some (continuous) $C^*$-seminorm $p$ on $A$:
\beq\label{DEF:C^*-okrestnost-kak-shar-C^*-polunormy}
U=\{x\in A:\ p(x)\le 1\}.
\eeq
The kernel
$$
\Ker U=\bigcap_{\lambda>0}\lambda\cdot U
$$
of any $C^*$-neighbourhood of zero $U$ in $A$ coincides with the kernel of the homomorphism $D$, hence it is a closed ideal in $A$. Consider the quotient algebra $A/\Ker U$ and endow it with the norm, for which $U+\Ker U$ is the unit ball. This algebra $A/\Ker U$ is a subalgebra in $B$ with the norm induced from $B$, or, what is the same, with the norm generated by the seminorm $p$ from \eqref{DEF:C^*-okrestnost-kak-shar-C^*-polunormy}. We denote by $A/U$ or by $A/p$ the completion of $A/\Ker U$ with respect to this norm,
\beq\label{A/U=(A/Ker U)^blacktriangledown}
A/U=A/p=(A/\Ker U)^\blacktriangledown
\eeq
and we call this the {\it quotient algebra of $A$ by the $C^*$-neighbourhood of zero $U$} or {\it by the  $C^*$-seminorm $p$}. Certainly, $A/U=A/p$ is a $C^*$-algebra (and we can think that this is a closed subalgebra in $B$). The corresponding mapping
\beq\label{DEF:pi_U}
\pi_U=\pi_p:A\to A/U=A/p
\eeq
is called the {\it quotient mapping} of the algebra $A$ by  the $C^*$-neighbourhood of zero $U$, or by the  $C^*$-seminorm $p$, or a {\it $C^*$-quotient mapping} of $A$.

\blm\label{LM:ph=ph_U-circ-pi_U-0}
For each homomorphism $\ph:A\to B$ into a $C^*$-algebra $B$ there is a $C^*$-neighbourhood of zero $U\subseteq A$ and a homomorphism $\ph_U:A/U\to B$ such that the following diagram is commutative:
\beq\label{ph=ph_U-circ-pi_U-0}
 \xymatrix @R=2pc @C=1.2pc
 {
  A\ar[rr]^{\ph}\ar[dr]_{\pi_U} & & B  \\
  & A/U\ar@{-->}[ur]_{\ph_U} &
  }
\eeq
\elm

\blm
If $U$ and $U'$ are two $C^*$-neighbourhoods of zero in $A$, and $U\supseteq U'$, then there is a unique homomorphism $\varkappa^{U'}_U:A/U\gets A/U'$ such that the following diagram is commutative:
\beq\label{varkappa^U'_U-0}
 \xymatrix @R=2pc @C=1.2pc
 {
  & A\ar[ld]_{\pi_U}\ar[rd]^{\pi_{U'}} &   \\
  A/U &  & A/ U'\ar@{-->}[ll]^{\varkappa^{U'}_U}
 }
\eeq
\elm

\blm\label{LM:diff-okr-nulya-uporyadocheny} The intersection $U\cap U'$ of any two $C^*$-neighbourhoods of zero $U$ and $U'$ in $A$ is a $C^*$-neighbourhood of zero.
\elm
\bpr
Indeed, suppose $D:A\to B$ and $D':A\to B'$ are homomorphisms that generate $U$ and $U'$.
$$
U=\{x\in A:\ \norm{D(x)}\le 1\},\qquad U'=\{x\in A:\ \norm{D'(x)}\le 1\}.
$$
Consider the $C^*$-algebra $B\oplus B'$ with the norm
$$
\norm{b\oplus b'}=\max\{\norm{b},\norm{b'}\},\qquad b\in B,\quad b'\in B'.
$$
The mapping
$$
D'':A\to B\oplus B'\quad\Big|\quad D''(x)=D(x)\oplus D'(x), \quad x\in A,
$$
is a homomorphism, and for each $x\in A$ we have
$$
x\in U\cap U'\quad\Longleftrightarrow\quad \norm{D(x)}\le 1\ \&\ \norm{D"(x)}\le 1\quad\Longleftrightarrow\quad
\norm{D''(x)}=\max\{\norm{D(x)},\norm{D"(x)}\}\le 1.
$$
\epr

\btm
The system $\pi_U:A\to A/U$ of $C^*$-quotient mapping is a net of epimorphisms\footnote{See definition on page \pageref{DEF:set-epimorf}.} in the category $\InvSteAlg$ of involutive stereotype algebras, i.e. has the following properties:
  \bit{
\item[(a)] each algebra $A$ has at least one $C^*$-neighbourhood of zero $U$, and the set of all $C^*$-neighbourhoods of zero in $A$ is directed by the semi-order
$$
U\le U'\quad\Longleftrightarrow\quad U\supseteq U',
$$

\item[(b)] for each algebra $A$ the system of morphisms $\varkappa_U^{U'}$ from \eqref{varkappa^U'_U-0} is covariant, i.e. for any three neighbourhoods of zero $U\supseteq U'\supseteq U''$ the following diagram is commutative,
$$
 \xymatrix @R=2pc @C=1.2pc
 {
  A/U &  & A/ U''\ar[ll]_{\varkappa^{U''}_U}\ar[dl]^{\varkappa^{U''}_{U'}}\\
  & A/U \ar[ul]^{\varkappa^{U'}_U} &
 }
$$
and this system $\varkappa_U^{U'}$ has a projective limit in $\InvSteAlg$;

\item[(c)] for any homomorphism $\alpha:A\gets A'$ in $\InvSteAlg$ and for any $C^*$-neighbourhood of zero $U$ in $A$ there is a $C^*$-neighbourhood of zero $U'$ in $A'$ and a homomorphism $\alpha_U^{U'}:A/U\gets A'/U'$ such that the following diagram is commutative:
 \beq\label{DIAGR:set-C^*} \xymatrix @R=2.5pc @C=4.0pc {
 A\ar[d]_{\pi_U} & A'\ar@{-->}[d]^{\pi_{U'}}\ar[l]_{\alpha} \\
A/U & A'/U'\ar@{-->}[l]^{\alpha_U^{U'}}
 } \eeq
 }\eit
\etm

By the condition (b) of this theorem, there is a projective limit $\leftlim_{0\gets U'} A/ U'$ of the system $\varkappa_U^{U'}$. As a corollary, there is a unique arrow $\pi:A\to \leftlim_{0\gets U'} A/ U'$ in $\InvSteAlg$, such that the following diagram is commutative:
\beq\label{A-to-leftlim-A/U-0}
 \xymatrix @R=2pc @C=1.2pc
 {
  & A\ar[ld]_{\pi_U}\ar@{-->}[rd]^{\pi_{U'}} &   \\
  A/U &  & \leftlim_{0\gets U'} A/ U'\ar[ll]^{\varkappa_U}
 }
\eeq
The image $\pi(A)$ of the mapping $\pi$ is (an involutive subalgebra and) a subspace in the stereotype space $\leftlim_{0\gets U'} A/ U'$. Hence it generates an immediate subspace in $\leftlim_{0\gets U'} A/ U'$, or the envelope $\Env\pi(A)$ \cite{Akbarov-env}, which is the biggest stereotype subspace in $\leftlim_{0\gets U'} A/ U'$ that has $\pi(A)$ as a dense subspace. Let us denote by $\rho:A\to \Env\pi(A)$ the lifting of the morphism $\pi$ in $\Env\pi(A)$.

\btm\label{TH:opisanie-nepr-obolochki} The morphism $\rho:A\to \Env\pi(A)$ is a continuous envelope of the algebra $A$:
$$
\Env\pi(A)=\Env_{\mathcal E}A.
$$
\etm
\bpr
The system of $C^*$-quotient mappings $\pi_U:A\to A/U$ generates on the inside the class $\varPhi$ of homomorphisms into  $C^*$-algebras by Lemma  \ref{LM:ph=ph_U-circ-pi_U-0} (see also \cite[Lemma 3.5]{Akbarov-env}):
    $$
    {\mathcal N}\subseteq\varPhi\subseteq\Mor(\SteAlg)\circ {\mathcal N}.
    $$
On the other hand, by Theorem \ref{TH:SMono-circledcirc-DEpi=InvSteAlg}, the class $\DEpi$ of all dense epimorphisms is monomorphically complemented in the category $\InvSteAlg$. Thus, by Theorem  \ref{TH:funktorialnost-pri-seti-Epi-i-dolonyaemosti}, the envelopes of $A$ with respect to the classes $\varPhi$ and ${\mathcal N}$ coincide with the morphism $\rho$.
\epr

\brem
The continuous envelope $\rho:A\to\Env_{\mathcal E}A$ is a composition of the elements $\red_\infty$ and $\coim_\infty$ of the nodal decomposition of the morphism $\pi:A\to\leftlim_{0\gets U'} A/ U'$ in the category $\tt Ste$ of stereotype spaces (not algebras!):
  \beq\label{C-envelope=im_infty-lim N_X}
\red_\infty\pi\circ\coim_\infty\pi=\env_{\mathcal C} A.
  \eeq
Visually this can be presented by the diagram
  \beq\label{DIAGR:C-envelope=im_infty-lim N_X}
\xymatrix @R=3.pc @C=6.0pc % @M=14pt
{
A\ar[d]_{\coim_\infty\pi}\ar[r]^{\pi=\leftlim_{0\gets U'}\pi_{U'}} & \leftlim_{0\gets U'} A/U' &   \\
\Coim_\infty\pi\ar[r]_{\red_\infty\pi} &  \Im_\infty \pi \ar[u]_{\im_\infty\pi}\ar@{=}[r] & \Env_{\mathcal C}A
}
 \eeq
And the algebra $\Env_{\mathcal C} A$ can be represented as an envelope (in the sense of \eqref{DEF:Env^XM}) of the set of values of the morphism $\pi$ in the stereotype space $\lim_{0\gets U'} A/ U'$:
  \beq\label{Env_EA=Env-pi(A)}
\Env_{\mathcal C} A=\Env\pi(A).
  \eeq

\erem

\subsection{Continuous algebras}

Let us say that involutive stereotype algebra $A$ is {\it continuous}, if it coincide with its continuous envelope, i.e. the continuous envelope $\env_{\mathcal C}A: A\to \Env_{\mathcal C}A$ is an isomorphism in the category $\InvSteAlg$ of involutive stereotype algebras. The class of all continuous algebras will be denoted by ${\mathcal C}\text{-}{\tt Alg}$. It forms a full subcategory in the category $\InvSteAlg$.

\paragraph{Continuous tensor product of involutive stereotype algebras.}

Let $\Env_{\mathcal C}$ be the functor of continuous envelope, defined in Corollary \ref{COR:Env_C-idemp-funktor}. For any two involutive sereotype algebras $A$ and $B$ let us define their {\it continuous tensor product} by the equality
\beq\label{C/circledast}
A\overset{\mathcal C}{\circledast} B=\Env_{\mathcal C}(A\circledast B)
\eeq
To each pair $\alpha:A\to A'$ and $\beta:B\to B'$ of morphisms of algebras we associate the morphism
\beq\label{DEF:alpha-C/circledast-beta}
\alpha\overset{\mathcal C}{\circledast}\beta=\Env_{\mathcal C}(\alpha\circledast\beta):
A\overset{\mathcal C}{\circledast} B\to A'\overset{\mathcal C}{\circledast} B'.
\eeq
Finally, to each pair of elements $a\in A$, $b\in B$ one can assign the elementary tensor
\beq\label{DEF:a-overset-C-circledast-b}
a\overset{\mathcal C}{\circledast}b=\env_{\mathcal C}(a\circledast b)
\eeq

\blm\label{LM:polnota-a-overset-C-circledast-b}
The elementary tensors $a\overset{\mathcal C}{\circledast}b$, $a\in A$, $b\in B$, are total\footnote{See definition on page \pageref{DEF:overline^X}.} in $A\overset{\mathcal C}{\circledast} B$.
\elm
\bpr
The tensors $a\circledast b$ are total in $A\circledast B$, and the image of $\env_{\mathcal C}$ is dense in $A\overset{\mathcal C}{\circledast} B$. Identity \eqref{eta(a-overset-C-circledast-b)=a-odot-b} follows from Diagram \eqref{DIAGR:eta_(A,B)}.
\epr

Further we need the following construction. For any two seminorms $q\in {\mathcal P}(A)$ and $r\in {\mathcal P}(B)$ let us consider the seminorm $q\otimes_{\max}r$ on $A\circledast B$ defined as the composition of the maps
\beq\label{DEF:q-otimes_max-r}
\xymatrix @R=1.pc @C=4.0pc
{
A\circledast B\ar@/_6ex/@{-->}[rrr]_{q\otimes_{\max}r}\ar[r]^(.4){\pi_q\circledast\pi_r} & A_q\circledast B_r=A_q\widehat{\otimes}B_r\ar[r]^{\tau} & A_q\underset{\max}{\otimes}B_r\ar[r]^(.6){\norm{\cdot}_{\max}}& \R_+
}
\eeq
where $\pi_q:A\to A_q$ and $\pi_r:B\to B_r$ are $C^*$-quotient maps defined in \eqref{DEF:pi_U}, $\pi_q\circledast\pi_r$ is their projective stereotype tensor product, $\tau$ the natural map of tensor products, and $\norm{\cdot}_{\max}$ the maximal tensor product of $C^*$-algebras.

\blm\label{LM:p(z)-le-q-otimes_max-r(z)}
Each $C^*$-seminorm $p$ on $A\circledast B$ is subrodinated to some $C^*$-seminorm $q\otimes_{\max}r$ on $A\circledast B$:
\beq\label{p(z)-le-q-otimes_max-r(z)}
p(x)\le (q\otimes_{\max}r)(x),\qquad x\in A\circledast B
\eeq
\elm
\bpr
Put $C=(A\circledast B)/p$. Then $p$ is the composition
$$
\xymatrix @R=1.pc @C=4.0pc
{
A\circledast B\ar@/_6ex/@{-->}[rr]_{p}\ar[r]^(.4){\pi_p} & (A\circledast B)/p=C\ar[r]^(.6){\norm{\cdot}_C}& \R_+
}
$$
where $\pi_p$ is the projection from \eqref{DEF:pi_U}, and $\norm{\cdot}_C$ the norm on $C$. Consider the seminorms
$$
q(a)=p(a\circledast 1_B),\qquad r(b)=p(1_A\circledast b),\qquad a\in A, \ b\in B.
$$
Let us show that the homomorphism $\pi_p:A\circledast B\to C$ is extended to some homomorphism $\pi_{q,r}:A_q\otimes_{\max} B_r\to C$:
\beq\label{PROOF:p(z)-le-q-otimes_max-r(z)}
\xymatrix @R=3.pc @C=4.0pc
{
A\circledast B\ar[dr]_{\pi_p}\ar[r]^(.4){\pi_q\circledast\pi_r} & A_q\circledast B_r\ar@{-->}[d]^{\rho}\ar[r]^{\tau}& A_q\otimes_{\max} B_r\ar@{-->}[dl]^{\sigma} \\
& C &
}
\eeq
For this we consider the homomorphisms
$$
\alpha:A\to C\quad\Big|\quad \alpha(a)=\pi_p(a\circledast 1_B),\qquad a\in A,
$$
$$
\beta:B\to C\quad\Big|\quad \beta(b)=\pi_p(1_A\circledast b),\qquad b\in B.
$$
Since $q$ and $r$ are restrictions of the seminorm $p$ under the maps $a\mapsto a\circledast 1_B$ and $b\mapsto 1_A\circledast b$, the homomorphisms $\alpha$ and $\beta$ can be extended to some homomorphisms $A_q\to C$ and $B_r\to C$:
$$
\xymatrix @R=2.pc @C=2.0pc
{
A\ar[dr]_{\alpha}\ar[rr]^{\pi_q} & & A_q\ar@{-->}[dl]^{\alpha_q} \\
& C &
}
\qquad
\xymatrix @R=2.pc @C=2.0pc
{
B\ar[dr]_{\beta}\ar[rr]^{\pi_r} & & B_r\ar@{-->}[dl]^{\beta_r} \\
& C &
}
$$

On the other hands, by Lemma \ref{LM:ph:A-circledast-B->C}, the images of the maps $\alpha$ and $\beta$ commute,
$$
\alpha(a)\cdot\beta(b)=\beta(b)\cdot\alpha(a),\qquad a\in A,\ b\in B,
$$
while $\pi_q$ and $\pi_r$ are dense epimorphisms. This implies that the images of $\alpha_q$ and $\beta_r$ also commute:
$$
\alpha_q(a')\cdot\beta_r(b')=\beta_r(b')\cdot\alpha_q(a'),\qquad a'\in A_q,\ b'\in B_r.
$$
As a corollary, again by Lemma \ref{LM:ph:A-circledast-B->C}, there is a homomorphism $\rho$ such that the left inner triangle in \eqref{PROOF:p(z)-le-q-otimes_max-r(z)} commutes. Then, since $A_q$ and $B_r$ are $C^*$-algebras, the homomorphism $\rho$ can be extended to a homomorphism $\sigma$ on $A_q\otimes_{\max}B_r$.

Finally, after constructing the map $\sigma$ in \eqref{PROOF:p(z)-le-q-otimes_max-r(z)}, we see that the seminorm  $z\mapsto \norm{\sigma(z)}_C$ is subordinated to the seminorm  $z\mapsto\norm{z}_{\max}$ (since any homomorphism of  $C^*$-algebras does not increase the norm, \cite[Theorem 2.1.7]{Murphy}):
$$
\norm{\sigma(z)}_C\le \norm{z}_{\max},\qquad z\in A_q\otimes_{\max}B_r.
$$
This implies \eqref{p(z)-le-q-otimes_max-r(z)}:
$$
p(x)=\norm{\pi_p(x)}_C=\norm{\sigma(\tau((\pi_q\circledast\pi_r)(x)))}_C
\le \norm{\tau((\pi_q\circledast\pi_r)(x))}_{\max}=(q\otimes_{\max}r)(x),\qquad z\in A\circledast B.
$$
\epr

\btm\label{TH:C-circledast->odot} For any two involutive stereotype algebras $A$ and $B$ thete is a unique linear continuous mapping $\eta_{A,B}:A\overset{\mathcal C}{\circledast}B\to A\odot B$ such that the following diagram is commutative,
    \beq\label{DIAGR:eta_(A,B)}
\xymatrix @R=2.pc @C=5.0pc % @M=14pt
{
A\circledast B\ar[dr]_{\env_{\mathcal C}{A\circledast B}\quad}\ar[rr]^{@_{A,B}} & & A\odot B \\
& A\overset{\mathcal C}{\circledast}B\ar[ur]_{\eta_{A,B}} & \\
}
\eeq
and the system of mappings $\eta_{A,B}:A\overset{\mathcal C}{\circledast}B\to A\odot B$ is a natural transformation of the functor $(A,B)\in{\tt InvAlg}^2\mapsto A\overset{\mathcal C}{\circledast}B\in{\mathcal C}\text{-}{\tt Alg}$ into the functor $(A,B)\in{\tt InvAlg}^2\mapsto A\odot B\in {\tt Ste}$.
\etm
\bpr
1. First, we prove that $\eta_{A,B}$ exists. Let us recall \cite[Corollary 31.15]{Fragoulopoulou}, that for $C^*$-algebras $A$ and $B$ there is a chain of morphisms
$$
\xymatrix @R=1.pc @C=1.0pc
{
A\otimes_{\pi}B\ar[r] & A\underset{\max}{\otimes}B\ar[r] & A\underset{\min}{\otimes}B\ar[r] & A\otimes_{\e}B.
}
$$
where $\underset{\max}{\otimes}$ and $\underset{\min}{\otimes}$ are maximal and minimal tensor product of $C^*$-algebras, and $\otimes_{\pi}$ and $\otimes_{\e}$ are their projective and injective tensor product as Banach spaces (the first three arrows are continuous homomorphisms of algebras, and the last one is the linear continuous mapping of Banach spaces). If we add the stereotype projective and injective tensor product we obtain the chain
$$
\xymatrix @R=1.pc @C=1.0pc
{
A\circledast B=A\otimes_{\pi}B\ar[r] & A\underset{\max}{\otimes}B\ar[r] & A\underset{\min}{\otimes}B\ar[r] & A\otimes_{\e}B\ar[r] & A\odot B.
}
$$
Denote by $\theta_{A,B}$ the composition of morphisms in this chain, that binds $A\underset{\max}{\otimes}B$ with $A\odot B$:
$$
\xymatrix @R=1.pc @C=1.0pc
{
A\circledast B=A\otimes_{\pi}B\ar[r] & A\underset{\max}{\otimes}B\ar[r]\ar@/_4ex/@{-->}[rrr]_{\theta_{A,B}} & A\underset{\min}{\otimes}B\ar[r] & A\otimes_{\e}B\ar[r] & A\odot B.
}
$$
When passing to the projective limit and using Lemma \ref{LM:p(z)-le-q-otimes_max-r(z)}, we obtain the equality
$$
\leftlim_{p\in{\mathcal P}(A\circledast B)}(A\circledast B)/p=(\text{Lemma \ref{LM:p(z)-le-q-otimes_max-r(z)}})=\leftlim_{q\in{\mathcal P}(A),\ r\in {\mathcal P}(B)} A/q\underset{\max}{\otimes}B/r,
$$
and after that, the chain of morphisms
$$
\xymatrix @R=1.pc @C=1.0pc
{
\leftlim_{p\in{\mathcal P}(A\circledast B)}\kern-11pt(A\circledast B)/p\ar@{=}[r]&\kern-20pt\leftlim_{q\in{\mathcal P}(A),\ r\in {\mathcal P}(B)}\kern-20pt A/q\underset{\max}{\otimes}B/r\ar@/^4ex/[rr]^{\leftlim\theta_{A/q,B/r}} & & \leftlim_{q\in{\mathcal P}(A),\ r\in {\mathcal P}(B)}\kern-20pt A/q\odot B/r\ar@{=}[r]& \leftlim_{q\in{\mathcal P}(A)} A/q\odot \leftlim_{r\in {\mathcal P}(B)} B/r\ar@{=}[r]& A\odot B
}
$$
(the last equality follows from the fact that $\odot$ commutes with the tensor products, see \cite[(2.53)]{Akbarov-env}). This gives the following chain of morphisms
$$
\xymatrix @R=1.pc @C=1.0pc
{
A\overset{\mathcal C}{\circledast} B\ar@{=}[r] & \Env_{\mathcal C}(A\circledast B)\ar@/^4ex/[rr]^{\leftlim_{p\in{\mathcal P}(A\circledast B)}\rho_p} & & \leftlim_{p\in{\mathcal P}(A\circledast B)}\kern-11pt(A\circledast B)/p\ar@{=}[r]&\kern-20pt\leftlim_{q\in{\mathcal P}(A),\ r\in {\mathcal P}(B)}\kern-20pt A/q\underset{\max}{\otimes}B/r\ar[rrrr]^(.6){\leftlim\theta_{A/q,B/r}} & & && A\odot B
}
$$
Their composition is $\eta_{A,B}$.

2. Now let us prove that the system of maps $\eta_{A,B}$ is a morphism of functors. Take two morphisms of algebras, $\alpha:A\to A'$ and $\beta:B\to B'$, and consider the diagram
\beq\label{PROOF:C-circledast->odot}
\xymatrix @R=2.pc @C=5.0pc % @M=14pt
{
A\circledast B\ar[dd]_{\alpha\circledast\beta}\ar[dr]_{\env_{\mathcal C}{A\circledast B}\quad}\ar[rr]^{@_{A,B}} & & A\odot B\ar[dd]_{\alpha\odot\beta} \\
& A\overset{\mathcal C}{\circledast}B\ar[dd]_(.2){\alpha\overset{\mathcal C}{\circledast}\beta}\ar[ur]_{\eta_{A,B}} & \\
A'\circledast B'\ar[dr]_{\env_{\mathcal C}{A'\circledast B'}\quad}\ar[rr]^(.7){@_{A',B'}}|!{[r]}\hole & & A'\odot B' \\
& A'\overset{\mathcal C}{\circledast}B'\ar[ur]_{\eta_{A',B'}} & \\
}
\eeq
Here the upper and the lower bases are commutative, since they are diagrams \eqref{DIAGR:eta_(A,B)}, the remote side is commutative, since this is diagram \eqref{@-preobr-functorov} of naturality of the Grothendieck transformation $@$, and the left near side is commutative, since this is the diagram of functoriality of the envelope \eqref{DIAGR:funktorialnost-env_varPhi^Epi-v-kat-s-uzl-razl-1-*}. On the other hand, from Lemma  \ref{LM:polnota-a-overset-C-circledast-b} it follows that the map $\env_{\mathcal C}{A\circledast B}$ is an epimorphism of stereotype spaces. Together this means that the right near side in diagram
\eqref{PROOF:C-circledast->odot} is commutative as well, and this is what we need.
\epr

\brem\label{REM:eta(a-overset-C-circledast-b)=a-odot-b}
Diagram \eqref{DIAGR:eta_(A,B)} implies that $\eta_{A,B}$ maps the elementary tensors $a\overset{\mathcal C}{\circledast}b$ into the elementary tensors $a\odot b$:
\beq\label{eta(a-overset-C-circledast-b)=a-odot-b}
\eta_{A,B}(a\overset{\mathcal C}{\circledast}b)=a\odot b.
\eeq
\erem

Theorems  \ref{TH:Env_C-reg-obolochka} and \ref{TH:sushestvovanie-tenz-proizv-v-L} imply

\btm\label{TH:C-obolochka=monoidalnyi-funktor} Formula \eqref{C/circledast} defines a tensor product in ${\mathcal C}\text{-}{\tt Alg}$ that turns ${\mathcal C}\text{-}{\tt Alg}$ into a monoidal category, and the functor of continuous envelope $E$ is a monoidal functor from the monoidal category $(\InvSteAlg,\circledast)$ of involutive stereotype algebras into the monoidal category $({\mathcal C}\text{-}{\tt Alg},\overset{\mathcal C}{\circledast})$ of continuous algebras. The corresponding morphisms of bifunctors
 $$
\Big((A,B)\mapsto \Env_{\mathcal C}(A)\overset{\mathcal C}{\circledast} \Env_{\mathcal C}(B)\Big)\overset{E^{\circledast}}{\rightarrowtail} \Big((A,B)\mapsto \Env_{\mathcal C}(A\circledast B)\Big)
 $$
is defined by the formula
$$
E^\circledast_{A,B}=\Env_{\mathcal C}(\env_{\mathcal C}A\circledast \env_{\mathcal C}B)^{-1}:\Env_{\mathcal C}(A)\overset{\mathcal C}{\circledast} \Env_{\mathcal C}(B)=\Env_{\mathcal C}(\Env_{\mathcal C}(A)\circledast \Env_{\mathcal C}(B))\to \Env_{\mathcal C}(A\circledast B),
$$
and the morphism $E^{\C}$ in ${\mathcal C}\text{-}{\tt Alg}$, that turns the identity object $\C$ in the category ${\mathcal C}\text{-}{\tt Alg}$ into the image $\Env_{\mathcal C}(\C)$ of the identity object $\C$ in the category $\InvSteAlg$, is the local identity:
$$
E^{\C}=1_{\C}:\C\to\C=\Env_{\mathcal C}(\C).
$$
\etm

\paragraph{Action of continuous envelope on bialgebras.}

\btm\label{LM:koalg-v-CAlg->koalg-v-odot}
If $A$ is a coalgebra in a monoidal category $({\mathcal C}\text{-}{\tt Alg},\overset{\mathcal C}{\circledast})$  of continuous algebras with the structure morphisms
$$
\varkappa:A\to A\overset{\mathcal C}{\circledast} A,\qquad \e:A\to \C,
$$
then $A$ is coalgebra in the monoidal category $({\tt Ste},\odot)$ of stereotype spaces with the structure morphisms
$$
\lambda=\eta_{A,A}\circ\varkappa:A\to A\odot A,\qquad \e:A\to\C.
$$
Every morphism $\ph:A\to B$ of coalgebras in $({\mathcal C}\text{-}{\tt Alg},\overset{\mathcal C}{\circledast})$ is a morphism of $A$ into $B$ as coalgebras in $({\tt Ste},\odot)$.
\etm
\bpr
1. Consider the associativity diagram for $\varkappa$
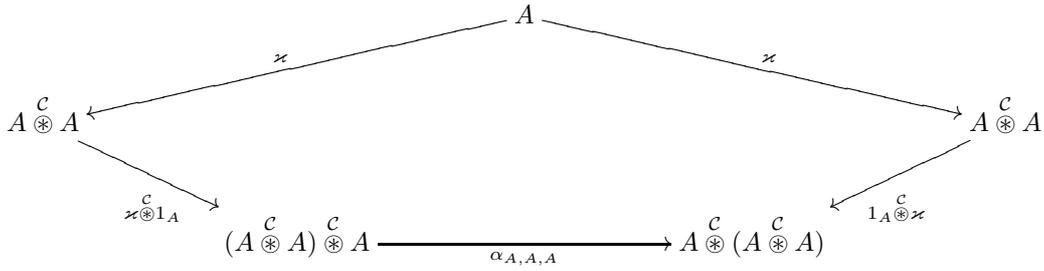
$$
\xymatrix @R=2.pc @C=4.0pc % @M=14pt
{
&& A\ar[dll]_{\varkappa}\ar[drr]^{\varkappa} && \\
A\overset{\mathcal C}{\circledast}A\ar[dr]_{\varkappa\overset{\mathcal C}{\circledast}1_A} && && A\overset{\mathcal C}{\circledast}A \ar[dl]^{1_A\overset{\mathcal C}{\circledast}\varkappa}\\
& (A\overset{\mathcal C}{\circledast}A)\overset{\mathcal C}{\circledast}A\ar[rr]_{\alpha_{A,A,A}} && A\overset{\mathcal C}{\circledast}(A\overset{\mathcal C}{\circledast}A)
}
$$
and add it to the diagram
$$
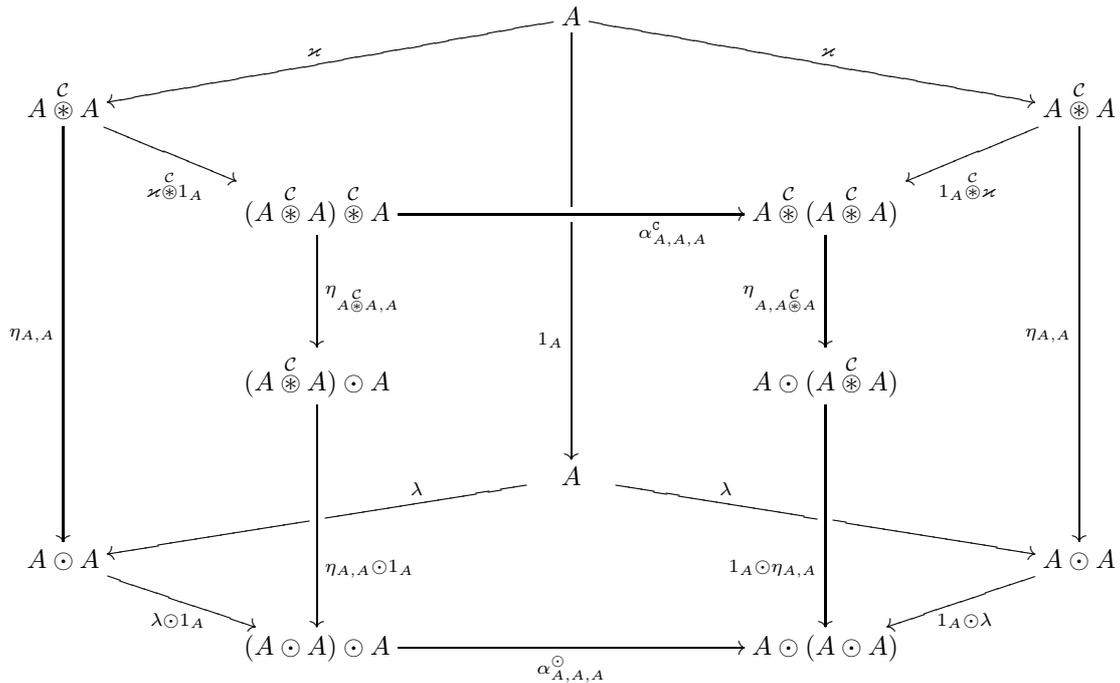
\xymatrix @R=1.5pc @C=4.0pc % @M=14pt
{
&& A\ar[dll]_{\varkappa}\ar[drr]^{\varkappa}\ar[ddddd]_(.7){1_A}|!{[dd]}\hole && \\
A\overset{\mathcal C}{\circledast}A\ar[dr]_{\varkappa\overset{\mathcal C}{\circledast}1_A}\ar[ddddd]_{\eta_{A,A}} && && A\overset{\mathcal C}{\circledast}A \ar[dl]^{1_A\overset{\mathcal C}{\circledast}\varkappa}\ar[ddddd]_{\eta_{A,A}}\\
& (A\overset{\mathcal C}{\circledast}A)\overset{\mathcal C}{\circledast}A\ar[rr]_(.7){\alpha^{\tt C}_{A,A,A}}\ar[dd]^{\eta_{A\overset{\mathcal C}{\circledast}A,A}} && A\overset{\mathcal C}{\circledast}(A\overset{\mathcal C}{\circledast}A)\ar[dd]_{\eta_{A,A\overset{\mathcal C}{\circledast}A}} & \\
&&&& \\
& (A\overset{\mathcal C}{\circledast}A)\odot A\ar[ddd]^(.7){\eta_{A,A}\odot 1_A} && A\odot (A\overset{\mathcal C}{\circledast}A)\ar[ddd]_(.7){1_A\odot\eta_{A,A}} & \\
&& \quad A\quad \ar[dll]_(.3){\lambda}|!{[d];[ll]}\hole\ar[drr]^(.3){\lambda}|!{[d];[rr]}\hole && \\
A\odot A\ar[dr]_{\lambda\odot 1_A} && && A\odot A \ar[dl]^{1_A\odot \lambda}\\
& (A\odot A)\odot A\ar[rr]_{\alpha^{\odot}_{A,A,A}} && A\odot (A\odot A) &
}
$$
Here the upper base of the prism is commutative and we have to prove the commutativity of the lower base. For this it is sufficient to verify the commutativity of the lateral faces. The two remote lateral faces are commutative just since they present the definition of morphism $\lambda$. The commutativity of the left nearby lateral face
\beq\label{DIAGR-1:koalg-v-CAlg->koalg-v-odot}
\xymatrix @R=2.pc @C=5.0pc % @M=14pt
{
A\overset{\mathcal C}{\circledast}A\ar[dd]_{\eta_{A,A}}\ar[r]^{\varkappa\overset{\mathcal C}{\circledast}1_A} & (A\overset{\mathcal C}{\circledast}A)\overset{\mathcal C}{\circledast}A
\ar[d]^{\eta_{A\overset{\mathcal C}{\circledast}A,A}} \\
 & (A\overset{\mathcal C}{\circledast}A)\odot A\ar[d]^{\eta_{A,A}\odot 1_A} \\
A\odot A\ar[r]_{\lambda\odot 1_A} & (A\odot A)\odot A
}
\eeq
can be verified on elementary tensors. Take $a,b\in A$ and let us represent $\varkappa(a)$ as a limit of a net of sums of elementary tensors (here is the first time when we use Lemma \ref{LM:polnota-a-overset-C-circledast-b}):
$$
\sum_{i\in I_s}x_i^s\overset{\mathcal C}{\circledast}y_i^s\underset{s\to\infty}{\longrightarrow}\varkappa(a)
$$
Then when moving by Diagram \eqref{DIAGR-1:koalg-v-CAlg->koalg-v-odot} the elementary tensor $a\overset{\mathcal C}{\circledast}b$ gives the following elements:
$$
\xymatrix @R=2.pc @C=2.0pc % @M=14pt
{
a\overset{\mathcal C}{\circledast}b\ar@{|->}[dd]_{\eta_{A,A}}\ar@{|->}[rrr]^{\varkappa\overset{\mathcal C}{\circledast}1_A} & & & \varkappa(a)\overset{\mathcal C}{\circledast}b\ar@{=}[r]
 & \lim\limits_{s\to\infty}(\sum\limits_{i\in I_s}x_i^s\overset{\mathcal C}{\circledast}y_i^s)\overset{\mathcal C}{\circledast}b\ar@{|->}[d]^{\eta_{A\overset{\mathcal C}{\circledast}A,A}} \\
& & & & \lim\limits_{s\to\infty}(\sum\limits_{i\in I_s}x_i^s\overset{\mathcal C}{\circledast}y_i^s)\odot b\ar@{|->}[d]^{\eta_{A,A}\odot 1_A} \\
a\odot b\ar@{|->}[rr]^(.4){\lambda\odot 1_A} & &\eta_{A,A}(\varkappa(a))\odot b\ar@{=}[r] & \eta_{A,A}\left(\lim\limits_{s\to\infty}(\sum\limits_{i\in I_s}x_i^s\overset{\mathcal C}{\circledast}y_i^s)\right)\odot b
\ar@{=}[r]& \lim\limits_{s\to\infty}(\sum\limits_{i\in I_s}x_i^s\odot y_i^s)\odot b
}
$$
Since tensors $a\overset{\mathcal C}{\circledast}b$ are total in $A\overset{\mathcal C}{\circledast}A$ (here we use Lemma \ref{LM:polnota-a-overset-C-circledast-b} second time), this proves the commutativity of \eqref{DIAGR-1:koalg-v-CAlg->koalg-v-odot}.

The same trick proves the commutativity of the right nearby lateral face
$$
\xymatrix @R=2.pc @C=5.0pc % @M=14pt
{
 A\overset{\mathcal C}{\circledast}(A\overset{\mathcal C}{\circledast}A)
\ar[d]_{\eta_{A,A\overset{\mathcal C}{\circledast}A}} & A\overset{\mathcal C}{\circledast}A\ar[dd]_{\eta_{A,A}}\ar[l]_{\varkappa\overset{\mathcal C}{\circledast}1_A} \\
A\odot(A\overset{\mathcal C}{\circledast}A)\ar[d]_{1_A\odot\eta_{A,A}} & \\
A\odot (A\odot A) & A\odot A\ar[l]_{\lambda\odot 1_A}
}
$$

For the central nearby lateral face
$$
\xymatrix @R=2.pc @C=5.0pc % @M=14pt
{
(A\overset{\mathcal C}{\circledast}A)\overset{\mathcal C}{\circledast}A\ar[r]^{\alpha^{\tt C}_{A,A,A}}
 \ar[d]_{\eta_{A\overset{\mathcal C}{\circledast}A,A}} & A\overset{\mathcal C}{\circledast}(A\overset{\mathcal C}{\circledast}A)\ar[d]_{\eta_{A,A\overset{\mathcal C}{\circledast}A}} \\
(A\overset{\mathcal C}{\circledast}A)\odot A\ar[d]_{\eta_{A,A}\odot1_A} & A\odot(A\overset{\mathcal C}{\circledast}A)\ar[d]_{1_A\odot\eta_{A,A}} \\
(A\odot A)\odot A\ar[r]^{\alpha^{\odot}_{A,A,A}} & A\odot (A\odot A)
}
$$
we have to consider a triple of elements $a,b,c\in A$. When moving by this diagram they give:
$$
\xymatrix @R=2.pc @C=5.0pc % @M=14pt
{
(a\overset{\mathcal C}{\circledast}b)\overset{\mathcal C}{\circledast}c\ar@{|->}[r]^{\alpha^{\tt C}_{A,A,A}}
 \ar@{|->}[d]_{\eta_{A\overset{\mathcal C}{\circledast}A,A}} & a\overset{\mathcal C}{\circledast}(b\overset{\mathcal C}{\circledast}c)\ar@{|->}[d]_{\eta_{A,A\overset{\mathcal C}{\circledast}A}} \\
(a\overset{\mathcal C}{\circledast}b)\odot c\ar@{|->}[d]_{\eta_{A,A}\odot1_A} & a\odot(b\overset{\mathcal C}{\circledast}c)\ar@{|->}[d]_{1_A\odot\eta_{A,A}} \\
(a\odot b)\odot c\ar@{|->}[r]^{\alpha^{\odot}_{A,A,A}} & a\odot (b\odot c)
}
$$

The diagrams for the counit are verified similarly.

2. Suppose $\ph:A\to B$ is a morphism of coalgebras in $({\mathcal C}\text{-}{\tt Alg},\overset{\mathcal C}{\circledast})$, i.e. a morphism of $A$ into $B$ as stereotype spaces, such that the following diagrams are commutative:
\beq\label{PROOF:LM:koalg-v-CAlg->koalg-v-odot}
\xymatrix @R=2.pc @C=5.0pc % @M=14pt
{
A\ar[d]_{\varkappa_A}\ar[r]^{\ph}  & B\ar[d]^{\varkappa_B} \\
A\overset{\mathcal C}{\circledast}A\ar[r]^{\ph\overset{\mathcal C}{\circledast}\ph} & B\overset{\mathcal C}{\circledast}B
}
\qquad
\xymatrix @R=2.pc @C=2.0pc % @M=14pt
{
A\ar[dr]_{\e_A}\ar[rr]^{\ph}  & & B\ar[dl]^{\e_B} \\
 & \C &
}
\eeq
where $\varkappa_A$, $\varkappa_B$, $\e_A$, $\e_B$ are structure morphisms. Then the left one of these diagrams can be complemented as follows:
$$
\xymatrix @R=2.pc @C=5.0pc % @M=14pt
{
A\ar[d]_{\varkappa_A}\ar[r]^{\ph}  & B\ar[d]^{\varkappa_B} \\
A\overset{\mathcal C}{\circledast}A
\ar[d]_{\eta_{A,A}}\ar[r]^{\ph\overset{\mathcal C}{\circledast}\ph} & B\overset{\mathcal C}{\circledast}B
\ar[d]_{\eta_{B,B}} \\
A\odot A\ar[r]_{\ph\odot\ph} & B\odot B,
}
$$
This diagram is commutative since by Theorem \ref{TH:C-circledast->odot} the morphisms $\eta$ are natural transformations of bifunctors. Together with the right diagram in \eqref{PROOF:LM:koalg-v-CAlg->koalg-v-odot} this means that $\ph:A\to B$ is a morphism of coalgebras in $(\Ste,\odot)$.
\epr

\btm\label{TH:C-obolochka-sohranyaet-Hopfov}
Let $H$ be a bialgebra in the category $({\tt Ste},\circledast)$ of stereotype spaces, or, what is the same, a coalgebra on the category ${\tt Ste}^{\circledast}$ of stereotype algebras, with the comultiplication $\varkappa$ and the counit $\e$. Then
 \bit{
\item[(i)] the continuous envelope $\Env_{\mathcal C}H$ is a coalgebra in the monoidal category $({\mathcal C}\text{-}{\tt Alg},\overset{\mathcal C}{\circledast})$ of continuous algebras with the comultiplication and the counit
    \beq\label{varkappa^E,e^E}
    \varkappa_E=\Env_{\mathcal C}(\env_{\mathcal C}H\circledast \env_{\mathcal C}H)\circ \Env_{\mathcal C}(\varkappa),\qquad \e_E=\Env_{\mathcal C}(\e),
    \eeq

\item[(ii)] the continuous envelope $\Env_{\mathcal C}H$ is a coalgebra in the monoidal category $({\tt Ste},\odot)$ of stereotype spaces with the comultiplication and the counit
    \beq\label{varkappa^odot,e^odot}
    \varkappa_\odot=\eta_{\Env_{\mathcal C}H,\Env_{\mathcal C}H}\circ \Env_{\mathcal C}(\env_{\mathcal C}H\circledast \env_{\mathcal C}H)\circ \Env_{\mathcal C}(\varkappa)=
    \eta_{\Env_{\mathcal C}H,\Env_{\mathcal C}H}\circ \varkappa_E,\qquad \e_\odot=\Env_{\mathcal C}(\e),
    \eeq

\item[(iii)] the morphism $\env_{\mathcal C}H^\star:H^\star\gets (\Env_{\mathcal C}H)^\star$, dual to the morphism of envelope $\env_{\mathcal C}H:H\to \Env_{\mathcal C}H$, is a morphism of stereotype algebras, if $(\Env_{\mathcal C}H)^\star$ is considered as an algebra with the multiplication and the unit, dual to \eqref{varkappa^odot,e^odot}, and $H^\star$ as the algebra with the multiplication and the unit
    $$
    \varkappa^\star\circ @_{H^\star,H^\star},\qquad \e^\star.
    $$

 }\eit
\etm
\bpr
An involutive bialgebra in the category $({\tt Ste},\circledast)$ is the same as a coalgebra in the category  $(\InvSteAlg,\circledast)$ of involutive stereotype algebras. Hence by theorem \ref{TH:C-obolochka=monoidalnyi-funktor} $\Env_{\mathcal C}H$ is a coalgebra in $({\mathcal C}\text{-}{\tt Alg},\overset{\mathcal C}{\circledast})$ with the comultiplication and the counit \eqref{varkappa^E,e^E}. After that we use Theorem \ref{LM:koalg-v-CAlg->koalg-v-odot}, and we obtain that $H$ is a coalgebra in the category $({\tt Ste},\odot)$ with the comultiplication and the counit \eqref{varkappa^odot,e^odot}. It remains to verify (iii). We have to note for this that the following diagram is commutative:
$$
\xymatrix @R=1.5pc @C=4.0pc % @M=14pt
{
\Env_{\mathcal C}H\ar[dr]_{\Env_{\mathcal C}(\varkappa)}\ar@/^8ex/[ddrr]^{\varkappa_E}\ar@/_30ex/[dddddr]_{\varkappa_\odot} & & \\
 H\ar[u]_{\env_{\mathcal C}H}\ar[d]^{\varkappa}  & \Env_{\mathcal C}(H\circledast H)\ar[dr]_{\Env_{\mathcal C}(\env_{\mathcal C}H\circledast \env_{\mathcal C}H)\quad\qquad}  & \\
 H\circledast H\ar[ur]_{\quad\env_{\mathcal C}{H\circledast H}}\ar[dr]^{\qquad\env_{\mathcal C}H\circledast \env_{\mathcal C}H}\ar[dd]_{@_{H,H}} & &  \Env_{\mathcal C}(\Env_{\mathcal C}H\circledast \Env_{\mathcal C}H)\ar@{=}[dd] \\
 & \Env_{\mathcal C}H\circledast \Env_{\mathcal C}H\ar[ur]_{\qquad \env_{\mathcal C}{\Env_{\mathcal C}H\circledast \Env_{\mathcal C}H}}\ar[dd]^{@_{\Env_{\mathcal C}H,\Env_{\mathcal C}H}} & \\
 H\odot H\ar[dr]^(.4){\qquad\env_{\mathcal C}H\odot \env_{\mathcal C}H} & & \Env_{\mathcal C}H\overset{\mathcal C}{\circledast} \Env_{\mathcal C}H\ar[dl]^{\qquad \eta_{\Env_{\mathcal C}H,\Env_{\mathcal C}H}} \\
 & \Env_{\mathcal C}H\odot \Env_{\mathcal C}H & \\
}
$$
After passing to the dual space we have
$$
\xymatrix @R=1.5pc @C=4.0pc % @M=14pt
{
(\Env_{\mathcal C}H)^\star\ar[d]^{\env_{\mathcal C}H^\star} & & \\
 H^\star  & \Env_{\mathcal C}(H\circledast H)^\star\ar[ul]^{\Env_{\mathcal C}(\varkappa)^\star}\ar[dl]^(.4){\env_{\mathcal C}{H\circledast H}^\star}  & \\
 H^\star\odot H^\star\ar[u]_{\varkappa^\star}  & &  \Env_{\mathcal C}(\Env_{\mathcal C}H\circledast \Env_{\mathcal C}H)^\star\ar[ul]^{\Env_{\mathcal C}(\env_{\mathcal C}H\circledast \env_{\mathcal C}H)^\star\quad\qquad} \ar@/_8ex/[uull]_{\varkappa_E^\star}\ar@{=}[dd]\ar[dl]^{\qquad \env_{\mathcal C}{\Env_{\mathcal C}H\circledast \Env_{\mathcal C}H}^\star} \\
 & (\Env_{\mathcal C}H)^\star\odot (\Env_{\mathcal C}H)^\star\ar[ul]_{\qquad\env_{\mathcal C}H^\star\odot \env_{\mathcal C}H^\star}  & \\
 H^\star\circledast H^\star\ar[uu]^{@_{H^\star,H^\star}} & & (\Env_{\mathcal C}H\overset{\mathcal C}{\circledast} \Env_{\mathcal C}H)^\star \\
 & (\Env_{\mathcal C}H)^\star\circledast (\Env_{\mathcal C}H)^\star\ar@/^30ex/[uuuuul]^{\varkappa_\odot^\star}\ar[ul]_{\qquad\env_{\mathcal C}H^\star\circledast \env_{\mathcal C}H^\star} \ar[uu]_{@_{(\Env_{\mathcal C}H)^\star,(\Env_{\mathcal C}H)^\star}}\ar[ur]_{\qquad \eta_{\Env_{\mathcal C}H,\Env_{\mathcal C}H}^\star} & \\
}
$$
The fragment in the left lower corner is exactly the diagram of coherence of the multiplications:
$$
\xymatrix @R=2.pc @C=6.0pc % @M=14pt
{
H^\star & (\Env_{\mathcal C}H)^\star\ar[l]_{\env_{\mathcal C}H^\star} \\
H^\star\odot H^\star\ar[u]^{\varkappa^\star} & \\
H^\star\circledast H^\star\ar[u]^{@_{H^\star,H^\star}} & (\Env_{\mathcal C}H)^\star\circledast (\Env_{\mathcal C}H)^\star
\ar[l]_{\env_{\mathcal C}H^\star\circledast \env_{\mathcal C}H^\star}\ar[uu]_{\varkappa_\odot^\star}
}
$$
\epr

\btm\label{TH:C-obolochka-sohranyaet-inv-Hopfov}
Let $H$ be an involutive Hopf algebra in the category $({\tt Ste},\circledast)$ of stereotype spaces. Then
 \bit{

\item[(i)] the continuous envelope $\Env_{\mathcal C}H$, as a coalgebra in the monoidal categories  $({\mathcal C}\text{-}{\tt Alg},\overset{\mathcal C}{\circledast})$ and $({\tt Ste},\odot)$, has interconsistent antipode $\Env_{\mathcal C}(\sigma)$ and involution $\Env_{\mathcal C}(\bullet)$, uniquely defined by diagrams in the category $\tt{Ste}$
    \beq\label{C-obolochka-sohranyaet-inv-Hopfov}
    \xymatrix @R=2.pc @C=4.0pc % @M=14pt
{
H\ar[d]_{\sigma}\ar[r]^{\env_{\mathcal C}H} & \Env_{\mathcal C}H\ar@{-->}[d]^{\Env_{\mathcal C}(\sigma)} \\
H\ar[r]^{\env_{\mathcal C}H} & \Env_{\mathcal C}H
}
\qquad
    \xymatrix @R=2.pc @C=4.0pc % @M=14pt
{
H\ar[d]_{\bullet}\ar[r]^{\env_{\mathcal C}H} & \Env_{\mathcal C}H\ar@{-->}[d]^{\Env_{\mathcal C}(\bullet)} \\
H\ar[r]^{\env_{\mathcal C}H} & \Env_{\mathcal C}H
}
\eeq

\item[(ii)] the morphism $\env_{\mathcal C}H^\star:H^\star\gets (\Env_{\mathcal C}H)^\star$, dual to the morphism of envelope $\env_{\mathcal C}H:H\to \Env_{\mathcal C}H$, is an involutive homomorphism of stereotype algebras over $\circledast$, if $H^\star$ and $(\Env_{\mathcal C}H)^\star$ are endowed with the structure of dual involutive algebras to the involutive coalgebras with antipode $H$ and $\Env_{\mathcal C}H$ by the property $4^\circ$ on page \pageref{4^0:inv-v-sopryazh-alg}.
 }\eit
\etm
\bpr
1. Denote by $H^{\op}$ the algebra $H$ with the opposite multiplication:
$$
\mu^{\op}=\mu\circ\br.
$$
Let $\op_H:H\to H^{\op}$ denote the identity mapping of $H$ into itself (we assume that the range here is the algebra with the opposite multiplication). This is an anti-homomorphism of algebras. The continuous envelopes of $H$ and $H^{\op}$ are also connected to each other through a natural anti-homomorphism, which we denote by $\Env_{\mathcal C}(\op)$:
\beq\label{E(op_H)}
    \xymatrix @R=2.pc @C=4.0pc % @M=14pt
{
H\ar[d]_{\op_H}\ar[r]^{\env_{\mathcal C}H} & \Env_{\mathcal C}H\ar@{-->}[d]^{\Env_{\mathcal C}(\op_H)} \\
H^{\op}\ar[r]^{\env_{\mathcal C}H} & \Env_{\mathcal C}(H^{\op})
}
\eeq
This follows from the fact  that for each $C^*$-seminorm $p$ the following diagram is commutative:
$$
    \xymatrix @R=2.pc @C=4.0pc % @M=14pt
{
H\ar[d]_{\op_H}\ar[r]^{\pi_p} & H/p\ar@{-->}[d]^{\op_{H/p}} \\
H^{\op}\ar[r]^{\pi_p} & (H/p)^{\op}
}
$$
As a corollary there is a unique anti-homomorphism between projective limits:
$$
    \xymatrix @R=2.pc @C=4.0pc % @M=14pt
{
H\ar[d]_{\op_H}\ar[r]^{\leftlim_p\pi_p} & \leftlim_p H/p\ar@{-->}[d] \\
H^{\op}\ar[r]^{\leftlim_p\pi_p} & \leftlim_p (H/p)^{\op}
}
$$
Then we pass to the immediate subspaces generated by the images of $H$ and $H^{\op}$, and we obtain the dotted arrow in  \eqref{E(op_H)}.

When $\Env_{\mathcal C}(\op_H)$ is already defined for each stereotype algebra $H$, the mapping $\Env_{\mathcal C}(\sigma)$ can be defined by the formula
$$
\Env_{\mathcal C}(\sigma)=\Env_{\mathcal C}(\op_{H^{\op}})\circ \Env_{\mathcal C}(\op_H\circ\sigma)
$$
or by the diagram
$$
    \xymatrix @R=2.pc @C=8.0pc % @M=14pt
{
H\ar[d]^{\sigma}\ar[r]^{\env_{\mathcal C}H}\ar@/_7ex/[ddd]_{\sigma} & \Env_{\mathcal C}H\ar[dd]_{\Env_{\mathcal C}(\op_H\circ\sigma)}\ar@{-->}@/^7ex/[ddd]^{\Env_{\mathcal C}(\sigma)} \\
H\ar[d]^{\op_H} &  \\
H^{\op}\ar[d]^{\op_{H^{\op}}}\ar[r]^{\env_{\mathcal C}{H^{\op}}} & \Env_{\mathcal C}(H^{\op})\ar[d]_{\Env_{\mathcal C}(\op_{H^{\op}})}\\
H\ar[r]^{\env_{\mathcal C}{H}} & \Env_{\mathcal C}H
}
$$

2. The existence of $\Env_{\mathcal C}(\bullet)$ is proved by the same trick as the existence of $\Env_{\mathcal C}(\op_H)$. First for arbitrary $C^*$-seminorm $p$ we notice the diagram
$$
    \xymatrix @R=2.pc @C=4.0pc % @M=14pt
{
H\ar[d]_{\bullet}\ar[r]^{\pi_p} & H/p\ar@{-->}[d]^{\bullet} \\
H\ar[r]^{\pi_p} & H/p
}
$$
And then we pass to the projective limit and to the immediate subspaces generated by the image of $H$, and we obtain the right diagram in \eqref{C-obolochka-sohranyaet-inv-Hopfov}.

3. For (ii) we have only to verify that $\env_{\mathcal C}H^\star$ preserves involution. For $a\in H$ and $f\in (\Env_{\mathcal C}H)^\star$ we have:
\begin{multline*}
\env_{\mathcal C}H^\star(f^\bullet)(a)=f^\bullet(\env_{\mathcal C}H(a))=\overline{f(\Env_{\mathcal C}(\sigma)(\env_{\mathcal C}H(a))^\bullet)}=
\overline{f((\bullet\circ \Env_{\mathcal C}(\sigma)\circ \env_{\mathcal C}H)(a))}=\\=
\overline{f((\bullet\circ \env_{\mathcal C}H\circ\sigma)(a))}=
\overline{f((\env_{\mathcal C}H\circ\bullet\circ\sigma)(a))}=\overline{f(\env_{\mathcal C}H(\sigma(a)^\bullet))}=
\overline{\env_{\mathcal C}H^\star(f)(\sigma(a)^\bullet)}=\env_{\mathcal C}H^\star(f)^\bullet(a).
\end{multline*}
\epr

\paragraph{Continuous tensor product with ${\mathcal C}(M)$.}

Let $X$ be a stereotype space, and $M$ a paracompact locally compact topological space. Consider the algebra ${\mathcal C}(M)$ of continuous functions on $M$ and the space ${\mathcal C}(M,X)$ of continuous functions mappings $M$ into $X$. We endow ${\mathcal C}(M)$ and ${\mathcal C}(M,X)$ with the standard topology of uniform convergence on compact sets in $M$
$$
u_i\overset{{\mathcal C}(M,X)}{\longrightarrow} 0\quad\Longleftrightarrow\quad \forall\ \text{compact}\ K\subseteq M\quad  u_i|_K\overset{{\mathcal C}(K,X)}{\longrightarrow} 0,\qquad u\in {\mathcal C}(M,X),\quad t\in M.
$$
and with the pointwise algebraic operations:
$$
(\lambda\cdot u)(t)=\lambda\cdot u(t)\qquad (u+v)(t)=u(t)+v(t),\qquad u,v\in {\mathcal C}(M,X),\quad \lambda\in\C,\quad t\in M.
$$

From \cite[Theorem 8.1]{Akbarov} we have

\bprop The space  ${\mathcal C}(M,X)$ is a stereotype module over ${\mathcal C}(M)$, and
\beq\label{C(M,X)-cong-C(M)-odot-X}
{\mathcal C}(M,X)\cong {\mathcal C}(M)\odot X
\eeq
\eprop

Further we shall be interested in the case when $X=A$ is a smooth (hence, stereotype) algebra. Then the space ${\mathcal C}(M,A)$ is also endowed with the structure of stereotype algebra with the pointwise multiplication
$$
(u\cdot v)(t)=u(t)\cdot v(t),\qquad u,v\in {\mathcal C}(M,A),\quad t\in M.
$$
From \eqref{C(M,X)-cong-C(M)-odot-X} it follows that ${\mathcal C}(M,A)$ is a stereotype module over $A$.

\btm\label{TH:C(M)-circledast-A->C(M,A)} For each continuous algebra $A$ and each paracompact locally compact space  $M$ the natural mapping
\beq\label{C(M)-circledast-A->C(M,A)}
\iota:{\mathcal C}(M)\circledast A\to {\mathcal C}(M,A) \quad\Big|\quad \iota(u\circledast a)(t)=u(t)\cdot a,\quad u\in {\mathcal C}(M),\ a\in A,\ t\in M,
\eeq
is a continuous envelope:
\beq\label{C(M)-circledast-A=C(M,A)}
{\mathcal C}(M)\overset{\mathcal C}{\circledast} A\cong{\mathcal C}(M,A).
\eeq
\etm

\blm\label{LM:iota-in-DEpi-C(M)}
The image of the mapping $\iota:{\mathcal C}(M)\circledast A\to {\mathcal C}(M,A)$ is dense in ${\mathcal C}(M,A)$.
\elm
\bpr
Take a mapping $f\in {\mathcal C}(M,A)$, a compact set $K\subseteq M$, and a convex neighbourhood of zero $U\subseteq A$. Since $f$ is uniformly continuous on $K$, there is an entourage $V$ in $K$ such that
$$
|s-t|<V\quad\Longrightarrow\quad f(s)-f(t)\in U.
$$
Take a finite sequence of points $t_1,...,t_n\in K$ such that the neighbourhoods $t_k+V$ form a covering of $K$. Find a partition of unity for them on $K$, i.e. functions $\eta_1,...,\eta_n\in{\mathcal C}(M)$ such that
$$
\eta_k\big|_{K\setminus(t_k+V)}=0,\qquad 0\le \eta_k\le 1,\qquad\sum_{k=1}^n\eta_k\big|_K=1.
$$
Put
$$
x=\sum_{k=1}^n\eta_k\circledast f(t_k)\in{\mathcal C}(M)\circledast A.
$$
Then for $s\in K$ we have
$$
f(s)-\iota(x)(s)=\sum_{k=1}^n\eta_k(s)\cdot f(s)-\sum_{k=1}^n\eta_k(s)\cdot f(t_k)=\sum_{k=1}^n\eta_k(s)\cdot \big(f(s)-f(t_k)\big)\in U.
$$
\epr

\blm\label{LM:J^0_C(M)C(M)-circledast-A-cong-J^0_C(M)C(M,A)}
The modules ${\mathcal C}(M)\circledast A$ and ${\mathcal C}(M,A)$ over the algebra ${\mathcal C}(M)$ have isomorphic value bundles:
\beq\label{J^0_C(M)C(M)-circledast-A-cong-J^0_C(M)C(M,A)}
\Jet^0_{{\mathcal C}(M)}{\mathcal C}(M)\circledast A\cong \Jet^0_{{\mathcal C}(M)}{\mathcal C}(M,A),\qquad n\in\N
\eeq
\elm
\bpr
For each point $t\in M$ the ideal $I_t$ has the codimention 1 in ${\mathcal C}(M)$, hence we can use Lemma \ref{PROP:[(X-circledast-Z)/(Y-circledast-Z)]^triangledown-cong-[(X-odot-Z)/(Y-odot-Z)]^triangledown}:
\begin{multline*}
\Jet^0_{{\mathcal C}(M)}{\mathcal C}(M)\circledast A=[({\mathcal C}(M)\circledast A)/ (I_t\circledast A)]^\vartriangle=\eqref{[(X-circledast-Z)/(Y-circledast-Z)]^triangledown-cong-[(X-odot-Z)/(Y-odot-Z)]^triangledown}=\\=
[({\mathcal C}(M)\odot A)/ (I_t\odot A)]^\vartriangle=\Jet^0_{{\mathcal C}(M)}{\mathcal C}(M)\odot A=\eqref{C(M,X)-cong-C(M)-odot-X}=
\Jet^0_{{\mathcal C}(M)}{\mathcal C}(M,A)
\end{multline*}
\epr

\blm\label{TH:diff-oper<-morfizm-rassl-znachenij} Let $M$ be a paracompact locally compact space, $F$ a $C^*$-algebra and $\ph:{\mathcal C}(M)\to F$ a homomorphism of stereotype algebras, and let $\ph({\mathcal C}(M))$ belong to the center of $F$:
$$
\ph({\mathcal C}(M))\subseteq Z(F).
$$
Then for any stereotype space $X$ each morphism of the value bundles
$$
\mu:\Jet_{{\mathcal C}(M)}^0({\mathcal C}(M,X))\to \Jet_{{\mathcal C}(M)}^0(F)
$$
defines a unique morphism of stereotype ${\mathcal C}(M)$-modules $D:{\mathcal C}(M,X)\to F$, satisfying the identities
\beq\label{j^0(Dx)=mu-circ-j^0(x)} \jet^0(\ph(x))=\mu\circ \jet^0(x),\qquad
x\in {{\mathcal C}(M,X)}.
 \eeq \elm
\bpr By Theorem \ref{TH:B-cong-Sec(val_AB)}, the mapping
$v:F\to\Sec(\pi^0_{{{\mathcal C}(M)},F})$, that turns $F$ into the algebra of continuous sections of the value bundle $\pi^0_{{{\mathcal C}(M)},F}: \Jet^0_{{\mathcal C}(M)}F\to\Spec({{\mathcal C}(M)})$ over the algebra  ${{\mathcal C}(M)}$, is an isomorphism of $C^*$-algebras:
$$
F\cong\Sec(\pi^0_{{{\mathcal C}(M)},F}).
$$
Consider the inverse isomorphism $v^{-1}:\Sec(\pi^0_{{{\mathcal C}(M)},F})\to F$:
\beq\label{lambda(j^0(b))=b} v^{-1}(\jet^0(b))=b,\qquad b\in F.
\eeq
Then to each morphism of the value bundles $\mu:\Jet_{{\mathcal C}(M)}^0({{\mathcal C}(M,X)})\to \Jet_{{\mathcal C}(M)}^0(F)$ one can assign an operator $\ph:{{\mathcal C}(M,X)}\to F$ by the formula
\beq\label{Da=lambda(mu-circ-j^0(a))} \ph(x)=v^{-1}\Big(\mu\circ \jet^0(x)\Big),\qquad
x\in {{\mathcal C}(M,X)}.
\eeq
Obviously, it satisfies
\eqref{j^0(Dx)=mu-circ-j^0(x)}.
\epr

\blm The mapping $\iota:{\mathcal C}(M)\circledast A\to {\mathcal C}(M,A)$ is a continuous extension.
\elm
\bpr
Suppose $\ph:{\mathcal C}(M)\circledast A\to B$ is a morphism into a $C^*$-algebra $B$. By Lemma \ref{LM:ph:A-circledast-B->C}, it is representable in the form
\beq\label{ph(u-circledast-a)=eta(u)-cdot-alpha(a)}
\ph(u\circledast a)=\eta(u)\cdot\alpha(a)=\alpha(a)\cdot\eta(u),\qquad u\in {\mathcal C}(M), \ a\in A,
\eeq
where $\eta:{\mathcal C}(M)\to B$, $\alpha:A\to B$ are some morphisms of stereotype algebras. Consider the operator $\eta$ and denote by $C$ its image in $B$:
$$
C=\overline{\eta({\mathcal C}(M))}.
$$
Let $F$ be the commutant of $C$ in $B$:
$$
F=C^!=\{x\in B:\quad \forall c\in C\quad x\cdot c=c\cdot x\}.
$$
Since the algebra $C$ is commutative, it lies in $F$, and moreover in the center of $F$:
$$
C\subseteq Z(F).
$$
Note that the image of the operator $\ph$ lies in $F$:
$$
\ph({\mathcal C}(M)\circledast A)\subseteq F,
$$
since
\begin{multline*}
\ph(v\circledast a)\cdot \eta(u)=\ph(v\circledast a)\cdot \ph(u\circledast 1)=\ph\big((v\circledast a)\cdot (u\circledast 1)\big)=\ph\big((v\cdot u)\circledast 1\big)=\ph\big((u\cdot v)\circledast 1\big)=\\=\ph\big((u\circledast 1)\cdot (v\circledast a)\big)=\ph(u\circledast 1)\cdot \ph(v\circledast a)=\eta(u)\cdot \ph(v\circledast a)
\end{multline*}

To verify that $\iota$ is a continuous extension, we have to show that there is a (unique) homomorphism $\ph':{\mathcal C}(M,A)\to F$ that extends $\ph$:
 \beq\label{prodolzhenie-D-na-C(M,A)}
 \xymatrix @R=2pc @C=1.2pc
 {
 {\mathcal C}(M)\circledast A\ar[rr]^{\iota}\ar[dr]_{\ph} & & {\mathcal C}(M,A)\ar@{-->}[dl]^{\ph'}\\
  & F &
 }
 \eeq
The homomorphism $\ph:{\mathcal C}(M)\circledast A\to F$ is a ${\mathcal C}(M)$-morphism, hense by Theorem
\ref{TH:morf-modul->morf-rassl-znach} there is a morphism of the value bundles
$\jet_0[\ph]:\Jet_{{\mathcal C}(M)}^0[{\mathcal C}(M)\circledast A]\to \Jet_{{\mathcal C}(M)}^0(F)=\pi^0_A F$, such that
$$
\jet^0(\ph(x))=\jet_n[\ph]\circ \jet^n(x),\qquad x\in {\mathcal C}(M)\circledast A.
$$
By Lemma \ref{LM:J^0_C(M)C(M)-circledast-A-cong-J^0_C(M)C(M,A)}, the value bundles of the algebras ${\mathcal C}(M)\circledast A$ and ${\mathcal C}(M,A)$ are isomorphic. Denote this isomorphism by $\mu:\Jet_{{\mathcal C}(M)}^0[{\mathcal C}(M)\circledast A]\gets \Jet^0_{{\mathcal C}(M)}[{\mathcal C}(M,A)]$. Consider the composition $\nu=\jet_0[\ph]\circ\mu:\Jet_{{\mathcal C}(M)}[{\mathcal C}(M,A)]\to \Jet_{{\mathcal C}(M)}^0(F)=\pi^0_A F$:
$$
 \xymatrix @R=2pc @C=1.2pc
 {
 \Jet_{{\mathcal C}(M)}^0[{\mathcal C}(M)\circledast A]\ar[dr]_{\jet_0[\ph]} & & \Jet^0_{{\mathcal C}(M)}[{\mathcal C}(M,A)]\ar@{-->}[dl]^{\quad\nu=\jet_0[\ph]\circ\mu}\ar[ll]_{\mu}\\
  & \Jet^0_{{\mathcal C}(M)}[F] &
 }
$$
By Lemma \ref{TH:diff-oper<-morfizm-rassl-znachenij}, this dotted arrow $\nu$ generates a morphism $\ph':{\mathcal C}(M,A)\to F$ (over the algebra ${\mathcal C}(M)$), such that
$$
\jet^0(\ph'f)=\jet_0[\ph]\circ \jet^0(f),\qquad  f\in {\mathcal C}(M,A).
$$
For each $x\in {\mathcal C}(M)\circledast A$ we have
\beq\label{PROOF:TH:C(M,A)-1}
\jet^0\big(\ph'(\iota(x))\big)=\jet_0[\ph]\circ \jet^0(\iota(x))=\jet_0[\ph]\circ
\jet^0(x)=\jet^0(\ph(x)).
\eeq
Note that by Theorem \ref{TH:B-cong-Sec(val_AB)}, the mapping
$\jet^0=v:F\to\Sec(\pi^0_{{\mathcal E}(M)}F)=\Sec(\Jet^0_{{\mathcal E}(M)}F)$, that turns $F$ into the algebra of continuous sections of the value bundles $\pi^0_{{\mathcal E}(M)}F: \Jet^0_{{\mathcal E}(M)}F\to\Spec({{\mathcal E}(M)})$ over the algebra ${{\mathcal E}(M)}$, is an isomorphism of $C^*$-algebras:
$$
F\cong\Sec(\pi^0_{{\mathcal E}(M)}F).
$$
Hence, we can apply the inverse operator to $\jet^0$ to the equality \eqref{PROOF:TH:C(M,A)-1}, and we obtain
$$
\ph'(\iota(x))=\ph(x).
$$
I.e., $\ph'$ extends $\ph$ in \eqref{prodolzhenie-D-na-C(M,A)}. By Lemma \ref{LM:iota-in-DEpi-C(M)}, the elements of the form $\iota(u\circledast a)$ are total in ${\mathcal C}(M,A)$, that is why this extension $\ph'$ is unique.

By definition, the operator $\ph'$ is a morphism with respect to the algebra ${\mathcal C}(M)$, but this is not sufficient: we have to prove that $\ph'$ is a homomorphism of algebras. To verify this let us take $u,v\in {\mathcal C}(M)$ and $a,b\in A$. Then
$$
\ph'(\iota(u\circledast a)\cdot\iota(v\circledast b))=\ph'(\iota(u\circledast a\cdot v\circledast b))=\ph(u\circledast a\cdot v\circledast b)=\ph(u\circledast a)\cdot\ph(v\circledast b)=\ph'(\iota(u\circledast a))\cdot\ph'(\iota(v\circledast b))
$$
Again by Lemma \ref{LM:iota-in-DEpi-C(M)} the elements of the form $\iota(u\circledast a)$ are total in ${\mathcal C}(M,A)$. As a corollary we can replace them with arbitrary vectors from ${\mathcal C}(M,A)$, hence $\ph'$ must be a homomorphism.
\epr

\blm The mapping $\iota:{\mathcal C}(M)\circledast A\to {\mathcal C}(M,A)$ is a continuous envelope.
\elm
\bpr
Suppose $\sigma:{\mathcal C}(M)\circledast A\to C$ is another continuous extension. We have to verify that there is a morphism $\upsilon$ such that the following diagram is commutative:
$$
 \xymatrix @R=2pc @C=1.2pc
 {
{\mathcal C}(M)\circledast A\ar[rr]^{\sigma}\ar[dr]_{\iota} & & C \ar@{-->}[dl]^{\upsilon}\\
  & {\mathcal C}(M,A) &
 }
$$

Take a compact set $K\subseteq M$ and a homomorphism $\eta:A\to B$ into a $C^*$-algebra, and put
$$
D(u\circledast a)(t)=u(t)\circledast \eta(a),\qquad u\in {\mathcal C}(M),\quad a\in A,\quad t\in K.
$$
The mapping $D$ is a homomorphism from ${\mathcal C}(M)\circledast A$ into ${\mathcal C}(K)\circledast B$. The last algebra is naturally mapped into the tensor product of $C^*$-algebras ${\mathcal C}(K)\underset{\max}{\otimes}B$, which is isomorphic to ${\mathcal C}(K)\odot B$ and to ${\mathcal C}(K,B)$:
$$
{\mathcal C}(K)\circledast B\to {\mathcal C}(K)\underset{\max}{\otimes}B\cong \eqref{A-min-B=A-odot-B}\cong
{\mathcal C}(K)\odot B\cong \eqref{C(M,X)-cong-C(M)-odot-X} \cong {\mathcal C}(K,B).
$$

This means that we can treat $D$ as a morphism into the $C^*$-algebra ${\mathcal C}(K,B)$:
$$
D:{\mathcal C}(M)\circledast A\to {\mathcal C}(K)\circledast B\to {\mathcal C}(K)\underset{\max}{\otimes}B \cong\eqref{A-min-B=A-odot-B}\cong{\mathcal C}(K)\odot B\cong\eqref{C(M,X)-cong-C(M)-odot-X}\cong {\mathcal C}(K,B)
$$
Since $\sigma:{\mathcal C}(M)\circledast A\to C$ is a continuous extension, the homomorphism $D:{\mathcal C}(M)\circledast A\to {\mathcal C}(K,B)$ can be uniquely extended to a homomorphism $D':C\to {\mathcal C}(K,B)$:
 \beq\label{DIAGR:env-C(M)-circledast-A-1}
 \xymatrix @R=2pc @C=1.2pc
 {
{\mathcal C}(M)\circledast A\ar[rr]^{\sigma}\ar[dr]_{D} & & C \ar@{-->}[dl]^{D'}\\
  & {\mathcal C}(K,B) &
 }
 \eeq
If we now freeze $c\in C$ and vary $K\subset M$, then the arising continuous functions $D'(c)$ on $K$ are coherent in the sense that they coinside on the intersections of their domains. As a corollary there is a continuous function $\iota_B'(c):M\to B$ such that its restriction on each compact set $K$ coincide with the corresponding function  $D'(c)$:
$$
\iota'(c)\big|_K=D'(c),\qquad K\subset M.
$$
In other words there is a mapping $\iota'_B:C\to{\mathcal C}(M,B)$ (which by construction is a homomorphism of algebras), such that the following diagram (that specifies \eqref{DIAGR:env-C(M)-circledast-A-1}) is commutative:
 \beq\label{DIAGR:env-C(M)-circledast-A-*-M}
 \xymatrix @R=2.5pc @C=2pc
 {
 {\mathcal C}(M)\circledast A\ar[rr]^{\sigma}\ar@/_5ex/[ddr]_{D}\ar[dr]_{\iota_B} & & C \ar@/^5ex/[ddl]^{D'}\ar@{-->}[dl]^{\iota'_B}\\
  & {\mathcal C}(M,B)\ar[d]_{\rho_K} & \\
  & {\mathcal C}(K,B) &
 }
 \eeq
(here $\rho_K$ is the mapping of restriction to $K$).

Let now $U$ be a $C^*$-neighbourhood of zero\footnote{$C^*$-neighbourhood of zero were defined on page \pageref{DEF:C*-neighbourhood-of-zero}.} in $A$, that corresponds to the homomorphism $\eta:A\to B$. Since $\sigma$ is a dense epimorphism, the upper inner triangle in \eqref{DIAGR:env-C(M)-circledast-A-*-M} can be added to the diagram
 \beq\label{DIAGR:env-E(M)-circledast-A-*-M-1}
 \xymatrix @R=2.5pc @C=2pc
 {
 {\mathcal C}(M)\circledast A\ar[rr]^{\sigma}\ar@/_5ex/[ddr]_{\iota_B}\ar[dr]_{\vartheta_U} & & C \ar@/^5ex/[ddl]^{\iota'_B}\ar@{-->}[dl]^{\vartheta'_U}\\
  & {\mathcal C}(M,A/U)\ar[d]_{\eta_U\oslash 1_M} & \\
  & {\mathcal C}(M,B) &
 }
 \eeq
where $\eta_U:A/U\to B$ is a morphism from \eqref{ph=ph_U-circ-pi_U-0}, and
$$
\vartheta_U(u\circledast a)(t)=u(t)\cdot\pi_U(a),\qquad u\in{\mathcal E}(M),\quad a\in A,\quad t\in M,
$$
$$
(\eta_U\oslash 1_M)(h)(t)=\eta_U(h(t)),\qquad h\in {\mathcal E}(M,A/U),\quad t\in M.
$$
From the definition of $\vartheta_U$ it follows immediatetly that if $U'\subseteq U$ is another $C^*$-neighbourhood of zero, then
\beq\label{vartheta_U=(varkappa^U'_U-oslash-1_M)-cdot-vartheta_U'-0}
\vartheta_U=(\varkappa^{U'}_U\oslash 1_M)\cdot\vartheta_{U'},\qquad U\supseteq U',
\eeq
where $\varkappa^{U'}_U$ is the morphism from \eqref{varkappa^U'_U-0}, and
$$
(\varkappa^{U'}_U\oslash 1_M)(h)(t)=\varkappa^{U'}_U(h(t)),\qquad h\in {\mathcal E}(M,A/U'),\quad t\in M.
$$
The equality \eqref{vartheta_U=(varkappa^U'_U-oslash-1_M)-cdot-vartheta_U'-0} is the left lower inner triangle in the diagram
$$
 
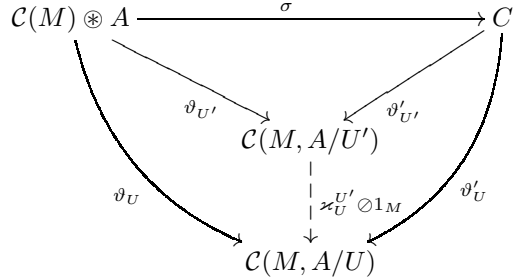
\xymatrix @R=2.5pc @C=3pc
 {
{\mathcal C}(M)\circledast A\ar[rr]^{\sigma}\ar[dr]_(.6){\vartheta_{U'}}\ar@/_5ex/[ddr]_{\vartheta_U} & & C \ar[dl]^(.6){\vartheta'_{U'}} \ar@/^5ex/[ddl]^{\vartheta'_U}\\
 & {\mathcal C}(M,A/U')\ar@{-->}[d]^{\varkappa^{U'}_U\oslash 1_M} & \\
  & {\mathcal C}(M,A/U) &
 }
$$
At the same time the perimeter and the upper inner triangle here are variants of the upper inner triangle in \eqref{DIAGR:env-C(M)-circledast-A-*-M}, and besides this $\sigma$ is an isomorphism. As a corollary, the remaining right lower inner triangle must be commutative as well.

This means that the morphisms $\vartheta'_U:C\to {\mathcal C}(M,A/U)$ form a projective cone of the system  $\varkappa^{U'}_U\oslash 1_M$, and thus there is a morphism $\vartheta'$ into the projective limit:
$$
 
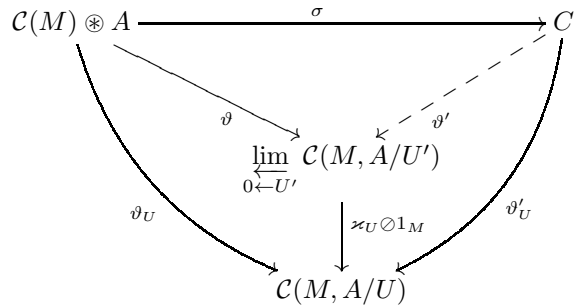
\xymatrix @R=2.5pc @C=3pc
 {
{\mathcal C}(M)\circledast A\ar[rr]^{\sigma}\ar[dr]_(.6){\vartheta}\ar@/_5ex/[ddr]_{\vartheta_U} & & C \ar@{-->}[dl]^(.6){\vartheta'} \ar@/^5ex/[ddl]^{\vartheta'_U}\\
 & \leftlim\limits_{0\gets U'}{\mathcal C}(M,A/U')\ar[d]^{\varkappa_U\oslash 1_M} & \\
  & {\mathcal C}(M,A/U) &
 }
$$
Let us note now the following chain:
$$
\leftlim\limits_{0\gets U'}{\mathcal C}(M,A/U')=\leftlim\limits_{0\gets U'}({\mathcal C}(M)\odot A/U')=\cite[(2.53)]{Akbarov-env}=
{\mathcal C}(M)\odot \leftlim\limits_{0\gets U'}A/U'={\mathcal C}(M,\leftlim\limits_{0\gets U'}A/U')
$$
And let us put the last space into our diagram:
$$
 
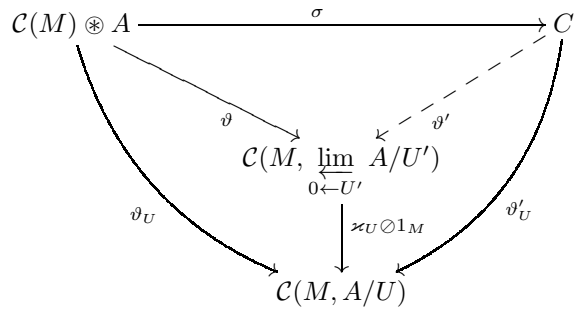
\xymatrix @R=2.5pc @C=3pc
 {
{\mathcal C}(M)\circledast A\ar[rr]^{\sigma}\ar[dr]_(.6){\vartheta}\ar@/_5ex/[ddr]_{\vartheta_U} & & C \ar@{-->}[dl]^(.6){\vartheta'} \ar@/^5ex/[ddl]^{\vartheta'_U}\\
 & {\mathcal C}(M,\leftlim\limits_{0\gets U'}A/U')\ar[d]^{\varkappa_U\oslash 1_M} & \\
  & {\mathcal C}(M,A/U) &
 }
$$

Recall again that $\sigma$ is a dense epimorphism. This implies that the arrow $\vartheta'$ can be lifted to an arrow $\upsilon$ with the values in the space ${\mathcal C}(M,\Im\pi)$ of functions, which  have images in the image of the mapping $\pi:A\to\leftlim\limits_{0\gets U'}A/U'$, or, what is the same, in the immediate subspace, generated by the set of values of the mapping $\pi$. This space coincide with $A$, since $A$ is a continuous algebra:
$$
\Im\pi\cong\Env_{\mathcal C}A\cong A
$$
We obtain the following diagram
$$
 
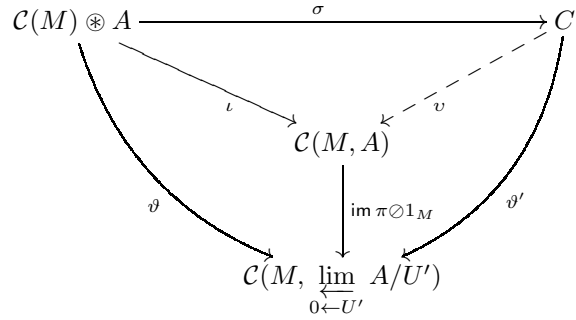
\xymatrix @R=2.5pc @C=3pc
 {
{\mathcal C}(M)\circledast A\ar[rr]^{\sigma}\ar[dr]_(.6){\iota}\ar@/_5ex/[ddr]_{\vartheta} & & C \ar@{-->}[dl]^(.6){\upsilon} \ar@/^5ex/[ddl]^{\vartheta'}\\
 & {\mathcal C}(M,A)\ar[d]^{\im\pi\oslash 1_M} & \\
  & {\mathcal C}(M,\leftlim\limits_{0\gets U'}A/U') &
 }
$$
where $\pi$ is the morphism from \eqref{A-to-leftlim-A/U-0}.
\epr

\subsection{$\mathcal{C}(M)$ as a continuous envelope of its subalgebras}

Let us call a continuous mapping of topological spaces $\e:X\to Y$ a {\it covering}, if each compact set $T\subseteq Y$ is contained in the image of some compact set $S\subseteq X$. If the space $Y$ is Hausdorff, then this implies automatically that $\e$ is surjective. If in addition $\e$ is injective, then we call it an {\it exact covering}.\label{DEF:nalozhenie} In an exact covering $\e:X\to Y$ the space $Y$ can be treated as a new, weaker topologization of $X$, which does not change the system of compact sets and the topology on each compact set.

\btm\label{C-obolochka-podalgebry-v-C(M)} Let $A$ be an involutive stereotype subalgebra in the algebra $\mathcal{C}(M)$ of continuous functions on a paracompact locally compact space $M$, i.e. there is a (continuous and unital) monomorphism of involutive stereotype algebras
$$
\iota:A\to \mathcal{C}(M).
$$
The continuous envelope of $A$ coincides with $\mathcal{C}(M)$
\beq\label{Env_C_A=C(M)}
\Env_{\mathcal C} A=\mathcal{C}(M)
\eeq
(i.e. $\iota$ is a continuous envelope of $A$), if and only if the dual mapping of spectra $\iota^{\Spec}:\Spec(A)\gets M$ is an exact covering. \etm

\bpr
1. First, we prove necessity. Suppose, \eqref{Env_C_A=C(M)} holds. Take a compact set $T\subseteq\Spec(A)$ and consider the mapping
$$
\ph_T:A\to{\mathcal C}(T)\quad\Big|\quad \ph_T(a)(t)=t(a),\qquad t\in T,\ a\in A.
$$
This is an involutive homomorphism into a $C^*$-algebra, hence it can be extended to the envelope $\Env_{\mathcal C} A=\mathcal{C}(M)$:
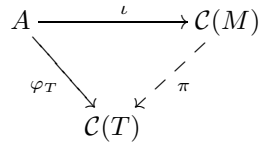
$$
 \xymatrix @R=2pc @C=1.2pc
 {
 A\ar[rr]^{\iota}\ar[dr]_{\ph_T} & & {\mathcal C}(M)\ar@{-->}[dl]^{\pi}\\
  & {\mathcal C}(T) &
 }
$$
where $\pi$ is an involutive homomorphism. The dual mapping of spectra $\pi^{\Spec}:\Spec({\mathcal C}(T))\to\Spec({\mathcal C}(M))$ turns $T=\Spec({\mathcal C}(T))$ into some compact set $S\subseteq \Spec({\mathcal C}(M))$. After that the mapping of spectra $\iota^{\Spec}:\Spec(A)\gets\Spec({\mathcal C}(M))$ turns $S$ exactly to the compact set $T$. This proves that $\iota^{\Spec}:\Spec(A)\gets\Spec({\mathcal C}(M))$ is a covering. Let us show that this mapping is injective. If this was not so, there would be two points $s\ne s'\in M$ that stick together under the mapping $\iota^{\Spec}$:
$$
s\circ\sigma=s'\circ\sigma=t\in\Spec(A)
$$
In other words, the character $t:A\to\C$ has two different extensions on ${\mathcal C}(M)$:
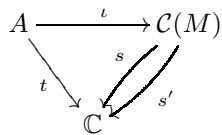
$$
 \xymatrix @R=2pc @C=1.2pc
 {
 A\ar[rr]^{\iota}\ar[dr]_{t} & & {\mathcal C}(M)\ar@/_1ex/[dl]_{s}\ar@/^1ex/[dl]^{s'}\\
  & \C &
 }
$$
But $\iota$ is a continuous extension, so the character $t:A\to\C$, being an involutive homomorphism into the $C^*$-algebra $\C$, must have unique extension.

2. Now we prove sufficiency. Suppose $\iota^{\Spec}$ is an exact covering. Then the algebra $A$ differs the points of $M$, and, since it contains the unit (and hence, all constants), by the Stone-Weierstrass theorem $A$ must be dense in $\mathcal{C}(M)$.

Let us show that $\iota:A\to \mathcal{C}(M)$ is a continuous extension. Suppose $\ph:A\to B$ is a morphism of $A$ into a $C^*$-algebra $B$. To construct a dotted arrow $\ph'$ for the diagram \eqref{DEF:diagr-nepr-rasshirenie},
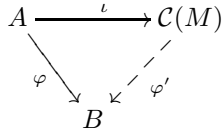
$$
 \xymatrix @R=2pc @C=1.2pc
 {
 A\ar[rr]^{\iota}\ar[dr]_{\ph} & & {\mathcal C}(M)\ar@{-->}[dl]^{\ph'}\\
  & B &
 }
$$
it is sufficient to think that $B$ is commutative and that $\ph(A)$ is dense in $B$ (since we always can replace $B$ by the closure $\overline{\ph(A)}$ in $B$, which is acommutative subalgebra in $B$). Then from the commutativity of  $B$ it follows that $B$ has the form ${\mathcal C}(T)$, and from the density of $\ph(A)$ in $B$ -- that the compact set $T$ is embedded (injectively) into $\Spec(A)$. By the premise of the theorem, $\iota^{\Spec}:\Spec(A)\gets M$ is a covering, hence $T\subseteq\Spec(A)$ is the image of a compact set $K\subseteq M$, $K\cong T$. Then the mapping $\ph$ is represented as a composition of the injection $\iota:A\to{\mathcal C}(M)$ and the mapping $\pi_T:{\mathcal C}(M)\to {\mathcal C}(T)$ of the restriction to the compact set $T$:
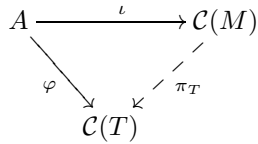
$$
 \xymatrix @R=2pc @C=1.2pc
 {
 A\ar[rr]^{\iota}\ar[dr]_{\ph} & & {\mathcal C}(M)\ar@{-->}[dl]^{\pi_T}\\
  & {\mathcal C}(T) &
 }
$$
We obtain that $\ph'=\pi_T$. And the dotted arrow is unique since $A$ is dense in ${\mathcal C}(M)$.

Let us check that $\iota:A\to {\mathcal C}(M)$ is a maximal extension, i.e. if we take another extension $\sigma:A\to C$, then there appears a morphism $\upsilon:C\to {\mathcal C}(M)$, such that the following  diagram is commutative:
 \beq
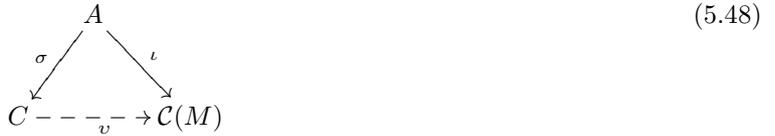
\label{TOP:C->C(M)}
 \xymatrix @R=2pc @C=1.2pc
 {
 & A\ar[rd]^{\iota}\ar[ld]_{\sigma} & \\
 C\ar@{-->}[rr]_{\upsilon} & & {\mathcal C}(M)
 }
 \eeq
For each compact set $T\subseteq M$ the homomorphism
$$
\iota_T:A\to{\mathcal C}(T)\quad\big|\quad \iota_T(a)(t)=t(a),\qquad t\in T\subseteq\Spec(A)
$$
is uniquely extended to a homomorphism $\iota'_T:A'\to{\mathcal C}(T)$
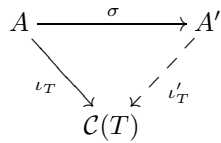
$$
 \xymatrix @R=2pc @C=1.2pc
 {
 A\ar[rr]^{\sigma}\ar[dr]_{\iota_T} & & A'\ar@{-->}[dl]^{\iota'_T}\\
  & {\mathcal C}(T) &
 }
$$
If $T\subseteq S$ are two compact sets in $M$, then this diagram is added to a diagram
$$
 
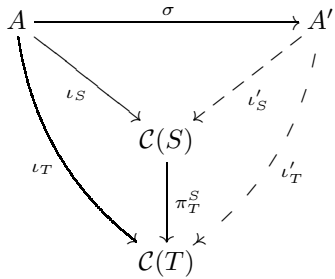
\xymatrix @R=2.5pc @C=3pc
 {
 A\ar[rr]^{\sigma}\ar@/_3ex/[ddr]_{\iota_T}\ar[dr]_{\iota_S} & & A'\ar@{-->}[dl]^{\iota'_S}\ar@/^3ex/@{-->}[ddl]^{\iota'_T}\\
  & {\mathcal C}(S)\ar[d]^{\pi_T^S} & \\
  & {\mathcal C}(T) &
 }
$$
where $\pi_T^S$ is the restriction to the compact $T$. The right lower triangle in this diagram means that the morphisms $\iota'_T:A'\to {\mathcal C}(T)$ form a projective cone in the contravariant system $\pi_T^S$. Hence there is an arrow $\iota'$, such that the right lower triangle is commutative in all diagrams
$$
 
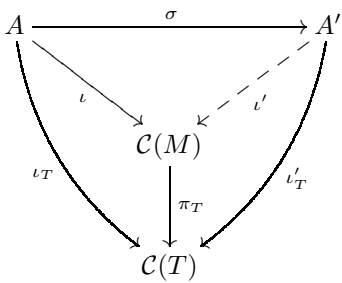
\xymatrix @R=2.5pc @C=3pc
 {
 A\ar[rr]^{\sigma}\ar@/_3ex/[ddr]_{\iota_T}\ar[dr]_{\iota} & & A'\ar@{-->}[dl]^{\iota'}\ar@/^3ex/[ddl]^{\iota'_T}\\
  & {\mathcal C}(M)\ar[d]^{\pi_T} & \\
  & {\mathcal C}(T) &
 }
$$
Since the perimeter and the left lower triangle in any such a diagram are commutative, we have
$$
\pi_T\circ\iota'\circ\sigma=\iota'_T\circ\sigma=\iota_T=\pi_T\circ\iota.
$$
This is true for all $T$, thus we have
$$
\iota'\circ\sigma=\iota.
$$
\epr

\paragraph{A counterexample.}

\bex\label{EX:C(8)} {\it
There is a dense involutive subalgebra $A$ in ${\mathcal C}(\R)$, such that the dual mapping of spectra $\iota^{\Spec}:\Spec(A)\gets M$ is a bijection, but the continuous envelope of $A$ does not coincide with ${\mathcal C}(\R)$:}
$$
\Spec(A)=M,\qquad \Env_{\mathcal C}A\ne {\mathcal C}(\R)
$$
\eex
\bpr
This is the algebra $A$ of continuous functions on $\R$, which have the limit in the infinity that coincides with the value in some point, for example in zero:
$$
u\in A\quad\Longleftrightarrow\quad u\in {\mathcal C}(\R)\quad\&\quad \lim_{t\to\infty}u(t)=u(0).
$$
The algebra $A$ is involutive, it contains constants and separates the points of $\R$, hence it is dense in  ${\mathcal C}(\R)$. We endow $A$ with the topology of uniform convergence on the whole line $\R$. Obvioulsy, the spectrum of $A$ is the line $\R$ with the new topology, where the basis neighbourhoods of the point 0 are the sets of the form
$$
(-\infty, A)\cup(a,b)\cup(B,+\infty)
$$
where $A<a<0<b<B$ are numbers in $\R$ (and the basis neighbourhoods of other points don't change). It is convenient to perceive this topology as the one induced on $\R$ by the embedding of $\R$ into the figure resembling in its shape the number 8, and that is why we denote it by 8 (in this embedding the ends of the line $\R$ incurve and approach to the point $0\in\R$). Certainly, the spectrum of $A$ is homeomorphic to 8,
$$
\Spec A\cong 8,
$$
and the algebra $A$ is isomorphic (as a stereotype algebra) to the algebra ${\mathcal C}(8)$ (of continuous functions on the compact space 8). Therefore the continuous envelope of $A$ coincides with $A$, and is not isomorphic  ${\mathcal C}(\R)$:
$$
\Env_{\mathcal C}A\cong A\cong{\mathcal C}(8)\not\cong {\mathcal C}(\R).
$$
\epr

\paragraph{The continuous envelope of the algebra $\Trig(G)=k(G)$ for a compact group $G$.}

By  \cite[(30.30)]{Hewitt-Ross-2}, $\Spec(\Trig(G))=G$. Together with the Theorem \ref{C-obolochka-podalgebry-v-C(M)} this gives the following result:

\btm The continuous envelope of the algebra $\Trig(G)=k(G)$ on a compact group $G$ coincides with the algebra ${\mathcal C}(G)$ of continuous functions on $G$:
\beq
\Env_{\mathcal C}\Trig(G)={\mathcal C}(G).
\eeq
\etm

\subsection{Continuous envelopes of group algebras}

\paragraph{Fourier transform on a commutative locally compact group.}

Let $C$ be a commutative locally compact group. Recall the algebra $\mathcal{C}(C)$ of continuous functions and the algebra $\mathcal{C}^\star(C)$ of measures with compact support on $C$, which we were talking about in \pageref{paragraph:C^star(G)}. The formula
\beq\label{Fourier-transform}
\overbrace{{\mathcal F}_C(\alpha)(\chi)}^{\scriptsize \begin{matrix}
\text{value of the function ${\mathcal F}_C(\alpha)\in {\mathcal C}(\widehat{C})$}\\
\text{in the point $\chi\in \widehat{C}$} \\ \downarrow \end{matrix}}\kern-35pt=\kern-50pt\underbrace{\alpha(\chi)}_{\scriptsize \begin{matrix}\uparrow \\
\text{action of the functional $\alpha\in{\mathcal C}^\star(C)$}\\
\text{at the function $\chi\in \widehat{C}\subseteq {\mathcal C}(C)$ }\end{matrix}}
\kern-50pt \qquad (\chi\in \widehat{C},\quad \alpha\in {\mathcal C}^\star(C))
\eeq
defines a mapping
$$
{\mathcal F}_C:{\mathcal C}^\star(C)\to {\mathcal C}(\widehat{C})
$$
which is a homomorphism of involutive stereotype algebras, and is called the {\it Fourier transform} on the group $C$.

In \cite[Theorem 2.11]{Kuznetsova} (see also \cite{Akbarov-env}) the following result is proved:

\bprop\label{PROP:Env_C-C^star(G)=C(widehat(G))}
The Fourier transform on a commutative locally compact group $C$ is a continuous envelope of the group algebra ${\mathcal C}^\star(C)$. As a corollary,
\beq\label{Env_C-C^star(G)=C(widehat(G))}
\Env_{\mathcal C}{\mathcal C}^\star(C)={\mathcal C}(\widehat{C}).
\eeq
\eprop

\paragraph{The continuous envelope of the group algebra of a compact group.}

\bprop\label{PROP:nepr-obol-komp-gruppy}
For a compact group $K$ the continuous envelope of its group algebra ${\mathcal C}^\star(K)$ is the Cartesian product of the algebras ${\mathcal B}(X_\pi)$, where $\pi$ runs over the dual object $\widehat{K}$, and $X_\pi$ is the space of the representation $\pi$:
\beq\label{nepr-obol-komp-gruppy}
\Env_{\mathcal C}{\mathcal C}^\star(K)=\prod_{\pi\in\widehat{K}}{\mathcal B}(X_\pi).
\eeq
\eprop
\bpr
1. Let us prove that the mapping $P=\prod_{\pi\in\widehat{K}}\pi:{\mathcal C}^\star(K)\to\prod_{\pi\in\widehat{K}}{\mathcal B}(X_\pi)$ is a continuous extension. Note first that $P$ is a dense epimorphism (this follows from the definition of the direct product). Further, suppose $\psi:{\mathcal C}^\star(K)\to B$ is an involutive homomorphism into a $C^*$-algebra $B$. Consider an inclusion of $C^*$-algebras $\eta:B\to {\mathcal B}(X)$. The composition $\psi=\eta\circ\ph$ generates a norm-continuous representation  $\rho=\psi\circ\delta:K\to {\mathcal B}(X)$, which by Theorem \ref{TH:nepr-po-norme-perdst-K} can be decomposed into a direct sum of unitary continuous representations, of which only finite number are not equivalent to each other. This means, in particular, that there exists a finite set $M\subseteq\widehat{K}$ such that $\psi$ can be represented as a composition $\psi'\circ P_M$, where $P_M$ is the natural projection of ${\mathcal C}^\star(K)$ into the direct product $\prod_{\pi\in M}{\mathcal B}(X_\pi)$. This in its turn implies that the homomorphism $\ph$ vanishes on the kernel of $P_M$: $\Ker P_M\subseteq\Ker\ph$. In addition, the algebra $\prod_{\pi\in M}{\mathcal B}(X_\pi)$ is finite dimensional and isomorphic to the quotient algebra ${\mathcal C}^\star(K)/\Ker P_M$. Hence $\ph$ is representable as a composition $\ph'\circ P_M$, and we obtain the diagram
$$
 
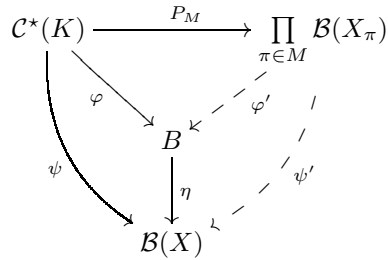
\xymatrix @R=2pc @C=1.2pc
 {
{\mathcal C}^\star(K)\ar[rr]^{P_M}\ar[dr]_{\ph}\ar@/_4ex/[ddr]_{\psi} & &
\prod\limits_{\pi\in M}{\mathcal B}(X_\pi)\ar@{-->}[dl]^{\ph'} \ar@{-->}@/^4ex/[ddl]^{\psi'}
\\
  & B\ar[d]^{\eta} & \\
  & {\mathcal B}(X) &
 }
$$
It implies the diagram
$$
 
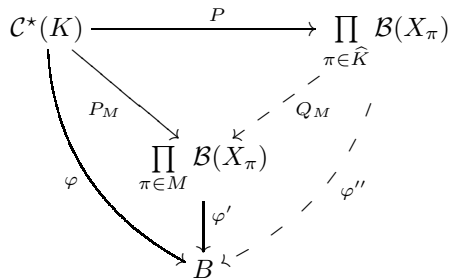
\xymatrix @R=2pc @C=1.2pc
 {
{\mathcal C}^\star(K)\ar[rr]^{P}\ar[dr]_{P_M}\ar@/_4ex/[ddr]_{\ph} & & \prod\limits_{\pi\in\widehat{K}}{\mathcal B}(X_\pi)\ar@{-->}[dl]^{Q_M} \ar@{-->}@/^4ex/[ddl]^{\ph''}
\\
  & \prod\limits_{\pi\in M}{\mathcal B}(X_\pi) \ar[d]^{\ph'} & \\
  & B &
 }
$$
where $Q_M$ is the natural projection of the direct product into its subproduct.

2. Further we prove that $P:{\mathcal C}^\star(K)\to\prod_{\pi\in\widehat{K}}{\mathcal B}(X_\pi)$ is a continuous envelope. Suppose $Q:{\mathcal C}^\star(K)\to A$ is another continuous extension. Then for each representation $\sigma\in\widehat{K}$ there is a (unique) morphism $\alpha_\sigma:A\to {\mathcal B}(X_\sigma)$ such that $\ph_\sigma=\alpha_\sigma\circ Q$. The family of morphisms $\{\alpha_\sigma:A\to {\mathcal B}(X_\sigma);\ \sigma\in\widetilde{K}\}$ generates a morphism $\upsilon:A\to \prod_{\pi\in\widehat{K}}{\mathcal B}(X_\pi)$ such that $\alpha_\sigma=\iota_\sigma\circ\upsilon$ for each $\sigma\in\widehat{K}$. We obtain a diagram
$$
 
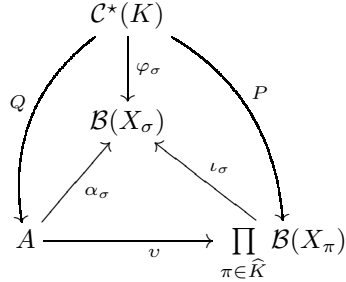
\xymatrix @R=2pc @C=1.2pc
 {
& {\mathcal C}^\star(K)\ar@/_4ex/[ddl]_{Q}\ar@/^4ex/[ddr]^{P}\ar[d]^{\ph_\sigma} &
\\
& {\mathcal B}(X_\sigma) &
\\
A\ar[rr]_{\upsilon}\ar[ur]_{\alpha_\sigma} & & \prod\limits_{\pi\in\widehat{K}}{\mathcal B}(X_\pi)\ar[ul]_{\iota_\sigma}
 }
$$
where the inner little triangles are commutative due to the properties of the mappings $P$, $Q$, $\upsilon$. As a corollary,
$$
\iota_\sigma\circ\upsilon\circ Q=\alpha_\sigma\circ Q=\ph_\sigma=\iota_\sigma\circ P,
$$
and since this is true for any $\sigma$, we have
$$
\upsilon\circ Q=P,
$$
i.e. in this diagram the perimeter is also commutative. The uniqueness of the morphism $\upsilon$ follows from the uniqueness of $\alpha_\sigma$.
\epr

\paragraph{The continuous envelope of the group algebra of the group $C\times K$.}

Let $C$ be an Abelian locally compact group, and $K$ a compact group (not necessarily Abelian).

\bprop\label{TH:env_C^star(R^n-times-K)}
the continuous envelope of the group algebra ${\mathcal C}^\star(C\times K)$ is the algebra ${\mathcal C}(\widehat{C},\prod_{\sigma\in\widehat{K}}{\mathcal B}(X_\sigma))$ of continuous mappings from the Pontryagin dual group $\widehat{C}$ into the Cartesian product of the algebras ${\mathcal B}(X_\sigma)$, where $\sigma$ runs over the dual object $\widehat{K}$, and $X_\sigma$ is the space of the representation $\sigma$:
\beq\label{env_C^star(R^n-times-K)}
\Env_{\mathcal C}{\mathcal C}^\star(C\times K)={\mathcal C}\Big(\widehat{C},\prod_{\sigma\in\widehat{K}}{\mathcal B}(X_\sigma)\Big)=
\prod_{\sigma\in\widehat{K}}{\mathcal C}\big(\widehat{C},{\mathcal B}(X_\sigma)\big)={\mathcal C}(\widehat{C})\odot \prod_{\sigma\in\widehat{K}}{\mathcal B}(X_\sigma)=\Env_{\mathcal C}{\mathcal C}^\star(C)\odot \Env_{\mathcal C}{\mathcal C}^\star(K).
\eeq
\eprop
\bpr The first equality is proved by the chain
\begin{multline*}
\Env_{\mathcal C}{\mathcal C}^\star(C\times K)=\Env_{\mathcal C}\Big({\mathcal C}^\star(C)\circledast{\mathcal C}^\star(K)\Big)=
\cite[(1.129)]{Akbarov-env}=
\Env_{\mathcal C}\Big(\Env_{\mathcal C}{\mathcal C}^\star(C)\circledast\Env_{\mathcal C}{\mathcal C}^\star(K)\Big)=\eqref{Env_C-C^star(G)=C(widehat(G))}=\\=
\Env_{\mathcal C}\Big({\mathcal C}(\widehat{C})\circledast\prod_{\sigma\in\widehat{K}}{\mathcal B}(X_\sigma)\Big)=\eqref{C/circledast}=
{\mathcal C}(\widehat{C})\overset{\mathcal C}{\circledast}\prod_{\sigma\in\widehat{K}}{\mathcal B}(X_\sigma)=\eqref{C(M)-circledast-A=C(M,A)}=
{\mathcal C}\Big(\widehat{C},\prod_{\sigma\in\widehat{K}}{\mathcal B}(X_\sigma)\Big)
\end{multline*}
The second one in \eqref{env_C^star(R^n-times-K)} is obvious, the third one follows from \cite[Theorem 8.1]{Akbarov}, and the last one from \eqref{Env_C-C^star(G)=C(widehat(G))} and \eqref{nepr-obol-komp-gruppy}.
\epr

\paragraph{The continuous evelope of the group algebra of a discrete group.}

For a discrete group $D$ its group algebra is the algebra of functions on $D$ with a finite support:
$$
{\mathcal C}^\star(D)=\C_D=\{\alpha=\{\alpha_x,x\in D\}:\quad \card\{x\in D:\alpha_x\ne 0\}<\infty\}.
$$
The convolution on $\C_D$ is defined by its action on delta-functionals \eqref{delta^a*delta^b=delta^(a-cdot-b)}.

Note that each $C^*$-seminorm $p$ on $\C_D$ turn the unit either into zero, or into unit:
$$
p(\delta^e)=p(\delta^e*\delta^e)=p(\delta^e*(\delta^e)^\bullet)=p(\delta^e)^2\quad\Longrightarrow\quad p(\delta^e)=0\quad\vee\quad p(\delta^e)=1.
$$
In the first case $p$ turns any element into zero (since $p$ is always submultiplicative). Hence if $p\ne 0$, then  $p(\delta^e)=1$. Moreover, in this case each delta functional is turned into unit:
$$
1=p(\delta^e)=p(\delta^a*\delta^{a^{-1}})=p(\delta^a*(\delta^a)^\bullet)=p(\delta^a)^2\quad\Longrightarrow\quad p(\delta^a)=1.
$$
This implies that each $C^*$-seminorm $p$ on $\C_D$ is subordinated to the $\ell_1$-norm:
\beq\label{||pi(alpha)||-le-||alpha||_1}
p(\alpha)\le \norm{\alpha}_1,\qquad \alpha\in\C_D,
\eeq
since
$$
p(\alpha)=p\Big(\sum_{x\in D}\alpha_x\cdot\delta^x\Big)\le
\sum_{x\in D}|\alpha_x|\cdot p(\delta^x)\le\sum_{x\in D}|\alpha_x|\cdot 1=\norm{\alpha}_1.
$$
From \eqref{||pi(alpha)||-le-||alpha||_1} it follows that for each $\alpha\in\C_D$ there exists the supremum by all  $C^*$-seminorms
\beq\label{||alpha||^*}
\norm{\alpha}_\bullet=\sup_{p\in {\mathcal P}(\C_D)}p(\alpha)\le \norm{\alpha}_1.
\eeq
This is a $C^*$-seminorm on $\C_D$, since
$$
\norm{\alpha*\alpha^\bullet}_\bullet=\sup_{\pi\in\widehat{G}}\norm{\pi(\alpha*\alpha^\bullet)} =\sup_{\pi\in\widehat{G}}\norm{\pi(\alpha)*\pi(\alpha)^\bullet}=\sup_{\pi\in\widehat{G}}\norm{\pi(\alpha)}^2=
\big(\sup_{\pi\in\widehat{G}}\norm{\pi(\alpha)}\big)^2=\norm{\alpha}_\bullet^2.
$$
Moreover, this is a norm on $\C_D$, since if $\alpha\ne0$, then in the left regular representation $\pi:D\to{\mathcal B}(L_2(D))$ it turns into a nonzero element, which is separated from the zero by the norm in ${\mathcal B}(L_2(D))$, and this norm defines the $C^*$-seminorm on $\C_D$. The completion of the algebra $\C_D$ with respect to this norm coincides with the completion of $\ell_1(D)$ with respect to this norm, and is called the {\it group $C^*$-algebra} of the group $G$ and is denoted by $C^*(D)$ \cite{Dixmier}.

\bprop\label{PROP:nepr-obol-diskr-gruppy}
For a discrete group $D$ the continuous envelope of its group algebra ${\mathcal C}^\star(D)=\C_D$ is the group  $C^*$-algebra $C^*(D)$:
$$
\Env_{\mathcal C}\C_D=C^*(D).
$$
\eprop
\bpr
Let $\rho:\C_D\to C^*(D)$ be the mapping of completion with respect to the norm $\norm{\cdot}_\bullet$. Let us show that it is a continuous extension. Suppose $\ph:\C_D\to B$ is an involutive homomorphism into a $C^*$-algebra $B$. To complete the diagram
\beq\label{diagr-ph-ph'}
 \xymatrix @R=2pc @C=1.2pc
 {
 \C_D\ar[rr]^{\rho}\ar[dr]_{\ph} & & C^*(D)\ar@{-->}[dl]^{\ph'} \\
  & B &
 }
\eeq
it is sufficient to assume that $\ph$ has dense image in $B$. Then one can treat $B$ as the completion of the algebra $\C_D$ with respect to a $C^*$-seminorm $p$ (after taking the quotient by the kernel of this norm). But $p$, being a  $C^*$-seminorm, must be subordinated to the norm $\norm{\cdot}_\bullet$. Hence $\Ker\norm{\cdot}_\bullet\subseteq\Ker p$. This implies that $\ph$ can be represented as a composition $\C_D\to \C_D/\Ker\norm{\cdot}_\bullet\to B$. This representation gives \eqref{diagr-ph-ph'}.

Now let us prove that $\rho:\C_D\to C^*(D)$ is a continuous envelope. Suppose $\sigma:\C_D\to A$ is another continuous extension. Then, since $\rho:\C_D\to C^*(D)$ is a homomorphism into a $C^*$-algebra, it can be factored through $\sigma$:
$$
 \xymatrix @R=2pc @C=1.2pc
 {
 \C_D\ar[rr]^{\sigma}\ar[dr]_{\rho} & & A\ar@{-->}[dl]^{\upsilon} \\
  & C^*(D) &
 }
$$
\epr

\paragraph{The continuous envelope of the group algebra of a SIN-group.}
Recall that the SIN-groups were defined on page \pageref{DEF:SIN-gruppa}. According to Theorem \ref{TH:stroenie-SIN},
each SIN-group $G$ is a SIN-group is a discrete extension of a group $\R^n\times K$, where $n\in\Z_+$, and $K$ is a compact group:
$$
1\to \R^n\times K=N\to G\to D\to 1
$$
($D$ is a discrete group).

We denote by  ${\mathcal P}(G)$ the set of all continuous $C^*$-seminorms on a group algebra ${\mathcal C}^\star(G)$ of a group $G$. Certainly, the group algebra ${\mathcal C}^\star(N)$ is embedded into the group algebra ${\mathcal C}^\star(G)$. We denote this embedding by $\theta:{\mathcal C}^\star(N)\to {\mathcal C}^\star(G)$.

 \vglue10pt
\centerline{\bf Properties of the seminorms extensions:\footnote{See Errata on page \pageref{Errata}.}}
 \vglue10pt

\bit{\it

\item[$1^\circ$.]\label{1^0:prodolzhenie-polunorm} Each seminorm $p\in{\mathcal P}(N)$ can be extended to a seminorm $q\in{\mathcal P}(G)$.
$$
 \xymatrix @R=2pc @C=1.2pc
 {
  {\mathcal C}^\star(N)\ar[rr]^{\theta}\ar[dr]_{p} & & {\mathcal C}^\star(G)\ar@{-->}[dl]^{q} \\
  & \R_+ &
 }
$$

\item[$2^\circ$.]\label{DEF:p^max} For each seminorm $p\in{\mathcal P}(N)$ the supremum of all these extensions
$$
p^{\max}(\alpha)=\sup_{q\in{\mathcal P}(G):\ q|_{{\mathcal C}^\star(N)}=p}q(\alpha)
$$
(is finite and) is a continuous $C^*$-seminorm on ${\mathcal C}^\star(G)$;

\item[$3^\circ$.]\label{PROP:p^max-3} If $p\in{\mathcal P}(N)$, $q\in{\mathcal P}(G)$, and $q\Big|_{{\mathcal C}^\star(N)}\le p$, then $q\le p^{\max}$.

\item[$4^\circ$.] For each $p_1,...,p_n\in{\mathcal P}(N)$
\beq\label{max(p_1,...,p_n)^max}
\max\{p_1^{\max},...,p_n^{\max}\}\le\Big(\max\{p_1,...,p_n\}\Big)^{\max}.
\eeq

}\eit
\bpr These propositions belong to Yu.~N.~Kuznetsova, and we give the proofs from \cite{Kuznetsova}.

1. Suppose $p:{\mathcal C}^\star(N)\to\R_+$ is a continuous $C^*$-seminorm. It is the norm of some norm-continuous representation $\pi:{\mathcal C}^\star(N)\to{\mathcal L}(X)$,
$$
p(\alpha)=\norm{\dot{\pi}(\alpha)}.
$$
By Theorem \ref{TH:Kuz-1}, the induced representation $\dot{T}:{\mathcal C}^\star(G)\to{\mathcal L}(L_2(D,X))$ is also norm-continuous. Hence, the norm
$$
q(\beta)=\norm{\dot{T}(\beta)},\qquad \beta\in{\mathcal C}^\star(G),
$$
is continuous on ${\mathcal C}^\star(G)$. And it extends $p$.

2. Put $D=G/N$. For each coset $R\in D=G/N$ we denote by $\eta_R$ its characteristic function on $G$
$$
\eta_R(s)=\begin{cases}1,& s\in R\\ 0,& s\notin R\end{cases},
$$
Since $N$ is open in $G$, $\eta_R\in{\mathcal C}(G)$. Besides this, in the notation of \eqref{eq10.5} we have:
$$
\eta_{tN}=\eta_N\cdot t^{-1}= t^{-1}\cdot\eta_N,\qquad t\in G.
$$
Foe each measure $\alpha\in{\mathcal C}^\star(G)$ we denote by $\alpha_R$ its ``restriction'' to $R$:
$$
\alpha_R(u)=\alpha(\eta_R\cdot u),\qquad u\in{\mathcal C}(G).
$$
The chain
\begin{multline*}
(\delta^{t^{-1}}*\alpha_{tN})(u)=(t^{-1}\cdot\alpha_{tN})(u)=(t^{-1}\cdot\alpha_{tN})(u)=
\alpha_{tN}(u\cdot t^{-1})=\alpha(\eta_{tN}\cdot (u\cdot t^{-1}))=\\=\alpha((\eta_N\cdot t^{-1})\cdot (u\cdot t^{-1}))=
\alpha((\eta_N\cdot u)\cdot t^{-1})=(t^{-1}\cdot \alpha)(\eta_N\cdot u)=(t^{-1}\cdot \alpha)_N(u)
\end{multline*}
implies $\delta^{t^{-1}}*\alpha_{tN}=(t^{-1}\cdot \alpha)_N$, and then
$$
\supp (\delta^{t^{-1}}*\alpha_{tN})\subseteq N
$$
Let us assign to each coset $R\in G/N$ a representative $t_R\in R$. Then $t_RN=R$.

Now for a given seminorm  $p\in{\mathcal P}(N)$ we consider the number
$$
Q(\alpha)=\sum_{R\in G/N}p\Big(\delta^{t_R^{-1}}\cdot\alpha_R\Big).
$$
The sum in the right side is finite, since $\alpha\in{\mathcal C}^\star(G)$ has compact support, and therefore among the shifts $t_R^{-1}\cdot\alpha$ by the elements of cosets $t_R\in R\in G/N$, where $N$ is an open subgroup, only finite number have the support that intersects with $N$:
$$
\card\{R\in G/N:\ (t_R^{-1}\cdot\alpha)_N\ne 0\}<\infty
$$
Thus, $Q(\alpha)\in\R_+$ for each $\alpha\in{\mathcal C}^\star(G)$, and we obtain that $Q$ is a seminorm on ${\mathcal C}^\star(G)$. On the other hand, for each measure $\alpha\in{\mathcal C}^\star(G)$ with the support in some coset,
$$
\supp \alpha\subseteq S\in G/N,
$$
we have
$$
Q(\alpha)=\sum_{R\in G/N}p\Big((t_R^{-1}\cdot\alpha)_N\Big)=p\Big((t_S^{-1}\cdot\alpha)_N\Big)=p\Big(t_S^{-1}\cdot\alpha\Big).
$$
This implies that the mapping $Q$ is continuous on each coset ${\mathcal C}^\star(R)$, since it is the composition of a shift $\alpha\mapsto t_R^{-1}\cdot\alpha$ and a continuous mapping $p$. On the other hand, ${\mathcal C}^\star(G)$, as a locally convex space, is a direct sum of the spaces ${\mathcal C}^\star(R)$, $R\in G/N$, hence $Q:{\mathcal C}^\star(G)\to\R_+$ is a continuous seminorm (not necessarily a $C^*$-seminorm).

Note that for a given seminorm $q\in{\mathcal P}(G)$ and for each unitary element $\upsilon\in{\mathcal C}^\star(G)$ we have $q(\upsilon*\alpha)=q(\alpha)$. In particular, $q((\delta^{t^{-1}}*\alpha_R))=q(\alpha_R)$. If in addition $q|_{{\mathcal C}^\star(N)}=p$, then
$$
q(\alpha)=q\l\sum_{R\in G/N}\alpha_R\r\le \sum_{R\in G/N}q(\alpha_R)=\sum_{R\in G/N}q\Big(\delta^{t^{-1}}*\alpha_R\Big)=\sum_{R\in G/N}p\Big(\delta^{t^{-1}}*\alpha_R\Big)=Q(\alpha)
$$
We can conclude that each continuous $C^*$-extension $q$ of the seminorm $p$ from ${\mathcal C}^\star(N)$ to ${\mathcal C}^\star(G)$ is subordinated to $Q$. This means that the supremum of all such extensions is also subordinated to $Q$,
$$
p^{\max}(\alpha)=\sup_{q\in{\mathcal P}(G):\ q|_{{\mathcal C}^\star(N)}=p}q(\alpha)\le Q(\alpha)
$$
and $Q$ here is a continuous seminorm. Thus $p^{\max}$ is a continuous seminorm on ${\mathcal C}^\star(G)$. Obviously, it is a $C^*$-seminorm.

3. Put $r=\max\{q,p^{\max}\}$. Clearly, $r\in{\mathcal P}(G)$, and
$$
r=\max\{q,p^{\max}\}\ge p^{\max}.
$$
On the other hand, $r|_{{\mathcal C}^\star(N)}=p$, and, by the already proven property (ii),
$$
r\le p^{\max}.
$$
We obtain that $r=\max\{q,p^{\max}\}=p^{\max}$, hence $q\le p^{\max}$.

4. Take $p_1,...,p_n\in{\mathcal P}(N)$, and put
$$
p=\max\{p_1,...,p_n\},\qquad q=\max\{p_1^{\max},...,p_n^{\max}\}.
$$
Then
$$
q|_{{\mathcal C}^\star(N)}=\max\{p_1^{\max},...,p_n^{\max}\}|_{{\mathcal C}^\star(N)}=\max\{p_1,...,p_n\}=p
$$
and, by $3^\circ$, $q\le p^{\max}$.
\epr

Consider the natural isomorphism of algebras \cite[Theorem 8.4]{Akbarov}
$$
{\mathcal C}^\star(\R^n)\circledast{\mathcal C}^\star(K)\cong {\mathcal C}^\star(\R^n\times K).
$$
For each character $\chi:\R^n\to\C^\times$ and for each representation $\sigma\in\widehat{K}$ the formula
$$
p_{\chi,\sigma}(\alpha\circledast\beta)=\abs{\alpha(\chi)}\cdot\norm{\dot{\sigma}(\beta)},\qquad \alpha\in
{\mathcal C}^\star(\R^n),\quad\beta\in{\mathcal C}^\star(K)
$$
uniquely defines a $C^*$-seminorm $p_{\chi,\sigma}:{\mathcal C}^\star(\R^n\times K)\to\R_+$.

If now $T$ is a compact set of characters on $\R^n$ (i.e. a compact set in the Pontryagin dual group $\widehat{\R^n}\cong\R^n$), then a $C^*$-seminorm is defined
$$
p_{T,\sigma}(\eta)=\sup_{\chi\in T}p_{\chi,\sigma}(\eta),\qquad \eta\in{\mathcal C}^\star(\R^n\times K).
$$
By the property $2^\circ$ on page \pageref{DEF:p^max}, it generates a seminorm
$$
p_{T,\sigma}^{\max}:{\mathcal C}^\star(G)\to\R_+.
$$
And from the property $3^\circ$ we have

\blm\footnote{See Errata on page \pageref{Errata}.}\label{LM:p_(T,S)-konf-sistema}
For each SIN-group $G$ the seminorms of the form
\beq\label{p_(T,S)=max}
p_{T,S}=\max_{\sigma\in S}p_{T,\sigma}^{\max},
\eeq
where $S$ runs over the system of finite subsets in $\widehat{K}$, and $T$ over the system compact subsets in $\widehat{\R^n}\cong\R^n$, form a cofinal system among all $C^*$-seminorms on ${\mathcal C}^\star(G)$, and the quotient algebras by these seminorms have the form
\beq\label{C^star(G)/p_(T,S)}
{\mathcal C}^\star(G)/\max_{\sigma\in S}p_{T,\sigma}^{\max}=\prod_{\sigma\in S}{\mathcal C}^\star(G)/p_{T,\sigma}^{\max}.
\eeq
\elm
\bpr
The first part if this proposition (the claim that $p_{T,S}$ form a cofinal system among all $C^*$-seminorms on ${\mathcal C}^\star(G)$) follows from Property $3^\circ$ (on page \pageref{PROP:p^max-3}), so we only need to prove the equality \eqref{C^star(G)/p_(T,S)}. Take a seminorm $p_{T,S}$ and consider the quotient map $\e_{T,S}:{\mathcal C}^\star(G)\to {\mathcal C}^\star(G)/p_{T,S}$. Since ${\mathcal C}^\star(G)/p_{T,S}$ is a $C^*$-algebra, it can be embedded (as a  $C^*$-algebra) into an algebra of the form ${\mathcal B}(X)$, where $X$ is a Hilbert space. Let $\rho:{\mathcal C}^\star(G)/p_{T,S}\to {\mathcal B}(X)$ be such an embedding. Consider the composition
$$
\xymatrix %@R=2.5pc @C=4.0pc
{
{\mathcal C}^\star(G)\ar[r]_{\e_{T,S}} \ar@/^3ex/[rr]^{\dot{\pi}} & {\mathcal C}^\star(G)/p_{T,S}\ar[r]_{\rho}& {\mathcal B}(X)
}
$$
By Property $3^\circ$ on page \pageref{X=oplus_pi-X_pi} the space $X$ is decomposed into a Hilbert sum  \eqref{X=oplus_pi-X_pi},
$$
X=\dot{\bigoplus_{\sigma\in\widehat{K}}}X_\sigma,
$$
and the intertwinner between the representations $\dot{\pi}$ and $\dot{\pi}_\sigma$ (where $\dot{\pi}_\sigma$ is defined in \eqref{DEF:T_pi}), is the operator $\varPhi_\sigma:X\to X_\sigma$ defined by the formula \eqref{DEF:varPhi_sigma}:
$$
\varPhi_\sigma=\dim X_\sigma\cdot\int_K \tr\sigma(s^{-1})\cdot \dot{\pi}(\delta^s)\cdot\mu_K(\d s)=\dot{\pi}(\nu_\sigma)
$$
($\nu_\sigma$ is defined in \eqref{DEF:nu_sigma}).

Note that since the seminorm $p_{T,S}^{\max}$ coincides with the seminorm $p_{T,S}$ on the subalgebra ${\mathcal C}^\star(K)\subseteq {\mathcal C}^\star(G)$, these seminorms coincide on the measures $\nu_\sigma\in {\mathcal C}^\star(K)$:
    \beq\label{p(nu_pi)=p^max(nu_pi)}
    p_{T,S}^{\max}(\nu_\sigma)=p_{T,S}(\nu_\sigma),\qquad \sigma\in\widehat{K}.
    \eeq
This implies that for $\tau\notin S$
$$
p_{T,S}^{\max}(\nu_\tau)=p_{T,S}(\nu_\tau)=0,
$$
hence
$$
\dot{\pi}(\nu_\tau)=0,\qquad \tau\notin S.
$$
As a corollary, all the spaces $X_\tau$ with the indices $\tau$ not lying in $S$, vanish,
$$
X_\tau=0\qquad \tau\notin S,
$$
and
$$
X=\dot{\bigoplus_{\sigma\in S}}X_\sigma.
$$
We see that the homomorphism $\dot{\pi}: {\mathcal C}^\star(G)\to {\mathcal B}(X)$ can be lifted to a homomorphism
$$
\dot{\pi}: {\mathcal C}^\star(G)\to\bigoplus_{\sigma\in S} {\mathcal B}(X_\sigma).
$$
We need to verify that its (closed) image coincides with the direct sum of the (closed) images of the homomorphisms  $\dot{\pi}_\sigma$:
$$
\overline{\dot{\pi}\Big({\mathcal C}^\star(G)\Big)}=\bigoplus_{\sigma\in S}\overline{\dot{\pi}_\sigma\Big({\mathcal C}^\star(G)\Big)}
$$
The direct inclusion is obvious
$$
\overline{\dot{\pi}\Big({\mathcal C}^\star(G)\Big)}\subseteq\bigoplus_{\sigma\in S}\overline{\dot{\pi}_\sigma\Big({\mathcal C}^\star(G)\Big)},
$$
and we need to prove the reverse one:
$$
\overline{\dot{\pi}\Big({\mathcal C}^\star(G)\Big)}\supseteq\bigoplus_{\sigma\in S}\overline{\dot{\pi}_\sigma\Big({\mathcal C}^\star(G)\Big)}.
$$
Take the element
$$
b\in \bigoplus_{\sigma\in S}\overline{\dot{\pi}_\sigma\Big({\mathcal C}^\star(G)\Big)}.
$$
and denote by $b_\sigma$ its component in ${\mathcal B}(X_\sigma)$:
$$
b=\sum_{\sigma\in S}b_\sigma.
$$

Take $\e>0$. Let $n=\card S$ be the cardinality of the (finite) set $S$. For each $\sigma\in S$ there exists a measure $\alpha_\sigma\in {\mathcal C}^\star(G)$ such that
\beq\label{PROOF:LM:p_(T,S)-konf-sistema}
p_{T,\sigma}(b_\sigma-\dot{\pi}_\sigma(\alpha_\sigma))<\frac{\e}{n}.
\eeq
Let us take such a family $\{\alpha_\sigma;\ \sigma\in S\}$ and put
$$
\alpha=\sum_{\sigma\in S}\nu_\sigma*\alpha_\sigma.
$$
Then
\beq\label{nu_pi*alpha=nu_pi*alpha_pi}
\nu_\tau*\alpha=\sum_{\sigma\in S}\nu_\tau*\nu_\sigma*\alpha_\sigma=\eqref{chi_pi-ortog-proj}=\nu_\tau*\alpha_\tau,
\qquad \tau\in S,
\eeq
and therefore
\begin{multline*}
p_{T,S}(b-\dot{\pi}(\alpha))=
p_{T,S}(\sum_{\sigma\in S}b_\sigma-\dot{\pi}(\sum_{\sigma\in S}\nu_\sigma*\alpha))=
p_{T,S}(\sum_{\sigma\in S}b_\sigma-\sum_{\sigma\in S}\dot{\pi}(\nu_\sigma*\alpha))=\eqref{nu_pi*alpha=nu_pi*alpha_pi}=\\=
p_{T,S}(\sum_{\sigma\in S}b_\sigma-\sum_{\sigma\in S}\dot{\pi}(\nu_\sigma*\alpha_\sigma))=
p_{T,S}\Big(\sum_{\sigma\in S}b_\sigma-\sum_{\sigma\in S}\dot{\pi}_\sigma(\alpha_\sigma)\Big)=
p_{T,S}\Big(\sum_{\sigma\in S}(b_\sigma-\dot{\pi}_\sigma(\alpha_\sigma)\big)=\\=
p_{T,S}\Big(\sum_{\sigma\in S}\big(b_\sigma-\dot{\pi}_\sigma(\alpha_\sigma)\big)\Big)\le
\sum_{\sigma\in S}p_{T,S}\Big(b_\sigma-\dot{\pi}_\sigma(\alpha_\sigma)\Big)=
\sum_{\sigma\in S}p_{T,\sigma}\Big(b_\sigma-\dot{\pi}_\sigma(\alpha_\sigma)\Big)=\\=
\sum_{\sigma\in S}p_{T,\sigma}\Big(b-\dot{\pi}_\sigma(\alpha_\sigma)\Big)
<\eqref{PROOF:LM:p_(T,S)-konf-sistema}<
\sum_{\sigma\in S}\frac{\e}{n}=\e.
\end{multline*}
\epr

\bprop\footnote{See Errata on page \pageref{Errata}.}\label{PROP:nepr-Env-SIN-gruppy}
For each representation \eqref{SIN-kak-rasshirenie} of a SIN-group $G$ as the discrete extension of some group $\R^n\times K$ the continuous envelope of the group algebra ${\mathcal C}^\star(G)$ as a stereotype algebra is a direct product
\beq\label{nepr-Env-SIN-gruppy}
\Env_{\mathcal C}{\mathcal C}^\star(G)=\prod_{\sigma\in\widehat{K}}{\mathcal C}_\sigma^\star(G),
\eeq
where the factors
$$
{\mathcal C}_\sigma^\star(G)=\leftlim_{T\subseteq\widehat{\R^n}}{\mathcal C}^\star(G)/p_{T,\sigma}^{\max}
$$
are Fr\'echet algebras.
\eprop
\bpr
Note first that the product in the right side of \eqref{nepr-Env-SIN-gruppy} coincides with the Kuznetsova envelope (see \cite{Akbarov-env}), i.e. with the projective product of the system of $C^*$-quotient mappings:
\begin{multline}\label{SIN:Ste-lim=LCS-lim}
\leftlim_{p\in{\mathcal P}(G)}{\mathcal C}^\star(G)/p=\text{(Lemma \ref{LM:p_(T,S)-konf-sistema})}=\leftlim_{T\subseteq\widehat{\R^n}}\leftlim_{S\subseteq\widehat{K}}{\mathcal C}^\star(G)/p_{T,S}=\eqref{p_(T,S)=max}=
\leftlim_{T\subseteq\widehat{\R^n}}\leftlim_{S\subseteq\widehat{K}}{\mathcal C}^\star(G)/\max_{\sigma\in S}p_{T,\sigma}^{\max}=\\=\eqref{C^star(G)/p_(T,S)}=
\leftlim_{T\subseteq\widehat{\R^n}}\leftlim_{S\subseteq\widehat{K}}\prod_{\sigma\in S}{\mathcal C}^\star(G)/p_{T,\sigma}^{\max}=\leftlim_{T\subseteq\widehat{\R^n}}\prod_{\sigma\in\widehat{K}}{\mathcal C}^\star(G)/p_{T,\sigma}^{\max}=\prod_{\sigma\in\widehat{K}}\leftlim_{T\subseteq\widehat{\R^n}}{\mathcal C}^\star(G)/p_{T,\sigma}^{\max}=\prod_{\sigma\in\widehat{K}}{\mathcal C}_\sigma^\star(G)
\end{multline}
Apart from this the system of compact sets $T\subseteq\widehat{\R^n}$ has a countable cofinal subsystem $T_k$, hence each factor
$$
{\mathcal C}_\sigma^\star(G)=\leftlim_{T\subseteq\widehat{\R^n}}{\mathcal C}^\star(G)/p_{T,\sigma}^{\max}=\leftlim_{\infty\gets k}{\mathcal C}^\star(G)/p_{T_k,\sigma}^{\max}
$$
is in fact a countable projective limit of Banach algebras, and therefore a Fr\'echet algebra.

Note finally that under each projection ${\mathcal C}^\star(G)\to{\mathcal C}_\sigma^\star(G)$ the (set-theoretic) image of the space ${\mathcal C}^\star(G)$ is dense in ${\mathcal C}_\sigma^\star(G)$. This implies that the image of ${\mathcal C}^\star(G)$ is dense in the product $\prod_{\sigma\in\widehat{K}}{\mathcal C}_\sigma^\star(G)$. Therefore the continuous envelope of the algebra ${\mathcal C}^\star(G)$ coincide with $\prod_{\sigma\in\widehat{K}}{\mathcal C}_\sigma^\star(G)$ (see \cite[(3.62)]{Akbarov-env}).
\epr

\bprop\footnote{See Errata on page \pageref{Errata}.}\label{PROP:nepr-Env-SIN-gruppy=LCS-lim} For SIN-groups $G$ the envelope of the group algebra ${\mathcal C}^\star(G)$ coincides with the Kuznetsova envelope, i.e. with the projective limit of its $C^*$-quotient algebras in the category of locally convex spaces (and in the category of topological algebras):
\beq\label{E(C*(G))=LCS-leftlim}
\Env_{\mathcal C}{\mathcal C}^\star(G)={\tt LCS}\text{-}\kern-10pt\leftlim_{p\in{\mathcal P}({\mathcal C}^\star(G))}{\mathcal C}^\star(G)/p
\eeq
\eprop
\bpr
The chain \eqref{SIN:Ste-lim=LCS-lim} holds in the category of locally convex spaces, and the factor in the end ${\mathcal C}_\sigma^\star(G)=\leftlim_{T\subseteq\widehat{\R^n}}{\mathcal C}^\star(G)/p_{T,\sigma}^{\max}$, being a projective limit of a sequence of Banach spaces, is a Fr\'echet space, hence a stereotype space. Therefore, nothing is changed after replacing the category $\tt{LCS}$ with the category $\tt{Ste}$.
\epr

\bprop\footnote{See Errata on page \pageref{Errata}.}\label{PROP:E(theta)}
The continuous envelope $\Env_{\mathcal C}(\theta):\Env_{\mathcal C}({\mathcal C}^\star(N))\to \Env_{\mathcal C}({\mathcal C}^\star(G))$ of a morphism $\theta:{\mathcal C}^\star(N)\to {\mathcal C}^\star(G)$ is an open and closed mapping of stereotype spaces (in the sence of \cite{Akbarov-env}).
\eprop
\bpr
The openness follows from the property $1^\circ$ on page \pageref{1^0:prodolzhenie-polunorm}. After that we have to note that this mapping turns each component ${\mathcal C}^\star_\sigma(N)$ into the component ${\mathcal C}^\star_\sigma(G)$, and this is an open mapping of Fr\'echet spaces. This implies the closedness of this map and the closedness of the whole $\Env_{\mathcal C}(\theta)$.
\epr

\paragraph{The continuous envelopes of the group algebra of distributions ${\mathcal E}^\star(G)$.}

If $G$ is a Lie group, then (apart from the group algebra ${\mathcal C}^\star(G)$ of measures with compact support)  we can consider the group algebra ${\mathcal E}^\star(G)$ of dictributions with compact support on $G$ (see \cite{Akbarov}). If we denote the natural inclusion ${\mathcal E}(G)\subseteq{\mathcal C}(G)$ by some symbol, say, by $\lambda$, we obtain the following diagram:
\beq\label{Env_C-C^star(G)-Env_C-E^star(G)}
 \xymatrix @R=2.5pc @C=4pc
 {
{\mathcal C}^\star(G)\ar[r]^{\lambda^\star}\ar[d]_{\env_{\mathcal C}{\mathcal C}^\star(G)} & {\mathcal E}^\star(G)\ar[d]^{\env_{\mathcal C}{\mathcal E}^\star(G)}
\\
\Env_{\mathcal C}{\mathcal C}^\star(G)\ar[r]^{\Env_{\mathcal C}(\lambda^\star)} & \Env_{\mathcal C}{\mathcal E}^\star(G)
 }
\eeq

\btm\label{TH:Env_C-C^star(G)=Env_C-E^star(G)}
For each real Lie group $G$ the continuous envelopes of group algebras ${\mathcal C}^\star(G)$ and ${\mathcal E}^\star(G)$ coincide:
\beq\label{Env_C-C^star(G)=Env_C-E^star(G)}
\Env_{\mathcal C}{\mathcal C}^\star(G)=\Env_{\mathcal C}{\mathcal E}^\star(G).
\eeq
\etm
\bpr
Consider the composition
$$
 \xymatrix @R=2.5pc @C=4pc
 {
{\mathcal C}^\star(G)\ar[r]^{\lambda^\star} & {\mathcal E}^\star(G)\ar[r]^{\env_{\mathcal C}{\mathcal E}^\star(G)}&
\Env_{\mathcal C}{\mathcal E}^\star(G)
 }
$$
and let us verify that it is a continuous extension of the algebra ${\mathcal C}^\star(G)$. Indeed, if $\ph:{\mathcal C}^\star(G)\to B$ is a morphism into a $C^*$-algebra $B$, then it generates a (norm) continuous representation  $\pi=\ph\circ\delta:G\to B$, which due to Example \ref{EX:gladkost-nepr-po-norme-predstavleniya}, is smooth, and therefore can be extended to some homomorphism $\ph'=\ddot{\pi}:{\mathcal E}^\star(G)\to B$. This homomorphism  $\ddot{\pi}$ extends $\ph$, since they coincide on delta-functionals $\delta^x$, $x\in G$, which are total in ${\mathcal C}^\star(G)$ and in ${\mathcal E}^\star(G)$. The homomorphism $\ph'$ can be uniquely extended to a homomorphism $\ph''$ on the algebra $\Env_{\mathcal C}{\mathcal E}^\star(G)$ (which is a continuous extension of  ${\mathcal E}^\star(G)$):
$$
 \xymatrix @R=3pc @C=4pc
 {
{\mathcal C}^\star(G)\ar@/_3ex/[dr]_{\ph}\ar[r]^{\lambda^\star} & {\mathcal E}^\star(G)\ar[r]^{\env_{\mathcal C}{\mathcal E}^\star(G)}\ar@{-->}[d]^{\ph'=\ddot{\pi}} &
\Env_{\mathcal C}{\mathcal E}^\star(G)\ar@/^3ex/@{-->}[dl]^{\ph''}\\
& B & \\
 }
$$
This proves that $\Env_{\mathcal C}{\mathcal E}^\star(G)$ is a continuous extension of the algebra ${\mathcal C}^\star(G)$. Hence, there exists a unique morphism $\upsilon$ from $\Env_{\mathcal C}{\mathcal E}^\star(G)$ into the continuous envelope of ${\mathcal C}^\star(G)$ such that the following diagram is commutative:
$$
 \xymatrix @R=2.5pc @C=4pc
 {
& {\mathcal C}^\star(G)\ar@/_3ex/[dl]_{\lambda^\star}\ar@/^3ex/[rd]^{\quad\env_{\mathcal C}{\mathcal C}^\star(G)\circ\lambda^\star} &
\\
\Env_{\mathcal C}{\mathcal C}^\star(G) & & \Env_{\mathcal C}{\mathcal E}^\star(G)\ar@{-->}[ll]_{\upsilon}
 }
$$
From the fact that ${\mathcal C}^\star(G)$ is dense in $\Env_{\mathcal C}{\mathcal C}^\star(G)$ and in $\Env_{\mathcal C}{\mathcal E}^\star(G)$, it follows that $\upsilon$ is the inverse mrophism for $\Env_{\mathcal C}(\lambda^\star)$ from \eqref{Env_C-C^star(G)-Env_C-E^star(G)}.
\epr

\subsection{The algebra ${\mathcal K}(G)$}

For each locally compact group $G$ its group algebra of measures ${\mathcal C}^\star(G)$ is an involutive Hopf algebra with respect to the stereotype tensor product $\circledast$. By Theorems \ref{TH:C-obolochka-sohranyaet-Hopfov} and \ref{TH:C-obolochka-sohranyaet-inv-Hopfov}, its continuous envelope $\Env_{\mathcal C}{\mathcal C}^\star(G)$ is a coalgebra with the interconsistent antipode and involution on the categories ${\mathcal C}\text{-}{\tt Alg}$ of continuous algebras and in the category $({\tt Ste},\odot)$ of stereotype spaces. Denote by ${\mathcal K}(G)$ the stereotype dual space to the space $\Env_{\mathcal C}{\mathcal C}^\star(G)$:
 \beq\label{DEF:K(G)}
{\mathcal K}(G):=\Big(\Env_{\mathcal C}{\mathcal C}^\star(G)\Big)^\star.
 \eeq
This is the dual space to a coalgebra in $({\tt Ste},\odot)$ with interconsistent antipode and involution, hence by Property $4^\circ$ on page \pageref{4^0:inv-v-sopryazh-alg}, the following theorem is true:

\btm For each locally compact group $G$ the space ${\mathcal K}(G)$ is an algebra in the category $({\tt Ste},\circledast)$ (i.a. a stereotype algebra) with the interconsistent antipode and involution.
\etm

By Theorem \ref{TH:C-obolochka-sohranyaet-inv-Hopfov}(ii) the morphism
\beq\label{DEF:env_C^star}
\big(\env_{\mathcal C}{{\mathcal C}^\star(G)}\big)^\star:{\mathcal K}(G)=\Big(\Env_{\mathcal C}{\mathcal C}^\star(G)\Big)^\star\to {\mathcal C}^\star(G)^\star={\mathcal C}(G),
\eeq
dual to the morphism of envelope, is an involutive homomorphism of algebras:
\beq\label{env_C^star-homomorphism}
\env_{\mathcal C}^\star1=1,\qquad  \env_{\mathcal C}^\star(u\cdot v)=\env_{\mathcal C}^\star(u)\cdot\env_{\mathcal C}^\star(v),\qquad\env_{\mathcal C}^\star \overline{u}=\overline{\env_{\mathcal C}^\star u},
\qquad u,v\in {\mathcal K}(G)
\eeq

Let us show that it has trivial kernel.\label{Ker-e_C-C^star(G)^star=0} Indeed, the group $G$ is embeded into the algebra $\Env_{\mathcal C}{\mathcal C}^\star(G)$ as a composition of the delta-mapping and the envelope
$$
G\overset{\delta}{\longrightarrow}{\mathcal C}^\star(G)\overset{\env_{\mathcal C}{{\mathcal C}^\star(G)}}{\longrightarrow}\Env_{\mathcal C}{\mathcal C}^\star(G).
$$
The image of $G$ is total\footnote{Definition on page \pageref{DEF:overline^X}.} in ${\mathcal C}^\star(G)$ (by \cite[Lemma 8.2]{Akbarov}), and the image of ${\mathcal C}^\star(G)$ is dense in $\Env_{\mathcal C}{\mathcal C}^\star(G)$. Hence the image of $G$ is total in $\Env_{\mathcal C}{\mathcal C}^\star(G)$. This implies that each element $u\in {\mathcal K}(G):=\Big(\Env_{\mathcal C}{\mathcal C}^\star(G)\Big)^\star$ is uniquely defined by the composition
$$
u\circ \env_{\mathcal C}{{\mathcal C}^\star(G)}\circ\delta:G\to\C,
$$
which can be treated as the restriction of $u$ on $G$. In particular, if this composition vanishes in ${\mathcal C}(G)$, then $u=0$ in ${\mathcal K}(G)$.

The important for us conclusion is that ${\mathcal K}(G)$ can be perceived as an involutive subalgebra in ${\mathcal C}(G)$:

\btm\label{TH:K(G)->C(G)} The mapping $u\mapsto u\circ \env_{\mathcal C}{{\mathcal C}^\star(G)}\circ\delta$
coincides with the mapping $\env_{\mathcal C}{{\mathcal C}^\star(G)}^\star$, dual to $\env_{\mathcal C}{{\mathcal C}^\star(G)}$:
\beq\label{K(G)->C(G)}
\env_{\mathcal C}{{\mathcal C}^\star(G)}^\star(u)=u\circ \env_{\mathcal C}{{\mathcal C}^\star(G)}\circ\delta
\eeq
and injectively and homomorphically embeds ${\mathcal K}(G)$ into ${\mathcal C}(G)$ as an involutive subalgebra (and therefore the operations of summing, multiplication and involution on ${\mathcal K}(G)$ are pointwise).
\etm

\btm\label{TH:K(G)=lim}
The algebra ${\mathcal K}(G)$ as a stereotype space can be represented as the nodal coimage (in the category of stereotype spaces)
 \beq\label{K(G)=lim}
{\mathcal K}(G)=\Coim_{\infty}\ph^\star
 \eeq
of the mapping $\ph^\star$, dual to the natural morphism of stereotype spaces
$$
\ph:{\mathcal C}^\star(G)\to \leftlim_{p\in{\mathcal P}({\mathcal C}^\star(G))}{\mathcal C}^\star(G)/p.
$$
\etm
\bpr
Consider the diagram \eqref{DIAGR:C-envelope=im_infty-lim N_X} for $A={\mathcal C}^\star(G)$
$$
\xymatrix @R=3.pc @C=6.0pc % @M=14pt
{
{\mathcal C}^\star(G)\ar[d]_{\coim_\infty\ph}\ar[r]^{\ph=\kern-10pt\leftlim_{p\in{\mathcal P}({\mathcal C}^\star(G))}\kern-10pt\rho_p} & \leftlim_{p\in{\mathcal P}({\mathcal C}^\star(G))} {\mathcal C}^\star(G)/p &   \\
\Coim_\infty\ph\ar[r]_{\red_\infty\ph} &  \Im_\infty \ph \ar[u]_{\im_\infty\ph}\ar@{=}[r] & \Env_{\mathcal C}{\mathcal C}^\star(G)
}
$$
The dual diagram is
\beq\label{DIAGR:K(G)}
\xymatrix @R=2.pc @C=6.0pc % @M=14pt
{
{\mathcal C}(G) & \Big(\leftlim_{p\in{\mathcal P}({\mathcal C}^\star(G))} {\mathcal C}^\star(G)/p\Big)^\star \ar[l]_{\ph^\star} \ar[d]^{\coim_\infty\ph^\star} & \rightlim_{p\in{\mathcal P}({\mathcal C}^\star(G))}\Big({\mathcal C}^\star(G)/p\Big)^\star \ar@{=}[l]\\
\Im_\infty\ph^\star\ar[u]^{\im_\infty\ph^\star}
 &  \Coim_\infty \ph^\star\ar[l]^{\red_\infty\ph^\star} & {\mathcal K}(G)\ar@{=}[l]
}
\eeq
\epr

Let us recall spaces $\Trig(G)$ and $k(G)$, defined on pages \pageref{DEF:k(G)} and \pageref{DEF:Trig(G)}.

\btm\footnote{See Errata on page \pageref{Errata}.}\label{TH:V-k-K-C} The following chain of set-theoretic inclusions hold,
 \beq\label{V-k-K-C}
\Trig(G)\subseteq k(G)\subseteq {\mathcal K}(G)\subseteq {\mathcal C}(G).
 \eeq
Therewith,
\bit{
\item[(i)] always
 \beq\label{Trig-plotno-v-K}
\overline{\Trig(G)}={\mathcal K}(G),
 \eeq
\item[(ii)] if $G$ is a SIN-group, then
 \beq\label{k=K}
k(G)={\mathcal K}(G)
 \eeq
and
 \beq\label{K->C}
\overline{\Trig(G)}=\overline{k(G)}=\overline{{\mathcal K}(G)}={\mathcal C}(G).
 \eeq

\item[(iii)] if $G$ is a compact group, then
 \beq\label{Trig=K}
\Trig(G)=k(G)={\mathcal K}(G)
 \eeq
 }\eit
\etm
\bpr
1. In the chain \eqref{V-k-K-C} the first inclusion is obvious, and the third one we already noticed in Theorem \ref{TH:K(G)->C(G)}. Let us prove the second one: $k(G)\subseteq {\mathcal K}(G)$. Take $u\in k(G)$, i.e. a function satisfying \eqref{DEF:k(G)}, where $\pi:G\to{\mathcal B}(X)$ is a norm continuous unitary representation. By Theorem \ref{TH:nepr-po-norme-predstavleniya} $\pi$ generates a (continuous) homomorphism of involutive algebras $\psi:{\mathcal C}^\star(G)\to{\mathcal B}(X)$. Here ${\mathcal B}(X)$ is a $C^*$-algebra, therefore $\psi$ has an extension $\psi':\Env_{\mathcal C}{\mathcal C}^\star(G)\to{\mathcal B}(X)$. Consider the functional
$$
f(\beta)=\langle \psi'(\beta)x,y\rangle,\qquad \beta\in \Env_{\mathcal C}{\mathcal C}^\star(G)
$$
on ${\mathcal K}^\star(G)$. It generates a function on $G$, coinciding with $u$:
$$
\xymatrix @R=3.pc @C=6.0pc % @M=14pt
{
G\ar@/_4ex/[ddr]_{u}\ar[r]^{\delta}\ar[dr]_{\pi} & {\mathcal C}^\star(G)\ar[r]^{\env_{\mathcal C}{{\mathcal C}^\star(G)}}
\ar@{-->}[d]_{\psi} & \Env_{\mathcal C}{\mathcal C}^\star(G)\ar@{-->}[dl]^{\psi'}\ar@/^4ex/@{-->}[ddl]^{f} \\
& {\mathcal B}(X)\ar[d] & \\
& \C &
}
$$
This means that $u\in {\mathcal K}(G)$.

2. Let us show that $\Trig(G)$ is dense in ${\mathcal K}(G)$. Each algebra ${\mathcal C}^\star(G)/p$ can be isometrically embedded into an algebra of the form ${\mathcal B}(X)$:
$$
{\mathcal C}^\star(G)/p\to {\mathcal B}(X).
$$
Note that, first, by the Hahn-Banach theorem each functional $f\in ({\mathcal C}^\star(G)/p)^\star$ can be extended to some functional $g\in{\mathcal B}(X)^\star$, and, second, each functional $g\in{\mathcal B}(X)^\star$ can be approximated in ${\mathcal B}(X)^\star$ by linear combinations of pure states, i.e. (Example \ref{EX:trig-monom}) by functions from $\Trig(G)$. This means that in the dual mapping
$$
({\mathcal C}^\star(G)/p)^\star\gets{\mathcal B}(X)^\star
$$
the functionals lying in $\Trig(G)$, turn into a dense subset in $({\mathcal C}^\star(G)/p)^\star$. This is true for each seminorm $p\in{\mathcal P}({\mathcal C}^\star(G))$, hence we obtain, that the functionals from
$$
\Big(\leftlim_{p\in{\mathcal P}({\mathcal C}^\star(G))} {\mathcal C}^\star(G)/p\Big)^\star=\rightlim_{p\in{\mathcal P}({\mathcal C}^\star(G))}\Big({\mathcal C}^\star(G)/p\Big)^\star,
$$
which generate functions in $\Trig(G)$, are dense in
$$
\Big(\leftlim_{p\in{\mathcal P}({\mathcal C}^\star(G))} {\mathcal C}^\star(G)/p\Big)^\star=\rightlim_{p\in{\mathcal P}({\mathcal C}^\star(G))}\Big({\mathcal C}^\star(G)/p\Big)^\star.
$$
On the other hand, in Diagram \eqref{DIAGR:K(G)} it is seen that
$$
\Big(\leftlim_{p\in{\mathcal P}({\mathcal C}^\star(G))} {\mathcal C}^\star(G)/p\Big)^\star=\rightlim_{p\in{\mathcal P}({\mathcal C}^\star(G))}\Big({\mathcal C}^\star(G)/p\Big)^\star
$$
is densely mapped into ${\mathcal K}(G)$. Together this means that $\Trig(G)$ is dense in ${\mathcal K}(G)$.

3. Let $G$ be a SIN-group. By Proposition \ref{PROP:nepr-Env-SIN-gruppy=LCS-lim}, the mapping $\Env_{\mathcal C}{\mathcal C}^\star(G)$ coincides with the locally convex projective limit of the quotient algebras ${\mathcal C}^\star(G)/p$:
$$
\Env_{\mathcal C}{\mathcal C}^\star(G)={\tt{LCS}}\text{-}\kern-10pt\leftlim_{p\in{\mathcal P}({\mathcal C}^\star(G))}{\mathcal C}^\star(G)/p.
$$
This implies that the dual space is a locally convex injective limit of the spaces $({\mathcal C}^\star(G)/p)^\star$
$$
{\mathcal K}(G)=\Env_{\mathcal C}{\mathcal C}^\star(G)^\star={\tt{LCS}}\text{-}\kern-10pt\rightlim_{p\in{\mathcal P}({\mathcal C}^\star(G))}({\mathcal C}^\star(G)/p)^\star.
$$
We obtain, that any function $u\in{\mathcal K}(G)$ is generated by some functional $f\in({\mathcal C}^\star(G)/p)^\star$. But on the other hand, each such a functional is a linear combination of a finite set of states on a $C^*$-algebra \cite[4.3.7]{Kad-Ring}, hence $u$ is a sum of positive definite functions. By Theorem  \ref{TH:u(t)=langle-pi(t)x,x-rangle} each positive definite function subordinated to a $C^*$-seminorm, lies in $k(G)$. This proves \eqref{k=K}. Further, if $G$ is a SIN-group, then the algebra $k(G)$ separates points of $G$. In addition it contains the unit, and hence, all constants. Thus, by the Stone-Weierstrass theorem, the restriction of $k(G)$ on each compact $K\subseteq G$ is dense in ${\mathcal C}(K)$. Therefore, $k(G)={\mathcal K}(G)$ is dense in ${\mathcal C}(G)$.

4. Proposition (iii) follows from (ii) and from Theorem \ref{TH:Trig(G)=k(G)-na-kompaktah}.
\epr

\paragraph{The mapping ${\mathcal K}(G)\circledast{\mathcal K}(H)\to {\mathcal K}(G\times H)$.}

Let $G$ and $H$ be two locally compact groups. To each pair of functions $u\in{\mathcal C}(G)$ and $v\in{\mathcal C}(H)$ we assign a function of the Cartesian product $G\times H$:
\beq\label{u-boxdot-v}
(u\boxdot v)(s,t)=u(s)\cdot v(t),\qquad s\in G,\ t\in H.
\eeq
As is known \cite[Theorem 8.4]{Akbarov-env}, there is a unique linear continuous mapping
$$
\iota:{\mathcal C}(G)\odot {\mathcal C}(H)\to {\mathcal C}(G\times H)
$$
that satisfies the identity
$$
\iota(u\odot v)=u\boxdot v, \qquad u\in{\mathcal C}(G),\quad v\in{\mathcal C}(H).
$$
($\iota$ is an isomorphism of stereotype algebras). Put
\beq\label{omega_G,H}
\omega_{G,H}=\eta_{{\mathcal K}^\star(G),{\mathcal K}^\star(H)}\circ \Env_{\mathcal C}(\env_{\mathcal C}{{\mathcal C}^\star(G)}\circledast \env_{\mathcal C}{{\mathcal C}^\star(H)})\circ \Env_{\mathcal C}(\iota^\star).
\eeq
The action of this mapping is illustrated by the following diagram:
{\scriptsize
$$
\xymatrix @R=2.pc @C=4.0pc % @M=14pt
{
{\mathcal C}^\star(G\times H)\ar[rr]^{\env_{\mathcal C}{{\mathcal C}^\star(G\times H)}}\ar[dd]_{\iota^\star} & &
\Env_{\mathcal C}{\mathcal C}^\star(G\times H)\ar[dd]^{\Env_{\mathcal C}(\iota^\star)} & {\mathcal K}^\star(G\times H)\ar@{=}[l]
\ar@{-->}[dddd]^{\omega_{G,H}} \\ &&& \\
{\mathcal C}^\star(G)\circledast{\mathcal C}^\star(H)\ar[rr]^{\env_{\mathcal C}\big({\mathcal C}^\star(G)\circledast {\mathcal C}^\star(G)\big)}\ar[dd]|{\env_{\mathcal C}{{\mathcal C}^\star(G)}\circledast \env_{\mathcal C}{{\mathcal C}^\star(H)}} & & \Env_{\mathcal C}\big({\mathcal C}^\star(G)\circledast{\mathcal C}^\star(H)\big)\ar[dd]|{\Env_{\mathcal C}\big(\env_{\mathcal C}{{\mathcal C}^\star(G)}\circledast \env_{\mathcal C}{{\mathcal C}^\star(H)}\big)} & \\
&&& \\
\Env_{\mathcal C}{\mathcal C}^\star(G)\circledast \Env_{\mathcal C}{\mathcal C}^\star(H) \ar[rr]^(.45){\env_{\mathcal C}\big(\Env_{\mathcal C}{\mathcal C}^\star(G)\circledast \Env_{\mathcal C}{\mathcal C}^\star(H)\big)}
& & \Env_{\mathcal C}(\Env_{\mathcal C}{\mathcal C}^\star(G)\circledast \Env_{\mathcal C}{\mathcal C}^\star(H))\ar[r]^(.6){\eta_{{\mathcal K}^\star(G),{\mathcal K}^\star(H)}} & {\mathcal K}^\star(G)\odot {\mathcal K}^\star(H) \\
&& {\mathcal K}^\star(G)\overset{\mathcal C}{\circledast} {\mathcal K}^\star(H)\ar@{=}[u] & \\
{\mathcal K}^\star(G)\circledast {\mathcal K}^\star(H)\ar@{=}[uu]\ar@/_18ex/[uurrr]_(.8){\quad @_{{\mathcal K}^\star(G),{\mathcal K}^\star(H)}} & &  & \\
 &  &
}
$$
}
Consider the dual mapping:
$$
\xymatrix @R=2.pc @C=8.0pc % @M=14pt
{
{\mathcal K}(G\times H) & {\mathcal K}(G)\circledast{\mathcal K}(H)
\ar[l]_{\omega_{G,H}^\star}
}
$$

\btm\label{TH:omega^star(u-circledast-v)=u-boxdot-v}
For each two locally compact groups $G$ and $H$ the following identity holds:
\beq\label{omega^star(u-circledast-v)=u-boxdot-v}
\omega_{G,H}^\star(u\circledast v)=u\boxdot v,\qquad u\in{\mathcal K}(G),\quad v\in{\mathcal K}(H).
\eeq
\etm
\bpr
Consider the diagram
{\scriptsize
$$
\xymatrix @R=4.pc @C=4.0pc % @M=14pt
{
{\mathcal C}(G\times H) & &
\big(\Env_{\mathcal C}{\mathcal C}^\star(G\times H)\big)^\star\ar[ll]_{\big(\env_{\mathcal C}{\mathcal C}^\star(G\times H)\big)^\star} & {\mathcal K}(G\times H)\ar@{=}[l]
 \\
{\mathcal C}(G)\odot{\mathcal C}(H)\ar[u]^{\iota}
& & \Env_{\mathcal C}\big({\mathcal C}^\star(G)\circledast{\mathcal C}^\star(H)\big)^\star\ar[ll]_{\Big(\env_{\mathcal C}\big({\mathcal C}^\star(G)\circledast {\mathcal C}^\star(G)\big)\Big)^\star}\ar[u]_{\Env_{\mathcal C}(\iota^\star)^\star}
 & \\
(\Env_{\mathcal C}{\mathcal C}^\star(G))^\star\odot (\Env_{\mathcal C}{\mathcal C}^\star(H))^\star
 \ar[u]|{\big(\env_{\mathcal C}{{\mathcal C}^\star(G)}\circledast \env_{\mathcal C}{{\mathcal C}^\star(H)}\big)^\star}
& & \Env_{\mathcal C}(\Env_{\mathcal C}{\mathcal C}^\star(G)\circledast \Env_{\mathcal C}{\mathcal C}^\star(H))^\star\ar[u]|{\Env_{\mathcal C}\big(\env_{\mathcal C}{{\mathcal C}^\star(G)}\circledast \env_{\mathcal C}{{\mathcal C}^\star(H)}\big)^\star}\ar[ll]_{\env_{\mathcal C}{\big(\Env_{\mathcal C}{\mathcal C}^\star(G)\circledast \Env_{\mathcal C}{\mathcal C}^\star(H)\big)}^\star}
 & {\mathcal K}(G)\circledast {\mathcal K}(H)\ar[l]_(.35){\eta_{{\mathcal K}^\star(G),{\mathcal K}^\star(H)}^\star} \ar@{-->}[uu]_{\omega_{G,H}^\star}\ar@/^10ex/[dlll]^{@_{{\mathcal K}(G),{\mathcal K}(H)}} \\
{\mathcal K}(G)\odot {\mathcal K}(H)\ar@{=}[u] & &  & \\
 &  &
}
$$
}
When moving from the right lower corner to the left upper one, the elementary tensor $u\circledast v$, on the one hand, follows the way
$$
\xymatrix @R=2.pc @C=8.0pc % @M=14pt
 {
 u\circledast v\ar[r]^{@} & u\odot v\ar[r]^{(\env_{\mathcal C}{{\mathcal C}^\star(G)}\circledast \env_{\mathcal C}{{\mathcal C}^\star(H)})^\star} & u\odot v\ar[r]^{\iota}& u\boxdot v
 }
$$
and on the other, the way
$$
\xymatrix @R=2.pc @C=8.0pc % @M=14pt
 {
 u\circledast v\ar[r]^{\omega_{G,H}^\star} & \omega_{G,H}^\star(u\odot v)\ar[r]^(.35){\big(\env_{\mathcal C}{\mathcal C}^\star(G\times H)\big)^\star} & \big(\env_{\mathcal C}{\mathcal C}^\star(G\times H)\big)^\star(\omega_{G,H}^\star(u\odot v))
 }
$$
Hence,
$$
\big(\env_{\mathcal C}{\mathcal C}^\star(G\times H)\big)^\star(\omega_{G,H}^\star(u\odot v))=u\boxdot v,
$$
and since by Theorem \ref{TH:K(G)->C(G)}, $\big(\env_{\mathcal C}{\mathcal C}^\star(G\times H)\big)^\star$ can be treated as a set-theoretic inclusion, we have
$$
\omega_{G,H}^\star(u\odot v)=u\boxdot v.
$$
\epr

\btm\label{TH:K(A-times-K)-cong-K(A)-circledast-K(K)}
If $A$ is an Abelian locally compact group, and $K$ a compact group, then the mapping $\omega_{A,K}:{\mathcal K}(A)\circledast {\mathcal K}(K)\to {\mathcal K}(A\times K)$ is an isomorphism:
\beq\label{K(A-times-K)-cong-K(A)-circledast-K(K)}
{\mathcal K}(A\times K)\cong{\mathcal K}(A)\circledast {\mathcal K}(K)
\eeq
\etm
\bpr
This follows from \eqref{env_C^star(R^n-times-K)}:
\begin{multline*}
{\mathcal K}(A\times K)=\Big(\Env_{\mathcal C}{\mathcal C}^\star(A\times K)\Big)^\star\cong \eqref{env_C^star(R^n-times-K)}\cong \Big(\Env_{\mathcal C}{\mathcal C}^\star(A)\odot\Env_{\mathcal C}{\mathcal C}^\star(K) \Big)^\star\cong\\ \cong \Big(\Env_{\mathcal C}{\mathcal C}^\star(A)\Big)^\star\circledast\Big(\Env_{\mathcal C}{\mathcal C}^\star(K)\Big)^\star\cong{\mathcal K}(A)\circledast {\mathcal K}(K)
\end{multline*}
\epr

\paragraph{The shift in ${\mathcal K}(G)$.}

\btm\label{LM:sdvig-v-K(G)}
The shift (either left or right) by an arbitrary element $a\in G$ is an isomorphism of the stereotype algebra ${\mathcal K}(G)$.
\etm
\bpr
We prove this for the operator of the left shift: take
$$
M_a:{\mathcal C}^\star(G)\to{\mathcal C}^\star(G)\quad\Big|\quad  M_a(\alpha)=\delta^a*\alpha,\qquad \alpha\in {\mathcal C}^\star(G).
$$
Denote by $\eta_a=\env_{\mathcal C}(\delta^a)$, and put
$$
N_a:\Env_{\mathcal C}{\mathcal C}^\star(G)\to
\Env_{\mathcal C}{\mathcal C}^\star(G)\quad\Big|\quad N_a(\omega)=\eta_a*\omega,\qquad \omega\in \Env_{\mathcal C}{\mathcal C}^\star(G)
$$
(here $*$ is the multiplication in the algebra $\Env_{\mathcal C}{\mathcal C}^\star(G)$). Then
$$
\env_{\mathcal C}(M_a(\alpha))=\env_{\mathcal C}(\delta^a*\alpha)=\env_{\mathcal C}(\delta^a)*\env_{\mathcal C}(\alpha)=N_a(\env_{\mathcal C}(\alpha)),
$$
i.e. the foolowing diagram is commutative:
\beq\label{PROOF:LM:sdvig-v-K(G)}
\xymatrix @R=2.pc @C=8.0pc % @M=14pt
{
{\mathcal C}^\star(G)\ar[r]^{\env_{\mathcal C}}\ar[d]_{M_a} & \Env_{\mathcal C}{\mathcal C}^\star(G)\ar[d]^{N_a}\\
{\mathcal C}^\star(G)\ar[r]^{\env_{\mathcal C}} & \Env_{\mathcal C}{\mathcal C}^\star(G)
}
\eeq
As a corollary, the dual diagram is commutative:
$$
\xymatrix @R=2.pc @C=8.0pc % @M=14pt
{
{\mathcal C}(G) & {\mathcal K}(G)\ar[l]_{\env_{\mathcal C}^\star}\\
{\mathcal C}(G)\ar[u]^{M_a^\star} & {\mathcal K}(G)\ar[l]_{\env_{\mathcal C}^\star}\ar[u]_{N_a^\star}
}
$$
It implies, first, that $N_a^\star$ is an operator of the left shift by the element $a$ in the functional algebra ${\mathcal K}(G)$, since under the homomorphical injection $\env_{\mathcal C}^\star$ (described in Theorem \ref{TH:K(G)->C(G)}) it turns into the operator $M_a^\star$ of the left shift by the element $a\in G$ in the algebra  ${\mathcal C}(G)$. And, second, $N_a^\star$ is a homomorphism of algebras:
\begin{multline*}
\env_{\mathcal C}^\star N_a^\star(u\cdot v)=M_a^\star\Big(\env_{\mathcal C}^\star(u\cdot v)\Big)=\eqref{env_C^star-homomorphism}=
M_a^\star\Big(\env_{\mathcal C}^\star(u)\cdot\env_{\mathcal C}^\star(v)\Big)=\\=
M_a^\star\env_{\mathcal C}^\star(u)\cdot M_a^\star\env_{\mathcal C}^\star(v)=
\env_{\mathcal C}^\star(N_a^\star u)\cdot\env_{\mathcal C}^\star(N_a^\star v)=\eqref{env_C^star-homomorphism}=
\env_{\mathcal C}^\star(N_a^\star u \cdot N_a^\star v)
\end{multline*}
and, since by Theorem \ref{TH:K(G)->C(G)}, $\env_{\mathcal C}^\star$ is an injective map,
$$
N_a^\star(u\cdot v)=N_a^\star u \cdot N_a^\star v.
$$
\epr

\subsection{Continuous duality for Moore groups}

Recall that we defined Moore groups on page \pageref{DEF:Moore-gruppa} as those for which all unitary irreducible representations $\pi:G\to{\mathcal L}(X)$ are finite dimensional. By Theorem \ref{TH:Moore->SIN}, every such a group is a SIN-group, therefore it has the representation \eqref{SIN-kak-rasshirenie},
$$
1\to \R^n\times K=N\to G\to D\to 1
$$
where $K$ is compact and $D$ discrete. By Corollary \ref{TH:G=Moore->D=Moore} $D$ is also a Moore group.

\paragraph{Density of the mapping $\omega_{G.H}^\star:{\mathcal K}(G)\circledast{\mathcal K}(H)\to{\mathcal K}(G\times H)$.}

\btm\label{TH:u-boxdot-v-polno-v-K(G-times-R^n)}
For each two Moore groups $G$ and $H$ the functions of the form $u\boxdot v$ (definied in \eqref{u-boxdot-v}), where $u\in{\mathcal K}(G)$ and $v\in{\mathcal K}(H)$, are dense in ${\mathcal K}(G\times H)$.
\etm
\bpr
Consider a function $w$ of the form \eqref{DEF:k(G)} on the group $G\times H$. Since $G\times H$ is also a Moore group, the irreducible representation $\pi:G\times H\to{\mathcal B}(X)$ must be finite dimensional. Let $e_1,...,e_n$ be an orthonormalized basis in $X$, and suppose $\rho:G\to{\mathcal B}(X)$ and $\sigma:H\to{\mathcal B}(X)$ are representations acting by formulas
$$
\rho(s)=\pi(s,1_H),\qquad \sigma(t)=\pi(1_G,t),\qquad s\in G,\quad t\in H.
$$
Then
\begin{multline*}
w(s,t)=\Big\langle\pi(s,t)x,y\Big\rangle=\Big\langle\pi\Big((s,1_H)\cdot(1_G,t)\Big)x,y\Big\rangle=
\Big\langle\Big(\pi(s,1_H)\cdot\pi(1_G,t)\Big)x,y\Big\rangle=\Big\langle\Big(\rho(s)\cdot\sigma(t)\Big)x,y\Big\rangle=\\=
\langle\sigma(t)x,\rho(s)^\bullet y\rangle=\sum_{i=1}^n \langle\sigma(t)x,e_i\rangle\cdot\langle e_i,\rho(s)^\bullet y\rangle=
\sum_{i=1}^n \langle\sigma(t)x,e_i\rangle\cdot\langle \rho(s)e_i,y\rangle=\sum_{i=1}^n (u_i\boxdot v_i)(s,t)
\end{multline*}
where
$$
u_i(s)=\langle \rho(s)e_i,y\rangle,\qquad v_i(t)=\langle\sigma(t)x,e_i\rangle.
$$
This implies that the space $\Trig(G\times H)$ of norm-continuous trigonometric polynomials on $G$ is contained in the linear span of the functions $u\boxdot v$, where $u\in{\mathcal K}(G)$ and $v\in{\mathcal K}(H)$. Hence,
$$
\overline{\sp\{u\boxdot v; \ u\in{\mathcal K}(G), \ v\in{\mathcal K}(H)\}}\supseteq\overline{\Trig(G\times H)}=\eqref{Trig-plotno-v-K}={\mathcal K}(G\times H).
$$
\epr

\paragraph{The spectrum and the continuous envelope of the algebra ${\mathcal K}(G)$ for Moore groups.}

The last inclusion in \eqref{V-k-K-C}
$$
{\mathcal K}(G)\subseteq{\mathcal C}(G)
$$
means by the way that the spectra of these algebras are connected by a natural continuous mapping
$$
\Spec{\mathcal K}(G)\gets\Spec{\mathcal C}(G)=G.
$$
In the case when $G$ is a Moore group, this mapping is a homeomorphism:

\btm\footnote{See Errata on page \pageref{Errata}.}\label{TH:Spec-K(G)=G}
If $G$ is a Moore group, then the involutive spectrum of the algebra ${\mathcal K}(G)$ is topologically isomorphic to $G$:
\beq\label{Spec-K(G)=G}
\Spec{\mathcal K}(G)=G
\eeq
\etm

We premise the proof with 8 lemmas.

\blm\label{LM:Spec-K(G)-dlya-diskretnoi-gruppy}
If $G$ is an amenable discrete group, then the mapping of spectra $G\to \Spec{\mathcal K}(G)$ is a bijection.
\elm
\bpr We use here the ideas from \cite[Theorem 6.3]{Kuznetsova}.
By Proposition \ref{PROP:nepr-obol-diskr-gruppy}, the continuous envelope of the group algebra for $G$ is the $C^*$-algebra of this group:
$$
\Env_{\mathcal C}{\mathcal C}^\star(G)=C^*(G)
$$
Its (Banach) dual space is the classical Fourier-Stiltjes algebra $B(G)$, hence ${\mathcal K}(G)$ coincides as a set with $B(G)$ \cite{Eymard,Renault}:
$$
{\mathcal K}(G)=\Env_{\mathcal C}{\mathcal C}^\star(G)^\star=C^*(G)^\star=B(G).
$$
(this is the equality of vector spaces, but the topology on $B(G)$ is stronger that on ${\mathcal K}(G)$). The algebra $B(G)$ contains the ideal $A(G)$ called the Fourier algebra \cite{Eymard,Renault}:
$$
{\mathcal K}(G)=B(G)\supseteq A(G).
$$
From the amenability of $G$ it follows that the annihilator of the algebra $A(G)$ in the space $C^*(G)$ vanishes \cite{Hulanicki}, \cite[Theorem 4.21]{Paterson}:
$$
A(G)^\perp=0
$$
As a corollary, $A(G)$ is dense in ${\mathcal K}(G)$:
$$
\overline{A(G)}={\mathcal K}(G).
$$
If $\chi:{\mathcal K}(G)\to\C$ is an involutive character, then it is an involutive character on $A(G)$, and by Eymard theorem \cite{Eymard}, it acts on functions $u\in A(G)$ as the delta-functional in a point $a\in G$:
$$
\chi(u)=u(a),\qquad A(G).
$$
Since $A(G)$ is dense in ${\mathcal K}(G)$, this is true for all functions $u\in{\mathcal K}(G)$.
\epr

\blm\label{LM:Spec-K(G)-dlya-kompaktnoi-gruppy}
If $G$ is a compact group, then the mapping of spectra $G\to \Spec{\mathcal K}(G)$ is a homeomorphism.
\elm
\bpr
Each involutive character $\chi:{\mathcal K}(G)\to\C$ is an involutive character on the algebra $\Trig(G)$ of trigonnometric polynomials on $G$, hence \cite[30.5]{Hewitt-Ross-2} there exists a point $a\in G$ such that
$$
\chi(u)=u(a),\qquad u\in\Trig(G).
$$
By \eqref{Trig-plotno-v-K} the algebra $\Trig(G)$ is dense in the algebra ${\mathcal K}(G)$. Hence,
$$
\chi(u)=u(a),\qquad u\in{\mathcal K}(G).
$$
We see that the mapping $G\to \Spec{\mathcal K}(G)$ is surjective. On the other hand, each two points $s\ne t$ on a compact group $G$ are separated by a trigonometric polynomial $u\in\Trig(G)$, and this means that the mapping $G\to \Spec{\mathcal K}(G)$ is injective. Thus, it is bijective and continuous. At the same time $G$ is compact. Therefore this mapping is a homeomorphism.
\epr

\blm\label{LM:Spec-K(G)-dlya-kommut-gruppy}
If $G$ is an Abelian locally compact group, then the mapping of spectra $G\to \Spec{\mathcal K}(G)$ is a homeomorphism.
\elm
\bpr
By Proposition \ref{PROP:Env_C-C^star(G)=C(widehat(G))}, $\Env_{\mathcal C}{\mathcal C}^\star(G)={\mathcal C}(\widehat{G})$.
Hence
$$
{\mathcal K}(G)={\mathcal C}^\star(\widehat{G}).
$$
Each character $f:{\mathcal C}^\star(\widehat{G})\to\C$ in the composition with the mapping $\delta:\widehat{G}\to {\mathcal C}^\star(\widehat{G})$ gives a complex character $f\circ\delta:\widehat{G}\to\C^\times$, and if $f$ is involutive, then $f\circ\delta$ acts into the circle $\T$. I.e., $f\circ\delta$ is a character on the group  $\widehat{G}$, hence
$$
(f\circ\delta)(\chi)=f(\delta^\chi)=\chi(a),\qquad \chi\in \widehat{G}
$$
for some point $a\in G$. Since linear combinations of delta-functions $\delta^\chi$ are dense in ${\mathcal C}^\star(\widehat{G})$, this identity can be extended to all elements ${\mathcal C}^\star(\widehat{G})$:
$$
f(u)=u(a),\qquad u\in{\mathcal K}(G)={\mathcal C}^\star(\widehat{G})
$$
We see that the mapping $G\to \Spec{\mathcal K}(G)=\Spec{\mathcal C}^\star(\widehat{G})$ is surjective. On the other hand, any two points $s\ne t$ on $G$, i.e. two characters on $\widehat{G}$ can be separated by some delta-functional $\delta_\chi\in{\mathcal C}^\star(\widehat{G})$, $\chi\in \widehat{G}$,
$$
\delta^\chi(s)=\chi(s)\ne\chi(t)=\delta^\chi(t)
$$
and this means that the mapping of spectra $G\to \Spec{\mathcal K}(G)=\Spec{\mathcal C}^\star(\widehat{G})$ is an injection. Thus, it is bijective and continuous. It remains to verify that it is open (i.e. continuous in the reverse direction). Suppose $f_i\to f$ in $\Spec{\mathcal C}^\star(\widehat{G})$, and $a_i$, $a$ are the corresponding points in $G$. If $K$  is compact in $\widehat{G}$, then the set $\{\delta^\chi;\ \chi\in K\}$ is compact in ${\mathcal C}^\star(\widehat{G})$, hence $f_i(\delta^\chi)$ tends to $f(\delta^\chi)$ uniformly on $\chi\in K$:
$$
\chi(a_i)=f_i(\delta^\chi)\underset{i\to\infty}{\underset{\chi\in K}{\rightrightarrows}}f(\delta^\chi)=\chi(a).
$$
This is true for each compact set $K\subseteq \widehat{G}$, hence $a_i\to a$ in $G$.
\epr

\blm\label{LM:Spec(R^n-times-K)}
For each compact group $K$ and each $n\in\N$ the mapping of spectra $\R^n\times K\to \Spec{\mathcal K}(\R^n\times K)$ is a homeomorphism.
\elm
\bpr
Let $\chi:{\mathcal K}(\R^n\times K)\to\C$ be an involutive character. Consider the embeddings
$$
\rho:{\mathcal K}(\R^n)\to {\mathcal K}(\R^n\times K),\qquad \rho(u)=u\boxdot 1_K
$$
and
$$
\sigma:{\mathcal K}(K)\to {\mathcal K}(\R^n\times K),\qquad \rho(v)=1_{\R^n}\boxdot v.
$$
The composition $\chi\circ\rho$ is an involutive character on the algebra ${\mathcal K}(\R^n)$. By Lemma  \ref{LM:Spec-K(G)-dlya-kommut-gruppy} $\Spec{\mathcal K}(\R^n)=\R^n$, hence there is a point $a\in\R^n$ such that
$$
\chi(u\boxdot 1_K)=(\chi\circ\rho)(u)=u(a), \qquad u\in {\mathcal K}(\R^n)
$$
On the other hand, the composition $\chi\circ\sigma$ is an involutive character on the algebra ${\mathcal K}(K)$, which contains the algebra $\Trig(G)$, and by Lemma \ref{LM:Spec-K(G)-dlya-kompaktnoi-gruppy} there is a point $b\in K$ such that
$$
\chi(1_{\R^n}\boxdot v)=(\chi\circ\sigma)(v)=v(b), \qquad v\in {\mathcal K}(K).
$$
Now for functions of the form $u\boxdot v$ we have
$$
\chi(u\boxdot v)=\chi(u\boxdot 1_K\cdot 1_{\R^n}\boxdot v)=\chi(u\boxdot 1_K)\cdot\chi(1_{\R^n}\boxdot v)=u(a)\cdot v(b)=(u\boxdot v)(a,b).
$$
This identity is extended to the linear combinations of the functions of the form $u\boxdot v$, and then by Theorem  \ref{TH:u-boxdot-v-polno-v-K(G-times-R^n)}, to the whole space ${\mathcal K}(\R^n\times K)$. This proves the surjectivity of the mapping $\R^n\times K\to \Spec{\mathcal K}(\R^n\times K)$.

Let us now prove its injectivity: take $(s,a)\ne(t,b)$ in $\R^n\times K$. then either $s\ne t$, or $a\ne b$. In the first case one can find the function $u\in{\mathcal K}(\R^n)$ such that $u(s)\ne u(t)$ (at this momemnt we use the Lemma \ref{LM:Spec-K(G)-dlya-kommut-gruppy}), and this means that the function $u\boxdot 1$ (that belongs to  ${\mathcal K}(\R^n\times K)$ by Theorem \ref{TH:omega^star(u-circledast-v)=u-boxdot-v}) separates the points $(s,a)$ and $(t,b)$:
$$
(u\boxdot 1)(s,a)=u(s)\cdot 1=u(s)\ne u(t)=u(t)\cdot 1=(u\boxdot 1)(t,b).
$$
And in the second case, when $a\ne b$, we can do the same, and we find (using Lemma \ref{LM:Spec-K(G)-dlya-kompaktnoi-gruppy}) a function from ${\mathcal K}(\R^n)$, that separates these points.

We see that the mapping of spectra $\R^n\times K\to \Spec{\mathcal K}(\R^n\times K)$ is bijective (and continuous).

Let us prove its openness (the continuity in the reverse direction). Suppose $(s_i,a_i)\to (s,a)$ in $\Spec{\mathcal K}(\R^n\times K)$. Take a compact set $T\subseteq {\mathcal K}(\R^n)$ and consider its image  $\rho(T)\subseteq{\mathcal K}(\R^n\times K)$ which is also compact, therefore we have the uniform convergence by $u\in T$:
$$
u(s_i)=u(s_i)\cdot 1= (u\boxdot 1)(s_i,a_i)=\rho(u)(s_i,a_i)\underset{i\to\infty}{\underset{u\in T}{\rightrightarrows}}\rho(u)(s,a)=(u\boxdot 1)(s,a)=u(s)\cdot 1=u(s).
$$
This is true for any compact set $T\subseteq {\mathcal K}(\R^n)$, therefore $s_i\to s$ in $\Spec{\mathcal K}(\R^n)$. By Lemma \ref{LM:Spec-K(G)-dlya-kommut-gruppy} this means that $s_i\to s$ in $\R^n$.

Similarly (using Lemma \ref{LM:Spec-K(G)-dlya-kompaktnoi-gruppy}) we prove that $a_i\to a$ in $K$.
\epr

\blm\label{LM:1_L-in-K(G)}
Let $G$ be a SIN-group. Then for each coset $L\in G/N$ its characteristic function  $1_L$ is an element of the space ${\mathcal K}(G)$:
\beq\label{1_L-in-K(G)}
1_L\in {\mathcal K}(G)
\eeq
\elm
\bpr
Consider the trivial representation $\pi:N\to\C$, $\pi(t)=1$. Its induced representation
$\pi':G\to {\mathcal L}(L_2(D))$ is defined by formula \eqref{ind-representation}, which in this case has the form
$$
\pi'(g)(\xi)(t)=\xi\big(\ph(\sigma(t)\cdot g)\big)=
\xi\big(t\cdot\ph(g)\big),\qquad \xi\in L_2(D),\quad t\in D,\quad g\in G.
$$
Replace $\xi$ by the characteristic function of the unit $1_D$ in $D$:
$$
\xi(t)=\begin{cases}1, & t=1_D\\ 0, & t\ne 1_D\end{cases}
$$
Then
$$
f(g)=\langle\pi'(g)(\xi),\xi\rangle=\sum_{t\in D}\pi'(g)(\xi)(t)\cdot\overline{\xi(t)}=
\sum_{t\in D}\xi\big(t\cdot\ph(g)\big)\cdot \xi(t)=\begin{Bmatrix}1, & \ph(g)=1\\ 0, & \ph(g)\ne 1\end{Bmatrix}=
\begin{Bmatrix}1, & g\in N\\ 0, & g\notin N\end{Bmatrix}.
$$
I.e. the characteristic function $1_N$ belongs to ${\mathcal K}(G)$. From Theorem \ref{LM:sdvig-v-K(G)} it follows that all its shifts also belong to ${\mathcal K}(G)$.
\epr

For each $L\in G/N$ we put
\beq\label{DEF:K_L(G)}
{\mathcal K}_L(G)=1_L\cdot{\mathcal K}(G),\qquad {\mathcal K}_{G\setminus L}(G)=(1-1_L)\cdot{\mathcal K}(G)
\eeq
(here 1 is the unit of the algebra ${\mathcal K}(G)$). From \eqref{1_L-in-K(G)} we have

\blm\label{LM:K_L(G)-K_G-L(G)} Let $G$ be a SIN-group. Then the spaces ${\mathcal K}_L(G)$ and ${\mathcal K}_{G\setminus L}(G)$ complement each other in ${\mathcal K}(G)$:
\beq\label{K_N+K_(G-N)=K(G)}
{\mathcal K}_L(G)\oplus {\mathcal K}_{G\setminus L}(G)={\mathcal K}(G)
\eeq
(i.e. ${\mathcal K}(G)$ is a direct sum in the category of stereotype spaces).
\elm

Let us consider the embedding of the group algebras $\theta:{\mathcal C}^\star(N)\to {\mathcal C}^\star(G)$, its envelope $\Env_{\mathcal C}(\theta):\Env_{\mathcal C}({\mathcal C}^\star(N))\to \Env_{\mathcal C}{\mathcal C}^\star(G)$, and its dual mapping $\psi=\Env_{\mathcal C}(\theta)^\star:{\mathcal K}(N)\gets {\mathcal K}(G)$.

\blm\footnote{See Errata on page \pageref{Errata}.}\label{LM:K(N)=K_N(G)} Suppose $G$ is a SIN-group. Then the morphism of stereotype spaces $\psi=\Env_{\mathcal C}(\theta)^\star:{\mathcal K}(N)\gets {\mathcal K}(G)$ has the following properties:
\bit{
\item[(i)] its kernel is the second component in the decomposition \eqref{K_N+K_(G-N)=K(G)} (with $L=N$):
\beq\label{ph(K_(G-N)(G))=0}
\Ker\psi={\mathcal K}_{G\setminus N}(G).
\eeq
\item[(ii)] the restriction $\psi|_{{\mathcal K}_N(G)}:{\mathcal K}_N(G)\to {\mathcal K}(N)$ is an isomorphism of stereotype algebras.
}\eit
\elm
\bpr
1. To prove (i) let us consider the diagram
$$
\xymatrix @R=2.pc @C=8.0pc % @M=14pt
{
{\mathcal C}^\star(N)\ar[r]^{\env_{\mathcal C}{{\mathcal C}^\star(N)}}\ar[d]_{\theta} & \Env_{\mathcal C}({\mathcal C}^\star(N))\ar[d]^{\Env_{\mathcal C}(\theta)}\\
{\mathcal C}^\star(G)\ar[r]^{\env_{\mathcal C}{{\mathcal C}^\star(G)}} & \Env_{\mathcal C}{\mathcal C}^\star(G)
}
$$
If $u\in \Ker\psi$, then we obtain the chain
\begin{multline*}
0=\psi(u)=u\circ \Env_{\mathcal C}(\theta)\quad\Longrightarrow\quad 0=u\circ \Env_{\mathcal C}(\theta)\circ \env_{\mathcal C}{{\mathcal C}^\star(N)}=
u\circ \env_{\mathcal C}{{\mathcal C}^\star(G)}\circ \theta \quad\Longrightarrow\\ \Longrightarrow\quad 0=u\circ \env_{\mathcal C}{{\mathcal C}^\star(G)}\circ \theta\circ\delta^N \quad\Longrightarrow\quad 0=u\Big|_N \quad\Longrightarrow\\ \Longrightarrow\quad 0=u\cdot 1_N
\quad\Longrightarrow\quad u=u\cdot (1-1_N) \quad\Longrightarrow\quad u\in {\mathcal K}_{G\setminus N}(G).
\end{multline*}
Conversely, if $u\in {\mathcal K}_{G\setminus N}(G)$, then
\begin{multline*}
0=\psi(u)=u\circ \Env_{\mathcal C}(\theta)\quad\overset{\env_{\mathcal C}{{\mathcal C}^\star(N)}\in\Epi}{\Longleftarrow}\quad 0=u\circ \Env_{\mathcal C}(\theta)\circ \env_{\mathcal C}{{\mathcal C}^\star(N)}=
u\circ \env_{\mathcal C}{{\mathcal C}^\star(G)}\circ \theta \quad\overset{\overline{\sp\delta^N}={\mathcal C}^\star(N)}{\Longleftarrow}\\ \Longleftarrow\quad 0=u\circ \env_{\mathcal C}{{\mathcal C}^\star(G)}\circ \theta\circ\delta^N \quad\Longleftarrow\quad 0=u\Big|_N \quad\Longleftarrow\\ \Longleftarrow\quad 0=u\cdot 1_N
\quad\Longleftarrow\quad \exists v\in {\mathcal K}(G)\quad u=v\cdot (1-1_N) \quad\Longleftarrow\quad u\in {\mathcal K}_{G\setminus N}(G).
\end{multline*}

2. Let us prove (ii). First, we have to show that the restriction $\psi|_{{\mathcal K}_N(G)}:{\mathcal K}_N(G)\to {\mathcal K}(N)$ is a bijection. Its injectivity is obvious, let us prove the surjectivity. Take $u\in{\mathcal K}(N)$, then by \eqref{k=K}
$$
u(t)=\langle\pi(t)x,y\rangle
$$
for some norm-continuous unitary representation $\pi:N\to{\mathcal B}(X)$, and $x,y\in X$. Consider the induced representation $\pi':N\to{\mathcal B}(L_2(D,X))$ \eqref{ind-representation} and put
\beq\label{xi-upsilon}
\xi(t)=\begin{cases}x,& t=1_D\\ 0,& t\ne 1_D\end{cases},\qquad \upsilon(t)=\begin{cases}y,& t=1_D\\ 0,& t\ne 1_D\end{cases}.
\eeq
Then the function
$$
w(g)=\langle\pi'(g)\xi,\upsilon\rangle
$$
coincides with $u$ on $N$ and vanishes on $G\setminus N$. Indeed, for $g\in N$ we have (here $\ph:G\to G/N=D$ is the quotient map)
\begin{multline*}
w(g)=\langle\pi'(g)(\xi),\upsilon\rangle=\sum_{t\in D}\langle\pi'(g)(\xi)(t),\upsilon(t)\rangle=\eqref{xi-upsilon}=
\langle\pi'(g)(\xi)(1_D),\upsilon(1_D)\rangle=\eqref{ind-representation}=\\=
\bigg\langle\pi\Big(\underbrace{\sigma(1_D)}_{\scriptsize\begin{matrix}\| \put(2,0){\eqref{sigma(1_D)=1_G}}\\ 1_G\end{matrix}}\cdot g\cdot\sigma\big(\ph(\underbrace{\sigma(1_D)}_{\scriptsize\begin{matrix}\| \put(2,0){\eqref{sigma(1_D)=1_G}}\\ 1_G\end{matrix}}\cdot g)\big)^{-1}\Big)\Big(\xi\big(\ph(\underbrace{\sigma(1_D)}_{\scriptsize\begin{matrix}\| \put(2,0){\eqref{sigma(1_D)=1_G}}\\ 1_G\end{matrix}}\cdot g)\big)\Big),y\bigg\rangle=\langle\pi(g\cdot\sigma\big(\underbrace{\ph(g)}_{\scriptsize\begin{matrix}\| \put(2,0){$(g\in N)$}\\ 1_D\end{matrix}}\big)^{-1})(\xi(\underbrace{\ph(g)}_{\scriptsize\begin{matrix}\| \put(2,0){$(g\in N)$}\\ 1_D\end{matrix}})),y\rangle=\\=
\langle\pi(g\cdot\underbrace{\sigma(1_D)^{-1}}_{\scriptsize\begin{matrix}\| \put(-32,0){\eqref{sigma(1_D)=1_G}}\\ 1_G\end{matrix}})(\underbrace{\xi(1_D)}_{\scriptsize\begin{matrix}\| \put(2,0){\eqref{xi-upsilon}}\\ x\end{matrix}}),y\rangle=\langle\pi(g)(x),y\rangle=u(g).
\end{multline*}
And if $g\notin N$, then
\begin{multline*}
w(g)=\langle\pi'(g)(\xi),\upsilon\rangle=\sum_{t\in D}\langle\pi'(g)(\xi)(t),\upsilon(t)\rangle=\eqref{xi-upsilon}=
\langle\pi'(g)(\xi)(1_D),\upsilon(1_D)\rangle=\eqref{ind-representation}=\\=
\bigg\langle\pi\Big(\underbrace{\sigma(1_D)}_{\scriptsize\begin{matrix}\| \put(2,0){\eqref{sigma(1_D)=1_G}}\\ 1_G\end{matrix}}\cdot g\cdot\sigma\big(\ph(\underbrace{\sigma(1_D)}_{\scriptsize\begin{matrix}\| \put(2,0){\eqref{sigma(1_D)=1_G}}\\ 1_G\end{matrix}}\cdot g)\big)^{-1}\Big)\Big(\xi\big(\ph(\underbrace{\sigma(1_D)}_{\scriptsize\begin{matrix}\| \put(2,0){\eqref{sigma(1_D)=1_G}}\\ 1_G\end{matrix}}\cdot g)\big)\Big),y\bigg\rangle=\\=
\bigg\langle\pi\Big(g\cdot\sigma\big(\ph(g)\big)^{-1}\Big)\Big(\underbrace{\xi\big(\overbrace{\ph(g)}^{\scriptsize\begin{matrix}1_D\\ \text{\rotatebox{90}{$\ne$}}\put(2,0){$(g\notin N)$}\end{matrix}}\big)}_{\scriptsize\begin{matrix}\| \put(2,0){\eqref{xi-upsilon}}\\ 0\end{matrix}}\Big),y\bigg\rangle=
\langle\pi(g\cdot\sigma(\ph(g))^{-1})(0),y\rangle=0.
\end{multline*}

3. Thus, we understood that the restriction $\psi|_{{\mathcal K}_N(G)}:{\mathcal K}_N(G)\to {\mathcal K}(N)$ is a (continuous and) bijective mapping, and it remains to verify that it is open. We should first note that the mapping $\psi=\Env_{\mathcal C}(\theta)^\star:{\mathcal K}(N)\gets {\mathcal K}(G)$ is open -- this follows from the fact that its predual mapping $\Env_{\mathcal C}(\theta):\Env_{\mathcal C}({\mathcal C}^\star(N))\to \Env_{\mathcal C}{\mathcal C}^\star(G)$ is closed by Proposition \ref{PROP:E(theta)}.

Suppose then that $U$ is a convex balanced neighbourhood of zero in ${\mathcal K}_N(G)$. Using  \eqref{K_N+K_(G-N)=K(G)} we can find a neighbourhood of zero $V$ in ${\mathcal K}(G)$ such that
\beq\label{V-cap-K_N(G)=U}
V\cap {\mathcal K}_N(G)=U,\qquad V+{\mathcal K}_{G\setminus N}(G)=V.
\eeq
Since $\psi$ is open, it maps the neighbourhood of zero $V$ into the neighbourhood of zero $W=\psi(V)$ in ${\mathcal K}(N)$. Let us show that $W=\psi(U)$. Since $V\supseteq U$, we have $W=\psi(V)\supseteq\psi(U)$. Let us prove the inverse inclusion. Take $w\in W=\psi(V)$, i.e. $w=\psi(v)$ for some $v\in V$. Since, as we have already noticed in 2, the restriction  $\psi|_{{\mathcal K}_N(G)}$ is bijective, there is an element $u\in {\mathcal K}_N(G)$ such that  $w=\psi(u)$. then $\psi(u-v)=\psi(u)-\psi(v)=w-w=0$, i.e.
$$
u-v\in \ker\psi=\eqref{ph(K_(G-N)(G))=0}={\mathcal K}_{G\setminus N}(G)
$$
hence
$$
u\in v+{\mathcal K}_{G\setminus N}(G)\subseteq V+{\mathcal K}_{G\setminus N}(G)=\eqref{V-cap-K_N(G)=U}=V
$$
and therefore,
$$
u\in V\cap {\mathcal K}_N(G)=\eqref{V-cap-K_N(G)=U}=U.
$$
\epr

Thus, the mapping $\psi|_{{\mathcal K}_N(G)}:{\mathcal K}_N(G)\to {\mathcal K}(N)$ is an isomorphism of stereotype spaces. Hence there is the inverse mapping $\psi:{\mathcal K}(N)\to {\mathcal K}_N(G)$, which is a morphism of stereotype spaces, and moreover, of stereotype algebras.

From \ref{LM:sdvig-v-K(G)} it follows that the shifts define morphisms of stereotype algebras
$$
\psi_L:{\mathcal K}(N)\to {\mathcal K}_L(G)=1_L\cdot{\mathcal K}(G)
$$

\bpr[Proof of Theorem \ref{TH:Spec-K(G)=G}.] Suppose $G$ is a Moore group.

1. Let us first show that the mapping of spectra $G\to \Spec{\mathcal K}(G)$ is a surjection. Suppose  $\chi:{\mathcal K}(G)\to\C$ is an involutive character. The homomorphism $G\to D$ from \eqref{SIN-kak-rasshirenie} generates a homomorphism ${\mathcal C}^\star(G)\to {\mathcal C}^\star(D)$, which then generates a homomorphism $\Env_{\mathcal C}{\mathcal C}^\star(G)\to \Env_{\mathcal C}({\mathcal C}^\star(D))$, and then a homomorphism ${\mathcal K}(G)\gets {\mathcal K}(D)$. We denote the last homomorphism by $\ph:{\mathcal K}(D)\to {\mathcal K}(G)$. The composition $\chi\circ\ph:{\mathcal K}(D)\to \C$ is an involutive character on ${\mathcal K}(D)$, and $D$ is a Moore group by Corollary \ref{TH:G=Moore->D=Moore}, hence an amenable group by Theorem \ref{TH:Moore->AM}. Thus by Lemma  \ref{LM:Spec-K(G)-dlya-diskretnoi-gruppy} $\chi$ is a delta-functional:
$$
(\chi\circ\ph)(u)=u(L),\qquad u\in {\mathcal K}(D),
$$
for some $L\in G/N$. Consider the space ${\mathcal K}_L(G)$ from \eqref{DEF:K_L(G)} and denote by $\rho_L$ its embedding into ${\mathcal K}(G)$. On the other hand, denote by $\sigma$ the embedding ${\mathcal K}(N)\to {\mathcal K}_N(G)$, i.e. the isomorphism defined in Lemma \ref{LM:K(N)=K_N(G)}. Take $b\in L$, i.e. $L=N\cdot b$, and let $\tau_b:{\mathcal K}(G)\to{\mathcal K}(G)$ be the shift by element $b^{-1}$ (that acts on ${\mathcal K}(G)$ by Theorem \ref{LM:sdvig-v-K(G)}):
$$
\tau_bu=b^{-1}\cdot u,\qquad {\mathcal K}(G).
$$
It turns the space ${\mathcal K}(N)$ into the space ${\mathcal K}_L(G)$, hence a mapping is defined $\sigma_L=\tau_b\circ\sigma:{\mathcal K}(N)\to {\mathcal K}_L(G)$. Put now
$$
\chi_L=\chi\circ\rho_L,\qquad \chi_N=\chi_L\circ\sigma_L
$$
we obtain a commutative diagram
$$
\xymatrix @R=2.pc @C=3.0pc % @M=14pt
{
{\mathcal K}(N)\ar[r]^{\sigma_L}\ar[dr]_{\chi_N} & {\mathcal K}_L(G)\ar[d]^{\chi_L}\ar[r]^{\rho_L} & {\mathcal K}(G)\ar[dl]^{\chi}\\
 & \C &
}
$$
Since $\chi_N$ is a character on ${\mathcal K}(N)$, by Lemma \ref{LM:Spec(R^n-times-K)} it must be a delta-functional:
\beq\label{chi_N(u)=u(a)}
\chi_N(u)=u(a),\qquad u\in{\mathcal K}(N)
\eeq
for some $a\in N$. Then
\begin{multline*}
\chi(u)=\chi(1_L)\cdot\chi(u)=\chi(1_L\cdot u)=\chi_L(1_L\cdot u)=\chi_N(\sigma_L^{-1}(1_L\cdot u))=\eqref{chi_N(u)=u(a)}=
\sigma_L^{-1}(1_L\cdot u)(a)=\sigma(\sigma_L^{-1}(1_L\cdot u))(a)=\\=
(\sigma\circ\sigma_L^{-1})(1_L\cdot u))(a)=(\sigma\circ(\tau_b\circ\sigma)^{-1})(1_L\cdot u))(a)=
(\sigma\circ\sigma^{-1}\circ\tau_b^{-1})(1_L\cdot u))(a)=(\sigma\circ\sigma^{-1}\circ\tau_b^{-1})(1_L\cdot u))(a)=\\=
\tau_{b^{-1}}(1_L\cdot u))(a)=(b\cdot(1_L\cdot u))(a)=(1_L\cdot u)(\underbrace{a\cdot b}_{\scriptsize\begin{matrix}\text{\rotatebox{90}{$\owns$}}\\ L\end{matrix}})=u(a\cdot b)=\delta^{a\cdot b}(u).
\end{multline*}

2. Now let us verify that the mapping of spectra $G\to \Spec{\mathcal K}(G)$ is an injection. Take $a\ne b\in G$. If  $a\cdot b^{-1}\notin N$, i.e. $a\notin b\cdot N$, then the characteristic function $1_L\in{\mathcal K}(G)$ of the class $L=b\cdot N$ from Lemma \ref{LM:1_L-in-K(G)} separates $a$ and $b$:
$$
1_L(a)=0\ne 1=1_L(b).
$$
Another possibility: suppose $a\in b\cdot N$, i.e. $a\cdot b^{-1}\in N$. Then by Lemma \ref{LM:Spec(R^n-times-K)} we can take a function $u\in{\mathcal K}(N)$ such that
$$
u(a\cdot b^{-1})\ne u(1)
$$
By Lemma \ref{LM:K(N)=K_N(G)} there is a function $v\in{\mathcal K}_N(G)$ such that $u\big|_N=v\big|_N$, hence
$$
v(a\cdot b^{-1})\ne v(1)
$$
By Theorem \ref{LM:sdvig-v-K(G)} the shift $b^{-1}\cdot v$ again belongs to ${\mathcal K}(G)$, hence for this function we have
$$
(b^{-1}\cdot v)(a)=v(a\cdot b^{-1})\ne v(1)=v(b\cdot b^{-1})=(b^{-1}\cdot v)(b).
$$

3. It remains to verify the openness of the mapping $G\to \Spec{\mathcal K}(G)$. Suppose $a_i\to a$ in $\Spec{\mathcal K}(G)$. Theorem \ref{LM:sdvig-v-K(G)} implies immediately that $a_i\cdot a^{-1}\to 1$ in $\Spec{\mathcal K}(G)$. For the characteristic function $1_N\in{\mathcal K}(G)$ of the subgroup $N$ we have $1_N(a_i\cdot a^{-1})\to 1_N(1)=1$, hence, starting from some index, all $a_i\cdot a^{-1}$ belong to $N$. Take a compact set $S\subseteq {\mathcal K}(N)$. By Lemma \ref{LM:K(N)=K_N(G)} we can find a compact set $T\subseteq {\mathcal K}(G)$, which consists of functions whose restrictions to $N$ belong to $K$, and a bijection between $T$ and $S$ appears. Since $a_i\cdot a^{-1}\to 1$ in $\Spec{\mathcal K}(G)$, we have
$$
v(a_i\cdot a^{-1})\underset{i\to\infty}{\underset{v\in T}{\rightrightarrows}}v(1)
$$
and this is equivalent to
$$
u(a_i\cdot a^{-1})\underset{i\to\infty}{\underset{u\in S}{\rightrightarrows}}u(1).
$$
This is true for any compact set $S\subseteq {\mathcal K}(N)$, hence $a_i\cdot a^{-1}\to 1$ in $\Spec{\mathcal K}(N)$. But we already proved in Lemma \ref{LM:Spec(R^n-times-K)} that $\Spec{\mathcal K}(N)=N$, therefore we obtain that $a_i\cdot a^{-1}\to 1$ in $N$, and thus, in $G$.
\epr

Theorems \ref{TH:Spec-K(G)=G} and \ref{C-obolochka-podalgebry-v-C(M)} imply

\btm\footnote{See Errata on page \pageref{Errata}.}\label{TH:E(K(G))=C(G)}
If $G$ is a Moore group, then the continuous envelope of the algebra ${\mathcal K}(G)$ is the algebra ${\mathcal C}(G)$:
\beq\label{E(K(G))=C(G)}
\Env_{\mathcal C}{\mathcal K}(G)={\mathcal C}(G)
\eeq
\etm

\paragraph{The structure of Hopf algebras on $\Env_{\mathcal C} {\mathcal C}^\star(G)$ an on ${\mathcal K}(G)$ in the case of Moore groups.}

For Moore groups it is possible to prove that the algebras $\Env_{\mathcal C}{\mathcal C}^\star(G)$ and ${\mathcal K}(G)$ are involutive Hopf algebras.

\btm\footnote{See Errata on page \pageref{Errata}.}\label{TH:E(H)-dlya-Moore} Suppose $G$ is a Moore group. Then
\bit{
\item[(i)] the continuous envelope $\Env_{\mathcal C}{\mathcal C}^\star(G)$ of its group algebra ${\mathcal C}^\star(G)$ is an involutive Hopf algebra in the category of stereotype spaces $(\tt{Ste},\odot)$.

\item[(ii)] the dual algebra ${\mathcal K}(G)$ is an involutive Hopf algebra in the category of stereotype spaces $(\tt{Ste},\circledast)$
}\eit
\etm

We premise the proof with tree lemmas.

\blm\label{LM:E(H)-dlya-Moore} Suppose $G$ is a Moore group. Then for each seminorm $p\in{\mathcal P}({\mathcal C}^\star(G))$ the quotient algebra ${\mathcal C}^\star(G)/p$ is a strict $C^*$-algebra.
\elm
\bpr
If $\sigma:{\mathcal C}^\star(G)/p\to{\mathcal B}(X)$ is a unitary irreducible representation, then its composition with the quotient mapping  $\pi_p:{\mathcal C}^\star(G)\to {\mathcal C}^\star(G)/p$ is a unitary irreducible representation of the algebra ${\mathcal C}^\star(G)$. The group $G$, being embedded into ${\mathcal C}^\star(G)$ by delta-functionals, is dense in ${\mathcal C}^\star(G)$. This means that the composition
$$
\xymatrix @R=2.pc @C=4.0pc % @M=14pt
{
G\ar[r]^{\delta}& {\mathcal C}^\star(G)\ar[r]^{\pi_p}& {\mathcal C}^\star(G)/p\ar[r]^{\sigma}& {\mathcal B}(X)
}
$$
is a unitary irreducible representation of the Moore group $G$, and therefore it must be finite dimensional: $\dim X<\infty$. The boundedness of these numbers (with freezed $p$ and varying $\tau$) are proved in several steps (we use here the ideas of \cite[Lemma 5.8]{Kuznetsova}).

1. First we consider the case when $G$ is a compact group. Then by Proposition \ref{PROP:nepr-obol-komp-gruppy} the continuous envelope of the algebra of mesures ${\mathcal C}^\star(G)$ is the algebra
$\prod_{\sigma\in\widehat{G}}{\mathcal B}(X_\sigma)$. Since the quotient algebra ${\mathcal C}^\star(G)/p$ is a  $C^*$-algebra, the quotient mapping $\pi_p$ can be extended to the continuous envelope $\prod_{\sigma\in\widehat{G}}{\mathcal B}(X_\sigma)$:
$$
\xymatrix @R=2.pc @C=1.7pc % @M=14pt
{
{\mathcal C}^\star(G)\ar@/_3ex/[rd]_{\pi_p}\ar[rr]^{\env_{\mathcal C}{\mathcal C}^\star(G)}& & \prod_{\sigma\in\widehat{G}}{\mathcal B}(X_\sigma)\ar@/^3ex/@{-->}[ld]^{\widetilde{\pi_p}} \\
& {\mathcal C}^\star(G)/p &
}
$$
Note that $\widetilde{\pi_p}$ is an epimorphism of locally convex spaces (since the composition $\widetilde{\pi_p}\circ\env_{\mathcal C}{\mathcal C}^\star(G)=\pi_p$ is an epimorphism of locally convex spaces). On the other hand, since all algebras ${\mathcal B}(X_\sigma)$ are finite dimensional, the algebra $\prod_{\sigma\in\widehat{G}}{\mathcal B}(X_\sigma)$, as a locally convex space, is a Cartesian power $\C^{\mathfrak m}$ of the field $\C$ (where ${\mathfrak m}$ is a cardinal number).

This implies that the mapping $\widetilde{\pi_p}$ has the kernel of finite co-dimension (since it acts to the Banach space ${\mathcal C}^\star(G)/p$). This means in its turn, that in the family $\{{\mathcal B}(X_\sigma);\ \sigma\in\widehat{G}\}$ there irs a finite subfamily of algebras ${\mathcal B}(X_1),...,{\mathcal B}(X_n)$, such that $\widetilde{\pi_p}$ is a projection on its product:
\beq\label{C^star(G)/p-cong-prod_1^n-B(X_i)}
\widetilde{\pi_p}: \prod_{\sigma\in\widehat{G}}{\mathcal B}(X_\sigma)\to \prod_{i=1}^n{\mathcal B}(X_i)\cong {\mathcal C}^\star(G)/p
\eeq
Certainly, ${\mathcal C}^\star(G)/p$ in this situation is a strict $C^*$-algebra.

2. Let then $G$ be a compact buildup of an Abelian group, i.e. $G=Z\cdot K$, where $Z$ is Abelian, $K$ is compact, and they commute (see definition on page \pageref{DEF:komp-nadstr-abelevoi-gruppy}).
Suppose $\pi:G\to{\mathcal B}(X)$ is a norm continuous unitary irreducible representation. Consider its restrictions $\rho=\pi\big|_K$ and $\sigma=\pi\big|_Z$ and denote by $C_\pi$, $C_\rho$ and $C_\sigma$ respectively the $C^*$-subalgebras in ${\mathcal B}(X)$, generated by the images of $\pi$, $\rho$ and $\sigma$. Since $C_\rho$ and $C_\sigma$ commute, $C_\pi$ is a continuous image of the maximal tensor product $C_\rho\underset{\max}{\otimes}C_\sigma$. Note that the theorems \ref{TH:nepr-predstavleniya} and \ref{TH:nepr-po-norme-predstavleniya} imply that $C_\rho$ and $C_\sigma$ are $C^*$-quotient algebras of ${\mathcal C}^\star(K)$ and ${\mathcal C}^\star(Z)$.
From \eqref{C^star(G)/p-cong-prod_1^n-B(X_i)} we have that $C_\rho=\prod_{i=1}^n{\mathcal B}(X_{\pi_i})$, where  $\pi_i\in\widehat{K}$ is a finite sequence of unitary irreducible representations of $K$. At the same time,  Proposition \ref{PROP:Env_C-C^star(G)=C(widehat(G))} imply that $C_\sigma={\mathcal C}(T)$ for some compact $T\subseteq\widehat{Z}$. Hence
$$
C_\rho\underset{\max}{\otimes}C_\sigma=\left(\prod_{i=1}^n{\mathcal B}(X_{\pi_i})\right)\underset{\max}{\otimes}{\mathcal C}(T)=\prod_{i=1}^n\left({\mathcal B}(X_{\pi_i})\underset{\max}{\otimes}{\mathcal C}(T)\right)=\prod_{i=1}^n{\mathcal C}\left(T,{\mathcal B}(X_{\pi_i})\right)={\mathcal C}\left(\bigsqcup_{i=1}^n T_i,{\mathcal B}(X_{\pi_i})\right),
$$
where $T_i$ are copies of the compact set $T$. Now by \cite[10.4.4]{Dixmier} each unitary irreducible representation of the last algebra is isomorphic to some $X_{\pi_i}$, hence it has the dimension not greater that $\max_{i=1,...,n}X_{\pi_i}$. The same is true for the uintary irreducible representations of the algebra $C_\pi$, since being moved to $C_\rho\underset{\max}{\otimes}C_\sigma$ they also become unitary irreducible representations.

3. Further, suppose $G$ is a Lie-Moore group. By Theorem \ref{TH:Lie-Moore} $G$ is a finite extension of a compact buildup $H=Z\times K$ of an Abelian group $Z$. Suppose $m=\card G/H$ is the index of $H$ in $G$. Take $p\in{\mathcal P}({\mathcal C}^\star(G))$ and a unitary irreducible representation $\tau$ of the $C^*$-algebra ${\mathcal C}^\star(G)/p$. By \cite[Theorem 1]{Clifford}, the restriction of $\tau$ to $H$ is decomposed into a sum of no more than $m$ unitary irreducible representations of $H$, hence of the algebra ${\mathcal C}^\star(H)/p$. But above we have already proved that the dimension of unitary irreducible representations of the algebra ${\mathcal C}^\star(H)/p$ is bounded by some number $n\in\N$. Thus the dimension of $\tau$ is not greater than $m\cdot n$.

4. Finally, suppose $G$ is an arbitrary Moore group and $p\in{\mathcal P}({\mathcal C}^\star(G))$. The homomorphism  $G\to {\mathcal C}^\star(G)/p$ is norm-continuous, hence by \cite[Theorem 1]{Shtern}, it is factored through some quotient mapping $G\to G/H$, where $G/H$ is a Lie group. By Theorem \ref{TH:Moore->quotient} $G/H$ is a Moore group. So we reduced the situation to the previous case.
\epr

\blm\footnote{See Errata on page \pageref{Errata}.}\label{LM:E(H)-dlya-Moore-ii} Suppose $G$ is a Moore group. Then the diagonal $\beta$ of the diagram
$$
\xymatrix @R=6.pc @C=8.pc % @M=14pt
{
{\mathcal C}^\star(G)\circledast {\mathcal C}^\star(G)\ar@{-->}[dr]^{\beta}
\ar[d]_{@}\ar[r]^(.4){\env_{\mathcal C}{{\mathcal C}^\star(G)}\circledast \env_{\mathcal C}{{\mathcal C}^\star(G)}} & \Env_{\mathcal C}{\mathcal C}^\star(G)\circledast \Env_{\mathcal C}{\mathcal C}^\star(G)\ar[d]_{@} \\
{\mathcal C}^\star(G)\odot {\mathcal C}^\star(G)\ar[r]^(.4){\env_{\mathcal C}{{\mathcal C}^\star(G)}\odot \env_{\mathcal C}{{\mathcal C}^\star(G)}} & \Env_{\mathcal C}{\mathcal C}^\star(G)\odot \Env_{\mathcal C}{\mathcal C}^\star(G)
}
$$
is a dense epimorphism (i.e. the images of elements from the domain are dense in the range).
\elm
\bpr Consider the representation \eqref{SIN-kak-rasshirenie} of the SIN-group $G$, and let us construct the following chain of morphisms:
$$
\xymatrix @R=2.pc @C=.5pc % @M=14pt
{
{\mathcal C}^\star(G)\circledast {\mathcal C}^\star(G)\ar@{=}[r]&
 \prod\limits_{\sigma\in\widehat{K}} {\mathcal C}_\sigma^\star(G)\circledast \prod\limits_{\tau\in\widehat{K}} {\mathcal C}_\tau^\star(G)\ar[d]_{\pi_\sigma\circledast\pi_\tau} & & \\
& {\mathcal C}_\sigma^\star(G)\circledast {\mathcal C}_\tau^\star(G)\ar@{=}[r]&
\leftlim_{\infty\gets k}{\mathcal C}^\star(G)/p_{T_k,\sigma}^{\max}\circledast
 \leftlim_{\infty\gets l}{\mathcal C}^\star(G)/p_{T_l,\tau}^{\max}
\ar[d]_(.4){\pi_{T_k,\sigma}\circledast\pi_{T_l,\tau}} & \\
   & &
   {\mathcal C}^\star(G)/p_{T_k,\sigma}^{\max}\circledast {\mathcal C}^\star(G)/p_{T_l,\tau}^{\max}\ar[d] & \\
   && {\mathcal C}^\star(G)/p_{T_k,\sigma}^{\max}\check{\otimes}\ {\mathcal C}^\star(G)/p_{T_l,\tau}^{\max}\ar@{=}[r]
   &   {\mathcal C}^\star(G)/p_{T_k,\sigma}^{\max}\odot {\mathcal C}^\star(G)/p_{T_l,\tau}^{\max}
}
$$
(we use the notations from Proposition \ref{PROP:nepr-Env-SIN-gruppy}, and the morphisms $\pi_{\sigma}$ and $\pi_{T_k,\sigma}$ are natural projections, the last arrow is the natural morphism of tensor products, and the last equality follows from Lemma \ref{LM:E(H)-dlya-Moore} and equalities \eqref{A-min-B=A-odot-B}). Each arrow here is a dense epimorphism, hence the composition
$$
\xymatrix @R=2.pc @C=5pc % @M=14pt
{
{\mathcal C}^\star(G)\circledast {\mathcal C}^\star(G)\ar[r]&
{\mathcal C}^\star(G)/p_{T_k,\sigma}^{\max}\odot {\mathcal C}^\star(G)/p_{T_l,\tau}^{\max}
}
$$
is also a dense epimorphism. When $k$ and $l$ tend to infinity, they give a natural morphism
$$
\xymatrix @R=2.pc @C=5pc % @M=14pt
{
{\mathcal C}^\star(G)\circledast {\mathcal C}^\star(G)\ar[r]&
\leftlim_{\infty\gets k,l}{\mathcal C}^\star(G)/p_{T_k,\sigma}^{\max}\odot {\mathcal C}^\star(G)/p_{T_l,\tau}^{\max}
}
$$
which is also a dense epimorphism (we use here the fact that a locally convex projective limit of a sequence of Banach spaces is a Fr\'echet space, and thereofre it coincides with the stereotype projective limit of this system).

After that we again use commutativity of the projective limit with the injective tensor product \cite[(2.53)]{Akbarov-env}, and we see that the morphism
$$
\xymatrix @R=2.pc @C=1pc % @M=14pt
{
{\mathcal C}^\star(G)\circledast {\mathcal C}^\star(G)\ar[r]&
\leftlim_{\infty\gets k,l}{\mathcal C}^\star(G)/p_{T_k,\sigma}^{\max}\odot {\mathcal C}^\star(G)/p_{T_l,\tau}^{\max}
=
\leftlim_{\infty\gets k}{\mathcal C}^\star(G)/p_{T_k,\sigma}^{\max}\odot \leftlim_{\infty\gets l}{\mathcal C}^\star(G)/p_{T_l,\tau}^{\max}= {\mathcal C}_\sigma^\star(G)\odot {\mathcal C}_\tau^\star(G)
}
$$
is also a dense epimorphism.

This is true for all $\sigma,\tau\in\widehat{K}$, hence the morphism into the direct product
$$
\xymatrix @R=2.pc @C=1pc % @M=14pt
{
{\mathcal C}^\star(G)\circledast {\mathcal C}^\star(G)\ar[r]&
\prod\limits_{\sigma,\tau\in\widehat{K}} {\mathcal C}_\sigma^\star(G)\odot {\mathcal C}_\tau^\star(G)=
\prod\limits_{\sigma\in\widehat{K}} {\mathcal C}_\sigma^\star(G)\odot
\prod\limits_{\tau\in\widehat{K}}{\mathcal C}_\tau^\star(G)=
\Env_{\mathcal C}{\mathcal C}^\star(G)\odot \Env_{\mathcal C}{\mathcal C}^\star(G)
}
$$
is again dense (we use here the fact that the locally convex direct product of stereotype spaces coincides with their stereotype direct product). This is the morphism that we need.
\epr

\bpr[Proof of Theorem \ref{TH:E(H)-dlya-Moore}] Clearly, propositions (i) and (ii) are equivalent, so it is sifficient to prove (i).

1. First, ${\mathcal C}^\star(G)$ is a stereotype algebra (i.e. an algebra in the category $(\tt{Ste},\circledast)$), and its envelope $\Env_{\mathcal C}{\mathcal C}^\star(G)$ is also a stereotype algebra (i.e. an algebra in  $(\tt{Ste},\circledast)$). Let us show that the multiplication in $\Env_{\mathcal C}{\mathcal C}^\star(G)$ can be extended to an operator on $\Env_{\mathcal C}{\mathcal C}^\star(G)\odot \Env_{\mathcal C}{\mathcal C}^\star(G)$.

By Lemma \ref{LM:E(H)-dlya-Moore}, ${\mathcal C}^\star(G)/p$ is a strict $C^*$-algebra. Hence, the multiplication $\mu_p$ is extended to an operator $\mu'_p:{\mathcal C}^\star(G)/p\odot {\mathcal C}^\star(G)/p\to {\mathcal C}^\star(G)/p$. Consider the diagram
$$
\xymatrix @R=2.pc @C=4.0pc % @M=14pt
{
{\mathcal C}^\star(G)\circledast{\mathcal C}^\star(G)\ar[d]^{\mu}\ar[r]^{\pi_p\circledast\pi_p}& {\mathcal C}^\star(G)/p\circledast{\mathcal C}^\star(G)/p\ar[d]^{\mu_p}\ar[r]^{@}& {\mathcal C}^\star(G)/p\odot{\mathcal C}^\star(G)/p\ar[d]^{\mu'_p} \\
{\mathcal C}^\star(G)\ar[r]^{\pi_p}& {\mathcal C}^\star(G)/p\ar[r]^{\id}& {\mathcal C}^\star(G)/p
}
$$
If we remove the middle column and pass to the projective limit in $\tt{Ste}$, we obtain:
$$
\xymatrix @R=2.pc @C=2.0pc % @M=14pt
{
{\mathcal C}^\star(G)\circledast{\mathcal C}^\star(G)\ar[d]^{\mu}\ar@{-->}[r]& \kern-13pt\leftlim_{p\in{\mathcal P({\mathcal C}^\star(G))}}\kern-13pt\Big({\mathcal C}^\star(G)/p\odot{\mathcal C}^\star(G)/p\Big)\ar@{-->}[d]\ar@{=}[r]&
\kern-13pt
\leftlim_{p\in{\mathcal P({\mathcal C}^\star(G))}}\kern-13pt{\mathcal C}^\star(G)/p\odot\kern-13pt \leftlim_{p\in{\mathcal P({\mathcal C}^\star(G))}}\kern-13pt{\mathcal C}^\star(G)/p\ar@{=}[r] &
\Env_{\mathcal C}{\mathcal C}^\star(G)\odot \Env_{\mathcal C}{\mathcal C}^\star(G)\ar@{-->}[d]^{\mu'}\\
{\mathcal C}^\star(G)\ar@{-->}[r]& \kern-13pt\leftlim_{p\in{\mathcal P({\mathcal C}^\star(G))}}\kern-13pt{\mathcal C}^\star(G)/p\ar@{=}[rr]& & \Env_{\mathcal C}{\mathcal C}^\star(G)
}
$$
Here the first equality in the upper line is the result of the commutativity of the projective limit with the injective stereotype tensor product \cite[(2.53)]{Akbarov-env}, and the second equality (together with the equality in the bottom line) is the formula \eqref{E(C*(G))=LCS-leftlim}. The operator $\mu'$ is the extension of the multiplication to the injective tensor square, what we need.

2. Further, ${\mathcal C}^\star(G)$ is a bialgebra in the category $(\tt{Ste},\circledast)$, hence by Theorem  \ref{TH:C-obolochka-sohranyaet-Hopfov} $\Env_{\mathcal C}{\mathcal C}^\star(G)$ is a coalgebra in the category
$(\tt{Ste},\odot)$. Thus, $\Env_{\mathcal C}{\mathcal C}^\star(G)$ is an algebra and a coalgebra in  $(\tt{Ste},\odot)$. Consider the diagram

{\scriptsize
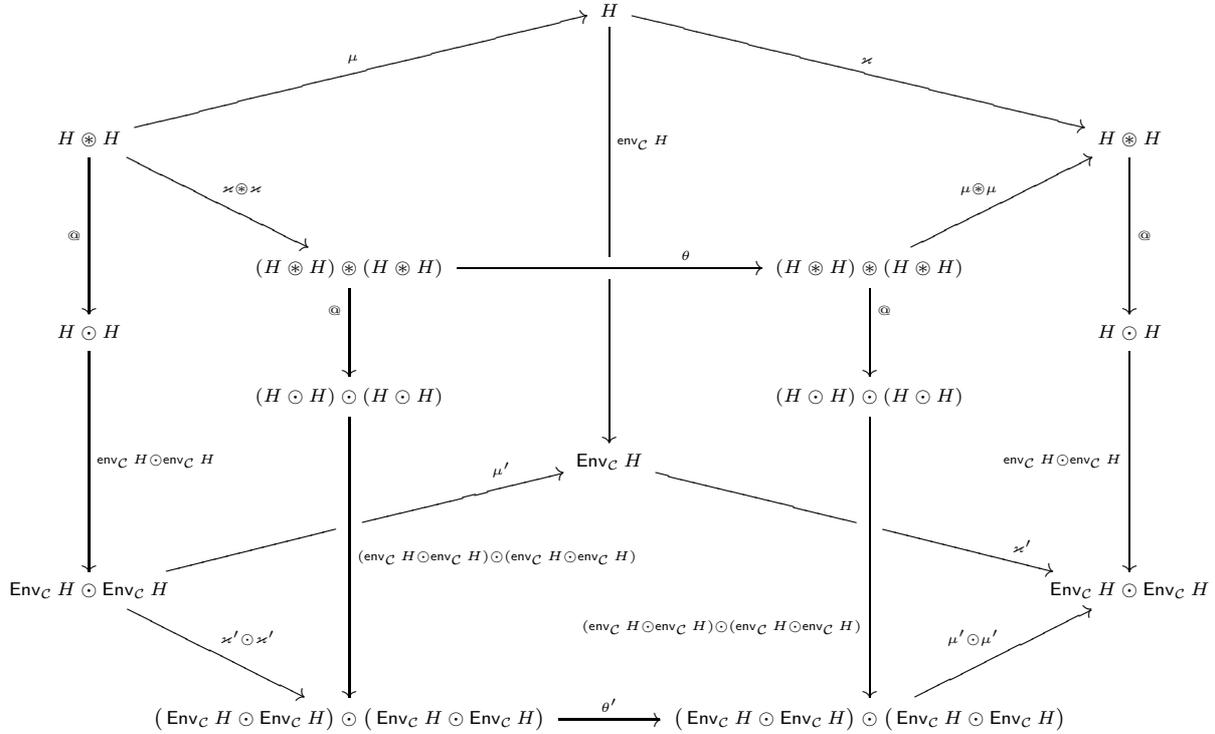
 $$\dgARROWLENGTH=-6em
\begin{diagram}
  \node[3]{H}\arrow[2]{se,l}{\varkappa}\arrow[4]{s,r,-}{\env_{\mathcal C}H}
 \\
 \\
 \node{H\circledast H}\arrow[3]{s,t}{@}
 \arrow{sse,l}{\varkappa\circledast\varkappa}
 \arrow[2]{ne,l}{\mu}
 \node[4]{H\circledast H}\arrow[3]{s,r}{@} \\
  \\
 \node[2]{(H\circledast H)\circledast (H\circledast H)}
 \arrow[2]{e,l,3}{\theta}
 \arrow[2]{s,l,1}{@}
 \node{}\arrow[3]{s}
 \node{(H\circledast H)\circledast (H\circledast H)}
 \arrow{nne,l}{\mu\circledast\mu}
 \arrow[2]{s,r,1}{@}
 \\
  \node{H\odot H}\arrow[4]{s,b}{\env_{\mathcal C}H\odot \env_{\mathcal C}H}
 \node[4]{H\odot H}\arrow[4]{s,l}{\env_{\mathcal C}H\odot \env_{\mathcal C}H}
 \\
 \node[2]{(H\odot H)\odot (H\odot H)}
 \arrow[5]{s,r}{(\env_{\mathcal C}H\odot \env_{\mathcal C}H)\odot(\env_{\mathcal C}H\odot \env_{\mathcal C}H)}
 \node[2]{(H\odot H)\odot (H\odot H)}
 \arrow[5]{s,l,3}{(\env_{\mathcal C}H\odot \env_{\mathcal C}H)\odot(\env_{\mathcal C}H\odot \env_{\mathcal C}H)}
 \\
  \node[3]{\Env_{\mathcal C}H}\arrow{se,-}
 \\
 \node[2]{}\arrow{ne,l,3}{\mu'}
 \node[2]{}\arrow{se,t,3}{\varkappa'}
 \\
 \node{\Env_{\mathcal C}H\odot \Env_{\mathcal C}H}
 \arrow{sse,l}{\varkappa'\odot\varkappa'}
 \arrow{ne,-}
 \node[4]{\Env_{\mathcal C}H\odot \Env_{\mathcal C}H}
  \\ \\
 \node[2]{\big(\Env_{\mathcal C}H\odot \Env_{\mathcal C}H\big)\odot \big(\Env_{\mathcal C}H \odot \Env_{\mathcal C}H\big)}
 \arrow[2]{e,l}{\theta'}
 \node[2]{\big(\Env_{\mathcal C}H\odot \Env_{\mathcal C}H\big)\odot \big(\Env_{\mathcal C}H \odot \Env_{\mathcal C}H\big)}
 \arrow{nne,l}{\mu'\odot\mu'}
\end{diagram}
 $$
 }\noindent
where $H={\mathcal C}^\star(G)$, and the sense of the other notations is obvious. Here the upper base is commutative, since $H$ is a Hopf algebra in $(\tt{Ste},\circledast)$, and the commutativity of the lateral faces is verified by the direct computations. In addition, by Lemma \ref{LM:E(H)-dlya-Moore-ii}, the very left morphism $\env_{\mathcal C}H\odot \env_{\mathcal C}H\circ @$ is an epimorphism in $\tt{Ste}$. This means that the lower base is also commutative.

In the same manner we prove the commutativity of the other diagrams in the definition of the Hopf algebra. For example, the diagram, that binds the comultiplication with the unit, is verified by building the following prism:
 $$\dgARROWLENGTH=1em
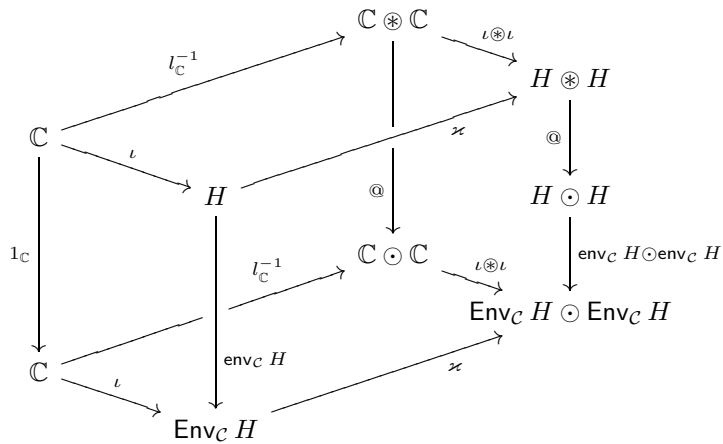
\begin{diagram}
  \node[3]{\C\circledast\C}\arrow{se,l}{\iota\circledast\iota}\arrow[2]{s,-}
 \\
   \node[4]{H\circledast H}\arrow[2]{s,l}{@}
 \\
 \node{\C}\arrow[2]{ne,l}{l_{\C}^{-1}}\arrow{se,l}{\iota}\arrow[4]{s,l}{1_{\C}}
 \node[2]{}\arrow[2]{s,l}{@}
 \\
 \node[2]{H}\arrow[2]{ne,r,3}{\varkappa}\arrow[4]{s,r,3}{\env_{\mathcal C}H}\node[2]{H\odot H}\arrow[2]{s,r}{\env_{\mathcal C}H\odot \env_{\mathcal C}H}
 \\
  \node[3]{\C\odot\C}\arrow{se,l}{\iota\circledast\iota}
 \\
  \node[2]{}\arrow{ne,l}{l_{\C}^{-1}} \node[2]{\Env_{\mathcal C}H\odot \Env_{\mathcal C}H}
 \\
 \node{\C}\arrow{ne,-}\arrow{se,l}{\iota}
 \\
 \node[2]{\Env_{\mathcal C}H}\arrow[2]{ne,r,3}{\varkappa}
 \\
\end{diagram}
 $$
Here the upper base is commutative, since $H$ is a Hopf algebra in $(\tt{Ste},\circledast)$, and the commutativity of the lateral faces follows from the properties of the functor $\Env_{\mathcal C}$. Therefore the lower base is commutative as well.

3. Now let us prove that the involution $\bullet$ in the Hopf algebra ${\mathcal C}^\star(G)$ generates the involution $\bullet'$ in the Hopf algebra $\Env_{\mathcal C}{\mathcal C}^\star(G)$. For each seminorm $p\in{\mathcal P}({\mathcal C}^\star(G))$ consider the natural projection $\pi_p:\Env_{\mathcal C}{\mathcal C}^\star(G)\to {\mathcal C}^\star(G)/p$. Let $\bullet/p$ be the involution on ${\mathcal C}^\star(G)/p$ generated by the involution $\bullet$ in ${\mathcal C}^\star(G)$. Put $\bullet_p=\bullet/p\ \circ\pi_p$:
$$
\xymatrix @R=2.pc @C=4.0pc % @M=14pt
{
\Env_{\mathcal C}{\mathcal C}^\star(G)\ar[r]^{\pi_p}\ar@/_2ex/[rd]_{\bullet_p} & {\mathcal C}^\star(G)/p\ar[d]^{\bullet/p} \\
& {\mathcal C}^\star(G)/p
}
$$
For each seminorms $p,q\in{\mathcal P}({\mathcal C}^\star(G))$, $p\le q$, in the category of stereotype spaces over the real field $\R$ the following diagram is commutative:
$$
\xymatrix @R=2.pc @C=2.0pc % @M=14pt
{
& \Env_{\mathcal C}{\mathcal C}^\star(G)\ar@/_2ex/[dl]_{\bullet_q}\ar@/^2ex/[dr]^{\bullet_p} &  \\
{\mathcal C}^\star(G)/q\ar[rr]_{\pi^q_p} & & {\mathcal C}^\star(G)/p
}
$$
where $\pi^q_p$ is the natural projection. This means that the family of morphisms $\bullet_p$ is a projective cone for the system $\pi^q_p$, hence there exists a unique morphism $\bullet'$ such that the following diagrams are commutative:
$$
\xymatrix @R=2.pc @C=2.0pc % @M=14pt
{
& & \Env_{\mathcal C}{\mathcal C}^\star(G)\ar@/_2ex/[dl]_{\bullet'}\ar@/^2ex/[dr]^{\bullet_p} &  \\
\Env_{\mathcal C}{\mathcal C}^\star(G)\ar@{=}[r]& \leftlim_{q}{\mathcal C}^\star(G)/q\ar[rr]_{\pi_p} & & {\mathcal C}^\star(G)/p
}
$$
Let us show that the morphism $\bullet'$ is the involution on $\Env_{\mathcal C}{\mathcal C}^\star(G)$ that we look for. First, it is connected with the initial involution $\bullet$ by the commutative diagram:
$$
\xymatrix @R=3.pc @C=6.0pc % @M=14pt
{
{\mathcal C}^\star(G)\ar[r]^{\env_{\mathcal C}{{\mathcal C}^\star(G)}}\ar[d]_{\bullet} & \Env_{\mathcal C}{\mathcal C}^\star(G)\ar[d]^{\bullet'} \\
{\mathcal C}^\star(G)\ar[r]^{\env_{\mathcal C}{{\mathcal C}^\star(G)}} & \Env_{\mathcal C}{\mathcal C}^\star(G) \\
}
$$
It implies the equality
$$
\bullet'\circ \bullet'\circ \env_{\mathcal C}{{\mathcal C}^\star(G)}=\env_{\mathcal C}{{\mathcal C}^\star(G)}\circ \bullet\circ \bullet=\env_{\mathcal C}{{\mathcal C}^\star(G)}=
\id_{\Env_{\mathcal C}{\mathcal C}^\star(G)}\circ \env_{\mathcal C}{{\mathcal C}^\star(G)}
$$
which, by the epimorphy of $\env_{\mathcal C}{{\mathcal C}^\star(G)}$, gives
$$
\bullet'\circ \bullet'=\id_{\Env_{\mathcal C}{\mathcal C}^\star(G)}.
$$

To prove \eqref{*-circ-mu=mu-circ-*-circledast-*-circ-op}, consider the diagram (in the category of stereotype spaces over $\R$), similar to the one from the paragraph 2 (here $H={\mathcal C}^\star(G)$, and $\br:x\otimes y\mapsto y\otimes x$ is the braiding morphism in the monoidal category):
 $$\dgARROWLENGTH=-2em
\begin{diagram}
  \node[3]{H}\arrow[2]{se,l}{\bullet}\arrow[4]{s,r,-}{\env_{\mathcal C}H}
 \\
 \\
 \node{H\circledast H}\arrow[3]{s,t}{@}
 \arrow{sse,l}{\br}
 \arrow[2]{ne,l}{\mu}
 \node[4]{H}\arrow[7]{s,r}{\env_{\mathcal C}H} \\
  \\
 \node[2]{H\circledast H}
 \arrow[2]{e,l,3}{\bullet\circledast\bullet}
 \arrow[2]{s,l,1}{@}
 \node{}\arrow[3]{s}
 \node{H\circledast H}
 \arrow{nne,l}{\mu}
 \arrow[2]{s,r,1}{@}
 \\
  \node{H\odot H}\arrow[4]{s,t}{\env_{\mathcal C}H\odot \env_{\mathcal C}H}
 \\
 \node[2]{H\odot H}
 \arrow[5]{s,r,3}{\env_{\mathcal C}H\odot \env_{\mathcal C}H}
 \node[2]{H\odot H}
 \arrow[5]{s,l,3}{\env_{\mathcal C}H\odot \env_{\mathcal C}H}
 \\
  \node[3]{\Env_{\mathcal C}H}\arrow{se,-}
 \\
 \node[2]{}\arrow{ne,r,3}{\mu'}
 \node[2]{}\arrow{se,l,3}{\bullet'}
 \\
 \node{\Env_{\mathcal C}H\odot \Env_{\mathcal C}H}
 \arrow{sse,r}{\br}
 \arrow{ne,-}
 \node[4]{\Env_{\mathcal C}H}
  \\ \\
 \node[2]{\Env_{\mathcal C}H\odot \Env_{\mathcal C}H}
 \arrow[2]{e,l}{\bullet'\odot\bullet'}
 \node[2]{\Env_{\mathcal C}H\odot \Env_{\mathcal C}H}
 \arrow{nne,r}{\mu'}
\end{diagram}
 $$
In this prism the upper base is commutative (as \eqref{*-circ-mu=mu-circ-*-circledast-*-circ-op} for $H$), and the commutativity of the lateral faces is checked by computation. Besides this, by Lemma \ref{LM:E(H)-dlya-Moore-ii} the left edge $\env_{\mathcal C}H\odot \env_{\mathcal C}H\circ @$ is an epimorphism. This implies that the lower base is also commutative:
$$
\bullet'\circ \mu=\mu\circ \bullet'\odot\bullet'\circ\br.
$$

To prove \eqref{varkappa-circ-*=*-circledast-*-circ-varkappa} let us consider the diagram (in the category of stereotype spaces over $\R$):
 $$\dgARROWLENGTH=1em
\begin{diagram}
  \node[3]{H}\arrow{se,l}{\varkappa}\arrow[2]{s,-}
 \\
   \node[4]{H\circledast H}\arrow{s,r}{@}
 \\
 \node{H}\arrow[2]{ne,l}{\bullet}\arrow{se,l}{\varkappa}\arrow[4]{s,l}{\env_{\mathcal C}H}
 \node[2]{}\arrow[2]{s,l}{\env_{\mathcal C}H}\node{H\odot H}\arrow[3]{s,r,3}{\env_{\mathcal C}H\odot \env_{\mathcal C}H}
 \\
 \node[2]{H\circledast H}\arrow[2]{ne,l,1}{\bullet\circledast\bullet}\arrow{s,r}{@}
 \\
 \node[2]{H\odot H}\arrow[3]{s,r,3}{\env_{\mathcal C}H\odot \env_{\mathcal C}H} \node{\Env_{\mathcal C}H}\arrow{se,l}{\varkappa'}
 \\
  \node[2]{}\arrow{ne,r}{\bullet'} \node[2]{\Env_{\mathcal C}H\odot \Env_{\mathcal C}H}
 \\
 \node{\Env_{\mathcal C}H}\arrow{ne,-}\arrow{se,l}{\varkappa'}
 \\
 \node[2]{\Env_{\mathcal C}H\odot \Env_{\mathcal C}H}\arrow[2]{ne,r,3}{\bullet'\odot\bullet'}
 \\
\end{diagram}
 $$
Here the upper base is commutative (as \eqref{varkappa-circ-*=*-circledast-*-circ-varkappa} for $H$), as well as the lateral faces, and the left morphism $\env_{\mathcal C}H$ is an epimorphism. This means that the similar equality is true in $\Env_{\mathcal C}H$:
$$
\varkappa'\circ\bullet'=\bullet'\odot\bullet'\circ\varkappa'.
$$
\epr

\paragraph{Reflexivity with respect to an envelope.}

Suppose $(\env,\Env)$ is an envelope in the category $\InvSteAlg$ of involutive stereotype algebras (in the sense of general definition \cite{Akbarov-env}).

Let us say that an involutive stereotype Hopf algebra $H$ with respect to the tensor product $\circledast$ is {\it reflexive with respect to the envelope}\label{DEF:reflexiv-otn-obolochki} $\Env$, if its envelope $\Env H$ has a structure of involutive stereotype Hopf algebra in the category $({\tt Ste},\odot)$ such that the following two conditions hold:
\bit{
\item[(i)] a morphism of the envelope $\env H:H\to \Env H$ is a homomorphism of Hopf algebras in the sense that the following diagrams are commutative:
\beq\label{DIAG:reflex-otn-obolochki-1}
 \xymatrix @R=2.pc @C=2.pc
{
& H\odot H\ar[dr]^{\quad \env H\odot \env H} & \\
H\circledast H\ar[ur]^{@}\ar[dr]^{\quad \env H\circledast \env H}\ar[dd]_{\mu} & & \Env H\odot \Env H\ar[dd]_{\mu_E} \\
& \Env H\circledast \Env H\ar[ur]^{@} & \\
H\ar[rr]^{\env H} && \Env H
}
\eeq
\beq\label{DIAG:reflex-otn-obolochki-2}
 \xymatrix @R=2.pc @C=2.pc
{
& H\odot H\ar[dr]^{\quad\env H\odot \env H} & \\
H\circledast H\ar[ur]^{@}\ar[dr]^{\quad \env H\circledast \env H} & & \Env H\odot \Env H \\
& \Env H\circledast \Env H\ar[ur]^{@} & \\
H\ar[rr]^{\env H}\ar[uu]^{\varkappa} && \Env H\ar[uu]^{\varkappa_E}
}
\eeq
\beq\label{DIAG:reflex-otn-obolochki-3}
 \xymatrix @R=2.pc @C=2.pc
{
H\ar[rr]^{\env H} & & \Env H\\
& \C\ar[ul]^{\iota}\ar[ur]_{\iota_E} &
}\qquad
 \xymatrix @R=2.pc @C=2.pc
{
H\ar[rr]^{\env H}\ar[dr]_{\e} & & \Env H\ar[dl]^{\e_E} \\
& \C &
}
\eeq
\beq\label{DIAG:reflex-otn-obolochki-4}
 \xymatrix @R=3.pc @C=4.pc
{
H\ar[r]^{\env H}\ar[d]_{\sigma} & \Env H\ar[d]^{\sigma_E} \\
H\ar[r]^{\env H} & \Env H
}\qquad
 \xymatrix @R=3.pc @C=4.pc
{
H\ar[r]^{\env H}\ar[d]_{\bullet} & \Env H\ar[d]^{\bullet_E} \\
H\ar[r]^{\env H} & \Env H
}
\eeq
-- here $@$ is the Grothendieck transform \cite{Akbarov}, $\mu$, $\iota$, $\varkappa$, $\e$, $\sigma$, $\bullet$ -- are the structure murphisms (multiplication, unit, comultiplication, counit, antipode, involution) in $H$, and $\mu_E$, $\iota_E$, $\varkappa_E$, $\e_E$, $\sigma_E$, $\bullet_E$ the corresponding structure morphisms in $\Env H$.

\item[(ii)]\label{(env-H)^star:H^star-gets-(Env-H)^star} the mapping $(\env H)^\star:H^\star\gets (\Env H)^\star$, dual to the morphism of envelope $\env H:H\to\Env H$, is an envelope in the same sense:
$$
(\env H)^\star=\env (\Env H)^\star
$$

}\eit

\brem
Suppose the envelope $\env: H\to \Env H$ and the morphism $\env H\odot\env H\circ @:H\circledast H\to \Env H\odot\Env H$ are epimorphisms of stereotype spaces (this means that the sets of values are dense in the ranges). Then the envelope $\Env H$ can have at most unique structure of involutive Hopf algebra in $({\tt Ste},\odot)$ satisfying the conditions  (i) and (ii).
\erem
\bpr
The morphism $\iota_E$ must be the composition of $\iota$ and $\env H$, thus its uniqueness is seen immediately. The epimorphy of $\env H$ implies the uniqueness of $\varkappa_E$, $\e_E$, $\sigma_E$, $\bullet_E$. And the epimorphy of  $\env H\odot\env H\circ @$ the uniqueness of $\mu_E$.
\epr

It is convenient to display the conditions (i) and (ii) as a diagram
 \beq\label{obshaya-diagramma-refleksivnosti}
 \xymatrix @R=1.pc @C=1.pc
 {
 H
 & \ar@{|->}[r]^{\env} & &
\Env H
 \\
 & & &
 \ar@{|->}[d]^{\star}
 \\
 \ar@{|->}[u]^{\star}
 & & &
 \\
 H^\star
 & &
 \ar@{|->}[l]_{\env}
 &
 (\Env H)^\star
 }
 \eeq
which we call the {\it reflexivity diagram}, and which we endow the following sense:
 \bit{
\item[1)] in the corners  of the square there are involutive Hopf algebras, and $H$ is the Hopf algebra in $({\tt Ste},\circledast)$, then it follows the Hopf algebra $\Env H$ in $({\tt Ste},\odot)$, and after that the categpries $({\tt Ste},\circledast)$ and $({\tt Ste},\odot)$ alternate,

\item[2)] the alternation of the operations $\env$ and $\star$ (no matter where we start) on the fourth step returns us back to the initial Hopf algebra (certainly, up to an isomorphism of functors).
 }\eit

The sense of the term ``reflexivity'' here is as follows. Denote the single successive application of the operations $\env$ and $\star$ by some symbol, for example, $\widehat{\ }\,\,$,
$$
\widehat H:=(\Env H)^\star
$$
Since $\Env H$ has a unique structure of Hopf algebra with respect to $\odot$, the dual space $\widehat H=(\Env H)^\star$ has a structure of involutive Hopf algebra with respect to $\circledast$. Moreover, $\widehat H=(\Env H)^\star$ is a Hopf algebra, reflexive with respect to $\Env$, since the application of $\star$ to the diagrams \eqref{DIAG:reflex-otn-obolochki-1}-\eqref{DIAG:reflex-otn-obolochki-4} gives the same diagrams with the replacement $H$ by $\widehat H=(\Env H)^\star$ (we use here the condition (ii) on page \pageref{(env-H)^star:H^star-gets-(Env-H)^star}).

Let us call $\widehat H=(\Env H)^\star$ the {\it dual Hopf algebra to
$H$ with respect to the envelope $\Env$}\index{Hopf algebra!dual with respect to an envelope}. The diagram \eqref{obshaya-diagramma-refleksivnosti} means that $H$ is naturally isomorphic to its second dual Hopf algebra in this sense:
 \beq\label{H-cong-(H^*)^*}
H\cong \widehat{\widehat H}
 \eeq

\paragraph{Continuous reflexivity.} In the special case when in the definition on page \pageref{DEF:reflexiv-otn-obolochki} $\Env$ means the continuous envelope $\Env_{\mathcal C}$, we call $H$ a {\it continuously reflexive Hopf algebra}, and the algebra $\widehat H=(\Env_{\mathcal C} H)^\star$ the {\it continuously dual Hopf algebra} to $H$.

Theorems \ref{TH:E(K(G))=C(G)} and \ref{TH:E(H)-dlya-Moore} imply the main result of \cite{Kuznetsova}:

 \btm\footnote{See Errata on page \pageref{Errata}.}\label{TH:nepr-dvoistvennost}
If $G$ is a Moore group, then the algebras ${\mathcal C}^\star(G)$ and ${\mathcal
K}(G)$ are continuously reflexive, and the reflexivity diagram for them is:
 \beq\label{chetyrehugolnik-C-C*}
 \xymatrix @R=1.pc @C=2.pc
 {
 {\mathcal C}^\star(G)
 & \ar@{|->}[r]^{\Env_{\mathcal C}} & &
 \Env_{\mathcal C}{\mathcal C}^\star(G)
 \\
 & & &
 \ar@{|->}[d]^{\star}
 \\
 \ar@{|->}[u]^{\star}
 & & &
 \\
 {\mathcal C}(G)
 & &
 \ar@{|->}[l]_{\Env_{\mathcal C}}
 &
 {\mathcal K}(G)
 }
 \eeq
 \etm

\bex Proposition \ref{PROP:Env_C-C^star(G)=C(widehat(G))} implies that for Abelian locally compact groups $A$ the reflexivity diagram has the form
 \beq\label{chetyrehugolnik-C-C*-F}
 \xymatrix @R=1.pc @C=2.pc
 {
 {\mathcal C}^\star(A)
 & \ar@{|->}[r]^{{\mathcal F}_A} & &
 {\mathcal C}(\widehat{A})
 \\
 & & &
 \ar@{|->}[d]^{\star}
 \\
 \ar@{|->}[u]^{\star}
 & & &
 \\
 {\mathcal C}(A)
 & &
 \ar@{|->}[l]_{{\mathcal F}_{\widehat{A}}}
 &
 {\mathcal C}^\star(\widehat{A})
  }
\eeq
(where $\widehat{A}$ is the Pontryagin dual group to $A$, and $\mathcal F$ the Fourier transform, defined in \eqref{Fourier-transform}).
\eex

\bex Proposition \ref{PROP:nepr-obol-diskr-gruppy} implies that for discrete Moore groups $D$ the reflexivity diagram is:
$$
 \xymatrix @R=1.pc @C=2.pc
 {
 \C_D
 & \ar@{|->}[r]^{\Env_{\mathcal C}} & &
 C^*(D)
 \\
 & & &
 \ar@{|->}[d]^{\star}
 \\
 \ar@{|->}[u]^{\star}
 & & &
 \\
 \C^D
 & &
 \ar@{|->}[l]_{\Env_{\mathcal C}}
 &
 {\mathcal K}(D)
 }
$$
Here $C^*(G)$ is the usual $C^*$-algebra of the group $G$ \cite{Dixmier}, and the algebra ${\mathcal K}(D)$ as a set  (and as an algebra) coincides with the Fourier-Stiltjes algebra $B(G)$ of the group $G$ \cite{Eymard}, and the difference is that ${\mathcal K}(D)$ has weaker topology (we already noticed this in the proof of Lemma \ref{LM:Spec-K(G)-dlya-diskretnoi-gruppy}).
\eex

\bex Proposition \ref{PROP:nepr-obol-komp-gruppy} implies that for compact groups $K$ the reflexivity diagram is:
$$
 \xymatrix @R=1.pc @C=2.pc
 {
 {\mathcal C}^\star(K)
 & \ar@{|->}[r]^{\Env_{\mathcal C}} & &
\prod_{\pi\in\widehat{K}}{\mathcal B}(X_\pi)
 \\
 & & &
 \ar@{|->}[d]^{\star}
 \\
 \ar@{|->}[u]^{\star}
 & & &
 \\
 {\mathcal C}(K)
 & &
 \ar@{|->}[l]_{\Env_{\mathcal C}}
 &
 \Trig(K)
 }
$$
\eex

\paragraph{Groups, discerned by $C^*$-algebras.}
Let us say that a locally compact group $G$ {\it is discerned by $C^*$-algebras}\label{DEF:gruppy-razlich-C^*-algebrami}, if (continuous involutive) homomorphisms of its measure algebra ${\mathcal C}^\star(G)\to B$ into various $C^*$-algebras $B$ separate elements of $G$ (with the injection of  $G$ into ${\mathcal C}^\star(G)$ by delta-functions). Clearly, if the group algebra of measures ${\mathcal C}^\star(G)$ is continuously reflexive, then the group $G$ is discerned by $C^*$-algebras, so this class of groups is interesing in estimating of how wide one can try to generalize Theorem \ref{TH:nepr-dvoistvennost}. The very Theorem \ref{TH:nepr-dvoistvennost} implies that all Moore groups are discerned by $C^*$-algebras. In the work by Yu.~N.~Kuznetsova \cite{Kuznetsova} it is shown that all SIN-groups have this property. On the other hand, by I.~M.~Singer's results \cite[Corollary 5]{Singer}, in the class of connected Lie groups only groups of the form $\R^n\times K$, where $n\in\N$, and $K$ is a compact Lie group, are discerned by $C^*$-algebras.

\btm[D.~Luminet, A.~Valette, \cite{Luminet-Valette}]\label{TH:Luminet-Valette} If a Lie group $G$ is discerned by  $C^*$-algebras, then $G$ is a linear group (i.e. $G$ can be embedded as a closed subgroup into a full linear group $GL_n(\C)$).
\etm

\section{Smooth envelopes and smooth duality}

\subsection{Joined self-adjoint elements and the system of partial derivatives}

\paragraph{Multi-indices.}

Let $d\in\N$ -- be a natural number\footnote{Everywhere natural numbers $\N$ are non-negative integers: $\N=\{d\in\Z:\ d\ge 0\}$.}. Let us call a {\it multi-index} of the length $d$ an arbitrary finite sequence of the length $d$ of natural numbers
$$
k=(k_1,...,k_d),\qquad k_i\in\N.
$$
For each two multi-indices $k,l\in\N^d$ the inequality $l\le k$ is defined coordinate-wise:
$$
l\le k\qquad\Longleftrightarrow\qquad \forall i=1,...,d\quad l_i\le k_i.
$$
A {\it sum} of two multi-indices $k,l\in\N^d$ is the multi-index
$$
k+l=(k_1+l_1,...,k_d+l_d).
$$
If $l\le k$, then the {\it subtraction} is again defined coordinate-wise:
$$
k-l=(k_1-l_1,...,k_d-l_d).
$$
The {\it order} and the {\it factorial} of a multi-index $k\in\N^d$ are defined by the equalities
$$
\abs{k}=k_1+...+k_d,\qquad k!=k_1!\cdot...\cdot k_d!.
$$
According to the last formula, the binomial coefficient is
$$
\begin{pmatrix}k \\ l\end{pmatrix}=\frac{k!}{l!\cdot (k-l)!}
$$

\paragraph{Algebras of power series with coefficients in a given algebra.}
Let $A$ be an arbitrary involutive stereotype algebra. Consider the algebra
$$
A[[d]]=A^{\N^d},
$$
consisting of all mappings $x:\N^d\to A$, or, what is the same, of families $x=\{x_k;\ k\in\N^d\}$ of elements from $A$, indexed by multi-indices of the length $d$. This set $A[[d]]$ is endowed with the topology of coordinate-wise convergence
$$
x^i\overset{A[[d]]}{\underset{i\to\infty}{\longrightarrow}}x\qquad\Longleftrightarrow\qquad
\forall k\in\N^d\quad
x_k^i\overset{A}{\underset{i\to\infty}{\longrightarrow}}x_k,
$$
and the algebraic operations on $A[[d]]$ -- involution, sum, multiplication by scalar and multiplication  -- are defined by formulas
 \begin{align}
& (x^\bullet)_k=(x_k)^\bullet, && x\in A[[d]],\ k\in \N^d \label{involutsiya-v-B[[d]]}\\
& (x+y)_k=x_k+y_k, && x,y\in A[[d]],\ k\in \N^d \label{summa-v-B[[d]]}\\
& (\lambda\cdot x)_k=\lambda\cdot x_k, && \lambda\in\C, \ x\in A[[d]],\  k\in \N^d  \label{umnopzh-na-skalar-v-B[[d]]}\\
& (x\cdot y)_k=\sum_{0\le l\le k}x_{k-l}\cdot y_l, &&  x,y\in A[[d]],\ k\in \N^d \label{proizv-v-B[[d]]}
 \end{align}
The unit in $A[[d]]$ is the family
\beq\label{1_(B[[d]])}
\underset{\scriptsize\begin{matrix}\text{\rotatebox{90}{$\owns$}} \\ A[[d]]\end{matrix}}{1}\kern-5pt{}_k=\begin{cases}
1,& k=0 \\ 0, & k\ne 0\end{cases}
\eeq

It is convenient to represent the elements of $A[[d]]$ as power series of $d$ variables $\tau_1,...,\tau_d$:
$$
x=\sum_{k\in\N^d}x_k\cdot\tau^k,
$$
where $\tau^k$ means the formal product
$$
\tau^k=\tau^{k_1}\cdot...\cdot\tau^{k_d},
$$
and the following identities are assumed to hold:
$$
(\tau^k)^\bullet=\tau^k,\qquad a\cdot\tau^k=\tau^k\cdot a,\qquad \tau^k\cdot\tau^l=\tau^{k+l},\qquad a\in A,\quad \tau\in\N^d.
$$
Then the sum, the multiplication by a scalar and the multiplication in $A[[d]]$ are defined by formulas for the power series:
 \begin{align}\label{alg-operatsii-kak-ryady}
& x+y=\sum_{k\in\N^d}(x_k+y_k)\cdot\tau^k,
&& \lambda\cdot x=\sum_{k\in\N^d}(\lambda\cdot x_k)\cdot\tau^k,
&& x\cdot y=\sum_{k\in\N^d}\l\sum_{0\le l\le k}x_{k-l}\cdot y_l\r\cdot\tau^k.
 \end{align}

Another way is to represent elements of $A[[d]]$ as the Taylor series of variables $\tau_1,...,\tau_d$:
\beq\label{predstavlenie-ryadom-Teilora}
x=\sum_{k\in\N^d}\frac{x^{(k)}}{k!}\cdot\tau^k,
\eeq
In this representation the coefficients $x^{(k)}$ of the series are connected with the usual coefficients $x_k$ by the formulas
$$
x^{(k)}=k!\cdot x_k,
$$
and the formulas for algebraic operatioosn \eqref{involutsiya-v-B[[d]]}-\eqref{proizv-v-B[[d]]} have the form
 \begin{align}\label{(xy)^(k)}
&(x^\bullet)^{(k)}=(x^{(k)})^\bullet &&
(\lambda\cdot x)^{(k)}=\lambda\cdot x^{(k)} &&
(x+y)^{(k)}=x^{(k)}+y^{(k)} &&
(x\cdot y)^{(k)}=\sum_{0\le l\le k}\begin{pmatrix}k \\ l\end{pmatrix}\cdot x^{(k-l)}\cdot y^{(l)}
 \end{align}

\paragraph{Algebras with the joined self-adjoint nilpotent elements.}

 \bit{
\item[$\bullet$] Let again $A$ be an involutive stereotype algebra, $d\in\N$ and $m\in\N^d$. Denote by
$I_m$ the closed ideal in the algebra $A[[d]]$ of power series with coefficients in
$A$, consisting of series whose coefficients with indices $k\le m$ vanish:
$$
I_m=\{x\in A[[d]]:\ \forall k\in\N^d\quad k\le m \quad\Longrightarrow\quad x_k=0\}.
$$
The quotient algebra
$$
A[m]:=A[[d]]/I_m
$$
is called the {\it algebra $A$ with joined self-adjoint nilpotent elements (of order $m$)}.
 }\eit

Denote by symbol $\N[m]$ the set of multi-indices not greater than $m$
$$
\N[m]=\{k\in\N^d:\ k\le m\}.
$$
Then $A[m]$ can be treated as the space of families $x=\{x_k;\ k\in \N[m]\}$ of elements from $A$, indexed by multi-indices $k\in\N[m]$. The involution, the sum, the multiplication by scalar and the multiplication in $A[m]$ are defined by formulas
 \begin{align}
& (x^\bullet)_k=(x_k)^\bullet, && x\in A[m],\ k\in \N[m] \label{involutsiya-v-B[m]}\\
& (x+y)_k=x_k+y_k, && x,y\in A[m],\ k\in \N[m] \label{summa-v-B[m]}\\
& (\lambda\cdot x)_k=\lambda\cdot x_k, && \lambda\in\C, \ x\in A[m],\  k\in \N[m]  \label{umnopzh-na-skalar-v-B[m]}\\
& (x\cdot y)_k=\sum_{0\le l\le k}x_{k-l}\cdot y_l, &&  x,y\in A[m],\ k\in
\N[m] \label{proizv-v-B[m]}
 \end{align}
The unit in $A[m]$ is, certainly, the family
\beq\label{1_(B[[d]])}
\underset{\scriptsize\begin{matrix}\text{\rotatebox{90}{$\owns$}} \\
A[m]\end{matrix}}{1_k}\kern-5pt=\begin{cases} 1,& k=0 \\ 0, & k\ne
0\end{cases}
\eeq

It is convenient to represent the elements of $A[m]$ as polynomials of degree $n$ of $d$ variables $\tau_1,...,\tau_d$:
\beq\label{predstavl-B[m]-ryadami}
x=\sum_{k\in\N[m]}x_k\cdot\tau^k=\sum_{k\in\N[m]}\frac{x^{(k)}}{k!}\cdot\tau^k,
\eeq
where $\tau^k$ is the formal product
$$
\tau^k=\tau^{k_1}\cdot...\cdot\tau^{k_d},
$$
and
\beq\label{tozhdestva-dlya-prisoedinennyh-elementov}
(\tau_i)^\bullet=\tau_i, \qquad a\cdot\tau_i=\tau_i\cdot a,\qquad \tau_i\cdot\tau_j=\tau_j\cdot\tau_i,
\qquad\tau_i^{m_i+1}=0,\qquad a\in A,\quad i=1,...,d.
\eeq
In this representation the variables $\tau_1,...,\tau_d$ can be treated as the joined elements to the algebra $A$ satisfying the conditions \eqref{tozhdestva-dlya-prisoedinennyh-elementov} (and this justifies the name that we gave to the algebra $A[m]$). Then the sum, the multiplication by scalar and the multiplication in $A[m]$ are defined by the very same formulas \eqref{alg-operatsii-kak-ryady} as for $A[[d]]$. In particular, in the representation of elements by the Taylor polynomial (the second equality in \eqref{predstavl-B[m]-ryadami}) the formulas for the algebraic operations \eqref{(xy)^(k)} are preserved.

\btm\label{TH:norm(x)=sum_(k-le-n)norm(x_k)_B} Let $B$ be an involutive Banach algebra with an involutive submultiplicative norm $\norm{\cdot}_B$. Then for each multi-index $n\in\N^d$ the algebra $B[n]$ is an involutive Banach algebra with involutive submultiplicative norm
\beq\label{norm(x)=sum_(k-le-n)norm(x_k)_B}
\norm{x}=\sum_{k\le n}\norm{x_k}_B,\qquad x\in B[n].
\eeq
\etm
\bpr
Put $M=\{(k,l)\subseteq\N[n]^2; \ l\le k\}$ and note that the mapping $(k,l)\in M\mapsto(k-l,l)\in\N[n]^2$ is injective:
$$
(k,l),(k',l')\in M\quad\&\quad (k-l,l)=(k'-l',l')\quad\Longrightarrow\quad l=l'\quad\&\quad k=k'.
$$
This implies the third inequality in the chain
\begin{multline*}
\norm{x\cdot y}=\sum_{k\le n}\norm{(x\cdot y)_k}_B=\sum_{k\le n}\norm{\sum_{l\le k}x_{k-l}\cdot y_l}_B\le
\sum_{k\le n}\sum_{l\le k}\norm{x_{k-l}\cdot y_l}_B\le
\sum_{k\le n}\sum_{l\le k}\norm{x_{k-l}}_B\cdot \norm{y_l}_B=\\=\sum_{(k,l)\in M}\norm{x_{k-l}}_B\cdot \norm{y_l}_B
\le \sum_{(m,l)\in\N[n]^2}\norm{x_m}_B\cdot \norm{y_l}_B=\sum_{m\in\N[n]}\norm{x_m}_B\cdot\sum_{l\in\N[n]} \norm{y_l}_B=\norm{x}\cdot\norm{y}.
\end{multline*}
\epr

Further we shall be interested almost exclusively in the case when $A$ is a $C^*$-algebra. The following example shows that the property of being a $C^*$-algebra is not inherited when passing from $A$ to $A[m]$.

\bex For non-vanishing $m$ and $d$ the algebra $A[m]$ can't be a $C^*$-algebra.
\eex
\bpr
Suppose the topology of $A[m]$ is generated by a $C^*$-norm. Consider a joined element $\tau_i$, $i\in\{1,...,n\}$. Let $B$ be a closed unital subalgebra in $A[m]$, generated by this element $\tau_i$. Then $B$ is a commutative  $C^*$-algebra, hence it is isomorphic to an algebra ${\mathcal C}(K)$ of functions on some compact space. On the other hand, the last condition in \eqref{tozhdestva-dlya-prisoedinennyh-elementov} means that $\tau_i$ is a nilpotent element:
$$
\tau_i^{m_i+1}=0.
$$
This is impossible since the algebras $B={\mathcal C}(K)$ don't have non-zero nilpotent elements.
\epr

For any two multi-indices $m\in\N^d$ and $n\in\N^{d'}$ (not necessarily $d=d'$) we define their {\it direct sum} $m\oplus n$ as a multi-index of the length $d+d'$, i.e. an element of $\N^{d+d'}$, by the formula
\beq\label{(k-oplus-l)_i}
(m\oplus n)_i=\begin{cases}m_i, & 1\le i\le d\\ n_{i-d}, & d<i\le d' \end{cases}.
\eeq
Then automatically
\beq\label{N[m]-times-N[n]=N[m-oplus-n]}
\N[m]\times\N[n]=\N[m\oplus n]
\eeq

\bprop\label{PROP:A[m]-circledast-B}
For any two involutive stereotype algebras $A$ and $B$ and any two multi-indices $m\in\N^d$ and $n\in\N^{d'}$ the following natural isomorphisms of involutive stereotype algebras hold:
\beq\label{A[m][n]=A[m-oplus-n]}
\big(A[m]\big)[n]\cong A[m\oplus n]\cong A[n\oplus m]\cong \big(A[n]\big)[m]
\eeq
\beq\label{A[m]-circledast-B}
A[m]\circledast B\cong (A\circledast B)[m]\cong A\circledast \big(B[m]\big)
\eeq
If $M$ is a paracompact locally compact space, then
\beq\label{C(M,B)[n]=C(M,B[n])}
{\mathcal C}(M,B)[n]\cong {\mathcal C}(M,B[n]).
\eeq
There exists also a natural homomorphism of stereotype involutive algebras
\beq\label{A[m]-oplus-B[n]-to-(A-oplus-B)[m-oplus-n]}
\ph:A[m]\oplus B[n]\to (A\oplus B)[m\oplus n]
\eeq
Suppose in addition that $A$ and $B$ are Banach algebras with submultiplicative norms, and each time the norm of the direct sum is defined as maximum
$$
\norm{x\oplus y}=\max\{\norm{x},\norm{y}\}
$$
and in the algebras with the joined self-adjoint elements as the sum of norms of components (by formula \eqref{norm(x)=sum_(k-le-n)norm(x_k)_B}). Then the norm of $\ph$ is estimated as follows:
\beq\label{1-le-norm(ph)-le-2}
1\le \norm{\ph}\le 2,
\eeq
\eprop
\bpr 1. Formula \eqref{A[m][n]=A[m-oplus-n]} is defined by the chain
$$
\big(A[m]\big)[n]=\big(A^{\N[m]}\big)^{\N[n]}=A^{\N[m]\times\N[n]}=\eqref{N[m]-times-N[n]=N[m-oplus-n]}=A^{\N[m\oplus n]}=A[m\oplus n]
$$

2. To prove \eqref{A[m]-circledast-B} let us define a map
$$
\gamma: A[m]\circledast B\to (A\circledast B)[m]
$$
by formula
\beq\label{DEF:A[m]-circledast-B->(A-circledast-B)[m]}
\gamma(x\circledast b)_k=x_k\circledast b,\qquad x\in A[m],\quad b\in B.
\eeq
(to each family $x=\{x_k:\ k\in\N[m]\}$ of elements in $A$ and to each element $b\in B$ the map $\gamma$ assigns the family $\gamma(x\circledast b)=\{\gamma(x\circledast b)_k;\ k\in\N[m]\}$ of elements in $A\circledast B$, defined by \eqref{DEF:A[m]-circledast-B->(A-circledast-B)[m]}). This map is an isomorphism of stereotype spaces, since the tensor product $\circledast$ is distributive with the operation of taking direct sum \cite[(2.52)]{Akbarov-env}. Let us check that it preserves multiplication: for each families $x=\{x_k:\ k\in\N[m]\}$ and $x'=\{x'_k:\ k\in\N[m]\}$ from $A$ and for each elements $b,b'\in B$ we have:
\begin{multline*}
\gamma\Big((x\circledast b)\cdot(x'\circledast b')\Big)_k=\gamma\big((x\cdot x')\circledast(b\cdot b')\big)_k=\eqref{DEF:A[m]-circledast-B->(A-circledast-B)[m]}=(x\cdot x')_k\circledast(b\cdot b')=\eqref{proizv-v-B[m]}=\\=\Big(\sum_{0\le l\le k}x_{k-l}\cdot x'_l\Big)\circledast (b\cdot b')=
\sum_{0\le l\le k}(x_{k-l}\cdot x'_l)\circledast (b\cdot b')=\sum_{0\le l\le k}(x_{k-l}\circledast b)\cdot (x'_l\circledast b')=\sum_{0\le l\le k}(x_{k-l}\circledast b)\cdot (x'_l\circledast b')=\\=\eqref{DEF:A[m]-circledast-B->(A-circledast-B)[m]}=\sum_{0\le l\le k}\gamma(x\circledast b)_{k-l}\cdot \gamma(x'\circledast b')_l=\eqref{proizv-v-B[m]}=\Big(\gamma(x\circledast b)\cdot\gamma(x'\circledast b')\Big)_k
\end{multline*}
This proves the first equality in \eqref{A[m]-circledast-B}. The second one is proved similarly.

3. The identity \eqref{C(M,B)[n]=C(M,B[n])} is evident.

4. Suppose $\sigma_1,...,\sigma_d$ is a sequence of joined self-adjoint nilpotent elements to $A$ in $A[m]$, and $\tau_1,...,\tau_{d'}$ the same sequence in $B[n]$. Consider the sequence of joined self-adjoint nilpotent elements to the algebra $A\oplus B$ in $(A\oplus B)[m\oplus n]$, and let us assign to its elements the notations $\overline{\sigma_i}$ and $\overline{\tau_j}$, arranging them in such a way that initially the elements $\overline{\sigma_i}$ appear in the growing order, and after them $\overline{\tau_j}$:
$$
\overline{\sigma_1},...,\overline{\sigma_d},\overline{\tau_1},...,\overline{\tau_{d'}}.
$$
Put
$$
\overline{\sigma}^k=\overline{\sigma_1}^{k_1}\cdot...\cdot \overline{\sigma_d}^{k_d},\qquad
\overline{\tau}^l=\overline{\tau_1}^{l_1}\cdot...\cdot \overline{\tau_{d'}}^{l_{d'}},\qquad k\in\N[m],\quad l\in\N[n].
$$
Then the homomorphism $\ph$ in the formula \eqref{A[m]-oplus-B[n]-to-(A-oplus-B)[m-oplus-n]} can be defined by the formula
$$
\ph(1_{A[m]}\oplus 1_{B[n]})=1_{(A\oplus B)[m\oplus n]},\quad \ph(\sigma_i\oplus 0_{B[n]})=(1_A\oplus 0_B)\cdot\overline{\sigma_i},\quad \ph(0_{A[m]}\oplus\tau_j)=(0_A\oplus 1_B)\cdot\overline{\tau_j}.
$$
or, equivalently, by the rule
$$
\ph\Big(\underbrace{\sum_{0\le k\le m}x_k\cdot\sigma^k}_{\scriptsize\begin{matrix}\text{\rotatebox{90}{$\owns$}}\\ A[m]\end{matrix}}\oplus\underbrace{\sum_{0\le l\le n} y_l\cdot\tau^l}_{\scriptsize\begin{matrix}\text{\rotatebox{90}{$\owns$}}\\ B[n]\end{matrix}}\Big)=\sum_{0\le k\le m}(\underset{\scriptsize\begin{matrix}\text{\rotatebox{90}{$\owns$}}\\ A\end{matrix}}{x_k}\oplus \underset{\scriptsize\begin{matrix}\text{\rotatebox{90}{$\owns$}}\\ B\end{matrix}}{0})\cdot\overline{\sigma}^k+ \sum_{0\le l\le n} (\underset{\scriptsize\begin{matrix}\text{\rotatebox{90}{$\owns$}}\\ A\end{matrix}}{0}\oplus \underset{\scriptsize\begin{matrix}\text{\rotatebox{90}{$\owns$}}\\ B\end{matrix}}{y_l})\cdot\overline{\tau}^l
$$
This mapping preserves the unit and is multiplicative:
\begin{multline*}
\ph\Big((x\oplus y)\cdot (x'\oplus y')\Big)=\ph\left(\Big(\sum_{0\le k\le m}x_k\cdot \sigma^k\oplus\sum_{0\le l\le n}y_l\cdot\tau^l\Big) \cdot \Big(\sum_{0\le k'\le m}x'_{k'}\cdot \sigma^{k'}\oplus\sum_{0\le l'\le n}y'_{l'}\cdot\tau^{l'}\Big)\right)=\\=
\ph\left(\Big(\sum_{0\le k\le m}x_k\cdot \sigma^k \cdot \sum_{0\le k'\le m}x'_{k'}\cdot \sigma^{k'}\Big)\oplus\Big(\sum_{0\le l\le n}y_l\cdot\tau^l\cdot\sum_{0\le l'\le n}y'_{l'}\cdot\tau^{l'}\Big)\right)=
\\=
\ph\left(\sum_{0\le k\le m}\Big(\sum_{0\le p\le k}x_{k-p}\cdot x'_{p}\Big)\cdot \sigma^k\oplus\sum_{0\le l\le n}\Big(\sum_{0\le q\le l}y_{l-q}\cdot y'_q\Big)\cdot\tau^l\right)=\\=
\sum_{0\le k\le m}\Big(\sum_{0\le p\le k}x_{k-p}\cdot x'_{p}\oplus 0\Big)\cdot \overline{\sigma}^k+ \sum_{0\le l\le n}\Big(0\oplus\sum_{0\le q\le l}y_{l-q}\cdot y'_q\Big)\cdot\overline{\tau}^l=\\=
\sum_{0\le k\le m}(x_k\oplus 0)\cdot \overline{\sigma}^k \cdot \sum_{0\le k'\le m}(x'_{k'}\oplus 0)\cdot \overline{\sigma}^{k'}+\sum_{0\le l\le n}(0\oplus y_l)\cdot\overline{\tau}^l\cdot\sum_{0\le l'\le n}(0\oplus y'_{l'})\cdot\overline{\tau}^{l'}=\\=
\sum_{0\le k\le m}(x_k\oplus 0)\cdot \overline{\sigma}^k \cdot \sum_{0\le k'\le m}(x'_{k'}\oplus 0)\cdot \overline{\sigma}^{k'}+\sum_{0\le l\le n}(0\oplus y_l)\cdot\overline{\tau}^l\cdot\sum_{0\le l'\le n}(0\oplus y'_{l'})\cdot\overline{\tau}^{l'}+\\+
\underbrace{\sum_{0\le k\le m}(x_k\oplus 0)\cdot \overline{\sigma}^k \cdot\sum_{0\le l'\le n}(0\oplus y'_{l'})\cdot\overline{\tau}^{l'}}_{\scriptsize\begin{matrix}\| \\ 0\end{matrix}}+\underbrace{\sum_{0\le l\le n}(0\oplus y_l)\cdot\overline{\tau}^l\cdot
\sum_{0\le k'\le m}(x'_{k'}\oplus 0)\cdot \overline{\sigma}^{k'}}_{\scriptsize\begin{matrix}\| \\ 0\end{matrix}}=\\=
\left(
\sum_{0\le k\le m}(x_k\oplus 0)\cdot \overline{\sigma}^k \oplus \sum_{0\le l\le n}(0\oplus y_l)\cdot\overline{\tau}^l
\right)\cdot \left(
\sum_{0\le k'\le m}(x'_{k'}\oplus 0)\cdot \overline{\sigma}^{k'} \oplus \sum_{0\le l'\le n}(0\oplus y'_{l'})\cdot\overline{\tau}^{l'}
\right)=\\=
\ph\Big(\sum_{0\le k\le m}x_k\cdot \sigma^k\oplus\sum_{0\le l\le n}y_l\cdot\tau^l\Big) \cdot\ph\Big(\sum_{0\le k'\le m}x'_{k'}\cdot \sigma^{k'}\oplus\sum_{0\le l'\le n}y'_{l'}\cdot\tau^{l'}\Big)=
\ph(x\oplus y)\cdot \ph(x'\oplus y')
\end{multline*}
In \eqref{1-le-norm(ph)-le-2} the first inequality is proved by the chain
\begin{multline*}
\norm{\ph\Big(\sum_{0\le k\le m}x_k\cdot \sigma^k\oplus\sum_{0\le l\le n}y_l\cdot\tau^l\Big)}_{(A\oplus B)[m\oplus n]}=\norm{\sum_{0\le k\le m}(x_k\oplus 0)\cdot \overline{\sigma}^k \oplus \sum_{0\le l\le n}(0\oplus y_l)\cdot\overline{\tau}^l}_{(A\oplus B)[m\oplus n]}=\\=
\sum_{0\le k\le m}\norm{x_k\oplus 0}_{A\oplus B}+\sum_{0\le l\le n}\norm{0\oplus y_l}_{A\oplus B}=
\sum_{0\le k\le m}\norm{x_k}_A+\sum_{0\le l\le n}\norm{y_l}_B=\\=\norm{\sum_{0\le k\le m}x_k\cdot \sigma^k}_{A[m]}+\norm{\sum_{0\le l\le n}y_l\cdot\tau^l}_{B[n]}\ge
\max\left\{\norm{\sum_{0\le k\le m}x_k\cdot \sigma^k}_{A[m]},\norm{\sum_{0\le l\le n}y_l\cdot\tau^l}_{B[n]}\right\}=\\=
\norm{\sum_{0\le k\le m}x_k\cdot \sigma^k\oplus\sum_{0\le l\le n}y_l\cdot\tau^l}_{A[m]\oplus B[n]}
\end{multline*}
and the second by the chain
\begin{multline*}
\norm{\ph\Big(\sum_{0\le k\le m}x_k\cdot \sigma^k\oplus\sum_{0\le l\le n}y_l\cdot\tau^l\Big)}_{(A\oplus B)[m\oplus n]}=\norm{\sum_{0\le k\le m}(x_k\oplus 0)\cdot \overline{\sigma}^k \oplus \sum_{0\le l\le n}(0\oplus y_l)\cdot\overline{\tau}^l}_{(A\oplus B)[m\oplus n]}=\\=
\sum_{0\le k\le m}\norm{x_k\oplus 0}_{A\oplus B}+\sum_{0\le l\le n}\norm{0\oplus y_l}_{A\oplus B}=
\sum_{0\le k\le m}\norm{x_k}_A+\sum_{0\le l\le n}\norm{y_l}_B=\\=\norm{\sum_{0\le k\le m}x_k\cdot \sigma^k}_{A[m]}+\norm{\sum_{0\le l\le n}y_l\cdot\tau^l}_{B[n]}\le
2\max\left\{\norm{\sum_{0\le k\le m}x_k\cdot \sigma^k}_{A[m]},\norm{\sum_{0\le l\le n}y_l\cdot\tau^l}_{B[n]}\right\}=\\=
2\norm{\sum_{0\le k\le m}x_k\cdot \sigma^k\oplus\sum_{0\le l\le n}y_l\cdot\tau^l}_{A[m]\oplus B[n]}
\end{multline*}
\epr

\paragraph[Systems of partial derivatives.]{Systems of partial derivatives and morphisms with values in the algebras of power series.}

Let $A$ be an involutive stereotype algebra, and $B$ a $C^*$-algebra. A system of operators $D_k:A\to B$, $k\in\N[m]$, is called a {\it system of partial derivatives} on the algebra $A$ with coefficients in the algebra $B$, if it satisfies the following conditions:
\begin{align}
& D_k(a^\bullet)=D_k(a)^\bullet, && k\in\N[m], \label{DEF:sist-chastn-proizv-*} \\
& D_k(1)=\begin{cases}1, & k=0
\\
0, & k\ne 0
\end{cases}, && k\in\N[m],\ 1\in A, \label{DEF:sist-chastn-proizv-1} \\
& D_k(a\cdot b)=\sum_{0\le l\le k}\begin{pmatrix}k \\ l\end{pmatrix}\cdot D_{k-l}(a)\cdot D_l(b), && k\in\N[m],\ a,b\in A. \label{DEF:sist-chastn-proizv-2}
\end{align}
In particular this means that the operator $D_0:A\to B$ is an involutive homomorphism of algebras,
$$
D_0(a^\bullet)=D_0(a)^\bullet,\qquad D_0(1)=1,\qquad  D_0(a\cdot b)=D_0(a)\cdot
D_0(b),\qquad a,b\in A.
$$
For $\abs{k}=1$ the operators $D_k:A\to B$ are derivatives with respect to the homomorphism $D_0$:
$$
 D_k(a\cdot b)=D_k(a)\cdot D_0(b)+D_0(a)\cdot D_k(b),\qquad a,b\in A.
$$

\btm\label{TH:homomorphizm-v-B[m]<->chast-proizv-v-B} For each involutive stereotype algebra $A$ and for each  $C^*$-algebra $B$ the formula
 \beq\label{homomorphizm-v-B[m]<->chast-proizv-v-B}
D_k(a)=D(a)^{(k)},\qquad k\in\N[m],\quad a\in A,
 \eeq
or, equivalently, the formula
 \beq\label{D(a)=sum-D_k(a)/k!-tau^k}
D(a)=\sum_{k\in\N[m]}\frac{D_k(a)}{k!}\cdot\tau^k,\qquad a\in A,
 \eeq
establishes a one-to-one correspondence between homomorphisms of involutive stereotype algebras $D:A\to B[m]$ and systems of partial derivatives $\{D_k;\ k\in\N[m]\}$ on $A$ with coefficients in $B$. \etm

\bpr 1. If $D:A\to B[m]$ is a homomorphism of injective stereotype algebras, then we can define the mappings
$\{D_k;\ k\in\N[m]\}:A\to B$ by formula \eqref{homomorphizm-v-B[m]<->chast-proizv-v-B}, and we get, first, that for each $a\in A$
$$
D_k(a^\bullet)=D(a^\bullet)^{(k)}=\big(D(a)^\bullet\big)^{(k)}=\eqref{(xy)^(k)}=\big(D(a)^{(k)}\big)^\bullet=\big(D_k(a)\big)^\bullet,
$$
second,
$$
D_k(1)=D(1)^{(k)}=1^{(k)}=1_k=\eqref{1_(B[[d]])}=\begin{cases}1, & k=0
\\
0, & k\ne 0
\end{cases},
$$
and, third, for any $a,b\in A$
\begin{multline*}
D_k(a\cdot b)=\eqref{homomorphizm-v-B[m]<->chast-proizv-v-B}=D(a\cdot b)^{(k)}=\big(D(a)\cdot D(b)\big)^{(k)}=\eqref{(xy)^(k)}=\\=
\sum_{0\le l\le k}\begin{pmatrix}k \\ l\end{pmatrix}\cdot  D(a)^{(k-l)}\cdot D(b)^{(l)}=\eqref{homomorphizm-v-B[m]<->chast-proizv-v-B}=\sum_{0\le l\le k}\begin{pmatrix}k \\ l\end{pmatrix}\cdot  D_{k-l}(a)\cdot D_l(b).
\end{multline*}
I.e. identities \eqref{DEF:sist-chastn-proizv-*}, \eqref{DEF:sist-chastn-proizv-1} and \eqref{DEF:sist-chastn-proizv-2} hold, and this means that the family $\{D_k;\ k\in\N[m]\}$ is a system of partial derivatives on $A$ with the coefficients in $B$.

2. Conversely, if $\{D_k;\ k\in\N[m]\}$ is a system of partial derivatives on $A$ with coefficients in $B$, then we define the mapping $D:A\to B[m]$ by formula \eqref{D(a)=sum-D_k(a)/k!-tau^k}, and we get that, first, for each $a\in A$
$$
D(a^\bullet)^{(k)}=D_k(a^\bullet)=\eqref{DEF:sist-chastn-proizv-*}=D_k(a)^\bullet=\big(D(a)^{(k)}\big)^\bullet,
$$
second,
$$
D(1)^{(k)}=D_k(1)=\eqref{DEF:sist-chastn-proizv-1}=\begin{cases}1, & k=0
\\
0, & k\ne 0
\end{cases}=\eqref{1_(B[[d]])}=1_k\qquad\Longrightarrow\qquad D(1)=1,
$$
and, third, for each $a,b\in A$
\begin{multline*}
D(a\cdot b)^{(k)}=\eqref{homomorphizm-v-B[m]<->chast-proizv-v-B}=D_k(a\cdot b)=\eqref{DEF:sist-chastn-proizv-2}=\sum_{0\le l\le k}\begin{pmatrix}k \\ l\end{pmatrix}\cdot D_{k-l}(a)\cdot D_l(b)=\eqref{homomorphizm-v-B[m]<->chast-proizv-v-B}=\\=\sum_{0\le l\le k}\begin{pmatrix}k \\ l\end{pmatrix}\cdot  D(a)^{(k-l)}\cdot D(b)^{(l)}=\eqref{(xy)^(k)}=\big(D(a)\cdot D(b)\big)^{(k)} \qquad\Longrightarrow\qquad D(a\cdot b)=D(a)\cdot D(b).
\end{multline*}
I.e. $D:A\to B[m]$ is an involutive homomorphism.
\epr

\bex\label{EX:sist-chast-proizv-na-C-infty(M)} Let $M$ be a smooth manifold of dimension $d\in\N$, $\ph:U\to V$ a  chart, where $U\subseteq M$, $V\subseteq\R^d$ and $K\subseteq U$ is a compact set. then the system of operators
\beq\label{sist-chast-proizv-na-C-infty(M)}
D_k:{\mathcal E}(M)\to\mathcal{C}(K)\quad\Big|\quad
D_k(a)=\frac{\partial^{\abs{k}}(a\circ\ph^{-1})}{\partial t_1^{k_1}...\partial
t_d^{k_d}}\circ\ph
\eeq
is a system of partial derivatives on ${\mathcal E}(M)$. The corresponding homomorphism of algebras
$D:{\mathcal E}(M)\to\mathcal{C}(K)[m]$ is a system of restrictions of partial derivatives to the compact set $K$:
$$
(D(a))_k=D_k(a)\big|_K
$$
\eex

\paragraph{Partial derivatives as differential operators.}

Again, suppose $A$ is an involutive stereotype algebra, $B$ a $C^*$-algebra, and $\{D_k;\ k\in\N[m]\}$ a system of partial derivatives on $A$ with coefficients in $B$. Then the homomorphism $\ph=D_0:A\to B$ turns $B$ into a module over $A$, and thus for each operator $P:A\to B$ and for each element $a\in A$ a commutator $[P,a]:A\to B$ is defined.

\bprop For each system of partial derivatives $\{D_k;\ k\in\N[m]\}$ on $A$ with the coefficients in $B$ the commutator of operators $D_k$ with an arbitrary element $a\in A$ with respect to the homomorphism $D_0:A\to B$ acts by the formula
\beq\label{[D_k,a]}
[D_k,a]=\sum_{0\le l< k}\begin{pmatrix}k \\ l\end{pmatrix}\cdot D_{k-l}(a)\cdot D_l
\eeq
For $k=0$ this formula has the form
\beq\label{[D_0,a]}
[D_0,a]=0.
\eeq
\eprop
\bpr Identity \eqref{[D_0,a]} holds trivially, since $D_0$ is a homomorphism. And \eqref{[D_k,a]} is proved by the following chain
$$
[D_k,a](x)=D_k(a\cdot x)-\ph(a)\cdot D_k(x)=\sum_{0\le l\le k}\begin{pmatrix}k \\ l\end{pmatrix}\cdot D_{k-l}(a)\cdot D_l(x)-D_0(a)\cdot D_k(x)=
\sum_{0\le l< k}\begin{pmatrix}k \\ l\end{pmatrix}\cdot D_{k-l}(a)\cdot D_l(x).
$$
\epr

\btm\label{TH:differentsialnye-chast-proizv} For each system of partial derivatives $\{D_k;\ k\in\N[m]\}$ on an involutive stereotype algebras $A$ with the coefficients in a $C^*$-algebra $B$ the following conditions are equivalent:
 \bit{
 \item[(i)] the operators $\{D_k;\ k\in\N[m]\}$ are differential operators from $A$ into $B$ with respect to the homomorphism $D_0:A\to B$ with the orders, not greater than the module of their multi-indices:
     \beq\label{D_k-in-D^|k|(D_0)}
     D_k\in \Diff^{\abs{k}}(D_0),
     \eeq

 \item[(ii)] for each multi-index $k>0$ the values of the operator $D_k$ lie in the space $Z^{|k|}(D_0)$:
     \beq\label{D_k(A)-subseteq-Z^|k|(D_0)}
     D_k(A)\subseteq Z^{|k|}(D_0),\qquad k>0.
     \eeq

\item[(iii)] for each multi-index $k>0$ the values of the operator $D_k$ lie in the space $Z^1(D_0)=D_0(A)^!$ (in the commutant of the operator $D_0$):
     \beq\label{D_k(A)-subseteq-Z^1(D_0)}
     D_k(A)\subseteq Z^1(D_0)=D_0(A)^!,\qquad k>0.
     \eeq
  }\eit
\etm

 \bit{
\item[$\bullet$]\label{DEF:differentsialnye-chast-proizv} A system of partial derivatives $\{D_k;\ k\in\N[m]\}$ on $A$ with the coefficients in $B$ is said to be {\it differential}, if it satisfies the equivalent conditions  (i)-(iii) of Theorem \ref{TH:differentsialnye-chast-proizv}.

\item[$\bullet$] A homomorphism of involutive stereotype algebras $D:A\to B[m]$ is said to be {\it differential}, if the system of partial derivatives $\{D_k;\ k\in\N[m]\}$ defined by \eqref{homomorphizm-v-B[m]<->chast-proizv-v-B} is differential. The class of all differential homomorphisms will be denoted by $\DiffMor$.

 }\eit

\bpr Note first that the equivalence of (ii) and (iii) is a corollary of \eqref{Z^n(ph)=Z^(n+1)(ph)}. Hence we have to verify the equivalence of (i) and (ii).

1. (i)$\Longrightarrow$(ii). Suppose (i) holds. We shall prove (ii) by induction.

1) Suppose $|k|=1$. Then for any $a,a_1\in A$ we have:
$$
[D_k,a]=\eqref{[D_k,a]}=\sum_{0\le l< k}\begin{pmatrix}k \\ l\end{pmatrix}\cdot D_{k-l}(a)\cdot D_l=D_k(a)\cdot D_0
$$
$$
\Downarrow
$$
$$
0=[[D_k,a],a_1]=\eqref{[b-cdot-P, a]}=D_k(a)\cdot \kern-6pt\underbrace{[D_0,a_1]}_{\tiny\begin{matrix}
\phantom{\eqref{[D_0,a]}}
\ \text{\rotatebox{90}{$=$}}\ \eqref{[D_0,a]} \\ 0\end{matrix}}\kern-6pt +[D_k(a),D_0(a_1)]\cdot D_0
$$
$$
\Downarrow
$$
$$
[D_k(a),D_0(a_1)]=0
$$
This is true for each $a_1$, hence $D_k(a)\in Z^1(D_0)$. This in its turn is true for any $a\in A$, thus $D_k(A)\in Z^1(D_0)$.

2) Suppose we have already proved \eqref{D_k(A)-subseteq-Z^|k|(D_0)} for all $k$ such that $|k|\le n$:
     \beq\label{D_k(A)-subseteq-Z^|k|(D_0),n}
     D_k(A)\subseteq Z^{|k|}(D_0),\qquad 0<|k|\le n.
     \eeq
Then for $|k|=n+1$ we have: for any $a\in A$
$$
\underbrace{[D_k,a]}_{\tiny\begin{matrix}
\phantom{\eqref{[D^n,A]-subseteq-D^(n-1)}}
\ \text{\rotatebox{90}{$\owns$}}\ \eqref{[D^n,A]-subseteq-D^(n-1)} \\ \Diff^{\abs{k}-1}\\ \text{\rotatebox{90}{$=$}} \\ \Diff^n \end{matrix}}\kern-5pt=\eqref{[D_k,a]}=\sum_{0\le l< k}\begin{pmatrix}k \\ l\end{pmatrix}\cdot D_{k-l}(a)\cdot D_l=D_k(a)\cdot D_0+
\sum_{0<l< k}\begin{pmatrix}k \\ l\end{pmatrix}\cdot\underbrace{\kern-5pt\overbrace{D_{k-l}(a)}^{\tiny\begin{matrix} Z^{|k-l|} \\
\eqref{D_k(A)-subseteq-Z^|k|(D_0),n}
\ \text{\rotatebox{90}{$\in$}}\ \phantom{\eqref{D_k(A)-subseteq-Z^|k|(D_0),n}} \end{matrix}}\kern-5pt\cdot\kern-15pt \overbrace{D_l}^{\tiny\begin{matrix} \Diff^{\abs{l}} \\
\phantom{\eqref{D_k-in-D^|k|(D_0)}}
\ \text{\rotatebox{90}{$\in$}}\ \eqref{D_k-in-D^|k|(D_0)} \end{matrix}}\kern-15pt}_{\tiny\begin{matrix}
\phantom{\eqref{Z^q.D^p->D^(q+p-1)}}
\ \text{\rotatebox{90}{$\owns$}}\ \eqref{Z^q.D^p->D^(q+p-1)} \\ \Diff^{\abs{k-l}+|l|-1}\\ \text{\rotatebox{90}{$=$}} \\ \Diff^n \end{matrix}}
$$
\vglue-40pt
$$
\Downarrow
$$
$$
D_k(a)\cdot D_0\in \Diff^n
$$
$$
\phantom{\tiny \eqref{b.ph-in-D^n(ph)<=>b-in-Z^(n+1)(ph)}}\quad \Downarrow \quad {\tiny \eqref{b.ph-in-D^n(ph)<=>b-in-Z^(n+1)(ph)}}
$$
$$
D_k(a)\in Z^{n+1}=Z^{|k|}
$$

2. (i)$\Longleftarrow$(ii). Conversely, suppose (ii) holds. Then (i) is also proved by induction.

0) For $k=0$ the proposition \eqref{D_k-in-D^|k|(D_0)} is always true, since the homomorphism $\ph=D_0$ is a differential operator of zero order by formula \eqref{[ph,a]=0}.

1) Suppose we have already proved \eqref{D_k-in-D^|k|(D_0)} for all $k$ such that $|k|\le n$:
     \beq\label{D_k-in-D^|k|(D_0),n}
     D_k\in \Diff^{\abs{k}}(D_0),\qquad |k|\le n.
     \eeq
Then for $|k|=n+1$ we have:
$$
\forall  a\in A\qquad
[D_k,a]=\eqref{[D_k,a]}=
\sum_{0\le l< k}\begin{pmatrix}k \\ l\end{pmatrix}\cdot\underbrace{\kern-5pt\overbrace{D_{k-l}(a)}^{\tiny\begin{matrix} Z^{|k-l|} \\
\eqref{D_k(A)-subseteq-Z^|k|(D_0)}
\ \text{\rotatebox{90}{$\in$}}\ \phantom{\eqref{D_k(A)-subseteq-Z^|k|(D_0)}} \end{matrix}}\kern-5pt\cdot\kern-15pt \overbrace{D_l}^{\tiny\begin{matrix} \Diff^{\abs{l}} \\
\phantom{\eqref{D_k-in-D^|k|(D_0),n}}
\ \text{\rotatebox{90}{$\in$}}\ \eqref{D_k-in-D^|k|(D_0),n} \end{matrix}}\kern-15pt}_{\tiny\begin{matrix}
\phantom{\eqref{Z^q.D^p->D^(q+p-1)}}
\ \text{\rotatebox{90}{$\owns$}}\ \eqref{Z^q.D^p->D^(q+p-1)} \\ \Diff^{\abs{k-l}+|l|-1}\\ \text{\rotatebox{90}{$=$}} \\ \Diff^n \end{matrix}},
$$
$$
\Downarrow
$$
$$
\forall  a\in A\qquad [D_k,a]\in \Diff^n,
$$
$$
\Downarrow
$$
$$
D_k\in \Diff^{n+1}.
$$
\epr

\subsection{Smooth envelopes}

\paragraph{Definition of the smooth envelope and functoriality.}

 \bit{
\item A {\it smooth envelope}\label{DEF:C^infty-obolochka} $\env_{\mathcal E}A:A\to\Env_{\mathcal E}A$ of an involutive stereotype algebra $A$ is its envelope in the class $\DEpi$ of dense epimorphisms in the category $\InvSteAlg$ of involutive stereotype algebras with respect to the class $\DiffMor$ of all differential homomorphisms into $C^*$-algebras $B[m]$ with the joined self-adjoint nilpotent elements:
  $$
  \Env_{\mathcal E} A=\Env_{\DiffMor}^{\DEpi}A
  $$
 }\eit

In detail, a {\it smooth extension} of an involutive stereotype algebra $A$ is defined as a dense epimorphism $\sigma:A\to A'$ of involutive stereotype algebras such that for each  $C^*$-algebra $B$, each multi-index $m\in\N^d$ and each differential involutive homomorphism $\ph:A\to B[m]$ there is a unique homomorphism of involutive stereotype algebras $\ph':A'\to B[m]$, such that the following diagram is commutative:
$$
 \xymatrix @R=2pc @C=1.2pc
 {
  A\ar[rr]^{\sigma}\ar[dr]_{\ph} & & A'\ar@{-->}[dl]^{\ph'} \\
  & B[m] &
 }
$$
And a {\it smooth envelope} of the algebra $A$ is a smooth extension $\rho:A\to \Env_{\mathcal E} A$, such that for each smooth extension $\sigma:A\to A'$ there is a unique homomorphism of involutive stereotype algebras $\upsilon:A'\to \Env_{\mathcal E} A$, such that the following diagram is commutative:
$$
 \xymatrix @R=2pc @C=1.2pc
 {
  & A\ar[ld]_{\sigma}\ar[rd]^{\rho} &   \\
  A'\ar@{-->}[rr]_{\upsilon} &  & \Env_{\mathcal E} A
 }
$$

\btm\label{TH:Env_C^infty-reg-obolochka} The smooth envelope $\Env_{\mathcal E}$ is regular and is coherent with the projective tensor product\footnote{Definitions on pages \pageref{DEF:reg-obolochka} and \pageref{DEF:obolochka-soglasovana-s-tenz-proizv}.} $\circledast$ in $\InvSteAlg$.
\etm
\bpr We use here the same tricks as in Theorem \ref{TH:Env_C-reg-obolochka}, and in the second part of the proposition -- the coherence with the tensor product -- the proof is exactly the same, while in the first part the differences appear in two last points.
\bit{
\item[R.4:] For each stereotype algebra $A$ there always exists a differential morphism $\ph:A\to B[m]$ (for instance one can take $B=0$, $m=0$ and $\ph=0$). This means (see definition on page \pageref{DEF:goes-from}), that the class $\DiffMor$ goes from $\InvSteAlg$. let us show that $\DiffMor$ is a right ideal in $\InvSteAlg$. Take  $D\in\DiffMor$, $D:B\to C[m]$, and $\sigma\in\Mor$, $\sigma:A\to B$. Then $D\circ\sigma:A\to C[m]$. For any $k>0$ we have
\begin{multline*}
(D\circ\sigma)_k(A)=D(\sigma(A))^{(k)}\subseteq D(B)^{(k)}=D_k(B)\overset{\eqref{D_k(A)-subseteq-Z^1(D_0)}}{\subseteq} D_0(B)^!\kern-30pt\overset{\scriptsize\begin{matrix}B\supseteq\sigma(A)\\ \Downarrow \\ D_0(B)\supseteq D_0(\sigma(A))\\ \Downarrow\end{matrix}}{\subseteq}\kern-30pt D_0(\sigma(A))^!=\\=
\Big(D(\sigma(A))^{(0)}\Big)^!=\Big((D\circ\sigma)(A)^{(0)}\Big)^!=
(D\circ\sigma)_0(A)^!
\end{multline*}
Thus, the morphism $D\circ\sigma:A\to C[m]$ satisfies the condition (iii) in Theorem \ref{TH:differentsialnye-chast-proizv}, i.e. $D\circ\sigma$ is a differential morphism.

\item[R.5:] It remains to check that the class $\DEpi$ pushes the class $\DiffMor$. Take the morphisms $D:B\to C[m]$ and $\sigma:A\to B$ such that $\sigma\in\DEpi$ and $D\circ\sigma\in\DiffMor$. Then for each $k>0$ we have:
$$
D_k(\sigma(A))=(D\circ\sigma)_k(A)\subseteq (D\circ\sigma)_0(A)^!=\Big((D\circ\sigma)(A)^{(0)}\Big)^!=\Big(D(\sigma(A))^{(0)}\Big)^! \kern-50pt\overset{\scriptsize\begin{matrix}\overline{\sigma(A)}=B\\ \Downarrow \\ D_0(\sigma(A))^!= D_0\Big(\overline{\sigma(A)}\Big)^!=D_0(B)^!\\ \Downarrow\end{matrix}}{\subseteq}\kern-50pt D_0(B)^!
$$
$$
\Downarrow
$$
$$
D_k(B)=D_k\Big(\overline{\sigma(A)}\Big)\subseteq\overline{D_k(\sigma(A))}\subseteq D_0(B)^!.
$$
Thus, the morphism $D:B\to C[m]$ satisfies the condition (iii) of Theorem \ref{TH:differentsialnye-chast-proizv}, i.e. $D$ is a differential morphism.
}\eit
\epr

\bcor\label{COR:Env_E-idemp-funktor}
The smooth envelope can be defined as an idempotent covariant functor from $\InvSteAlg$ into $\InvSteAlg$: there exist
\bit{
\item[1)] a map $A\mapsto (\Env_{\mathcal E}A,\env_{\mathcal E}A)$, that assigns to each involutive stereotype algebra $A$ an involutive stereotype algebra $\Env_{\mathcal E}A$ and a morphism of involutive stereotype algebras $\env_{\mathcal E}A:A\to\Env_{\mathcal E}A$, which is a smooth envelope of $A$, and

\item[2)] a map $\ph\mapsto \Env_{\mathcal E}(\ph)$, thatassigns to each morphism of involutive stereotype algebras $\ph:A\to B$ a morphism of involutive stereotype algebras $\Env_{\mathcal E}(\ph):\Env_{\mathcal E}A\to \Env_{\mathcal E}B$, such that the following diagram is commutative:
\beq\label{DIAGR:funktorialnost-env_E}
\xymatrix @R=2.pc @C=5.0pc % @M=14pt
{
A\ar[d]^{\ph}\ar[r]^{\env_{\mathcal E}A} & \Env_{\mathcal E}A\ar@{-->}[d]^{\Env_{\mathcal E}(\ph)} \\
B\ar[r]^{\env_{\mathcal E}B} & \Env_{\mathcal E}B \\
}
\eeq
}\eit
and the following identities hold:
\beq\label{funktorialnost-env_E-v-Ste^circledast}
\Env_{\mathcal E}(1_A)=1_{\Env_{\mathcal E}A},\quad \Env_{\mathcal E}(\beta\circ\alpha)=\Env_{\mathcal E}(\beta)\circ \Env_{\mathcal E}(\alpha),
\eeq
\beq\label{funktorialnost-env_E-v-Ste^circledast-2}
\Env_{\mathcal E}(\Env_{\mathcal E}A)=\Env_{\mathcal E}A,\quad \env_{\mathcal E}\Env_{\mathcal E}A=1_{\Env_{\mathcal E}A},
\eeq
\beq\label{Env_E(C)=C}
\Env_{\mathcal E}\C=\C
\eeq
\ecor

\btm[connection with the continuous envelope]\label{7:Env_infty->Env} The smooth envelope is mapped into the ccontinuous envelope: for any involutive stereotype algebra $A$ there exists a unique morphism of involutive stereotype algebras $\zeta_A:\Env_{\mathcal E}A\to\Env_{\mathcal C} A$ such that the following diagram is commutative
\beq\label{Env_infty->Env}
\begin{diagram}
\node[2]{A} \arrow{sw,t}{\env_{\mathcal E} A} \arrow{se,t}{\env_{\mathcal C} A}\\
\node{\Env_{\mathcal E} A}\arrow[2]{e,b,--}{\zeta_A}   \node[2]{\Env_{\mathcal C} A}
\end{diagram}
\eeq
This morphism $\zeta_A$ is always a dense epimorphism.
\etm
\bpr
This follows from \cite[(3.10)]{Akbarov-env}.
\epr

\paragraph{A net of differential quotient mappings.}
The construction of smooth envelope can be described a bit more visual in the following way. A neighbourhood of zero $U$ in an involutive stereotype algebra $A$ will be called {\it differential}, if it is a reverse image of the unit ball under some differential homomorphism $D:A\to B[m]$ into a $C^*$-algebra $B$ with joined adjoined nilpotent elements:
$$
U=\{x\in A:\ \norm{D(x)}\le 1\}
$$
(the norm $\norm{\cdot}$ is defined in \eqref{norm(x)=sum_(k-le-n)norm(x_k)_B}). For each differential neighbourhood of zero $U$ in $A$ its kernel
$$
\Ker U=\bigcap_{\lambda>0}\lambda\cdot U
$$
coincides with the kernel of the homomorphism $D$, hence it is a closed ideal in $A$. Consider the quotient algebra $A/\Ker U$ and endow it with the norm, for which the set of classes $U+\Ker U$ is a unit ball. This algebra $A/\Ker U$ can be treated as a subalgebra in $B[m]$ with the norm induced from this space. Its completion
\beq\label{A/U=(A/Ker U)^blacktriangledown}
A/U=(A/\Ker U)^\blacktriangledown
\eeq
is a Banach algebra (and we can consider it as a closed subalgebra in $B[m]$). We call $A/U$ a {\it quotient algebra of the algebra $A$ by a differential neighbourhood of zero $U$}, and the mapping
$$
\pi_U:A\to A/U
$$
we call the {\it quotient mapping} of the algebra $A$ by the differential neighbourhoods of zero $U$, or a {\it differential quotient mapping} of the algebra $A$.

\blm
For each differential homomorphism $\ph:A\to B[n]$ there is a unique homomorphism $\ph_U:A/U\to B[n]$, such that the following diagram is commutative:
\beq\label{ph=ph_U-circ-pi_U}
 \xymatrix @R=2pc @C=1.2pc
 {
  A\ar[rr]^{\ph}\ar[dr]_{\pi_U} & & B[n]  \\
  & A/U\ar@{-->}[ur]_{\ph_U} &
  }
\eeq
\elm

\blm
If $U$ and $U'$ are two differential neighbourhoods of zero in $A$, and $U\supseteq U'$, then there is a unique homomorphism $\varkappa^{U'}_U:A/U\gets A/U'$, such that the following diagram is commutative:
\beq\label{varkappa^U'_U}
 \xymatrix @R=2pc @C=1.2pc
 {
  & A\ar[ld]_{\pi_U}\ar[rd]^{\pi_{U'}} &   \\
  A/U &  & A/ U'\ar@{-->}[ll]^{\varkappa^{U'}_U}
 }
\eeq
\elm

\blm\label{LM:diff-okr-nulya-uporyadocheny} Any two differential neighbourhoods of zero $U$ and $U'$ in $A$ contain a common differential neighbourhood of zero $U''$:
$$
U\cap U'\supseteq U''.
$$
\elm
\bpr
Indeed, let $D:A\to B[m]$ and $D':A\to B'[m']$ be differential homomorphisms generating $U$ and $U'$.
$$
U=\{x\in A:\ \norm{D(x)}\le 1\},\qquad U'=\{x\in A:\ \norm{D'(x)}\le 1\}.
$$
Consider the homomorphism
$$
\ph:B[m]\oplus B'[m']\to (B\oplus B')[m\oplus m'],
$$
described in \eqref{A[m]-oplus-B[n]-to-(A-oplus-B)[m-oplus-n]}. The mapping
$$
D'':A\to (B\oplus B')[m\oplus m']\quad\Big|\quad D''(x)=\ph(D(x)\oplus D'(x)), \quad x\in A,
$$
is a homomorphism, and if we put
$$
U''=\{x\in A:\ \norm{D''(x)}\le 1\},
$$
then for each  $x\in U''$ we have
$$
1\ge \norm{D''(x)}=\norm{\ph(D(x)\oplus D'(x))}\ge \eqref{1-le-norm(ph)-le-2}\ge \norm{D(x)\oplus D'(x)}=
\max\{\norm{D(x)},\norm{D'(x)}\}
$$
i.e. $x\in U$ and $x\in U'$.
\epr

\btm
The system $\pi_U:A\to A/U$ of differential quotient mappings forms a net of epimorphisms\footnote{See definition on page \pageref{DEF:set-epimorf}.} in the category $\InvSteAlg$ of involutive stereotype algebras, i.e. has the following properties:
  \bit{
\item[(a)] each algebra $A$ has at least one differential neighbourhood of zero $U$, and the set of all differential neighbourhoods of zero in $A$ is directed by the pre-order
$$
U\le U'\quad\Longleftrightarrow\quad U\supseteq U',
$$

\item[(b)] for each algebra $A$ the system of morphisms $\varkappa_U^{U'}$ from \eqref{varkappa^U'_U}
is covariant, i.e. for each three neighbourhoods of zero $U\supseteq U'\supseteq U''$ the following diagram is commutative:
$$
 
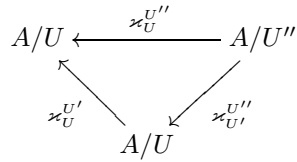
\xymatrix @R=2pc @C=1.2pc
 {
  A/U &  & A/ U''\ar[ll]_{\varkappa^{U''}_U}\ar[dl]^{\varkappa^{U''}_{U'}}\\
  & A/U \ar[ul]^{\varkappa^{U'}_U} &
 }
$$
and this system $\varkappa_U^{U'}$ has a projective limit in $\InvSteAlg$;

\item[(c)] for each morphism $\alpha:A\gets A'$ in $\InvSteAlg$ and for each differential neighbourhood of zero $U$ in $A$ there is a differential neighbourhood of zero $U'$ in $A'$ and a morphism $\alpha_U^{U'}:A/U\gets A'/U'$ such that the following diagram is commutative:
 \beq
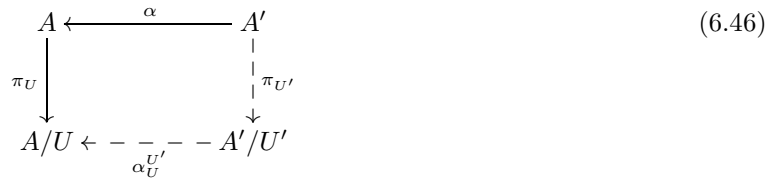
\label{DIAGR:set-differ} \xymatrix @R=2.5pc @C=4.0pc {
 A\ar[d]_{\pi_U} & A'\ar@{-->}[d]^{\pi_{U'}}\ar[l]_{\alpha} \\
A/U & A'/U'\ar@{-->}[l]^{\alpha_U^{U'}}
 } \eeq
 }\eit
\etm

By condition (b) of this theorem, there exists a projective limit $\leftlim_{0\gets U'} A/ U'$ of the system  $\varkappa_U^{U'}$. As a corollary, there exists a unique arrow $\pi:A\to \leftlim_{0\gets U'} A/ U'$ in $\InvSteAlg$, such that the following diagrams are commutative:
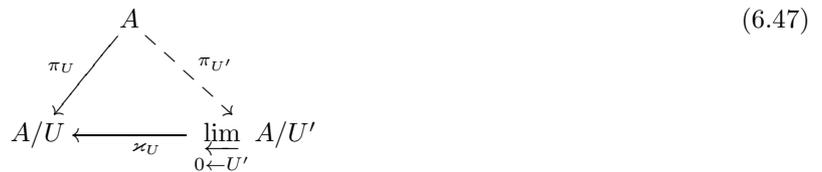
\beq\label{A-to-leftlim-A/U}
 \xymatrix @R=2pc @C=1.2pc
 {
  & A\ar[ld]_{\pi_U}\ar@{-->}[rd]^{\pi_{U'}} &   \\
  A/U &  & \leftlim_{0\gets U'} A/ U'\ar[ll]^{\varkappa_U}
 }
\eeq
The set of values $\pi(A)$ of the mapping $\pi$ is (an involutive subalgebra and) a subspace in the stereotype space  $\leftlim_{0\gets U'} A/ U'$. Hence it generates an immediate subspace in $\leftlim_{0\gets U'} A/ U'$, or, an envelope $\Env\pi(A)$ \cite{Akbarov-env}, i.e. a maximal stereotype subspace in $\leftlim_{0\gets U'} A/ U'$ having  $\pi(A)$ as its dense subspace. Denote by $\rho:A\to \Env\pi(A)$ the lifting of the morphism $\pi$ to $\Env\pi(A)$.

\btm The morphism $\rho:A\to \Env\pi(A)$ is a smooth envelope of the algebra $A$:
$$
\Env\pi(A)=\Env_{\mathcal E}A.
$$
\etm
\bpr
Here we ahve to follow the rpoof of Theorem 3.42 in \cite{Akbarov-env}: a net of epimorphisms $\mathcal N$ constructed there differs from the system $\pi_U:A\to A/U$ of differential quotient map that we constructed here in the detail that elements of $\mathcal N$ are finite sets of maps $\pi_U:A\to A/U$ (this is important in Theorem 3.42 in \cite{Akbarov-env} for ${\mathcal N}^X$ being directed with respect to the preorder in the class of epimorphisms, but in our case the system of quotient maps $\pi_U:A\to A/U$ is alreday directed by Lemma  \ref{LM:diff-okr-nulya-uporyadocheny}). Certainly, the local limits of these nets (i.e. the projective limits for each given $X$) coincide. As a corollary, all the conclusions for $\mathcal N$ are true for the net $\pi_U:A\to A/U$, in particular, the conclusion that the envelope with respect to these nets coincide with the envelope with respect to the initial class of morphisms $\varPhi=\DiffMor$.
\epr

\brem
As in the case of the continuous envelope, the smooth envelope $\rho:A\to\Env_{\mathcal E}A$ is the composition of the elements $\red_\infty$ and $\coim_\infty$ of the nodal decomposition of the morphism $\pi:A\to\leftlim_{0\gets U'} A/ U'$ in the category $\tt Ste$ of stereotype spaces (not algebras!):
  \beq\label{C^infty-envelope=im_infty-lim N_X}
\red_\infty\pi\circ\coim_\infty\pi=\env_{\mathcal E} A,
  \eeq
Visually this can be illustrated by the diagram
  \beq\label{DIAGR:E-envelope=im_infty-lim N_X}
\xymatrix @R=3.pc @C=6.0pc % @M=14pt
{
A\ar[d]_{\coim_\infty\pi}\ar[r]^{\pi=\leftlim_{0\gets U'}\pi_{U'}} & \leftlim_{0\gets U'} A/U' &   \\
\Coim_\infty\pi\ar[r]_{\red_\infty\pi} &  \Im_\infty \pi \ar[u]_{\im_\infty\pi}\ar@{=}[r] & \Env_{\mathcal C}A
}
 \eeq
The algebra $\Env_{\mathcal E} A$ can be understood as the envelope (in the sense of \eqref{DEF:Env^XM}) of the set of values of $\pi$ in the stereotype space $\lim_{0\gets U'} A/ U'$:
  \beq\label{Env_EA=Env-pi(A)}
\Env_{\mathcal E} A=\Env\pi(A).
  \eeq
\erem

\subsection{Smooth algebras}

We call an involutive stereotype algebra $A$ a {\it smooth algebra}, if it coincides with its smooth envelope, or, equivalently, if its smooth envelope $\env_{\mathcal E}A: A\to \Env_{\mathcal E}A$ is an isomorphism in the category  $\InvSteAlg$ of involutive stereotype algebras. The class of all smooth algebras is denoted by ${\mathcal E}\text{-}{\tt Alg}$. It forms a full subcategory in $\InvSteAlg$.

\paragraph{Smooth tensor product of involutive stereotype algebras.}

Let $\Env_{\mathcal E}$ be the functor of smooth envelope, defined in Corollary \ref{COR:Env_C^infty-idemp-funktor}.
For each two involutive stereotype algebras $A$ and $B$ we call its {\it smooth tensor product} the algebra
\beq\label{E/circledast}
A\overset{\mathcal E}{\circledast} B=\Env_{\mathcal E}(A\circledast B)
\eeq

\btm\label{TH:C^infty-circledast->odot} For each two involutive stereotype algebras $A$ and $B$ there is a unique linear continuous mapping $\eta^\infty_{A,B}:A\overset{\mathcal E}{\circledast}B\to A\odot B$ such that the following diagraam is commutative,
    \beq\label{DIAGR:eta^infty_(A,B)}
\xymatrix @R=2.pc @C=5.0pc % @M=14pt
{
A\circledast B\ar[dr]_{\env_{\mathcal E}{A\circledast B}\quad}\ar[rr]^{@_{A,B}} & & A\odot B \\
& A\overset{\mathcal E}{\circledast}B\ar[ur]_{\eta^\infty_{A,B}} & \\
}
\eeq
and the system of mappings $\eta^\infty_{A,B}:A\overset{\mathcal E}{\circledast}B\to A\odot B$ is a natural transformation of the functor $(A,B)\in{\mathcal E}\text{-}{\tt Alg}^2\mapsto A\overset{\mathcal E}{\circledast}B\in{\mathcal E}\text{-}{\tt Alg}$ into the functor $(A,B)\in{\mathcal E}\text{-}{\tt Alg}^2\mapsto A\odot B\in {\tt Ste}$.
\etm
\bpr Consider the diagram
$$
\xymatrix @R=3.pc @C=8.0pc % @M=14pt
{
A\circledast B\ar[d]_{\env_{\mathcal E}{A\circledast B}}\ar[dr]_{\env_{\mathcal C}{A\circledast B}\quad}\ar[r]^{@_{A,B}} &  A\odot B \\
A\overset{\mathcal E}{\circledast}B\ar[r]_{\zeta_{A\circledast B}} & A\overset{\mathcal C}{\circledast}B\ar[u]_{\eta_{A,B}}  \\
}
$$
Here the left lower triangle is Diagram \eqref{Env_infty->Env} for the algebra $A\circledast B$, and the right upper triangle is Diagram \eqref{DIAGR:eta_(A,B)}. The morphism
$$
\eta^\infty_{A,B}=\eta_{A,B}\circ\alpha_{A,B}
$$
is the one we need.
 \epr

\paragraph{Smooth tensor product of smooth algebras.}

From Theorem \ref{TH:Env_C^infty-reg-obolochka} and \ref{TH:sushestvovanie-tenz-proizv-v-L} it follows

\btm\label{TH:E-obolochka=monoidalnyi-funktor} Formula \eqref{C/circledast} defines in ${\mathcal E}\text{-}{\tt Alg}$ a tensor product, which turns ${\mathcal E}\text{-}{\tt Alg}$ into a monoidal category, and the functor  $\Env_{\mathcal E}$ is a monoidal functor from the monoidal category $(\InvSteAlg,\circledast)$ of involutive stereotype algebras into the monoidal category $({\mathcal E}\text{-}{\tt Alg},\overset{\mathcal E}{\circledast})$ of smooth algebras. The corresponding morphism of bifunctors
 $$
\Big((A,B)\mapsto \Env_{\mathcal E}(A)\overset{\mathcal E}{\circledast} \Env_{\mathcal E}(B)\Big)\overset{E^{\circledast}}{\rightarrowtail} \Big((A,B)\mapsto \Env_{\mathcal E}(A\circledast B)\Big)
 $$
is defined by the formula
$$
E^\circledast_{A,B}=\Env_{\mathcal E}(\env_{\mathcal E}A\circledast \env_{\mathcal E}B)^{-1}:\Env_{\mathcal E}(A)\overset{\mathcal E}{\circledast} \Env_{\mathcal E}(B)=\Env_{\mathcal E}(\Env_{\mathcal E}(A)\circledast \Env_{\mathcal E}(B))\to \Env_{\mathcal E}(A\circledast B),
$$
and the local identity
$$
E^{\C}=1_{\C}:\C\to\C=\Env_{\mathcal E}(\C).
$$
is the morphism in ${\mathcal C}\text{-}{\tt Alg}$, that turns the identity object $\C$ in ${\mathcal C}\text{-}{\tt Alg}$ into the image $\Env_{\mathcal E}(\C)$ of the identity object $\C$ in $\InvSteAlg$.
\etm

To each pair of elements $a\in A$, $b\in B$ one can assign the elementary tensor
\beq\label{DEF:a-overset-C^infty-circledast-b}
a\overset{\mathcal E}{\circledast}b=\env_{\mathcal E}(a\circledast b)
\eeq

\blm\label{LM:polnota-a-overset-C^infty-circledast-b}
The elementary tensors $a\overset{\mathcal E}{\circledast}b$, $a\in A$, $b\in B$, are total in $A\overset{\mathcal E}{\circledast} B$ and the mapping $\eta_{A,B}$ turn them into the elementary tensors $a\odot b$:
\beq\label{eta(a-overset-C^infty-circledast-b)=a-odot-b}
\eta_{A,B}(a\overset{\mathcal E}{\circledast}b)=a\odot b.
\eeq
\elm
\bpr
The tensors $a\circledast b$ are total in $A\circledast B$, and the set of values of $\env_{\mathcal E}$ is dense in  $A\overset{\mathcal E}{\circledast} B$. The identity \eqref{eta(a-overset-C^infty-circledast-b)=a-odot-b} follows from the diagram \eqref{DIAGR:eta^infty_(A,B)}.
\epr

\paragraph{Action of smooth envelope on bialgebras.}

The following three propositions are analogues of Teorems \ref{LM:koalg-v-CAlg->koalg-v-odot},  \ref{TH:C-obolochka-sohranyaet-Hopfov} and \ref{TH:C-obolochka-sohranyaet-inv-Hopfov}, and are proved similarly.

\blm\label{LM:koalg-v-C^inftyAlg->koalg-v-odot}
If $A$ is a coalgebra in the monoidal category $({\mathcal E}\text{-}{\tt Alg},\overset{\mathcal E}{\circledast})$  of smooth algebras with the structure morphisms
$$
\varkappa:A\to A\overset{\mathcal E}{\circledast} A,\qquad \e:A\to\C,
$$
then $A$ is a coalgebra in the moniodal category $({\tt Ste},\odot)$ of stereotype spaces with the structure morphisms $$
\lambda=\eta_{A,A}\circ\varkappa:A\to A\odot A,\qquad \e:A\to\C.
$$
\elm

\btm\label{TH:C^infty-obolochka-sohranyaet-Hopfov}
Let $H$ be a bialgebra in the category $({\tt Ste},\circledast)$ of stereotype spaces, or, equivalently, a coalgebra in the category ${\tt Ste}^{\circledast}$ of stereotype algebras with the comultiplication $\varkappa$ and the counit $\e$. Then
 \bit{
\item[(i)] the smooth envelope $\Env_{\mathcal E}H$ is a coalgebra in the monoidal category $({\mathcal C}\text{-}{\tt Alg},\overset{\mathcal E}{\circledast})$ of smooth algebras with the comultiplication and the counit
    \beq\label{varkappa^E^infty,e^E^infty}
    \varkappa_{\Env_{\mathcal E}}=\Env_{\mathcal E}(\env_{\mathcal E}H\circledast \env_{\mathcal E}H)\circ \Env_{\mathcal E}(\varkappa),\qquad \e_{\Env_{\mathcal E}}=\Env_{\mathcal E}(\e),
    \eeq

\item[(ii)] the smooth envelope $\Env_{\mathcal E}H$ is a coalgebra in the monoidal category $({\tt Ste},\odot)$ of stereotype spaces with the comultiplication and the counit
    \beq\label{varkappa^odot,e^odot-infty}
    \varkappa_\odot=\eta_{\Env_{\mathcal E}H,\Env_{\mathcal E}H}\circ \Env_{\mathcal E}(\env_{\mathcal E}H\circledast \env_{\mathcal E}H)\circ \Env_{\mathcal E}(\varkappa)=
    \eta_{\Env_{\mathcal E}H,\Env_{\mathcal E}H}\circ \varkappa_{\Env_{\mathcal E}},\qquad \e_\odot=\Env_{\mathcal E}(\e),
    \eeq

\item[(iii)] the morphism $(\env_{\mathcal E}H)^\star:H^\star\gets \Env_{\mathcal E}H^\star$, dual to the morphism of the envelope  $\env_{\mathcal E}H:H\to \Env_{\mathcal E}H$, is a morphism of stereotype algebras, if  $\Env_{\mathcal E}H^\star$ is considered as an algebra with the multiplication and the unit, dual to \eqref{varkappa^odot,e^odot-infty}, and $H^\star$ the algebra with the multiplication and the unit
    $$
    \varkappa^\star\circ @_{H^\star,H^\star},\qquad \e^\star.
    $$

 }\eit
\etm

\btm\label{TH:C^infty-obolochka-sohranyaet-inv-Hopfov}
Let $H$ be an involutive Hopf algebra in the category $({\tt Ste},\circledast)$ of stereotype spaces. Then
 \bit{
\item[(i)] the smooth envelope $\Env_{\mathcal E}H$, as a coalgebra in the monoidal categories $({\mathcal E}\text{-}{\tt Alg},\overset{\mathcal E}{\circledast})$ and $({\tt Ste},\odot)$, has interconsistent antipode $\Env_{\mathcal E}(\sigma)$ and involution $\Env_{\mathcal E}(\bullet)$, defined by the diagrams in the category $\tt{Ste}$
    \beq\label{C^infty-obolochka-sohranyaet-inv-Hopfov}
    \xymatrix @R=2.pc @C=4.0pc % @M=14pt
{
H\ar[d]_{\sigma}\ar[r]^{\env_{\mathcal E}H} & \Env_{\mathcal E}H\ar@{-->}[d]^{\Env_{\mathcal E}(\sigma)} \\
H\ar[r]^{\env_{\mathcal E}H} & \Env_{\mathcal E}H
}
\qquad
    \xymatrix @R=2.pc @C=4.0pc % @M=14pt
{
H\ar[d]_{\bullet}\ar[r]^{\env_{\mathcal E}H} & \Env_{\mathcal E}H\ar@{-->}[d]^{\Env_{\mathcal E}(\bullet)} \\
H\ar[r]^{\env_{\mathcal E}H} & \Env_{\mathcal E}H
}
\eeq

\item[(ii)] the morphism $(\env_{\mathcal E}H)^\star:H^\star\gets \Env_{\mathcal E}H^\star$, dual to the morphism of envelope $\env_{\mathcal E}H:H\to \Env_{\mathcal E}H$, is an involutive homomorphism of stereotype algebras over $\circledast$, if $H^\star$ and $\Env_{\mathcal E}H^\star$ are endowed with the structure of dual involutive algebras to the involutive coalgebras with the antipode $H$ and $\Env_{\mathcal E}H$ by Property $4^\circ$ on page \pageref{4^0:inv-v-sopryazh-alg}.
 }\eit
\etm

\paragraph{Smooth tensor product with ${\mathcal E}(M)$.}

Let $X$ be a stereotype space, and $M$ a smooth (locally euclidean) manifold. Consider the algebra ${\mathcal E}(M)$ of smooth functions on $M$ and the space ${\mathcal E}(M,X)$ of smooth mappings on $M$ with values in $X$. We endow  ${\mathcal E}(M)$ and ${\mathcal E}(M,X)$ by the standard topology of uniform convergence on compact sets by any partial derivative
$$
u_i\overset{{\mathcal E}(M,X)}{\longrightarrow} 0\quad\Longleftrightarrow\quad \forall U\subseteq M\quad \forall k\in\N^d \qquad u_i^{(k)}|_U\overset{{\mathcal C}(U,X)}{\longrightarrow} 0,\qquad u\in {\mathcal E}(M,X),\quad t\in M.
$$
(where $u^{(k)}$ is the partial derivative with respect to a chart on an open set $U\subseteq M$) and the pointwise operations:
$$
(\lambda\cdot u)(t)=\lambda\cdot u(t)\qquad (u+v)(t)=u(t)+v(t),\qquad u,v\in {\mathcal C}(M,X),\quad \lambda\in\C,\quad t\in M.
$$

From \cite[Theorem 8.9]{Akbarov} we have

\bprop The following identity holds:
\beq\label{E(M,X)-cong-E(M)-odot-X}
{\mathcal E}(M,X)\cong {\mathcal E}(M)\odot X
\eeq
\eprop

Further we shall be interested in the case when $A$ is a smooth (and thus, a stereotype) algebra. We endow the space  ${\mathcal E}(M,A)$ with the structure of stereotype algebra with the pointwise operations
$$
(u\cdot v)(t)=u(t)\cdot v(t),\qquad u,v\in {\mathcal C}(M,A),\quad t\in M.
$$
From \eqref{E(M,X)-cong-E(M)-odot-X} we see that ${\mathcal E}(M,A)$ is a stereotype $A$-module.

\btm\label{TH:C^infty(M)-circledast-A->C(M,A)} For each smooth algebra $A$ and for each smooth manifold $M$ the natural mapping
\beq\label{C^infty(M)-circledast-A->C(M,A)}
\iota:{\mathcal E}(M)\circledast A\to {\mathcal E}(M,A) \quad\Big|\quad \iota(u\circledast a)(t)=u(t)\cdot a,\quad u\in {\mathcal E}(M),\ a\in A,\ t\in M,
\eeq
is a smooth envelope and generates an isomorphism of stereotype algebras:
\beq\label{E(M)-circledast-A=E(M,A)}
{\mathcal E}(M)\overset{\mathcal E}{\circledast} A\cong{\mathcal E}(M,A).
\eeq
\etm

We split the proof into 5 lemmas.

\blm\label{LM:iota-in-DEpi-E(M)}
The mapping $\iota:{\mathcal E}(M)\circledast A\to {\mathcal E}(M,A)$ is a dense epimorphism.
\elm
\bpr
This is proved by analogy with Lemma \ref{LM:iota-in-DEpi-C(M)}.
\epr

\blm\label{LM:J^n_E(M)E(M)-circledast-A-cong-J^n_E(M)E(M,A)}
The modules ${\mathcal E}(M)\circledast A$ and ${\mathcal E}(M,A)$ over the algebra ${\mathcal E}(M)$ have isomorphic jet bundles:
\beq\label{J^n_E(M)E(M)-circledast-A-cong-J^n_E(M)E(M,A)}
\Jet^n_{{\mathcal E}(M)}{\mathcal E}(M)\circledast A\cong \Jet^n_{{\mathcal E}(M)}{\mathcal E}(M,A),\qquad n\in\N
\eeq
\elm
\bpr
In each point $t\in M$ the ideal $I_t^{n+1}$ has finite co-dimension in ${\mathcal E}(M)$, hence we can use Lemma \ref{PROP:[(X-circledast-Z)/(Y-circledast-Z)]^triangledown-cong-[(X-odot-Z)/(Y-odot-Z)]^triangledown}:
\begin{multline*}
\Jet^n_{{\mathcal E}(M)}{\mathcal E}(M)\circledast A=[({\mathcal E}(M)\circledast A)/ (I_t^{n+1}\circledast A)]^\vartriangle=\eqref{[(X-circledast-Z)/(Y-circledast-Z)]^triangledown-cong-[(X-odot-Z)/(Y-odot-Z)]^triangledown}=\\=
[({\mathcal E}(M)\odot A)/ (I_t^{n+1}\odot A)]^\vartriangle=\Jet^n_{{\mathcal E}(M)}{\mathcal E}(M)\odot A=\eqref{E(M,X)-cong-E(M)-odot-X}=
\Jet^n_{{\mathcal E}(M)}{\mathcal E}(M,A)
\end{multline*}
\epr

\blm\label{TH:diff-oper<-morfizm-rassl-struuj} Let $M$ be a smooth manifold, $F$ a $C^*$-algebra, and $\ph:{\mathcal E}(M)\to F$ a homomorphism of involutive stereotype algebras, with $\ph({\mathcal E}(M))$ lying in the center of  $F$:
$$
\ph({\mathcal E}(M))\subseteq Z(F).
$$
Then for any stereotype space $X$ each morphism of the jet bundle $\nu:\Jet_{{\mathcal E}(M)}^n({\mathcal E}(M,X))\to \Jet_{{\mathcal E}(M)}^0(F)$ defines a unique differential operator of order $n$ between stereotype  ${\mathcal E}(M)$-modules $D:{\mathcal E}(M)\to F$, such that
\beq\label{j^0(Dx)=mu-circ-j^n(x)}
 \xymatrix @R=2pc @C=1.2pc
 {
 \Jet_{{\mathcal E}(M)}^0[{\mathcal E}(M,X)]\ar[dd]_{\nu} & \\
 & M \ar[ul]_{\jet^n(u)}\ar[dl]^{\jet^n(Du)} \\
 \Jet_{{\mathcal E}(M)}^0[F] &
 }
\qquad \jet^0(Du)=\nu\circ \jet^n(u),\qquad
u\in  {\mathcal E}(M,X).
 \eeq
In other words, $\nu$ is the morphism of the jet bundles, generated by the differential operator $D$ by Theorem  \ref{TH:diff-oper->rassl-struuj}:
 $$
 \nu=\jet_n[D].
 $$
 \elm
\bpr By Theorem \ref{TH:B-cong-Sec(val_AB)}, the mapping
$v:F\to\Sec(\pi^0_{{{\mathcal E}(M)},F})$, that turns $F$ into the algebra of continuous sections of the value bundle $\pi^0_{{{\mathcal E}(M)},F}: \Jet^0_{{\mathcal E}(M)}F\to\Spec({{\mathcal E}(M)})$ over the algebra  ${{\mathcal E}(M)}$, is an isomorphism of $C^*$-algebras:
$$
F\cong\Sec(\pi^0_{{{\mathcal E}(M)},F}).
$$
Consider the reverse isomorphism $v^{-1}:\Sec(\pi^0_{{{\mathcal E}(M)},F})\to F$:
\beq\label{lambda(j^0(b))=b}
v^{-1}(\jet^0(b))=b,\qquad b\in F.
\eeq
To each morphism of jet bundles $\nu:\Jet_{{\mathcal E}(M)}^n({{\mathcal E}(M,X)})\to \Jet_{{\mathcal E}(M)}^0(F)$ one can assign an operator $D:{{\mathcal E}(M,X)}\to F$ by the formula
\beq\label{Da=lambda(mu-circ-j^n(a))}
Du=v^{-1}\Big(\nu\circ \jet^n(u)\Big),\qquad u\in {{\mathcal E}(M,X)}.
\eeq
Evidently, $D$ satisfies the identity \eqref{j^0(Dx)=mu-circ-j^n(x)}. On the other hand, in the space of smooth functions on $M$ with values in an arbitrary stereotype space $X$ (as well as in the usual space of functions with values in $\C$) the Newton-Leibnitz formula is true, hence for a given chart, the Hadamard lemma \cite{Petrovsky} and the Taylor expansion with the remainder in $\overline{I_t^{n+1}\cdot X}$ hold (with some $n\in\N$). Since the action of $D$ on element $x$ can be factored through the jet $\jet^n(x)$, $D$ can be linearly expressed through the coefficients of the Taylor decomposition of the element $x$ in the neighbourhood of given $t\in M$ (with the choice of a chart). These Taylor coefficients are differential operators over ${\mathcal E}(M)$. As a corollary, $D$ is also a differential operator over ${\mathcal E}(M)$.
\epr

\blm The mapping $\iota:{\mathcal E}(M)\circledast A\to {\mathcal E}(M,A)$ is a smooth extension.
\elm
\bpr Let $D:{\mathcal E}(M)\circledast A\to B[m]$ be a morphism into a $C^*$-algebra $B$ with joined self-adjoint nilpotent elements. Consider the family of partial derivatives $D_k:{\mathcal E}(M)\circledast A\to B$, generated by $D$, and put
$$
\eta_k(u)=D_k(u\circledast 1),\qquad \alpha_k(a)=D_k(1\circledast a),\qquad u\in {\mathcal E}(M), \ a\in A.
$$
Then $\eta:{\mathcal E}(M)\to B[m]$ and $\alpha:A\to B[m]$ are morphisms of involutive stereotype algebras, and by Lemma \ref{LM:ph:A-circledast-B->C},
\beq\label{D(u-circledast-a)=eta(u)-cdot-alpha(a)}
D(u\circledast a)=\eta(u)\cdot\alpha(a)=\alpha(a)\cdot\eta(u),\qquad u\in {\mathcal E}(M), \ a\in A,
\eeq
In particular,
\beq\label{D_0(u-circledast-a)=eta(u)-cdot-alpha(a)}
D_0(u\circledast a)=\eta_0(u)\cdot\alpha_0(a)=\alpha_0(a)\cdot\eta_0(u),\qquad u\in {\mathcal E}(M), \ a\in A.
\eeq

Consider the operator $\eta_0$ and denote by $C$ its image in $B$:
$$
C=\overline{\eta_0({\mathcal E}(M))}.
$$
Let $F$ be the commutant of $C$ in $B$:
$$
F=C^!=\{x\in B:\quad \forall c\in C\quad x\cdot c=c\cdot x\}.
$$
Since $C$ is a commutative algebra, it also lies in $F$, and moreover, in the center of $F$:
$$
C\subseteq Z(F).
$$
Note also that the images of all operators $D_k$ lie in $F$:
\beq\label{Im-D_k-subseteq-F}
D_k({\mathcal E}(M)\circledast A)\subseteq F.
\eeq
For $k=0$ this can be proved directly:
\begin{multline*}
D_0(v\circledast a)\cdot \eta_0(u)=D_0(v\circledast a)\cdot D_0(u\circledast 1)=D_0\big((v\circledast a)\cdot (u\circledast 1)\big)=D_0\big((v\cdot u)\circledast 1\big)=D_0\big((u\cdot v)\circledast 1\big)=\\=D_0\big((u\circledast 1)\cdot (v\circledast a)\big)=D_0(u\circledast 1)\cdot D_0(v\circledast a)=\eta_0(u)\cdot D_0(v\circledast a)
\end{multline*}
And for $k>0$ we have to apply \eqref{D_k(A)-subseteq-Z^1(D_0)}: since for $k>0$ the values of the operators $D_k$ and $D_0$ commute, we obtain
$$
D_k(v\circledast a)\cdot\eta_0(u)=D_k(v\circledast a)\cdot D_0(u\circledast 1)=D_0(u\circledast 1)\cdot D_k(v\circledast a)=\eta_0(u)\cdot D_k(v\circledast a).
$$

To show that $\iota$ is a smooth extension, we have to construct the system of differential partial derivatives $\{D_k';\ k\in\N^d\}$ (over ${\mathcal E}(M,A)$!), which extend the operators $D_k$ from ${\mathcal E}(M)\circledast A$ to ${\mathcal E}(M,A)$ with values in $F$:
 \beq\label{prodolzhenie-D_k-na-E(M,A)}
 \xymatrix @R=2pc @C=1.2pc
 {
 {\mathcal E}(M)\circledast A\ar[rr]^{\iota}\ar[dr]_{D_k} & & {\mathcal E}(M,A)\ar@{-->}[dl]^{D_k'}\\
  & F &
 }
 \eeq

Each differential operator $D_k:{\mathcal E}(M)\circledast A\to F$ by Theorem
\ref{TH:diff-oper->rassl-struuj} generates a morphism of the jet bundles
$\jet_n[D_k]:\Jet_{{\mathcal E}(M)}^n[{\mathcal E}(M)\circledast A]\to \Jet_{{\mathcal E}(M)}^0(F)=\pi^0_A F$, where  $n=|k|$, such that
$$
\jet^0(D_k x)=\jet_n[D_k]\circ \jet^n(x),\qquad x\in{\mathcal E}(M)\circledast A.
$$
By Lemma \ref{LM:J^n_E(M)E(M)-circledast-A-cong-J^n_E(M)E(M,A)}, the jet bundles of algebras  ${\mathcal E}(M)\circledast A$ and ${\mathcal E}(M,A)$ are isomorphic. Denote this isomorphism by $\mu:\Jet_{{\mathcal E}(M)}^n[{\mathcal E}(M)\circledast A]\gets \Jet^n_{{\mathcal E}(M)}[{\mathcal E}(M,A)]$. Consider the composition $\nu=\jet_n[D_k]\circ\mu:\Jet_{{\mathcal E}(M)}[{\mathcal E}(M,A)]\to \Jet_{{\mathcal E}(M)}^0(F)=\pi^0_A F$:
$$
 \xymatrix @R=2pc @C=1.2pc
 {
 \Jet_{{\mathcal E}(M)}^n[{\mathcal E}(M)\circledast A]\ar[dr]_{\jet_n[D_k]} & & \Jet^n_{{\mathcal E}(M)}[{\mathcal E}(M,A)]\ar@{-->}[dl]^{\quad\nu=\jet_n[D_k]\circ\mu}\ar[ll]_{\mu}\\
  & \Jet^0_{{\mathcal E}(M)}[F] &
 }
$$
By Lemma \ref{TH:diff-oper<-morfizm-rassl-struuj} this dotted arrow $\nu$ generates a differential operator $D_k':{\mathcal E}(M,A)\to F$ (over the algebra ${\mathcal E}(M)$), such that
$$
\jet^0(D_k'f)=\jet_n[D_k]\circ \jet^n(f),\qquad  f\in {\mathcal E}(M,A).
$$
For each $x\in {\mathcal E}(M)\circledast A$ we have
\beq\label{PROOF:TH:E(M,A)-1}
\jet^0\big(D_k'\iota(x)\big)=\jet_n[D_k]\circ \jet^n(\iota(x))=\jet_n[D_k]\circ
\jet^n(x)=\jet^0(D_kx).
\eeq
Note then that by Theorem \ref{TH:B-cong-Sec(val_AB)}, the mapping
$\jet^0=v:F\to\Sec(\pi^0_{{\mathcal E}(M)}F)=\Sec(\Jet^0_{{\mathcal E}(M)}F)$, that turns $F$ into the algebra of continuous sections of the value bundle $\pi^0_{{\mathcal E}(M)}F: \Jet^0_{{\mathcal E}(M)}F\to\Spec({{\mathcal E}(M)})$ over the algebra ${{\mathcal E}(M)}$, is an isomorphism of $C^*$-algebras:
$$
F\cong\Sec(\pi^0_{{\mathcal E}(M)}F).
$$
Hence we can apply the operator reverse to $\jet^0$ to \eqref{PROOF:TH:E(M,A)-1}, and we get the equality
$$
D_k'\iota(x)=D_kx.
$$
I.e. $D_k'$ extends $D_k$ in the diagram \eqref{prodolzhenie-D_k-na-E(M,A)}. Besides this, from the fact that $\iota$ maps ${\mathcal E}(M)\circledast A$ densely into ${\mathcal E}(M,A)$ it follows that the conditions
\eqref{DEF:sist-chastn-proizv-*}-\eqref{DEF:sist-chastn-proizv-2} are inherited from $D_k$ to $D_k'$.

By construction, the operators $D_k'$ are differential with respect to the algebra ${\mathcal E}(M)$, but this is not enough for us: each $D_k'$ must be a differential operator of order $|k|$ with respect to the algebra ${\mathcal E}(M,A)$, and the operators $D_k'$ must form a system of partial derivatives on ${\mathcal E}(M,A)$.

Both propositions follow from the fact that the operators $D_k$ form a differential system of partial derivatives on ${\mathcal E}(M)\circledast A$. First, each $D_k$ is a differential operator of order $|k|$ with respect ${\mathcal E}(M)\circledast A$, hence for each $u_0,u_1,...,u_{|k|}\in {\mathcal E}(M)$ we have
$$
[...[[D_k,u_0\circledast a_0],u_1\circledast a_1],... u_{|k|}\circledast a_{|k|}]=0.
$$
This implies
$$
[...[[D_k',\iota(u_0\circledast a_0)],\iota(u_1\circledast a_1)],... \iota(u_{|k|}\circledast a_{|k|})]=0,
$$
and, since elements of the form $\iota(u\circledast a)$ are total in ${\mathcal E}(M,A)$, this means that they can be replaced by arbitrary vectors from ${\mathcal E}(M,A)$, and we obtain that $D_k'$ is a differential operator of order $|k|$ over ${\mathcal E}(M,A)$. Second, formulas \eqref{DEF:sist-chastn-proizv-*}-\eqref{DEF:sist-chastn-proizv-2} can be transferred from $D_k$ to $D_k'$. For instance, \eqref{DEF:sist-chastn-proizv-2} for $D_k$,
$$
D_k(x\cdot y)=\sum_{0\le l\le k}\begin{pmatrix}k \\ l\end{pmatrix}\cdot D_{k-l}(x)\cdot D_l(y),\qquad x,y\in {\mathcal E}(M)\circledast A,
$$
implies
$$
D_k'(\iota(x)\cdot\iota(y))=D_k'(\iota(x\cdot y))=D_k(x\cdot y)=\sum_{0\le l\le k}\begin{pmatrix}k \\ l\end{pmatrix}\cdot D_{k-l}(x)\cdot D_l(y)=\sum_{0\le l\le k}\begin{pmatrix}k \\ l\end{pmatrix}\cdot D_{k-l}(\iota(x))\cdot D_l(\iota(y)),
$$
This is true for each $x,y\in {\mathcal E}(M)\circledast A$. Since the image of $\iota$ is dense in  ${\mathcal E}(M,A)$ (Lemma \ref{LM:iota-in-DEpi-E(M)}), we see that $\iota(x)$ and $\iota(y)$ can be replaced by arbitrary vectors from ${\mathcal E}(M,A)$, i.e. \eqref{DEF:sist-chastn-proizv-2} holds for operators $D_k'$ as well.
\epr

\blm The mapping $\iota:{\mathcal E}(M)\circledast A\to {\mathcal E}(M,A)$ is a smooth envelope.
\elm
\bpr
Suppose $\sigma:{\mathcal E}(M)\circledast A\to C$ is another smooth extension. We need to verify that there exists a morphism $\upsilon$, such that the following diagram is commutative:
$$
 \xymatrix @R=2pc @C=1.2pc
 {
{\mathcal E}(M)\circledast A\ar[rr]^{\sigma}\ar[dr]_{\iota} & & C \ar@{-->}[dl]^{\upsilon}\\
  & {\mathcal E}(M,A) &
 }
$$

Take a chart $\ph:U\to V$, where $U\subseteq M$, $V\subseteq\R^d$, and let $K\subseteq U$ be a compact set that coincides with the closure of its interior:
$\overline{\Int(K)}=K$. The operators from Example \ref{EX:sist-chast-proizv-na-C-infty(M)}
$$
\varPhi_k(u)=\frac{\partial^{\abs{k}}(u\circ\ph^{-1})}{\partial t_1^{k_1}...\partial
t_d^{k_d}}\circ\ph,\qquad u\in {\mathcal E}(M),\quad k\in\N^d,
$$
form a system of partial derivatives on ${\mathcal E}(M)$ with values in ${\mathcal C}(K)$.
Let $\varPhi:{\mathcal E}(M)\to \mathcal{C}(K)[m]$ be a differential homomorphism, corresponding to the system $\{\varPhi_k\}$.

Take an arbitrary homomorphism $\eta:A\to B[n]$ into a $C^*$-algebra with joined self-adjoint nilpotent elements $B[n]$, and put
$$
D(u\circledast a)=\varPhi(u)\circledast \eta(a),\qquad u\in {\mathcal E}(M),\quad a\in A.
$$
The mapping $D$ is a homomorphism from ${\mathcal E}(M)\circledast A$ into the algebra ${\mathcal C}(K)[m]\circledast B[n]$, which by \eqref{A[m]-circledast-B}, is isomorphic to $({\mathcal C}(K)\circledast B)[m\oplus n]$. The stereotype tensor product ${\mathcal C}(K)\circledast B$ is naturally mapped into the maximal tensor product of $C^*$-algebras ${\mathcal C}(K)\underset{\max}{\otimes}B$, which in its turn is isomorphic to ${\mathcal C}(K)\odot B$ and ${\mathcal C}(K,B)$:
$$
{\mathcal C}(K)\circledast B\to {\mathcal C}(K)\underset{\max}{\otimes}B\cong \eqref{A-min-B=A-odot-B}\cong
{\mathcal C}(K)\odot B\cong \eqref{C(M,X)-cong-C(M)-odot-X} \cong {\mathcal C}(K,B).
$$

Thus we can treat $D$ as a morphism into a $C^*$-algebra with the joined self-addjoined nilpotent elements ${\mathcal C}(K,B)[m\oplus n]$,
\begin{multline*}
D:{\mathcal E}(M)\circledast A\to {\mathcal C}(K)[m]\circledast B[n]\cong\eqref{A[m]-circledast-B},\eqref{A[m][n]=A[m-oplus-n]}\cong ({\mathcal C}(K)\circledast B)[m\oplus n]\to   ({\mathcal C}(K)\underset{\max}{\otimes}B)[m\oplus n]\cong\\ \cong\eqref{A-min-B=A-odot-B}\cong({\mathcal C}(K)\odot B)[m\oplus n]\cong\eqref{C(M,X)-cong-C(M)-odot-X}\cong {\mathcal C}(K,B)[m\oplus n]
\end{multline*}
Since $\sigma:{\mathcal E}(M)\circledast A\to C$ is a smooth envelope, the homomorphism $D:{\mathcal E}(M)\circledast A\to {\mathcal C}(K,B)[m\oplus n]$ is uniquely extended to a homomorphism $D':C\to {\mathcal C}(K,B)[m\oplus n]$:
 \beq\label{DIAGR:env-E(M)-circledast-A-0}
 \xymatrix @R=2pc @C=1.2pc
 {
{\mathcal E}(M)\circledast A\ar[rr]^{\sigma}\ar[dr]_{D} & & C \ar@{-->}[dl]^{D'}\\
  & {\mathcal C}(K,B)[m\oplus n] &
 }
 \eeq
Note that ${\mathcal C}(K,B)[m\oplus n]$ is isomorphic to ${\mathcal C}\big(K, B[n]\big)[m]$,
$$
{\mathcal C}(K,B)[m\oplus n]\cong\eqref{A[m][n]=A[m-oplus-n]}\cong {\mathcal C}(K,B)[n][m]\cong\eqref{C(M,B)[n]=C(M,B[n])}\cong {\mathcal C}\big(K, B[n]\big)[m],
$$
so diagram \eqref{DIAGR:env-E(M)-circledast-A-0} can be changed as follows:
 \beq\label{DIAGR:env-E(M)-circledast-A-1}
 \xymatrix @R=2pc @C=1.2pc
 {
{\mathcal E}(M)\circledast A\ar[rr]^{\sigma}\ar[dr]_{D} & & C \ar@{-->}[dl]^{D'}\\
  & {\mathcal C}(K,B[n])[m] &
 }
 \eeq
Now we can return back to the system of partial derivatives $D_k$, and for each index $k\in\N[m]$ we have
 $$
  \xymatrix @R=2pc @C=1.2pc
 {
 {\mathcal E}(M)\circledast A\ar[rr]^{\sigma}\ar[dr]_{D_k} & & C \ar@{-->}[dl]^{D_k'}\\
  & {\mathcal C}(K,B[n]) &
 }
 $$

Take an arbitrary element $c\in C$. Since $\sigma:{\mathcal E}(M)\circledast A\to C$ is a dense epimorphism, there is a net of elements $x_i\in {\mathcal E}(M)\circledast A$ such that
$$
\sigma(x_i)\overset{C}{\underset{i\to\infty}{\longrightarrow}} c.
$$
For each index $k\in\N[m]$ we have
 \beq\label{D_k(alpha_i)->D_k'(c)}
D_k(x_i)=D_k'(\sigma(x_i))\overset{\mathcal{C}(K)}{\underset{i\to\infty}{\longrightarrow}}
D_k'(c).
 \eeq
Now let us take a smooth curve in $K$, i.e. a smooth mapping $\gamma:[0,1]\to K$. For each index $k\in\N[m]$ of order $|k|=1$ and for each $t\in[0,1]$ let $\gamma^k(t)$ denote the $k$-th component of the derivative $\gamma'(t)$ in the expansion by local coordinates on $U$. For each function $u\in{\mathcal E}(M)$ we have by the Newton-Leibnitz theorem
$$
\varPhi_0(u)(\gamma(1))-\varPhi_0(u)(\gamma(0))=\sum_{|k|=1}\int_0^1\gamma^k(t)\cdot
\varPhi_k(u)(\gamma(t))\d t.
$$
If we multiply $u$ by an arbitrary element $a\in A$, we get
\begin{multline*}
D_0(u\circledast a)(\gamma(1))-D_0(u\circledast a)(\gamma(0))=\varPhi_0(u)(\gamma(1))\circledast\eta(a) -\varPhi_0(u)(\gamma(0))\circledast\eta(a)=\\=\sum_{|k|=1}\int_0^1\gamma^k(t)\cdot
\varPhi_k(u)(\gamma(t))\circledast\eta(a)\d t=\sum_{|k|=1}\int_0^1\gamma^k(t)\cdot
D_k(u\circledast a)(\gamma(t))\d t.
\end{multline*}
Since elements of the form $u\circledast a$ are total in ${\mathcal E}(M)\circledast A$, we can replace them by an arbitrary element $x\in{\mathcal E}(M)\circledast A$.

Together with \eqref{D_k(alpha_i)->D_k'(c)} this gives
 \begin{multline*}
D_0'(c)(\gamma(1))-D_0'(c)(\gamma(0)) \underset{\infty\gets i}{\longleftarrow}
D_0(x_i)(\gamma(1))-D_0(x_i)(\gamma(0))=\\
=\sum_{|k|=1}\int_0^1\gamma^k(t)\cdot
D_k(x_i)(\gamma(t))\d t \underset{i\to\infty}{\longrightarrow}
\sum_{|k|=1}\int_0^1\gamma^k(t)\cdot D_k'(c)(\gamma(t))\d t,
 \end{multline*}
therefore.
$$
D_0'(c)(\gamma(1))-D_0'(c)(\gamma(0))=\sum_{|k|=1}\int_0^1\gamma^k(t)\cdot
D_k'(c)(\gamma(t))\d t.
$$
This connection between the function $D_0'(c)\in {\mathcal C}(K,B)$ and the functions $D_k'(c)\in
{\mathcal C}(K,B)$, $|k|=1$, means that $D_0'(c)$ is continuously differentiable on $K$, and its partial derivatives (in our local coordinates) are the functions $D_k'(c)$, $|k|=1$.

After that we can take one of the derivatives $D_k'(c)$, $|k|=1$, and consider the indices of order 2. Using the same trick, we see that $D_k'(c)$ is also continuously differentiable. The induction over indices shows that the functions $D_k'(c)$ are infinitely differentiable, and are related to each other as partial derivatives of the function  $D_0'(c)$ (with respect to the chosen local coordinates). This means that diagram \eqref{DIAGR:env-E(M)-circledast-A-1} is commutative:
 \beq\label{DIAGR:env-E(M)-circledast-A-*-K}
 \xymatrix @R=2.5pc @C=2pc
 {
 {\mathcal E}(M)\circledast A\ar[rr]^{\sigma}\ar@/_5ex/[ddr]_{D}\ar[dr]_{\iota_{K,B}} & & C \ar@/^5ex/[ddl]^{D'}\ar@{-->}[dl]^{\iota_{K,B}'}\\
  & {\mathcal E}(K,B[n])\ar[d]_{\varPhi_B} & \\
  & \mathcal{C}(K,B[n])[m] &
 }
 \eeq
where
$$
\iota_{K,B}(u\circledast a)=u(t)\cdot\eta(a),\qquad u\in {\mathcal E}(M),\quad a\in A,
$$
$$
(\varPhi_B)_k(f)=\frac{\partial^{\abs{k}}(f\circ\ph^{-1})}{\partial t_1^{k_1}...\partial
t_d^{k_d}}\circ\ph,\qquad f\in {\mathcal E}(M,B[n]),\quad k\in\N^d,
$$
From this diagram it follows that $\iota_{K,B}'$ is continuous, since if $c_i\to c$, then this condition is preserved under the action of the operator $D_k'$, i.e. $D_k'(c_i)\to D_k'(c)$, and this is exactly the convergence in the space ${\mathcal E}(K,B[n])$.

If now we change the compact set $K\subset U$ and the open set $U\subseteq M$, then the arising smooth functions $D_0'(c)$ on $K$ are compatible to each other, since on the intersections of their domains they coincide. Thus a common smooth function is defined $\iota_B'(c):M\to B[n]$, such that its restriction to each compact set $K$ coincides with the corresponding function $D_0'(c)$:
$$
\iota'(c)\big|_K=D_0'(c),\qquad K\subset U\subseteq M.
$$
and the partial derivatives under the chosen system of coordinates coincide with the values of the operators $D_k'$ on $c$. In other words, a mapping is defined $\iota'_B:C\to{\mathcal E}(M,B[n])$ (by construction this is a homomorphism of algebras), such that the following diagram specifying \eqref{DIAGR:env-E(M)-circledast-A-*-K} is commutative:
 \beq\label{DIAGR:env-E(M)-circledast-A-*-M-0}
 \xymatrix @R=2.5pc @C=2pc
 {
 {\mathcal E}(M)\circledast A\ar[rr]^{\sigma}\ar@/_5ex/[ddr]_{\iota_{K,B}}\ar[dr]_{\iota_B} & & C \ar@/^5ex/[ddl]^{\iota'_{K,B}}\ar@{-->}[dl]^{\iota'_B}\\
  & {\mathcal E}(M,B[n])\ar[d]_{\rho_K} & \\
  & {\mathcal E}(K,B[n]) &
 }
 \eeq
(here $\rho_K$ is the mapping of restriction to $K$).

Let now $U$ be a differential neighbourhood of zero in $A$ generated by a homomorphism $\eta:A\to B[n]$. Since $\sigma$ is a dense epimorphism, the upper inner triangle in \eqref{DIAGR:env-E(M)-circledast-A-*-M-0} can be completed to the diagram
 \beq\label{DIAGR:env-E(M)-circledast-A-*-M}
 \xymatrix @R=2.5pc @C=2pc
 {
 {\mathcal E}(M)\circledast A\ar[rr]^{\sigma}\ar@/_5ex/[ddr]_{\iota_B}\ar[dr]_{\vartheta_U} & & C \ar@/^5ex/[ddl]^{\iota'_B}\ar@{-->}[dl]^{\vartheta'_U}\\
  & {\mathcal E}(M,A/U)\ar[d]_{\eta_U\oslash 1_M} & \\
  & {\mathcal E}(M,B[n]) &
 }
 \eeq
where $\eta_U:A/U\to B[n]$ is the morphism from \eqref{ph=ph_U-circ-pi_U}, and
$$
\vartheta_U(u\circledast a)(t)=u(t)\cdot\pi_U(a),\qquad u\in{\mathcal E}(M),\quad a\in A,\quad t\in M,
$$
$$
(\eta_U\oslash 1_M)(h)(t)=\eta_U(h(t)),\qquad h\in {\mathcal E}(M,A/U),\quad t\in M.
$$
From the definition of $\vartheta_U$ it follows that of $U'\subseteq U$ is another differential neighbourhood of zero, then
\beq\label{vartheta_U=(varkappa^U'_U-oslash-1_M)-cdot-vartheta_U'}
\vartheta_U=(\varkappa^{U'}_U\oslash 1_M)\cdot\vartheta_{U'},\qquad U\supseteq U',
\eeq
where $\varkappa^{U'}_U$ is a morphism from \eqref{varkappa^U'_U}, and
$$
(\varkappa^{U'}_U\oslash 1_M)(h)(t)=\varkappa^{U'}_U(h(t)),\qquad h\in {\mathcal E}(M,A/U'),\quad t\in M.
$$
The equality \eqref{vartheta_U=(varkappa^U'_U-oslash-1_M)-cdot-vartheta_U'} is the left lower inner triangle in the diagram
 $$
 
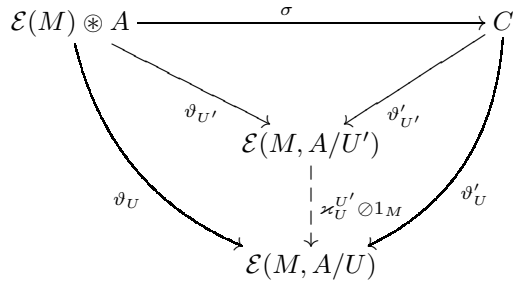
\xymatrix @R=2.5pc @C=3pc
 {
{\mathcal E}(M)\circledast A\ar[rr]^{\sigma}\ar[dr]_(.6){\vartheta_{U'}}\ar@/_5ex/[ddr]_{\vartheta_U} & & C \ar[dl]^(.6){\vartheta'_{U'}} \ar@/^5ex/[ddl]^{\vartheta'_U}\\
 & {\mathcal E}(M,A/U')\ar@{-->}[d]^{\varkappa^{U'}_U\oslash 1_M} & \\
  & {\mathcal E}(M,A/U) &
 }
$$
Here the perimeter and the upper inner triangle are variants of the upper inner triangle in \eqref{DIAGR:env-E(M)-circledast-A-*-M}, and in addition $\sigma$ is an epimorphism. As a corollary, the remaining right lower triangle also must be commutative.

This means that the morphisms $\vartheta'_U:C\to {\mathcal E}(M,A/U)$ form a projective cone of the system  $\varkappa^{U'}_U\oslash 1_M$, and thus there is a morphism $\vartheta'$ into the projective limit:
$$
 
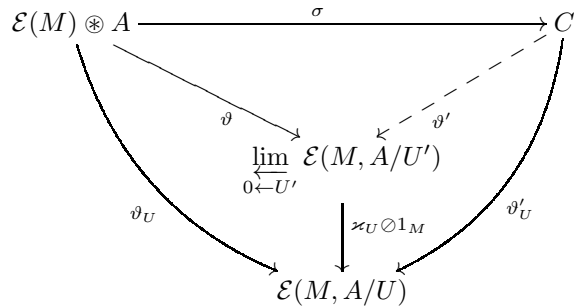
\xymatrix @R=2.5pc @C=3pc
 {
{\mathcal E}(M)\circledast A\ar[rr]^{\sigma}\ar[dr]_(.6){\vartheta}\ar@/_5ex/[ddr]_{\vartheta_U} & & C \ar@{-->}[dl]^(.6){\vartheta'} \ar@/^5ex/[ddl]^{\vartheta'_U}\\
 & \leftlim\limits_{0\gets U'}{\mathcal E}(M,A/U')\ar[d]^{\varkappa_U\oslash 1_M} & \\
  & {\mathcal E}(M,A/U) &
 }
$$
Now we note the chain
$$
\leftlim\limits_{0\gets U'}{\mathcal E}(M,A/U')=\leftlim\limits_{0\gets U'}({\mathcal E}(M)\odot A/U')=\cite[(2.53)]{Akbarov-env}=
{\mathcal E}(M)\odot \leftlim\limits_{0\gets U'}A/U'={\mathcal E}(M,\leftlim\limits_{0\gets U'}A/U')
$$
and place the last space into our diagram:
$$
 
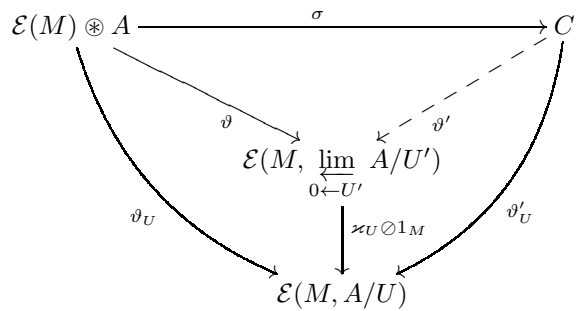
\xymatrix @R=2.5pc @C=3pc
 {
{\mathcal E}(M)\circledast A\ar[rr]^{\sigma}\ar[dr]_(.6){\vartheta}\ar@/_5ex/[ddr]_{\vartheta_U} & & C \ar@{-->}[dl]^(.6){\vartheta'} \ar@/^5ex/[ddl]^{\vartheta'_U}\\
 & {\mathcal E}(M,\leftlim\limits_{0\gets U'}A/U')\ar[d]^{\varkappa_U\oslash 1_M} & \\
  & {\mathcal E}(M,A/U) &
 }
$$

Again let us recall that $\sigma$ is a dense epimorphism. This implies that the arrow $\vartheta'$ can be lifted to an arrow $\upsilon$ with values in the space ${\mathcal E}(M,\Im\pi)$ of functions with values in the image of the mapping $\pi:A\to\leftlim\limits_{0\gets U'}A/U'$, or, what is the same, in the immediate subspace, generated by the set of values of the mapping $\pi$, and this space coincides with $A$, since $A$ is a smooth algebra:
$$
\Im\pi\cong\Env_{\mathcal E}A\cong A
$$
We obtain the diagram
$$
 
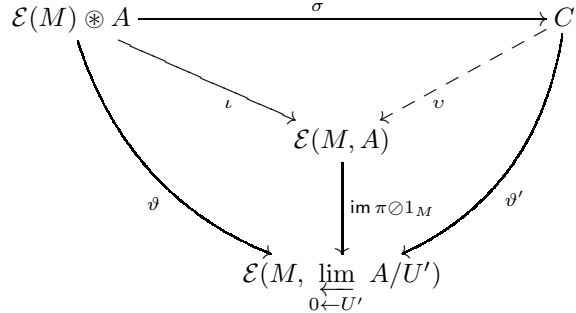
\xymatrix @R=2.5pc @C=3pc
 {
{\mathcal E}(M)\circledast A\ar[rr]^{\sigma}\ar[dr]_(.6){\iota}\ar@/_5ex/[ddr]_{\vartheta} & & C \ar@{-->}[dl]^(.6){\upsilon} \ar@/^5ex/[ddl]^{\vartheta'}\\
 & {\mathcal E}(M,A)\ar[d]^{\im\pi\oslash 1_M} & \\
  & {\mathcal E}(M,\leftlim\limits_{0\gets U'}A/U') &
 }
$$
where $\pi$ is the morphism from \eqref{A-to-leftlim-A/U}.
\epr

\subsection{${\mathcal E}(M)$, as a smooth envelope of its subalgebras}

Again, let $M$ be a smooth (locally euclidean) manifold, and ${\mathcal E}(M)$ the algebra of smooth functions on  $M$. We denote by $T_s(M)$, $T_s^\star(M)$ and $\Jet^n_s(M)$ the usual tangent space, cotangent space and the algebra of jets in a given point $s\in M$. These notations are related with those we introduced before by the equations
$$
T_s(M)=T_s[{\mathcal E}(M)],\quad T_s^\star(M)=T_s^\star[{\mathcal E}(M)],\quad \Jet_s^n(M)=\Jet_s^n[{\mathcal E}(M)].
$$

\btm\label{C^infty-obolochka-podalgebry-v-C^infty}
Let $A$ be an involutive stereotype subalgebra in the algebra $\mathcal{E}(M)$ of smooth functions on a smooth (locally euclidean) manifold $M$, i.e. there is a (continuous and unital) monomorphism of involutive stereotype algebras
$$
\iota:A\to \mathcal{E}(M).
$$
The smooth envelope of $A$ coincides with ${\mathcal E}(M)$
\beq\label{Env_C^infty_A=C^infty(M)}
\Env_{\mathcal E} A={\mathcal E}(M)
\eeq
(i.e. $\iota$ is a smooth envelope of $A$), if and only if the following two conditions hold:
 \bit{
\item[(i)] the dual mapping of spectra
$$
\Spec(A)\gets M
$$
is an exact covering\footnote{In the sense of definition on page \pageref{DEF:nalozhenie}.};

\item[(ii)] for each point $s\in M$ the natural mapping of tangent spaces
 $$
T_s[A]\gets T_s(M)
 $$
is an isomorphism (of finite dimensional vector spaces).
 }\eit
\etm

We split the proof into 5 lemmas.

\blm\label{LM:neobhodimost-v-TH:C^infty-obolochka-podalgebry-v-C^infty} Conditions (i) and (ii) are necessary for the enclosure $A\subseteq{\mathcal E}(M)$ to be a smooth extension of $A$.
\elm
 \bpr
Assume that $\iota$ is a smooth extension, i.e. an extension in the class $\DEpi$ of dense epimorphisms with respect to the class of differential involutive homomorphisms into  $C^*$-algebras with joined self-adjoint nilpotent elements. Then $\iota$ is an extension in $\DEpi$ with respect to the class of homomorphisms into $C^*$-algebras, since a $C^*$-algebra $B$ can be considered as a $C^*$-algebra with joined empty set of nilpotent elements: $B=B[0]$. I.e., $\iota$ is a continuous extension. The same reasoning as in Theorem \eqref{C-obolochka-podalgebry-v-C(M)} proves that the mapping of spectra $\iota^{\Spec}:\Spec(A)\gets M$ is an exact covering.

So we have to prove the condition (ii). Consider the algebra $\C_1[[1]]$ of polynomials of degree 1 of one variable. As a vector space it is isomorphic to the direct sum $\C\oplus\C$, and the morphisms $A\to \C_1[[1]]$ are the pairs
$(s,\sigma)$, where $s\in\Spec(A)$ is the point of the spectrum, and $\sigma\in T_s[A]$ the tangent vector in this point. For the point $t\in M$ in (ii) we have a commutative diagram
$$
 \xymatrix @R=2pc @C=1.2pc
 {
 A\ar[rr]^{\iota}\ar[dr]_{(t\circ\iota,\tau\circ\iota)} & & {\mathcal E}(M)\ar[dl]^{(t,\tau)}\\
  & \C_1[[1]] &
 }
$$
Since $\iota:A\to {\mathcal E}(M)$ is a smooth extension, each arrow $A\to\C_1[[1]]$ generates a unique arrow ${\mathcal E}(M)\to\C_1[[1]]$, and this means that the mapping of tangent spaces
$\tau\mapsto\iota^\star(\tau)=\tau\circ\iota$ is a bijection, and hence an isomorphism of vector spaces.  \epr

\blm\label{LM:Ker-rho_t^(T^star)=0} Let $\iota:A\to B$ be a homomorphism of involutive stereotype algebras, and for a point $t\in\Spec(B)$ the corresponding morphism of cotangent spaces $\iota_t^{\C
T^\star}:{\C}T_{t\circ\iota}^\star[A]\to{\C}T_t^\star[B]$ is injective:
 \beq\label{Ker-rho_t^(T^star)=0}
\Ker\iota_t^{T^\star}=0
 \eeq
Then
 \beq\label{overline(I_s^2)[A]=rho^(-1)(overline(I_t^2)[B])}
\overline{I_{t\circ\iota}^2}[A]=\iota^{-1}\Big(\overline{I_t^2}[B]\Big).
 \eeq
\elm \bpr The direct inclusion follows from the homomorphy (and continuity) of $\iota$ (we use the notations on page \pageref{DEF:M-cdot-N}):
 \begin{multline*}
\iota\Big(I_{t\circ\iota}[A]\Big)\subseteq I_t[B]\quad\Longrightarrow\\
\Longrightarrow\quad
\iota\Big(I_{t\circ\iota}^2[A]\Big)=\iota\Big(I_{t\circ\iota}[A]\cdot
I_{t\circ\iota}[A]\Big)=
\iota\Big(I_{t\circ\iota}[A]\Big)\cdot\iota\Big(I_{t\circ\iota}[A]\Big) \subseteq
I_s[B]\cdot I_s[B]=
I_s^2[B]\subseteq \overline{I_t^2}[B] \quad\Longrightarrow\\
\Longrightarrow\quad
I_{t\circ\iota}^2[A]\subseteq\iota^{-1}\Big(\overline{I_t^2}[B]\Big)
\quad\Longrightarrow\quad
\overline{I_{t\circ\iota}^2}[A]\subseteq\iota^{-1}\Big(\overline{I_t^2}[B]\Big).
 \end{multline*}
This inclusion can be presented in the form
$$
\iota\Big(\overline{I_{t\circ\iota}^2}[A]\Big)\subseteq \overline{I_t^2}[B],
$$
and we can conclude that the following diagram is commutative:
$$
 \xymatrix @R=2pc @C=3pc
 {
 & \overline{I_{t\circ\iota}^2}[A]\ar[r]^{\iota}\ar[d]_{\sigma_A} & \overline{I_t^2}[B]\ar[d]^{\sigma_B} &\\
 & I_{t\circ\iota}[A]\ar[r]^{\iota}\ar[d]_{\pi_A} & I_t[B]\ar[d]^{\pi_B} & \\
 {\C}T_{t\circ\iota}[A]\ar@{=}[r] & \Big(I_{t\circ\iota}[A]\Big/\overline{I_{t\circ\iota}^2}[A]\Big)^\triangledown\ar[r]^{\iota_t^{T^\star}} &
 \Big(I_t[B]\Big/\overline{I_t^2}[B]\Big)^\triangledown \ar@{=}[r] &  {\C}T_t[B]
 }
$$
It implies the inverse chain, that is necessary for the proof of
\eqref{overline(I_s^2)[A]=rho^(-1)(overline(I_t^2)[B])}:
 \begin{multline*}
a\in\iota^{-1}\Big(\overline{I_t^2}[B]\Big) \quad\Longrightarrow\quad \iota(a)\in
\overline{I_t^2}[B] \quad\Longrightarrow\quad
\iota_t^{T^\star}(\pi_A(a))=\pi_B(\iota(a))=0
\quad\Longrightarrow\\
\Longrightarrow\quad
\pi_A(a)\in\Ker\iota_t^{T^\star}=\eqref{Ker-rho_t^(T^star)=0}=0
\quad\Longrightarrow\quad a\in\Ker(\pi_A)=I_{t\circ\iota}[A].
 \end{multline*}
\epr

\blm\label{LM-3-0-dlya-TH:C^infty-obolochka-podalgebry-v-C^infty}
Under the conditions (i) and (ii) of Theorem \ref{C^infty-obolochka-podalgebry-v-C^infty} the natural morphism of bundles
$$
\Jet^n_{A}[{\mathcal E}(M)]\longleftarrow \Jet^n_{{\mathcal E}(M)}[{\mathcal E}(M)]
$$
is a fiber-wise isomorphism.
\elm
\bpr
By the Nachbin theorem \ref{TH:Nachbin} the algebra $A$ is dense in the algebra ${\mathcal E}(M)$. Hence by Lemma  \ref{LM:o-plotnom-ideale}, the ideal $I_t(A)=\{a\in A: a(t)=0\}$ is dense in the ideal $I_t({\mathcal E}(M))=\{u\in {\mathcal E}(M): u(t)=0\}$. Therefore,
$$
\overline{I_t^n(A)}=\overline{I_t^n({\mathcal E}(M))}=I_t^n({\mathcal E}(M)),
$$
and thus,
$$
{\mathcal E}(M)/\overline{I_t^n(A)}={\mathcal E}(M))/\overline{I_t^n({\mathcal E}(M))}={\mathcal E}(M)/I_t^n({\mathcal E}(M)).
$$
\epr

\blm\label{LM-3-dlya-TH:C^infty-obolochka-podalgebry-v-C^infty}
If $M$ is a smooth manifold and $\iota:A\to{\mathcal E}(M)$ a monomorphism of involutive stereotype algebras, then (ndependently on (i)) the condition (ii) of Theorem \ref{C^infty-obolochka-podalgebry-v-C^infty} is equivalent to each of the conditions
\bit{

\item[(iii)] for each point $s\in M$ the natural mapping of cotangent spaces
 $$
T_s^\star[A]\to T_s^\star(M)
 $$
is an isomorphism (of finite-dimensional vector spaces),

\item[(iv)] for each point $s\in M$ and for each number $n\in\N$ the natural mapping of the spaces of jets
 $$
\Jet_s^n[A]\to \Jet_s^n(M)
 $$
is an isomorphism (of finite-dimensional vector spaces).
 }\eit
and implies the condition
\bit{

\item[(v)] for each $n\in\N$ there is a continuous mapping $\mu:\Jet^n_A[A]\gets\Jet^n_{{\mathcal E}(M)}[{\mathcal E}(M)]$, which together with the mapping of spectra $\iota^{\Spec}:\Spec(A)\gets M$ forms a morphism of bundles,
$$
 \xymatrix @R=2pc @C=4pc
 {
\Jet^n_A[A]\ar[d]^{\pi^n_{A,A}} &  \Jet^n_{{\mathcal E}(M)}[{\mathcal E}(M)]\ar[d]^{\pi^n_{{\mathcal E}(M),{\mathcal E}(M)}}\ar[l]_{\mu} \\
\Spec(A)  & M\ar[l]_{\iota^{\Spec}}
 }
$$
    satisfying the identity
    \beq
    \mu(\jet^n(a\circ\iota^{\Spec})(t))=\jet^n(a)(\iota^{\Spec}(t)),\qquad a\in A,\quad t\in M.
    \eeq
 }\eit\noindent
If in addition the condition (i) of Theorem \ref{C^infty-obolochka-podalgebry-v-C^infty} holds, then the morphism of bundles $\mu$ is a bijection.

\elm
\bpr
1. Let us first show that the conditions (ii)-(iv) are equivalent. The fact that (ii) and (iii) are equivalent follows immediately from Theorem \ref{TH:svyaz-T(M)-i-T^*(M)}: if $T_s[A]$ and $T_s(M)$ are isomorphic as (finite dimensional) vector spaces, then their dual spaces $T_s^\star[A]$ and $T_s^\star(M)$ are isomorphic as well. And vice versa.

On the other hand (iv) implies (iii), since the isomorphism
$$
\Jet_s^1[A]\cong \Jet_s^1(M)
$$
maps the ideal $I_s[A]\cap \Jet_s^1[A]$ into the ideal $I_s(M)\cap \Jet_s^1(M)$, and this means the isomorphism
 \begin{multline*}
T_s^\star[A]=\Big(I_s[A]\Big/\overline{I_s^2}[A]\Big)^\triangledown=I_s[A]\Big/\overline{I_s^2}[A]
=I_s[A]\cap \Jet_s^1[A]\cong I_s(M)\cap \Jet_s^1(M)=\\=
I_s(M)\Big/\overline{I_s^2}(M)=
\Big(I_s(M)\Big/\overline{I_s^2}(M)\Big)^\triangledown=T_s^\star(M)
 \end{multline*}

Thus, in the equivalences $(ii)\Longleftrightarrow (iii)\Longleftrightarrow (iv)$ the only unclear link is the implication $(iii)\Longrightarrow (iv)$. Let us prove it. Suppose (iii) holds. Then first, we have equality \eqref{Ker-rho_t^(T^star)=0}, which in this case can be written as follows:
 \beq\label{overline(I_s^2)[A]=A-cap-I_s^2(M)}
\overline{I_s^2}[A]=A\cap I_s^2(M)
 \eeq
Further, since the space $T_s^\star[A]=\Big(I_s[A]\Big/\overline{I_s^2}[A]\Big)^\triangledown$ is finite dimensional, it coincides with the space $I_s[A]\Big/\overline{I_s^2}[A]$, which is finite dimensional as well, and therefore we can choose a finite basis there. In other words, there exists a sequence of vectors $e^1,...,e^d\in I_s[A]$ such that the cosets  $e^1+\overline{I_s^2}[A],...,e^d+\overline{I_s^2}[A]$ form a basis in $I_s[A]\Big/\overline{I_s^2}[A]$.
Take a neighbourhood $U$ of $s$, where the functions $e^1,...,e^d$ form a local chart of the manifold $M$. After that for each point $t\in U$ and for each multiindex $k\in\N^d$ we put
 \beq\label{e^k=(e^1)^(k_1)...(e^d)^(k_d)}
e_t^k=(e^1-e^1(t))^{k_1}\cdot...\cdot (e^d-e^d(t))^{k_d}.
 \eeq
For any function $f\in{\mathcal E}(M)$ and for any number $n\in\N$ let the symbol $E_s^n[f]$ denote the linear combination of functions $e^k$, $\abs{k}\le n$,
$$
E_s^n[f]=\sum_{|k|\le n}\lambda_k\cdot e_s^k.
$$
which in the point $s$ has the same jet of order $n$, as the function $f$:
 \beq\label{f-stackrel(I_s^(n+1)(M))(equiv)E_s^n[f]}
f\stackrel{I_s^{n+1}(M)}{\equiv} E_s^n[f]
 \eeq
(such a function exists and is unique, since $e_s^1,...,e_s^d$ form a local chart in a neighbourhood $V_s$ of $s$, and this chart diffeomorphically maps $V_s$ onto a neighbourhood of zero in $\R^d$, and in this diffeomorphism the functions $E_s^n[f]$ turn exactly into the Taylor polynomials of the function $f$ in the point 0).

2. Note that the operation $f\mapsto E_s^n[f]$ is multiplicative
 \beq\label{E_s^n[f-cdot-g]=E_s^n[f]-cdot-E_s^n[g]}
E_s^n[f\cdot g]=E_s^n[f]\cdot E_s^n[g],\qquad f,g\in{\mathcal E}(M)
 \eeq
satisfies a variant of the idempotency property
 \beq\label{idempotentnost-E_s^n[f]}
E_s^q\big[E_s^p[f]\big]=E_s^p[f]=E_s^p\big[E_s^q[f]\big],\qquad p\le
q\in\N,\quad f\in{\mathcal E}(M),
 \eeq
and is continuous in $A$:
 \beq\label{nepreryvnost-E_s^n}
a_i\overset{A}{\underset{i\to\infty}{\longrightarrow}}a
\quad\Longrightarrow\quad
E_s^n[a_i]\overset{A}{\underset{i\to\infty}{\longrightarrow}} E_s^n[a].
 \eeq
The first two properties follow from the representation of $E_s^n[f]$ as Taylor polynomials under the diffeomorphism formed by the local chart $e_1,...,e_d$, and the third one -- from the fact that the algebra of polynomials of a given degree is finite dimensional: if
$a_i\overset{A}{\underset{i\to\infty}{\longrightarrow}}a$, then
$a_i\overset{{\mathcal E}(M)}{\underset{i\to\infty}{\longrightarrow}}a$,
hence
$E_s^n[a_i]\overset{{\mathcal E}(M)}{\underset{i\to\infty}{\longrightarrow}}E_s^n[a]$,
and since the left and the right side belong to a finite dimensional space
$P=\sp\{e^k;\ |k|\le n\}$, the convergence in ${\mathcal E}(M)$ can be replaced by the convergence in $P$, and then by the convergence in $A$.

3. Now let us prove the formula
 \beq\label{a-stackrel(overline(I_s^(n+1))[A])(equiv)E_s^n[a]}
a\stackrel{\overline{I_s^{n+1}}[A]}{\equiv} E_s^n[a],\qquad a\in
\overline{I_s^n}[A].
 \eeq
This is done by induction. First, this formula is true for $n=0$,
 \beq\label{a-stackrel(I_s[A])(equiv)E_s^0[a]}
a\stackrel{I_s[A]}{\equiv} E_s^0[a],\qquad a\in A,
 \eeq
since
$$
a-E_s^0[a]=a-a(s)\cdot 1\in I_s[A].
$$
And for $n=1$,
 \beq\label{a-stackrel(overline(I_s^2)[A])(equiv)E_s^2[a]}
a\stackrel{\overline{I_s^2}[A]}{\equiv} E_s^1[a],\qquad a\in\overline{I_s}[A],
 \eeq
since \eqref{f-stackrel(I_s^(n+1)(M))(equiv)E_s^n[f]} implies
$a-E_s^1[a]\in I_s^2(M)$, while on the other hand, $a-E_s^1[a]\in A$, and together this means that $a-E_s^1[a]\in A\cap I_s^2(M)=\eqref{overline(I_s^2)[A]=A-cap-I_s^2(M)}=\overline{I_s^2}[A]$.

Assume then that the formula \eqref{a-stackrel(overline(I_s^(n+1))[A])(equiv)E_s^n[a]} is true for some integer  $n-1$:
 \beq\label{a-stackrel(overline(I_s^n)[A])(equiv)E_s^(n-1)[a]}
a\stackrel{\overline{I_s^n}[A]}{\equiv} E_s^{n-1}[a],\qquad a\in
\overline{I_s^{n-1}}[A].
 \eeq
Take $a\in\overline{I_s^n}[A]$. This means that there is a net $a_i\in I_s^n[A]$, tending to $a$ in $A$:
$$
a_i\overset{A}{\underset{i\to\infty}{\longrightarrow}}a.
$$
Each element $a_i$ belongs to $I_s^n[A]$, and thus, has the form
$$
a_i=\sum_{j=1}^p b_i^j\cdot c_i^j,\qquad b_i^j\in I_s^{n-1}[A],\quad c_i^j\in
I_s[A].
$$
From  \eqref{a-stackrel(overline(I_s^n)[A])(equiv)E_s^(n-1)[a]} and
\eqref{a-stackrel(overline(I_s^2)[A])(equiv)E_s^2[a]} we have the following chain:
$$
b_i^j\stackrel{\overline{I_s^n}[A]}{\equiv} E_s^{n-1}[b_i^j],\quad
c_i^j\stackrel{\overline{I_s^2}[A]}{\equiv} E_s^1[c_i^j],
$$
$$
\Downarrow
$$
$$
b_i^j\cdot c_i^j\stackrel{\overline{I_s^{n+1}}[A]}{\equiv}
E_s^{n-1}[b_i^j]\cdot E_s^1[c_i^j]=E_s^n[b_i^j\cdot c_i^j],
$$
$$
\Downarrow
$$
$$
a \overset{A}{\underset{\infty\gets i}{\longleftarrow}} a_i=\sum_{j=1}^p
b_i^j\cdot c_i^j\stackrel{\overline{I_s^{n+1}}[A]}{\equiv} \sum_{j=1}^p
E_s^n[b_i^j\cdot
c_i^j]=E_s^n[a_i]\overset{\eqref{nepreryvnost-E_s^n}}{\overset{A}{\underset{i\to\infty}{\longrightarrow}}}E_s^n[a],
$$
$$
\Downarrow
$$
$$
a\stackrel{\overline{I_s^{n+1}}[A]}{\equiv}E_s^n[a].
$$

When the formula \eqref{a-stackrel(overline(I_s^(n+1))[A])(equiv)E_s^n[a]} is proved, we can replace there  $a\in \overline{I_s^{n-1}}[A]$ by $a\in A$:
 \beq\label{a-stackrel(overline(I_s^(n+1))[A])(equiv)E_s^n[a],a-in-A}
a\stackrel{\overline{I_s^{n+1}}[A]}{\equiv} E_s^n[a],\qquad a\in A.
 \eeq
Indeed, for $n=0$ this is turned into a formula that we already noticed, \eqref{a-stackrel(I_s[A])(equiv)E_s^0[a]}. Further, if \eqref{a-stackrel(overline(I_s^(n+1))[A])(equiv)E_s^n[a],a-in-A} is proved for some $n-1$,
$$
a\stackrel{\overline{I_s^n}[A]}{\equiv} E_s^{n-1}[a],\qquad a\in A,
$$
then we can write it in the form
$$
a-E_s^{n-1}[a]\in \overline{I_s^n}[A],
$$
and by \eqref{a-stackrel(overline(I_s^(n+1))[A])(equiv)E_s^n[a]} we have:
$$
a-E_s^{n-1}[a]\stackrel{\overline{I_s^{n+1}}[A]}{\equiv}
E_s^n\big[a-E_s^{n-1}[a]\big]=\eqref{idempotentnost-E_s^n[f]}=E_s^n[a]-E_s^{n-1}[a],
$$
and thus,
$$
a\stackrel{\overline{I_s^{n+1}}[A]}{\equiv} E_s^n[a].
$$

4. Finally, \eqref{a-stackrel(overline(I_s^(n+1))[A])(equiv)E_s^n[a],a-in-A} implies that the quotient mapping $a\mapsto a+\overline{I_s^{n+1}}[A]$ can be factored through the mapping $a\mapsto E_s^n[a]$, and therefore the quotient algebra $\Jet_s^n[A]$ is isomorphic to the image $E_s^n[A]$ of the mapping $a\mapsto E_s^n[a]$:
$$
\Jet_s^n[A]\cong E_s^n[A].
$$
On the other hand, the algebra of jets $\Jet_s^n(M)$ on the manifold $M$ is also isomorphic to
$E_s^n[A]$:
$$
\Jet_s^n(M)\cong E_s^n[A].
$$
These equalities together mean (iv).

5. Let us prove that the conditions (ii), (iii), (iv) imply (v). We understood below that in the neighbourhood $U$ of the point $t\in M$ the bundles $\Jet_A^n[A]$ and $\Jet_{C^\infty(M)}^n[C^\infty(M)]$ coincide as sets. Suppose now that
$$
\xi_i\overset{\Jet_{C^\infty(M)}^n[C^\infty(M)]}{\underset{i\to\infty}{\longrightarrow}}\xi,
$$
Put $s_i=\pi(\xi_i)$, $s=\pi(\xi)$, and for each $i$ consider the linear combinations
$$
h_i=\sum_{|k|\le n}\lambda_k\cdot e_{s_i}^k: \qquad h_i=E_{s_i}[h_i],
$$
and
$$
h=\sum_{|k|\le n}\lambda_k\cdot e_s^k: \qquad h_i=E_s[h].
$$
The mapping of spectra $\Spec(A)\gets\Spec(C^\infty(M))=M$ is continuous, hence we have the following chain:
$$
s_i=\pi(\xi_i)\overset{M}{\underset{i\to\infty}{\longrightarrow}}\pi(\xi)=s\in U
$$
$$
\Downarrow
$$
$$
s_i=\pi(\xi_i)\overset{\Spec(A)}{\underset{i\to\infty}{\longrightarrow}}\pi(\xi)=s
$$
$$
\Downarrow
$$
$$
e_{s_i}^k=(e^1-e^1(s_i))^{k_1}\cdot...\cdot (e^d-e^d(s_i))^{k_d}\overset{A}{\underset{i\to\infty}{\longrightarrow}}
(e^1-e^1(s))^{k_1}\cdot...\cdot (e^d-e^d(s))^{k_d}=e_s^k
$$
$$
\Downarrow
$$
$$
h_i=\sum_{|k|\le n}\lambda_k\cdot e_{s_i}^k\overset{A}{\underset{i\to\infty}{\longrightarrow}}h=\sum_{|k|\le n}\lambda_k\cdot e_{s}^k
$$
$$
\Downarrow
$$
$$
\xi_i=\jet^n(h_i)\overset{\Jet_A^n[A]}{\underset{i\to\infty}{\longrightarrow}}\jet^n(h)=\xi.
$$

6. If the condition (i) of Theorem \ref{C^infty-obolochka-podalgebry-v-C^infty} holds, i.e. the mapping of spectra $\iota^{\Spec}:M\to\Spec(A)$ is an exact covering, then the morphism of bundles $\mu$ is a bijection, since it is a bijection between fibers and a bijection inside each bundle.
\epr

\blm\label{LM:dostatochnost-v-TH:C^infty-obolochka-podalgebry-v-C^infty} The conditions (i) and (ii) are sufficient for the morphism $\iota:A\to{\mathcal E}(M)$ to be a smooth envelope for $A$.
\elm
\bpr 1. First, let us check that under (i) and (ii) the morphism $\iota:A\to
{\mathcal E}(M)$ is a smooth extension. From the Nachbin theorem \ref{TH:Nachbin} it follows that $\iota:A\to {\mathcal E}(M)$ is a dense epimorphism. Let $B$ be a $C^*$-algebra and $\ph:A\to B[m]$ a differential morphism, i.e. such that the corresponding system of partial derivatives $\{D_k;\ k\in\N[m]\}$ consists of differential operators $D_k:A\to B$ of the orders $\ord D_k\le |k|$. The existence of the differential homomorphism $\ph'$ in the diagram
$$
 \xymatrix @R=2pc @C=1.2pc
 {
 A\ar[rr]^{\iota}\ar[dr]_{\ph} & & {\mathcal E}(M)\ar@{-->}[dl]^{\ph'}\\
  & B[m] &
 }
$$
is equivalent to the existence of the system of differential partial derivatives
$\{D_k';\ k\in\N[m]\}$ in diagrams
$$
 \xymatrix @R=2pc @C=1.2pc
 {
 A\ar[rr]^{\iota}\ar[dr]_{D_k} & & {\mathcal E}(M)\ar@{-->}[dl]^{D_k'}\\
  & B &
 }
$$
Put
$$
C=\overline{D_0(A)},\qquad
F=Z^1(D_0)=\{b\in B:\ \forall a\in A\quad [b,D_0(a)]=0\}=C^!.
$$
Since algebra $A$ is commutative, $C$ is also commutative, therefore
$$
C\subseteq C^!=F.
$$
This means that $D_0$ has values in $F$. On the other hand, by Theorem \ref{TH:differentsialnye-chast-proizv}, all operators $\{D_k;\ k\in\N[m], \ k>0\}$ also have values in $F$. Hence it is sufficient for us to find a system of differential partial derivatives $\{D_k';\
k\in\N[m]\}$, extending $D_k$ from $A$ to ${\mathcal E}(M)$ and having values in  $F$:
 \beq\label{prodolzhenie-D_k-na-C^infty}
 \xymatrix @R=2pc @C=1.2pc
 {
 A\ar[rr]^{\iota}\ar[dr]_{D_k} & & {\mathcal E}(M)\ar@{-->}[dl]^{D_k'}\\
  & F &
 }
 \eeq

2. First, we prove this for $k=0$:
 \beq\label{prodolzhenie-D_0-na-C^infty}
 \xymatrix @R=2pc @C=1.2pc
 {
 A\ar[rr]^{\iota}\ar[dr]_{D_0} & & {\mathcal E}(M)\ar@{-->}[dl]^{D_0'}\\
  & F &
 }
 \eeq
Let us embed the algebra ${\mathcal E}(M)$ into the algebra ${\mathcal C}(M)$. Then $\iota$ can be treated as a homomorphism from $A$ into ${\mathcal C}(M)$, and the condition (i) will mean that $\iota$ satisfies the premise of Theorem \ref{C-obolochka-podalgebry-v-C(M)}. Hence $D_0$ can be continuously extended to ${\mathcal C}(M)$. After that we restrict it to ${\mathcal E}(M)$, and we obtain $D_0'$.

3. When the homomorphism $D_0':{\mathcal E}(M)\to F$ is built, the algebra $F$ becomes a module over ${\mathcal E}(M)$. If we denote by $I_t^A$ and $I_t^{\mathcal E}$ the ideals in $A$ and ${\mathcal E}(M)$, consisting of functions vanishing in $t\in M$, then
$$
I_t^A\subseteq I_t^{\mathcal E}
$$
and by Lemma \ref{LM:o-plotnom-ideale} $I_t^A$ is dense in $I_t^{\mathcal E}$. Hence
$$
\overline{I_t^A\cdot F}=\overline{I_t^{\mathcal E}\cdot F}.
$$
Put $T=\Spec(C)$. Then we can think that $C={\mathcal C}(T)$. Consider two cases.
\bit{
\item[---] Suppose $t\notin T$. Then the functions from $I_t^{\mathcal E}$ approximate the characteristic function $\chi_T$ of the set $T$, hence the same is true for the functions from $I_t^A$:
$$
\chi_T\in\overline{I_t^A}=\overline{I_t^{\mathcal E}}\subseteq {\mathcal E}(M).
$$
Therefore
$$
\overline{I_t^A\cdot F}=\overline{I_t^{\mathcal E}\cdot F}=\chi_T\cdot F=1\cdot F=F,\qquad t\notin T.
$$
This means that the value bundles of the module $F$ over the algebras $A$ and ${\mathcal E}(M)$ in the points outside of $T$ vanish:
\beq\label{F/IF-t-notin-T}
F/\overline{I_t^A\cdot F}=F/\overline{I_t^{\mathcal E}\cdot F}=0,\qquad t\notin T=\Spec(C).
\eeq

\item[---] Suppose $t\in T$. Then, since ${\mathcal E}(M)$ (being restricted to $T$) is dense in $C={\mathcal C}(T)$, again by Lemma \ref{LM:o-plotnom-ideale} we obtain, that the ideal $I_t^{\mathcal E}$ is dense in the ideal  $I_t^C$ of functions from $C={\mathcal C}(T)$, vanishing in $t$:
$$
\overline{I_t^{\mathcal E}}=I_t^C.
$$
Therefore,
$$
\overline{I_t^A\cdot F}=\overline{I_t^{\mathcal E}\cdot F}=\overline{I_t^C\cdot F},\qquad t\in T,
$$
and this means that over each point $t\in T=\Spec(C)$ the value fibers of the module $F$ over the algebras $A$, ${\mathcal E}(M)$ and $C$ coincide:
\beq\label{F/IF-t-in-T}
F/\overline{I_t^A\cdot F}=F/\overline{I_t^{\mathcal E}\cdot F}=F/\overline{I_t^C\cdot F},\qquad t\in T=\Spec(C).
\eeq
}\eit
Besides this, due to condition (i), the compact set $T$ is continuously (injectively and with the same topology) embedded into $\Spec(A)$ and into $M$. Together this gives a description of the value bundles $\Jet^0_A[F]$ and $\Jet^0_{{\mathcal E}(M)}[F]$: they are trivial outside of the compact set $T$, while on $T$ they coincide with the value bundle $\Jet^0_C[F]$. We can conclude that the natural mapping of the jet bundles
$$
\Jet^0_A[F]=\bigsqcup_{t\in M} (F/\overline{I_t^A\cdot F})^\vartriangle\overset{\lambda}{\longleftarrow} \bigsqcup_{t\in M} (F/\overline{I_t^{\mathcal E}\cdot F})^\vartriangle=\Jet^0_{{\mathcal E}(M)}[F]
$$
is an isomorphism.

4. To build $D'_k$ for $k>0$, recall that by Theorem
\ref{TH:diff-oper->rassl-struuj} each differential operator $D_k$ generates a morphism of jet bundles
$\jet_n[D_k]:\Jet_A^n[A]\to \Jet_A^0(F)=\pi^0_A F$, where $n=|k|$, satisfying the identity
$$
\jet^0(D_ka)=\jet_n[D_k]\circ \jet^n(a),\qquad a\in A.
$$
By condition (v) of Lemma \ref{LM-3-dlya-TH:C^infty-obolochka-podalgebry-v-C^infty}, the jet bundles of the algebras  $A$ and ${\mathcal E}(M)$ are connected to each other by a natural morphism $\mu:\Jet_A^n[A]\gets \Jet^n_{{\mathcal E}(M)}[{\mathcal E}(M)]$, which in addition is a bijection. We obtain the diagram
$$
 \xymatrix @R=2pc @C=4pc
 {
 \Jet^n_A[A]\ar[dd]_{\jet_n[D_k]}\ar[dr]^{\pi^n_{A,A}} & & & \Jet^n_{{\mathcal E}(M)}[{\mathcal E}(M)]\ar@{-->}[dd]^{\nu=\jet_n[D_k']}\ar[lll]_{\mu}\ar[dl]_{\pi^n_{{\mathcal E}(M),{\mathcal E}(M)}}  \\
  &  \Spec(A) &  M\ar[l]_{\iota^{\Spec}} &  \\
 \Jet^0_A[F]\ar[rrr]^{\lambda^{-1}}\ar[ur]_{\pi^0_{A,F}} & & & \Jet^0_{{\mathcal E}(M)}[F]  \ar[ul]^{\pi^0_{{\mathcal E}(M),F}}
 }
$$
Let us consider the composition (i.e. the dotted arrow)
$$
\nu=\lambda^{-1}\circ\jet_n[D_k]\circ\mu
$$
It is a morphism of bundles, hence by Lemma \ref{TH:diff-oper<-morfizm-rassl-struuj} it generates a differential operator $D_k':{\mathcal E}(M)\to F$ (over the algebra ${\mathcal E}(M)$), such that
$$
\jet^0(D_k'u)=\jet_n[D_k]\circ \jet^n(u),\qquad u\in {\mathcal E}(M).
$$
For each $a\in A$ we have
\beq\label{PROOF:TH:C^infty-1}
\jet^0\big(D_k'\iota(a)\big)=\jet_n[D_k]\circ \jet^n(\iota(a))=\jet_n[D_k]\circ
\jet^n(a)=\jet^0(D_ka).
\eeq
Note that by Theorem \ref{TH:B-cong-Sec(val_AB)} the mapping $\jet^0=v:F\to\Sec(\pi^0_AF)=\Sec(\Jet^0_AF)$, that turns $F$ into the algebra of continuous sections of the jet bundles $\pi^0_AF: \Jet^0_AF\to\Spec(A)$ over $A$, is an isomorphism of $C^*$-algebras:
$$
F\cong\Sec(\pi^0_AF).
$$
Hence we can apply to \eqref{PROOF:TH:C^infty-1} the operator, inverse to $\jet^0$, and we obtain the equality
$$
D_k'\iota(a)=D_ka.
$$
In other words, $D_k'$ extends $D_k$ in \eqref{prodolzhenie-D_k-na-C^infty}. Besides this, since $\iota$ maps
$A$ densely into ${\mathcal E}(M)$, the conditions
\eqref{DEF:sist-chastn-proizv-*}-\eqref{DEF:sist-chastn-proizv-2} are transferred from $D_k$ to $D_k'$. Thus, the family $\{D_k',\ k\in\N[m]\}$ is a system of partial derivatives on ${\mathcal E}(M)$.

5. Now let us verify that from (i) and (ii) it follows, that the extension $\iota:A\to
{\mathcal E}(M)$ is a smooth envelope. Let $\sigma:A\to C$ be another smooth extension. Take a local chart $\ph:U\to V$, where $U\subseteq M$, $V\subseteq\R^d$, and let $K\subseteq U$ be a compact set that coincides with the closure of its interior: $\overline{\Int(K)}=K$. Then the operators from the example \ref{EX:sist-chast-proizv-na-C-infty(M)}
$$
D_k:{\mathcal E}(M)\to\mathcal{C}(K)\quad\Big|\quad
D_k(a)=\frac{\partial^{\abs{k}}(a\circ\ph^{-1})}{\partial t_1^{k_1}...\partial
t_d^{k_d}}\circ\ph,\qquad a\in A,\quad k\in\N[m],
$$
form a system of partial derivatives from $A$ into $\mathcal{C}(K)$. Suppose $D:A\to
\mathcal{C}(K)[m]$ is a differential homomorphism, generated by $\{D_k\}$. Since $\sigma:A\to C$ is a smooth extension, the homomorphism $D:A\to C(K)[[d]]$ is uniquely extended to some homomorphism $D':C\to C(K)[[d]]$:
 \beq\label{DIAGR:env-A=E(M)}
 \xymatrix @R=2pc @C=1.2pc
 {
 A\ar[rr]^{\sigma}\ar[dr]_{D} & & C \ar@{-->}[dl]^{D'}\\
  & \mathcal{C}(K)[m] &
 }
 \eeq
If we return back to the system of partial derivatives $D_k$, we obtain for each index $k\in\N[m]$ a diagram
 $$
  \xymatrix @R=2pc @C=1.2pc
 {
 A\ar[rr]^{\sigma}\ar[dr]_{D_k} & & C \ar@{-->}[dl]^{D_k'}\\
  & \mathcal{C}(K) &
 }
 $$

Take an arbitrary element $c\in C$. Since $\sigma:A\to C$ is a dense epimorphism, there is a net of elements $a_i\in A$ such that
$$
\sigma(a_i)\overset{C}{\underset{i\to\infty}{\longrightarrow}} c.
$$
For each index $k\in\N[m]$ we have
 \beq\label{D_k(a_i)->D_k'(b)}
D_k(a_i)=D_k'(\sigma(a_i))\overset{\mathcal{C}(K)}{\underset{i\to\infty}{\longrightarrow}}
D_k'(c).
 \eeq
Now consider a smooth curve in $K$, i.e. a smooth mapping $\gamma:[0,1]\to K$. For each index $k\in\N[m]$ of the order $|k|=1$ and for each point $t\in[0,1]$ the symbol $\gamma^k(t)$ denotes the $k$-th component of the derivative $\gamma'(t)$ in the expansion by local coordinates on $U$. For each function $a\in A$ by the Newton-Leibnitz theorem we have
$$
D_0(a)(\gamma(1))-D_0(a)(\gamma(0))=\sum_{|k|=1}\int_0^1\gamma^k(t)\cdot
D_k(a)(\gamma(t))\d t.
$$
Together with \eqref{D_k(a_i)->D_k'(b)} this gives
 \begin{multline*}
D_0'(c)(\gamma(1))-D_0'(c)(\gamma(0)) \underset{\infty\gets i}{\longleftarrow}
D_0(a_i)(\gamma(1))-D_0(a_i)(\gamma(0))=\\
=\sum_{|k|=1}\int_0^1\gamma^k(t)\cdot
D_k(a_i)(\gamma(t))\d t \underset{i\to\infty}{\longrightarrow}
\sum_{|k|=1}\int_0^1\gamma^k(t)\cdot D_k'(c)(\gamma(t))\d t,
 \end{multline*}
and thus,
$$
D_0'(c)(\gamma(1))-D_0'(c)(\gamma(0))=\sum_{|k|=1}\int_0^1\gamma^k(t)\cdot
D_k'(c)(\gamma(t))\d t
$$
This connection between the function $D_0'(c)\in {\mathcal C}(K)$ and the functions $D_k'(c)\in
{\mathcal C}(K)$, $|k|=1$, means that $D_0'(c)$ is continuously differentiable on $K$, and its partial derivatives in the chosen local coordinates are the functions $D_k'(c)$, $|k|=1$.

After that we take arbitrary $D_k'(c)$, $|k|=1$, and consider indices of order 2. The same reasoning shows that $D_k'(c)$ is also continuously differentiable. The induction over the indices gives that all functions $D_k'(c)$ are smooth and connected to each other as partial derivatives of the function $D_0'(c)$ (with respect to the chosen local coordinates).

If we change the compact set $K\subset U$ and the open set $U\subseteq M$, then the arising smooth functions  $D_0'(c)$ on $K$ are compatible to each other, since on the intersections of their domains they coincide. This means that there is a common smooth function $\iota'(c):M\to\C$ such that its restriction to each compact set $K$ coincides with the corresponding function $D_0'(c)$:
$$
\iota'(c)\big|_K=D_0'(c),\qquad K\subset U\subseteq M.
$$
and the partial derivatives with respect to the chosen local coordinates coincide with the values of the operators  $D_k'$ on $c$. In other words, there is a mapping $\iota':C\to{\mathcal E}(M)$ (by construction this is a homomorphism of algebras), such that the following diagram specifying \eqref{DIAGR:env-A=E(M)} is commutative:
 \beq\label{DIAGR:env-A=E(M)-*}
 \xymatrix @R=2pc @C=1.2pc
 {
 A\ar[rr]^{\sigma}\ar@/_4ex/[ddr]_{D}\ar[dr]_{\iota} & & C \ar@/^4ex/[ddl]^{D'}\ar@{-->}[dl]^{\iota'}\\
  & {\mathcal E}(M)\ar[d]_{\lambda} & \\
  & \mathcal{C}(K)[m] &
 }
 \eeq
(here the mapping $\lambda$ is the expansion of the smooth function in the Taylor polynomial on the compact set  $K\subset U\subseteq M$ with respect to the chosen local coordinates on $U$). From this diagram we have that $\iota'$ is continuous, since if $c_i\to c$, then this condition is preserved under the action of each operator $D_k'$, i.e.  $D_k'(c_i)\to D_k'(c)$, and this is exactly the convergence in the space ${\mathcal E}(M)$.

The upper inner triangle in \eqref{DIAGR:env-A=E(M)-*} is exactly the diagram that we need.
 \epr

\paragraph{Counterexamples.}

The following example shows that the weakening of the condition (i) in Theorem \ref{C^infty-obolochka-podalgebry-v-C^infty} makes this proposition false.

\bex\label{EX:E(8)} {\it
There is a dense involutive subalgebra $A$ in ${\mathcal E}(\R)$ such that
\bit{
\item[(i)] the dual mapping of spectra $\iota^{\Spec}:\Spec(A)\gets M$ is a bijection (but not a covering),

\item[(ii)] for each point $s\in\R$ the natural mapping of tangent spaces $T_s[A]\gets T_s(\R)$ is an isomorphism  (of finite dimensional vector spaces),

\item[(iii)] the smooth envelope of $A$ is not isomorphic to ${\mathcal E}(\R)$:
$$
\Env_{\mathcal E}A\not\cong {\mathcal E}(\R)
$$
}\eit
}
\eex
\bpr This is a modification of Example \ref{EX:C(8)}. Consider the circle $\T=\R/\Z$ and the algebra ${\mathcal E}(\T)$ of smooth functions on it. In the Cartesian square ${\mathcal E}(\T)^2$ consider the subalgebra ${\mathcal E}(8)$, consisting of pairs of functions $(u,v)\in{\mathcal E}(\T)^2$ with the same derivatives in the point 0:
$$
(u,v)\in{\mathcal E}(8)\quad\Longleftrightarrow\quad u,v\in{\mathcal E}(\T)\quad\&\quad \forall k\ge 0\quad u^{(k)}(0)=v^{(k)}(0).
$$
This is a closed subalgebra in ${\mathcal E}(\T)^2$, and its spectrum is the space which in Example \ref{EX:C(8)} was denoted by the symbol $8$, and which can be described as the result of glueing of two copies of the circle $\T$ in the point 0:
\beq\label{Spec(E(8))}
\Spec{\mathcal E}(8)=8.
\eeq
Each point $s\in 8$ either belongs to the first copy of the circle $\T$, which we denote by $\T_1$, or to the second, $\T_2$ (or, if $s=0$, then we consider that $s$ belongs to both circles $\T_1$ and $\T_2$). Let us show that in each point $s\in 8$ the tangent space to ${\mathcal E}(8)$ is isomorphic to the tangent space to ${\mathcal E}(\T)$ in this point:
\beq\label{T_s(E(8))}
T_s({\mathcal E}(8))\cong\begin{cases} T_s({\mathcal E}(\T_1)), & s\in\T_1\\  T_s({\mathcal E}(\T_2)), & s\in\T_2 \end{cases},\qquad s\in 8.
\eeq

Take a function $\ph\in{\mathcal E}(\T)$, which in some neighbourhood $U$ of 0 in $\R$ is the inverse mapping for the projection $\pi:\R\to\R/\Z=\T$:
$$
\ph(\pi(t))=t,\qquad t\in U.
$$
Let $I_s(\T)$ be the ideal in ${\mathcal E}(\T)$, consisting of functions, vanishing in $s$. Then by the Hadamard lemma \cite{Petrovsky} each function $u\in{\mathcal E}(\T)$ differs from $u(s)+u'(s)\ph$ by the function from $I_s^2(\T)$:
\beq\label{Hadamard-dlya-T}
u\overset{I_0^2(\T)}{\equiv} u(0)+u'(0)\ph.
\eeq
Consider two cases.
\bit{

\item[1)] Take a point $s\in 8\setminus\{0\}$. As an element of the space 8, $s$  either belongs to $\T_1$, or to  $\T_2$. Let $\ph_s$ denote the shift of the function $\ph$ by $s$ in $\T$:
$$
\ph_s(t)=\ph(t-s),\quad t\in\T,
$$
and let $\psi_s$ be a function with zero germ in 0, and the same germ with $\ph_s$ in $s$:
$$
\psi_s(t)\equiv\begin{cases} 0, & \mod 0\\ \ph_s, & \mod s\end{cases}.
$$
From \eqref{Hadamard-dlya-T} we have
$$
(u,v)\overset{I_s^2(\T)}{\equiv} \begin{cases}u(s)+u'(s)\cdot(\psi_s,0),& s\in\T_1\\ v(s)+v'(s)\cdot(0,\psi_s),& s\in\T_2 \end{cases},
$$
As a corollary, the action of the tangent vector $\tau\in T_s({\mathcal E}(8))$ at the element $(u,v)$ is described by the formula
$$
\tau(u,v)=\begin{cases}u(s)+u'(s)\cdot\tau(\psi_s,0),& s\in\T_1\\ v(s)+v'(s)\cdot\tau(0,\psi_s),& s\in\T_2 \end{cases}
$$
This implies that (for $s\ne 0$) the tangent space $T_s({\mathcal E}(8))$ is isomorphic to $\R$ and to $T_s(\T)$.
$$
T_s({\mathcal E}(8))\cong\R\cong T_s({\mathcal E}(\T)).
$$

\item[2)] Consider the point $s=0$. Then from \eqref{Hadamard-dlya-T} we have
$$
(u,v)\overset{I_0^2(\T)}{\equiv} u(s)\cdot(1,1)+u'(s)\cdot(\ph,\ph)=v(s)\cdot(1,1)+v'(s)\cdot(\ph,\ph),
$$
so the action of the tangent vector $\tau\in T_0({\mathcal E}(8))$ on this element is described by the formula
$$
\tau(u,v)=u'(s)\cdot\tau(\ph,\ph)=v'(s)\cdot\tau(\ph,\ph).
$$
This implies that (for $s=0$) the tangent space $T_0({\mathcal E}(8))$ is isomorphic to $\R$ and to $T_0(\T)$.
$$
T_0({\mathcal E}(8))\cong\R\cong T_0({\mathcal E}(\T)).
$$

}\eit

When \eqref{Spec(E(8))} and \eqref{T_s(E(8))} are proved, let us take a smooth bijective mapping $\omega:\R\to 8$ (such a mapping always exists, and it is not unique). It generates a homomorphism of algebras ${\mathcal E}(\R)\gets{\mathcal E}(8)$. If we denote by $A$ the image of this homomorphism (with the algebraic operations and the topology induced from ${\mathcal E}(8)$), we obtain the algebra with the properties (i) and (ii), and therefore with (iii) as well.
\epr

The following two examples show the mutual independence of the conditions (i) and (ii) in Theorem \ref{C^infty-obolochka-podalgebry-v-C^infty}.

\bex
There exists a closed subalgebra $A$ in the algebra ${\mathcal E}(\R)$ with the following properties:
 \bit{
\item[(i)] the spectrum of $A$ does not coincide with $\R$ as a set,
$$
\Spec(A)\ne \R,
$$

\item[(ii)] for each point $s\in\R$ the natural mapping of tangent spaces
 $$
T_s[A]\gets T_s(\R)
 $$
is an isomorphism (of vector spaces),
 }\eit
\eex \bpr The algebra of periodic smooth functions on $\R$ with a given period, say, $1$, has these properties. \epr

\bex\label{EX:kontrprimer-podalgebra-E(R)}
There exists a closed subalgebra $A$ in the algebra ${\mathcal E}(\R)$ with the following properties:
 \bit{
\item[(i)] the spectrum of $A$ coincides with $\R$,
$$
\Spec(A)=\R,
$$

\item[(ii)] in the point $s=0$ the tangent space to $A$ is trivial:
$$
T_0[A]=0
$$
 }\eit
\eex

To prove this we need the following

\blm\label{LM:ploskaya-v-nule-funktsija} Suppose a function $x$ in ${\mathcal E}(\R)$ has zero derivatives of all orders in the point 0,
$$
\forall n\ge 0\qquad x^n(0)=0.
$$
Then it can be approximated in ${\mathcal E}(\R)$ by functions of zero germs in the point $0$.
\elm
\bpr
Take $n\in\N$, $T>1$, $\e>0$ and find $\delta>0$ such that
$$
\forall k\in\{0,...,n\}\qquad \sup_{|t|\le\delta}\abs{x^{(k)}(t)}\le\e.
$$
Take a function $\eta\in{\mathcal E}(\R)$,
$$
0\le\eta(s)\le 1,\qquad \eta(s)=\begin{cases}0,& |s|\le\frac{\delta}{2}\\ 1,& |s|\ge\delta \end{cases}
$$
and put
$$
y_n(t)=\int_0^t \eta(s)\cdot x^{(n+1)}(s)\d s,\quad y_{n-1}(t)=\int_0^t y_n(s)\d s,\quad ...,\quad y_0(t)=\int_0^t y_1(s)\d s.
$$
Then
\begin{multline*}
\sup_{|t|\le T}\abs{x^n(t)-y_n(t)}=\sup_{|t|\le T}\abs{\int_0^t x^{(n+1)}(s)\d s-\int_0^t \eta(s)\cdot x^{(n+1)}(s)\d s}\le\\
\le\sup_{|t|\le T}\int_0^t (1-\eta(s))\cdot \abs{x^n(s)}\d s\le \sup_{|s|\le
T}\abs{x^n(s)}\cdot \int_0^t (1-\eta(s))\d s\le \sup_{|s|\le
T}\abs{x^{(n+1)}(s)}\cdot\delta.
\end{multline*}
$$
\Downarrow
$$
\begin{multline*}
\sup_{|t|\le T}\abs{x^{(n-1)}(t)-y_{n-1}(t)}=\sup_{|t|\le T}\abs{\int_0^t
x^n(s)\d s-\int_0^t y_n(s)\d s}\le\\ \le\sup_{|t|\le T}\int_0^t
\abs{x^n(s)-y_{n}(s)}\d s\le \sup_{|s|\le T}\abs{x^{(n+1)}(s)}\cdot
 \delta\cdot T.
\end{multline*}
$$
\Downarrow
$$
$$
...
$$
$$
\Downarrow
$$
\begin{multline*}
\sup_{|t|\le T}\abs{x^{(0)}(t)-y_0(t)}=\sup_{|t|\le T}\abs{\int_0^t x^1(s)\d
s-\int_0^t y_1(s)\d s}\le\\ \le\sup_{|t|\le T}\int_0^t \abs{x^1(s)-y_{1}(s)}\d
s\le \sup_{|s|\le T}\abs{x^{(n+1)}(s)}\cdot
 \delta\cdot T^n.
\end{multline*}
We see that for given $n\in\N$, $T>1$, $\e>0$ one can choose $\delta>0$ such that
$$
\delta <\frac{\e}{T^n\cdot \sup_{|s|\le T}\abs{x^{(n+1)}(s)}},
$$
and then the function $y_0$ will differ from $x$ less than $\e$ uniformly on all derivatives of order less than $n$ on the interval $[-T,T]$:
$$
\max_{0\le k\le n}\sup_{|t|\le T}\abs{x^{(k)}(t)-y^{(k)}(t)}=\max_{0\le k\le n}\sup_{|t|\le T}\abs{x^{(k)}(t)-y_k(t)}\le T^n\cdot\sup_{|s|\le T}\abs{x^{(n+1)}(s)}<\e.
$$
\epr

\bpr[Proof of Example \ref{EX:kontrprimer-podalgebra-E(R)}.] This is the algebra of smooth functions which in the point $s=0$ have vanishing derivatives of every positive order:
$$
A=\{x\in {\mathcal E}(M):\quad \forall n\ge 1\quad x^n(0)=0\}.
$$
1. Let us first show that $\Spec(A)=\R$. Take $s\in\Spec(A)$, i.e. $s$ is an involutive, continuous, linear, multiplicative, preserving unit functional on $A$:
$$
s(x^\bullet)=\overline{s(x)},\quad s(\lambda\cdot x+y)=\lambda\cdot s(x)+s(y),\quad s(x\cdot y)=s(x)\cdot s(y),\quad s(1)=1.
$$
Consider its kernel $\Ker s=\{x\in A:\ s(x)=0\}$. Let us show that there is a point $t\in\R$ such that if we consider it as a functional on $A$, then its kernel contains $\Ker(s)$:
\beq\label{Ker s-subseteq-Ker t}
\Ker s\subseteq\Ker t=\{x\in A:\ x(t)=0\}.
\eeq
Suppose that this is not true, i.e.
\beq\label{forall-t-in-R-exists-x_t-in-Ker(s)-x_t(t)-ne-0}
\forall t\in\R\quad \exists x_t\in\Ker s\quad x_t(t)\ne 0.
\eeq
Consider the family of sets
$$
U_t=\{ r\in\R: \ x_t(r)\ne 0\}.
$$
They are open in $\R$, and the condition \eqref{forall-t-in-R-exists-x_t-in-Ker(s)-x_t(t)-ne-0} means that they form a covering of $\R$. Hence we can refine in it a smooth locally finite partition of unity:
$$
\eta_t\in {\mathcal E}(M),\quad \supp\eta_t\subseteq U_t,\quad \eta_t\ge 0,\quad \sum_{t\in\R}\eta_t=1.
$$
We can choose this partition such that the following supplementary conditions hold:
$$
\forall n\ge 1\quad \eta_0^n(0)=0\quad\&\quad \forall t\ne 0\quad
0\notin\supp\eta_t.
$$
They automatically imply that all functions $\eta_t$ belong to $A$, and the series $\sum_{t\in\R}\eta_t\cdot x_t\cdot x_t^\bullet$ converges in $A$ (in the topology induced from ${\mathcal E}(\R)$). Denote its sum by $y$:
$$
y=\sum_{t\in\R}\eta_t\cdot x_t\cdot x_t^\bullet,
$$
and note that since all $x_t$ belong to the closed ideal $\Ker s$, the element $y$ also lies in it:
$$
y\in\Ker s.
$$
On the other hand, the function $y$ nowhere vanishes. Indeed, for each point $r\in\R$ there is a function  $\eta_{t_r}$, non vanishing in this point, hence $\eta_{t_r}(r)>0$. Herewith, the condition $\supp\eta_{t_r}\subseteq U_{t_r}=\{q\in\R:\ x_{t_r}(q)\ne 0\}$ imply the condition $\abs{x_{t_r}(r)}>0$, and we get
$$
y(r)=\sum_{t\in\R}\eta_t(r)\cdot \abs{x_t(r)}^2\ge \underbrace{\eta_{t_r}(r)}_{\tiny\begin{matrix}
\text{\rotatebox{90}{$<$}}\\ 0 \end{matrix}}\cdot \underbrace{\abs{x_{t_r}(r)}^2}_{\tiny\begin{matrix}
\text{\rotatebox{90}{$<$}}\\ 0 \end{matrix}}>0.
$$
Thus, $y$ is invertible, as an element of the algebra ${\mathcal E}(\R)$. In addition, by construction, all the derivatives of the function $y$ in the point $0$ vanish, therefore its inverse function $\frac{1}{y}$ also has vanishing derivatives in $0$. Hence,
$$
\frac{1}{y}\in A.
$$
I.e. $y$ is invertible as an element of $A$.

Thus, the function $y$ belongs to the ideal $\Ker s$ of $A$, and on the other hand it is invertible in $A$. Hence, $\Ker s=A$, and this contradicts to the condition $s(1)=1$. Thus, our assumption \eqref{forall-t-in-R-exists-x_t-in-Ker(s)-x_t(t)-ne-0} turned out to be not true, and for some $t\in\R$ we have \eqref{Ker s-subseteq-Ker t}. We see that there are two non-zero functionals $s$ and $t$ on $A$, such that, first,  $\Ker s\subseteq\Ker t$, and, second, $s(1)=t(1)=1$. This is possible only if they coincide everywhere on $A$: $s=t$.

2. Let us now prove (ii). It is sufficient to show that the kernel $I_0[A]=\{x\in A:\ x(0)=0\}$ of the character $x\mapsto x(0)$ has the following property:
\beq\label{overline(I_0^2)[A]=I_0[A]}
\overline{I_0^2}[A]=I_0[A].
\eeq
-- then \eqref{T_s[A]-cong-Real(I_s/overline(I_s^2)^star)} will imply
$$
T_0[A]=\Real\big(I_0[A]/\overline{I_0^2}[A]\big)^\triangledown=0.
$$
Indeed, take $x\in I_0[A]$, i.e. $x$ is a function from ${\mathcal E}(\R)$, that vanishes in 0 with all its derivatives,
$$
\forall n\ge 0\qquad x^n(0)=0.
$$
By Lemma \ref{LM:ploskaya-v-nule-funktsija} it can be approximated in ${\mathcal E}(\R)$ by functions with the zero germ in the point $0$. On the other hand, each function $y$ with a zero germ in $0$ belongs to the ideal $I_0^2[A]$  (we can multiply $y$ by another function with a zero germ in 0, but having the value 1 on the support of $y$). We see, that $x$ can be approximated in ${\mathcal E}(\R)$ (hence, in $A$) by functions from $I_0^2[A]$.
\epr

The following example, suggested by M.~B\"achtold, shows that in Lemma \ref{LM:Ker-rho_t^(T^star)=0} the injectivity of the mapping $\iota_t^{\C T^\star}$ is necessary.

\bex
If for a homomorphism of involutive stereotype algebras $\iota:A\to B$ in some point $t\in\Spec(B)$ the morphism of cotangent spaces $\iota_t^{\C T^\star}:{\C}T_{t\circ\iota}^\star[A]\to{\C}T_t^\star[B]$ is not injective,
$$
\Ker\iota_t^{T^\star}\ne 0
$$
then formula \eqref{overline(I_s^2)[A]=rho^(-1)(overline(I_t^2)[B])} is not necessarily true:
$$
\overline{I_{t\circ\iota}^2}[A]\ne \iota^{-1}\Big(\overline{I_t^2}[B]\Big).
$$
 \eex
\bpr
Consider the algebra ${\mathcal E}(\R)$ of smooth functions on the line $\R$ and a subalgebra $A$ in it, that is generated (as a pure algebra) by functions $x^2$ and $x\cdot e^x$. Thus, $A$ is a linear span in ${\mathcal E}(\R)$ of the functions
$$
x^{2p+q}\cdot e^{qx},\qquad p,q\in\N=\{0,1,2,...\}
$$
(these functions form an algebraic basis in $A$). The algebra $A$ separates the points of $\R$
$$
s\ne t\in\R\qquad\Longrightarrow\qquad \exists a\in A\quad a(s)\ne a(t),
$$
and the tangent vectors of ${\mathcal E}(\R)$
$$
\forall s\in\R\quad\forall\tau\in T_s[{\mathcal E}(\R)]\quad \exists a\in A\quad \tau(a)\ne 0.
$$
Hence, by the Nachbin theorem \ref{TH:Nachbin}, $A$ is dense in ${\mathcal E}(\R)$. If we endow $A$ with the strongest locally convex topology, $A$ becomes a stereotype algebra (by Example \ref{EX:siln-loc-vyp-topol}), and the natural embedding into ${\mathcal E}(\R)$ is an epimorphism of stereotype algebras $\iota:A\to {\mathcal E}(\R)$. Consider the ideal $I_0[A]$ of functions in $A$, vanishing in the point $0\in\R$. It is a linear span of the functions
$$
x^{2p+q}\cdot e^{qx},\qquad p,q\in\N=\{0,1,2,...\},\qquad p+q\ne 0.
$$
Its square $I_0^2[A]$ is a linear span of the functions
$$
x^{2p+q}\cdot e^{qx},\qquad p,q\in\N=\{0,1,2,...\},\qquad p>0.\qquad q>0.
$$
In this list the number of functions is less than in the list for $I_0[A]$: two functions, $x^2$ and $x\cdot e^x$, are absent. This means that the quotient space $\C T_0^\star[A]=I_0[A]/I_0^2[A]$ has the dimension 2. The same is true for its real part:
$$
\dim_{\R}T_0^\star[A]=2.
$$
Hence the mapping of the cotangent spaces $\iota_0^{T^\star}:T_0^\star[A]\to T_0^\star[{\mathcal E}(\R)]$ cannot be injective. At the same time the function $x^2$, that lies in $A$ and in $I_0^2[{\mathcal E}(\R)]$ (hence, in  $\iota^{-1}\Big(\overline{I_t^2}[B]\Big)$ as well), does not belong to $I_0^2[A]$.
\epr

\paragraph{Smooth envelope of the algebra $k(G)$ on a compact Lie group $G$.}

Recall the algebra $k(G)$ of matrix elements defined on page \pageref{DEF:k(G)}.

\btm The smooth envelope of the algebra $k(G)$ of matrix elements on a compact Lie group $G$ coincides with the algebra ${\mathcal E}(G)$ of smooth functions on $G$:
\beq
\Env_{\mathcal E}k(G)={\mathcal E}(G).
\eeq
\etm
\bpr
The algebra $k(G)$ is embedded into ${\mathcal E}(G)$, their spectra coincide, and by Corollary  \ref{COR:kasat-prostr-k-Trig(G)}, the tangent spaces also coincide in each point. Then Theorem \ref{C^infty-obolochka-podalgebry-v-C^infty} works.
\epr

\subsection{Smooth envelopes of group algebras}

\paragraph{Coincidence of the smooth envelopes of ${\mathcal C}^\star(G)$ and ${\mathcal E}^\star(G)$.}

We have already noticed on page \pageref{Env_C-C^star(G)-Env_C-E^star(G)} that a real Lie group $G$ has two group algebras: the algebra ${\mathcal C}^\star(G)$ of measures with compact support and the algebra ${\mathcal E}^\star(G)$ of distributions with compact support (these constructions are described in detail in \cite{Akbarov}). If we denote by $\lambda$ the natural inclusion ${\mathcal E}(G)\subseteq{\mathcal C}(G)$, then we obtain the diagram
\beq\label{Env_E-C^star(G)-Env_E-E^star(G)}
 \xymatrix @R=2.5pc @C=4pc
 {
{\mathcal C}^\star(G)\ar[r]^{\lambda^\star}\ar[d]_{\env_{\mathcal E}{\mathcal C}^\star(G)} & {\mathcal E}^\star(G)\ar[d]^{\env_{\mathcal E}{\mathcal E}^\star(G)}
\\
\Env_{\mathcal E}{\mathcal C}^\star(G)\ar[r]^{\Env_{\mathcal C}(\lambda^\star)} & \Env_{\mathcal E}{\mathcal E}^\star(G)
 }
\eeq

By analogy with Theorem \ref{TH:Env_C-C^star(G)=Env_C-E^star(G)} one can prove

\btm\label{TH:Env_E-C^star(G)=Env_E-E^star(G)}
For each real Lie group $G$ the smooth envelopes of group algebras ${\mathcal C}^\star(G)$ and ${\mathcal E}^\star(G)$ coincide:
\beq\label{Env_E-C^star(G)=Env_E-E^star(G)}
\Env_{\mathcal E}{\mathcal C}^\star(G)=\Env_{\mathcal E}{\mathcal E}^\star(G).
\eeq
\etm

\paragraph{Fourier transform on a commutative Lie group.}
Recall the Fourier transform on a commutative locally compact group $G$, defined by formula \eqref{Fourier-transform}. If $G$ is a compactly generated commutative Lie group, then the dual group $\widehat{G}$ is also a (compactly generated) Lie group. Hence in this case the algebras of smooth function ${\mathcal E}(G)$, ${\mathcal E}(\widehat{G})$ are defined, as well as the algebras of distributions ${\mathcal E}^\star(G)$, ${\mathcal E}^\star(\widehat{G})$. The very same formula \eqref{Fourier-transform} defines a mapping
$$
{\mathcal F}_G:{\mathcal E}^\star(G)\to {\mathcal C}(\widehat{G}).
$$
which we also call a {\it Fourier transform} on the group $G$.

\blm For a compactly generated commutative Lie group $G$ the Fourier transform continuously maps the algebra of measures ${\mathcal C}^\star(G)$ and the algebra of distributions ${\mathcal E}^\star(G)$ into the algebra ${\mathcal E}(\widehat{G})$ of smooth functions on $\widehat{G}$.
$$
{\mathcal F}_G:{\mathcal C}^\star(G)\to {\mathcal E}(\widehat{G}),\qquad {\mathcal F}_G:{\mathcal E}^\star(G)\to {\mathcal E}(\widehat{G}).
$$
\elm
\bpr
Take $\alpha\in{\mathcal E}^\star(G)$. Each real character (a continuous homomorphism) $r:G\to\R$ defines a one-parametric subgroup (a continuous homomorphism)
$$
h:\R\to\widehat{G}\quad\Big|\quad h(t)=e^{itr},
$$
and vice versa, each one-parametric subgroup $h:\R\to\widehat{G}$ has this representation \cite[(24.43)]{Hewitt-Ross}. Hence
\begin{multline*}
\frac{{\mathcal F}_G(\alpha)(\chi\cdot h(t))-{\mathcal F}_G(\alpha)(\chi)}{t}=
\frac{\alpha(\chi\cdot e^{itr})-\alpha(\chi)}{t}=\alpha\left(\frac{\chi\cdot e^{itr}-\chi}{t}\right)=
\alpha\left(\chi\cdot \frac{e^{itr}-1}{t}\right)\underset{t\to 0}{\longrightarrow}
\alpha(\chi\cdot ir)
\end{multline*}
and from this relation we see that ${\mathcal F}_G(\alpha)$ is continuously differentiable along each one-parameter subgroup in $\widehat{G}$. Similarly one can check the continuous differentiability of arbitrary order, and we get ${\mathcal F}_G(\alpha)\in{\mathcal E}(\widehat{G})$, and the partial derivatives are expressed in the formulas
$$
\partial_{h_k}...\partial_{h_1}{\mathcal F}_G(\alpha)(\chi)=\alpha(\chi\cdot (ir_1)\cdot...\cdot (ir_k)).
$$
Further, if $\alpha_\nu\to 0$ in ${\mathcal C}^\star(G)$, then for each compact set $K\subseteq\widehat{G}$ and for each  $r_1,...,r_k$ the set $\{\chi\cdot (ir_1)\cdot...\cdot (ir_k);\ \chi\in K \}$ is compact in ${\mathcal C}(\widehat{G})$, hence
$$
\partial_{h_k}...\partial_{h_1}{\mathcal F}_G(\alpha_\nu)(\chi)=\alpha(\chi\cdot (ir_1)\cdot (ir_k))\underset{\nu\to \infty}{\longrightarrow}0
$$
uniformly by $\chi\in K$. This means that $\partial_{h_k}...\partial_{h_1}{\mathcal F}_G(\alpha_\nu)\to 0$ in the space ${\mathcal C}(\widehat{G})$. This is true for any $h_1,...,h_k$, therefore, ${\mathcal F}_G(\alpha_\nu)\to 0$ in the space ${\mathcal E}(\widehat{G})$. Similarly we consider the case of $\alpha\in{\mathcal C}^\star(G)$.
\epr

\btm\label{TH:Fourier-Lie=obolochka} For a compactly generated commutative Lie group $G$ both Fourier transforms
$$
{\mathcal F}_G:{\mathcal C}^\star(G)\to {\mathcal E}(\widehat{G}),\qquad {\mathcal F}_G:{\mathcal E}^\star(G)\to {\mathcal E}(\widehat{G}).
$$
are smooth envelopes of group algebras:
\beq\label{Env_E-C^star(G)-commut}
\Env_{\mathcal E} {\mathcal C}^\star(G)=\Env_{\mathcal E} {\mathcal E}^\star(G)={\mathcal E}(\widehat{G}).
\eeq
 \etm
\bpr By Theorem \ref{TH:Env_E-C^star(G)=Env_E-E^star(G)} it is sufficient to consider the case of ${\mathcal E}^\star(G)$. And for it we have just check the conditions (i) and (ii) of Theorem \ref{C^infty-obolochka-podalgebry-v-C^infty}. Let $\delta:G\to {\mathcal E}^\star(G)$ be the inclusion of the group  $G$ into its group algebra ${\mathcal E}^\star(G)$ as delta-functions:
 $$
\delta_a(u)=u(a),\qquad u\in {\mathcal E}(G).
 $$

Then, first, by \cite[Theorem 10.12]{Akbarov}, the formula
$$
\chi=\ph\circ\delta
$$
defines a bijection between the characters $\ph:{\mathcal E}^\star(G)\to\C$ on the algebra ${\mathcal E}^\star(G)$ and the complex characters $\chi:G\to\C^\times$ on the group $G$. Herewith the involutive characters $\ph:{\mathcal E}^\star(G)\to\C$ correspond to the usual characters $\chi:G\to\T$ (with values in the circle $\T$). This correspondence $\ph\leftrightarrow\chi$ is continuous in both directions, hence it defines an isomorphism
$$
\Spec{\mathcal E}^\star(G)\cong \widehat{G}.
$$
This is the condition (i) in Theorem \ref{C^infty-obolochka-podalgebry-v-C^infty}.

Second, let $\tau:{\mathcal E}^\star(G)\to\C$ be a tangent vector in the point $\e\in\Spec\Big({\mathcal E}^\star(G)\Big)$, that corresponds to the unit character $1(a)=1$, $a\in G$. For it the Leibnitz identity
\eqref{DEF:kasatelnyj-vektor} has the form
$$
\tau(\alpha *
\beta)=\tau(\alpha)\cdot\e(\beta)+\e(\alpha)\cdot\tau(\beta),\qquad
\alpha,\beta\in {\mathcal E}^\star(G),
$$
and if we replace $\alpha$ by $\delta_a$, and $\beta$ by $\delta_b$, we obtain
$$
\tau(\delta_{a\cdot b})= \tau(\delta_a *\delta_b)
=\tau(\delta_a)\cdot\e(\delta_b)+\e(\delta_a)\cdot\tau(\delta_b)=
\tau(\delta_a)+\tau(\delta_b),\qquad a,b\in G.
$$
This means that the mapping $a\mapsto\tau(\delta_a)$ is a homomorphism from
$G$ into the additive group of $\C$. If we claim in addition that $\tau$ is an involutive tangent vector, then the numbers $\tau(\delta_a)$ become real, and therefore the mapping $a\mapsto\tau(\delta_a)$ turns into a
(continuous) homomorphism of groups $G\to\R$, i.e. a real character. Thie establishes a bijection between the tangent space $T_{\e}\Big({\mathcal E}^\star(G)\Big)$ to the algebra ${\mathcal E}^\star(G)$ and the group $\Hom(G,\R)$ of real characters on $G$. But $\Hom(G,\R)$ is isomorphic to the group of one-parameter subgroups in $\widehat{G}$, i.e. to the group $\Hom(\R,\widehat{G})$ of (continuous) homomorphisms $\R\to \widehat{G}$ (see e.g.
\cite[(24.42)]{Hewitt-Ross}). As a corollary, $\Hom(G,\R)$ is isomorphic to the tangent space to the group  $\widehat{G}$ in the point $1\in \widehat{G}$, and we get
$$
T_{\e}\Big({\mathcal E}^\star(G)\Big)\cong
\Hom(G,\R)\cong\Hom(\R,\widehat{G})\cong T_1(\widehat{G}).
$$
If we omit the intermediate equations, then the point 1 can be replaces by any other point  $\chi\in \widehat{G}$, and we come to an isomorphism
$$
T_{\chi}\Big({\mathcal E}^\star(G)\Big)\cong T_{\chi}(\widehat{G})\cong
T_{\chi}\Big({\mathcal E}(\widehat{G})\Big).
$$
This is the condition (ii) of Theorem \ref{C^infty-obolochka-podalgebry-v-C^infty}.
 \epr

\paragraph{Smooth envelope of the group algebra of a compact group.}

\btm\label{PROP:glad-obol-komp-gruppy}
For each compact group $K$ the smooth envelope of its group algebra of measures ${\mathcal C}^\star(K)$ is the Cartesian product of the algebras ${\mathcal B}(X_\pi)$, where $\pi$ runs over the dual object $\widehat{K}$, and $X_\pi$ is the space of the representation $\pi$:
\beq\label{glad-obol-komp-gruppy}
\Env_{\mathcal E}{\mathcal C}^\star(K)=\prod_{\pi\in\widehat{K}}{\mathcal B}(X_\pi).
\eeq
If in addition $K$ is a Lie group, then the same algebra is the smooth envelope of the group algebra of distributions ${\mathcal E}^\star(K)$:
\beq\label{glad-obol-komp-gruppy-E^star(G)}
\Env_{\mathcal E}{\mathcal E}^\star(K)=\Env_{\mathcal E}{\mathcal C}^\star(K)=\prod_{\pi\in\widehat{K}}{\mathcal B}(X_\pi).
\eeq
\etm
We need the following

\blm\label{LM:D_k|_K=0} Suppose $K$ is a compsct group and $B$ a $C^*$-algebra. Then in each system of partial derivatives $D_k:{\mathcal C}^\star(G)\to B$, $k\in\N[m]$, only the operator $D_0$ can be non-zero:
\beq\label{D_k|_K=0}
\forall k>0\qquad D_k=0.
\eeq
\elm
\bpr Consider the mapping
$$
\ph_k=D_k\circ\delta:K\to B.
$$
We need to verify that
\beq\label{ph_k=0,k-ne-0}
\forall k\ne 0\qquad \ph_k=0
\eeq
For each $a\in K$ we have
$$
\ph_0(a)^{-1}=\ph_0(a^{-1})=D_0(\delta^{a^{-1}})=D_0((\delta^a)^\bullet)=D_0(\delta^a)^\bullet=\ph_0(a)^\bullet,
$$
i.e. $\ph_0(a)$ is a unitary element in $B$. Hence,
\beq\label{||ph_0(a)||=1}
||\ph_0(a)||=1,\qquad a\in K.
\eeq
Suppose now that $k$ is a multiindex of order 1. Then, first,
\beq\label{norm(ph_k(x))=C<infty}
\sup_{x\in K}\norm{\ph_k(x)}=C<\infty
\eeq
(since $\ph_k:K\to B$ is a continuous mapping on the compact space $K$).
And, second,
\beq\label{ph_k(a-cdot-b)=ph_0(a)-cdot-ph_k(b)+ph_k(a)-cdot-ph_0(b)}
\ph_0(a)\cdot\ph_k(b)=\ph_k(b)\cdot\ph_0(a),\qquad
\ph_k(a\cdot b)=\ph_0(a)\cdot\ph_k(b)+\ph_k(a)\cdot\ph_0(b),\qquad a,b\in K
\eeq
$$
\Downarrow
$$
$$
\ph_k(a^p)=p\cdot \ph_0(a)^{p-1}\cdot\ph_k(a),\qquad a\in K,\quad p\in\N.
$$
$$
\Downarrow
$$
$$
\ph_k(a)=\frac{1}{p}\cdot \ph_0(a^{-1})^{p-1}\cdot\ph_k(a^p),\qquad a\in K,\quad p\in\N.
$$
$$
\Downarrow
$$
$$
\norm{\ph_k(a)}=\frac{1}{p}\cdot\norm{\ph_0(a^{-1})^{p-1}\cdot\ph_k(a^p)}\le
\frac{1}{p}\cdot\underbrace{\norm{\ph_0(a^{-1})}^{p-1}}_{\scriptsize \begin{matrix}\eqref{||ph_0(a)||=1}\ \|\ \phantom{\eqref{||ph_0(a)||=1}} \\ 1\end{matrix}}\cdot\kern-10pt\underbrace{\norm{\ph_k(a^p)}}_{\scriptsize \begin{matrix} \eqref{norm(ph_k(x))=C<infty}\ \text{\rotatebox{90}{$\ge$}} \ \phantom{\eqref{norm(ph_k(x))=C<infty}} \\ C\end{matrix}}\kern-10pt
 \le \frac{1}{p}\cdot 1\cdot C\underset{p\to\infty}{\longrightarrow}0
$$
$$
\Downarrow
$$
$$
\ph_k(a)=0,\qquad a\in K.
$$

We proved \eqref{ph_k=0,k-ne-0} for multi-indices $k$ of order 1. If now $|k|=2$, then we have
\begin{multline*}
\ph_k(a\cdot b)=\sum_{0\le l\le k}\begin{pmatrix}k \\ l\end{pmatrix}\cdot \ph_{k-l}(a)\cdot \ph_l(b)
=\\=\ph_k(a)\cdot \ph_0(b)+\sum_{|l|=1}\begin{pmatrix}k \\ l\end{pmatrix}\cdot \ph_{k-l}(a)\cdot \underbrace{\ph_l(b)}_{\scriptsize\begin{matrix}\| \\ 0\end{matrix}}+
\ph_0(a)\cdot \ph_k(b)=\ph_k(a)\cdot \ph_0(b)+\ph_0(a)\cdot \ph_k(b)
\end{multline*}
I.e. $\ph_k$ satisfies the right identity in \eqref{ph_k(a-cdot-b)=ph_0(a)-cdot-ph_k(b)+ph_k(a)-cdot-ph_0(b)}
(and therefore, both identities), and by the same reasons $\ph_k=0$. In general case we have to organize induction by $k$. \epr

\bpr[Proof of Theorem \ref{PROP:glad-obol-komp-gruppy}]
It is sufficient to prove \eqref{glad-obol-komp-gruppy}, since \eqref{glad-obol-komp-gruppy-E^star(G)} will follow from Theorem \ref{TH:Env_E-C^star(G)=Env_E-E^star(G)}. The identity \eqref{D_k|_K=0} imply that, when we compute the smooth envelope of the algebra ${\mathcal C}^\star(G)$, the class of test morphisms coincides with the class of morphisms into $C^*$-algebras. As a corollary, the smooth envelope of the algebra ${\mathcal C}^\star(G)$ coincides with its continuous envelope, and thus, by Proposition \ref{PROP:nepr-obol-komp-gruppy},
$$
\Env_{\mathcal E}{\mathcal C}^\star(K)=\prod_{\pi\in\widehat{K}}{\mathcal B}(X_\pi).
$$
\epr

\paragraph{Smooth envelope of the group algebra of the group $C\times K$.}

\btm\label{TH:env_E-C(R^n-times-K)}
Suppose $C$ is an Abelian compactly generated Lie group, and $K$ a compact group. Then the smooth envelope of the group algebra of measures ${\mathcal C}^\star(C\times K)$ is the algebra ${\mathcal E}(\widehat{C},\prod_{\sigma\in\widehat{K}}{\mathcal B}(X_\sigma))$ of smooth mappings from the Pontryagin dual group $\widehat{C}$ to the Cartesian product of the algebras ${\mathcal B}(X_\sigma)$, where $\sigma$ runs over the dual object $\widehat{K}$, and $X_\sigma$ is the space of representation $\sigma$:
\beq\label{env_E-C(R^n-times-K)}
\Env_{\mathcal E}{\mathcal E}^\star(C\times K)={\mathcal E}\Big(\widehat{C},\prod_{\sigma\in\widehat{K}}{\mathcal B}(X_\sigma)\Big)=
\prod_{\sigma\in\widehat{K}}{\mathcal E}\big(\widehat{C},{\mathcal B}(X_\sigma)\big)={\mathcal E}(\widehat{C})\odot \prod_{\sigma\in\widehat{K}}{\mathcal B}(X_\sigma)=
\Env_{\mathcal E}{\mathcal C}^\star(C)\odot \Env_{\mathcal E}{\mathcal C}^\star(K).
\eeq
If in addition $K$ is a Lie group, then the same algebra is a smooth envelope of the algebra of distributions ${\mathcal E}^\star(C\times K)$:
\begin{multline}\label{env_E-C^star(R^n-times-K)}
\Env_{\mathcal E}{\mathcal C}^\star(C\times K)=\Env_{\mathcal C}{\mathcal E}^\star(C\times K)={\mathcal E}\Big(\widehat{C},\prod_{\sigma\in\widehat{K}}{\mathcal B}(X_\sigma)\Big)=\\=
\prod_{\sigma\in\widehat{K}}{\mathcal E}\big(\widehat{C},{\mathcal B}(X_\sigma)\big)={\mathcal E}(\widehat{C})\odot \prod_{\sigma\in\widehat{K}}{\mathcal B}(X_\sigma)=
\Env_{\mathcal E}{\mathcal E}^\star(C)\odot \Env_{\mathcal E}{\mathcal E}^\star(K).
\end{multline}
\etm
\bpr Both propositions follow from Theorems \ref{TH:Fourier-Lie=obolochka} and \ref{PROP:glad-obol-komp-gruppy}. For example, for ${\mathcal E}^\star(C\times K)$ the chain of reasoning is as follows:
\begin{multline*}
\Env_{\mathcal E}{\mathcal E}^\star(C\times K)=\Env_{\mathcal E}\Big({\mathcal E}^\star(C)\circledast{\mathcal E}^\star(K)\Big)=
\cite[(1.129)]{Akbarov-env}=
\Env_{\mathcal E}\Big(\Env_{\mathcal E}{\mathcal E}^\star(C)\circledast\Env_{\mathcal E}{\mathcal E}^\star(K)\Big)=\\=\eqref{Env_E-C^star(G)-commut},\eqref{glad-obol-komp-gruppy}=
\Env_{\mathcal E}\Big({\mathcal E}(\widehat{C})\circledast\prod_{\sigma\in\widehat{K}}{\mathcal B}(X_\sigma)\Big)=\eqref{E/circledast}=
{\mathcal E}(\widehat{C})\overset{\mathcal E}{\circledast}\prod_{\sigma\in\widehat{K}}{\mathcal B}(X_\sigma)=\eqref{E(M)-circledast-A=E(M,A)}=
{\mathcal E}\Big(\widehat{C},\prod_{\sigma\in\widehat{K}}{\mathcal B}(X_\sigma)\Big)
\end{multline*}
This proves the second equality in \eqref{env_E-C^star(R^n-times-K)}. Similarly we prove the first equality in \eqref{env_E-C(R^n-times-K)}, and together these equalities give the first equality in  \eqref{env_E-C^star(R^n-times-K)}:
$$
\Env_{\mathcal E}{\mathcal C}^\star(C\times K)={\mathcal E}\Big(\widehat{C},\prod_{\sigma\in\widehat{K}}{\mathcal B}(X_\sigma)\Big)=\Env_{\mathcal E}{\mathcal E}^\star(C\times K).
$$
The third equality in \eqref{env_E-C^star(R^n-times-K)} is obvious. The fourth follows from \cite[Theorem 8.9]{Akbarov}. Finally, the last equality in \eqref{env_E-C^star(R^n-times-K)} follows from \eqref{Env_E-C^star(G)-commut} and \eqref{glad-obol-komp-gruppy}.
\epr

\subsection{The algebra ${\mathcal K}_\infty(G)$}

For each Lie group $G$ its group algebra of distributions ${\mathcal E}^\star(G)$ is an involutive Hopf algebra with respect to the projective stereotype tensor product $\circledast$. Hence, by Theorems \ref{TH:C^infty-obolochka-sohranyaet-Hopfov} and \ref{TH:C^infty-obolochka-sohranyaet-inv-Hopfov}, its smooth envelope $\Env_{\mathcal E}({\mathcal E}^\star(G))$ is a coalgebra with interconsistent antipode and involution on the categories ${\tt E}\text{-}{\tt Alg}$ of smooth algebras and $({\tt Ste},\odot)$ of stereotype spaces. Let us denote by ${\mathcal K}_\infty(G)$ the stereotype dual space to the space $\Env_{\mathcal E}{\mathcal E}^\star(G)$:
 \beq\label{DEF:K_infty(G)}
{\mathcal K}_\infty(G):=\Big(\Env_{\mathcal E}{\mathcal E}^\star(G)\Big)^\star.
 \eeq
This is a dual space to a coalgebra in $({\tt Ste},\odot)$ with interconsistent antipode and involution, hence by property $4^\circ$ on page \pageref{4^0:inv-v-sopryazh-alg}, we have

\btm For each Lie group $G$ the space ${\mathcal K}_\infty(G)$ is an algebra in the category $({\tt Ste},\circledast)$ (i.e. a stereotype algebra) with interconsistent antipode and involution.
\etm

By Theorem \ref{TH:C^infty-obolochka-sohranyaet-inv-Hopfov}(ii) the morphism
$$
\big(\env_{\mathcal E}{{\mathcal E}^\star(G)}\big)^\star:{\mathcal K}_\infty(G)=\Big(\Env_{\mathcal E}{\mathcal E}^\star(G)\Big)^\star\to {\mathcal E}^\star(G)^\star={\mathcal E}(G),
$$
dual to the morphism of envelope, is an involutive homomorphism of algebras. The same reasoning as for the algebra ${\mathcal K}(G)$ on page \pageref{Ker-e_C-C^star(G)^star=0} show that the morphism $\big(\env_{\mathcal E}{{\mathcal E}^\star(G)}\big)^\star$ has zero kernel.

As a corollary, the algebra ${\mathcal K}_\infty(G)$ can be treated as an involutive subalgebra in ${\mathcal E}(G)$:

\btm\label{TH:K_infty(G)->E(G)} The mapping $u\mapsto u\circ \env_{\mathcal E}{{\mathcal E}^\star(G)}\circ\delta$
coincides with the mapping $\big(\env_{\mathcal E}{{\mathcal E}^\star(G)}\big)^\star$, dual to $\env_{\mathcal E}{{\mathcal E}^\star(G)}$:
\beq\label{K_infty(G)->E(G)}
\env_{\mathcal E}{{\mathcal E}^\star(G)}^\star(u)=u\circ \env_{\mathcal E}{{\mathcal E}^\star(G)}\circ\delta
\eeq
and injectively embeds ${\mathcal K}_\infty(G)$ into ${\mathcal E}(G)$ as an involutive subalgebra (hence the operations of summing, multiplication and involution on ${\mathcal K}_\infty(G)$ are pointwise).
\etm

The following proposition is analogous to Theorem \ref{TH:K(G)=lim}:

\btm\label{TH:K_infty(G)=lim}
The algebra ${\mathcal K}_\infty(G)$ as a stereotype space is a nodal coimage (in the category of stereotype spaces)
 \beq\label{K_infty(G)=lim}
{\mathcal K}_\infty(G)=\Coim_{\infty}\ph^\star
 \eeq
of the mapping $\ph^\star$, dual to the natural morphism of stereotype spaces
$$
\ph:{\mathcal E}^\star(G)\to \leftlim_{U}{\mathcal E}^\star(G)/U.
$$
where $U$ runs over the system of differential neighbourhoods of zero in ${\mathcal E}^\star(G)$.
\etm

Recall diagram \eqref{Env_E-C^star(G)-Env_E-E^star(G)} and let us complete it to the diagram
\beq\label{Env_E-C^star(G)-Env_E-E^star(G)-NEW}
 \xymatrix @R=4pc @C=4pc
 {
{\mathcal C}^\star(G)\ar[r]^{\lambda^\star}\ar@{..>}@/_8ex/[dd]_{\env_{\mathcal C}{\mathcal C}^\star(G)}\ar@{..>}[d]^{\env_{\mathcal E}{\mathcal C}^\star(G)} & {\mathcal E}^\star(G)\ar@{..>}[d]_{\env_{\mathcal E}{\mathcal E}^\star(G)}\ar@{..>}@/^8ex/[dd]^{\env_{\mathcal C}{\mathcal E}^\star(G)}
\\
\Env_{\mathcal E}{\mathcal C}^\star(G)\ar@{..>}[d]^{\zeta_{{\mathcal C}^\star(G)}}\ar@{=>}[r]^{\Env_{\mathcal C}(\lambda^\star)} & \Env_{\mathcal E}{\mathcal E}^\star(G)\ar@{..>}[d]_{\zeta_{{\mathcal E}^\star(G)}}
\\
\Env_{\mathcal C}{\mathcal C}^\star(G)\ar@{=>}[r]^{\Env_{\mathcal C}(\lambda^\star)} & \Env_{\mathcal C}{\mathcal E}^\star(G)
 }
\eeq
Here the perimeter is diagram \eqref{Env_C-C^star(G)-Env_C-E^star(G)}, and the triangles on each side are diagrams  \eqref{Env_infty->Env}. The upper ("continuous") arrow $\lambda^\star$ is a dense injection, the two horizontal  ("double") arrows are isomorphisms by \eqref{Env_E-C^star(G)=Env_E-E^star(G)} and \eqref{Env_C-C^star(G)=Env_C-E^star(G)}, and the rest ("dotted") arrows are dense morphisms by the definitions of envelopes.

The dual diagram is
\beq\label{Env_E-C^star(G)-Env_E-E^star(G)-NEW-*}
 \xymatrix @R=3pc @C=3pc
 {
& {\mathcal C}(G) & {\mathcal E}(G)\ar[l]^{\lambda} &
\\
& \Big(\Env_{\mathcal E}{\mathcal C}^\star(G)\Big)^\star\ar@{..>}[u] & \Big(\Env_{\mathcal E}{\mathcal E}^\star(G)\Big)^\star\ar@{=>}[l]\ar@{=}[r]\ar@{..>}[u] & {\mathcal K}_\infty(G)
\\
{\mathcal K}(G)\ar@{=}[r] & \Big(\Env_{\mathcal C}{\mathcal C}^\star(G)\Big)^\star\ar@{..>}[u] & \Big(\Env_{\mathcal C}{\mathcal E}^\star(G)\Big)^\star \ar@{=>}[l]\ar@{..>}[u] &
 }
\eeq
and the upper arrow $\lambda$ is a dense injection, two horizontal ("double") arrows are isomorphisms, and the rest  ("dotted") arrows are injections. Thus, we can conclude that the following chain of injections hold:
$$
{\mathcal K}(G)\subseteq {\mathcal K}_\infty(G)\subseteq {\mathcal E}(G)\subseteq {\mathcal C}(G).
$$
It can be completed by the chain \eqref{V-k-K-C}, and we obtain the following

\btm\label{TH:V-k-K-K_infty-E-C}
For a real Lie group $G$ the following chain of set-theoretic inclusions hold:
 \beq\label{V-k-K-K_infty-E-C}
\Trig(G)\subseteq k(G)\subseteq {\mathcal K}(G)\subseteq {\mathcal K}_\infty(G)\subseteq {\mathcal E}(G)\subseteq {\mathcal C}(G),
 \eeq
And
\bit{

\item[(i)] always
$$
\overline{\Trig(G)}={\mathcal K}(G),\qquad \overline{{\mathcal E}(G)}={\mathcal C}(G)
$$

\item[(ii)] if $G=C\times K$, where $C$ is an abelian compactly generated Lie group, and $K$ a compact Lie group, then
\beq\label{overline-K(G)=K_infty(G)}
\overline{{\mathcal K}(G)}={\mathcal K}_\infty(G),
\eeq

\item[(iii)] if $G$ is a SIN-group, then
\beq\label{overline-K(G)=E(G)}
\overline{{\mathcal K}(G)}={\mathcal E}(G).
\eeq
}\eit
\etm

We need the following

\blm\label{LM:T_a[K(G)]}
If $G$ is a SIN-group, which is at the same time a Lie group, then the tangent space to the algebra ${\mathcal K}(G)$ at any point $a\in G$ coincides with the tangent spec to $G$ at this point:
\beq\label{T_a[K(G)]=T_a(G)}
T_a[{\mathcal K}(G)]=T_a(G).
\eeq
\elm
\bpr
1. First let $G$ be an Abeilan group. Since by Theorem \ref{LM:sdvig-v-K(G)}, ${\mathcal K}(G)$ is invarint with respect to shifts, we can take $a=0\in G$. Then
\begin{multline*}
T_a[{\mathcal K}(G)]=T_0\Big[\big(\Env_{\mathcal C}{\mathcal C}^\star(G)\big)^\star\Big]=\eqref{Env_C-C^star(G)=C(widehat(G))}=
T_0\Big[\big({\mathcal C}(\widehat{G})\big)^\star\Big]=\\=T_0\Big[{\mathcal C}^\star(\widehat{G})\Big]=\Hom(\widehat{G},\R)=\Hom(\R,\widehat{\widehat{G}})=\Hom(\R,G)=T_0(G)=T_a(G).
\end{multline*}

2. Further, let $G$ be a compact group. Then
$$
T_a[{\mathcal K}(G)]=\eqref{k=K}=T_a[k(G)]=\eqref{T_a[k(G)]=T_a(G)}=T_a(G).
$$

3. Take $G=\R^n\times K$. For any $s\in\R^n$ and $t\in K$ we have
\begin{multline*}
T_{s,t}[{\mathcal K}(\R^n\times K)]=\eqref{K(A-times-K)-cong-K(A)-circledast-K(K)}=T_{s,t}\Big[{\mathcal K}(\R^n)\circledast {\mathcal K}(K)\Big]=\eqref{T_s[A]-oplus-T_t[B]}=\\=
T_s[{\mathcal K}(\R^n)]\oplus T_t[{\mathcal K}(K)]=
T_s(\R^n)\times T_t(K)=T_{s,t}(\R^n\times K).
\end{multline*}

4. Finally, let $G$ be an arbitrary SIN-group (and a Lie group). Let us represent $G$ as a discrete extension \eqref{SIN-kak-rasshirenie} of some group $N=\R^n\times K$. Take a point $a\in G$, and choose a coset $L$ with respect to the subgroup $N$ that contains $a$, and consider the algebra ${\mathcal K}_L(G)$ defined in \eqref{DEF:K_L(G)}. If $e$ is the unit of the group $G$, then
\begin{multline*}
T_a[{\mathcal K}(G)]=T_a[{\mathcal K}_L(G)]=(\text{Theorem \ref{LM:sdvig-v-K(G)}})=
T_e[{\mathcal K}_N(G)]=(\text{Lemma \ref{LM:K(N)=K_N(G)}})=\\=T_e[{\mathcal K}(N)]=(\text{already proven})=T_e(N)=T_e(G)=T_a(G)
\end{multline*}
(the equalities mean isomorphisms in the obvious transformations).
\epr

\bpr[Proof of Theorem \ref{TH:V-k-K-K_infty-E-C}.]
1. The first formula in (i) is already proven in Theorem \ref{TH:V-k-K-C}, and the second is the standard relation between the space of smooth and continuous functions on a smooth manifold.

2. If $G=C\times K$, then, on the one hand,
$$
\Env_{\mathcal C}{\mathcal C}^\star(G)=\Env_{\mathcal C}{\mathcal C}^\star(C\times K)=
\eqref{env_C^star(R^n-times-K)}={\mathcal C}\Big(\widehat{C},\prod_{\sigma\in\widehat{K}}{\mathcal B}(X_\sigma)\Big)
$$
and on the other,
$$
\Env_{\mathcal E}{\mathcal E}^\star(G)=\Env_{\mathcal E}{\mathcal E}^\star(C\times K)=\eqref{env_E-C^star(R^n-times-K)}={\mathcal E}\Big(\widehat{C},\prod_{\sigma\in\widehat{K}}{\mathcal B}(X_\sigma)\Big)
$$
and thus the second space is embedded injectively into the first one. This implies that the dual spaces are densely mapped one into another:
$$
\overline{{\mathcal K}(G)}={\mathcal K}_\infty(G).
$$

3. Suppose $G$ is a SIN-group. Let us represent it as a discrete extension \eqref{SIN-kak-rasshirenie} of some group $N=\R^n\times K$. By Lemma \ref{LM:K(N)=K_N(G)}, the restriction of the space ${\mathcal K}(G)$ on $N$ is isomorphic to ${\mathcal K}(N)$, and by Lemma \ref{LM:Spec(R^n-times-K)}, the spectrum of the algebra ${\mathcal K}(N)$ coincides with $N$:
$$
\Spec{\mathcal K}(N)=N.
$$
We can conclude that the algebra ${\mathcal K}(G)$ separates the points of $N$. By Theorem \ref{LM:sdvig-v-K(G)}, the shifts are isomorphisms of ${\mathcal K}(G)$, hence ${\mathcal K}(G)$ separates the points of each coset $L\in G/N$. Besides this by Lemma \ref{LM:1_L-in-K(G)}, the characteristic function $1_L$ of each such a class $L$ belongs to ${\mathcal K}(G)$. This implies that ${\mathcal K}(G)$ separates the points not only inside each coset $L$, but also in different cosets $L,M\in G/N$. Thus ${\mathcal K}(G)$ separates points on the whole group $G$.

On the other hand, \eqref{T_a[K(G)]=T_a(G)} holds. Together this means that ${\mathcal K}(G)$ satisfies the conditions of the Nahbin theorem \ref{TH:Nachbin}, hence the algebra ${\mathcal K}(G)$ is dense in ${\mathcal E}(G)$.
\epr

\paragraph{The mapping ${\mathcal K}_\infty(G)\circledast{\mathcal K}_\infty(H)\to {\mathcal K}_\infty(G\times H)$.}

Let $G$ and $H$ be Lie groups. By analogy with the mapping $\omega_{G,H}:{\mathcal K}(G)\circledast{\mathcal K}(H)\to {\mathcal K}(G\times H)$ (defined on page \pageref{omega_G,H}) we define the mapping
${\mathcal K}_\infty(G)\circledast{\mathcal K}_\infty(H)\to {\mathcal K}_\infty(G\times H)$.
And by analogy with Theorem \ref{TH:K(A-times-K)-cong-K(A)-circledast-K(K)} we prove

\btm\label{TH:K_infty(A-times-K)-cong-K_infty(A)-circledast-K_infty(K)}
If $C$ is an Abelian compactly generated locally compact group, and $K$ a compact group, then the mapping ${\mathcal K}_\infty(C)\circledast {\mathcal K}_\infty(K)\to {\mathcal K}_\infty(C\times K)$ is an isomorphism:
\beq\label{K_infty(C-times-K)-cong-K_infty(A)-circledast-K_infty(K)}
{\mathcal K}_\infty(C\times K)\cong{\mathcal K}_\infty(C)\circledast {\mathcal K}_\infty(K)
\cong{\mathcal K}_\infty(C)\circledast {\mathcal K}(K)
\eeq
\etm
\bpr
This follows from \eqref{env_E-C^star(R^n-times-K)}:
\begin{multline*}
{\mathcal K}_\infty(C\times K)=\Big(\Env_{\mathcal E}{\mathcal E}^\star(C\times K)\Big)^\star\cong \eqref{env_E-C^star(R^n-times-K)}\cong \Big(\Env_{\mathcal E}{\mathcal E}^\star(C)\odot\Env_{\mathcal E}{\mathcal E}^\star(K) \Big)^\star\cong\\ \cong \Big(\Env_{\mathcal E}{\mathcal E}^\star(C)\Big)^\star\circledast\Big(\Env_{\mathcal E}{\mathcal E}^\star(K)\Big)^\star\cong{\mathcal K}_\infty(A
)\circledast {\mathcal K}_\infty(K)
\end{multline*}
\epr

\subsection{Smooth duality for groups $C\times K$}

\paragraph{Smooth envelope of the algebra ${\mathcal K}_\infty(C\times K)$.}

\btm\label{TH:E_E(K_infty(G))=E(G)}
Let $C$ be an Abelian Lie group, $K$ a compact Lie group, and $G=C\times K$. Then
\bit{

\item[(i)] the spectrum of the algebra ${\mathcal K}_\infty(G)$ is topologically isomorphic to $G$:
\beq\label{Spec-K_infty(G)=G}
\Spec{\mathcal K}_\infty(G)=G
\eeq

\item[(ii)] in each point $a\in G$ the tangent spaces of ${\mathcal K}_\infty(G)$ and $G$ are isomorphic:
\beq\label{T_a[K_infty(G)]=T_a(G)}
T_a[{\mathcal K}_\infty(G)]=T_a(G)
\eeq

\item[(iii)] the smooth envelope of the algebra ${\mathcal K}_\infty(G)$ coincides with ${\mathcal E}(G)$:
\beq\label{Env_E-K_infty(G)=E(G)}
\Env_{\mathcal E}{\mathcal K}_\infty(G)={\mathcal E}(G).
\eeq

}\eit
\etm
\bpr
1. From the chain \eqref{V-k-K-K_infty-E-C} let us extract the fragment
\beq\label{K(G)->K_infty(G)->E(G)}
{\mathcal K}(G)\subseteq {\mathcal K}_\infty(G)\subseteq {\mathcal E}(G).
\eeq
From \eqref{overline-K(G)=K_infty(G)} and \eqref{overline-K(G)=E(G)} it follows that this is a chain of (continuous and) dense injections. Thus after passing to spectra we obtain a chain of (continuous) injections
$$
G=\Spec{\mathcal K}(G)\gets\Spec{\mathcal K}_\infty(G)\gets\Spec{\mathcal E}(G)=G
$$
(here the first equality is proved in Theorem \ref{TH:Spec-K(G)=G}, and the second one is obvious). Certainly, these are homeomorphisms.

2. Since the injections in \eqref{K(G)->K_infty(G)->E(G)} are dense, the chain of mappings of tangent spaces also consists of injections
$$
T_a(G)=T_a[{\mathcal K}(G)]\gets T_a[{\mathcal K}_\infty(G)]\gets T_a[{\mathcal E}(G)]=T_a(G),
$$
(here the first equality follows from Corollary \ref{COR:kasat-prostr-k-Trig(G)}, and the last one is again obvious). And again it is clear that these mappings are isomorphisms (of finite dimensional vector spaces).

3. We proved (i) and (ii), and by Theorem \ref{C^infty-obolochka-podalgebry-v-C^infty} this implies (iii).
\epr

\paragraph{Structure of Hopf algebras on $\Env_{\mathcal E} {\mathcal E}^\star(G)$ and ${\mathcal K}_\infty(G)$ for  $G=C\times K$.}

\btm\label{TH:Env_E(E*(G))-Hopf-dlya-R^n-times-K} Suppose $C$ is an Abealin compactly generated Lie group, $K$ a compact Lie group, and $G=C\times K$. Then
\bit{
\item[(i)] the smooth envelope $\Env_{\mathcal E}{\mathcal E}^\star(G)$ of a group algebra ${\mathcal E}^\star(G)$ is an involutive Hopf algebra in the category of stereotype spaces $(\tt{Ste},\odot)$.

\item[(ii)] the dual algebra ${\mathcal K}_\infty(G)$ is an involutive Hopf algebra in the category of stereotype spaces $(\tt{Ste},\circledast)$
}\eit
\etm
\bpr It is sufficient to prove (i). First,
$$
\Env_{\mathcal E}{\mathcal E}^\star(C)=\eqref{Env_E-C^star(G)-commut}={\mathcal E}(\widehat{C}),
$$
and this is a Hopf algebra in the category $(\tt{Ste},\odot)$, for example, by \cite[Example 10.25]{Akbarov}. On the other hand,
$$
\Env_{\mathcal E}{\mathcal E}^\star(K)=\eqref{glad-obol-komp-gruppy-E^star(G)}=\prod_{\pi\in\widehat{K}}{\mathcal B}(X_\pi)=\eqref{nepr-obol-komp-gruppy}=\Env_{\mathcal C}{\mathcal C}^\star(K),
$$
and this is a Hopf algebra in the category $(\tt{Ste},\odot)$ by Theorem \ref{TH:E(H)-dlya-Moore}. Thus the space
\beq\label{Env_E(E*(G))-Hopf-dlya-R^n-times-K}
\Env_{\mathcal E}{\mathcal E}^\star(C\times K)=\eqref{env_E-C^star(R^n-times-K)}=
\Env_{\mathcal E}{\mathcal E}^\star(C)\odot\Env_{\mathcal E}{\mathcal E}^\star(K).
\eeq
is a Hopf algebra in $(\tt{Ste},\odot)$ as a tensor product of Hopf algebras.
\epr

\paragraph{Smoothly reflexive Hopf algebras.}
Let $H$ be an involutive stereotype Hopf algebra with respect to the tensor product $\circledast$. We say that $H$ is {\it smoothly reflexive}, if it is reflexive with respect to the smooth envelope $\Env_{\mathcal E}$ (in the sense of definition on page \pageref{DEF:reflexiv-otn-obolochki}).

Theorems \ref{TH:E_E(K_infty(G))=E(G)} and \ref{TH:Env_E(E*(G))-Hopf-dlya-R^n-times-K} imply

 \btm\label{TH:glad-dvoistvennost}
Let $C$ be an Abelian compactly generated Lie group, $K$ a compact Lie group, and $G=C\times K$. Then the algebras ${\mathcal E}^\star(G)$ and ${\mathcal K}_\infty(G)$ are smoothly reflexive, and the reflexivity diagram for them is:
 \beq\label{chetyrehugolnik-E-E*}
 \xymatrix @R=1.pc @C=2.pc
 {
 {\mathcal E}^\star(G)
 & \ar@{|->}[r]^{\Env_{\mathcal E}} & &
 \Env_{\mathcal E}{\mathcal E}^\star(G)
 \\
 & & &
 \ar@{|->}[d]^{\star}
 \\
 \ar@{|->}[u]^{\star}
 & & &
 \\
 {\mathcal E}(G)
 & &
 \ar@{|->}[l]_{\Env_{\mathcal E}}
 &
 {\mathcal K}_\infty(G)
 }
 \eeq
 \etm

\bex Theorem \ref{TH:Fourier-Lie=obolochka} implies that for a compactly generated Lie group $C$ the reflexivity diagram is
 \beq\label{chetyrehugolnik-E-E*-F}
  \xymatrix @R=1.pc @C=2.pc
 {
 {\mathcal E}^\star(C)
 & \ar@{|->}[r]^{{\mathcal F}_C} & &
 {\mathcal E}(\widehat{C})
 \\
 & & &
 \ar@{|->}[d]^{\star}
 \\
 \ar@{|->}[u]^{\star}
 & & &
 \\
 {\mathcal E}(C)
 & &
 \ar@{|->}[l]_{{\mathcal F}_{\widehat{C}}}
 &
 {\mathcal E}^\star(\widehat{C})
  }
\eeq
(here $\widehat{C}$ is the Pontryagin dual group for $C$, ${\mathcal F}_C$ the Fourier transform, defined in \eqref{Fourier-transform}).
\eex

\bex Theorem \ref{PROP:glad-obol-komp-gruppy} implies that for a compact Lie group $K$ the reflexivity diagram is
$$
 \xymatrix @R=1.pc @C=2.pc
 {
 {\mathcal E}^\star(K)
 & \ar@{|->}[r]^{\Env_{\mathcal E}} & &
\prod_{\pi\in\widehat{K}}{\mathcal B}(X_\pi)
 \\
 & & &
 \ar@{|->}[d]^{\star}
 \\
 \ar@{|->}[u]^{\star}
 & & &
 \\
 {\mathcal E}(K)
 & &
 \ar@{|->}[l]_{\Env_{\mathcal E}}
 &
 \Trig(K)
 }
$$
\eex

\paragraph{Groups, discerned by $C^*$-algebras with joined self-adjoint nilpotent elements.} By analogy with the definitions on page \pageref{DEF:gruppy-razlich-C^*-algebrami}, let us say that a locally compact group $G$ {\it is discerned by $C^*$-algebras with joined self-adjoint nilpotent elements}\label{DEF:gruppy-razlich-C^*-algebrami-s-nilp},
if (continuous involutive) homomorphisms of its measure algebra ${\mathcal C}^\star(G)\to B[m]$ into various algebras of the form $B[m]$, where $B$ is a $C^*$-algebra, and $m\in\N[n]$, separate elements of $G$ (with the injection of  $G$ into ${\mathcal C}^\star(G)$ by delta-functions).

\btm\label{TH:gruppy-razlich-C*-s-pris-nilp}
If a Lie group $G$ is discerned by $C^*$-algebras with joined self-adjoint nilpotent elements, then $G$ is discerned by $C^*$-algebras (and therefore is a linear group by Theorem \ref{TH:Luminet-Valette}).
\etm

To prove this we need

\blm\label{LM:x->e^x}
For each $C^*$-algebra $A$ the exponential mapping
$$
x\mapsto e^x=\sum_{n=0}^\infty\frac{x^n}{n!}\quad : A\to A
$$
is injective on self-adjoint elements:
$$
x\ne y\in\Real A\quad\Longrightarrow\quad e^x\ne e^y.
$$
\elm
\bpr
Let $\Real_+ A$ be the set of positive self-adjoint elements in $A$.
For each $z\in\Real_+ A$ and for each $n\in\N$ a square $\sqrt[n]{z}\in\Real_+ A$ is uniquely defined (this follows from the spectral theorem \cite{Kad-Ring}). On the other hand, for each $x\in\Real A$ its exponent $e^x$ belongs to $\Real_+A$, since
$$
e^x=e^{\frac{x}{2}}\cdot e^{\frac{x}{2}}=e^{\frac{x}{2}}\cdot\left( e^{\frac{x}{2}}\right)^\bullet\ge 0.
$$
Hence, the roots $\sqrt[n]{e^x}=e^{\frac{x}{n}}$ are uniquely defined, and we obtain
$$
x=\lim_{n\to\infty}n\left(e^{\frac{x}{n}}-1\right)=\lim_{n\to\infty}n\left(\sqrt[n]{e^x}-1\right),
$$
i.e. $x$ is uniquely defined by $e^x$.
\epr

\bpr[Proof of Theorem \ref{TH:gruppy-razlich-C*-s-pris-nilp}.]
Suppose that $G$ is not discerned by $C^*$-algebras, i.e. all homomorphisms $\ph:G\to B$ into various $C^*$-algebras $B$ have a common non-trivial kernel $N$, $\{e\}\ne N\subseteq G$. Consider a homomorphism $D:G\to B[m]$. For each multiindex $k\in\N[m]$ of first order, $|k|=1$, and for each $x,y\in N$ we have:
$$
D_k(x\cdot y)=C_k(x)\cdot\underbrace{D_0(y)}_{\scriptsize\begin{matrix}\|\\ 1\end{matrix}}+\underbrace{D_0(x)}_{\scriptsize\begin{matrix}\|\\ 1\end{matrix}}\cdot D_k(y)=D_k(x)+D_k(y).
$$
Hence, $D_k:N\to B$ is a logarithm (i.e. it turns multiplication in $N$ into summing in $B$). Besides this,
$$
D_k(x)^\bullet=D_k(x^\bullet)=D_k(x^{-1})=-D_k(x),
$$
and this means that each element $iD_k(x)$, $x\in N$, is self-adjoint:
$$
(iD_k(x))^\bullet=iD_k(x).
$$
Consider the mapping
$$
\ph(x)=e^{iD_k(x)}=\sum_{l=0}^\infty\frac{i^{|l|}}{l!}\cdot D_k(x)^l.
$$
It is an involutive homomorphism, i.e. a representation of $N$ in the $C^*$-algebra $B$. One can consider its induced representation $\psi:G\to{\mathcal B}(Z)$ in a Hilbert space $Z$. This is a homomorphism of $G$ into a $C^*$-algebra ${\mathcal B}(Z)$, hence on the subgroup $N$ the mapping $\psi$ must be trivial:
$$
\psi(x)=1,\qquad x\in N.
$$
This means that the initial homomorphism $\ph:N\to B$ must be trivial as well:
$$
\ph(x)=e^{iD_k(x)}=1_B.
$$
Hence each element $iD_k(x)$ is self-adjoint. By Lemma \ref{LM:x->e^x} this means that
$$
D_k(x)=0,\qquad x\in N.
$$
We see that on the subgroup $N$ all partial derivatives of order 1 vanish. If now $k$ has order 2, then for $x,y\in N$ we have
\begin{multline*}
D_k(x\cdot y)=\sum_{0\le l\le k}\begin{pmatrix}k \\ l\end{pmatrix}\cdot D_{k-l}(x)\cdot D_l(y)=\\=
\overbrace{D_k(x)\cdot \underbrace{D_0(y)}_{\scriptsize\begin{matrix}\|\\ 1\end{matrix}}}^{\scriptsize |l|=0}+\sum_{|l|=1}\begin{pmatrix}k \\ l\end{pmatrix}\cdot D_{k-l}(x)\cdot \underbrace{D_l(y)}_{\scriptsize\begin{matrix}\|\\ 0\end{matrix}}+\overbrace{\underbrace{D_0(x)}_{\scriptsize\begin{matrix}\|\\ 1\end{matrix}}\cdot D_k(y)}^{\scriptsize |l|=2}=D_k(x)+D_k(y)
\end{multline*}
I.e. $D_k$ is again a logarithm. By the same reason, $D_k=0$ on $N$. And so on.
\epr

\newpage

\renewcommand{\thesection}{}

\section{Errata}\label{Errata}

After sending this paper to the journal the author found that Section 5 contains a mistake in the proof of one statement, and this implies several gaps in the proofs of some next statements in this section. The author apologizes for this fault. Namely, in the proof of Property $1^\circ$ on page \pageref{1^0:prodolzhenie-polunorm} the reasoning that the induced representation preserves the norm is false. The author thanks Fan Zheng for the following counterexample.

\bex\label{EX:Fan-Zheng}
Let $G=S_3$ be the group of permutations of 3 elements. Let us set $x=(2\ 3\ 1)$ and $y=(2\ 1\ 3)$. Then
$$
y\cdot x\cdot y^{-1}=x^{-1}=x^2.
$$
Consider the subgroup $N=A_3$ of even permutations in $G=S_3$. (The group $N=A_3$ consists of three permutations, $(1\ 2\ 3)$, $(2\ 3\ 1)$ and $(3\ 1\ 2)$, and obviously, $x\in A_3$, and $y\notin A_3$.) Let $\pi:N\to\C$ be the action of $N$ on $\C$, defined by the rule
$$
\pi(1)\lambda=\lambda,\quad \pi(x)\lambda=e^{\frac{2\pi i}{3}}\cdot\lambda,\quad \pi(x^2)\lambda=e^{-\frac{2\pi i}{3}}\cdot\lambda,\qquad\lambda\in\C.
$$

Consider the element of the group algebra ${\mathcal C}^\star(G)$
$$
\alpha=\delta^1+i\cdot \delta^x.
$$
The extension $\dot{\pi}:{\mathcal C}^\star(N)\to\C$ of the representation $\pi$ to the group algebra ${\mathcal C}^\star(N)$ turns $\alpha$ into the number
$$
\dot{\pi}(\alpha)=1+i\cdot e^{\frac{2\pi i}{3}}
$$
with the absolute value
$$
\norm{\dot{\pi}(\alpha)}=\abs{1+i\cdot e^{\frac{2\pi i}{3}}}=\sqrt{2-\sqrt{3}}.
$$

Consider the quotient group $F=G/N$, and the quotient map $\ph:G\to G/N=F$. Choose an arbitrary retraction $\sigma:F\to G$ and put $z=\ph(y)\in F$. Then consider the induced representation $\pi':G\to{\mathcal B}(L_2(F))$ and its extension to the group algebra $\dot{\pi}':{\mathcal C}^\star(G)\to{\mathcal B}(L_2(F))$. Let $\xi\in L_2(F)$ be a function defined by the rule
$$
\xi(t)=\begin{cases}1,& t=z\\ 0,& t=1\end{cases}.
$$
Then
$$
\pi'(x)(\xi)(z)=\pi(\sigma(z)\cdot x\cdot\sigma(z)^{-1})(\xi(z))=\pi(y\cdot x\cdot y^{-1})(1)=\pi(x^2)\cdot 1=
e^{-\frac{2\pi i}{3}},
$$
and
$$
\norm{\dot{\pi}'(\alpha)}\ge\norm{\dot{\pi}'(\alpha)(\xi)}=\abs{\dot{\pi}'(\alpha)(\xi)(z)}=\abs{1+i\cdot e^{-\frac{2\pi i}{3}}}=\sqrt{2+\sqrt{3}}>
\sqrt{2-\sqrt{3}}=\norm{\dot{\pi}(\alpha)}.
$$
\eex

Property $1^\circ$ on page \pageref{1^0:prodolzhenie-polunorm} can be proved in a narrower case when the subgroup $N$ is a direct summand in the group $G$:

\bprop\label{PROP:prodolzhenie-polunorm} Suppose $G=N\times D$, where $N$ is a locally compact group, and $D$ a discrete group. Then each seminorm $p\in{\mathcal P}(N)$ can be extended to a seminorm $q\in{\mathcal P}(G)$.
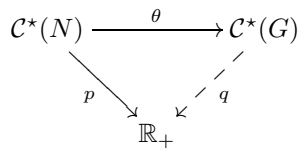
$$
 \xymatrix @R=2pc @C=1.2pc
 {
  {\mathcal C}^\star(N)\ar[rr]^{\theta}\ar[dr]_{p} & & {\mathcal C}^\star(G)\ar@{-->}[dl]^{q} \\
  & \R_+ &
 }
$$
\eprop
\bpr
Let $p:{\mathcal C}^\star(N)\to\R_+$ be a continuous $C^*$-seminorm. It is a norm of some norm-continuous representation $\dot{\pi}:{\mathcal C}^\star(N)\to{\mathcal L}(X)$,
$$
p(\alpha)=\norm{\dot{\pi}(\alpha)}.
$$
By Theorem \ref{TH:Kuz-1}, the induced representation $\dot{\pi}':{\mathcal C}^\star(G)\to{\mathcal L}(L_2(D,X))$ is also norm-continuous. Hence, the norm
$$
q(\beta)=\norm{\dot{\pi}'(\beta)},\qquad \beta\in{\mathcal C}^\star(G),
$$
is continuous on ${\mathcal C}^\star(G)$. Let us look at the construction of the induced representation \eqref{ind-representation}. In our case the quotient map $\ph:G\to D$ is just a projection to the second component:
$$
\ph(a,t)=t,\qquad a\in N,\ t\in D.
$$
Let $\rho(t)$ be the component of the element $\sigma(t)\in G$ in the group $N$:
$$
\sigma(t)=(\rho(t),t),\qquad \rho(t)\in N,\ t\in D.
$$
Then for $a\in N$, $t\in D$ we have:
\begin{multline*}
\pi'(a,1)(\xi)(t)=\pi(\sigma(t)\cdot (a,1)\cdot \sigma(t\cdot\underbrace{\ph(a,1)}_{\scriptsize\begin{matrix}\|\\ 1\end{matrix}})^{-1})\Big(\xi(t\cdot\underbrace{\ph(a,1)}_{\scriptsize\begin{matrix}\|\\ 1\end{matrix}})\Big)=
\pi(\sigma(t)\cdot (a,1)\cdot \sigma(t)^{-1})(\xi(t))=\\=
\pi\big((\rho(t),t)\cdot (a,1)\cdot (\rho(t),t)^{-1}\big)(\xi(t))=
\pi\big((\rho(t),t)\cdot (a,1)\cdot (\rho(t)^{-1},t^{-1})\big)(\xi(t))=\\=
\pi\big((\rho(t)\cdot a\cdot \rho(t)^{-1},t\cdot t^{-1})\big)(\xi(t))=
\pi\big(\underbrace{(\rho(t)\cdot a\cdot \rho(t)^{-1},1)}_{\scriptsize\begin{matrix}\text{\rotatebox{90}{$\owns$}}\\ N\times{1} \end{matrix}}\big)(\xi(t))=\pi(\rho(t)\cdot a\cdot \rho(t)^{-1})(\xi(t))=\\=
\Big(\pi(\rho(t))\circ \pi(a)\circ \pi(\rho(t))^{-1}\Big)(\xi(t))
\end{multline*}
From this for $\alpha\in {\mathcal C}^\star(N)$ we have:
$$
\dot{\pi}'(\alpha)(\xi)(t)=\Big(\pi(\rho(t))\circ \dot{\pi}(\alpha)\circ\pi(\rho(t))^{-1}\Big)(\xi(t))
$$
and, taking into account that $\pi(\rho(t))$ is a unitary operator,
\begin{multline*}
\norm{\dot{\pi}'(\alpha)}^2=\sup_{\norm{\xi}\le 1}\norm{\dot{\pi}'(\alpha)(\xi)}^2=\sup_{\norm{\xi}\le 1}\sum_{t\in D}\norm{\dot{\pi}'(\alpha)(\xi)(t)}^2=\\=\sup_{\norm{\xi}\le 1}
\sum_{t\in D}\norm{\Big(\pi(\rho(t))\circ \dot{\pi}(\alpha)\circ\pi(\rho(t))^*\Big)(\xi(t))}^2=
\sup_{\norm{\xi}\le 1}
\sum_{t\in D}\norm{\dot{\pi}(\alpha)(\xi(t))}^2=\norm{\dot{\pi}(\alpha)}^2
\end{multline*}
\epr

An immediate corollary from these observations is that the most part of the propositions of Section 5 related to the envelopes of group algebras of SIN-groups can be treated as proven statements in this text only for the special case when the group $G$ is a cartesian product of an Euclidean space $\R^n$, a compact group and a discrete group. In detail,
 \bit{
 \item[1)] Properties $1^\circ$-$4^\circ$ on page \pageref{1^0:prodolzhenie-polunorm}, Lemma  \ref{LM:p_(T,S)-konf-sistema}, Proposition \ref{PROP:nepr-Env-SIN-gruppy}, Proposition  \ref{PROP:nepr-Env-SIN-gruppy=LCS-lim}, Proposition \ref{PROP:E(theta)}, Theorem \ref{TH:V-k-K-C}(ii), Lemma  \ref{LM:K(N)=K_N(G)}(ii) -- are proved for the case of $G=\R^n\times K\times D$, where $K$ is a compact group, and $D$ a discrete group,
 \item[2)] Theorem \ref{TH:Spec-K(G)=G}, Theorem \ref{TH:E(K(G))=C(G)}, Theorem \ref{TH:E(H)-dlya-Moore}, Lemma  \ref{LM:E(H)-dlya-Moore-ii}, Theorem \ref{TH:nepr-dvoistvennost} -- are proved for the case when the Moore group $G$ can be represented in the form $G=\R^n\times K\times D$, where $K$ is a compact group, and $D$ a discrete Moore group.
 }\eit

\tableofcontents

\end{document}